\documentclass{article}

\usepackage{longtable}
\usepackage{booktabs}
\def\tightlist{}

\usepackage{charter}
\usepackage{fancyhdr}
\usepackage[margin=1.6in]{geometry}

\usepackage{amsmath,amssymb}
\usepackage[colorlinks=true]{hyperref}
\usepackage{xcolor}
\definecolor{myurlcolor}{rgb}{0.6,0,0}
\definecolor{mycitecolor}{rgb}{0,0,0.8}
\definecolor{myrefcolor}{rgb}{0,0,0.8}
\hypersetup{linkcolor=myrefcolor,citecolor=mycitecolor,urlcolor=myurlcolor}

\renewcommand{\texttt}[1]{%
  \begingroup
  \ttfamily
  \begingroup\lccode`~=`/\lowercase{\endgroup\def~}{/\discretionary{}{}{}}%
  \begingroup\lccode`~=`[\lowercase{\endgroup\def~}{[\discretionary{}{}{}}%
  \begingroup\lccode`~=`.\lowercase{\endgroup\def~}{.\discretionary{}{}{}}%
  \catcode`/=\active\catcode`[=\active\catcode`.=\active
  \scantokens{#1\noexpand}%
  \endgroup
}

\usepackage{titlesec}

\renewcommand{\thesection}{Week~\arabic{section}}
\titleformat{\section}[display]{\normalfont}{\Large\bfseries\thesection}{1em}{\large\normalfont}

\usepackage{titletoc}
\titlecontents{section}[0em]{\normalfont}{\bfseries\thecontentslabel\hspace{1em}\normalfont\small}{}{\titlerule*[0.3pc]{.}\small\itshape\thecontentspage}[\vspace{0.5em}]

\usepackage[toc]{multitoc}

\title{This Week's Finds in Mathematical Physics (51--100)}
\author{John Baez}
\date{April 23, 1995 to March 23, 1997}

\usepackage{tikz}
\usetikzlibrary{knots}
\usetikzlibrary{arrows}

\usepackage{etoolbox}
\AtBeginEnvironment{quote}{\itshape}

\usepackage{tikz-cd}

\usepackage{graphicx}
\usepackage[export]{adjustbox}

\makeatletter
\newcommand{\xRightarrow}[2][]{\ext@arrow 0359\Rightarrowfill@{#1}{#2}}
\makeatother

\usepackage[main]{embedall}

\begin{document}

\begin{titlepage}
  \begin{center}
    {\Huge\textbf{This Week's Finds in}}
  \\[0.7em]{\Huge\textbf{Mathematical Physics}}
  \\[1em]{\huge\textit{Weeks 51 to 100}}
  \\[4em]{\LARGE \textit{April 23, 1995} to \textit{March 23, 1997}}
  \\[4em]{\huge by John Baez}
  \\[0.5em]{\Large{Typeset by Tim Hosgood}}
  \end{center}
  
  \vskip 5em
   \[
  \begin{tikzpicture}
    \begin{knot}
      \strand[thick] (0,0.5)
        to (0,0)
        to [out=down,in=down,looseness=2] (1,0)
        to (1,0.5);
    \end{knot}
    \node[fill=white] at (0,0) {$x^*$};
    \node[fill=white] at (1,0) {$x$};
    \node[label=below:{$$}] at (0.5,-0.6) {$$};
    \begin{knot}
      \strand[thick] (-0.75,0.5)
        to (-0.75,0)
        to [out=down,in=up] (0.125,-1.75)
        to (0.125,-2.5);
      \strand[thick] (1.75,0.5)
        to (1.75,0)
        to [out=down,in=up] (0.875,-1.75)
        to (0.875,-2.5);
    \end{knot}
    \node[fill=white] at (-0.75,0) {$x$};
    \node[fill=white] at (1.75,0) {$x^*$};
    \node[fill=white] at (0,-2) {$x$};
    \node[fill=white] at (1,-2) {$x^*$};
    \begin{scope}[shift={(0.875,-3)}]
      \begin{knot}
        \strand[thick] (0,0.5)
          to (0,0)
          to [out=down,in=down,looseness=2] (1,0)
          to (1,0.5);
      \end{knot}
      \node[fill=white] at (0,0) {$x^*$};
      \node[fill=white] at (1,0) {$x$};
      \node[label=below:{$$}] at (0.5,-0.6) {$$};
      \begin{knot}
        \strand[thick] (-0.75,0.5)
          to (-0.75,0)
          to [out=down,in=up] (0.125,-1.75)
          to (0.125,-2.5);
        \strand[thick] (1.75,0.5)
          to (1.75,0)
          to [out=down,in=up] (0.875,-1.75)
          to (0.875,-2.5);
      \end{knot}
      \node[fill=white] at (-0.75,0) {$x$};
      \node[fill=white] at (1.75,0) {$x^*$};
      \node[fill=white] at (0,-2) {$x$};
      \node[fill=white] at (1,-2) {$x^*$};
    \end{scope}
    \begin{scope}[shift={(1.875,-2.5)}]
      \begin{knot}
        \strand[thick] (0,0)
          to (0,0.5)
          to [out=up,in=down,looseness=0.75] (1,2.5)
          to (1,3);
        \strand[thick] (0.75,0)
          to (0.75,0.5)
          to [out=up,in=down,looseness=0.75] (1.75,2.5)
          to (1.75,3);
      \end{knot}
      \node[fill=white] at (0,0.5) {$x$};
      \node[fill=white] at (0.75,0.5) {$x^*$};
      \node[fill=white] at (1,2.5) {$x$};
      \node[fill=white] at (1.75,2.5) {$x^*$};
    \end{scope}
    \node at (4.5,-2.5) {$=$};
    \begin{scope}[xscale=-1,shift={(-9,0)}]
      \begin{knot}
        \strand[thick] (0,0.5)
          to (0,0)
          to [out=down,in=down,looseness=2] (1,0)
          to (1,0.5);
      \end{knot}
      \node[fill=white] at (0,0) {$x$};
      \node[fill=white] at (1,0) {$x^*$};
      \node[label=below:{$$}] at (0.5,-0.6) {$$};
      \begin{knot}
        \strand[thick] (-0.75,0.5)
          to (-0.75,0)
          to [out=down,in=up] (0.125,-1.75)
          to (0.125,-2.5);
        \strand[thick] (1.75,0.5)
          to (1.75,0)
          to [out=down,in=up] (0.875,-1.75)
          to (0.875,-2.5);
      \end{knot}
      \node[fill=white] at (-0.75,0) {$x^*$};
      \node[fill=white] at (1.75,0) {$x$};
      \node[fill=white] at (0,-2) {$x^*$};
      \node[fill=white] at (1,-2) {$x$};
      \begin{scope}[shift={(0.875,-3)}]
        \begin{knot}
          \strand[thick] (0,0.5)
            to (0,0)
            to [out=down,in=down,looseness=2] (1,0)
            to (1,0.5);
        \end{knot}
        \node[fill=white] at (0,0) {$x$};
        \node[fill=white] at (1,0) {$x^*$};
        \node[label=below:{$$}] at (0.5,-0.6) {$$};
        \begin{knot}
          \strand[thick] (-0.75,0.5)
            to (-0.75,0)
            to [out=down,in=up] (0.125,-1.75)
            to (0.125,-2.5);
          \strand[thick] (1.75,0.5)
            to (1.75,0)
            to [out=down,in=up] (0.875,-1.75)
            to (0.875,-2.5);
        \end{knot}
        \node[fill=white] at (-0.75,0) {$x^*$};
        \node[fill=white] at (1.75,0) {$x$};
        \node[fill=white] at (0,-2) {$x^*$};
        \node[fill=white] at (1,-2) {$x$};
      \end{scope}
      \begin{scope}[shift={(1.875,-2.5)}]
        \begin{knot}
          \strand[thick] (0,0)
            to (0,0.5)
            to [out=up,in=down,looseness=0.75] (1,2.5)
            to (1,3);
          \strand[thick] (0.75,0)
            to (0.75,0.5)
            to [out=up,in=down,looseness=0.75] (1.75,2.5)
            to (1.75,3);
        \end{knot}
        \node[fill=white] at (0,0.5) {$x^*$};
        \node[fill=white] at (0.75,0.5) {$x$};
        \node[fill=white] at (1,2.5) {$x^*$};
        \node[fill=white] at (1.75,2.5) {$x$};
      \end{scope}
    \end{scope}
  \end{tikzpicture}
\] 
\end{titlepage}

\begin{center}
{\huge \bf Preface}
\end{center}

\vskip 1em

\noindent
These are issues 51 to 100 of \emph{This Week's Finds of Mathematical Physics}.
This series has sometimes been called the world's first blog, though it was originally 
posted on a ``usenet newsgroup'' called sci.physics.research --- a form of communication 
that predated the world-wide web.  I began writing this series as a way to talk about papers I was reading and writing, and I continued that in issues 51--100, but these issues also contain two expository series.   One is on ``ADE classifications'':
\begin{itemize}
\item \hyperlink{week62_ade}{Week 62} - The ubiquity of ADE classifications.  The
classification of finite reflection groups using Coxeter diagrams.
\item \hyperlink{week63_ade}{Week 63} - The classification of crystallographic 
finite reflection groups and semisimple Lie algebras.
\item \hyperlink{week64_ade}{Week 64} - Semisimple Lie algebras.   Affine Lie
algebras, Wess--Zumino--Witten models and quantum groups.
\item \hyperlink{week65_ade}{Week 65} - The A, D, and E lattices.  The McKay
correpondence and the classification of
minimal models.
\end{itemize}
The other, called the ``The Tale of \(n\)-Categories'', was an introduction to categories
and higher categories --- mainly just \(2\)-categories.  It can be found in these issues:
\begin{itemize}
\item \hyperlink{week73_tale}{Week 73} - The category of sets and the 2-category of categories.
\item \hyperlink{week74_tale}{Week 74} - Kinds of categories: monoids, groups, and groupoids.   The periodic table of \(n\)-categories.
\item \hyperlink{week75_tale}{Week 75} - The fundamental groupoid of a topological space.   The classifying space of a groupoid.
\item \hyperlink{week76_tale}{Week 76} - Equations, isomorphisms, and equivalences.   Adjoint functors.
\item \hyperlink{week77_tale}{Week 77} - Adjoint functors.
\item \hyperlink{week78_tale}{Week 78} - Adjoint functors and adjoint linear operators.
\item \hyperlink{week79_tale}{Week 79} - The unit and counit of an adjunction.
\item \hyperlink{week80_tale}{Week 80} - The definition of 2-category.
\item \hyperlink{week83_tale}{Week 83} - Adjunctions in 2-categories and dual objects in monoidal categories.
\item \hyperlink{week84_tale}{Week 84} - Review.  Monads and monoids.
\item \hyperlink{week89_tale}{Week 89} - Monads in 2-categories.  Monoids and monoidal categories as monads.
\item \hyperlink{week92_tale}{Week 92} - Monads from adjunctions.
\item \hyperlink{week99_tale}{Week 99} - 2-Hilbert spaces.  Coproducts.
\item \hyperlink{week100_tale}{Week 100} - Definitions of \(n\)-category.
\end{itemize}

Tim Hosgood kindly typeset all 300 issues of \emph{This Week's Finds} in 2020. They
will be released in six installments of 50 issues each.  I have edited the issues here
to make the style a bit more uniform and also to change some references to preprints,
technical reports, etc.\ into more useful arXiv links.  This accounts for some anachronisms
where I discuss a paper that only appeared on the arXiv later.   

I thank Fridrich Valach for helping me fix many mistakes.  
There are undoubtedly many still remaining. 
If you find some, please contact me and I will try to fix them.

\tableofcontents

hypertarget{week51}{%
\section{April 23, 1995}\label{week51}}

For people in theoretical physics, Trieste is a kind of mecca. It's an
Italian town on the Adriatic quite near the border with Slovenia, and
it's quite charming, especially the castle of Maximilian near the coast,
built when parts of northern Italy were under Hapsburg rule. Maximilian
later took his architect with him to Mexico when he became Emperor
there, who built another castle for him in Mexico City. (The Mexicans,
apparently unimpressed, overthrew and killed Maximilian.) These days,
physicists visit Trieste partially for the charm of the area, but mainly
to go to the ICTP and SISSA, two physics institutes, the latter of which
has grad students, the former of which is purely for research. There are
lots of conferences and workshops at Trieste, and I was lucky enough to
be invited to Trieste while one I found interesting was going on.

As I described to some extent in \protect\hyperlink{week44}{``Week 44''}
and \protect\hyperlink{week45}{``Week 45''}, Seiberg and Witten have
recently shaken up the subject of Donaldson theory by using some
physical reasoning to radically simplify the computations involved.
Donaldson theory has always had a lot to do with physics, since it uses
the special features of of gauge theory in 4 dimensions to obtain
invariants of \(4\)-dimensional manifolds. So perhaps it is not
surprising that physicists have had a lot to say about Donaldson theory
all along, even before the recent Seiberg--Witten revolution. And indeed,
at Trieste lots of mathematicians and physicists were busy talking to
each other about Donaldson theory, trying to catch up with the latest
stuff and trying to see what to do next.

Now I don't know much about Donaldson theory, but I have a vague
interest in it for various reasons. First, it's \emph{supposed} to be a
4-dimensional topological quantum field theory, or TQFT. Indeed, the
very first paper on TQFTs was about Donaldson theory in 4 dimensions:

\begin{enumerate}
\def\labelenumi{\arabic{enumi})}
\tightlist
\item
 Edward Witten,  ``Topological quantum field theory'', \emph{Comm.
  Math. Phys.} \textbf{117} (1988) 353.
\end{enumerate}

\noindent
Only later did Witten turn to the comparatively easier case of
Chern--Simons theory, which is a \(3\)-dimensional TQFT:

\begin{enumerate}
\def\labelenumi{\arabic{enumi})}
\setcounter{enumi}{1}
\tightlist
\item
   Edward Witten, ``Quantum field theory and the Jones polynomial'',
  \emph{Commun. Math. Phys.} \textbf{121} (1989) 351.
\end{enumerate}

However, when \emph{mathematicians} talk about TQFTs they usually mean
something satisfying Atiyah's axioms for a TQFT (which are nicely
presented in his book --- see \protect\hyperlink{week39}{``Week 39''}).
Now it turns out that Chern--Simons theory can be rigorously constructed
as a TQFT satisfying these axioms most efficiently using braided
monoidal categories, which play a big role in 3d topology. So it makes
quite a bit of sense in a \emph{general} sort of way that Crane and
Frenkel are trying to construct Donaldson theory using braided monoidal
\(2\)-categories, which seem to play a comparable role in 4d topology.
In the paper which I cite in \protect\hyperlink{week50}{``Week 50''},
they begin to construct a braided monoidal \(2\)-category related to the
group \(\mathrm{SU}(2)\), which they conjecture gives a TQFT related to
Donaldson theory. That also makes some \emph{general} sense, because
Donaldson theory, at least ``old'' Donaldson theory, is closely related
to gauge theory with gauge group \(\mathrm{SU}(2)\). Still, I've always
wanted to see a more \emph{specific} reason why Donaldson theory should
be related to the Crane-Frenkel ideas, not necessarily a proof, but at
least a good heuristic argument.

Luckily George Thompson, who invited me to Trieste, knows a bunch about
TQFTs. Unluckily he was sick and I never really got to talk to him very
much! But luckily his collaborator Matthias Blau was also there, so I
took the opportunity to pester him with questions. I learned a bit, most
of which is in their paper:

\begin{enumerate}
\def\labelenumi{\arabic{enumi})}
\setcounter{enumi}{2}
\tightlist
\item
  Matthias Blau and George Thompson, ``\(N = 2\) topological gauge theory, 
  the Euler characteristic of moduli spaces, and the Casson invariant'', 
  \emph{Comm.\ Math.\ Phys.} \textbf{152} (1993), 41--71.
\end{enumerate}
This paper helped me a lot in understanding Crane and Frenkel's ideas.
But so that this ``week'' doesn't get too long, I'll just focus on one
basic aspect of the paper, which is the importance of supersymmetric
quantum theory for TQFTs. Then next week I'll say a bit more about the
Donaldson theory business.

If you look at Witten's paper on Donaldson theory above, you'll see he
writes down the Lagrangian for a ``supersymmetric'' field theory, which
is supposed to be a TQFT, namely, Donaldson theory. Supersymmetric field
theories treat bosons and fermions in an even-handed way. But why does
supersymmetry show up here? The connection with TQFTs is actually pretty
simple and beautiful, at least in essence.

Suppose we are doing quantum field theory, and ``space'' (as opposed to
``spacetime'') is some manifold \(M\). Then we have some Hilbert space
of states \(Z(M)\) and some Hamiltonian \(H\), which is a self-adjoint
operator on \(Z(M)\). To evolve a state (a vector in \(Z(M)\)) in time,
we hit it with the unitary operator \(\exp(-itH)\), where \(t\) is the
amount of time we want to evolve by, and the minus sign is just a
convention designed to confuse you.

We can think of this geometrically as follows. We are taking spacetime
to be \([0,t] \times M\). You can visualize spacetime as a kind of pipe,
if you want, and then imagine sticking in the state \(\psi\) at one end
and having \(\exp(-itH)\psi\) pop out at the other end.

Now say we bend the pipe around and connect the input end to the output
end! Then we get the spacetime \(S^1\times M\), where \(S^1\) is the
circle of circumference \(t\), formed by gluing the two ends of the
interval \([0,t]\) together. For this kind of ``closed'' spacetime, or
compact manifold, a quantum field theory should give us not an operator
like \(\exp(-itH)\), but a number, the ``partition function'', which in
this case is just the \emph{trace} \(\operatorname{tr}(\exp(-itH))\).

The deep reason for this is that taking the trace of an operator ---
remember, that means taking the sum of the diagonal entries, when you
think of it as a matrix --- is really very much like as taking a pipe
and bending it around, connecting the input end to the output end,
forming a closed loop. This may seem bizarre, but observe that taking
the sum of the diagonal entries really is just a quantitative measure of
how much the ``output constructively interferes with the input''. (And a
very nice one, since it winds up not depending on the basis in which we
write the matrix!) This sort of idea is basic in the Bohm-Aharonov
effect, where we take a particle in an electromagnetic field around a
loop and see how much it interferes with itself, and it is also the
basic idea of a ``Wilson loop'', where we do the same thing for a
particle in a gauge field. In other words, the trace measures the amount
of ``positive feedback''. If this still seems bizarre, or just vague,
take a look at:

\begin{enumerate}
\def\labelenumi{\arabic{enumi})}
\setcounter{enumi}{3}
\tightlist
\item
   Louis Kauffman, \emph{Knots and Physics}, World Scientific,
  Singapore, 1991.
\end{enumerate}

You will see that the same idea shows up in knot theory, where taking a
trace corresponds to taking something (like a braid or tangle) and
folding it over to connect the input and output. In a later ``week''
I'll talk a bit about a new paper by Joyal, Street and Verity that
studies the notion of ``trace'', ``feedback'' and ``folding over'' in a
really general context, the context of category theory.

Anyway, the partition function \(\operatorname{tr}(\exp(-itH))\)
typically depends on \(t\), or in other words, it depends on the
circumference of our circle \(S^1\), not just on the topology of the
manifold \(S^1\times M\). In a TQFT, the partition function is only
supposed to depend on the topology of spacetime! So, how can we get
\(\operatorname{tr}(\exp(-itH))\) to be independent of \(t\)?

There is a banal answer and a clever answer. The banal answer is to take
\(H = 0\)! Then \(\operatorname{tr}(\exp(-itH)) = \operatorname{tr}(1)\)
is just the \emph{dimension} of the Hilbert space:
\[\operatorname{tr}(\exp(-itH)) = \dim(Z(M)).\]  Actually this isn't quite
as banal as it may sound; indeed, the basic equation of quantum gravity
is the Wheeler--DeWitt equation, \[H \psi = 0,\] which must hold for all
physical states. In other words, in quantum gravity there is a big space
of ``kinematical states'' on which \(H\) is an operator, but the really
meaningful ``physical states'' are just those in the subspace
\[Z(M) = \{\psi: H \psi = 0\}.\] Read \protect\hyperlink{week11}{``Week
11''} for more on this.

But there is a clever answer involving supersymmetry! You might hope
that there were some more interesting self-adjoint operators \(H\) such
that \(\operatorname{tr}(\exp(-itH))\) is time-independent, but there
aren't. So we seem stuck. This reminds me of a course I took from Raoul
Bott. He said ``So, we think about the problem... and think some more... 
and still we are stuck, so what should we do? \emph{Superthink!}''

Recall that the Hamiltonian of a free particle in quantum mechanics is
--- up to boring constants --- just minus the Laplacian on configuration
space which is some Riemannian manifold that the particle roams around
on. For this Hamiltonian, \(\operatorname{tr}(\exp(-itH))\) doesn't
quite make sense, since the Hilbert space is infinite-dimensional and
the sum of the diagonal matrix entries diverges. But
\(\operatorname{tr}(\exp(-tH))\) often \emph{does} converge. This is why
folks often replace \(t\) by \(-it\) in formulas, which is called
``going to imaginary time'' or a ``Wick transform''; it really amounts
to replacing Schrodinger's equation by the heat equation: i.e., instead
of a quantum particle, we have a particle undergoing Brownian motion! In
any event, \(\operatorname{tr}(\exp(-tH))\) certainly depends on \(t\)
in these situations, but there is something very similar that does \emph{not}.

Namely, let's replace the Laplacian on \emph{functions} by the Laplacian
on \emph{differential forms}. I won't try to remind you what these are;
I'll simply note that functions are 0-forms, but there are also
\(1\)-forms, 2-forms, and so on --- tensor fields of various sorts ---
and the Laplacian of a \(j\)-form is another \(j\)-form. So for each
\(j\) we get a kind of Hamiltonian \(H_j\), which is just minus the
Laplacian on \(j\)-forms. We can also consider the space of \emph{all}
forms, never mind the \(j\), and on this space there is a Hamiltonian
\(H\), which is just minus the Laplacian on \emph{all} forms. Now, we
could try to take the trace of \(\exp(-tH)\), but it's more interesting
to take the ``supertrace'':
\[\operatorname{str}(\exp(-tH)) = \operatorname{tr}(\exp(-tH_0)) - \operatorname{tr}(\exp(-tH_1)) + \operatorname{tr}(\exp(-tH_2)) - \cdots\]
in other words, the trace of \(\exp(-tH)\) acting on even forms,
\emph{minus} the trace on odd forms.

Why? Well, odd forms are sort of ``fermionic'' in nature, while even
forms are sort of ``bosonic''. The idea of supersymmetry is to throw in
minus signs when you've got ``odd things'', because they are like
fermions, and physicists know that lots formulas for fermions are just
like formulas for bosons, which are ``even things'', except for these
signs. That's the rough idea. It's all related to how, when you
interchange two identical bosons, their wavefunction remains unchanged,
while for fermions it picks up a phase of \(-1\).

Now the amazing cool thing is that \(\operatorname{str}(\exp(-tH))\) is
independent of \(t\). This follows from some stuff called Hodge theory,
or, if you want to really show off, index theory. Basically it works
like this. If you have an operator \(A\) with eigenvalues \(\lambda_i\),
then \[\operatorname{tr}(\exp(-tA)) = \sum_i \exp(-t \lambda_i)\] if the
sum makes sense. We can use this formula to write out
\(\operatorname{str}(\exp(-tH))\) in terms of eigenvalues of the
Laplacians \(H_j\), and it turns out that all the terms coming from
nonzero eigenvalues exactly cancel! So all that's left is the part
coming from the zero eigenvalues, which is independent of \(t\). If you
believe this for a second, it means we can compute the supertrace by
taking the limit as \(t\to\infty\). The eigenvalues are all nonnegative,
so all the quantities \(\exp(-t \lambda_i)\) go to zero except for the
zero eigenvalues, and we're left with \(\operatorname{str}(\exp(-tH))\)
being equal to the alternating sum of the dimensions of the spaces
\[\{\psi : H_j \psi = 0\}\] Now in fact, Hodge theory tells us that
these spaces are really just the ``cohomology groups'' of our
configuration space, so the answer we get for
\(\operatorname{str}(\exp(-tH))\) is what folks call the ``Euler
characteristic'' of our configuration space... an important
topological invariant.

So, generalizing the heck out of this idea, we can hope to get TQFTs
from supersymmetric quantum field theories as follows. Start with some
recipe for associating to each choice of ``space'' \(M\) a
``configuration space'' \(C(M)\) --- some space of fields on \(M\),
typically. Let \(Z(M)\) be the space of all forms on \(C(M)\), and let
\(H\) be the minus the Laplacian, as an operator on \(Z(M)\). Then we
expect that the partition function \(\operatorname{str}(\exp(-tH))\)
will be independent of \(t\). This is just what one wants in a TQFT.
Moreover, the partition function will be the Euler characteristic of the
configuration space \(C(M)\).

But what if we want to get a TQFT out of this trick, and avoid reference
to the Laplacian? Then we can just do the following equivalent thing (at
least it's morally equivalent: there will usually be things to check).
Let \(Z_+(M)\) be the direct sum of all the even cohomology groups of
\(C(M)\), and let \(Z_-(M)\) be the direct sum of all the odd ones. Then
\[\operatorname{str}(\exp(-tH)) = \dim(Z_+(M))- \dim(Z_-(M))\] so what we
expect is, not quite a TQFT in the Atiyah sense, but a ``superTQFT''
whose space of states has an ``even'' part equal to \(Z_+(M)\) and an
``odd'' part equal to \(Z_-(M)\); the right hand side is then the
``superdimension'' of the space of states this ``superTQFT'' assigns to
\(M\).

Now actually in real life things get tricky because the configuration
space \(C(M)\) might be infinite-dimensional, or a singular variety. If
\(C(M)\) is too weird, it gets hard to say what its Euler characteristic
should be! But as Blau and Thompson's paper and the references in it
point out, one can often still make it make sense, with enough work. In
particular, when we are dealing with Donaldson theory, \(C(M)\) is just
the moduli space of flat \(\mathrm{SU}(2)\) connections on \(M\). This
means that the partition function of \(S^1\times M\) should be the Euler
characteristic of moduli space, better known as the Casson invariant.
And what is the vector space our superTQFT assigns to \(M\)? Well, it's
called Floer homology. Now actually there are a lot of subtleties here
I'm deliberately sloughing over. Read Blau and Thompson's paper for some
of them --- and read the references for more!

\hypertarget{week52}{%
\section{May 9, 1995}\label{week52}}

So, last ``week'', I said a bit about how supersymmetry could be handy
for constructing topological quantum field theories, and this week I
want to say a bit more about what that has to do with getting a purely
combinatorial description of Donaldson theory.

But first, I want to lighten things up a bit by mentioning a good
science fiction novel!

\begin{enumerate}
\def\labelenumi{\arabic{enumi})}
\tightlist
\item
  Greg Egan, \emph{Permutation City}, Millenium, 1994.
\end{enumerate}

\noindent
There is a lot of popular interest these days in the anthropic
principle. Roughly, this claims to explain certain features of the
universe by noting that if the universe didn't have those features,
there would be no intelligent life. So, presumably, the very fact that
we are here and asking certain questions guarantees that the questions
will have certain answers.

Of course, the anthropic principle is controversial. Suppose one could
really show that if the universe didn't have property \(X\), there would
be no intelligent life. Does this really count as an ``explanation'' of
property \(X\)? People like arguing about this. But this question is
much too subtle for a simple-minded soul such as myself. I'm still stuck
on more basic things!

For example, are there any examples where we \emph{can} really show that
if the universe didn't have property \(X\), there would be no
intelligent life? It seems that to answer this, we need to have some
idea about what we're counting as ``all possible universes'', and what
counts as ``intelligent life''. So far we only know \emph{one} example of a
universe and \emph{one} example of intelligent life, so it is difficult to
become an expert on these subjects! It'd be all too easy for us to
unthinkingly assume that all intelligent life is carbon-based,
metabolizes using oxidation, and eats pizza, just because folks around
here do.

Our unthinking parochialism is probably all the worse as far as
different universes are concerned! What counts as a possible universe,
anyway? Rather depressingly, we must admit we don't even know the laws
of \emph{this} universe, so we don't really know what it takes for a
universe to be possible, in the strong sense of capable of actually
existing as a universe. We are a little bit better off if we consider
all \emph{logically possible} universes, but not a whole lot better.
Certainly every axiom system counts as a logically possible set of laws
of a universe --- every set of axioms in every possible formal system. But
who is to say that universes must have laws of this form? We don't even
know for sure that \emph{ours} does!

So this whole topic will remain a hopeless quagmire until one takes a
small, carefully limited piece of it and studies that. People studying
artificial life are addressing one of these bite-sized pieces, and
getting some interesting results. I hope everyone has heard about Thomas
Ray's program Tierra, for example: he created an artificial ecosystem -
one could call it a ``possible universe'' --- and found, after seeding it
with one self-reproducing program, a rapid evolution of parasites, etc.,
following many of the patterns of ecology here. But so far, perhaps
merely due to time and memory limitations, no intelligence!

\emph{One} of the cool things about \emph{Permutation City} is an imagined
cellular automaton, the ``Autoverse'', which is complicated enough to
allow life. But something much cooler is the main theme of the book.
Egan calls it the ``Dust Theory''. It's an absolutely outrageous theory,
but if you think about it carefully, you'll see that it's rather hard to
spot a flaw. It depends on the tricky puzzles concealed in the issue of
``isomorphism''.

Being a mathematician, one thing that always puzzled me about the
notions of ``intelligent life'' and ``all possible universes'' was the
question of isomorphisms between universes. Certainly we all agree that,
say, the Heisenberg ``matrix mechanics'' and Schrodinger ``wave
mechanics'' formulations of quantum mechanics are isomorphic. In both of
them, the space of states is a Hilbert space, but in one the states are
described as sequences of numbers, while in the other they are described
as wavefunctions. At first they look like quite different theories. But
in a while people realized that there was a unitary operator from
Heisenberg's space of states to Schrodinger's, and that via this
correspondence all of matrix mechanics is equivalent to wave mechanics.

So does Heisenberg's universe count as the same one as Schrodinger's, or
a different one? It seems clear that they're the same. But say we had
two quantum-mechanical systems whose Hamiltonians have the same
eigenvalues (or spectrum); does that mean they are the ``same'' system,
really? Is that all there is to a physical system, a list of
eigenvalues??? If we are going to go around talking about ``all possible
universes'', it would probably pay to think a little about this sort of
thing!

Say we had two candidates for ``laws of the universe'', written down as
axioms in different formal systems. How would we decide if these were
describing different universes, or were simply different ways of talking
about the same universe? Pretty soon it becomes clear that the issue is
not a black-and-white one of ``same'' versus ``different'' universes.
Instead, laws of physics, or universes satisfying these laws, can turn
out to be isomorphic or not depending on how much structure you want the
isomorphism to preserve. And even if they are isomorphic, there may not
be a ``unique'' isomorphism or a ``canonical'' isomorphism. (Very
roughly speaking, a canonical isomorphism is a ``God-given best one'',
but one can use some category theory to make this precise.) If you think
about this carefully you'll see that our universe could be isomorphic to
some very different-seeming ones, or could have some very
different-seeming ones `embedded' in it.

Greg Egan takes this issue and runs with it -- in a very interesting
direction. Everyone interested in cellular automata, artificial life,
virtual reality, or other issues of simulation should read this, as well
as anyone who likes philosophy or just a good story.

Okay, back to business here....

\begin{enumerate}
\def\labelenumi{\arabic{enumi})}
\setcounter{enumi}{1}
\item
  Alberto Cattaneo, ``Teorie topologiche di tipo BF ed invarianti dei
  nodi'', doctoral thesis, department of physics, University of Milan.

  Alberto Cattaneo, Paolo Cotta-Ramusino, J\"urg Fr\"ohlich, and Maurizio
  Martellini, ``Topological BF theories in 3 and 4 dimensions'',
  available as
  \href{https://arxiv.org/abs/hep-th/9505027}{\texttt{hep-th/9505027}}.
\end{enumerate}
So, last week I said a wee bit about Blau and Thompson's paper on
supersymmetry and the Casson invariant. All I said was that suitably
concocted supersymmetric field theories could be used to compute the
Euler characteristics of your favorite spaces, and that Blau and
Thompson talked about one which computed the Casson invariant, which is
(in a rather subtle sense) the Euler characteristic of the moduli space
of flat connections on a trivial \(\mathrm{SU}(2)\) bundle over a
3-manifold. Traditionally one requires that the 3-manifold be a homology
3-sphere, but Kevin Walker showed how to do it for rational homology
spheres, and Blau and Thompson mention other work in which the Casson
invariant is generalized still further.

But I didn't say \emph{which} supersymmetric field theory computes the
Casson invariant for you. The answer is, \(N = 2\) supersymmetric \(BF\)
theory with gauge group \(\mathrm{SU}(2)\). So now I should say a little
about \(BF\) theory. Actually I have already mentioned it here and
there, especially in \protect\hyperlink{week36}{``Week 36''}. But I
should say a bit more. This is going to be pretty technical, though, so
fasten your seatbelts.

The people I know who are the most excited about \(BF\) theory are the
folks I was visiting at Milan, namely Cotta-Ramusino, Martellini and his
student Cattaneo. They are working on \(BF\) theory in 3 and 4
dimensions. Let me talk about \(BF\) theory in 3 dimensions, which is
what's most directly relevant here. Well, in \emph{any} dimension, say
\(n\), the fields in \(BF\) theory are a connection \(A\) on a trivial
bundle (take your favorite gauge group \(G\)), whose curvature \(F\)
we'll think of as a \(2\)-form taking values in the Lie algebra of
\(G\), and Lie-algebra-valued \((n-2)\)-form \(B\). Then the Lagrangian
of the theory is \[L(B,F) = \operatorname{tr}(B \wedge F)\] where in the
``trace'' we're using something like the Killing form to get an honest
\(n\)-form which we can integrate over spacetime.

But in 3 dimensions, since \(B\) is a \(1\)-form, you can add an extra
``cosmological constant'' term and take as your Lagrangian
\[L(B,F,c) = \operatorname{tr}(B \wedge F + (c^2/3) B \wedge B \wedge B)\]
where I have put in ``\(c^2/3\)'' as my ``cosmological constant'' for
insidious reasons to become clear momentarily. Now what the article by
Cattaneo, Cotta-Ramusino, Froehlich and Martellini makes really clear is
how \(BF\) theory is related to Chern--Simons theory. This is implicit in
Witten's work on 3d gravity (see \protect\hyperlink{week16}{``Week
16''}), which is just the special case where \(G\) is
\(\mathrm{SO}(2,1)\) or \(\mathrm{SO}(3)\), and where the cosmological
constant really is the usual cosmological constant. But I'd never
noticed it. Recall that the Chern--Simons action is
\[L(A) = \operatorname{tr}(A \wedge dA + (2/3)A \wedge A \wedge A)\]
Thus if we have \(1\)-form B around as well, we can set \[
  \begin{aligned}
    A' &= A + cB,
  \\A'' &= A - cB
  \end{aligned}
\] so we get two different Chern--Simons theories with actions \(L(A')\)
and \(L(A'')\), respectively. \emph{Or}, we can form a theory whose action is
the difference of these two, and, lo and behold:
\[L(A') - L(A'') = 4cL(B,F,c).\] In other words, \(BF\) theory with
cosmological constant is just a ``difference of two Chern--Simons
theories''. Fans of topological quantum field theory may perhaps be more
familiar with this if I point out that the Turaev--Viro theory is just
\(BF\) theory with gauge group \(\mathrm{SU}(2)\), and the fact that the
partition function for this theory is the absolute value squared of that
for Chern--Simons theory is a special case of what I'm talking about. The
nice thing about all this is that the funny phases coming from framings
in Chern--Simons theory precisely cancel out when you form this
``difference of two Chern--Simons theories''.

Now the Casson invariant is related to \(BF\) theory in 3 dimensions
\emph{without} cosmological constant, i.e., taking \(c = 0\). We might
get worried by the equation above, which we can't solve for \(L(B,F)\)
when \(c = 0\), but as Cattaneo and company point out,
\[L(B,F) = \lim_{c\to0}\frac{L(A')-L(A'')}{4c}\] so \(BF\) theory
without cosmological constant is just a limiting case, actually a kind
of \emph{derivative} of Chern--Simons theory. They use this to make
clearer the relation between the vacuum expectation values of Wilson
loops in Chern--Simons theory --- which give you the HOMFLY polynomial
for \(G = \mathrm{SU}(N)\) --- and the corresponding vacuum expectation
values in \(BF\) theory without cosmological constant --- which give you
the Alexander polynomial! Very pretty stuff.

Now back to the Casson invariant and some flagrant speculation on my
part concerning Crane and Frenkel's ideas on Donaldson theory. (I said
last week that this is where I was heading, and now I'm almost there!)
Okay: we know how to define Chern--Simons theory in a purely
combinatorial way using quantum groups. I.e., we can compute the
partition function of Chern--Simons theory with gauge group \(G\) using
the quantum version of the group \(G\)... let me just call it
``quantum \(G\)''. If we take \(c\) to be imaginary above, one can show
that \(BF\) theory with cosmological constant can be computed in a very
similar way starting with the quantum group corresponding to the
\emph{complexification} of \(G\), i.e.~``quantum \(\mathbb{C}G\)''. The
point is that \(A+cB\) can then be thought of as a connection on a
bundle with gauge group \(\mathbb{C}G\). So far this is not flagrant
speculation. Slightly more flagrantly, but not really very much at all,
the formulas above hint that \(BF\) theory without cosmological constant
can be computed in a similar way starting with the quantum group
corresponding to the \emph{tangent bundle} of \(G\), or ``quantum
\(TG\)''. (The tangent bundle of a Lie group is again a Lie group, and
as we let \(c \to 0\) what we are really doing is taking a limit in
which \(\mathbb{C}G\) approaches \(TG\); folks call this a
``contraction'', and in the \(\mathrm{SU}(2)\) case many of the details
appear in Witten's paper on 3d quantum gravity; the tangent bundle of
\(\mathrm{SO}(2,1)\) being just the Poincar\'e group in 3 dimensions.) If
anyone knows whether folks have worked out the quantization of these
tangent bundle groups, let me know! I think some examples have been
worked out.

Okay, but Blau and Thompson say that to compute the Casson invariant you
need to use, not \(BF\) theory with gauge group \(\mathrm{SU}(2)\), but
\emph{supersymmetric} \(BF\) theory with gauge group \(\mathrm{SU}(2)\).
Well, no problemo --- just compute it with ``quantum
super-\(T(\mathrm{SU}(2))\)''! Here I'm getting a bit flagrant; there
\emph{are} theories of quantum supergroups, but I don't know much about
them, especially ``quantum super-\(TG\)'' for \(G\) compact semisimple.
Again, if anybody does, please let me know! (Actually Blau told me to
check out a paper by Saleur and somebody on this, but I never
did....)

Okay, but now let's get seriously flagrant. Recall that the Casson
invariant is really the Euler characteristic of something, just a
number, but this number is just the superdimension of a
super-vector-space, namely the Floer cohomology. From numbers to vector
spaces: this is a typical sort of ``categorification'' process that one
would expect as one goes from 3d to 4d TQFTs. And indeed, folks suspect
that the Floer cohomology is the space of states for a 4d TQFT, or
something like a 4d TQFT, namely Donaldson theory. (``Something like
it'' because of many quirky twists that one wouldn't expect of a
full-fledged TQFT satisfying the Atiyah axioms.) So, just as the Casson
invariant is associated to a certain Hopf algebra, namely ``quantum
super-\(T(\mathrm{SU}(2))\)'', we'd expect Donaldson theory to be
associated to a certain Hopf \emph{category}, the ``categorification of
quantum super-\(T(\mathrm{SU}(2))\)''. So all we need to do is figure
out how to categorify quantum super-\(T(\mathrm{SU}(2))\) and we've got
a purely combinatorial definition of Donaldson theory!

Well, that's not quite so easy, of course. And I may have made, not only
the inevitable errors involved in painting a simplified sketch of what
is bound to be a rather big task, but also other worse errors. Still,
this business should clarify, if only a wee bit, what Crane and Frenkel
are up to when they are categorifying \(\mathrm{SU}(2)\). In fact, it's
likely that working with \(\mathrm{SU}(2)\) rather than
\(T(\mathrm{SU}(2))\) will remove some of the divergences from the state
sum, since, being compact, \(\mathrm{SU}(2)\) has a discrete set of
representations (and quantum \(\mathrm{SU}(2)\) has finitely many
interesting ones, at roots of unity). So they may get a theory that's
allied to but not exactly the same as Donaldson theory, yet
better-behaved as far as the TQFT axioms go.

If anyone actually does anything interesting with these ideas I'd very
much appreciate hearing about it.

\hypertarget{week53}{%
\section{May 18, 1995}\label{week53}}

Near the end of April I was invited by Ronnie Brown to Bangor, Wales for
a very exciting get-together. Readers of ``This Week's Finds'' will know
I'm interested in \(n\)-categories and higher-dimensional algebra these
days. Brown is the originator of the term ``higher-dimensional algebra''
and has been sort of a prophet of the subject for many years. Tim Porter
at Bangor also works on this subject; I'll try to say a bit more about
his stuff next week. And visiting Bangor at the time were John Power and
Ross Street, two category theorists who do a lot of work on
\(n\)-categories. So I had a chance to learn some more
higher-dimensional algebra and category theory and see what these folks
thought of my crazy ideas.

\begin{enumerate}
\def\labelenumi{\arabic{enumi})}
\tightlist
\item
  Ronald Brown, ``Out of line'', \emph{Royal Institution Proceedings}
  \textbf{64}, 207--243.
\end{enumerate}

Brown is very interested in explaining mathematics to the public, and
this paper is based on a talk he gave to a general audience. It is a
very accessible introduction to what mathematics is really all about,
and what higher-dimensional algebra is about in particular. ``Out of
line'' is a pun, of course, because not only does higher-dimensional
algebra seek to burst free of certain habits of ``linear thinking'' that
tend to go along with symbol string manipulation, it also has been a bit
outside the mainstream of mathematics until recently.

Now, when I speak of ``linear thinking'' I am not indulging in some
vague new-agey complaint against rationality. I mean something very
precise: the tendency to focus ones energy on operations that are easily
modelled by the juxtaposition of symbols in a line. The primordial
example is addition: we have a bunch of sticks in a row:
\[\vert\vert\vert\vert\vert\] and we say there are ``5'' sticks and
write \[1+1+1+1+1=5.\] Fine. But when we have a \(2\)-dimensional array
of sticks:
\[\begin{gathered}\vert\vert\vert\vert\\\vert\vert\vert\vert\end{gathered}\]
we are in a hurry to bring the situation to linear form by making up a
new operation, ``multiplication'', and saying we have \(2 \times 4\)
sticks. This isn't so bad for plenty of purposes; it's not as if I'm
against times tables! But certain things, particular in topology, can
get obscured by neglecting operations that correspond most naturally to
higher-dimensional forms of juxtaposition, and Brown's paper explains
some of these, and how to deal with these problems. The point is not to
avoid linear notation; it's to avoid falling into certain mental traps
it can lead you into if you're not careful!

\begin{enumerate}
\def\labelenumi{\arabic{enumi})}
\setcounter{enumi}{1}
\tightlist
\item
  A.\ J.\ Power, ``Why tricategories?'', \emph{Information and Computation} 
  \textbf{120} (1995), 251--262.  
\end{enumerate}
\noindent
When I mentioned this paper to a friend, she puzzledly asked ``\,`Why
try categories?'?'' And indeed, one must have tried categories and
enjoyed them before moving on to bicategories, tricategories and that
great beckoning terra incognita of mathematics, \(n\)-category theory.

In a sense I already know ``why tricategories''. I think they're great,
and in a paper with James Dolan --- summarized in
\protect\hyperlink{week49}{``Week 49''} --- I did my best to get
everyone else interested in general \(n\)-categories. For me, the great
thing about \(n\)-category theory is that it strives to formalize the
notion of ``process'' in all its recursive splendor. An \(n\)-category
is a mathematical structure containing not only objects, which one might
think of as ``things'', and morphisms, which one might think of as
``processes leading from one thing to another'', but also
\(2\)-morphisms, which are ``processes leading from one process to
another'', and then \(3\)-morphisms, etc., on up to \(n\)-morphisms.

In physics and topology applications, the \(k\)-morphisms can often be
thought of as \(k\)-dimensional geometrical objects, since (as the
Greeks knew) the process of motion of a point traces out a
\(1\)-dimensional figure, and similarly the motion of a
\(1\)-dimensional figure traces out a \(2\)-dimensional surface... and
\(n\)-dimensional spacetime is in some rough sense ``traced out'' by the
motion of \((n-1)\)-dimensional spacelike slices through time. If you
think this is vague, you're right --- and that's why one needs
\(n\)-category theory, to make it precise! When one understands
\(n\)-categories (which so far we really do only up to \(n = 3\)) one
sees that the possibilities inherent in \(n\)-dimensional topology match
up very nicely with one might have stumbled on not knowing topology at
all, but just playing around with this iterated notion of processes
between processes between processes.... This ``natural
correspondence'' between purely algebraic concepts and topological ones
is what makes topological quantum field theory tick, and I can't help
but feel that it will have quite a bit to say about physics eventually.

Power, however, gives a quite different set of reasons for being
interested in tricategories. These are drawn from computer science and
logic, and they make me realize yet again how poor and outdated my
education in logic was, and how much interesting stuff there is going on
in the subject!

At a completely naive level, one might already expect that relation
between ``processes'' and ``things'' is subtle and interesting in
computer science. For after all, we can think of a program either as a
process that turns one thing into another, or as data, a sort of thing,
which other programs can act on. Power does not really emphasize this
issue explicitly, but I can't help remarking on it, especially because
it reminds me of the curious fact that in mathematical physics, ``time
is the last dimension''.

That is, in topological quantum field theory, the \(n\)-morphisms in an
\(n\)-category, which are the processes having no further processes
going between them, represent the passage of time. And indeed, for many
practical purposes it appears that \(n = 4\) is where things leave off,
since spacetime appears \(4\)-dimensional. On the other hand, nobody
knows any mathematical reason why one has to stop at any given \(n\). So
ultimately we should try to understand ``\(\omega\)-categories'', having
\(n\)-morphisms for all \(n\) greater than or equal to zero (0-morphisms
being simply objects, and \(1\)-morphisms being morphisms). This
corresponds philosophically to the notion that every ``process'' can
also be regarded as a ``thing'' which other processes can transform.
Moreover, we should also try to understand
``\(\mathbb{Z}\)-categories'', having \(n\)-morphisms for all integers
\(n\), even negative ones! In this world, where there is no bottom as
well as no top, every ``thing'' can also be regarded as a ``process''.

But I digress. Power is actually more interested in a different (though
perhaps related) hierarchy. Sometimes people like to say computers just
do stuff with bunches of numbers, but that's pretty misleading. Of
course computers \emph{can} do things with numbers, but that's one of
the simpler mathematical things they can do. A number is an element of a
set (the set of real numbers, or some set of more computer-manageable
numbers.) And computers have no problem dealing with elements of sets.
But computers can also deal with sets themselves --- and more fancy
mathematical objects.

Many mathematical objects are sets, or bunches of sets, equipped with
operations satisfying equational laws. For example, a group is a set
equipped with a product and inverse operation satisfying various laws.
Sometimes these operations are only defined if certain conditions hold,
of course. For example, a category is a set of ``objects'' and a set of
``morphisms'', together with various operations like composition of
morphisms, but one can only compose two morphisms \(f\colon x\to y\) and
\(g\colon w\to z\) if \(y = w\). Other examples might include graphs,
partially ordered sets... and all sorts of things computer
scientists know and love.

We could call all of these ``sets with essentially algebraic
structure.'' Mathematically sophisticated computer scientists want to be
able define data types corresponding to arbitrary sorts of sets with
essentially algebraic structure, and to play around with them easily. So
they need to ponder such things in considerable generality.

Note that in all cases, there is not just a bunch of objects to play
with --- like ``groups'' or ``partially ordered sets'' --- but a
\emph{category} in which the morphisms are structure-preserving maps
between the objects in question. For example, there is a category
\(\mathsf{Grp}\) whose objects are groups and whose morphisms are group
homomorphisms.

The categories one gets this way are of a certain sort. Power calls them
``categories of models of finite limit theories''. To define this
requires a bit of know-how, but it's basically simple. For example,
suppose I were trying to explain the definition of a category to a
computer scientist; I might say, every category has a set
\(\mathrm{ob}\) of objects and a set \(\mathrm{mor}\) of morphisms;
every morphism has an object called its source (or domain), so there is
a function \[\operatorname{source}\colon\mathrm{mor}\to\mathrm{ob}\] and
similarly every morphism has an object called its target (or codomain)
so there is a function
\[\operatorname{target}\colon\mathrm{mor}\to\mathrm{ob}.\] Now, we can
compose a morphism \(f\) and a morphism \(g\) to get \(fg\) if
\(\operatorname{target}(f) = \operatorname{source}(g)\), so we have a
composition function
\[\operatorname{composition}\colon C\to\mathrm{mor}\] defined only on
the subset \(C\) of \(\mathrm{mor}\times\mathrm{mor}\) that is the
biggest subset making the following diagram commute: \[
  \begin{tikzcd}
    C \ar[r,"p_1"]
      \ar[d,swap,"p_2"]
    &\mathrm{mor} \ar[d,"\operatorname{target}"]
  \\\mathrm{mor} \ar[r,"\operatorname{source}"]
    &\mathrm{ob}
  \end{tikzcd}
\] where \(p_1\colon(f,g)\mapsto f\) and \(p_2\colon(f,g)\mapsto g\).

Now category theorists have a slick way of dealing with these functions
defined only a subset satisfying equational conditions; instead of
talking about the ``biggest subset'' \(C\) they would say that \(C\) is
the ``limit'' of the diagram \[
  \begin{tikzcd}
    &\mathrm{mor} \ar[d,"\operatorname{target}"]
  \\\mathrm{mor} \ar[r,"\operatorname{source}"]
    &\mathrm{ob}
  \end{tikzcd}
\] If you don't get this, don't worry; in a sense it's just another way
of talking about the same thing, with the advantage of being infinitely
more general, since one can talk about the limit of any diagram, though
here we will only need limits of \emph{finite} diagrams.

Then, after having lined up these ingredients (and I have left some
out!), I could go ahead and say what equational laws they need to
satisfy, like associativity of composition; and if I wanted I could
write all these laws out using commutative diagrams, too! Then I would
have laid out the ``theory of categories'' --- a complete description of
the operations in a category and the laws they obey.

The theory of categories is a typical example of a ``finite limit
theory'', because what I really did above, in describing the ``theory of
categories'', is describe a \emph{category}, say \(\mathsf{Th}\), having
objects called \(\mathrm{ob}\) and \(\mathrm{mor}\), and morphisms
called \(\operatorname{source}\), \(\operatorname{target}\),
\(\operatorname{composition}\), and so on, such that various diagrams
commute! Moreover, we should think of \(\mathsf{Th}\) as a category with
all finite limits, that is, one in which all finite diagrams have
limits. That allows us to deal with things like the object \(C\) above,
which are defined as limits of finite diagrams.

So we have this thing \(\mathsf{Th}\), the ``theory of categories''. And
then, any \emph{particular} category is a ``model'' of this theory
\(\mathsf{Th}\). A ``model'' assigns to each object in \(\mathsf{Th}\) a
particular set --- for example, ``mor'' above gets assigned the set of
morphisms in some particular category \(\mathcal{C}\) --- and assigns to
each morphism in \(\mathsf{Th}\) a particular function --- for example,
``composition'' above gets assigned the function representing
composition in \(\mathcal{C}\). Moreover, this assignment satisfies a
bunch of utterly obvious consistency conditions which one summarizes by
saying that a ``model of the theory \(\mathsf{Th}\) is a functor from
\(\mathsf{Th}\) to \(\mathsf{Set}\) that preserves finite limits''. In
logic, you know, a model of a theory is something that assigns to each
thingie in the theory an actual thingie, in such a way that all the
stuff the theory \emph{says} is true about these thingies, \emph{is} true!

Now if you are with me thus far you either know this stuff better than I
do, or else I congratulate you, because the example I picked was
deliberately self-referential and confusing --- I was using category
theory to describe the theory of categories, and also, the theory
\(\mathsf{Th}\) itself was a category! But the world of thought does
have a funny way of wrapping back on itself like that... so it's
good to get used to it.

In fact there is a big literature on ``sets with essentially algebraic
structure'' and ``categories of models of finite limit
theories''... this is a branch of logic they never taught me about
in school, but it definitely exists, and Power gives some references to
it:

\begin{enumerate}
\def\labelenumi{\arabic{enumi})}
\setcounter{enumi}{2}
\item
  Pierre Gabriel and Friedrich Ulmer, \emph{Lokal pr\"asentierbare Kategorien}, 
  Lecture Notes in Mathematics \textbf{221}, Springer, Berlin, 1971.

  G.\ Max Kelly, ``Structures defined by finite limits in the enriched context
  I'', \emph{Cahiers de Top.\ et Geom.\ Diff.} \textbf{23} (1982), 3--41.

  Michael Makkai and Robert Pare, ``Accessible categories: the
  foundations of categorical model theory'', \emph{Contemp.\ Math.} 
  \textbf{104}, AMS, Providence, Rhode Island, 1989.
\end{enumerate}
\noindent
But let's dig in a bit further, since really the fun is just starting.
Now, I told you what a model of one of these finite limit theories Th
was, but not what a morphism between models is! Well, if a model is a
sort of functor, a morphism between them should be a sort of natural
transformation between functors; that's how it usually goes. So there is
really a category \(\mathsf{Mod}(\mathsf{Th})\) of models of one of
these theories \(\mathsf{Th}\). If \(\mathsf{Th}\) were the theory of
categories as above, \(\mathsf{Mod}(\mathsf{Th})\) would be the category
of (small) categories, which is called \(\mathsf{Cat}\). To take a less
fiendish example, if \(\mathsf{Th}\) were the theory of groups,
\(\mathsf{Mod}(\mathsf{Th})\) would the category \(\mathsf{Grp}\).

But now suppose one wanted to build a computer language that could not
only deal with all sorts of data types corresponding to different ``sets
with essentially algebraic structure'', but also various ``categories
with essentially algebraic structure''. For if one likes category theory
well enough to do a lot of computer science using it, it makes sense to
let the computer itself join the fun, by creating a language in which
it's easy to talk about categories. After all, our eventual goal with
computers is to have them completely replace computer scientists, right?

Well, in a way ``categories with essentially algebraic structure''
aren't terribly different from sets with essentially algebraic
structure. Roughly, the idea is that instead of having a data type that
consists of a bunch of sets with functions between them satisfying some
equational laws, we have a data type consisting of a bunch of
categories, functors between them, and natural transformations between
THEM satisfying equational laws. What this means is that if we try to
copy the above stuff, instead of a ``theory'' we will have a
``2-theory'' \(\mathsf{Th}\), which is some sort of \(2\)-category, and
then a model of this would be a 2-functor from \(\mathsf{Th}\) to
\(\mathsf{Cat}\). We want to wind up getting a \(2\)-category
\(\mathsf{Mod}(\mathsf{Th})\) of models of \(\mathsf{Th}\).

But actually carrying this out is a bit tricky, and much of Power's
paper goes into the details of various proposed schemes. Of course there
is no reason in principle to stop here, other than our limited
understanding of \(n\)-categories, sheer bewilderment, or boredom.
Reasoning about \(n\)-categories always tends to drag in
\((n+1)\)-categories, because the collection of all \(n\)-categories
with some particular structure (such as the ``essentially algebraic
structures'' I've focussed on here, but also other sorts) typically
forms an \((n+1)\)-category. This is how Power motivates tricategories.
Right now we are stuck at \(n = 3\), but there are good reasons to
expect that pretty soon we'll go beyond that. In fact, Power and Street
showed me a sketch of their ideas on tetracategories....

\hypertarget{week54}{%
\section{June 2, 1995}\label{week54}}

I just got back from a quantum gravity conference in Warsaw, and I'm
dying to talk about some of the stuff I heard there, but first I should
describe some work on topology and higher-dimensional algebra that I
have been meaning to discuss for some time now.

\begin{enumerate}
\def\labelenumi{\arabic{enumi})}
\tightlist
\item
  Timothy Porter, ``Abstract homotopy theory: the interaction of category
  theory and homotopy theory''.  Available at \href{https://ncatlab.org/nlab/files/Abstract-Homotopy.pdf}{https://ncatlab.org/nlab/files/Abstract-Homotopy.pdf}.
\end{enumerate}

Timothy Porter is another expert on higher-dimensional algebra whom I
met in Bangor, Wales, where he teaches. As paper 3) below makes clear,
he is very interested in the relationship between traditional themes in
topology and the new-fangled topological quantum field theories (TQFTs)
people have been coming up with these days. The above paper does not
mention TQFTs; instead, it is an overview of various approaches that
people have used to study homotopy theory in an algebraic way. But
anyone seriously interested in the intersection of physics and topology
would do well to get ahold of it, since it's a pleasant way to get
acquainted with some of the beautiful techniques algebraic topologists
have been developing, which many physicists are just starting to catch
up with.

What's homotopy theory? Well, roughly, it's the study of the properties
of spaces that are preserved by a wide class of stretchings and
squashings, called ``homotopies''.

For example, a closed disc \(D\) and a one-point set \(\{p\}\) are quite
different as topological spaces, in that there is no continuous map from
one to the other having a continuous inverse. (This is obvious because
they have a different number of points!) But there is clearly something
similar about them, because you can squash a disc down to a point
without crushing any holes in the process (since the disc has no holes).
To formalize this, note that we can find continuous functions
\[f\colon D\to\{p\}\] and \[g\colon\{p\}\to D\] that are inverses ``up
to homotopy''. For example, let \(f\) be the only possible function from
\(D\) to \(\{p\}\), taking every point in \(D\) to \(p\), and let \(g\)
be the map that sends \(p\) to the point \(0\), where we think of \(D\)
as the unit disc in the plane. Now if we first do \(g\) and then do
\(f\) we are back where we started from, so \(gf\) is the identity on
\(\{p\}\). But if we first do \(f\) and then \(g\) we are NOT
necessarily back where we started from: instead, the function \(fg\)
takes every point in \(D\) to the point \(0\) in \(D\). So \(fg\) is not
the identity. But it is ``homotopic'' to the identity, by which I mean
that there is a continuously varying family of continuous functions
\(F_t\) from \(D\) to itself, such that \(F_0 = fg\) and \(F_1\) is the
identity on \(D\). Simply let \(F_t\) be scalar multiplication by \(t\)!
As \(t\) goes from \(1\) to \(0\), we see that \(F_t\) squashes the disc
down to a point.

A bit more precisely, and more generally too, if we have two topological
spaces \(X\) and \(Y\) we say that two continuous functions
\(f,g\colon X \to Y\) are homotopic if there is a continuous function
\[F\colon[0,1]\times X\to Y\] such that \[F(0,x)=f(x)\] and
\[F(1,x) = g(x).\] Intuitively, this means that \(f\) can be
``continuously deformed'' into \(g\). Then we say that two spaces \(X\)
and \(Y\) are homotopic if there are continuous functions
\(f\colon X\to Y\), \(g\colon Y \to X\) which are inverse up to
homotopy, i.e., such that \(gf\) and \(fg\) are homotopic to the
identity on \(X\) and \(Y\), respectively.

The main goal in homotopy theory is to understand when functions are
homotopic and when spaces are homotopic. This is incredibly hard in
\emph{general}, but in special cases a huge amount is known. To take a
random (but important) example, people know that all maps from the
sphere to the circle are homotopic. Remember that algebraists call the
sphere \(S^2\) since its surface is \(2\)-dimensional, and call the
circle \(S^1\); in general the unit sphere in \(\mathbb{R}^{n+1}\) is
called \(S^n\). So for short, one says that all maps from \(S^2\) to
\(S^1\) are homotopic. But: there are infinitely many different
nonhomotopic maps from \(S^3\) to \(S^2\)! In fact there is a nice way
to label all these ``homotopy classes'' of maps by integers. And then:
there are only two homotopy classes of maps from \(S^4\) to \(S^3\).
There are also only two homotopy classes of maps from \(S^5\) to
\(S^4\), and from \(S^6\) to \(S^5\), and so on.

Now, the famous topologist J. H. C. Whitehead put forth an important
program in 1950, as follows: ``The ultimate aim of \emph{algebraic
homotopy} is to construct a purely algebraic theory, which is equivalent
to homotopy theory in the same way that `analytic' is equivalent to
`pure' projective geometry.'' Since then a lot of people have approached
this program from various angles, and Porter's paper tours some of the
key ideas involved.

Part of the reason for pursuing this program is simply to get good at
computing things, in a manner similar to how analytic geometry helps you
solve problems in ``pure'' geometry. This is not my main interest; if I
want to know how many homotopy classes of maps there are from \(S^9\) to
\(S^6\), or something, I know where to look it up, or whom to ask ---
which is infinitely more efficient than trying to figure it out myself!
And indeed, there is a formidable collection of tools out there for
solving various sorts of specific homotopy-theoretic problems, not all
of which rely crucially on a \emph{general} purely algebraic theory of
homotopy.

I'm more interested in this program for another reason, which is simply
to find an algebraic language for talking about things being true ``up
to homotopy''. As I've tried to explain in recent ``weeks'', there are
many situations where equations should be replaced by some weaker form
of equivalence. Taking this seriously leads naturally to the study of
\(n\)-categories, in which equations between \(j\)-morphisms can be
replaced by specified \((j+1)\)-morphisms. But Porter describes a host
of different (though related) formalisms set up to handle this sort of
issue. A few of the main ones are: simplicial sets, simplicial objects
in more general categories, Kan complexes, Quillen's ``model
categories'', \(\mathsf{Cat}^n\) groups, and homotopy coherent diagrams.
Understanding how all these formalisms are related and what they are
good for is quite a job, but this paper helps one get started.

So far everything I've actually said is quite elementary --- I've made
reference to some impressive buzzwords without explaining them, but that
doesn't count. So I should put in something for the folks who want more!
Let me say a word or two about \(\mathsf{Cat}^n\) groups. The definition
of these is a typical mind-blowing piece of higher-dimensional algebra,
so I can't resist explaining it. (After a while these definitions stop
seeming so mind-boggling, and then one is presumably beginning
understand the point of the subject!) In
\protect\hyperlink{week53}{``Week 53''} I gave a definition of a
category using category theory. This might seem completely circular and
useless, but of course I was illustrating quite generally how one could
define a ``model'' of a ``finite limit theory'' using category theory.
The idea was that a category is a \emph{set} of objects, a \emph{set} of
morphisms, together with various \emph{functions} like the source and
target functions which assign to any morphism (or ``arrow'') its source
and target (or ``tail'' and ``tip''). These sets and functions needed to
satisfy various axioms, of course.

Now \emph{sets} and \emph{functions} are the objects and morphisms in
the category of sets, which folks call Set. So in
\protect\hyperlink{week53}{``Week 53''} I cooked up a little category
\(\mathsf{Th}\) called ``the theory of categories'', which has objects
called ``\(\mathrm{ob}\)'' and ``\(\mathrm{mor}\)'', morphisms called
``\(s\)'' and ``\(t\)'', etc.  These were completely abstract gizmos,
not actual sets and functions. But we required them to satisfy the exact
same laws that the sets of objects and morphisms, and the source and
target functions, and so on, satisfy in an actual category. Then a
functor from \(\mathsf{Th}\) to \(\mathsf{Set}\) which preserves finite
limits is called a ``model'' of the theory of categories, because it
assigns to the completely abstract gizmos actual sets and functions
satisfying the same laws. In other words, if we have a functor
\[F\colon\mathsf{Th}\to\mathsf{Set}\] we have an actual set
\(F(\mathrm{ob})\) of objects, an actual set \(F(\mathrm{mor})\) of
morphisms, an actual function \(F(s)\) from \(F(\mathrm{ob})\) to
\(F(\mathrm{mor})\), and so on. In short, we have an actual category!

Now to get this trick to work we didn't need much to be true about the
category Set: all we needed was that it had finite limits. (Ignore this
technical stuff about limits if you don't get it; you can still get the
basic idea here.) And there are lots of categories with this property,
like the category \(\mathsf{Grp}\) of groups. So we can also talk about
a model of the theory of categories in the category of groups! What is
this? Well, it's just a functor from \(\mathsf{Th}\) to \(\mathsf{Grp}\)
that preserves finite limits. More concretely, it's exactly like a
category, except everywhere in the definition of category where you see
the word ``set'', scratch that out and write in ``group'', and
everywhere you see the word ``function'', scratch that out and write in
``homomorphism''. So you have a \emph{group} of objects, a \emph{group}
of morphisms, together with various \emph{homomorphisms} like the source
and target, and so on... satisfying laws perfectly analogous to
those in the definition of a category!

Folks call this kind of thing a ``categorical group'', a ``category
object in \(\mathsf{Grp}\)'' or an ``internal category in
\(\mathsf{Grp}\)''. From the point of view of sheer audacity alone, it's
a wonderful concept: we have taken the definition of a category and
transplanted it from the soil in which it was born, namely the category
\(\mathsf{Set}\), into new soil, namely the category \(\mathsf{Grp}\)!
But more remarkably still, the study of these ``categorical groups'' is
equivalent to the study of ``homotopy 2-types'' --- that is, topological
spaces, but where you only care about them up to homotopy, and even more
drastically, where nothing above dimension 2 concerns you. This result
is due (as far as I can tell) to Ronnie Brown and C. B. Spencer,
building on earlier work of Mac Lane and Whitehead.

But why stop here? Consider the category \(\mathsf{Cat}(\mathsf{Grp})\)
of these category objects in \(\mathsf{Grp}\). Take my word for it, such
a thing exists and it has finite limits. That means we can look for
models of the theory of categories in \(\mathsf{Cat}(\mathsf{Grp})\) ---
i.e., functors from \(\mathsf{Th}\) to \(\mathsf{Cat}(\mathsf{Grp})\),
preserving finite limits. In fact, \emph{there} things form a category,
say \(\mathsf{Cat}^2(\mathsf{Grp})\), and \emph{this} category has
finite limits, so we can look for models of the theory of categories in
\emph{this} category, and \emph{these} form a category
\(\mathsf{Cat}^3(\mathsf{Grp})\), which also has finite limits...
etc. So we can construct an insanely recursive hierarchy:

\begin{itemize}
\tightlist
\item
  groups
\item
  category objects in the the category of groups
\item
  category objects in the category of (category objects in the category
  of groups)
\item
  etc....
\end{itemize}
\noindent
Now, truly wonderfully, Jean-Louis Loday showed that the study of
\(\mathsf{Cat}^n(\mathsf{Grp})\) is equivalent (in a certain precise
sense) to the study of homotopy \(n\)-types --- i.e., homotopy theory
where you don't care about phenomena above dimension \(n\):

\begin{enumerate}
\def\labelenumi{\arabic{enumi})}
\setcounter{enumi}{1}
\tightlist
\item
  Jean-Louis Loday, ``Spaces with finitely many non-trivial homotopy groups'',
  \emph{Jour. Pure Appl. Algebra} \textbf{24} (1982), 179--202.
\end{enumerate}
\noindent
Subsequently, Ronnie Brown, Loday and others have done some interesting
topology using this fact. But the most remarkable thing, in a way, is
how taking a perfectly basic concept, the concept of \emph{group}, and then
doing category theory ``internally'' in the category of groups in an
iterated fashion, winds up being very much the same as doing topology ---
or at least homotopy theory. This suggests that there is something
deeply algebraic about homotopy theory in the first place.

\begin{enumerate}
\def\labelenumi{\arabic{enumi})}
\setcounter{enumi}{2}
\tightlist
\item
  Timothy Porter, ``Interpretations of Yetter's notion of \(G\)-coloring:
  simplicial fibre bundles and non-abelian cohomology'', available at
  \url{https://citeseer.ist.psu.edu/viewdoc/summary?doi=10.1.1.50.3718}.
\end{enumerate}
\noindent
Physicists know and love the Dijkgraaf--Witten model, a 2+1-dimensional
TQFT that depends on a finite group \(G\). It's easy to compute the
Hilbert space of states for any compact oriented 2-manifold in this
TQFT. Just triangulate your 2-manifold and let your Hilbert space have
as a basis the set of all possible ways of labelling the edges with
elements of \(G\) such that \(g_1g_2g_3 = 1\) whenever we have 3 edges
going counterclockwise around any triangle. Yetter figured out how to
generalize this to get an interesting TQFT from any finite categorical
group:

\begin{enumerate}
\def\labelenumi{\arabic{enumi})}
\setcounter{enumi}{3}
\item
  David N. Yetter, ``Topological quantum field theories associated to
  finite groups and crossed G-sets'', \emph{Journal of Knot Theory and
  its Ramifications} \textbf{1} (1992), 1--20.

  ``TQFTs from homotopy 2-types'', \emph{Journal of Knot Theory and its
  Ramifications} \textbf{2} (1993), 113--123.
\end{enumerate}
\noindent
This should be the beginning of some bigger pattern relating homotopy
theory and TQFTs. Jim Dolan and I have our own theories as to how this
pattern should work (see \protect\hyperlink{week49}{``Week 49''}) but
they are still a rather long ways from being theorems. Porter, who is an
expert in simplicial methods, makes the relationship (or \emph{one} of the
relationships) very clear in this special case.

\begin{enumerate}
\def\labelenumi{\arabic{enumi})}
\setcounter{enumi}{4}
\item
  Justin Roberts, ``Skein theory and Turaev--Viro invariants'', \emph{Topology},
  \textbf{34}, 771--787.

  ``Refined state-sum invariants of 3- and 4-manifolds'', in \emph{Geometric Topology 
    (Athens, GA, 1993)}, ed.\ Will Kazez, AMS/IP Stud. Adv. Math \textbf{2},
   1997, pp.\ 217--234.

  ``Skeins and mapping class groups'', \emph{Math.\ Proc.\ Camb.\ Phil.\
  Soc.} \textbf{115} (1994), 53--77.

  Gregor Masbaum and Justin Roberts, ``On central extensions of mapping
  class groups'', \emph{Math.\ Annalen} \textbf{302} (1995), 131--150.
\end{enumerate}
\noindent
The first two papers here might be the most interesting for physicists.
They both deal with 3d and 4d TQFTs constructed using quantum
\(\mathrm{SU}(2)\): in particular, the Turaev--Viro theory in dimension
3, and the Crane--Yetter--Broda theory in dimension 4. The first theory is
interesting physically because it corresponds to 3d Euclidean quantum
gravity with cosmological constant. The second theory is interesting
mainly because it's one of the few 4d TQFTs for which the Atiyah axioms
have been shown; sometime I will explain the Lagrangian for this theory,
which seems not to be well-known.

In Roberts' first paper, which was already discussed in
\protect\hyperlink{week14}{``Week 14''}, he gave a simple proof that the
partition function for the Turaev--Viro theory was the absolute value
squared of that for Chern--Simons theory, perhaps the most famous of
TQFTs. He also showed that the partition function for the
Crane--Yetter--Broda theory was a function of the signature and Euler
characteristic (classical invariants of 4-manifolds). In the second
paper, he considers observables for the Turaev--Viro and Crane--Yetter--Broda
theories depending on
cohomology classes in the manifold. I wish I had energy now to explain a
bit more about these observables, since they are very interesting, but I
don't!

\begin{enumerate}
\def\labelenumi{\arabic{enumi})}
\setcounter{enumi}{5}
\tightlist
\item
  Lawrence Breen, \emph{On the Classification of 2-Gerbes and 2-Stacks},
  \emph{Asterisque} \textbf{225}, 1994.
\end{enumerate}
\noindent
Suffice it to say that if gerbes and stacks --- which are, very roughly,
sort of like sheaves of categories --- are too simple to interest you,
it's time to read about 2-gerbes and 2-stacks --- which are roughly like
sheaves of \(2\)-categories.

\hypertarget{week55}{%
\section{June 4, 1995}\label{week55}}

I recently went to a workshop on canonical quantum gravity in Warsaw,
organized by Jerzy Kijowski and Jerzy Lewandowski, and I learned some
interesting things. I'll talk about some of them in this issue, and some
in the next.

Conferences are a funny thing. On science newsgroups on the net, there
is very little talk about conferences. This is probably because the
people who really understand conferences are too busy flying from one
conference to the next to post to newsgroups very often. Academic
success is in part measured by the number of conference invitations one
receives, the prestige of the conferences, and the type of invitation.
For example, a big plenary lecture on an impressive stage, preceded by a
little warmup where someone explains how great you are, counts for
infinitely many talks in those parallel sessions where dozens of people
get 10 minutes each to explain their work before the moderator begins to
make little coughs indicating that it's time for the next one, while all
the while people drift in and out in a feeble attempt to find the really
interesting talks. Still, giving any sort of talk is regarded as better
than giving none, so academics spend a lot of time doing this sort of
thing.

One of the great dangers of being a successful academic, in fact, is
that one may get invited to so many conferences that one never has time
to think. Winning the Nobel prize is purported to be the kiss of death
in this respect. Of course, it's a universal platitude that the real
thinking at conferences gets done not during the talks, but informally
in small groups. But the funny thing is that at most conferences people
are so worn out after going to a day's worth of talks that they have
limited energy for serious conversation afterwards: they usually seem
more interested in finding the good local restaurants and scenic
attractions. If people could have conferences with no lectures
whatsoever, or maybe one a day, it would probably be more productive.
But the idea that a bunch of people could figure something out just by
sitting around and chatting informally is absolutely foreign to our
conception of ``work''. People expect to receive money from bureaucrats
to go to conferences, but to convince a bureaucrat that you are deserve
the money, you need to give a lecture, so of course all conferences have
too many lectures.

Turning back towards Warsaw, a city with a marvelous mathematical
history, I am reminded of Gian-Carlo Rota's biographical sketch of
Stanislaw Ulam, in which (as a master of irony) he talks about how lazy
Ulam was: all he wanted to do was sit around in cafes and come up with
interesting conjectures and research programs, and leave it to others to
work them out. And this in turn reminds me of the Scottish Cafe, where
Polish mathematicians used to hang out and write on the tablecloths,
until the owner provided them with a notebook, in which many famous
conjectures were formulated, and I believe prizes like bottles of wine
were offered for their solutions. Was the Scottish Cafe in Warsaw?
{[}No, Lwow.{]} Does it still exist? I completely forgot to check while
I was there. The Banach Center, in which the conference participants
stayed, comes from a later stratum of Polish mathematical history; it
was built after the war, and one room still contains a portrait of
Lenin. I know that because a film crew used it to shoot a scene for a
historical movie!

Anyway, I enjoyed this conference in Warsaw quite a bit, because a lot
of people working on the loop representation of quantum gravity were
there, and I managed to have a fair number of serious conversations.
Before going into what I learned there, I should say that I just found a
fun thing for people to read who are interested in quantum gravity, but
are not necessarily specialists:

\begin{enumerate}
\def\labelenumi{\arabic{enumi})}
\tightlist
\item
  Gary Au, ``The quest for quantum gravity'', available as
  \href{https://arxiv.org/abs/gr-qc/9506001}{\texttt{gr-qc/9506001}}.
\end{enumerate}

This consists mainly of interviews with Chris Isham, Abhay Ashtekar and
Edward Witten. What's nice is that the interviews are conducted by
someone who knows physics. The questions and answers are technical
enough to convey some of the real substance of the subject, while still
(I hope) non-technical enough so that you don't have to be an expert to
get a lot out of them. Isham talks mainly about the ``problem of time''
in quantum gravity, Ashtekar talks mainly about the loop representation
of quantum gravity, and Witten talks about string theory.

Anyway, Ashtekar and a bunch of other good people were at this Warsaw
conference, which is why I went. The main topics of conversation were
spin networks and their use in studying the area and volume operators in
quantum gravity. As I explained earlier in
\protect\hyperlink{week43}{``Week 43''}, one may very roughly think of a
spin network as a graph whose edges are labelled with ``spins''
\(0\),\(1/2\),\(1\),\(3/2\), and so on, and whose vertices are labelled
with certain gadgets called ``intertwining operators'' (which in the
simplest case are just the Clebsch-Gordon coefficients you learn about
when studying angular momentum in quantum mechanics). Penrose introduced
these as abstract graphs (see \protect\hyperlink{week22}{``Week 22''}
and \protect\hyperlink{week41}{``Week 41''}), as a kind of substitute
for thinking of space as a manifold, but more recently Rovelli and
Smolin started thinking of them as graphs embedded into 3d space, and
saw that these were a really natural way to describe states of quantum
gravity: even better than loops, because they form an orthonormal basis!
Actually, it was mainly me who proved in a really rigorous way that they
form an orthonormal basis, but Rovelli and Smolin had already been doing
calculations using this idea for a while. One thing they computed was
the eigenvalues of the observables in quantum gravity corresponding to
the area of a surface in space, or the volume of a region.

Now there are all sorts of technical caveats and subtleties that I don't
want to get into here, but in a really rough sort of sense, what their
answers suggest is that IF the loop representation of quantum gravity is
right, and we are on the right track about how it works, then the area
of surfaces comes in certain (not integer, but discrete) multiples of
the Planck length squared, and the volume of regions comes in multiples
of the Planck length cubed. Note: that was a big ``IF''. This is
especially interesting because it doesn't arise by assuming from the
start that spacetime has a discrete structure. In fact, their
computations assume spacetime is a continuous manifold. Nonetheless this
discreteness pops out. It's not completely surprising: after all,
Schrodinger's equation for the hydrogen atom is a perfectly
``continuous'' sort of thing, a partial differential equation, but the
energy of the bound states winds up being a discrete sort of thing.
Still, it's sort of exciting and new.

An interesting thing happened at the conference. Renate Loll, who works
on the loop representation of gauge theories and also lattice gauge
theory, has recently developed a lattice formulation of quantum gravity
closely modelled after the loop representation:

\begin{enumerate}
\def\labelenumi{\arabic{enumi})}
\setcounter{enumi}{1}
\tightlist
\item
  Renate Loll, ``Nonperturbative solutions for lattice quantum
  gravity'' \emph{Nucl.\ Phys.\ B} \textbf{444} (1995), 619--639.  Also available as
  \href{https://arxiv.org/abs/gr-qc/9502006}{\texttt{gr-qc/9502006}}.
\end{enumerate}
\noindent
This has the wonderful feature that it's perfectly rigorous and also one
can start using computers to start calculating things with it. For
example, the most subtle aspect of the loop representation of quantum
gravity is the Wheeler--DeWitt equation \[H\psi=0\] where \(H\) is an
operator called the ``Hamiltonian constraint''. More on this later; my
point here is just that physical states of quantum gravity need to
satisfy this equation. Getting \(H\) to be well-defined is tricky when
space is a continuum, but in Loll's lattice version of theory (which is
an approximation to the full continuum theory) she has already done
this, so one can now start trying numerically to find solutions and see
what they look like. She has also found some explicit solutions.

\emph{Also}, she did some work on the volume operator in her lattice
approach, and came up with a result in contradiction to Rovelli and
Smolin's paper on the subject (cited in
\protect\hyperlink{week43}{``Week 43''}). They had said that states
corresponding to trivalent spin networks --- spin networks with only 3
edges at each vertex --- could have nonzero volume. But using her
version of the theory she computed that trivalent states --- states with
only 3 nonzero spins at the edges of the lattice incident to any vertex
--- all had zero volume, and that she needed at least 4 nonzero spins to
get volume! The volume operator, in case you're wondering, acts as a
certain sum over vertices: each one winds up contributing a certain
finite amount of volume, which the theory allows you to compute.

This led to a whole lot of discussion and scribbling on the blackboards
of the Banach center. I found it truly delightful to see all these
physicists drawing pictures of spin networks and doing graphical
computations just the way a knot theorist like Kauffman does all the
time. It was as if the universe had this spin network aspect to it, and
everyone was finally starting to catch on. Either that or mass delusion!
I hadn't quite gotten the hang of how to compute these volume operators
before, so it was a great chance to learn: one person would do a
computation, then someone else would do it a different way and get a
different answer, then someone else would do it yet another way and get
yet another answer, and so on, so you could ask lots of questions
without seeming too dumb. Even I did a computation after a while, and
got zero volume for at least a certain class of trivalent vertices. The
votes in favor of trivalent vertices having zero volume kept piling up.
Finally Smolin noticed that he and Rovelli had made a sign mistake. This
is incredibly easy to do, since there are lots of tricky sign
conventions in spin network theory. Fundamentally these are due to the
fact that spin-\(1/2\) particles are fermions... but I don't think
people fully understand the physical implications of this. (There is
also a marvelous category-theoretic explanation of it, but I fear that
if I go into that all the physicists will stop reading. Maybe some other
time.) Rovelli and Smolin got pretty depressed about this for a while,
but I tried to reassure them that only people who write really
interesting papers ever get anybody to find the mistakes.

So perhaps we know a little more about the meaning of volume in a
quantum theory of spacetime.

Spin networks are very beautiful and simple things. To learn about them,
in addition to the various papers listed in the ``weeks'' above, one can
now turn to Rovelli and Smolin's paper:

\begin{enumerate}
\def\labelenumi{\arabic{enumi})}
\setcounter{enumi}{2}
\tightlist
\item
  Carlo Rovelli and Lee Smolin, ``Spin networks in quantum gravity'',
  \emph{Phys.\ Rev.\ D} \textbf{52} (1995), 5743--5759.
   Also available as 
   \href{https://arxiv.org/abs/gr-qc/9505006}{\texttt{gr-qc/9505006}}.
\end{enumerate}
\noindent
If you are more of a mathematician, or less of an expert on quantum
gravity, you might also try a review article I wrote about them, which
starts with a quick summary of what the heck canonical quantum gravity
is about, why it's hard to do, and why the loop representation seems to
help:

\begin{enumerate}
\def\labelenumi{\arabic{enumi})}
\setcounter{enumi}{3}
\tightlist
\item
  John Baez, ``Spin networks in nonperturbative canonical quantum
  gravity'', in \emph{The Interface of Knots and Physics}, ed.\ Louis Kauffman, 
  AMS, Providence, 1996, pp.\ 167--203.   Also available as 
  \href{https://arxiv.org/abs/gr-qc/9504036}{\texttt{gr-qc/9504036}}.
\end{enumerate}

Now so far I have been trying to make things sound simple, but here I
should point out that when one talks about ``states of quantum gravity''
there are at least three quite different things one might mean. This is
because the loop representation follows Dirac's general philosophy of
quantizing systems with constraints, with some extra twists here and
there. As I've repeatedly explained
(e.g.~\protect\hyperlink{week43}{``Week 43''}), Einstein's equation for
general relativity has 10 components, and if you split spacetime up into
space and time (more or less arbitrarily --- there's no ``best'' way), 4
of these can be seen as constraints that the metric on space and its
first time derivative must satisfy (at any given time), while the
remaining 6 describe how the metric on space evolves in time (which
makes sense, because the metric has 6 components). When you follow
Dirac's procedure for quantizing the equations what you do is this.
First you forget about the constraint and get a big space of states, the
``kinematical state space''. There are lots of mathematical choices
involved here, but Ashtekar and Lewandowski came up with a particular
nice way of doing this rigorously, and one calls this space of states
``\(L^2\) of the space of \(\mathrm{SU}(2)\) connections modulo gauge
transformations with respect to the Ashtekar-Lewandowski generalized
measure''. Spin networks form an orthonormal basis of this Hilbert
space. All the stuff about area and volume operators above refers to
operators on this space.

Then, however, you need to deal with the constraints. Now 3 of the 4
constraints simply amount to requiring that your states be invariant
under all diffeomorphisms of space, so these are usually dealt with
first, and called the ``diffeomorphism constraint''. Imposing these
constraints is a bit tricky; naively one would first guess that this
``diffeomorphism-invariant state space'' is just a subspace of the
original kinematical state space, but actually it's not quite so simple.
In any event, there are also spin network states at the
diffeomorphism-invariant level, corresponding not to \emph{particular}
embeddings of graphs in space, but to diffeomorphism equivalence classes
thereof. This again has been used by Rovelli, Smolin and others for a
while now, but it was first rigorously shown in the following paper:

\begin{enumerate}
\def\labelenumi{\arabic{enumi})}
\setcounter{enumi}{4}
\tightlist
\item
  Abhay Ashtekar, Jerzy Lewandowski, Don Marolf, Jose Mourao, and Thomas
  Thiemann, ``Quantization of diffeomorphism invariant theories of
  connections with local degrees of freedom'', \emph{Jour.\ Math.\ Phys.} \textbf{36}
  (1995), 6456--6493. Also available as
  \href{https://arxiv.org/abs/gr-qc/9504018}{\texttt{gr-qc/9504018}}.
\end{enumerate}

\noindent
This paper is nice in part because it doesn't assume you already have
read every previous paper about this stuff; instead, it describes the
general plan of the loop representation before constructing the
diffeomorphism-invariant spin network states. Also, buried in an
appendix somewhere, it gives nice conceptual formulas for the area and
volume operators, which serve as a complement to Rovelli and Smolin's
explicit computations of their matrix elements in terms of the spin
network basis.

Anyway, after taking care of the diffeomorphism constraint, one finally
needs to take care of the Hamiltonian constraint, meaning one needs to
find states satisfying the Wheeler--DeWitt equation. This is the hardest
thing to make rigorous, and the most exciting aspect of the whole
subject, because it expresses the fact that ``physical states'' of
quantum gravity are invariant under diffeomorphisms of space-TIME, not
just space. There is much more to say about this, but I won't go into it
here.

Now besides Loll and Rovelli and Smolin, all the authors of the above
paper except Mourao were at the conference in Warsaw, so there was a
large contingent of spin network fans around, not even counting some
other folks whose work I will get to in a while. This is why I was so
eager to go there, especially because my own talk was on a rather
esoteric subject which only these experts could be expected to give a
darn about. Namely...

The breakthrough of Ashtekar and Lewandowski, when it came to making the
loop representation rigorous, involved working with piecewise
real-analytic loops rather than piecewise smooth loops. (Actually
Penrose suggested this idea.) This is because piecewise smooth loops can
intersect in crazy ways, like in a Cantor set, which nobody could figure
out how to handle. But the price of this breakthrough was that one had
to assume the 3-manifold representing space was real-analytic, and
things then only work nicely for real-analytic diffeomorphisms, as
opposed to smooth ones. This always bugged me, so I have been working
away for about a year trying to deal with smooth loops, and finally I
got smart and teamed up with Steve Sawin, and we recently figured out
how to get things to work with smooth loops (at least a bunch of things,
like the Ashtekar-Lewandowski generalized measure). Our paper will be
out pretty soon, but for now anyone who wants a taste of the
mathematical technology involved should look at:

\begin{enumerate}
\def\labelenumi{\arabic{enumi})}
\setcounter{enumi}{5}
\tightlist
\item
  Steve Sawin, ``Path integration in two-dimensional topological quantum
  field theory'',  \emph{Jour.\ Math.\ Phys.} \textbf{36}
  (1995), 6130--6136.  Also available as \href{http://arxiv.org/abs/gr-qc/9505040}{\texttt{gr-qc/9505040}}.
\end{enumerate}

Loop representation ideas are applicable not only to canonical quantum
gravity but also to path integrals in gauge theory, because in both
cases one is doing integrals over a space of connections mod gauge
transformations. Here Sawin uses them to give a rigorous formulation of
2d TQFTs in terms of path integrals. There aren't many unitary 2d TQFTs,
and all of them are isomorphic to \(2\)-dimensional quantum gravity with
the usual Einstein-Hilbert action, with different values of the coupling
constant, or else direct sums of such theories.

Next ``week'' I'll talk about cool new idea Smolin has about TQFTs,
quantum gravity, and Bekenstein's bound on the entropy of a physical
system in terms of its surface area.

\hypertarget{week56}{%
\section{June 16, 1995}\label{week56}}

I got a copy of the following paper when I showed up in Warsaw:

\begin{enumerate}
\def\labelenumi{\arabic{enumi})}
\tightlist
\item
  Lee Smolin, ``Linking topological quantum field theory and
  nonperturbative quantum gravity'', available as
  \href{https://arxiv.org/abs/gr-qc/9505028}{\texttt{gr-qc/9505028}}.
\end{enumerate}

\noindent
and then I spent a fair amount of time reading it and thinking about it
throughout the conference. If the big hypothesis formulated in this
paper is correct, I think we are on the verge of having a really
beautiful theory of \(4\)-dimensional quantum gravity, at least given
certain boundary conditions. Mind you, I just mean a really beautiful
theory, not necessarily a physically correct theory. But beautiful
theories have a certain tendency to be right, or at least close, so let
me explain this hypothesis.

First of all, we have to remember that Ashtekar reformulated Einstein's
equation so that the configuration space for general relativity on the
spacetime \(\mathbb{R}\times S\), instead of being the space of
\emph{metrics} on a 3-manifold \(S\), is a space of \emph{connections}
on \(S\). A connection is just what a physicist often calls a vector
potential, but for any old gauge theory, not just electromagnetism.
Different gauge theories have different gauge groups, so I had better
tell you the gauge group of Ashtekar's version of general relativity:
it's \(\mathrm{SL}(2,\mathbb{C})\), the group of \(2\times2\) complex
matrices with determinant equal to \(1\). And I should probably tell you
which bundle over \(S\) we have an \(\mathrm{SL}(2,\mathbb{C})\)
connection on... but luckily, all \(\mathrm{SL}(2,\mathbb{C})\)
bundles over 3-manifolds are trivial, so I can cut corners by saying
it's the trivial bundle. We can think of a connection \(A\) on the
trivial \(\mathrm{SL}(2,\mathbb{C})\) bundle over \(S\) as \(1\)-forms
taking values in the Lie algebra \(\mathfrak{sl}(2,\mathbb{C})\),
consisting of \(2\times2\) complex matrices with trace zero.

Okay, so naively you might think a state in the \emph{quantum} version
of general relativity a la Ashtekar is just a wavefunction \(\psi(A)\).
That's not too far wrong and I won't bother about certain nitpicky
technicalities here (again, for the full story try my paper 
``\href{https://arxiv.org/abs/gr-qc/9504036}{Spin networks
in nonperturbative quantum gravity}''). But
there's one very important catch I can't ignore: general relativity has
\emph{constraint} equations, meaning that \(\psi\) has to satisfy some
equations. The first constraint, the Gauss law, just says that we must
have \[\psi(A) = \psi(A')\] whenever \(A'\) is the result of doing a
gauge transformation to \(A\). Or at the very least, this should hold up
to a phase; the point is that \(\psi\) is only supposed to record
physically significant information about the state of the universe, and
two connections are physically equivalent if they differ by a gauge
transformation. The second constraint, the diffeomorphism constraint,
says we need to have \[\psi(A) = \psi(A')\] when \(A'\) is the result of
applying a diffeomorphism of space, \(S\), to \(A\). Again, the point is
that two solutions of general relativity are physically the same if they
differ only by a coordinate transformation, or --- \emph{roughly} the
same thing --- a diffeomorphism. The third constraint is the real
killer. It's meaning is that \(\psi(A)\) doesn't change when we do a
diffeomorphism of space\emph{time} to the connection \(A\), but it's usually
formulated `infinitesimally' as the Wheeler--DeWitt equation
\[H \psi = 0\] meaning roughly that the time derivative of \(\psi\) is
zero. Think of it as a screwy quantum gravity version of Schrodinger's
equation, where \(d\psi/dt\) had better be zero!

It's hard to find explicit solutions of these equations. Indeed, it's
hard to know what the heck these equations \emph{mean} in a sufficiently
precise way to recognize a solution if we found one! However, things
were even worse back in the old days. Back in the old days when we
thought of states as wavefunctions on the space of metrics, we didn't
know \emph{any} solutions of these equations. But nowadays we are very happy,
because we know infinitely many times as many solutions! To be precise,
we now know \emph{one} solution. This is called the Chern--Simons state, and it
was discovered by Kodama:

\begin{enumerate}
\def\labelenumi{\arabic{enumi})}
\setcounter{enumi}{1}
\tightlist
\item
  H. Kodama, ``Holomorphic wavefunction of the universe'', \emph{Phys.\
  Rev.\ D} \textbf{42} (1990), 2548--2565.
\end{enumerate}

\noindent
Now actually people have proposed other explicit solutions, but
personally I have certain worries about all those other solutions, and I
am not alone in this, whereas everyone seems to agree that, no matter
how you slice it, the Chern--Simons state is a solution.

Now there is a slight catch: the Chern--Simons state is a solution of
quantum gravity \emph{with cosmological constant}. This is an extra term
that Einstein threw into his equations so that they wouldn't make the
obviously ridiculous prediction that the universe is expanding. When
Hubble took a look and saw galactic redshifts all over, Einstein called
this extra term the biggest blunder in his life. That kind of remark,
coming from that kind of person, might make you a little bit reluctant
to get too excited about having found a state of quantum gravity with
this extra term thrown in! Luckily it turns out that you can take the
limit as the cosmological constant goes to zero, and get a state of the
theory where the cosmological constant is zero. I like to call this the
``flat state'', because it's zero except where the connection \(A\) is
flat.

(In fact, if the space \(S\) is not simply connected, there are lots of
different flat states, because there is what experts call a moduli space
of flat connections, i.e., lots of different flat connections modulo
gauge transformations. Not many people talk too much about these flat
states, but I do so in my paper
``\href{https://arxiv.org/abs/gr-qc/9504036}{Spin networks
in nonperturbative quantum gravity}'' and also
the harder one ``\href{https://arxiv.org/abs/gr-qc/9410018}{Knots and quantum gravity: progress and prospects}''.)

Now what is this Chern--Simons state? Well, there is a wonderful thing
you can compute from a connection \(A\) on a (compact oriented)
3-manifold \(S\), called the Chern--Simons action:
\[CS(A) = \int_S \operatorname{tr}(A \wedge dA + (2/3)A \wedge A \wedge A)\]
which looks weird when you first see it, but gradually starts seeming
very sensible and nice. The reason why folks like it is that it doesn't
change when you do a small gauge transformation --- i.e., one you can
get to following a continuous path from the identity --- and it changes
only by an integral multiple of \(8\pi^2\) if you do a large gauge
transformation. Plus, it's diffeomorphism-invariant. It's incredibly
hard to write down many functions of \(A\) with these properties, so
they are precious. There are deeper reasons why it's so nice, but let's
leave it at that for now.

So, the Chern--Simons state is \[\psi(A) = \exp(-6 CS(A)/\Lambda)\] where
\(\Lambda\) is the cosmological constant. Don't worry about the factor
of 6 too much; depending on how you set up various things you might get
different numbers, and I can never keep certain factors of 2 straight in
this particular calculation. Note however that it looks as if things go
completely haywire as \(\Lambda\) approaches zero, which is why my
earlier remark about the `flat state' is a bit nontrivial.

Why does this satisfy the constraints? Well, I just said above that the
Chern--Simons action was hand-tailored to have the gauge-invariance and
diffeomorphism-invariance we want, so the only big surprise is that we
\emph{also} have a solution of the Wheeler--DeWitt equation. Well, we do:
a two-line computation shows it.

But clearly nature, or at least the goddess of mathematics, is trying to
tell us something if this Chern--Simons state, which has all sorts of
wonderful properties relating to \emph{3-dimensional} geometry, is also
a solution of the Wheeler--DeWitt equation, which is all about
\emph{4-dimensional} geometry, since it expresses the invariance of
\(\psi\) under evolution in \emph{time}. I have been thinking about this for
several years now and I think I finally really understand it. There are
probably people out there to whom it's perfectly obvious, because it's
not really all that complicated, but unfortunately none of these people
has ever explained it. Let me briefly say, for those who know about such
things, that it all comes down to the fact that the Chern--Simons action
was \emph{born} as a \(3\)-dimensional spinoff of a \(4\)-dimensional
thing called the 2nd Chern class. (If you want more details, I explained
this as well as I could at the time in ``\href{https://arxiv.org/abs/gr-qc/9410018}{Knots 
and quantum gravity: progress and prospects}''.)

What is the physical meaning of the Chern--Simons state? As far as I know
Kodama's paper hasn't been vastly surpassed in explaining this. He shows
that in the classical limit this state becomes something called the
anti-deSitter universe, a solution of Einstein's equation describing a
(roughly) exponentially expanding universe. If you are wondering what
this has to do with Einstein's introduction of the constant to \emph{keep} the
universe from expanding, let me just say this. In our current big bang
theory the universe is expanding, but the presence of matter, or any
sort of positive energy density, tends to pull it back in, and if there
is enough matter it'll collapse again. Einstein's stuck in a
cosmological constant term to give the vacuum some negative energy
density, exactly enough to counteract the positive energy density of
matter, so things would neither collapse nor expand, but instead remain
in an (unstable, alas) equilibrium. In the deSitter universe there's no
matter, just a cosmological constant of the opposite sign, so that the
vacuum has positive energy density. In the anti-deSitter universe
(invented by deSitter's nemesis anti-deSitter) there's no matter either,
but the cosmological constant has the sign giving the vacuum negative
energy density, which pushes the universe to keep expanding faster and
faster.

Now in addition to this physical interpretation, the Chern--Simons state
also has some remarkable relationships to knot theory, which were
discovered by Witten and, since then, studied intensively by lots of
people. I have written a lot in This Week's Finds about this! But
briefly, there should be an invariant of knots and links associated to
any state of quantum gravity, and the one associated to the Chern--Simons
state is the Kauffman bracket, a close relative of the Jones polynomial,
which is distinguished by having a very simple, beautiful definition,
and also lots of wonderful relationships to an algebraic structure, the
quantum group \(SU_q(2)\). I should add that in addition to an invariant
of knots and links, a state of quantum gravity should also give an
invariant of \emph{spin networks}, and indeed the Kauffman bracket
extends to a wonderful invariant of spin networks. One can read about
this in many places, but perhaps the most detailed account is Kauffman
and Lins' book referred to in \protect\hyperlink{week30}{``Week 30''}.

So the question arises: are all these wonderful features of the
Chern--Simons state of quantum gravity very special things that tell us
very little about quantum gravity in general, or are they important
clues that, if we understood them, would reveal a lot about the nature
of \emph{all} states of quantum gravity?

Of course, everyone who has fallen in love with the beauty of
Chern--Simons theory would \emph{like} the answer to be the latter. Louis
Crane, for example, is deeply convinced that Chern--Simons theory is
indeed the key to the whole business. He has written many papers on the
subject, most of which I've gone over in earlier Finds, and also one
brand new one, which is actually a very readable introduction to the
grand scheme he has in mind:

\begin{enumerate}
\def\labelenumi{\arabic{enumi})}
\setcounter{enumi}{2}
\tightlist
\item
  Louis Crane: ``Clock and category: is quantum gravity algebraic?'', \emph{Jour.\
  Math.\ Phys.\ }\textbf{36} (1995), 6180--6193.  Also available as
  \href{https://arxiv.org/abs/gr-qc/9504038}{\texttt{gr-qc/9504038}}.
\end{enumerate}
\noindent
This grand scheme involves a thorough refashioning of quantum gravity in
terms of category theory, and uses some of the very beautiful
mathematics of \(n\)-categories, but neglecting all these subtleties,
let us simply say that he argues that if the universe is \emph{in} the
Chern--Simons state, there is no need to understand any other states! He
doesn't really think all there is in the universe is gravity, of course,
so he envisages a souped-up theory containing other forces and
particles, but he argues that a generalization of quantum gravity to
include all these other forces and particles will still have a special
state, and that that's the state of the universe.

Being a conservative fellow myself, I prefer to remain agnostic on this
issue, but I did write a paper showing how, if you assumed that space,
the manifold above I called \(S\), is a \(3\)-dimensional sphere --- a
so-called \(S^3\) --- then if quantum gravity was in the Chern--Simons
state, there would be very nice Hilbert spaces of ``relative states'' on
each ``half'' of space. The relative state notion is often associated
with Everett, who made a big deal out of the fact that, even if a
two-part system was in a single state (a ``pure state'' in the language
of quantum theory), each part could be regarded as being in a
probabilistic mixture of lots of states (a ``mixed state''). Whether or
not you like the ``many-worlds interpretation'' of quantum theory which
Everett's thesis spawned, it is true that a single pure state on a
two-part system specifies a whole \emph{space} of states on each half.
So my idea was to take \(S^3\), arbitrarily split it into two balls
(\(D^3\)'s in math jargon), and work out the space of states on each
half. If you want to wax rhapsodic of this you can call one half the
``observer'' and the other the ``observed'', though it's crucial that
there is a symmetry interchanging the two --- there's not any way to
tell them apart.

On the more technical side, there is a lot of nice (though well-understood) 
knot theory involved. For example, a special property of the
quantum group \(\mathrm{SU}_q(2)\) --- called the ``positivity of the Markov
trace'', and discovered by Vaughan Jones of Jones polynomial fame --- equips each
space of states with an inner product, even in this situation where we
have no idea of an inner product to begin with. This paper is:

\begin{enumerate}
\def\labelenumi{\arabic{enumi})}
\setcounter{enumi}{3}
\tightlist
\item
  John Baez, ``Quantum gravity and the algebra of tangles'', 
  \emph{Class.\ Quant.\ Grav.}\textbf{10} (1993), 673--694.  Also available
  (without the all-important pictures!) as \href{https://arxiv.org/abs/hep-th/9205007}
  {\texttt{hep-th/9205007}}.
\end{enumerate}

So what has Lee Smolin done? Actually I have spent so much time leading
up to it that now I'm too tired to do it justice! So I'll explain it
next time. But let me just say, in order to tantalize you into tuning in
to the next episode, that he carefully analyzes quantum gravity on a
ball, imposing boundary conditions that let you see why relative states
of Chern--Simons theory give states of quantum gravity. And then he makes
the hypothesis that I mentioned at the beginning of this article. This
is that \emph{all} states of quantum gravity with these boundary
conditions come from relative states of Chern--Simons theory. He even
gives some good evidence for this hypothesis, coming originally from
Hawking's work on the thermal radiation produced by black holes! (To be
continued.)

\hypertarget{week57}{%
\section{July 3, 1995}\label{week57}}

This week I'll start by finishing up my introduction to the following
paper:

\begin{enumerate}
\def\labelenumi{\arabic{enumi})}
\tightlist
\item
  Lee Smolin, ``Linking topological quantum field theory and
  nonperturbative quantum gravity'', available as
  \href{https://arxiv.org/abs/gr-qc/9505028}{\texttt{gr-qc/9505028}}.
\end{enumerate}

\noindent
So: recall where we were. Let me not repeat the details, but simply note
that we were playing around with quantum gravity on a \(4\)-dimensional
spacetime, using the Ashtekar `new variables' formalism, and we'd
noticed that in the theory with nonzero cosmological constant
\(\Lambda\), there is an explicit solution of the theory, the
`Chern--Simons' state. Actually I hadn't really shown that this state
satisfies the key equation, the Wheeler--DeWitt equation, but if you know
the formulas it's easy to check.

Now one might think that one solution isn't all that much, apart from it
being a whole lot better than none, which was the situation before these
discoveries. However, as I began to explain last time, one can get a
whole slew of states if one takes as ones space S, not a closed
3-dimensional manifold (as we were doing at first) but a 3-manifold with
boundary. This is where Lee Smolin starts. He considers quantum gravity
with certain `self-dual boundary conditions' that are specially
compatible with Chern--Simons theory.

There is presumably an enormous space of states of quantum gravity
satisfying these boundary conditions, although we don't really know what
they look like. Say we want to understand these states as much as
possible. What do they look like? Well, first of all, the loop
representation gives us a nice picture of the `kinematical states' ---
i.e., states not necessarily satisfying the diffeomorphism constraint or
the Wheeler--DeWitt equation. (This picture may be wrong, of course, but
let me throw caution to the winds and just explain the picture.) Every
kinematical state is a linear combination of `spin network states'. For
more on spin networks, check out \protect\hyperlink{week55}{``Week 55''}
and the references in there, but let me remind you what spin networks
look like in this case.

For simplicity and ease of visualization, you can pretend \(S\) is a
ball, so its boundary is a sphere. Think of a spin network state as a
graph embedded in this ball, possibly with some edges ending on the the
boundary, with all the edges labelled by spins
\(j = 0,1/2,1,3/2,\ldots\), and with the vertices labelled by some extra
numbers that we won't worry about here. Let's call the points where
edges end on the boundary `punctures', because the idea is that they
really poke through the boundary and keep on going. Physically, these
edges graph represent `flux tubes of area': if we measure the area of
some surface in this state (or at least a surface that doesn't intersect
the vertices), the area is just the quantity \[L^2  \sqrt{j(j+1)}\]
summed over all edges that poke through the surface, where \(L\) is the
Planck length and \(j\) is the spin labelling that edge. Gauge theories
often have ``flux tube'' solutions when you quantize them: for example,
type II superconductors admit flux tubes of the magnetic field, while
superfluids admit flux tubes of angular momentum (vortices). The idea
behind spin networks in quantum gravity, physically speaking, is that
gravity is a gauge field which at the Planck scale is organized into
branching flux tubes of area.

But we want to understand, not the kinematical states in general, but
the actual physical states, which satisfy the diffeomorphism constraint
and the Wheeler--DeWitt equation. We can start by measuring everything we
measure by doing experiments right at the boundary of \(S\). More
precisely, we can try to find a maximal set of commuting observables
that ``live on the boundary'' in this sense. For example, the area of any
patch of \(S\) counts as one of these observables, and all these
`surface patch area' observables commute. If we measure all of them, we
know everything there is to know about the area of all regions on the
boundary of \(S\). Thanks to spin network technology, as described
above, specifying all their eigenvalues amounts to specifying the
location of a bunch of punctures on the boundary of \(S\), together with
the spins labelling the edges ending there.

Now Chern--Simons theory gives an obvious candidate for the space of
physical states of quantum gravity for which these ``surface patch area''
observables have specified eigenvalues. In fact, if you hand
Chern--Simons theory a surface like the boundary of \(S\), together with
a bunch of punctures labelled by spins, it gives you a \emph{finite-dimensional} 
state space. Let's not explain just now how it gives
you this state space; let's simply mumble that it gives you this space
by virtue of being an ``extended topological quantum field theory.'' If
you want, you can think of these states as being the ``relative states;' I
discussed in last week's Finds, but not all of them: only those for
which the ``surface patch area'' observables have specified eigenvalues.
There is a wonderfully simple combinatorial recipe for describing all
these states in terms of spin networks living in \(S\), having edges
that end at the punctures, with the right spins at these ends.

Smolin's hypothesis is that this finite-dimensional space of states
coming from Chern--Simons theory \emph{is} the space of all physical
states of quantum gravity on \(S\) that

\begin{enumerate}
\def\labelenumi{\arabic{enumi})}
\tightlist
\item
  satisfy the self-dual boundary conditions, and
\item
  have the specified values of the surface patch area observables.
\end{enumerate}
\noindent
Now if this hypothesis is true, it means we have a wonderfully simple
description of all the physical states on \(S\) satisfying the self-dual
boundary conditions!

So why should such a wonderful thing be true? I wish I knew! In fact,
I'm busily trying to figure it out. Smolin doesn't give any direct
evidence that it \emph{is} true, so it might not be. But he does give
some very interesting indirect evidence, coming from thermodynamics.

Thanks to work by Hawking, Bekenstein and others, there is a lot of
evidence that if one takes quantum gravity into account, the maximal
entropy of any system contained in a region with surface area \(A\)
should be proportional to \(A\). The basic idea is this. For various
reasons, one expects that the entropy of a black hole is proportional to
the area of its event horizon. For example, when you smash some black
holes together it turns out that the total area of the event horizons
goes up --- this is called the `second law of black hole
thermodynamics'. This and many more fancy thought experiments suggest
that when you have some black holes around the right notion of entropy
should include a term proportional to the total area of their event
horizons. Now suppose you had some other system which had even \emph{more}
entropy than this, but the same surface area. Then you could dump in
extra matter until it became a black hole, which would therefore have
less entropy, violating the second law.

This is a hand-waving argument, all right! It's not the sort of thing
that would convince a mathematician. But it does suggest an intriguing
connection between the vast literature on black hole thermodynamics and
the more mathematical problem of relating quantum gravity and
Chern--Simons theory.

Now the maximum entropy of a system is proportional to the logarithm of
the total number of states it can assume. So if the `Bekenstein bound'
holds, the dimension of the space of states of a system contained in a
region with surface area \(A\) is proportional to \(\exp(A/c)\) for some
constant \(c\) (which should be about the Planck length squared). Now
the remarkable thing about Smolin's hypothesis is that if it's true,
this is what one gets, because the dimension of the space given by
Chern--Simons theory does grow like this.

There is another approach leading to this conclusion that the space of
states of a bounded region should have dimensional proportional to
\(\exp(A/c)\), called the `t Hooft--Susskind holographic hypothesis. I
was going to bone up on this for This Week's Finds, but I have been too
busy! It's getting late and I'm getting bleary-eyed, so I'll stop here.
I will simply give the references to this ``holographic hypothesis''; if
anyone wants to explain it, please post to \texttt{sci.physics.research}
--- I'd be immensely grateful.

\begin{enumerate}
\def\labelenumi{\arabic{enumi})}
\setcounter{enumi}{1}
\item
  Gerhard 't Hooft, ``Dimensional reduction in quantum gravity'', 
  available as
  \href{https://arxiv.org/abs/gr-qc/9310026}{\texttt{gr-qc/9310026}}.
\item
  Leonard Susskind, ``The world as a hologram'', \emph{Jour.\
   Math.\ Phys.} \textbf{36} (1995), 6377--6396.  Also   Also available as
  \href{https://arxiv.org/abs/hep-th/9409089}{\texttt{hep-th/9409089}}.

  Leonard Susskind, ``Strings, black holes and Lorentz contractions'',
  available as
  \href{https://arxiv.org/abs/hep-th/9308139}{\texttt{hep-th/9308139}}.
\end{enumerate}

\hypertarget{week58}{%
\section{July 12, 1995}\label{week58}}

A few weeks ago I went to the IVth Porto Meeting on Knot Theory and
Physics, to which I had been kindly invited by Jose Mourao. Quite a few
of the (rather few) believers in the relevance of \(n\)-categories to
physics were there. I spoke on higher-dimensional algebra and
topological quantum field theory, and also a bit on spin networks. Louis
Crane spoke on his ideas, especially the idea of getting
\(4\)-dimensional TQFTs out of state-sum models. And John Barrett spoke
on

\begin{enumerate}
\def\labelenumi{\arabic{enumi})}
\tightlist
\item
  John Barrett, ``Quantum gravity as topological quantum field theory'',
  \emph{Jour.\ Math.\ Physics}, \textbf{36} (1995), 6161--6179.
  Also available as
  \href{https://arxiv.org/abs/gr-qc/9506070}{\texttt{gr-qc/9506070}}.
\end{enumerate}
\noindent
This is a nice introduction to the concepts of topological quantum field
theory (TQFT) that doesn't get bogged down in the (still substantial)
technicalities. In particular, it pays more emphasis than usual to the
physical interpretation of TQFTs, and how this meshes with more
traditional issues in the interpretation of quantum mechanics. One of
the main things I got out of the conference, in fact, was a sense that
there is a budding field along these lines, just crying out to be
developed. As Barrett notes, Atiyah's axioms for a TQFT can really be
seen as coming from combining

\begin{enumerate}
\def\labelenumi{\alph{enumi})}
\tightlist
\item
  The rules of quantum mechanics for composing amplitudes
\end{enumerate}
\noindent
and

\begin{enumerate}
\def\labelenumi{\alph{enumi})}
\setcounter{enumi}{1}
\tightlist
\item
  Functoriality, or the correct behavior under diffeomorphisms of
  manifolds.
\end{enumerate}
\noindent
Indeed, he convincingly recovers the TQFT axioms from these two
principles. And of course these two principles could be roughly called
``basic quantum mechanics'' and ``general covariance''... lending
credence to the idea that whatever the theory of quantum gravity turns
out to be, it should be something closely related to a TQFT. (I should
emphasize, though, that this question is one of the big puzzles in the
subject.)

The richness inherent in b) makes the business of erecting a formalism
to interpret topological quantum field theory much more interesting than
the (by now) rather stale discussions that only treat a), or ``basic
quantum mechanics''. In particular, in a TQFT, every way of combining
manifolds --- spaces or spacetimes --- yields a corresponding rule for
composing amplitudes. For example, if we have two spacetimes that look
like \[
  \begin{tikzpicture}[scale=0.5]
    \draw[thick] (0,0) ellipse (2cm and 1cm);
    \draw[thick] (-2,0) to (-2,-6);
    \draw[thick] (2,0) to (2,-6);
    \begin{scope}[shift={(0,-6)}]
      \draw[thick,dashed] (0:2) arc (0:180:2cm and 1cm);
      \draw[thick] (180:2) arc (180:360:2cm and 1cm);
    \end{scope}
  \end{tikzpicture}
\] (that's supposed to look like a pipe!) and \[
  \begin{tikzpicture}[scale=0.5]
    \draw[thick] (-3,0) ellipse (2cm and 1cm);
    \draw[thick] (3,0) ellipse (2cm and 1cm);
    \draw[thick] (-5,0) .. controls (-5,-2) and (-2,-4) .. (-2,-6);
    \draw[thick] (5,0) .. controls (5,-2) and (2,-4) .. (2,-6);
    \draw[thick] (-1,0) .. controls (-1,-1) .. (0,-2);
    \draw[thick] (1,0) .. controls (1,-1) .. (0,-2);
    \begin{scope}[shift={(0,-6)}]
      \draw[thick,dashed] (0:2) arc (0:180:2cm and 1cm);
      \draw[thick] (180:2) arc (180:360:2cm and 1cm);
    \end{scope}
  \end{tikzpicture}
\] --- that is, a cylinder and a ``trinion'' (or upside-down pair of
pants) --- we can combine them either ``horizontally'' like this: \[
  \begin{tikzpicture}[scale=0.5]
    \begin{scope}[shift={(8,0)}]
      \draw[thick] (0,0) ellipse (2cm and 1cm);
      \draw[thick] (-2,0) to (-2,-6);
      \draw[thick] (2,0) to (2,-6);
      \begin{scope}[shift={(0,-6)}]
        \draw[thick,dashed] (0:2) arc (0:180:2cm and 1cm);
        \draw[thick] (180:2) arc (180:360:2cm and 1cm);
      \end{scope}
    \end{scope}
    \draw[thick] (-3,0) ellipse (2cm and 1cm);
    \draw[thick] (3,0) ellipse (2cm and 1cm);
    \draw[thick] (-5,0) .. controls (-5,-2) and (-2,-4) .. (-2,-6);
    \draw[thick] (5,0) .. controls (5,-2) and (2,-4) .. (2,-6);
    \draw[thick] (-1,0) .. controls (-1,-1) .. (0,-2);
    \draw[thick] (1,0) .. controls (1,-1) .. (0,-2);
    \begin{scope}[shift={(0,-6)}]
      \draw[thick,dashed] (0:2) arc (0:180:2cm and 1cm);
      \draw[thick] (180:2) arc (180:360:2cm and 1cm);
    \end{scope}
  \end{tikzpicture}
\] or ``vertically'' like this: \[
  \begin{tikzpicture}[scale=0.5]
    \draw[thick] (-3,0) ellipse (2cm and 1cm);
    \draw[thick] (3,0) ellipse (2cm and 1cm);
    \draw[thick] (-5,0) .. controls (-5,-2) and (-2,-4) .. (-2,-6);
    \draw[thick] (5,0) .. controls (5,-2) and (2,-4) .. (2,-6);
    \draw[thick] (-1,0) .. controls (-1,-1) .. (0,-2);
    \draw[thick] (1,0) .. controls (1,-1) .. (0,-2);
    \begin{scope}[shift={(0,-6)}]
      \draw[thick,dashed] (0:2) arc (0:180:2cm and 1cm);
      \draw[thick] (180:2) arc (180:360:2cm and 1cm);
    \end{scope}
    \draw[thick] (-2,-6) to (-2,-10);
    \draw[thick] (2,-6) to (2,-10);
    \begin{scope}[shift={(0,-10)}]
      \draw[thick,dashed] (0:2) arc (0:180:2cm and 1cm);
      \draw[thick] (180:2) arc (180:360:2cm and 1cm);
    \end{scope}
  \end{tikzpicture}
\]

Corresponding to each spacetime we have a ``time evolution operator''
--- a linear operator that describes how states going in one end pop out
the other, ``evolved in time''. And corresponding to horizontal and
vertical composition of spacetimes we have two ways to compose
operators: horizontal composition usually being called ``tensor
product'', and vertical composition being called simply ``composition''.
These two ways satisfy some compatibility conditions, as well.

Now if one has read a bit about \(n\)-categories and/or ``extended''
topological quantum field theories, one already knows that this is just
the tip of the iceberg. If we allow ourselves to cut spacetimes into
smaller bits --- e.g., pieces with ``corners'', such as tetrahedra or
their higher-dimensional kin --- one gets more possible ways of
composing operators, and more compatibility conditions. These become
algebraically rather sophisticated, but luckily, there is a huge amount
of evidence that existing TQFTs extend to more sophisticated structures
of this sort, through a miraculous harmony between algebra and topology.

This leads to some interesting new concepts when it comes to the
physical interpretation of extended TQFTs. As Crane described in his
talk (see also his papers listed in \protect\hyperlink{week2}{``Week
2''}, \protect\hyperlink{week23}{``Week 23''} and
\protect\hyperlink{week56}{``Week 56''}), in a 4-dimensional extended TQFT
one expects the following sort of thing. If we think of an ``observer''
as a 3-manifold with boundary --- imagine a person being the 3-manifold
and his skin being the boundary, if one likes --- the extended TQFT
should assign to his boundary a ``Hilbert category'' or ``2-Hilbert
space''. This is the categorical analog of a Hilbert space. In other
words, just as a Hilbert space is a \emph{set} in which you can
\emph{sum} things and \emph{multiply} them by \emph{complex numbers},
and get \emph{complex numbers} by taking \emph{inner products} of
things, a 2-Hilbert space is an analogous structure in which every term
I just put in italics is replaced by its analog one step up the
categorical ladder. This means: 
\[
  \begin{aligned}
    \text{set} &\to \text{category}
  \\\text{sum} &\to \text{direct sum}
  \\\text{multiply} &\to \text{tensor}
  \\\text{complex numbers} &\to \text{vector spaces}
  \\\text{inner products} &\to \text{homs}
  \end{aligned}
\]

There's a good chance that you know the analogy between numbers and
vector spaces: just as you can add numbers and multiply them, you can
take direct sums and tensor products of vector spaces, and many of the
same rules still apply (in a somewhat more sophisticated form, because
laws that were equations are now isomorphisms). A little less familiar
is the analogy between inner products and ``homs''. Given two vectors
\(v\) and \(w\) in a Hilbert space you can take the inner product
\(\langle v,w\rangle\) and get a number; similarly, given two
(finite-dimensional) Hilbert spaces \(V\) and \(W\) you can form
\(\mathrm{hom}(V,W)\) --- that is, the set of all linear maps from \(V\)
to \(W\) --- and get a Hilbert space. The same thing works in any
``2-Hilbert space''.

The most basic example of a 2-Hilbert space would be Hilb, the category
of finite-dimensional Hilbert spaces, but also \(\mathsf{Reps}(G)\), the
category of finite-dimensional unitary representations of a finite
group. (Similar remarks hold for quantum groups at root of unity.) Just
as the inner product is linear in one argument and conjugate-linear in
the other, ``\(\mathrm{hom}\)'' behaves nicely under direct sums in each
argument, but each argument behaves a bit differently under tensor
product, so one can say it's ``linear'' in one and ``conjugate-linear''
in the other.

So anyway, just as in a 4d TQFT a 3-manifold \(M\) determines a Hilbert
space \(Z(M)\), and a 4-manifold \(N\) with boundary equal to \(M\)
determines a vector \(Z(N)\) in \(Z(M)\), something similar happens in
an extended TQFT. (For experts, here I'm really talking about
``unitary'' TQFTs and extended TQFTs --- these are the physically
sensible ones.) Namely, a ``skin of observation'' or 2-manifold \(S\)
determines a 2-Hilbert space \(Z(S)\), and an ``observer'' or 3-manifold
\(M\) with boundary equal to \(S\) determines an object in \(Z(S)\).
Now, given two observers \(M\) and \(M'\) with the same ``skin'' --- for
example, the observer ``you'' and the observer ``everything in the world
except you'' --- one gets two objects \(Z(M)\) and \(Z(M')\) in
\(Z(S)\), so one can form the ``inner product''
\(\mathrm{hom}(Z(M),Z(M'))\), which is a Hilbert space. This is
\emph{your} Hilbert space for describing states of \emph{everything in
the world except you}. Note that we are using the term ``observer'' here
in a somewhat whimsical sense; in particular, every region of space
counts as an observer in this game, so we can flip things around and
form the inner product \(\mathrm{hom}(Z(M'),Z(M))\), which is the
Hilbert space that \emph{everything in the world except you} can use to
describe states of \emph{you}. These two Hilbert spaces, with roles
reversed, are conjugate to each other (using an obvious but perhaps
slightly unfamiliar definition of ``conjugate'' Hilbert space), so
they're pretty much the same.

Now this may at first seem weird, but if you think about it, it becomes
a bit less so. Of course, all of this stuff simply follows from the
notion of a unitary extended TQFT, and whether the actual laws of
physics are given by such a structure is a separate issue. But there is
clearly a lot of relevance to the ``holographic hypothesis'' and Lee
Smolin's more mathematical version of that hypothesis, as sketched in
\protect\hyperlink{week57}{``Week 57''}. The basic idea, there as here,
is that we are concentrating our attention on the things about a system
that can be measured at its boundary, and what we measure there can be
either thought of describing the state of the ``inside'' or dually the
``outside''.

I think if I go out on a limb here, and rhapsodize a bit, the point
might be clearer: but don't take this too seriously. Namely: all of the
stuff you see, hear, and otherwise observe about the world --- which
seems to be ``information about the outside'' --- is also stuff going on
in your brain, hence ``information about the inside''. What this stuff
really is, of course, is \emph{correlations} between the inside and the
outside. This is the reason for the duality between observer and
observed mentioned above. Note: we need not worry here whether or not
there's ``really'' a lot more going on outside than what you observe.
The point is simply that everything \emph{you} observe about what's
going on in the world outside is correlated to stuff that the world
could observe about what is going on in you. (Maybe.)

I should perhaps also add that the mathematicians are getting a bit
behind on the job of developing the ``higher linear algebra'' needed to
support this sort of physics. So it's a bit hard to point to a good
reference for all this 2-Hilbert space stuff. I'm slowly writing a paper
on it, but for now the best sources seem to be Kapranov and Voevodsky's
work on 2-vector spaces:

\begin{enumerate}
\def\labelenumi{\arabic{enumi})}
\setcounter{enumi}{1}
\tightlist
\item
  M. Kapranov and V. Voevodsky, ``2-Categories and Zamolodchikov
  tetrahedra equations'', in \emph{Proc.\ Symp.\ Pure Math.} \textbf{56},
  Part 2 (1994), AMS, Providence, pp.~177--260.
\end{enumerate}

\noindent
(see also \protect\hyperlink{week4}{``Week 4''}), along with Dan Freed's work on
higher algebraic structures in gauge theory
(\protect\hyperlink{week12}{``Week 12''},
\protect\hyperlink{week48}{``Week 48''}), and David Yetter's new paper:

\begin{enumerate}
\def\labelenumi{\arabic{enumi})}
\setcounter{enumi}{2}
\tightlist
\item
  David Yetter, ``Categorical linear algebra: a setting for questions
  from physics and low-dimensional topology''.
  Available as \href{http://math.ucr.edu/home/baez/yetter.pdf}{\texttt{http://math.ucr.edu/home/baez/yetter.pdf}} and \href{http://math.ucr.edu/home/baez/yetter.ps}{\texttt{http://math.ucr.edu/home/baez/yetter.ps}}.
\end{enumerate}

\noindent
This treats 2-vector spaces in a very beautiful way, but not 2-Hilbert
spaces. Definitely worth reading for anyone interested in this sort of
thing!

While visiting Porto, I managed somehow to miss talking to Eugenia Cesar
de Sa, which was really a pity because she was the one who developed the
way of describing 4-manifolds that Broda (see
\protect\hyperlink{week9}{``Week 9''}, \protect\hyperlink{week10}{``Week
10''}) used to construct a 4-dimensional TQFT. This TQFT was later shown
by Roberts (see \protect\hyperlink{week14}{``Week 14''}) to be
isomorphic to that described by Crane and Yetter using a state sum model
--- i.e., by a discrete analog of a path integral in which one chops
spacetime up into \(4\)-dimensional ``hypertetrahedra'' (better known as
\(4\)-simplices!), labels their 2d and 3d faces by spins, and sums over
labellings. Her work is cited in the Broda reference in
\protect\hyperlink{week17}{``Week 17''}, but I managed luckily to get a
copy of her thesis:

\begin{enumerate}
\def\labelenumi{\arabic{enumi})}
\setcounter{enumi}{3}
\tightlist
\item
  Eugenia Cesar de Sa, \emph{Automorphisms of 3-Manifolds and
  Representations of 4-Manifolds}, Ph.D.~thesis, University of Warwick,
  1977.
\end{enumerate}

\noindent
This should let me learn more about \(4\)-dimensional topology, a
fascinating subject on which I'm woefully ignorant.

One reason why Broda's work, and thus de Sa's, is interesting to me, is
that people have suspected for a while that the Crane--Yetter--Broda
theory, which is constructed purely combinatorially, is isomorphic to \(BF\)
theory with cosmological term. \(BF\) theory (see
\protect\hyperlink{week36}{``Week 36''}) is a \(4\)-dimensional field
theory that can be described starting from a Lagrangian in the
traditional manner of physics. The theory ``with cosmological term'' can
be regarded as a baby version of quantum gravity with nonzero
cosmological constant, a baby version having only one state, the
``Chern--Simons state''. As I discussed in
\protect\hyperlink{week56}{``Week 56''}, it's this Chern--Simons state
that plays a key role in Smolin's attempt to ``exactly solve'' quantum
gravity. Thus I suspect that \(BF\) theory is a good thing to understand
really well. Recently I showed in

\begin{enumerate}
\def\labelenumi{\arabic{enumi})}
\setcounter{enumi}{4}
\tightlist
\item
  John Baez, ``4-dimensional \(BF\) theory with cosmological term as a
  topological quantum field theory'', \emph{Lett.\ Math.\ Phys.\ }\textbf{38} 
  (1996), 129--143.  Available as
  \href{https://arxiv.org/abs/q-alg/9507006}{\texttt{q-alg/9507006}}.
\end{enumerate}

\noindent
that the Crane--Yetter--Broda theory is indeed isomorphic as a TQFT to a
certain \(BF\) theory. With a bit more work, this should give us a state
sum model for the \(BF\) theory that's a baby version of quantum gravity
in 4 dimensions. This should come in handy for studying Smolin's
hypothesis and its ramifications.

\begin{enumerate}
\def\labelenumi{\arabic{enumi})}
\setcounter{enumi}{5}
\tightlist
\item
  Timothy Porter, ``TQFTs from homotopy \(n\)-types'', 
  \emph{Jour.\ London Math.\ Soc.} \textbf{58} (1998), 723--732.
\end{enumerate}

\noindent
The Dijkgraaf--Witten model is an \(n\)-dimensional TQFT one gets from a
finite group \(G\). It's given by a really simple state sum model. Chop
your manifold into simplices; then the allowed ``states'' are just
labellings of the edges with elements of \(G\) subject to the constraint
that the product around any triangle is \(1\). You can think of a
labelling as a kind of ``connection'' that tells you how to parallel
transport along the edges, and the constraint says the connection is
flat. Expectation values of physical observables are then computed as
sums over these states. In fact, this TQFT is a baby version of \(BF\)
theory \emph{without} cosmological constant. A toy model of a toy model
of quantum gravity, in other words: the classical solutions of \(BF\)
theory without cosmological constant are just flat connections on some
\(G\)-bundle where \(G\) is a Lie group, while the Dijkgraaf--Witten model does
something similar for a finite group.

In a previous paper (see \protect\hyperlink{week54}{``Week 54''}) Porter
studied the Dijkgraaf--Witten model and a generalization of it due to
Yetter that allows one to label faces with things too... one can
think of this generalization as allowing ``curvature'', because the
product of elements of \(G\) around a triangle need no longer be \(1\);
instead, it's something determined by the labelling of the face.

\begin{enumerate}
\def\labelenumi{\arabic{enumi})}
\setcounter{enumi}{6}
\tightlist
\item
  David Yetter, ``TQFTs from homotopy 2-types'', \emph{Journal of Knot
  Theory and its Ramifications} \textbf{2} (1993), 113--123.
\end{enumerate}

\noindent
In his new paper Porter takes this idea to its logical conclusion and
constructs analogous theories that allow labellings of simplices in any
dimension. Technically, the input data is no longer just a finite group,
but a finite simplicial group \(G\).

What's a simplicial group? It's a wonderful thing; using the
``internalization'' trick I've referred to in some previous Finds, all I
need to say is that it's a simplicial object in the category of groups.
A simplicial set is a bunch of sets, one for each natural number,
together with a bunch of ``face'' and ``degeneracy'' maps satisfying the
same laws that the face and degeneracy maps do for a simplex. (Students
of singular or simplicial homology will know what I'm talking about.)
Similarly, a simplicial group is a bunch of \emph{groups}, together with
a bunch of of ``face'' and ``degeneracy'' \emph{homomorphisms}
satisfying the same laws.

A triangulated manifold gives a simplicial set in an obvious way, and
from any simplicial set one can obtain a simplicial groupoid (like a
simplicial group, but with groupoids instead!) called its ``loop
groupoid''. The sort of labellings Porter considers are homomorphisms
from this simplicial groupoid to the given simplicial group G.

I will refrain from trying to say what all this has to do with homotopy
\(n\)-types. Nonetheless, from a pure mathematics point of view, that's
the most exciting aspect of the whole business! Part of the puzzle about
TQFTs is their relation to traditional algebraic topology (and
not-so-traditional algebraic topology like nonabelian cohomology,
\(n\)-stacks, etc.), and this work serves as a big clue about that
relationship.

\hypertarget{week59}{%
\section{August 3, 1995}\label{week59}}

\begin{quote}
As you crack your eyes one morning your reason is assaulted by a
strange sight. Over your head, humming quietly, there floats a
monitor, an ethereal otherworldly screen on which is written a
curious message. ``I am the Screen of ultimate Truth. I am bulging
with information and ask nothing better than to be allowed to
impart it."
\end{quote}

\noindent
It would be nice if more math books started with something
attention-grabbing like this. In fact, it appears near the beginning of

\begin{enumerate}
\def\labelenumi{\arabic{enumi})}
\tightlist
\item
  Geoffrey M.\ Dixon, \emph{Division Algebras: Octonions, Quaternions,
  Complex Numbers and the Algebraic Design of Physics}, Springer, Berlin,
  1994.
\end{enumerate}

Dixon is convinced that the details of the Standard Model of particle
interactions can be understood better by taking certain mathematical
structures very seriously. There are very few algebras over the reals
where we can divide by nonzero elements: if we demand associativity and
commutativity, just the reals themselves and the complex numbers. If we
drop the demand for commutativity, we also get a \(4\)-dimensional
algebra called the quaternions, invented by Hamilton. If in addition we
drop the demand for associativity, and ask only that our algebra be
``alternative'', we also get an \(8\)-dimensional algebra called the
octonions, or Cayley numbers. (I'll say what ``alternative'' means in
\protect\hyperlink{week61}{``Week 61''}.) Clearly these are very special
structures, and also clearly they play an important role in
physics... or do they?

Well, few people doubt that the real numbers are fundamental to physics
(though some advocates of the discrete might prefer the integers), and
with emergence of quantum theory, if not sooner, the basic role of the
complex numbers also became clear. Hamilton discovered the quaternions
in the 1800s, and used them to formulate a beautiful theory of rotations
in \(3\)-dimensional space. They fell out of favor somewhat when the
vectors of Gibbs proved simpler for many purposes, but their deeper
importance became clear when people started studying spin: indeed, the
Pauli matrices so important in physics are closely related to the
quaternions, and it is the group of unit quaternions,
\(\mathrm{SU}(2)\), rather than the group of rotations in 3d space,
\(\mathrm{SO}(3)\), which turns out to be the symmetry group whose
different representations correspond to particles of different spin. But
what about the octonions?

Well, there are not too many places in physics yet where the octonions
reach out and grab one with the force the reals, complexes, and
quaternions do. But they are certainly out there, they have a certain
beauty to them, and they are the natural stopping-point of a certain
finite sequence of structures, so it is natural for people of a certain
temperament to believe that they are there for a reason. Dixon makes a
passionate case for this in the beginning of his book.

Suppose you were confronted with the Screen of Truth. What would you ask
it? Dixon, being a physicist, naturally fantasizes asking it why the
universe is the way it is! What kind of answer could this possibly have?
Perhaps there is only one consistent way for things to be, and
mathematics, with its unique and beautiful structures that are pure
expressions of logical necessity, is trying to tell us something about
this?

On the one hand this seems outrageous... especially to the
hard-nosed pragmatist or empiricist in us. It seems old-fashioned,
naive, and too good to be true. On the other hand, the universe
\emph{is} outrageous! It's outrageous that it exists in the first place,
and doubly outrageous that it has the particular physical laws it does
and no others. It has only been through the old-fashioned, naive belief
that we can understand it using mathematics that we discovered what we
have of its physical laws. So maybe eventually we \emph{will} see that
the basic structures of mathematics determine, in some mysterious sense,
all the basic laws of physics. Or maybe we won't. In either case, there
is a long way yet to go. As Dixon's Screen of Truth comments, before it
departs:

\begin{quote}
``Do you believe that were I to explain as much of what I know as you
could comprehend that you would recognize it, that you would say, oh
yes, this is but an extension of what we have already done, and
though the mathematics is broader, the principles deeper, I am not
surprised?  Do you think you have asked even a fraction of the
questions you need to ask?''
\end{quote}

\noindent
Anyway, it is at least worth considering all the beautiful mathematical
structures one runs into for their potential importance. For example,
the octonions.

In order to write this week's Finds, I needed to learn a little about
the octonions. I wanted some good descriptions of the octonions, that
hopefully would ``explain'' them or at least make them easy to remember.
So I asked for help on \texttt{sci.physics.research}, and I got some help from
Greg Kuperberg, Ezra Getzler, Matthew Wiener, and Alexander Vlasov.
After a while Geoffrey Dixon got wind of this and referred me to his
work! I'll probably talk to him later this summer when I go back to
Cambridge Massachusetts, and hopefully I'll learn more about octonions
and the like.

But for now let me just give a quick beginner's introduction to the
octonions. A lot of this appears in

\begin{enumerate}
\def\labelenumi{\arabic{enumi})}
\setcounter{enumi}{1}
\tightlist
\item
  William Fulton and Joe Harris, \emph{Representation Theory --- a First
  Course}, Springer, Berlin, 1991.
\end{enumerate}

\noindent
I should add that this book is a very good place to learn about Lie
groups, Lie algebras, and their representations.  I wish I had
taken a course based on this book when I was in grad school!

Let's start with the real numbers. Then the complex number \[a+bi\] can
be thought of as a pair \[(a,b)\] of real numbers. Addition is done
component-wise, and multiplication goes like this:
\[(a,b)(c,d) = (ac - db,da + bc)\] We can also define the conjugate of a
complex number by \[(a,b)^* = (a,-b).\] Now that we have the complex
numbers, we can define the quaternions in a similar way. A quaternion
can be thought of as a pair \[(a,b)\] of complex numbers. Addition is
component-wise and multiplication goes like this
\[(a,b)(c,d) = (ac - d^*b, da + bc^*)\] This is just like how we defined
multiplication of complex numbers, but with a couple of conjugates
(\({}^*\)'s) thrown in. To emphasize how similar the two multiplications
are, we could have included the conjugates in the first formula, since
the conjugate of a real number is just itself.

We can also define the conjugate of a quaternion by
\[(a,b)^* = (a^*,-b).\] The game continues! Now we can define an
octonion to be a pair of quaternions; as before, we add these
component-wise and multiply them as follows:
\[(a,b)(c,d) = (ac - d^*b, da + bc^*).\] One can also define the
conjugate of an octonion by \[(a,b)^* = (a^*,-b).\] Why do the real
numbers, complex numbers, quaternions and octonions have multiplicative
inverses? I take it as obvious for the real numbers. For the complex
numbers, you can check that \[(a,b)^* (a,b) = (a,b) (a,b)^* = K (1,0)\]
where \(K\) is a real number called the ``norm squared'' of \((a,b)\).
The multiplicative identity for the complex numbers is \((1,0)\). This
means that the multiplicative inverse of \((a,b)\) is \((a,b)^*/K\). You
can check that the same holds for the quaternions and octonions!

This game of getting new algebras from old is called the
``Cayley--Dickson'' construction. Of course, the fun we've had so far
should make you want to keep playing this game and develop a
16-dimensional algebra, the ``hexadecanions,'' consisting of pairs of
octonions equipped with the same sort of multiplication law. What do you
get? Why aren't there multiplicative inverses anymore? I haven't
checked, because this is my summer vacation! I am learning about
octonions just for fun, since I just finished writing some rather
technical papers, and my idea of fun does not presently include
multiplying two hexadecanions together to see why the norm-squared law
\((a,b) (a,b)^* = (a,b)^* (a,b) = K (1,0)\) breaks down. But I'm sure
someone out there will enjoy doing this... and I'm sure someone
else out there has already done it! So they should let me know what
happens. There is something out there called ``Pfister forms'', which I
think might be related.

{[}Toby Bartels did some nice work on hexadecanions in response to the
above challenge, which appears at the end of this article.{]}

Now if we unravel the above definition of quaternions, by writing the
quaternion \((a+bi,c+di)\) as \(a+bi+cj+dk\), we see that the
multiplication law is \[i^2 = j^2 = k^2 = -1,\] and
\[ij = -ji = k, \quad jk = -kj = i, \quad ki = -ik = j.\]

For more about the inner meaning of these rules, see
\protect\hyperlink{week5}{``Week 5''}. Similarly, we can unravel the
above definition of octonions by writing the octonion
\((a+bi+cj+dk,e+fi+gj+hk)\) as
\[a + b e_1 + c e_2 + d e_3 + e e_4 + f e_5 + g e_6 + h e_7.\] Note:
since mathematicians are very impersonal, they usually call these seven
dwarves \(e_1,\ldots,e_7\) instead of Sleepy, Grumpy, etc., as in the
Disney movie. Any one of these 7 guys times himself is \(-1\). Also, any
two distinct ones anticommute; for example, \(e_3 e_7 = -e_7 e_3\).
There is a nice way to remember how to multiply them using the ``Fano
plane''. This is a projective plane with 7 points, where by a
``projective plane'' I mean that any two points determine an abstract
sort of ``line'', which in this case consists of just 3 points, and any
two lines intersect in a point. It looks like this:
\[\includegraphics[max width=0.65\linewidth]{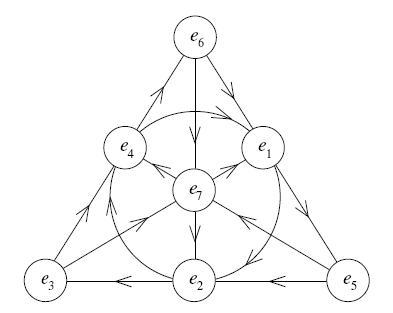}\]

The ``lines'' are the 3 edges of the big triangle, the 3 lines going
through a vertex, the center and the midpoint of the opposite edge, and
the circle including \(e_1\), \(e_2\), and \(e_4\). All the ``lines''
are cyclically ordered, and that tells you how to multiply the seven
dwarves. For example, the line that's actually a circle goes clockwise,
so \(e_1 e_2 = e_4\), \(e_2 e_4 = e_1\), and \(e_4 e_1 = e_2\). The
lines that are edges of the big triangle also point clockwise, so for
example \(e_5 e_2 = e_3\), and cyclic permutations thereof, and
\(e_6 e_3 = e_4\). The lines that go through the center point from the
vertex to the midpoint of the opposite edge, so for example
\(e_3 e_7 = e_1\). I hope that made sense; you can work it out yourself,
of course.

My convention for numbering the seven dwarves in the picture above is
\emph{completely arbitrary}, so don't bother remembering it --- make up
your own if you prefer! The convention I chose looks sort of weird at
first, but it has a couple of endearing features:

\begin{itemize}
\tightlist
\item
  Index cycling: if \(e_i e_j = e_k\), then
  \(e_{i+1} e_{j+1} = e_{k+1}\).
\item
  Index doubling: if \(e_i e_j = e_k\), then \(e_{2i} e_{2j} = e_{2k}\).
\end{itemize}

\noindent
Here we add and multiply mod \(7\). Index doubling corresponds to
rotating the Fano plane.

So those are the octonions in a nutshell. I should say a bit about how
they relate to triality for \(\mathrm{SO}(8)\), the exceptional Lie
group \(\mathrm{G}_2\), the group \(\mathrm{SU}(3)\) which is so
important in the study of the strong force, and to lattices like
\(\mathrm{E}_8\), \(\Lambda 16\) and the Leech lattice. But I will
postpone that; for now you can consult Fulton and Harris, and also
various papers by Dixon:

\begin{enumerate}
\def\labelenumi{\arabic{enumi})}
\setcounter{enumi}{2}
\item
  Geoffrey Dixon, ``Octonion X-product orbits'', available as
  \href{https://arxiv.org/abs/hep-th/9410202}{\texttt{hep-th/9410202}}.

  ``Octonion X-product and \(\mathrm{E}_8\) lattices'', 
  available as
  \href{https://arxiv.org/abs/hep-th/9411063}{\texttt{hep-th/9411063}}.

  ``Octonions: \(\mathrm{E}_8\) lattice to \(\Lambda 16\)'', 
  available as
  \href{https://arxiv.org/abs/hep-th/9501007}{\texttt{hep-th/9501007}}.

  ``Octonions: invariant representation of the Leech lattice'', 
  available as
  \href{https://arxiv.org/abs/hep-th/9504040}{\texttt{hep-th/9504040}}.

  ``Octonions: invariant Leech lattice exposed'', available as
  \href{https://arxiv.org/abs/hep-th/9506080}{\texttt{hep-th/9506080}}.
\end{enumerate}

\noindent
I am not presently in a position to assess these papers or Dixon's work
relating division algebras and the Standard Model, but hopefully
sometime I will be able to say a bit more.

Let me wrap up by saying a bit about the Leech lattice. As described in
my review of Conway and Sloane's book (\protect\hyperlink{week20}{``Week
20''}), there is a wonderful branch of mathematics that studies the
densest ways of packing spheres in n dimensions. Most of the results so
far concern lattice packings, packings in which the centers of the
spheres form a subset of \(\mathbb{R}^n\) closed under addition and
scalar multiplication by integers. When \(n = 8\), the densest known
packing is given by the so-called \(\mathrm{E}_8\) lattice. In
\protect\hyperlink{week20}{``Week 20''} I described how to get this
lattice using the quaternions and the icosahedron. Briefly, it goes as
follows. The group of rotational symmetries of the icosahedron (not
counting reflections) is a subgroup of the rotation group
\(\mathrm{SO}(3)\) containing 60 elements. As mentioned above,
\(\mathrm{SO}(3)\) has as a double cover the group \(\mathrm{SU}(2)\) of
unit quaternions. So there is a 120-element subgroup of
\(\mathrm{SU}(2)\) consisting of elements that map to elements of
\(\mathrm{SO}(3)\) that are symmetries of the icosahedron. Now form all
integer linear combinations of these 120 special elements of
\(\mathrm{SU}(2)\). We get a subring of the quaternions known as the
``icosians''.

We can think of icosians as special quaternions, but we can also think
of them as special vectors in \(\mathbb{R}^8\), as follows. Every
icosian is of the form
\[(a + \sqrt{5} b) + (c + \sqrt{5} d)i + (e + \sqrt{5} f)j + (g + \sqrt{5} h)k\]
with \(a,b,c,d,e,f,g,h\) rational --- but not all rational values of
\(a,\ldots,h\) give icosians. The set of all vectors
\(x = (a,b,c,d,e,f,g,h)\) in \(\mathbb{R}^8\) that correspond to
icosians in this way is the \(\mathrm{E}_8\) lattice!

The Leech lattice is the densest known packing in 24 dimensions. It has
all sorts of remarkable properties. Here is an easy way to get ones
hands on it. First consider triples of icosians \((x,y,z)\). Let \(L\)
be the set of such triples with \[x = y = z \mod h\] and
\[x + y + z = 0 \mod h^*\] where \(h\) is the quaternion
\((-\sqrt{5} + i + j + k)/2\). Since we can think of an icosian as a
vector in \(\mathbb{R}^8\), we can think of \(L\) as a subset of
\(\mathbb{R}^{24}\). It is a lattice, and in fact, it's the Leech
lattice! I have a bit more to say about the Leech lattice in
\protect\hyperlink{week20}{``Week 20''}, but the real place to go for
information on this beast is Conway and Sloane's book. It turns out to
be related to all sorts of other "exceptional'\,' algebraic structures.
People have found uses for many of these in string theory, so if string
theory is right, maybe they are important in physics. Personally, I want
to understand them more deeply as pure mathematics before worrying too
much about their applications to physics.

\begin{center}\rule{0.5\linewidth}{0.5pt}\end{center}

Here is what Toby Bartels wrote about the hexadecanions.  (By the way,
these are usually called the ``sedenions''.)

\begin{quote}
From: Toby Bartels  \hfill \break
Subject: Re: why hexadecanions have no inverses \hfill \break
To: John Baez \hfill \break
Date: Sun, 20 Aug 1995 \hfill
\end{quote}

\begin{quote}
I spent a couple days thinking about why hexadecanions have no inverses,
and the first thing I want to say about it is that they do. However,
these inverses are of limited applicability, because the hexadecanions
are not a division algebra. A division algebra allows you to conclude,
given \(x y = 0\), that \(x\) or \(y\) is \(0\). If your algebra has
inverses, you might try to multiply this equation by the inverse of
\(x\) or \(y\) (whichever one isn't \(0\)) to prove the other is \(0\),
but this only works if the algebra is associative. Since the octonions
and hexadecanions aren't associative, there's no reason (yet) to think
either of these is a division algebra. It turns out that the octonions
are a division algebra, despite not being associative, but the
hexadecanions aren't.

Why aren't the hexadecanions a division algebra? Because the real
numbers aren't of characteristic 2. Allow me to explain.

I will prove below that the \(2^n\) onions are a division algebra only
if the \(2^{n-1}\) onions are associative. So, the question becomes: why
aren't the octonions associative? Well, I've found a proof that \(2^n\)
onions are associative only if \(2^{n-1}\) onions are commutative. So,
why aren't the quaternions commutative? Again, I have a proof that
\(2^n\) onions are commutative only if \(2^{n-1}\) onions equal their
own conjugates. So, why don't the complex numbers equal their own
conjugates? I have a proof that \(2^n\) onions do equal their own
conjugates, but it works only if the \(2^{n-1}\) onions are of
characteristic 2. The real numbers are not of characteristic 2, so the
complex numbers don't equal their own conjugates, so the quaternions
aren't commutative, so the octonions aren't associative, so the
hexadecanions aren't a division algebra.

I require a few identities about conjugates that hold for all \(2^n\)
onions: \((x^*)^* = x\), \((x + y)^* = x^* + y^*\), and
\((x y)^* = y^* x^*\). (If these identities are reminiscent of
identities for transposes of matrices, it is no coincidence.) I will
prove these by induction. That is, if an identity holds for \(2^{n-1}\)
onions, I show it holds for \(2^n\) onions. Since they hold trivially
for the reals (\(n = 0\)), they hold for all.

\[((a, b)^*)^* = (a^*, -b)^* = ((a^*)^*, -(-b)).\] By the induction
hypothesis and the nature of addition (an Abelian group),
\[((a^*)^*, -(-b)) = (a, b).\]
\[((a, b) + (c, d))^* = (a + c, b + d)^* = ((a + c)^*, -(b + d)).\] By
the induction hypothesis and addition again,
\[((a + c)^*, -(b + d)) = (a^* + c^*, -b + -d) = (a^*, -b) + (c^*, -d) = (a, b)^* + (c, d)^*.\]

The next proof needs the distribution of multiplication over addition.
\[(a, b) ((c, d) + (e, f)) = (a, b) (c + e, d + f) = (a (c + e) - (d + f)^* b, (d + f) a + b (c + e)^*).\]
By the induction hypothesis and the identity immediately above, \[
  \begin{gathered}
    (a (c + e) - (d + f)^* b, (d + f) a + b (c + e)^*)
  \\= (a c + a e - d^* b - f^* b, d a + f a + b c^* + b e^*)
  \\= (a c - d^* b, d a + b c^*) + (a e - f^* b, f a + b e^*)
  \\= (a, b) (c, d) + (a, b) (e, f).
  \end{gathered}
\] Also, \[
  \begin{gathered}
    ((a, b) + (c, d)) (e, f)
  \\= (a + c, b + d) (e, f)
  \\= ((a + c) e - f^* (b + d), f (a + c) + (b + d) e^*).
  \end{gathered}
\] By the induction hypothesis again, \[
  \begin{gathered}
    ((a + c) e - f^* (b + d), f (a + c) + (b + d) e^*)
  \\= (a e + c e - f^* b - f^* d, f a + f c + b e^* + d e^*)
  \\= (a e - f^* b, f a + b e^*) + (c e - f^* d, f c + d e^*)
  \\= (a, b) (e, f) + (c, d) (e, f).
  \end{gathered}
\]

\[((a, b) (c, d))^* = (a c - d^* b, d a + b c^*)^* = ((a c - d^* b)^*, -(d a + b c^*)).\]
Using the induction hypothesis and each of the above identities, \[
  \begin{gathered}
    ((a c - d^* b)^*, -(d a + b c^*))
  \\= (c^* a^* - (-b)^* (-d), -d a + (-b) c^*)
  \\= (c^*, -d) (a^*, -b)
  \\= (c, d)^* (a, b)^*.
  \end{gathered}
\]

In light of the above identities, if I ever come across, say,
\((x y^* + z)^*\), I'll simply write \(y x^* + z^*\) without a moment's
hesitation.

Since inductive proofs have been so useful, I'll use one to prove that
\(2^n\) onions always have inverses. First, I'll extend the method in
John's article, beginning with an inductive proof that \(x x^* = x^* x\)
is real. \[(a, b) (a, b)^* = (a, b) (a^*, -b) = (a a^* + b^* b, 0),\]
and \[(a, b)^* (a, b) = (a^*, -b) (a, b) = (a^* a + b^* b, 0).\] The
inductive hypothesis states that both \(a^* a = a a^*\) and \(b^* b\)
are real, so \((a, b) (a, b)^* = (a, b)^* (a, b)\) is real. Since the
sum of a positive real and a nonnegative real is positive, I can take
this as a proof by induction that \(x x^* = x^* x\) is not only real,
but is also positive unless \(x = 0\) (which will be important). All you
have to do now is check that these things are true of the \(2^0\)
onions, and they are, quite trivially (since the \(2^0\) onions are the
reals).

Since the \(2^n\) onions are always a vector space over the reals (as
mentioned in John's article),
\[x (x^* / (x x^*)) = (x x^*) / (x x^*) = 1.\] Since one can always
divide by the real \(x x^*\), the inverse of \(x\) is \(x^* / (x x^*)\)
in any \(2^n\) onion algebra.

To continue with the streak of inductive proofs, I will now try to prove
that the \(2^n\) onions are always a division algebra. (I will fail.)
Assume \[0 = (0, 0) = (a, b) (c, d) = (a c - d^* b, d a + b c^*).\] This
gives the system of equations \[a c - d^* b = 0 = d a + b c^*.\]
Multiplying,
\[(a c) c^* - (d^* b) c^* = 0 c^* = 0 = d^* 0 = d^* (d a) + d^* (b c^*).\]
If \(2^{n-1}\) onions are associative, I can add the equations to get
\[a (c c^*) + (d^* d) a = 0.\] Since \(c c^*\) and \(d^* d\) are real,
they commute with \(a\), and the division algebra nature of \(2^{n-1}\)
onions allows me to conclude that \(c c^* + d^* d = 0\) (which implies
\(c = d = 0\) in light of positive definiteness) or that \(a = 0\) (from
which the original equation gives \(b = 0\)). Thus, the octonions are a
division algebra (since the quaternions are associative), but the
hexadecanions aren't (since the octonions aren't associative).

(If you're reading carefully, you realize that I haven't really proved
that the hexadecanions aren't a division algebra. I've failed to prove
that they are, but that's not the same thing. When I first wrote this, I
wasn't reading carefully; I will return to plug this hole later.)

Thus, the \(2^n\) onions are a division algebra iff the \(2^{n-1}\)
onions are a division algebra and are associative. So, let's try to
prove associativity of \(2^n\) onions by induction. \[
  \begin{gathered}
    ((a, b) (c, d)) (e, f)
  \\= (a c - d^* b, d a + b c^*) (e, f)
  \\= ((a c - d^* b) e - f^* (d a + b c^*), f (a c - d^* b) + (d a + b c^*) e^*)
  \\=((ac)e - (d^* b)e - f^* (da) - f^* (b c^*), f(ac) - f(d^* b) + (da) e^* + (b c^*) e^*).
  \end{gathered}
\] On the other hand, \[
  \begin{gathered}
    (a, b) ((c, d) (e, f))
  \\= (a, b) (c e - f^* d, f c + d e^*)
  \\= (a (c e - f^* d) - (f c + d e^*)^* b, (f c + d e^*) a + b (c e - f^* d)^*)
  \\= (a(ce) - a(f^* d) - (c^* f^*)b - (e d^*)b, (fc)a + (d e^*)a + b(e^* c^*) - b(d^* f)).
  \end{gathered}
\] Assuming the induction hypothesis that \(2^{n-1}\) onions are
associative, these are equal in general iff \(2^{n-1}\) onions also are
commutative.

Thus, \(2^n\) onions are associative iff \(2^{n-1}\) onions are
associative and are commutative. So, let's try to prove commutativity of
\(2^n\) onions by induction.
\[(a, b) (c, d) = (a c - d^* b, d a + b c^*).\] On the other hand,
\[(c, d) (a, b) = (c a - b^* d, b c + d a^*).\] Assuming the induction
hypothesis that \(2^{n-1}\) onions are commutative, these are equal in
general iff \(2^{n-1}\) onions also equal their own conjugates.

Thus, \(2^n\) onions are commutative iff \(2^{n-1}\) onions are
commutative and equal their own conjugates. So, let's try to prove
conjugate equality of \(2^n\) onions by induction. \[(a, b) = (a, b).\]
On the other hand, \[(a, b)^* = (a^*, -b).\] Assuming the induction
hypothesis that \(2^{n-1}\) onions equal their own conjugates, these are
equal in general iff \(2^{n-1}\) onions also have characteristic 2.
(\(b = -b\) means \(0 = b + b = 1 b + 1 b = (1 + 1) b = 2 b\); this is
true in general iff \(0 = 2\), which is what characteristic 2 means.)

Thus, \(2^n\) onions equal their own conjugates iff \(2^{n-1}\) onions
equal their own conjugates and have characteristic 2. Since the reals
don't have characteristic 2, there's no point in trying to prove
anything about that by induction. However, it's a general result that
any algebra has characteristic 2 if it has a superalgebra of
characteristic 2. Since the \(2^n\) onions are all superalgebras of the
reals (which means the reals are always isomorphic to a subset of the
\(2^n\) onions), none of the \(2^n\) onions can have characteristic 2 if
the reals don't.

In summary, the definition of the reals as the complete ordered field,
along with an initial definition that \(x^* = x\) in the reals, allows
trivial proofs that: they form a division algebra, they are associative,
they are commutative, and they equal their own conjugates, but they
don't have characteristic 2. (All of these, in fact, are true of any
ordered field with this definition of conjugate, complete or not.) From
this and the above considerations, the complex numbers form a division
algebra, are associative, and are commutative, but they neither equal
their own conjugates nor have characteristic 2. From this, the
quaternions form a division algebra and are associative, but they
neither are commutative, equal their own conjugates, nor have
characteristic 2. From this, the octonions form a division algebra but
they neither associative, are commutative, equal their own conjugates,
nor have characteristic 2. Finally, the hexadecanions neither form a
division algebra, are associative, are commutative, equal their own
conjugates, nor have characteristic 2.

At this point, I must return to the logical hole I mentioned earlier.
But I want to work with a different algebraic concept than a division
algebra; instead I will use (inspired by Doug Merrit's post to
{\rm \texttt{sci.physics.research}}) what I guess is called `alternativity',
which says \(x (x y) = (x x) y\). I don't like putting alternativity
into the pattern, since associativity implies alternativity. All the
other properties (commutativity, conjugate equality, characteristic) are
logically independent in general. I'd like to prove that every
associative \(2^n\) onion algebra is alternative, just as I proved every
commutative one was associative, without its having been obvious to
begin with. Well, I will be disappointed even more badly later on.

Taking the conjugate of \(x (x y) = (x x) y\),
\[(y^* x^*) x^* = y^* (x^* x^*),\] so left alternativity implies right
alternativity, for \(2^n\) onions.

I require an additional general identity of \(2^n\) onions. Earlier, I
proved by induction that \(x x^*\) was real, but now I need the reality
of \(x + x^*\). Like everything else, this is proved by induction.
\[(a, b) + (a, b)^* = (a, b) + (a^*, -b) = (a + a^*, 0).\] Thus, if
\(a + a^*\) is real, \((a, b) + (a, b)^*\) is real. Since \(x + x^*\) is
real when \(x\) is real, \(x + x^*\) is real when \(x\) is any \(2^n\)
onion.

Now suppose we're working in an alternative \(2^n\) onion algebra.
\[x (x y) + x^* (x y) = (x + x^*) (x y).\] Since \(x + x^*\) is real, it
associates, so
\[x (x y) + x^* (x y) = ((x + x^*) x) y = (x x) y + (x^* x) y.\] Since
\(x (x y) = (x x) y\), \[x^* (x y) = (x^* x) y,\] which will be needed.

Let's attempt to prove by induction that \(2^n\) onions are always
alternative. \[
  \begin{gathered}
    (a, b) ((a, b) (c, d))
  \\= (a, b) (a c - d^* b, d a + b c^*)
  \\= (a (a c - d^* b) - (d a + b c^*)^* b, (d a + b c^*) a + b (a c - d^* b)^*)
  \\= (a(ac) - a(d^* b) - (a^* d^*)b - (c b^*)b, (da)a + (b c^*)a + b(c^* a^*) - b(b^* d)).
  \end{gathered}
\] Meanwhile, \[
  \begin{gathered}
    ((a, b) (a, b)) (c, d)
  \\= (a a - b^* b, b a + b a^*) (c, d)
  \\= ((aa)c - (b^* b)c - d^* (ba) - d^* (b a^*),d(aa) - d(b^* b) + (ba) c^* + (b a^*) c^*).
  \end{gathered}
\] These are indeed equal in general iff \(2^{n-1}\) onions are
associative.

The last sentence may not be immediately obvious. The induction
hypothesis and its corollaries leave us with
\(x (y z) + (x^* y) z = y (z x) + y (z x^*)\) as a necessary and
sufficient condition. It may not be clear that associativity implies
this, much less vice versa. However, the reality of \(x + x^*\) once
more enters the picture.
\[y (z x) + y (z x^*) = y (z (x + x^*)) = (x + x^*) (y z) = x (y z) + x^* (y z).\]
Thus, the condition becomes
\[x (y z) + (x^* y) z = x (y z) + x^* (y z),\] which is equivalent, in
the general case, to associativity.

To sum up the findings so far: For any n, the \(2^n\) onions form a
vector space over the reals. \(x + x^*\) and \(x x^*\) are real if \(x\)
is any \(2^n\) onion; additionally, \(x x^* = x^* x.\) Every \(2^n\)
onion has an inverse, which is a real multiple of its conjugate.
Conjugation is analogous to matrix transposition in that
\[(x^*)^* = x, (x + y)^* = x^* + y^*, \textit{ and } (x y)^* = y^* x^*.\]
Multiplication distributes over addition every time. For no n do all
\(2^n\) onions equal their own negatives. \(2^{n+1}\) onions equal their
own conjugates iff \(2^n\) onions equal their own conjugates and their
own negatives. all \(2^{n+1}\) onions commute iff all \(2^n\) onions
commute and equal their own conjugates. \(2^{n+1}\) onions are
associative iff \(2^n\) onions are associative and commutative.
\(2^{n+1}\) onions are alternative iff \(2^n\) onions are alternative
and associative. The \(2^n\) onions form a division algebra if they are
alternative.

I will be satisfied if I can prove the converse of the last statement.
In light of the results about alternativity, my original attempt to
prove that division of \(2^n\) onions requires associativity of
\(2^{n-1}\) onions looks even more convincing, (since alternativity of
\(2^{n-1}\) onions can be included in the induction hypothesis), but
it's still not valid. I still haven't shown that, if \(2^{n-1}\) onions
aren't alternative, there must be nonzero \(2^n\) onions \(x\) and \(y\)
such that \(x y = 0\). There doesn't seem to be any reason why there
shouldn't be, but there just might happen not to be any. So, despite the
inelegance of it all, in order to prove that the hexadecanions aren't a
division algebra, I'm forced to exhibit nonzero \(x\) and \(y\) such
that \(x y = 0\).

Just playing around, I found \[
  \begin{gathered}
    (e_1, e_4) (-1, e_5)
  \\= (e_1 (-1) - (e_5)^* e_4, e_5 e_1 + e_4 (-1)^*)
  \\= (-e_1 + e_5 e_4, e_5 e_1 - e_4).
  \end{gathered}
\] Since \(e_5 e_4 = (0, i) (0, 1) = (i, 0) = e_1\) and
\(e_5 e_1 = (0, i) (i, 0) = (0, i^* i) = (0, 1) = e_4\),
\[(e_1, e_4) (-1, e_5) = (0, 0) = 0.\]

The \(2^n\) onions can't be a division algebra if the \(2^{n-1}\) onions
aren't. If \(x y = 0\) in the \(2^{n-1}\) onions,
\((x, 0) (y, 0) = (x y, 0) = (0, 0) = 0\). Thus, the octonions and below
are the only \(2^n\) onions to be division algebras. Still, I wish I had
a proof of this that didn't require the ugly brute force use of a
specific counterexample. (This is the interested reader's cue....)

-- Toby
\end{quote}

By the way, in a post to \texttt{sci.physics.research} on November 2,
1999, Ralph Hartley pointed out that even if we start with a field of
characteristic 2, repeatedly applying the Cayley--Dickson construction
will \emph{not} lead to an infinite sequence of division algebras,
because it's not true in this case that if \(x\) is nonzero, \(xx^*\) is
nonzero. The problem is that a field of characteristic 2 can't be an
ordered field.

\hypertarget{week60}{%
\section{August 8, 1995}\label{week60}}

The end of a sabbatical is a somewhat sad affair: so many plans
one had, and so few accomplished! As I pack my bags to return from
Cambridge England to Cambridge Massachusetts, and then wing my way back
to Riverside, I find I have quite a stack of preprints that I meant to
include in This Week's Finds, but haven't gotten around to yet.

\begin{enumerate}
\def\labelenumi{\arabic{enumi})}
\tightlist
\item
  N. P. Landsman, ``Rieffel induction as generalized quantum
  Marsden--Weinstein reduction'', \emph{Journal of Geometry and Physics}
  \textbf{15} (1995), 285--319.  Also available as \href{https://arxiv.org/abs/hep-th/9305088}{\texttt{hep-th/9305088}}.
\end{enumerate}

\noindent
Marsden--Weinstein reduction, also called symplectic reduction, is the
modern way to deal with constraints in classical mechanics problems.
It's a two-step procedure where first one looks at the subspace of your
phase space on which the constraints vanish, and then a quotient of this
by a certain equivalence relation. For example, if you have a particle
in a plane, its phase space is \(\mathbb{R}^4\), with coordinates
\((x,y,p_x,p_y)\) representing the \(x\) and \(y\) components of the
position and the \(x\) and \(y\) components of the momentum. If we have
a constraint \(x = 0\), Marsden--Weinstein reduction tells us first to
form the subspace of our phase space on which \(x = 0\), and then
quotient by the equivalence relation where two points are equivalent if
they have the same value of \(p_x\). We get down to \(\mathbb{R}^2\),
with coordinates \((y,p_y)\). But Marsden- Weinstein reduction works in
great generality and has become a basic part of the toolkit of
mathematical physics.

What's the quantum analog of Marsden--Weinstein reduction? That's what
this paper is about. Quantum mechanics in the presence of constraints is
a tricky and important business, and there are lots of theories about
how to do it. Gauge theories all have constraints, and so does general
relativity... and in quantizing general relativity, the presence of
constraints is what gives rise to the ``problem of time''. (See
\protect\hyperlink{week43}{``Week 43''} for more on this.) What
attracted my attention to this paper is a two-stage procedure for
dealing with contraints, quite analogous to Marsden--Weinstein reduction.
This should shed some interesting light on the problem of time.

\begin{enumerate}
\def\labelenumi{\arabic{enumi})}
\setcounter{enumi}{1}
\item
  Tomotada Ohtsuki, ``Finite type invariants of integral homology 3-spheres'',
  \emph{Journal of Knot Theory and its Ramifications} \textbf{5} (1996), 101--116.

  Lev Rozansky, ``The trivial connection contribution to Witten's
  invariant and finite type invariants of rational homology spheres'',
   available as
  \href{https://arxiv.org/abs/q-alg/9504015}{\texttt{q-alg/9505015}}.

  Stavros Garoufalidis, ``On finite type 3-manifold invariants I'', \emph{Journal 
  of Knot Theory and its Ramifications} \textbf{5}, 441--461.

  Stavros Garoufalidis and Jerome Levine, ``On finite type 3-manifold
  invariants II'', available as \href{https://arxiv.org/abs/q-alg/9506012}{\texttt{q-alg/9506012}}.

  Ruth J.\ Lawrence, ``Asymptotic expansions of Witten--Reshetikhin--Turaev
  invariants for some simple 3-manifolds'', \emph{Jour.\ Math.\ Phys.} \textbf{36} 
  (1995), 6106--6129.
\end{enumerate}
\noindent
Chern--Simons theory gives invariant of links in \(\mathbb{R}^3\), which
are functions of Planck's constant \(\hbar\), and if one expands them as
power series in \(\hbar\), the coefficients are link invariants with special
properties, which one summarizes by calling them ``Vassiliev
invariants'' or ``invariants of finite type''. (See
\protect\hyperlink{week3}{``Week 3''} for more.) But the partition
function of Chern--Simons theory on a compact oriented 3-manifold is also
interesting; it's an invariant of the 3-manifold defined for certain
values of \(\hbar\). (Often instead one thinks of it as a
function of a quantity \(q\), the limit \(q \to 1\) corresponding to the
limit \(\hbar \to 0\).)

Recently people have studied the partition function of special
3-manifolds called homology spheres, which have the same homology as
\(S^3\). (People have looked at both integral and rational homology
spheres.) After a bit of subtle fiddling, one can extract from the
partition function of a homology sphere a power series in
\[\hbar' = q - 1,\] and the coefficients of the powers of \(\hbar'\)
have been conjectured by Rozansky to have nice properties which one may
summarize by calling them ``finite type'' invariants, in analogy to the
link invariant case. (Namely, that they transform in nice ways under
Dehn surgery.) For example, the coefficient of \(\hbar'\) itself is 6
times the Casson invariant of the (integral) homology 3-sphere. So there
appears to be a budding branch of ``perturbative 3-manifold invariant
theory''. I just wish I understood better what's really going on behind
all this!

\begin{enumerate}
\def\labelenumi{\arabic{enumi})}
\setcounter{enumi}{2}
\tightlist
\item
  Thomas Friedrich, ``Neue Invarianten der \(4\)-dimensionalen
  Mannigfaltigkeiten'', Berlin preprint.
\end{enumerate}

\noindent
This is a nice introduction to the new Seiberg--Witten approach to
Donaldson theory, which does not assume you already know the old stuff
by heart. Very pretty mathematics!

\begin{enumerate}
\def\labelenumi{\arabic{enumi})}
\setcounter{enumi}{3}
\tightlist
\item
  Andr\'e Joyal, Ross Street, and Dominic Verity, ``Traced monoidal
  categories'', \emph{Math.\ Proc.\ Camb.\ Phil.\ Soc.} \textbf{119} (1996), 447--468.
\end{enumerate}
\noindent
This is an abstract characterization of monoidal categories (categories
with tensor products) which have a good notion of the ``trace'' of a
morphism. Many abstract treatments of traces assume that your category
is ``rigid symmetric'' or ``balanced'', meaning that your objects have
duals and you can switch around objects in order to define the trace of
a morphism \(f\colon V \to V\) in a manner analogous to how one usually
does it in linear algebra, as a certain composite:
\[I\to V\otimes V^* \xrightarrow{f\otimes1}V\otimes V^*\to I\] where
\(I\) is the ``unit object'' for the tensor product (e.g.~the complex
numbers, when we're working in the category of vector spaces.) But one
does not really need all this extra structure if all one wants is a good
notion of ``trace''. This paper isolates the bare minimum. As one might
expect if one knows the relation between knot theory and category
theory, there are lots of nice pictures of tangles in this paper!

\begin{enumerate}
\def\labelenumi{\arabic{enumi})}
\setcounter{enumi}{4}
\tightlist
\item
  Michael Reisenberger, ``Worldsheet formulations of gauge theories and
  gravity'', available as
  \href{https://arxiv.org/abs/gr-qc/9412035}{\texttt{gr-qc/9412035}}.
\end{enumerate}
\noindent
The loop representation of a gauge theory describes states as linear
combinations of loops in space, or more generally, ``spin networks''.
What's the spacetime picture of which this is a spacelike slice? The
obvious thing that comes to mind is a two-dimensional surface of some
sort. I've advocated this point of view myself in an attempt to relate
the loop representation of gravity to string theory (see
\protect\hyperlink{week18}{``Week 18''}). Here Reisenberger makes some
progress in making this precise for some simpler theories analogous to
gravity --- for example, \(BF\) theory.

And now for some things I \emph{did} manage to finish up on my
sabbatical:

\begin{enumerate}
\def\labelenumi{\arabic{enumi})}
\setcounter{enumi}{5}
\tightlist
\item
  John Baez and Stephen Sawin, ``Functional integration on spaces of
  connections'', \emph{Jour. Functional Analysis} \textbf{150} (1997), 1--27. 
  Also available as
  \href{https://arxiv.org/abs/q-alg/9507023}{\texttt{q-alg/9507023}}.
\end{enumerate}

\noindent
As I described in \protect\hyperlink{week55}{``Week 55''}, it's now
possible to set up a rigorous version of the loop representation without
assuming (as had earlier been required) that ones manifold is
real-analytic and the loops are all analytic. This means that one can do
things in a manner invariant under all diffeomorphisms, not just
analytic ones. To achieve this, one needs to ponder rather carefully the
complicated ways smooth paths, even embedded ones, can intersect (for
example, they can intersect in a Cantor set).

\begin{enumerate}
\def\labelenumi{\arabic{enumi})}
\setcounter{enumi}{6}
\tightlist
\item
  John Baez, Javier P.\ Muniain and Dardo Piriz, ``Quantum gravity
  hamiltonian for manifolds with boundary'', \emph{Phys. Rev. D} \textbf{52} 
  (1995), 6840--6845. 
  Also available as
  \href{https://arxiv.org/abs/gr-qc/9501016}{\texttt{gr-qc/9501016}}.
\end{enumerate}

\noindent
When space is a compact manifold with no boundary, there is no Hamiltonian
in quantum gravity, just a Hamiltonian constraint (see
\protect\hyperlink{week43}{``Week 43''}). This makes it tricky to
understand time evolution in the theory --- the ``problem of time''. But
with a boundary, there is a Hamiltonian, given by a surface integral
over the boundary. (The reason is that, at least when the equations of
motion hold, the Hamiltonian is a total divergence, so you can use
Gauss' theorem to express it as an integral over the boundary, which of
course is zero when there is no boundary.)

Rovelli and Smolin (see \protect\hyperlink{week42}{``Week 42''}) worked
out a loop representation of quantum gravity --- in a heuristic sort of
way which various slower sorts like myself have been struggling to make
rigorous in the subsequent years --- and a key step in this was
expressing the Hamiltonian constraint in terms of loops. In this paper
we do the same sort of thing for the Hamiltonian, when there is a
boundary. This requires considering not only loops but also paths that
start and end at the boundary.

Remarkably, the Hamiltonian acts on paths that start and end at the
boundary in a manner very similar to the Hamiltonian constraint for
quantum gravity coupled to massless chiral spinors (e.g.~neutrinos, if
neutrinos are really massless and have a ``handedness'' as they appear
to). This suggests that on a manifold with boundary, the degrees of
freedom ``living on the boundary'' are described by a chiral spinor
field. Steve Carlip has already shown something very similar for quantum
gravity in 2+1 dimensional spacetime, a more tractable simplified model
--- see \protect\hyperlink{week41}{``Week 41''}. Moreover, he used this
to explain why the entropy of a black hole is proportional to its area
(or length in 2+1 dimensions). The idea is that the entropy is really
accounted for by the degrees of freedom of the event horizon itself. It
would be nice to do something similar in 3+1-dimensional spacetime.

\hypertarget{week61}{%
\section{August 24, 1995}\label{week61}}

I'd like to return to the theme of octonions, which I began to explore
in \protect\hyperlink{week59}{``Week 59''}. The recipe I described
there, which starts with the real numbers, and then builds up the
complex numbers, quaternions, octonions, hexadecanions etc. by a
recursive process, is called the ``Cayley--Dickson process''. Now let me
describe a way to obtain the octonions using a special property of
rotations in 8-dimensional space, called ``triality''. I'll start with a
gentle introduction to the theory of rotation groups; for this, a nice
reference is the book by Fulton and Harris that I mentioned in
\protect\hyperlink{week59}{``Week 59''}. Then I will turn up the heat a
bit and describe triality and how to use it to get the octonions. I
learned some of this stuff from:

\begin{enumerate}
\def\labelenumi{\arabic{enumi})}
\tightlist
\item
  Alex J. Feingold, Igor B. Frenkel, and John F. X. Rees, \emph{Spinor
  construction of vertex operator algebras, triality, and
  \(\mathrm{E}_8^{(1)}\)}, \emph{Contemp.\ Math.} \textbf{121}, AMS, Providence
  Rhode Island, 1991.
\end{enumerate}
I should emphasize, however, that what I will talk about is older, while
the above book starts with triality and then does far more sophisticated
things. An older reference for what I'll talk about is

\begin{enumerate}
\def\labelenumi{\arabic{enumi})}
\setcounter{enumi}{1}
\tightlist
\item
  Claude Chevalley, \emph{The Algebraic Theory of Spinors}, Columbia U.\
  Press, New York, 1954.
\end{enumerate}

\noindent
I think the concept of triality goes back to Cartan, but I don't really
know the history. By the way, I'd really appreciate any corrections to
what I say below.

Okay, so, how should we start? Well, probably we should start with the
group of rotations in \(n\)-dimensional Euclidean space. This group is
called \(\mathrm{SO}(n)\). It is not simply connected if \(n > 1\),
meaning that there are loops in it which cannot be continuously shrunk
to a point. This is easy to see for \(\mathrm{SO}(2)\), which is just
the circle --- or, if you prefer, the unit complex numbers. It's a bit
trickier to see for \(\mathrm{SO}(3)\), but it is easy enough to
demonstrate --- either mathematically or via the famous ``belt trick''
--- that the loop consisting of a 360 degree rotation around an axis
cannot be continuously shrunk to a point, while the loop consisting of a
720 degree rotation around an axis can.

This ``doubly connected'' property of \(\mathrm{SO}(3)\) implies that it
has an interesting ``double cover'', a group \(G\) in which all loops
\emph{can} be contracted to a point, together with a two-to-one function
\(F\colon G \to \mathrm{SO}(3)\) with \(F(gh) = F(g)F(h)\). (This sort
of function, the nice kind of function between groups, is called a
``homomorphism''.) And this double cover \(G\) is just
\(\mathrm{SU}(2)\), the group of \(2\times2\) complex matrices which are
unitary and have determinant \(1\). Better yet --- if we are warming up
for the octonions --- we can think of \(\mathrm{SU}(2)\) as the unit
quaternions!

Now elements of \(\mathrm{SO}(n)\) are just \(n\times n\) real matrices
which are orthogonal and have determinant \(1\), so given an element
\(g\) of \(\mathrm{SO}(n)\) and a vector \(v\) in \(\mathbb{R}^n\), we
can do matrix multiplication to get a new vector \(gv\) in
\(\mathbb{R}^n\), which of course is just the result of rotating \(v\)
by the rotation \(g\). This makes \(\mathbb{R}^n\) into a
``representation'' of \(\mathrm{SO}(n)\), meaning simply that
\[(gh)v = g(hv)\] and \[1v = v.\] We call \(\mathbb{R}^n\) the
``vector'' representation of the rotation group \(\mathrm{SO}(n)\), for
obvious reasons.

Now \(\mathrm{SO}(n)\) has lots of other representations, too. If we
consider \(\mathrm{SO}(3)\), for example, there is in addition to the
vector representation (which is 3-dimensional) also the trivial
\(1\)-dimensional representation (where the group element \(g\) acts on
a complex number \(v\) by leaving it alone!) and also interesting
representations of dimensions 5, 7, 9, etc. The interesting
representation of dimension \(2j+1\) is called the ``spin-\(j\)''
representation by physicists. All representations of \(\mathrm{SO}(3)\)
can be built up from these representations, and none of these
representations can be broken down into smaller ones --- one says they
are irreducible.

But the double cover of \(\mathrm{SO}(3)\), namely \(\mathrm{SU}(2)\),
has more representations! Using the two-to-one homomorphism
\(F\colon \mathrm{SU}(2) \to \mathrm{SO}(3)\) we can convert any
representation of \(\mathrm{SO}(3)\) into one of \(\mathrm{SU}(2)\), but
not vice versa. For example, since \(\mathrm{SU}(2)\) consists of
\(2\times2\) complex matrices, it has a representation on
\(\mathbb{C}^2\), given by the obvious matrix multiplication. This is
called the ``spinor'' or ``spin-\(1/2\)'' representation of
\(\mathrm{SU}(2)\). It doesn't come from a representation of
\(\mathrm{SO}(3)\).

To digress a bit, the reason physicists got so interested in
\(\mathrm{SU}(2)\) is that to describe what happens when you rotate a
particle (in the framework of quantum theory) it turns out you need, not
just the representations of \(\mathrm{SO}(3)\), but of its double cover,
\(\mathrm{SU}(2)\). E.g., an electron, proton or neutron is described by
the spin-\(1/2\) representation. This implies that when you turn an
electron around 360 degrees about an axis, its wavefunction changes
sign, but when you rotate it another 360 degrees, its wavefunction is
back to where it started. You can't describe this behavior using
representations of \(\mathrm{SO}(3)\), but you can using
\(\mathrm{SU}(2)\). In general, for any
\(j = 0, 1/2, 1, 3/2, 2, \ldots\), there is an irreducible
representation of \(\mathrm{SU}(2)\), called the ``spin-\(j\)''
representation, which is \((2j+1)\)-dimensional. Only when the spin is
an integer does the representation come from one of \(\mathrm{SO}(3)\).

Things get more complicated when we consider rotations in higher
dimensional space. For any \(n\) greater than or equal to 3, the group
\(\mathrm{SO}(n)\) is doubly connected, and has a simply connected
double cover, which in general is called \(\mathrm{Spin}(n)\). Folks
have figured out all the representations of \(\mathrm{Spin}(n)\) and
which of these come from representations of \(\mathrm{SO}(n)\). It is
more complicated for \(n > 3\) than for \(n = 3\) (in particular, they
aren't just classified by ``spin''), but it is still quite
comprehensible and charming. Just to head off any confusions that might
occur, let me emphasize that it's sort of a lucky coincidence that
\(\mathrm{Spin}(3) = \mathrm{SU}(2)\). In general, the spin groups don't
have too much to do with the groups \(\mathrm{SU}(n)\) of \(n\times n\)
unitary complex matrices with determinant \(1\).

There is, however, a doubly lucky coincidence in dimension 4; namely,
\(\mathrm{Spin}(4) = \mathrm{SU}(2) \times \mathrm{SU}(2)\). In other
words, an element of \(\mathrm{Spin}(4)\) can be thought of as a pair of
\(\mathrm{SU}(2)\) matrices, and we multiply these pairs like
\((g,g')(h,h') = (gh,g'h')\). This implies that the irreducible
representations of \(\mathrm{Spin}(4)\) are given by a ``tensor
product'' of two irreducible representations of \(\mathrm{SU}(2)\), so
we can classify them by pairs of spins, say \((j,j')\). The dimension of
the \((j,j')\) representation is \((2j+1)(2j'+1)\), since the dimension
of a tensor product is the product of the dimensions. In particular, we
call the \((1/2,0)\) representation the ``left-handed'' spinor
representation and the \((0,1/2)\) representation the ``right-handed''
spinor representation. The reason is that a reflection transforms one
into the other. Since spacetime is \(4\)-dimensional, representations of
\(\mathrm{Spin}(4)\) are important in physics, although really one
should keep track of the fact that time works a bit differently than
space, which \(\mathrm{Spin}(4)\) fails to do. In any event, ignoring
the subtleties about how time works differently than space, we can
roughly say that the existence of left-handed and right-handed spinor
representations of \(\mathrm{Spin}(4)\) is the reason why massless
spin-\(1/2\) particles such as neutrinos can have a ``handedness'' or
``chirality''.

More generally, it turns out that the representation theory of
\(\mathrm{Spin}(n)\) depends strongly on whether \(n\) is even or odd.
When \(n\) is even (and bigger than 2), it turns out that
\(\mathrm{Spin}(n)\) has left-handed and right-handed spinor
representations, each of dimension \(2^{n/2-1}\). When \(n\) is odd
there is just one spinor representation. Of course, there is always the
representation of \(\mathrm{Spin}(n)\) coming from the vector
representation of \(\mathrm{SO}(n)\), which is \(n\)-dimensional.

This leads to something very curious. If you are an ordinary
4-dimensional physicist you undoubtedly tend to think of spinors as
``smaller'' than vectors, since the spinor representations are
2-dimensional, while the vector representation is \(3\)-dimensional.
However, in general, when the dimension \(n\) of space (or spacetime) is
even, the dimension of the spinor representations is \(2^{n/2-1}\),
while that of the vector representation is \(n\), so after a while the
spinor representation catches up with the vector representation and
becomes bigger!

This is a little bit curious, or at least it may seem so at first, but
what's \emph{really} curious is what happens exactly when the spinor
representation catches up with the vector representation. That's when
\(2^{n/2-1} = n\), or \(n = 8\). The group \(\mathrm{Spin}(8)\) has
three \(8\)-dimensional irreducible representations: the vector,
left-handed spinor, and right-handed spinor representation. While they
are not equivalent to each other, they are darn close; they are related
by a symmetry of \(\mathrm{Spin}(8)\) called ``triality''. And, to top
it off, the octonions can be seen as a kind of spin-off of this triality
symmetry... as one might have guessed, from all this
\(8\)-dimensional stuff. One can, in fact, describe the product of
octonions in these terms.

So now let's dig in a bit deeper and describe how this triality business
works. For this, unfortunately, I will need to assume some vague
familiarity with exterior algebras, Clifford algebras, and their
relation to the spin group. But we will have a fair amount of fun
getting our hands on a \(24\)-dimensional beast called the Chevalley
algebra, which contains the vector and spinor representations of
\(\mathrm{Spin}(8)\)!

Start with an \(8\)-dimensional \emph{complex} vector space \(V\) with a
nondegenerate symmetric bilinear form on it. We can think of \(V\) as
the representation of \(\mathrm{SO}(8)\), hence \(\mathrm{Spin}(8)\),
where now I've switched notation and write \(\mathrm{SO}(8)\) to mean
\(\mathrm{SO}(8,\mathbb{C})\), and \(\mathrm{Spin}(8)\) to mean
\(\mathrm{Spin}(8,\mathbb{C})\). We can split \(V\) into two
\(4\)-dimensional subspaces \(V_+\) and \(V_-\) such that
\(\langle v,w\rangle = 0\) whenever \(v\) and \(w\) are either both in
\(V_+\), or both in \(V_-\). Let \(\mathrm{Cliff}\) be the Clifford
algebra over \(V\). Note that as a vector space, there is a natural
identification of \(\mathrm{Cliff}\) with
\[\Lambda V_+ \otimes \Lambda V_-\] where \(\Lambda\) means
``exterior algebra'' and \(\otimes\) means ``tensor product''. (If you
are physicist you may enjoy following Dirac and thinking of
\(\Lambda V_-\) as a Fock space for ``particles'', and \(\Lambda V_+\)
as a Fock space for ``holes''. If you don't enjoy this, ignore it!
We will now to proceed to fill all the holes.) Pick a nonzero vector
\(v\) in \(\Lambda^4 V_+\), the top exterior power of \(V_+\). Let
\(S\) denote the subspace of \(\mathrm{Cliff}\) consisting of all
elements of the form \(uv\) with \(u\) in \(\mathrm{Cliff}\). Note that
\(\mathrm{Cliff}\) and \(S\) are representations of \(\mathrm{Cliff}\)
by left multiplication, and therefore are representations of
\(\mathrm{Spin}(8)\) --- because \(\mathrm{Spin}(8)\) sits inside
\(\mathrm{Cliff}\). (This is a standard way to get one's hands on the
spin groups.)

Note that \(\Lambda V_+\) and \(\Lambda V_-\) both have dimension
\(2^4 = 16\). We can think of both of these as subspaces of
\(\mathrm{Cliff}\); for example, we can think of the vector \(u\) in
\(\Lambda V_-\) as the vector \(1 \otimes u\) in \(\mathrm{Cliff}\).
Note that \(uv = 0\) when \(u\) is in \(\Lambda V_+\). (For
physicists: since the sea of holes is filled, you can't create another.)
Thus \(S\) consists of vectors of the form \(uv\) where \(u\) lies in
\(\Lambda V_-\), and if you think a bit, it follows that \(S\) is
16-dimensional.

So we have our hands on a \(16\)-dimensional representation of
\(\mathrm{Spin}(8)\), namely \(S\). However, we can split it into two
\(8\)-dimensional representations, the left- and right-handed spinor
representations, as follows. Let \[\Lambda^\text{even} V_-\] denote
the part of the exterior algebra consisting of stuff with even degree,
and \[\Lambda^\text{odd} V_-\] the part with odd degree. Then we can
write \(S = S_+ \oplus S_-\), where \(\oplus\) means ``direct sum'' and
\[S_+ = (\Lambda^\text{even} V_-)v , \quad  S_- = (\Lambda^\text{odd} V_-)v.\]
Now, since any element of \(\mathrm{Cliff}\) that's in
\(\mathrm{Spin}(8)\) has even degree in \(\mathrm{Cliff}\), and since
even times even is even, while even times odd is odd, it follows that as
a representation of \(\mathrm{Spin}(8)\), \(S\) splits into \(S_+\) and
\(S_-\), which we call the left-handed and right-handed spinors,
respectively. (Actually I don't know which one is called which, but
being left-handed myself I think the positive one should obviously be
called the left-handed one.)

Note, by the way, that everything so far works quite generally for
\(\mathrm{Spin}(n)\) when \(n\) is even, and it's only in the last
paragraph that I used the fact that \(n\) was even. I certainly haven't
done anything weird using the fact that \(n\) is 8. So as a bonus we're
learning some quite general stuff about spinors.

Now let's do something weird using the fact that \(n\) is 8. We've got
these three \(8\)-dimensional representations of \(\mathrm{Spin}(8)\) on
our hands, namely \(V\), \(S_+\), and \(S_-\). How do they relate?
Recall that \(S_+ + S_- = S\) is a representation of \(\mathrm{Cliff}\),
and since \(V\) sits inside \(\mathrm{Cliff}\) as the elements of degree
\(1\), we have for any \(a\) in \(V\),
\[\mbox{$ab$ is in $S_-$ if $b$ is in $S_+$}\] and
\[\mbox{$ab$ is in $S_+$ if $a$ is in $S_-$}\] If we are in the mood,
this might tempt us to lump \(V\), \(S_+\), and \(S_-\) together to form
a \(24\)-dimensional algebra! Let's call this the Chevalley algebra and
write \[\mathrm{Chev} = V + S_+ + S_-\]

Of course, we need to figure out how to multiply any two guys in
\(\mathrm{Chev}\). We define the product of any two guys in \(V\) to be
zero, and ditto for \(S_+\) or \(S_-\). But we can find an interesting
way to multiply a guy in \(S_+\) by a guy in \(S_-\) to get a guy in
\(V\). I think the vector representation always sits inside the tensor
product of the left- and right-handed spinor representations when space
is even-dimensional, and that this is what we're looking for. But
explicitly, here's what we do in this case. There is a kind of \({}^*\)
operation on \(\mathrm{Cliff}\) given by
\[(abc \cdots def)^* = fed\cdots cba\] where \(a,b,c,\cdots,d,e,f\) lie
in \(V\). This lets us define a symmetric bilinear form on \(S\) by
\[\langle x,y\rangle v = x^* y\] Together with the symmetric bilinear
form we started with on \(V\), this gives us a symmetric bilinear form
on all of \(\mathrm{Chev}\), defining \(\langle a,b \rangle\) to be \(0\) if
\(a\) is in \(V\) and \(b\) is in \(S_+\) or \(S_-\). This bilinear form
on \(\mathrm{Chev}\) turns out to be nondegenerate, and
\(\langle a,b \rangle = 0\) whenever \(a\) and \(b\) lie in different
ones of three summands of \(\mathrm{Chev}\).

So now \(\mathrm{Chev}\) has a nondegenerate symmetric bilinear form on it.
This lets us define a cubic form on \(\mathrm{Chev}\)! Say we have
\((a,b,c)\) in \(V \oplus S_+ \oplus S_- = \mathrm{Chev}\). Then we
define our cubic form \(F\) by \[F(a,b,c) = \langle ab,c \rangle\] using
the fact that we already know how to multiply a guy in \(V\) with a guy
in \(S_+\), and get a guy in \(S_-\).

You probably know --- if you've survived this far! --- that from a
quadratic form you can get a symmetric bilinear form by
``polarization''. Well, similarly, we can get a symmetric trilinear form
\(f\) on \(\mathrm{Chev}\) by polarizing \(F\). Explicitly, for any
\(u_1,u_2,u_3\) in \(\mathrm{Chev}\), we have
\[f(u_1,u_2,u_3) =  F(u_1 + u_2 + u_3) -F(u_1 + u_2) -F(u_2 + u_3) -F(u_1 + u_3) 
+F(u_1) +F(u_2) +F(u_3).\]
Then, since we have a nondegenerate symmetric bilinear form on
\(\mathrm{Chev}\), we can turn \(f\) into a product on
\(\mathrm{Chev}\), by setting
\[\langle u_1 u_2, u_3 \rangle = f(u_1,u_2,u_3).\] 
The assiduous reader can
check that this product on \(\mathrm{Chev}\) agrees with the product we
had partially defined so far; the only new thing it does is define the
product of a guy in \(S_+\) and a guy in \(S_-\), obtaining something in
\(V\). This product turns out to be commutative, but not associative.

Now, if I were really gung-ho about describing triality, I would
describe how the group of permutations of 3 letters, \(S_3\), acts as
automorphisms of \(\mathrm{Chev}\) in a way that lets one scramble the
summands \(V\), \(S_+\), and \(S_-\) at will. In fact, \(S_3\) acts as
automorphisms of \(\mathrm{Spin}(8)\) in a way that gives rise to this
action on \(\mathrm{Chev}\). But right now I'm running out of steam, so
I think I'll just say how to get the octonions out of the Chevalley
algebra!

It's simple: pick a vector \(v\) in \(V\) with
\(\langle v,v \rangle = 1\), and a vector \(s_+\) in \(S_+\) with
\(\langle s_+,s_+ \rangle = 1\). Then \(s_- = vs_+\) is a vector in
\(S_-\) with \(\langle s_-,s_- \rangle = 1\). We now turn \(V\) into the
octonions as follows. Given \(v\) and \(w\) in \(V\), define their
octonion product \(v \ast w\) to be \[v \ast w = (v s_-) (w s_+)\] where the
product on the right hand side is the product in the Chevalley algebra.
In other words: take \(v\) and turn it into something in \(S_+\) by
forming \(v s_-\). Take \(w\) and turn it into something in \(S_-\) by
forming \(w s_+\). The product of these is then something in \(V\). In
short, we form the octonions from the three \(8\)-dimensional
representations of \(\mathrm{Spin}(8)\) by a kind of
ring-around-the-rosie process using triality!

Note: what we just obtained was a \emph{complex} \(8\)-dimensional
algebra, which is the complexification of the octonions. Using the fact
that the vector representation of \(\mathrm{SO}(8,\mathbb{C})\) on
\(\mathbb{C}^8\) contains the vector representation of
\(\mathrm{SO}(8,\mathbb{R})\) on \(\mathbb{R}^8\) as a ``real part'', we
should be able to get the octonions themselves.

One can work out the details following the book of Fulton and Harris,
and the references therein. I should add that they do a lot more fun
stuff involving the exceptional Lie groups and the \(27\)-dimensional
exceptional Jordan algebra... all of this ``exceptional'' stuff
seems to form a unified whole! There is a lot more fun stuff along these
lines in

\begin{enumerate}
\def\labelenumi{\arabic{enumi})}
\setcounter{enumi}{2}
\tightlist
\item
  Ian R. Porteous, \emph{Topological Geometry}, Cambridge U. Press,
  Cambridge, 1981.
\end{enumerate}
In particular, to correct a widespread misimpression, it's worth noting
that there are a lot of nonassociative division algebras over the reals
besides the octonions; Porteous describes one of dimension 4 in his
book. However, all division algebras over \(\mathbb{R}\) are of
dimension 1,2,4, or 8. Also, all normed division algebras over
\(\mathbb{R}\) are the reals, complexes, quaternions, or octonions, and
these are also all the alternative division algebras over
\(\mathbb{R}\), as well... where an ``alternative'' algebra is one
for which any two elements generate an associative algebra. Nota bene:
here a division algebra is one such that for all nonzero \(x\), the map
\(y \mapsto xy\) is invertible. In the finite-dimensional case, this
implies that every element has a left and right inverse. If assume
associativity, the converse is true, but in the nonassociative case it
ain't. Whew! Nonassociative algebras are tricky, if you're used to
associative ones, so if you're interested, you might try this:

\begin{enumerate}
\def\labelenumi{\arabic{enumi})}
\setcounter{enumi}{3}
\tightlist
\item
  R. D. Schafer, \emph{An Introduction to Non-Associative Algebras},
  Dover, New York, 1995.
\end{enumerate}

In addition to the people listed in \protect\hyperlink{week59}{``Week
59''}, I should thank Dan Asimov, Michael Kinyon, Frank Smith, and Dave
Rusin for help with this post. I also thank Doug Merritt for reminding
me about the following nice book on quaternions, octonions, and all
sorts of similar algebras:

\begin{enumerate}
\def\labelenumi{\arabic{enumi})}
\setcounter{enumi}{4}
\tightlist
\item
  I.\ L.\ Kantor and A.\ S.\ \emph{Solodovnikov, Hypercomplex Numbers --- an
  Elementary Introduction to Algebras}, Springer, Berlin, 1989;
  translation of ``Giperkompleksnye chisla'', Moscow, 1973.
\end{enumerate}

\noindent
Back in the old days when there weren't too many algebras around besides
the reals, complexes and quaternions, people called algebras
``hypercomplex numbers''.

\hypertarget{week62}{%
\section{August 28, 1995}\label{week62}}

\hypertarget{week62_ade}{
Now I'd like to talk about a fascinating subject of importance in both
mathematics and physics, the subject of ``ADE classifications''.} Here A,
D, and E aren't abbreviations for anything; they are just names for
certain diagrams. But these diagrams show up all over the place when you
start trying to classify beautiful and symmetrical things.

Let's start with something nice and simple: the Platonic solids. It's
not terribly hard to classify all the regular polyhedra in
\(3\)-dimensional Euclidean space. Roughly, it goes like this. The faces
could all be equilateral triangles. Obviously there need to be at least
3 faces meeting at each vertex to get a polyhedron. If there are exactly
3, you have a tetrahedron. If there are 4, you have an octahedron. If
there are 5, you have an icosahedron. There can't be 6 or more, since
when you have 6 they lie flat in the plane, and more is even worse. The
faces could also be squares. If there are 3 squares meeting at each
vertex you have a cube. There can't be 4 or more, since when you have 4
they lie flat in the plane. The faces could also be regular pentagons.
If there are 3 pentagons meeting at each vertex you have a dodecahedron.
There can't be 4 or more, since when you have 4 you already have more
than 360 degree's worth of angles.

So, there we are: the 5 regular polyhedra are the tetrahedron,
octahedron, icosahedron, cube, and dodecahedron! Of course, we haven't
shown these solids actually exist. Sometimes people forget that you
really need to check that all these possibilities are realized! But the
Greeks did that a while back. This is perhaps the first example of an
ADE classification.

This had such beauty that in his ``Timaeus'' dialog, Plato suggested
that the 4 elements were made of these solids, not counting for the
dodecahedron. Interestingly, Plato considered decomposing the faces of
these solids into ``elementary triangles'', in order to explain how one
element could turn into another. This is presumably why he left out the
dodecahedron: one can't chop up a regular pentagon into 30-60-90
triangles. In a passage that's notoriously hard to translate, he
suggested that the dodecahedron corresponding to some sort of
``quintessence'', or perhaps the zodiac. It's worth pointing out, also,
that Plato explicitly says it's okay if someone comes up with a better
scheme. He makes it clear that he is just trying to lay out an
\emph{example} of a mathematical scheme for explaining the elements, to
get people interested.

Later, of course, Kepler suggested that the 5 Platonic solids
corresponded to the orbits of the 5 planets:
\[\includegraphics[max width=0.65\linewidth]{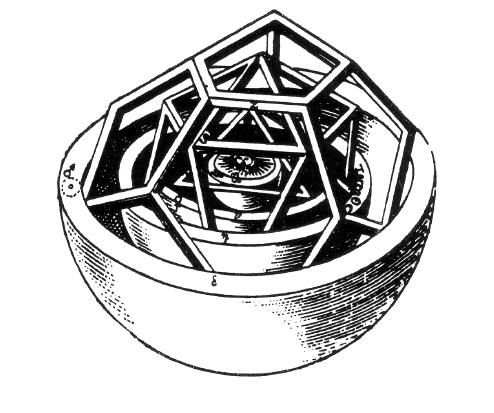}\]
As it turns out, Plato and Kepler were in the right ball-park, but not
really right. Both the solar system and atoms are described pretty well
by similar laws --- the inverse-square force laws for gravity and
electrostatics. And solving this problem (in either the classical or
quantum case) does indeed require a deep understanding of rotations in
3-dimensional space. It's sort of amusing, however, that the Platonic
solids have as their symmetries finite subgroups of the rotation group
in 3 dimensions, while the study of quantum-mechanical atoms instead
involves the theory of ``representations'' of this group, which are in
some sense dual. The rotation group in \(n\) dimensions, by the way, is
called \(\mathrm{SO}(n)\). See \protect\hyperlink{week61}{``Week 61''}
for a bit more about it. For a grand tour of the inverse square law,
both classical and quantum, read:

\begin{enumerate}
\def\labelenumi{\arabic{enumi})}
\tightlist
\item
  Victor Guillemin and Shlomo Sternberg, \emph{Variations on a Theme by
  Kepler}, AMS, Providence, Rhode Island,
  1990.
\end{enumerate}
You will see, among other things, that the real reason the inverse
square force law problem is exactly solvable is that it has a hidden
symmetry under \(\mathrm{SO}(4)\), not just \(\mathrm{SO}(3)\).

But I digress! Recall how I said that ``obviously'' a regular polyhedron
has to have 3 faces meeting at each vertex? What would happen if you
relaxed the definition a little bit, and let there be just 2 faces
meeting at a vertex? Well, then any regular polygon would give a
regular polyhedron, I guess. Then we would have an infinite series of
regular polyhedra with only two faces, together with 5 exceptions, the
Platonic solids. That's actually typical of ADE-type classifications:
often, when you are classifying really symmetrical things, you find some
infinite series of ``obvious'' or ``classical'' cases, together with
finitely many weird ``exceptional'' cases.

Before I get further into ADE classifications, let me note that the
\emph{problem} of why there are so many ADE classifications, and how
they are all related, was explicitly raised by the famous mathematical
physicist V. I. Arnol'd:

\begin{enumerate}
\def\labelenumi{\arabic{enumi})}
\setcounter{enumi}{1}
\tightlist
\item
  Vladimir I.\ Arnol'd, ``Problems of Present Day Mathematics'' in \emph{Mathematical
  Developments Arising from Hilbert's Problems}, ed.~F. E. Browder,
  Proc. Symp. Pure Math. \textbf{28}, AMS,
  Providence, Rhode Island, 1976.
\end{enumerate}
This lists a lot of important math problems, following up on Hilbert's
famous turn-of-the-century listing of problems. Problem VIII in this
book is the ``ubiquity of ADE classifications''. Arnol'd lists the
following examples:

\begin{itemize}
\tightlist
\item
  Platonic solids
\item
  Finite groups generated by reflections
\item
  Weyl groups with roots of equal length
\item
  Representations of quivers
\item
  Singularities of algebraic hypersurfaces with definite intersection
  form
\item
  Critical points of functions having no moduli
\end{itemize}
Don't worry if you don't know what those are except for the first one!
I'll try to explain some of them. Later I'll also explain two new ones
that came out of string theory:

\begin{itemize}
\tightlist
\item
  Minimal models
\item
  Certain ``quantum categories''
\end{itemize}
Perhaps the best single place to start learning about ADE
classifications is:

\begin{enumerate}
\def\labelenumi{\arabic{enumi})}
\setcounter{enumi}{2}
\tightlist
\item
  M. Hazewinkel, W. Hesselink, D. Siermsa, and F. D. Veldkamp, ``The
  ubiquity of Coxeter--Dynkin diagrams (an introduction to the ADE
  problem)'', \emph{Niew.\ Arch.\ Wisk.} \textbf{25} (1977), 257-307.
  Also available at \href{http://math.ucr.edu/home/baez/hazewinkel_et_al.pdf}{\texttt{http://math.ucr.edu/home/baez/hazewinkel\_et\_al.pdf}}
\end{enumerate}

Okay, so what the heck is an ADE classification, after all? It's
probably good to start by looking at ``finite reflection groups.'' Say
we are in \(n\)-dimensional Euclidean space. Then given any unit vector
\(v\), there is a reflection that takes \(v\) to \(-v\), and doesn't do
anything to the vectors orthogonal to \(v\). Let's call this a
``reflection through \(v\)''. A finite reflection group is a finite
group of transformations of Euclidean space such that every element is a
product of reflections. For example, the group of symmetries of a
regular \(n\)-gon is a finite reflection group. (This is a useful
exercise if you don't see it right off the bat.)

Note that if we do two reflections, we get a rotation. In particular,
suppose we have vectors \(v\) and \(w\) at an angle \(\theta\) from each
other, and let \(r\) and \(s\) be the reflections through \(v\) and
\(w\), respectively. Then \(rs\) is a rotation by the angle \(2\theta\). Draw
a picture and check it! This means that if \(\theta = \pi / n\), then
\((rs)^n\) is a rotation by the angle \(2\pi\), which is the same as no
rotation at all, so \((rs)^n = 1\). On the other hand, if \(\theta\) is not a
rational number times \(\pi\), we never have \((rs)^n = 1\), so \(r\)
and \(s\) can not both be in some \emph{finite} reflection group.

With a little more work, we can convince ourselves that any finite
reflection group is captured by a ``Coxeter diagram''. The idea is that
the group is generated by reflections through unit vectors that are all
at angles of \(\pi/n\) from each other. To keep track of things, we draw
a dot for each one of these vectors. Suppose two of the vectors are at
an angle \(\pi/n\) from each other. If \(n = 2\), we don't bother
drawing a line between the two dots. Otherwise, we draw a line between
them, and label it with the number \(n\). Typically, if \(n = 3\) people
don't bother writing the number; they just draw that line. That's what
I'll do. (People also sometimes draw \(n - 2\) lines instead of writing
the number \(n\).)

Algebraically speaking, if someone hands us a Coxeter diagram like \[
  \begin{tikzpicture}
    \draw[thick] (0,0) node{$\bullet$} to (1,0) node{$\bullet$} to node[label=above:{7}]{} (2,0) node {$\bullet$};
  \end{tikzpicture}
\] we get a group having one generator for each dot, and with one
relation \(r^2 = 1\) for each generator \(r\) (since that's what
reflections do), and one relation of the form \((rs)^n = 1\) for each
line connecting dots, or \((rs)^2 = 1\) if there is no line connecting
two dots. It turns out that if a Coxeter diagram yields a \emph{finite}
group this way, it's a finite reflection group.

However, not every diagram we draw yields a finite group! Here are all
the possible Coxeter diagrams giving finite groups. They have names.
First there is \(\mathrm{A}_n\), which has \(n\) dots in a row, like this: \[
  \begin{tikzpicture}
    \draw[thick] (0,0) node{$\bullet$} to (1,0) node{$\bullet$} to (2,0) node {$\bullet$} to (3,0) node{$\bullet$};
  \end{tikzpicture}
\] For example, the group of symmetries of the equilateral triangle is
\(\mathrm{A}_2\). The two dots can correspond to the reflections \(r\) and \(s\)
through two of the altitudes of the triangle, which are at an angle of
\(\pi/3\) from each other. Thus they satisfy \((rs)^3 = 1\). More
generally, \(\mathrm{A}_n\) corresponds to the group of symmetries of an
\(n\)-dimensional simplex --- which is just the group of permutations of
the \(n+1\) vertices.

Then there is \(\mathrm{B}_n\), which has \(n\) dots, where \(n > 1\):
\[
  \begin{tikzpicture}
    \draw[thick] (0,0) node{$\bullet$} to (1,0) node{$\bullet$} to (2,0) node{$\bullet$} to node[label=above:{$4$}]{} (3,0) node {$\bullet$};
  \end{tikzpicture}
\] It has just one edge labelled with a 4. \(\mathrm{B}_n\) turns out to
be the group of symmetries of a hypercube or hyperoctahedron in \(n\)
dimensions.

Then there is \(\mathrm{D}_n\), where \(n > 3\).  This has $n$ dots arranged in a row that
branches at the end:
\[
  \begin{tikzpicture}
    \draw[thick] (0,0) node{$\bullet$} to (1,0) node{$\bullet$} to (2,0) node{$\bullet$} to (3,0) node {$\bullet$};
    \draw[thick] (3,0) to (4,1) node {$\bullet$};
    \draw[thick] (3,0) to (4,-1) node {$\bullet$};
  \end{tikzpicture}
\] 
Then there are \(\mathrm{E}_6\), \(\mathrm{E}_7\), and
\(\mathrm{E}_8\): \[
  \begin{gathered}
    \begin{tikzpicture}
      \draw[thick] (0,0) node{$\bullet$} to (1,0) node{$\bullet$} to (2,0) node{$\bullet$} to (3,0) node {$\bullet$} to (4,0) node {$\bullet$};
      \draw[thick] (2,0) to (2,-1) node{$\bullet$};
    \end{tikzpicture}
\\\begin{tikzpicture}
      \draw[thick] (0,0) node{$\bullet$} to (1,0) node{$\bullet$} to (2,0) node{$\bullet$} to (3,0) node {$\bullet$} to (4,0) node {$\bullet$} to (5,0) node {$\bullet$};
      \draw[thick] (3,0) to (3,-1) node{$\bullet$};
    \end{tikzpicture}
  \\\begin{tikzpicture}
      \draw[thick] (0,0) node{$\bullet$} to (1,0) node{$\bullet$} to (2,0) node{$\bullet$} to (3,0) node {$\bullet$} to (4,0) node {$\bullet$} to (5,0) node {$\bullet$} to (6,0) node {$\bullet$};
      \draw[thick] (4,0) to (4,-1) node{$\bullet$};
    \end{tikzpicture}
  \end{gathered}
\] 
Interestingly, this series does \emph{not} go on. That's what I meant
about ``classical'' versus ``exceptional'' structures.

Then there is \(\mathrm{F}_4\): \[
  \begin{tikzpicture}
    \draw[thick] (0,0) node{$\bullet$} to (1,0) node{$\bullet$} to node[label=above:{4}]{} (2,0) node{$\bullet$} to (3,0) node {$\bullet$};
  \end{tikzpicture}
\] Then there's \(\mathrm{G}_2\): \[
  \begin{tikzpicture}
    \draw[thick] (0,0) node{$\bullet$} to node[label=above:{6}]{} (1,0) node{$\bullet$};
  \end{tikzpicture}
\] and \(\mathrm{H}_3\) and \(\mathrm{H}_4\): \[
  \begin{tikzpicture}
    \draw[thick] (0,0) node{$\bullet$} to node[label=above:{5}]{} (1,0) node{$\bullet$} to (2,0) node{$\bullet$};
  \end{tikzpicture}
  \qquad
  \begin{tikzpicture}
  \draw[thick] (0,0) node{$\bullet$} to node[label=above:{5}]{} (1,0) node{$\bullet$} to (2,0) node{$\bullet$} to (3,0) node {$\bullet$};
\end{tikzpicture}
\] \(\mathrm{H}_3\) is the group of symmetries of the dodecahedron or
icosahedron. \(\mathrm{H}_4\) is the group of symmetries of a regular solid in 4
dimensions which I talked about in \protect\hyperlink{week20}{``Week
20''}. This regular solid is also called the ``unit icosians'' --- it
has 120 vertices, and is a close relative of the icosahedron and
dodecahedron. One amazing thing is that it itself \emph{is} a group in a
very natural way. There are no ``hypericosahedra'' or
``hyperdodecahedra'' in dimensions greater than 4, which is related to
the fact that the \(\mathrm{H}\) series quits at this point.

Finally, there is another infinite series, \(\mathrm{I}_m\): \[
  \begin{tikzpicture}
    \draw[thick] (0,0) node{$\bullet$} to node[label=above:{$m$}]{} (1,0) node{$\bullet$};
  \end{tikzpicture}
\] This corresponds to the symmetry group of the \(2m\)-gon in the
plane, and people usually require \(m = 5\) or \(m > 6\), so as to not
count twice some Coxeter diagrams that we've already run into.

\emph{That's all!}

So, we have an ``\(\mathrm{ABDEFGHI}\) classification'' of finite
reflection groups. (In some future week I had better say what happened
to ``\(\mathrm{C}\)''.) Note that the symmetry groups of the Platonic
solids and some of their higher-dimensional relatives fit in nicely into
this classification, so that's one sense in which the Greeks' discovery
of these solids counts as the first ``\(\mathrm{ADE}\) classification''.
But there is at least one another, deeper, way to fit the Platonic
solids themselves into an \(\mathrm{ADE}\) classification. I'll try to
say more about this in future weeks.

You may still be wondering what's so special about \(\mathrm{A}\),
\(\mathrm{D}\), and \(\mathrm{E}\). I'll have to get to that, too.

\hypertarget{week63}{%
\section{September 14, 1995}\label{week63}}

\hypertarget{week63_ade}{
Let me continue the tale of ``ADE classifications''.} Last week I
described an ``ABDEFGHI classification'' of all finite reflection groups
- that is, finite symmetry groups of Euclidean space, every element of
which is a product of reflections. Now we'll build on that to get other
related classifications.

So, recall:

Every element of a finite reflection group is a product of reflections
through certain special vectors, which people call ``roots''. These
roots are all at angles \(\pi/n\) from each other, where \(n > 1\) is an
integer. To describe the group, we draw a diagram with one dot for each
root. If two roots are perpendicular we don't draw a line between them,
but otherwise, if they are at an angle \(\pi/n\) from each other, we
draw a line and label it with the integer \(n\). Actually, the integer
\(n = 3\) comes up so often that we don't bother labelling the line in
this case.

Now, not all of these diagrams correspond to finite reflection groups.
The following ones, together with disjoint unions of them, are all the
possibilities.

\begin{quote}
\(\mathrm{A}_n\), which has \(n\) dots like this: \[
  \begin{tikzpicture}
    \draw[thick] (0,0) node{$\bullet$} to (1,0) node{$\bullet$} to (2,0) node {$\bullet$} to (3,0) node{$\bullet$};
  \end{tikzpicture}
\] \(\mathrm{B}_n\), which has \(n\) dots, where \(n > 1\): \[
  \begin{tikzpicture}
    \draw[thick] (0,0) node{$\bullet$} to (1,0) node{$\bullet$} to (2,0) node{$\bullet$} to node[label=above:{$4$}]{} (3,0) node {$\bullet$};
  \end{tikzpicture}
\] \(\mathrm{D}_n\), which has \(n\) dots, where \(n > 3\): \[
  \begin{tikzpicture}
    \draw[thick] (0,0) node{$\bullet$} to (1,0) node{$\bullet$} to (2,0) node{$\bullet$} to (3,0) node {$\bullet$};
    \draw[thick] (3,0) to (4,1) node {$\bullet$};
    \draw[thick] (3,0) to (4,-1) node {$\bullet$};
  \end{tikzpicture}
\]
 \(\mathrm{E}_6\), \(\mathrm{E}_7\), and \(\mathrm{E}_8\): 
\[
  \begin{gathered}
    \begin{tikzpicture}
      \draw[thick] (0,0) node{$\bullet$} to (1,0) node{$\bullet$} to (2,0) node{$\bullet$} to (3,0) node {$\bullet$} to (4,0) node {$\bullet$};
      \draw[thick] (2,0) to (2,-1) node{$\bullet$};
    \end{tikzpicture}
\\\begin{tikzpicture}
      \draw[thick] (0,0) node{$\bullet$} to (1,0) node{$\bullet$} to (2,0) node{$\bullet$} to (3,0) node {$\bullet$} to (4,0) node {$\bullet$} to (5,0) node {$\bullet$};
      \draw[thick] (3,0) to (3,-1) node{$\bullet$};
    \end{tikzpicture}
  \\\begin{tikzpicture}
      \draw[thick] (0,0) node{$\bullet$} to (1,0) node{$\bullet$} to (2,0) node{$\bullet$} to (3,0) node {$\bullet$} to (4,0) node {$\bullet$} to (5,0) node {$\bullet$} to (6,0) node {$\bullet$};
      \draw[thick] (4,0) to (4,-1) node{$\bullet$};
    \end{tikzpicture}
  \end{gathered}
\] \(\mathrm{F}_4\): \[
  \begin{tikzpicture}
    \draw[thick] (0,0) node{$\bullet$} to (1,0) node{$\bullet$} to node[label=above:{$4$}]{} (2,0) node{$\bullet$} to (3,0) node {$\bullet$};
  \end{tikzpicture}
\] \(\mathrm{G}_2\): \[
  \begin{tikzpicture}
    \draw[thick] (0,0) node{$\bullet$} to node[label=above:{$6$}]{} (1,0) node{$\bullet$};
  \end{tikzpicture}
\] \(\mathrm{H}_3\) and \(\mathrm{H}_4\): \[
  \begin{tikzpicture}
    \draw[thick] (0,0) node{$\bullet$} to node[label=above:{$5$}]{} (1,0) node{$\bullet$} to (2,0) node{$\bullet$};
  \end{tikzpicture}
  \qquad
  \begin{tikzpicture}
  \draw[thick] (0,0) node{$\bullet$} to node[label=above:{$5$}]{} (1,0) node{$\bullet$} to (2,0) node{$\bullet$} to (3,0) node {$\bullet$};
\end{tikzpicture}
\] \(\mathrm{I}_m\), where \(m = 5\) or \(m > 6\): \[
  \begin{tikzpicture}
    \draw[thick] (0,0) node{$\bullet$} to node[label=above:{$m$}]{} (1,0) node{$\bullet$};
  \end{tikzpicture}
\]
\end{quote}

Recall that \(\mathrm{I}_m\) is the symmetry group of the of regular
\(m\)-gon, while others of these are the symmetry groups of Platonic
solids, and still others are symmetry groups of regular polytopes in
\(n\)-dimensional space. For example, the symmetry group of the
dodecahedron is \(\mathrm{H}_3\), while that of its \(\mathrm{4}\)-dimensional relative is
\(\mathrm{H}_4\).

Now you may know that there are no perfect crystals in the shape of a
regular dodecahedron. However, iron pyrite comes close. In his wonderful
book:

\begin{enumerate}
\def\labelenumi{\arabic{enumi})}
\tightlist
\item
  Hermann Weyl, \emph{Symmetry}, Princeton U.\ Press, Princeton,
  New Jersey, 1989.
\end{enumerate}
\noindent
Weyl suggests that this is how people discovered the regular
dodecahedron:

\begin{quote}
... the discovery of the last two {[}Platonic solids{]} is certainly
one of the most beautiful and singular discoveries made in the whole
history of mathematics. With a fair amount of certainty, it can be
traced to the colonial Greeks in southern Italy. The suggestion has been
made that they abstracted the regular dodecahedron from the crystals of
pyrite, a sulfurous mineral abundant in Sicily.
\end{quote}
\noindent
Thus while iron pyrite is nothing but ``fool's gold'' to the miner, it
may have done a good deed by fooling the Greeks into discovering the
regular dodecahedron. Could this be why the ratio of the diagonal to the
side of a regular pentagon, \((\sqrt{5} + 1)/2\), is called the golden
ratio? Or am I just getting carried away? One is tempted to call the
shape of pyrite crystals the ``fool's dodecahedron,'' but in fact it's
called a ``pyritohedron''. (All this information on pyrite, as well as
the puns, I owe to Michael Weiss.)

More recently, I think people have discovered ``quasicrystals'' (related
to Penrose tiles) having true dodecahedral symmetry. But no perfectly
repetitive crystals form dodecahedra! And the reason is that there is no
lattice having \(H_3\) as its symmetries.

Recall that we get a ``lattice'' by taking \(n\) linearly independent
vectors in \(n\)-dimensional Euclidean space and forming all linear
combinations with integer coefficients. If someone hands us a finite
reflection group, we can look for a lattice having it as symmetries. If
one exists, we say the group satisfies the ``crystallographic
condition''. The only ones that do are
\[\mbox{$\mathrm{A}_n$, $\mathrm{B}_n$, $\mathrm{D}_n$, $\mathrm{E}_6$, $\mathrm{E}_7$, $\mathrm{E}_8$, $\mathrm{F}_4$, and $\mathrm{G}_2$}\]
(and those corresponding to disjoint unions of these diagrams). In other
words, the symmetry groups of the pentagon (\(\mathrm{I}_5\)), the heptagon and
so on (\(\mathrm{I}_m\) with \(m > 6\)), and the dodecahedron and its
\(4\)-dimensional relative (\(\mathrm{H}_3\) and \(\mathrm{H}_4\)) are ruled out.

Now let us turn to the theory of Lie groups. Lie groups are the most
important ``continuous'' (as opposed to discrete) symmetry groups.
Examples include the real line (with addition as the group operation),
the circle (with addition mod \(2\pi\)), and the groups
\(\mathrm{SO}(n)\) and \(\mathrm{SU}(n)\) discussed in
\protect\hyperlink{week61}{``Week 61''}. These groups are incredibly
important in both physics and mathematics. Thus it is wonderful, and
charmingly ironic, that the same diagrams that classify the
oh-so-discrete finite reflection groups also classify some of the most
beautiful of Lie groups: the ``simple'' Lie groups. It turns out that
the simple Lie groups correspond to the diagrams of forms
\(\mathrm{A}\), \(\mathrm{B}\), \(\mathrm{D}\), \(\mathrm{E}\), 
\(\mathrm{F}\), and \(\mathrm{G}\). Oh yes, and \(\mathrm{C}\). I have to
tell you what happened to \(\mathrm{C}\).

There is a vast amount known about semisimple Lie groups, and everyone
really serious about mathematics winds up needing to learn some of this
stuff. I took courses on Lie groups and their Lie algebras in grad
school, but it was only later that I really came to appreciate the
beauty of the simple Lie groups. Back then I found it mystifying because
the work involved in the classification was so algebraic, and I
preferred the more geometrical aspects of Lie groups. Part of the reason
is that the treatment I learned emphasized the Lie algebras and
downplayed the groups. A nice treatment that emphasizes the groups is:

\begin{enumerate}
\def\labelenumi{\arabic{enumi})}
\setcounter{enumi}{1}
\tightlist
\item
  John Frank Adams, \emph{Lectures on Lie Groups}, Benjamin, New York,
  1969.
\end{enumerate}

So what's the basic idea? Let me summarize two semesters of grad school,
and tell you the basic stuff about Lie groups and the classification of
simple Lie groups. Forgive me if it's a bit rushed, sketchy, and even
mildly inaccurate: hopefully the main ideas will shine through the murk
this way.

A Lie group is a group that's also a manifold, for which the group
operations (multiplication and taking inverses) are smooth functions.
This lets you form the tangent space to any point in the group, and the
tangent space at the identity plays a special role. It's called the Lie
algebra of the group. If we have any element \(x\) in the Lie algebra,
we can exponentiate it to get an element \(\exp(x)\) in the group, and
we can keep track of the noncommutativity of the group by forming the
``bracket'' of elements \(x\) and \(y\) in the Lie algebra:

\[[x,y] = \frac{d}{dt}\frac{d}{ds} \exp(sx) \exp(ty) \exp(-sx) \exp(-ty)\]

where \(s\) and \(t\) are real numbers, and we evaluate the derivative
at \(s,t = 0\). Note that \([x,y] = 0\) if the group is commutative.
This bracket operation satisfies some axioms, and algebraists call
anything a Lie algebra that satisfies those axioms. For example, you
could take \(n \times n\) matrices and let \([x,y] = xy - yx\).

Now a Lie algebra is called ``semisimple'' if for any \(z\), there are
\(x\) and \(y\) with \(z = [x,y]\). This is sort of the opposite of an
abelian, or commutative, Lie algebra, where \([x,y] = 0\) for all \(x\)
and \(y\). It turns out that we can take direct sums of Lie algebras by
defining operations componentwise, and it turns out that if you have a
\emph{compact} Lie group, its Lie algebra is always the direct sum of a
semisimple Lie algebra and an abelian one. The abelian ones are pretty
trivial, so all the hard works lies in understanding the semisimple
ones. Any semisimple one is the direct sum of a bunch of semisimple ones
that aren't sums of anything else, and these basic building blocks are
called the ``simple'' ones. They are like the prime numbers of Lie
algebra theory. Unlike the prime numbers, though, we can completely
classify all of them!

Now how does one classify the simple Lie algebras? Basically, it goes
like this. We'll assume our simple Lie algebra is the Lie algebra of a
compact Lie group \(G\) --- it turns out that they all are. Now, sitting
inside \(G\) there is a maximal commutative subgroup \(T\) that's a
torus: a product of a bunch of circles. Let \(\mathrm{Lie}(T)\) stand
for the Lie algebra of this torus \(T\). Now, sitting inside
\(\mathrm{Lie}(T)\) there is a lattice, consisting of all elements \(x\)
with \(\exp(x) = 1\). This is how lattices sneak into the picture!

Moreover, for some elements \(g\) in \(G\), if we ``conjugate'' \(T\) by
\(g\), that is, form the set of all elements \(gtg^{-1}\) where \(t\) is
in \(T\), we get \(T\) back. In other words, these elements of \(g\) act
as symmetries of the torus \(T\). Now, if something acts as symmetries
of something else, it also acts as symmetries of everything naturally
cooked up from that something else. (Roughly speaking, ``naturally''
means "without dependence on arbitrary choices.) For this reason, these
special elements of \(G\) also act as symmetries of \(\mathrm{Lie}(T)\)
and of the lattice sitting inside \(\mathrm{Lie}(T)\). So we have a
lattice together with a group of symmetries, which by the way is called
the Weyl group of \(G\). Now the cool part is that the Weyl group is
actually a finite reflection group, so it must correspond to one of the
diagrams in the list above! Even better, it turns out that the Lie
algebra of \(G\) is determined by the lattice together with its Weyl
group.

The upshot is that to classify semisimple Lie algebras, all we need is
the classification of finite reflection groups satisfying the
crystallographic condition --- which we've done already using diagrams
--- together with a classification of lattices having such finite
reflection groups as symmetries. It turns out that the operation of
taking direct sums of semisimple Lie algebras corresponds to taking
disjoint unions of diagrams, so to get the ``building blocks'' --- the
\emph{simple} Lie algebras --- we only need to worry about the diagrams
we've drawn above, not disjoint unions of them. Now it turns out that
for every type except \(\mathrm{B}\), there is (up to isomorphism) only
\emph{one} lattice having that group of symmetries, but for \(\mathrm{B}\) there
are two. Recall the diagram \(\mathrm{B}_n\) looks like this:
\[
  \begin{tikzpicture}
    \draw[thick] (0,0) node{$\bullet$} to (1,0) node{$\bullet$} to (2,0) node{$\bullet$} to node[label=above:{$4$}]{} (3,0) node {$\bullet$};
  \end{tikzpicture}
\] with \(n\) dots. And recall that the dots correspond to ``roots'',
which in the present context are vectors in \(\mathrm{Lie}(T)\). Now it
turns out that whenever we have a finite reflection group satisfying the
crystallographic condition, we can get a lattice having it as symmetries
by taking integer linear combinations of the roots, but \emph{not}
necessarily roots that are unit vectors; the lengths of the roots
matter. In all cases except \(\mathrm{B}\), there is basically just one way to
get the lengths right, but for \(\mathrm{B}\) there are two. We can keep track of
the root lengths with some extra markings on our diagrams, and then we
get two diagrams, which we call \(\mathrm{B}_n\) and \(\mathrm{C}_n\).
One of them has the root at the right of the diagram being longer, and
one has the root right next to it being longer. This makes no difference
when \(n = 2\), since then we just have \[
  \begin{tikzpicture}
    \draw[thick] (0,0) node{$\bullet$} to node[label=above:{$4$}]{} (1,0) node {$\bullet$};
  \end{tikzpicture}
\] which is perfectly symmetrical. So folks usually consider
\(\mathrm{C}_n\) only for \(n > 2\), to avoid double counting.

In other words, all the simple Lie algebras are of the form:

\begin{itemize}
\tightlist
\item
  \(\mathrm{A}_n\), \(n > 0\)
\item
  \(\mathrm{B}_n\), \(n > 1\)
\item
  \(\mathrm{C}_n\), \(n > 2\)
\item
  \(\mathrm{D}_n\), \(n > 3\)
\item
  \(\mathrm{E}_6\), \(\mathrm{E}_7\), \(\mathrm{E}_8\)
\item
  \(\mathrm{F}_4\)
\item
  \(\mathrm{G}_2\)
\end{itemize}

Okay, so what \emph{are} these things, really? What do they \emph{mean},
and what are the implications of the fact that the symmetries of the
forces of nature are given by the some of the corresponding Lie groups?
Why are 4 infinite series of them and 5 ``exceptional'' Lie algebras?
What's so special about A, D, and E, that makes people keep talking
about ``ADE classifications''? What do the exceptional Lie algebras (and
their corresponding Lie groups) have to do with octonions? Why do some
string theorists think the symmetry group of nature is \(\mathrm{E}_8\),
the largest exceptional Lie group???

Well, I'm afraid that I'm going camping in a couple of hours, so I'll
have to leave you hanging, even though I do know the answers to
\emph{some} of these questions. I'll try to finish talking about ADE
classifications in the next couple of issues.

\begin{center}\rule{0.5\linewidth}{0.5pt}\end{center}

\begin{quote}
\emph{... without fantasy one would never become a mathematician,
and what gave me a place among the mathematicians of our day, despite my
lack of knowledge and form, was the audacity of my thinking.}

--- Sophus Lie
\end{quote}

\hypertarget{week64}{%
\section{September 23, 1995}\label{week64}}

\hypertarget{week64_ade}{
I have been talking about different ``ADE classifications''.} This time
I'll start by continuing the theme of last Week, namely simple Lie
algebras, and then move on to discuss affine Lie algebras and quantum
groups. These are algebraic structures that describe the symmetries
appearing in quantum field theory in 2 and 3 dimensions. They are very
important in string theory and topological quantum field theory, and
it's largely physics that has gotten people interested in them.

Remember, we began by classifying finite reflection groups. A finite
reflection group is simply a finite group of linear transformations of
\(\mathbb{R}^n\), every element of which is a product of reflections.
Finite reflection groups are in 1-1 correspondence with the following
``Coxeter diagrams'', together with disjoint unions of such diagrams:

\begin{quote}
\(\mathrm{A}_n\), which has \(n\) dots like this:  \[
  \begin{tikzpicture}
    \draw[thick] (0,0) node{$\bullet$} to (1,0) node{$\bullet$} to (2,0) node {$\bullet$} to (3,0) node{$\bullet$};
  \end{tikzpicture}
\]\(\mathrm{B}_n\), which has \(n\) dots, where \(n > 1\): \[
  \begin{tikzpicture}
    \draw[thick] (0,0) node{$\bullet$} to (1,0) node{$\bullet$} to (2,0) node{$\bullet$} to node[label=above:{$4$}]{} (3,0) node {$\bullet$};
  \end{tikzpicture}
\] \(\mathrm{D}_n\), which has \(n\) dots, where \(n > 3\): \[
  \begin{tikzpicture}
    \draw[thick] (0,0) node{$\bullet$} to (1,0) node{$\bullet$} to (2,0) node{$\bullet$} to (3,0) node {$\bullet$};
    \draw[thick] (3,0) to (4,1) node {$\bullet$};
    \draw[thick] (3,0) to (4,-1) node {$\bullet$};
  \end{tikzpicture}
\] \(\mathrm{E}_6\), \(\mathrm{E}_7\), and \(\mathrm{E}_8\): 
\[
  \begin{gathered}
    \begin{tikzpicture}
      \draw[thick] (0,0) node{$\bullet$} to (1,0) node{$\bullet$} to (2,0) node{$\bullet$} to (3,0) node {$\bullet$} to (4,0) node {$\bullet$};
      \draw[thick] (2,0) to (2,-1) node{$\bullet$};
    \end{tikzpicture}
\\\begin{tikzpicture}
      \draw[thick] (0,0) node{$\bullet$} to (1,0) node{$\bullet$} to (2,0) node{$\bullet$} to (3,0) node {$\bullet$} to (4,0) node {$\bullet$} to (5,0) node {$\bullet$};
      \draw[thick] (3,0) to (3,-1) node{$\bullet$};
    \end{tikzpicture}
  \\\begin{tikzpicture}
      \draw[thick] (0,0) node{$\bullet$} to (1,0) node{$\bullet$} to (2,0) node{$\bullet$} to (3,0) node {$\bullet$} to (4,0) node {$\bullet$} to (5,0) node {$\bullet$} to (6,0) node {$\bullet$};
      \draw[thick] (4,0) to (4,-1) node{$\bullet$};
    \end{tikzpicture}
  \end{gathered}
\]  \(\mathrm{F}_4\): \[
  \begin{tikzpicture}
    \draw[thick] (0,0) node{$\bullet$} to (1,0) node{$\bullet$} to node[label=above:{$4$}]{} (2,0) node{$\bullet$} to (3,0) node {$\bullet$};
  \end{tikzpicture}
\] \(\mathrm{G}_2\): \[
  \begin{tikzpicture}
    \draw[thick] (0,0) node{$\bullet$} to node[label=above:{$6$}]{} (1,0) node{$\bullet$};
  \end{tikzpicture}
\] \(\mathrm{H}_3\) and \(\mathrm{H}_4\): \[
  \begin{tikzpicture}
    \draw[thick] (0,0) node{$\bullet$} to node[label=above:{$5$}]{} (1,0) node{$\bullet$} to (2,0) node{$\bullet$};
  \end{tikzpicture}
  \qquad
  \begin{tikzpicture}
  \draw[thick] (0,0) node{$\bullet$} to node[label=above:{$5$}]{} (1,0) node{$\bullet$} to (2,0) node{$\bullet$} to (3,0) node {$\bullet$};
\end{tikzpicture}
\] \(\mathrm{I}_m\), where \(m = 5\) or \(m > 6\): \[
  \begin{tikzpicture}
    \draw[thick] (0,0) node{$\bullet$} to node[label=above:{$m$}]{} (1,0) node{$\bullet$};
  \end{tikzpicture}
\]
\end{quote}

Not all of these finite reflection groups satisfy the ``crystallographic
condition'', namely that they act as symmetries of some lattice. The
ones that do are of types A,B,D,E,F, and G, and disjoint unions thereof
--- but I'm going to stop reminding you about disjoint unions all the
time!

Now, if we have a finite reflection group that's the symmetries of some
lattice, we can take the dimension of the lattice to be the number of
dots in the Coxeter diagram. In fact, the dots correspond to a basis of
the lattice, and the lines between them (and their numberings) keep
track of the angles between the basis vectors. These basis vectors are
called ``roots''. To describe the lattice completely, in principle we
need to know the lengths of the roots as well as the angles between
them. But it turns out that except for type B, there is just one choice
of lengths that works (up to overall scale). For type B there are two
choices, which people call \(\mathrm{B}_n\) and \(\mathrm{C}_n\),
respectively. People keep track of the lengths with a ``Dynkin diagram''
like this:

\begin{itemize}
\tightlist
\item
  \(\mathrm{B}_n\) has \(n\) dots, where \(n>1\): 
\[
      \begin{tikzpicture}
        \draw[thick] (0,0) node{$\bullet$} to (1,0) node{$\bullet$} to (2,0) node{$\bullet$} to node[label=above:{$4$}]{\textgreater} (3,0) node {$\bullet$};
      \end{tikzpicture}
    \]
\item
  \(\mathrm{C}_n\) has \(n\) dots, where \(n>2\): \[
      \begin{tikzpicture}
        \draw[thick] (0,0) node{$\bullet$} to (1,0) node{$\bullet$} to (2,0) node{$\bullet$} to node[label=above:{$4$}]{\textless} (3,0) node {$\bullet$};
      \end{tikzpicture}
    \]
\end{itemize}

The arrow points to the shorter root; for \(\mathrm{B}_n\) all the roots
except the last one are \(\sqrt{2}\) times as long as the last one,
while for \(\mathrm{C}_n\) all the roots except the last one are
\(1/\sqrt{2}\) as long. (In fact, the lattices corresponding to
\(\mathrm{B}_n\) and \(\mathrm{C}_n\) are ``dual'', in the hopefully
obvious sense.) The only reason why we require \(n > 2\) for
\(\mathrm{C}_n\) is that \(\mathrm{C}_2\) is basically the same as 
\(\mathrm{B}_2\)!

Now last Week, I also sketched how the Lie algebras of the compact
simple Lie groups were \emph{also} classified by the same Dynkin
diagrams of types \(\mathrm{A}, \mathrm{B}, \mathrm{C}, \mathrm{D}, 
\mathrm{E}, \mathrm{F},\) and \(\mathrm{G}\). These were real Lie algebras;
we can also switch viewpoint and work with complex Lie algebras if we
like, in which case we simply say we're studying the complex simple Lie
algebras, as opposed to their ``compact real forms''.

Unfortunately, I didn't say much about what these Lie algebras actually
are! Basically, they go like this:

\(\mathrm{A}_n\) --- The Lie algebra \(\mathrm{A}_n\) is just
\(\mathfrak{sl}_{n+1}(\mathbb{C})\), the \((n+1) \times (n+1)\) complex
matrices with vanishing trace, which form a Lie algebra with the usual
bracket \([x,y] = xy -yx\). The compact real form of
\(\mathfrak{sl}_n(\mathbb{C})\) is \(\mathfrak{su}_n\), and the
corresponding compact Lie group is \(\mathrm{SU}(n)\), the \(n\times n\)
unitary matrices with determinant \(1\). The symmetry group of the
electroweak force is \(\mathrm{U}(1) \times \mathrm{SU}(2)\), where
\(\mathrm{U}(1)\) is the \(1 \times 1\) unitary matrices. The symmetry
group of the strong force is \(\mathrm{SU}(3)\). The study of
\(\mathrm{A}_n\) is thus a big deal in particle physics. People have
also considered ``grand unified theories'' with symmetry groups like
\(\mathrm{SU}(5)\).

\(\mathrm{B}_n\) --- The Lie algebra \(\mathrm{B}_n\) is
\(\mathfrak{so}_{2n+1}(\mathbb{C})\), the \((2n+1) \times (2n+1)\)
skew-symmetric complex matrices with vanishing trace. The compact real
form of \(\mathfrak{so}_n(\mathbb{C})\) is \(\mathfrak{so}_n\), and the
corresponding compact Lie group is \(\mathrm{SO}(n)\), the
\(n \times n\) real orthogonal matrices with determinant \(1\), that is,
the rotation group in Euclidean \(n\)-space. For some basic cool facts
about \(\mathrm{SO}(n)\), check out \protect\hyperlink{week61}{``Week
61''}.

\(\mathrm{C}_n\) --- The Lie algebra \(\mathrm{C}_n\) is
\(\mathfrak{sp}_n(\mathbb{C})\), the \(2n \times 2n\) complex matrices
of the form \[
  \left(
    \begin{array}{cc}
      A&B\\C&D
    \end{array}
  \right)
\] where \(B\) and \(C\) are symmetric, and \(D\) is minus the transpose
of \(A\). The compact real form of \(\mathfrak{sp}_n(\mathbb{C})\) is
\(\mathfrak{sp}_n\), and the corresponding compact Lie group is called
\(\mathrm{Sp}(n)\). This is the group of \(n \times n\) quaternionic
matrices which preserve the usual inner product on the space
\(\mathbb{H}^n\) of \(n\)-tuples of quaternions.

\(\mathrm{D}_n\) --- The Lie algebra \(\mathrm{D}_n\) is
\(\mathfrak{so}_{2n}(\mathbb{C})\), the \(2n \times 2n\) skew-symmetric
complex matrices with vanishing trace. See \(\mathrm{B}_n\) above for
more about this. It may seem weird that \(\mathrm{SO}(n)\) acts so
differently depending on whether \(n\) is even or odd, but it's true:
for example, there are ``left-handed'' and ``right-handed'' spinors in
even dimensions, but not in odd dimensions. Some clues as to why are
given in \protect\hyperlink{week61}{``Week 61''}.

Those are the ``classical'' Lie algebras, and they are things that are
pretty easy to reinvent for yourself, and to get interested in for all
sorts of reasons. As you can see, they are all about ``rotations'' in
real, complex, and quaternionic vector spaces.

The remaining ones are called ``exceptional'', and they are much more
mysterious. They were only discovered when people figured out the
classification of simple Lie algebras. As it turns out, they tend to be
related to the octonions! Some other week I will say more about them,
but for now, let me just say:

\(\mathrm{F}_4\) --- This is a 52-dimensional Lie algebra. Its smallest
representation is \(26\)-dimensional, consisting of the traceless
\(3\times3\) hermitian matrices over the octonions, on which it
preserves a trilinear form.

\(\mathrm{G}_2\) --- This is a \(14\)-dimensional Lie algebra, and the
compact Lie group corresponding to its compact real form is also often
called \(\mathrm{G}_2\). This group is just the group of symmetries
(automorphisms) of the octonions! In fact, the smallest representation
of this Lie algebra is 7-dimensional, corresponding to the purely
imaginary octonions.

\(\mathrm{E}_6\) --- This is a 78-dimensional Lie algebra. Its smallest
representation is \(27\)-dimensional, consisting of all the \(3\times3\)
hermitian matrices over the octonions this time, on which it preserves
the anticommutator.

\(\mathrm{E}_7\) --- This is a 133-dimensional Lie algebra. Its smallest
representation is 56-dimensional, on which it preserves a tetralinear
form.

\(\mathrm{E}_8\) --- This is a 248-dimensional Lie algebra, the biggest
of the exceptional Lie algebras. Its smallest representation is
248-dimensional, the so-called ``adjoint'' representation, in which it
acts on itself. Thus in some vague sense, the simplest way to understand
the Lie group corresponding to \(\mathrm{E}_8\) is as the symmetries of
itself! (Thanks go to Dick Gross for this charming information; I think
he said all the other exceptional Lie algebras have representations
smaller than themselves, but I forget the sizes.) In
\protect\hyperlink{week20}{``Week 20''} I described a way to get its
root lattice (the \(8\)-dimensional lattice spanned by the roots) by
playing around with the icosahedron and the quaternions.

People have studied simple Lie algebras a lot this century, basically
studied the hell out of them, and in fact some people were getting a
teeny bit sick of it recently, when along came some new physics that put
a lot of new life into the subject. A lot of this new physics is related
to string theory and quantum gravity. Some of this physics is
``conformal field theory'', the study of quantum fields in 2-dimensional
spacetime that are invariant under all conformal (angle-preserving)
transformations. This is important in string theory because the string
worldsheet is \(2\)-dimensional. Some other hunks of this physics go by
the name of ``topological quantum field theory'', which is the study of
quantum fields, usually in 3 dimensions so far, that are invariant under
\emph{all} transformations (or more precisely, all diffeomorphisms).

Every simple Lie algebra gives rise to an ``affine Lie algebra'' and a
``quantum group''. The symmetries of conformal field theories are
closely related to affine Lie algebras, and the symmetries of
topological quantum field theories are quantum groups. I won't tell you
what affine Lie algebras and quantum groups ARE, since it would take
quite a while. Instead I'll refer you to a good good introduction to
this stuff:

\begin{enumerate}
\def\labelenumi{\arabic{enumi})}
\tightlist
\item
  J\"urgen Fuchs, \emph{Affine Lie Algebras and Quantum Groups},
   Cambridge U.\ Press, Cambridge 1992.
\end{enumerate}

Let me whiz through his table of contents and very roughly sketch what
it's all about.

\begin{enumerate}
\def\labelenumi{\arabic{enumi}.}
\item
  \textbf{Semisimple Lie algebras}

  This is a nice summary of the theory of semisimple Lie algebras
  (remember, those are just direct sums of simple Lie algebras) and
  their representations. Especially if you are a physicist, a slick
  summary like this might be a better way to start learning the subject
  than a big fat textbook on the subject. He covers the Dynkin diagram
  stuff and much, much more.
\item
  \textbf{Affine Lie algebras}

  This starts by describing Kac--Moody algebras, which are certain
  \emph{infinite-dimensional} analogs of the simple Lie algebras. Fuchs
  concentrates on a special class of these, the affine Lie algebras, and
  describes the classification of these using Dynkin diagrams. The main
  nice thing about the affine Lie algebras is that their corresponding
  infinite-dimensional Lie groups are very nice: they are almost ``loop
  groups''. If we have a Lie group \(G\), the loop group \(LG\) is just
  the set of all smooth functions from the circle to \(G\), which we
  make into a group by pointwise multiplication. If you're a physicist,
  this is obviously relevant to string theory, because at each time, a
  string is just a circle (or bunch of circles), and if you are doing
  gauge theory on the string, with symmetry group \(G\), the gauge group
  is then just the loop group \(LG\). So you'd expect the representation
  theory of loop groups and their Lie algebras to be really important.

  You'd \emph{almost} be right, but there is a slight catch. In quantum
  theory, what counts are the ``projective'' representations of a group,
  that is, representations that satisfy the rule \(g(h(v)) = (gh)(v)\)
  \emph{up to a phase}. (This is because ``phases are unobservable in
  quantum theory'' --- one of those mottos that needs to be carefully
  interpreted to be correct.) The projective representations of the loop
  group \(LG\) correspond to the honest-to-goodness representations of a
  ``central extension'' of \(LG\), a slightly fancier group than \(LG\)
  itself. And the Lie algebra of \emph{this} group is called an affine
  Lie algebra.

  So, people who like gauge theory and string theory need to know a lot
  about affine Lie algebras and their representations, and that's what
  this chapter covers. A real heavy-duty string theorist will need to
  know more about Kac--Moody algebras, so if you're planning on becoming
  one of those, you'd better also try:

  \begin{enumerate}
  \def\labelenumii{\arabic{enumii})}
  \setcounter{enumii}{1}
  \tightlist
  \item
    Victor Kac, \emph{Infinite Dimensional Lie Algebras},
    Cambridge U.\ Press, Cambridge, 1990.
  \end{enumerate}

  You'll also need to know more about loop groups, so try:

  \begin{enumerate}
  \def\labelenumii{\arabic{enumii})}
  \setcounter{enumii}{2}
  \tightlist
  \item
    Andrew Pressley and Graeme Segal, \emph{Loop Groups}, Oxford
    U.\ Press, Oxford, 1986.
  \end{enumerate}
\item
  \textbf{WZW theories}

  Well, I just said that physicists liked affine Lie algebras because
  they were the symmetries of conformal field theories that were also
  gauge theories. Guess what: a Wess--Zumino--Witten, or WZW, theory, is a
  conformal field theory that's also a gauge theory! You can think of it
  as the natural generalization of the wave equation in 2 dimensions
  (which is conformally invariant, btw) from the case of real-valued
  fields to general \(G\)-valued fields, where \(G\) is our favorite
  Lie group.
\item
  \textbf{Quantum groups}

  When you quantize a WZW theory whose symmetry group \(G\) is some
  simple Lie group, something funny happens. In a sense, the group
  itself also gets quantized! In other words, the algebraic structure of
  the group, or its Lie algebra, gets ``deformed'' in a way that depends
  on the parameter \(\hbar\) (Planck's constant). I have muttered much
  about quantum groups on This Week's Finds, especially concerning their
  relevance to topological quantum field theory, and I will not try to
  explain them any better here! Eventually I will discuss a bunch of
  books that have come out on quantum groups, and I hope to give a
  mini-introduction to the subject in the process.
\item
  \textbf{Duality, fusion rules, and modular invariance}

  The previous chapter described quantum groups as abstract algebraic
  structures, showing how you can get one from any simple Lie algebra.
  Here Fuchs really shows how you get them from quantizing a WZW theory.
  WZW theories are invariant under conformal transformations, and
  quantum groups inherit lots of cool properties from this fact. For
  example, suppose you form a torus by taking the complex plane and
  identifying two points if they differ by any number of the form
  \(n z_1 + m z_2\), where \(z_1\) and \(z_2\) are fixed complex numbers
  and \(n\), \(m\) are arbitrary integers. For example, we might
  identify all these points: \[
     \begin{tikzpicture}[scale=0.7]
       \draw[->] (-3,0) to (4,0) node[label=below:{$\mathrm{Re}(z)$}]{};
       \draw[->] (0,-3) to (0,4) node[label=left:{$\mathrm{Im}(z)$}]{};
       \foreach \m in {-1,0,1,2}
       {
         \foreach \n in {-1,0,1,2}
         {
           \node at ({\m*1.5-\n/3-0.2},{1.5*\n+\m-0.5}) {$\bullet$};
         }
       }
     \end{tikzpicture}
   \] The resulting torus is a ``Riemann surface'' and it has lots of
  transformations, called ``modular transformations''. The group of
  modular transformations is the discrete group
  \(\mathrm{SL}(2,\mathbb{Z})\) of \(2\times2\) integer matrices with
  determinant \(1\); I leave it as an easy exercise to guess how these
  give transformations of the torus. (This is an example of a ``mapping
  class group'' as discussed in \protect\hyperlink{week28}{``Week
  28''}.) In any event, the way the the WZW theory transforms under
  modular transformations translates into some cool properties of the
  corresponding quantum group, which Fuchs discusses. That's roughly
  what ``modular invariance'' means.

  Similarly, ``fusion rules'' have to do with the thrice-punctured
  sphere, or ``trinion'', which is another Riemann surface. And
  ``duality'' has to do with the sphere with four punctures, which can
  be viewed as built up from trinions in either of two ``dual'' ways: \[
     \begin{tikzpicture}[scale=0.3,rotate=90]
       \begin{scope}
         \draw[thick] (-3,0) ellipse (2cm and 1cm);
         \draw[thick] (3,0) ellipse (2cm and 1cm);
         \draw[thick] (-5,0) .. controls (-5,-2) and (-2,-4) .. (-2,-6);
         \draw[thick] (5,0) .. controls (5,-2) and (2,-4) .. (2,-6);
         \draw[thick] (-1,0) .. controls (-1,-1) .. (0,-2);
         \draw[thick] (1,0) .. controls (1,-1) .. (0,-2);
         \draw[thick] (-2,-6) to (-2,-7);
         \draw[thick] (2,-6) to (2,-7);
       \end{scope}
       \begin{scope}[rotate=180,shift={(0,14)}]
         \begin{scope}[shift={(-3,0)},rotate=180]
           \draw[thick,dashed] (0:2) arc (0:180:2cm and 1cm);
           \draw[thick] (180:2) arc (180:360:2cm and 1cm);
         \end{scope}
         \begin{scope}[shift={(3,0)},rotate=180]
           \draw[thick,dashed] (0:2) arc (0:180:2cm and 1cm);
           \draw[thick] (180:2) arc (180:360:2cm and 1cm);
         \end{scope}
         \draw[thick] (-5,0) .. controls (-5,-2) and (-2,-4) .. (-2,-6);
         \draw[thick] (5,0) .. controls (5,-2) and (2,-4) .. (2,-6);
         \draw[thick] (-1,0) .. controls (-1,-1) .. (0,-2);
         \draw[thick] (1,0) .. controls (1,-1) .. (0,-2);
         \draw[thick] (-2,-6) to (-2,-7);
         \draw[thick] (2,-6) to (2,-7);
       \end{scope}
     \end{tikzpicture}
   \] or \[
     \begin{tikzpicture}[scale=0.3]
       \begin{scope}
         \draw[thick] (-3,0) ellipse (2cm and 1cm);
         \draw[thick] (3,0) ellipse (2cm and 1cm);
         \draw[thick] (-5,0) .. controls (-5,-2) and (-2,-4) .. (-2,-6);
         \draw[thick] (5,0) .. controls (5,-2) and (2,-4) .. (2,-6);
         \draw[thick] (-1,0) .. controls (-1,-1) .. (0,-2);
         \draw[thick] (1,0) .. controls (1,-1) .. (0,-2);
         \draw[thick] (-2,-6) to (-2,-7);
         \draw[thick] (2,-6) to (2,-7);
       \end{scope}
       \begin{scope}[rotate=180,shift={(0,14)}]
         \begin{scope}[shift={(-3,0)},rotate=180]
           \draw[thick,dashed] (0:2) arc (0:180:2cm and 1cm);
           \draw[thick] (180:2) arc (180:360:2cm and 1cm);
         \end{scope}
         \begin{scope}[shift={(3,0)},rotate=180]
           \draw[thick,dashed] (0:2) arc (0:180:2cm and 1cm);
           \draw[thick] (180:2) arc (180:360:2cm and 1cm);
         \end{scope}
         \draw[thick] (-5,0) .. controls (-5,-2) and (-2,-4) .. (-2,-6);
         \draw[thick] (5,0) .. controls (5,-2) and (2,-4) .. (2,-6);
         \draw[thick] (-1,0) .. controls (-1,-1) .. (0,-2);
         \draw[thick] (1,0) .. controls (1,-1) .. (0,-2);
         \draw[thick] (-2,-6) to (-2,-7);
         \draw[thick] (2,-6) to (2,-7);
       \end{scope}
     \end{tikzpicture}
   \] This is one of the reasons string theory was first discovered ---
  we can think of the above pictures as two Feynman diagrams for
  interacting strings, and the fact that they are really just distorted
  versions of each other gives rise to identities among Feynman
  diagrams. Similarly, this fact gives rise to identities satisfied by
  the fusion rules of quantum groups.
\end{enumerate}
\noindent
So --- Fuchs' book should make clear how, starting from the austere
beauty of the Dynkin diagrams, we get not only simple Lie groups, but
also a wealth of more complicated structures that people find important
in modern theoretical physics.

\begin{center}\rule{0.5\linewidth}{0.5pt}\end{center}

\begin{quote}
\emph{Mathematics, rightly viewed, possesses not only truth, but supreme
beauty --- a beauty cold and austere, like that of sculpture, without
appeal to any part of our weaker nature, without the gorgeous trappings
of painting or music, yet sublimely pure, and capable of a stern
perfection such as only the greatest art can show.}

--- Bertrand Russell.
\end{quote}

\hypertarget{week65}{%
\section{October 3, 1995}\label{week65}}

\hypertarget{week65_ade}{
This time I'll finish up talking about ``ADE classifications'' for a
while, although there is certainly more to say.} Recall where we were:
the following diagrams correspond to the simple Lie algebras, and they
also define certain lattices, the ``root lattices'' of those Lie
algebras:

\begin{quote}
\(\mathrm{A}_n\), which has \(n\) dots like this:
 \[
  \begin{tikzpicture}
    \draw[thick] (0,0) node{$\bullet$} to (1,0) node{$\bullet$} to (2,0) node {$\bullet$} to (3,0) node{$\bullet$};
  \end{tikzpicture}
\] 
\(\mathrm{B}_n\), which has \(n\) dots, where \(n > 1\): 
\[
      \begin{tikzpicture}
        \draw[thick] (0,0) node{$\bullet$} to (1,0) node{$\bullet$} to (2,0) node{$\bullet$} to node[label=above:{$4$}]{\textgreater} (3,0) node {$\bullet$};
      \end{tikzpicture}
    \]
\(\mathrm{C}_n\), which has \(n\) dots, where \(n > 2\): 
\[
      \begin{tikzpicture}
        \draw[thick] (0,0) node{$\bullet$} to (1,0) node{$\bullet$} to (2,0) node{$\bullet$} to node[label=above:{$4$}]{\textless} (3,0) node {$\bullet$};
      \end{tikzpicture}
    \]
 \(\mathrm{D}_n\), which has \(n\) dots, where \(n > 3\): \[
  \begin{tikzpicture}
    \draw[thick] (0,0) node{$\bullet$} to (1,0) node{$\bullet$} to (2,0) node{$\bullet$} to (3,0) node {$\bullet$};
    \draw[thick] (3,0) to (4,1) node {$\bullet$};
    \draw[thick] (3,0) to (4,-1) node {$\bullet$};
  \end{tikzpicture}
\]  \(\mathrm{E}_6\), \(\mathrm{E}_7\), and \(\mathrm{E}_8\): 
 \[
  \begin{gathered}
    \begin{tikzpicture}
      \draw[thick] (0,0) node{$\bullet$} to (1,0) node{$\bullet$} to (2,0) node{$\bullet$} to (3,0) node {$\bullet$} to (4,0) node {$\bullet$};
      \draw[thick] (2,0) to (2,-1) node{$\bullet$};
    \end{tikzpicture}
\\\begin{tikzpicture}
      \draw[thick] (0,0) node{$\bullet$} to (1,0) node{$\bullet$} to (2,0) node{$\bullet$} to (3,0) node {$\bullet$} to (4,0) node {$\bullet$} to (5,0) node {$\bullet$};
      \draw[thick] (3,0) to (3,-1) node{$\bullet$};
    \end{tikzpicture}
  \\\begin{tikzpicture}
      \draw[thick] (0,0) node{$\bullet$} to (1,0) node{$\bullet$} to (2,0) node{$\bullet$} to (3,0) node {$\bullet$} to (4,0) node {$\bullet$} to (5,0) node {$\bullet$} to (6,0) node {$\bullet$};
      \draw[thick] (4,0) to (4,-1) node{$\bullet$};
    \end{tikzpicture}
  \end{gathered}
\]   \(\mathrm{F}_4\): \[
  \begin{tikzpicture}
    \draw[thick] (0,0) node{$\bullet$} to (1,0) node{$\bullet$} to node[label=above:{$4$}]{\textgreater} (2,0) node{$\bullet$} to (3,0) node {$\bullet$};
  \end{tikzpicture}
\] \(\mathrm{G}_2\): \[
  \begin{tikzpicture}
    \draw[thick] (0,0) node{$\bullet$} to node[label=above:{$6$}]{\textgreater} (1,0) node{$\bullet$};
  \end{tikzpicture}
\]
\end{quote}

The dots in one of these ``Dynkin diagrams'' correspond to certain set
of basis vectors, or ``roots'', of the lattice. The lines, with their
decorative numbers and arrows, give enough information to recover the
lattice from the diagram. In particular, two dots that are not connected
by a line correspond to roots that are at a 90 degree angle from each
other, while two dots connected by an unnumbered line correspond to
roots that are at a 60 degree angle from each other. Numbered lines mean
the angle between roots is something else, and the arrows point from the
longer to the shorter root in this case, as partially explained in
\protect\hyperlink{week63}{``Week 63''}.

However, we will now concentrate on the cases  \(\mathrm{A}, \mathrm{D},\) 
and  \(\mathrm{E}\), where all the
roots are 90 or 60 degrees from each other, and they are all the same
length --- usually taken to be length 2. These are the ``simply laced''
Dynkin diagrams. I want to explain what's so special about them! But
first, I should describe the corresponding lattices more explicitly, to
make it clear how simple they really are.

The following information, and more, can be found in Chapter 4 of:

\begin{enumerate}
\def\labelenumi{\arabic{enumi})}
\tightlist
\item
  J. H. Conway and N. J. A.
  Sloane, \emph{Sphere Packings, Lattices and Groups}, 
 Springer, Berlin, 1993.
\end{enumerate}
\noindent
which I described in more detail in \protect\hyperlink{week20}{``Week
20''}.

So, what are \(\mathrm{A}, \mathrm{D},\) and \(\mathrm{E}\) like?

\textbf{A}. We can describe the lattice \(\mathrm{A}_n\) as the set of
all \((n+1)\)-tuples of integers \((x_1,\ldots,x_{n+1})\) such that
\[x_1+\cdots+x_{n+1}=0.\] It's a fun exercise to show that \(\mathrm{A}_2\) is a
\(2\)-dimensional hexagonal lattice, the sort of lattice you use to pack
pennies as densely as possible. Similarly, \(\mathrm{A}_3\) gives a standard way
of packing grapefruit, which is in fact the densest lattice packing of
spheres in 3 dimensions. (Hsiang has claimed to have shown it's the
densest packing, lattice or not, but this remains controversial.)

\textbf{D}. We can describe \(\mathrm{D}_n\) as the set of all
\(n\)-tuples of integers \((x_1,\ldots,x_n)\) such that
\[x_1+\cdots+x_n\quad\text{is even}.\] Or, if you like, you can imagine
taking an \(n\)-dimensional checkerboard, coloring the cubes alternately
red and black, and taking the center of each red cube. In four
dimensions, \(\mathrm{D}_4\) gives a denser packing of spheres than 
\(\mathrm{A}_4\); in
fact, it gives the densest lattice packing possible. Moreover, \(\mathrm{D}_5\)
gives the densest lattice packing of in dimension 5. However, in
dimensions 6, 7, and 8, the \(\mathrm{E}_n\) lattices are the best!

\textbf{E}. We can describe \(\mathrm{E}_8\) as the set of 8-tuples
\((x_1,\ldots,x_8)\) such that the \(x_i\) are either all integers or all
half-integers --- a half-integer being an integer plus \(1/2\) --- and
\[x_1+\cdots+x_8\quad\text{is even}.\] Each point has 240 closest
neighbors. Alternatively, as described in
\protect\hyperlink{week20}{``Week 20''}, we can describe
\(\mathrm{E}_8\) in a slick way in terms of the quaternions. Another
neat way to think of \(\mathrm{E}_8\) is in terms of the octonions! Not
too surprising, perhaps, since the octonions and \(\mathrm{E}_8\) are
both \(8\)-dimensional, but it's still sorta neat. For the details,
check out

\begin{enumerate}
\def\labelenumi{\arabic{enumi})}
\setcounter{enumi}{1}
\tightlist
\item
  Geoffrey Dixon, ``Octonion X-product and \(\mathrm{E}_8\) lattices'',
  available as
  \href{https://arxiv.org/abs/hep-th/9411063}{\texttt{hep-th/9411063}}.
\end{enumerate}

Briefly, this goes as follows. In \protect\hyperlink{week59}{``Week
59''} we described a multiplication table for the ``seven dwarves'' ---
a basis of the imaginary octonions --- but there are lots of other
multiplication tables that would also give an algebra isomorphic to the
octonions. Given any unit octonion \(a\), we can define an ``octonion
\(\times\)-product'' as follows: \[b \times c = (b a)(a^* c)\] where
\(a^*\) is the conjugate of \(a\) (as defined in
\protect\hyperlink{week59}{``Week 59''}) and the product on the
right-hand side is the usual octonion product, parenthesized because it
ain't associative. For exactly 480 choices of the unit octonion \(a\),
the \(\times\)-product gives us a new multiplication table for the seven
dwarves, such that we get an algebra isomorphic to the octonions again!
240 of these choices have all rational coordinates (in terms of the
seven dwarves), and these are precisely the 240 closest neighbors of the
origin in a copy of the \(\mathrm{E}_8\) lattice! The other 240 have all
irrational coordinates, and these are the closest neighbors to the
origin of a \emph{different} copy of the \(\mathrm{E}_8\) lattice. (Here
we've rescaled the \(\mathrm{E}_8\) lattice so the nearest neighbors
have distance \(1\) from the origin, instead of \(\sqrt{2}\) as above.)

Once we have \(\mathrm{E}_8\) in hand, we can get its little pals
\(\mathrm{E}_7\) and \(\mathrm{E}_6\) as follows. To get
\(\mathrm{E}_7\), just take all the vectors in \(\mathrm{E}_8\) that are
perpendicular to some closest neighbor of the origin. To get
\(\mathrm{E}_6\), find a copy of the lattice \(\mathrm{A}_2\) in \(\mathrm{E}_8\)
(there are lots) and then take all the vectors in \(\mathrm{E}_8\)
perpendicular to everything in that copy of \(\mathrm{A}_2\).

So, now that we have a nodding acquaintance with \(\mathrm{A}, 
\mathrm{D},\)  and \(\mathrm{E}\), let me
describe some of the many places they show up. First, what's so great
about these lattices, apart from the fact that they're the root lattices
of simple Lie algebras, with a special ``simply-laced'' property? I
don't think I really understand the answer to this in a deep way, but I
know various things to say!

First, Witt's theorem says that the \(\mathrm{A}, \mathrm{D},\) and 
\(\mathrm{E}\) lattices and their
direct sums are the only integral lattices having a basis consisting of
vectors \(v\) with \(\|v\|^2 = 2\). Here a lattice is ``integral'' if
the dot product of any two vectors in it is an integer. In fact, any
integral lattice having a basis consisting of vectors with \(\|v\|^2\)
equal to \(1\) or \(2\) is a direct sum of copies of  \(\mathrm{A}, \mathrm{D},\) 
and  \(\mathrm{E}\) lattices and the integers, thought of as a \(1\)-dimensional 
lattice.

This makes ADE classifications pop up in various places in math and
physics. For example, there is a cool relationship between the ADE
diagrams and the symmetry groups of the Platonic solids, called the
McKay correspondence. Briefly, this is what you do to get it. First,
learn about \(\mathrm{SO}(3)\) and \(\mathrm{SU}(2)\) from
\protect\hyperlink{week61}{``Week 61''} or somewhere. Then, take the
symmetry group of a Platonic solid, or more generally any finite
subgroup \(G\) of \(\mathrm{SO}(3)\). Since \(\mathrm{SO}(3)\) has
\(\mathrm{SU}(2)\) as a double cover, you can get a double cover of
\(G\), say \(\widetilde{G}\), sitting inside \(\mathrm{SU}(2)\). For
example, if \(G\) was the symmetry group of the icosahedron,
\(\widetilde{G}\) would be the icosians (see
\protect\hyperlink{week24}{``Week 24''}).

Since \(\widetilde{G}\) is finite, it has finitely many irreducible
representations. Draw a dot for each of the irreducible representations.
One of these will be \(2\)-dimensional, coming from the spin-\(1/2\)
representation of \(\mathrm{SU}(2)\). Now, when you tensor this 2d rep
with any other irreducible rep \(R\), you get a direct sum of
irreducible reps; draw one line from the dot for \(R\) to each other dot
for each time that other irreducible rep appears in this direct sum.
What do you get? Well, you get an ``affine Dynkin diagram'' of type 
\(\mathrm{A}, \mathrm{D},\) or \(\mathrm{E}\), 
which is like the usual Dynkin diagram but with an extra dot
thrown in (corresponding to the trivial rep of \(\widetilde{G}\)). And,
you get all of them this way!

In fact, playing around with this stuff some more, you can get the
affine Dynkin diagrams of the other simple Lie algebras. There is a lot
more to this... you should probably look at:

\begin{enumerate}
\def\labelenumi{\arabic{enumi})}
\setcounter{enumi}{2}
\item
  John McKay, ``Graphs, singularities and finite groups'', in
  \emph{Proc.\ Symp.\ Pure Math.} vol \textbf{37}, AMS, Providence, 
  Rhode Island,  1980, pp.\ 183--186.

  D.\ Ford and John McKay, ``Representations and Coxeter graphs'', in \emph{The
  Geometric Vein} Coxeter Festschrift (1982), Springer, Berlin,
  pp.\ 549--554.

  John McKay, ``A rapid introduction to ADE theory'', available at
  \href{http://math.ucr.edu/home/baez/ADE.html}{\texttt{http://math.ucr.edu/home/baez/ADE.html}}.
\item
  Pavel Etinghof and Mikhail Khovanov, ``Representations of tensor
  categories and Dynkin diagrams'', available as
  \href{https://arxiv.org/abs/hep-th/9408078}{\texttt{hep-th/9408078}}.
\end{enumerate}

McKay does lots of other mindblowing group theory. He's clearly in tune
with the symmetries of the universe... and occasionally he deigns
to post to the net! A beautiful way of thinking about the McKay
correspondence in terms of category theory appears in the paper by
Etinghof and Khovanov; what we are really doing, it turns out, is
classifying the representations of the tensor category of unitary
representations of \(\mathrm{SU}(2)\). This tensor category is generated
by one object, the spin-\(1/2\) representation, meaning that every other
representation sits inside some tensor power of that one. This way of
thinking of it is important in

\begin{enumerate}
\def\labelenumi{\arabic{enumi})}
\setcounter{enumi}{4}
\tightlist
\item
  J\"urg Fr\"ohlich and Thomas Kerler, \emph{Quantum Groups, Quantum
  Categories, and Quantum Field Theory}, Lecture Notes in
  Mathematics \textbf{1542}, Springer, Berlin, 1991.
\end{enumerate}
Here Fr\"ohlich and Kerler give a classification of certain ``quantum
categories'' that show up in conformal field theory and topological
quantum field theory. These are certain braided tensor categories with
properties like those of the categories of representations of quantum
groups at roots of unity. In such categories, every object has a
``quantum dimension'', which need not be integral, and Fr\"ohlich and
Kerler concentrate on those categories which are generated by a single
object of quantum dimension less than \(2\), getting an ADE-type
classification of them. The category of representations of
\(\mathrm{SU}(2)\), on the other hand, is generated by a single object
of dimension equal to \(2\) --- the spin-\(1/2\) representation --- so
Fr\"ohlich and Kerler are basically generalizing the McKay stuff to
certain quantum groups related to \(\mathrm{SU}(2)\).

Where else do ADE diagrams show up? Well, here I won't try to say
anything about their role in the representation theory of ``quivers'',
or in singularity theory; these are covered pretty well in here:

\begin{enumerate}
\def\labelenumi{\arabic{enumi})}
\setcounter{enumi}{5}
\tightlist
\item
  M. Hazewinkel, W. Hesselink, D. Siermsa, and F. D. Veldkamp, ``The
  ubiquity of Coxeter--Dynkin diagrams (an introduction to the ADE
  problem)'', \emph{Niew. Arch. Wisk.}, \textbf{25} (1977), 257--307.
   Also available at \href{https://ir.cwi.nl/pub/10039/10039D.pdf}{\texttt{https://ir.cwi.nl/pub/10039/10039D.pdf}}.
\end{enumerate}
\noindent
Instead, I'll mention something more recent. In string theory, there is
a Lie algebra called the Virasoro algebra that plays a crucial role; its
almost just the Lie algebra of the group of diffeomorphisms of the
circle, but it's actually just one dimension bigger, being a ``central
extension'' thereof; projective representations of the Lie algebra of
the group of diffeomorphisms of the circle correspond to honest
representations of the Virasoro algebra. An important task in string
theory was to classify the unitary representations of the Virasoro
algebra having a given ``central charge'' \(c\) (this describes the
action of that one extra dimension) and ``conformal weight'' \(h\) (this
describes the action of dilations). It turns out that to get unitary
reps one needs \(c\) and \(h\) to be nonnegative. The representations
with \(c\) between \(0\) and \(1\) are especially nice, for reasons I
don't really understand, and they are called ``minimal models''. An ADE
classification of these was conjectured by Capelli and Zuber, and proved
by

\begin{enumerate}
\def\labelenumi{\arabic{enumi})}
\setcounter{enumi}{6}
\item
  Capelli and Zuber, \emph{Commun. Math. Phys.} \textbf{113} (1987) 1.
\item
  Kato, \emph{Mod. Phys. Lett. A} \textbf{2} (1987) 585.
\end{enumerate}
\noindent
Friedan, Qiu, and Shenker also played a big role in this, in part by
figuring out the allowed values of \(c\). For a good introduction to
this stuff and what it has to do with honest \emph{physics} (which I
admit I've been slacking off on here), try:

\begin{enumerate}
\def\labelenumi{\arabic{enumi})}
\setcounter{enumi}{8}
\tightlist
\item
  Claude Itzykson and Jean-Michel Drouffe, \emph{Statistical Field
  Theory, 1: From Brownian Motion to Renormalization and Lattice Gauge
  Theory}, and \emph{2: Strong Coupling, Monte Carlo Methods, Conformal
  Field Theory and Random Systems.} Cambridge U. Press, 1989.
\end{enumerate}

I will probably come back to this ADE stuff as time goes by, since I'm
sort of fascinated by it, and hopefully folks can refer back to the last
few weeks when I do, so they'll remember what I'm talking about. But in
the next few Weeks I want to catch up with some new developments in math
and physics that have happened in the last few months....

\begin{center}\rule{0.5\linewidth}{0.5pt}\end{center}

\begin{quote}
\emph{A mathematician, like a painter or poet, is a maker of patterns.
If his patterns are more permanent than theirs, it is because they are
made with ideas.} --- G.\ H.\ Hardy
\end{quote}

\hypertarget{week66}{%
\section{October 10, 1995}\label{week66}}

Well, I want to get back to talking about some honest physics, but I
think this week I won't get around to it, since I can't resist
mentioning two tidbits of a more mathematical sort. The first one is
about \(\pi\), and the second one is about the Monster. The second one
\emph{does} have a lot to do with string theory, if only indirectly.

First, thanks to my friend Steven Finch, I just found out that Simon
Plouffe, Peter Borwein and David Bailey have computed the ten billionth
digit in the hexadecimal (i.e., base 16) expansion of \(\pi\). They use
a wonderful formula which lets one compute a given digit of \(\pi\) in
base 16 without needing to compute all the preceding digits! Namely,
 \[  \pi = \sum_{n = 0}^\infty
  \left[
    \frac{4}{8n+1} -\frac{2}{8n+4} -\frac{1}{8n+5} -\frac{1}{8n+6}
  \right] \frac{1}{16^n}
\] 
Since the quantity in square brackets is not an integer, it requires
cleverness to use this formula to get a given digit of \(\pi\), but they
figured out a way. Moreover, their method generalizes to a variety of
other constants.  You can see more information here:

\begin{enumerate}
\def\labelenumi{\arabic{enumi})}

\item
  David Bailey, Peter Borwein and Simon Plouffe, ``On the rapid
  computation of various polylogarithmic constants'', available at 
  \href{http://wayback.cecm.sfu.ca/~pborwein/Z_OLD/PISTUFF/Apistuff.html}{
  \texttt{http://wayback.cecm.sfu.ca/$\sim$pborwein/Z$\underline{\;\;}$OLD/PISTUFF/Apistuff.html}}.
\item
   Steven
  Finch, ``The miraculous Bailey-Borwein-Plouffe \(\pi\) algorithm'',
  \href{http://plouffe.fr/simon/articles/Miraculous.pdf}{\texttt{http://plouffe.fr/simon/articles/Miraculous.pdf}}.
\end{enumerate}

\noindent
The first one has details about how
the billionth digit of \(\pi\), \(\pi^2\), \(\ln(2)\), and some other
constant were first computed.  The second gives a good
overview of what's up.

Can we hope for a similar formula in base 10? More importantly, could
these ideas let us prove that \(\pi\) is ``normal'', that is, that every
possible string of digits appears in it with the frequency one would
expect of a ``random'' number? It seems that this would be a natural
avenue of attack.

Next, a tidbit of a more erudite sort concerning the elusive Monster
manifold. Recall from \protect\hyperlink{week63}{``Week 63''} and
\protect\hyperlink{week64}{``Week 64''} that the compact simple Lie
groups can classified into 4 infinite families and 5 exceptions. I have
always been puzzled by these ``exceptional Lie groups'', so I tried to
explain a bit about how they are related to some other ``exceptional
structures'' in mathematics, such as the icosahedron and the octonions.
In physics, Witten has suggested that the correct theory of our universe
might also be an exceptional structure of some sort. This idea has found
some support in string theory, which uses the exceptional Lie group
\(\mathrm{E}_8\) and other structures I'll mention a bit later. In a
more hand-waving way, one may argue that the theory of our universe must
be incredibly special, since out of all the theories we can write down,
just this \emph{one} describes the universe that actually \emph{exists}.
All sorts of simpler universes apparently don't exist. So maybe the
theory of the universe needs to use special, ``exceptional'' mathematics
for some reason, even though it's complicated.

Anyway, as a hard-nosed mathematician, vague musings along these lines
get tiresome to me rather quickly. Instead, what interests me most about
this stuff is the whole idea of ``exceptional structures'' ---
symmetrical structures that don't fit into the neat regular families in
classification theorems. The remarkable fact is that many of them are
deeply related. As Geoffrey Dixon put it to me, they seem to have a
``holographic'' quality, meaning that each one contains in encoded form
some of the information needed to construct all the rest! It thus seems
pointless to hope that one is ``the explanation'' for the rest, but I
would still like some conceptual ``explanation'' for the whole
collection of them --- though I have no idea what it should be.

Surely a clue must lie in the theory of finite simple groups. Just as
the simple Lie groups are the building blocks of the theory of
continuous symmetries, these are the building blocks of the theory of
discrete --- indeed finite --- symmetries. More precisely ``finite
simple'' group is a group \(G\) with finitely many elements and no
normal subgroups, that is, no nontrivial subgroups \(H\) such that \(h\)
in \(H\) implies \(ghg^{-1}\) in \(H\) for all \(g\) in \(G\). This
condition means that you cannot form the ``quotient group'' \(G/H\),
which one can think of as a ``more simplified'' version of \(G\).

The classification of the finite simple groups is one of remarkable
achievements of 20th-century mathematics. The entire proof of the
classification theorem is estimated to take 10,000 pages, done by many
mathematicians. To start learning about it, try:

\begin{enumerate}
\def\labelenumi{\arabic{enumi})}
\setcounter{enumi}{3}
\tightlist
\item
  Ron Solomon, ``On finite simple groups and their classification'',
  \emph{Notices Amer. Math. Soc.} \textbf{45}, February 1995, 231--239.
\end{enumerate}

\noindent
and the references therein. Again, there are some infinite families and
26 exceptions called the ``sporadic'' groups. The biggest of these is
the Monster, with \[
  \begin{gathered}
    246\cdot 320\cdot 59\cdot 76\cdot 112\cdot 133\cdot 17\cdot 19\cdot 23\cdot 29\cdot 31\cdot 41\cdot 47\cdot 59\cdot 71
    \\= 808017424794512875886459904961710757005754368000000000
  \end{gathered}
\] elements. It is a kind of Mt. Everest of the sporadic groups, and all
the routes to it I know involve a tough climb through all sorts of
exceptional structures: \(\mathrm{E}_8\) (see
\protect\hyperlink{week65}{``Week 65''}), the Leech lattice (see
\protect\hyperlink{week20}{``Week 20''}), the Golay code, the Parker
loop, the Griess algebra, and more. I certainly don't understand this
stuff....

Even before the Monster was proved to exist, it started casting its
enormous shadow over mathematics. For example, consider the theory of
modular functions. What are those? Well, consider a lattice in the
complex plane, like this: \[
  \begin{tikzpicture}[scale=0.7]
    \draw[->] (-3,0) to (4,0) node[label=below:{$\mathrm{Re}(z)$}]{};
    \draw[->] (0,-3) to (0,4) node[label=left:{$\mathrm{Im}(z)$}]{};
    \foreach \m in {-1,0,1,2}
    {
      \foreach \n in {-1,0,1,2}
      {
        \node at ({\m*1.5-\n/3},{1.5*\n+\m}) {$\bullet$};
      }
    }
  \end{tikzpicture}
\] These are important in complex analysis, as described in
\protect\hyperlink{week13}{``Week 13''}. To describe one of these you
can specify two ``periods'' \(\omega_1\) and \(\omega_2\), complex
numbers such that every point in the lattice of the form
\[n \omega_1 + m \omega_2.\] But this description is redundant, because if
we choose instead to use \[
  \begin{aligned}
    \omega'_1 &= a\omega_1+b\omega_2
  \\\omega'_2 &= c\omega_1+b\omega_2
  \end{aligned}
\] we'll get the same lattice as long as the matrix of integers \[
  \left(
    \begin{array}{cc}
      a&b\\c&d
    \end{array}
  \right)
\] is invertible and its inverse also consists of integers. These
transformations preserve the ``handedness'' of the basis \(\omega_1\),
\(\omega_2\) if they have determinant \(1\), and that's generally a good
thing to require. The group of \(2\times2\) invertible matrices over the
integers with determinant \(1\) is called \(\mathrm{SL}(2,\mathbb{Z})\),
or the ``modular group'' in this context. I said a bit about it and its
role in string theory in \protect\hyperlink{week64}{``Week 64''}.

Now, if we are only interested in parametrizing the different
\emph{shapes} of lattices, where two rotated or dilated versions of the
same lattice count as having the same shape, it suffices to use one
complex number, the ratio \[\tau=\frac{\omega_1}{\omega_2}.\] We might
as well assume \(\tau\) is in the upper halfplane, \(H\), while we're at
it. But for the reason given above, this description is redundant; if we
have a lattice described by \(\tau\), and a matrix in
\(\mathrm{SL}(2,\mathbb{Z})\), we get a new guy \(\tau'\) which really
describes the same shaped lattice. If you work it out,
\[\tau' = \frac{a\tau + b}{c\tau + d}.\] So the space of different
possible shapes of lattices in the complex plane is really the quotient
\[H/\mathrm{SL}(2,\mathbb{Z}).\] Now, a function on this space is just a
function of \(\tau\) that doesn't change when you replace \(\tau\) by
\(\tau'\) as above. In other words, it's basically just a function
depending only on the shape of a 2d lattice. Now it turns out that there
is essentially just one ``nice'' function of this sort, called \(j\);
all other ``nice'' functions of this sort are functions of \(j\). (For
experts, what I mean is that the meromorphic
\(\mathrm{SL}(2,\mathbb{Z})\)-invariant functions on \(H\) union the
point at infinity are all rational functions of this function \(j\).)

It looks like this:
\[j(\tau) = q^{-1} + 744 + 196884 q + 21493706 q^2 + \cdots\] where
\(q = \exp(2\pi i\tau)\). In fact, starting from a simple situation, we
have quickly gotten into quite deep waters. The simplest explicit
formula I know for \(j\) involves lattices in \(24\)-dimensional space!
This could easily be due to my limited knowledge of this stuff, but it
suits my present purpose: first, we get a vague glimpse of where
\(\mathrm{E}_8\) and the Leech lattice come in, and second, we get a
vague glimpse of the mysterious significance of the numbers 24 and 26 in
string theory.

So what is this \(j\) function, anyway? Well, it turns out we can define
it as follows. First form the Dedekind eta function
\[\eta(q) = q^{\frac{1}{24}}\prod_{n=1}^\infty(1-q^n).\] This is not
invariant under the modular group, but it transforms in a pretty simple
way. Then take the \(\mathrm{E}_8\) lattice --- remember, that's a very
nice lattice in 8 dimensions, in fact the only ``even unimodular''
lattice in 8 dimensions, meaning that the inner product of any two
vectors in the lattice is even, and the volume of each fundamental
domain in it equals \(1\). Now take the direct sum of 3 copies of
\(\mathrm{E}_8\) to get an even unimodular lattice \(L\) in 24
dimensions. Then form the theta function
\[\theta(q) = \sum_{x\in L}q^{\langle x,x\rangle/2}.\] In other words,
we take all lattice points \(x\) and sum \(q\) to the power of their
norm squared over \(2\). Now we have
\[j(\tau) = \frac{\theta(q)}{\eta(q)^{24}}.\]

Quite a witches' brew of a formula, no? If someone could explain to me
the deep inner reason for \emph{why} this works, I'd be delighted, but
right now I am clueless. I will say this, though: we could replace \(L\)
with any other even unimodular lattice in 24 dimensions and get a
function differing from \(j\) only by a constant. Guess how many even
unimodular lattices there are in 24 dimensions? Why, 24, of course!
These ``Niemeier lattices'' were classified by Niemeier in 1968. All but
one of them have vectors with length squared equal to \(2\), but there
is one whose shortest vector has length squared equal to \(4\), and
that's the Leech lattice. This one has a very charming relation to
\(26\)-dimensional spacetime, described in
\protect\hyperlink{week20}{``Week 20''}.

Since the constant term in \(j\) can be changed by picking different
lattices in 24 dimensions, and constant functions aren't very
interesting anyway, we can say that the first interesting coefficient in
the above power series for \(j\) is 196884. Then, right around when the
Monster was being dreamt up, McKay noticed that the dimension of its
smallest nontrivial representation, namely 196883, was suspiciously
similar. Coincidence? No.~It turns out that all the coefficients of
\(j\) can be computed from the dimensions of the irreducible
representations of the Monster! Similarly, Ogg noticed in the study of
the modular group, the primes 2, 3, 5, 7, 11, 13, 17, 19, 23, 29, 31,
41, 47, 59 and 71 play a special role. He went to a talk on the Monster
and noticed the ``coincidence''. Then he wrote a paper offering a bottle
of Jack Daniels to anyone who could explain it. This was the beginning
of a subject called ``Monstrous Moonshine''... the mysterious
relation between the Monster and the modular group.

Well, as it eventually turned out, one way to get ahold of the Monster
is as a group of symmetries of a certain algebra of observables for a
string theory, or more precisely, a ``vertex operator algebra'':

\begin{enumerate}
\def\labelenumi{\arabic{enumi})}
\setcounter{enumi}{4}
\tightlist
\item
  Igor Frenkel, James Lepowsky, and Arne Meurman, \emph{Vertex Operator
  Algebras and the Monster}, Academic Press, Boston, 1988.
\end{enumerate}
\noindent
The relation of string theory to modular invariance and 26 dimensional
spacetime then ``explains'' some of the mysterious stuff mentioned
above. (By the way, the authors of the above book say the fact that
there are 26 sporadic finite simple groups is just a coincidence. I'm
not so sure myself\ldots{} not that I understand any of this stuff, but
it's just too spooky how the number 26 keeps coming up all over!)

Anway, now let me fast-forward to some recent news. I vaguely gather
that people would like to explain the relation between the Monster and
string theory more deeply, by finding a \(24\)-dimensional manifold
having the Monster as symmetries, and cooking up a field theory of maps
from the string worldsheet to this ``Monster manifold'', so that the
associated vertex operator algebra would have a good reason for having
the Monster as symmetries. Apparently Hirzebruch has offered a prize for
anyone who could do this in a nice way, by finding a ``24-manifold with
\(p_1=0\) whose Witten genus is \((j-744)\Delta\)'' on which the Monster
acts. Recently, Mike Hopkins at MIT and Mark Mahowald at Northwestern
have succeeded in doing the first part, the part in quotes above. They
haven't gotten a Monster action yet. Their construction uses a lot of
homotopy theory.

I don't have much of a clue about any of this stuff, but Allen Knutson
suggests that I read

\begin{enumerate}
\def\labelenumi{\arabic{enumi})}
\setcounter{enumi}{5}
\tightlist
\item
  Friedrich Hirzebruch, Thomas Berger, and Rainer Jung, \emph{Manifolds
  and Modular Forms}, translated by Peter S.\ Landweber,
  Braunschweig, Vieweg, 1992.
\end{enumerate}
\noindent
for more about this ``Witten genus'' stuff. He also has referred me to
the following articles by Borcherds:

\begin{enumerate}
\def\labelenumi{\arabic{enumi})}
\setcounter{enumi}{6}
\item
  Richard E. Borcherds, ``The Monster Lie-algebra'', \emph{Adv. Math.}
  \textbf{83} (1990), 30--47.

  Richard E. Borcherds, ``Monstrous Moonshine and monstrous
  Lie-superalgebras'', \emph{Invent. Math.} \textbf{109} (1992),
  405--444.
\end{enumerate}

\noindent
For your entertainment and edification, I include the abstract of the
second one below:

\begin{quote}
We prove Conway and Norton's moonshine conjectures for the infinite
dimensional representation of the monster simple group constructed by
Frenkel, Lepowsky and Meurman. To do this we use the no-ghost theorem
from string theory to construct a family of generalized Kac--Moody
superalgebras of rank 2, which are closely related to the monster and
several of the other sporadic simple groups. The denominator formulas of
these superalgebras imply relations between the Thompson functions of
elements of the monster (i.e.~the traces of elements of the monster on
Frenkel, Lepowsky, and Meurman's representation), which are the
replication formulas conjectured by Conway and Norton. These replication
formulas are strong enough to verify that the Thompson functions have
most of the ``moonshine'' properties conjectured by Conway and Norton,
and in particular they are modular functions of genus 0. We also
construct a second family of Kac--Moody superalgebras related to elements
of Conway's sporadic simple group \(\mathrm{Co}_1\). These superalgebras have even
rank between 2 and 26; for example two of the Lie algebras we get have
ranks 26 and 18, and one of the superalgebras has rank 10. The
denominator formulas of these algebras give some new infinite product
identities, in the same way that the denominator formulas of the affine
Kac--Moody algebras give the Macdonald identities.
\end{quote}

\hypertarget{week67}{%
\section{October 23, 1995}\label{week67}}

I'm pretty darn busy now, so the forthcoming Weeks will probably be
pretty hastily written. This time I'll mainly write about quantum
gravity.

\begin{enumerate}
\def\labelenumi{\arabic{enumi})}
\tightlist
\item
  Margaret Wertheim, \emph{Pythagoras' Trousers: God, Physics, and the
  Gender Wars}, Times Books/Random House, New York, 1995.
\end{enumerate}

\noindent
I enjoyed this book, despite or perhaps because of the fact that it may
annoy lots of physicists. It notes that, starting with Pythagoras,
theoretical physics has always had a crypto-religious aspect. With
Pythagoras it was obvious; he seems to have been the leader of a special
sort of religious cult. With people like Kepler, Newton and Einstein it
is only slightly less obvious. The operative mythology in every case is
that of the mage. Think of Einstein, stereotypically with long white
hair (though most of best work was actually done before his hair got
white), a powerful but benign figure devoted to finding out ``the
thoughts of God''. The mage, of course, is typically male, and Wertheim
argues that this makes it harder for women to become physicists. (A lot
of the same comments would apply to mathematics.) It is not a very
scholarly book, but I wouldn't dismiss it.

\begin{enumerate}
\def\labelenumi{\arabic{enumi})}
\setcounter{enumi}{1}
\tightlist
\item
  Stephen W.\ Hawking, ``Virtual black holes'', \emph{Phys.\ Rev.\ D} 
  \textbf{53} (1996), 3099--3107.   Also available as
  \href{https://arxiv.org/abs/hep-th/9510029}{\texttt{hep-th/9510029}}.
\end{enumerate}
\noindent
Hawking likes the ``Euclidean path-integral approach'' to quantum
gravity. The word ``Euclidean'' is a horrible misnomer here, but it
seems to have stuck. It should really read ``Riemannian'', the idea
being to replace the Lorentzian metric on spacetime by one in which time
is on the same footing as space. One thus attempts to compute answers to
quantum gravity problems by integrating over all Riemannian metrics on
some 4-manifold, possibly with some boundary conditions. Of course, this
is tough --- impossible so far --- to make rigorous. But Hawking isn't
scared; he also wants to sum over all 4-manifolds (possibly having a
fixed boundary). Of course, to do this one needs to have some idea of
what ``all 4-manifolds'' are. Lots of people like to consider wormholes,
which means considering 4-manifolds that aren't simply connected. Here,
however, Hawking argues against wormholes, and concentrates on
simply-connected 4-manifolds. He writes: 
\begin{quote}
Barring some pure
mathematical details, it seems that the topology of simply connected
four-manifolds can be essentially represented by gluing together three
elementary units, which I call bubbles. The three elementary units are
\(S^2 \times S^2\), \(\mathbb{C}\mathrm{P}^2\), and \(\mathrm{K}3\). 
The latter two have
orientation reversed versions, \(-\mathbb{C}\mathrm{P}^2\) and \(-\mathrm{K}3\).
\end{quote}
\(S^2 \times S^2\) is just the product of the 2-dimensional sphere with
itself, and he argues that this sort of bubble corresponds to a virtual
black hole pair. He considers the effect on the Euclidean path integral
when you have lots of these around (i.e., when you take the connected
sum of \(S^4\) with lots of these). He argues that particles scattering
off these lose quantum coherence, i.e., pure states turn to mixed
states. And he argues that this effect is very small at low energies
\emph{except} for scalar fields, leading him to predict that we may
never observe the Higgs particle! Yes, a real honest particle physics
prediction from quantum gravity! As he notes, ``unless quantum gravity
can make contact with observation, it will become as academic as
arguments about how many angels can dance on the head of a pin''. I
suspect he also realizes that he'll never get a Nobel prize unless he
goes out on a limb like this. He also gives an argument for why
the ``\(\theta\) angle" measuring CP violation by the strong force may be
zero. This parameter sits in front of a term in the Standard Model
Lagrangian; there seems to be no good reason for it to be zero, but
measurements of the neutron electric dipole moment show that it has to
be less than \(10^{-9}\), according to the following book:

\begin{enumerate}
\def\labelenumi{\arabic{enumi})}
\setcounter{enumi}{2}
\tightlist
\item
  Kerson Huang, \emph{Quarks, Leptons, and Gauge Fields}, World
  Scientific, Singapore, 1982.
\end{enumerate}

\noindent
Perhaps there are better bounds now, but this book is certainly one of
the nicest introductions to the Standard Model, and if you want to learn
about this ``\(\theta\) angle'' stuff, it's a good place to start.

Anyway, rather than going further into the physics, let me say a bit
about the ``pure mathematical details''. Here I got some help from Greg
Kuperberg, Misha Verbitsky, and Zhenghan Wang (via Xiao-Song Lin, a
topologist who is now here at Riverside). Needless to say, the mistakes
are mine alone, and corrections and comments are welcome!

First of all, Hawking must be talking about homeomorphism classes of
compact oriented simply-connected smooth 4-manifolds, rather than
diffeomorphism classes, because if we take the connected sum of 9 copies
of \(\mathbb{C}\mathrm{P}^2\) and one of \(-\mathbb{C}\mathrm{P}^2\), that has infinitely
many different smooth structures. Why the physics depends only on the
homeomorphism class is beyond me... maybe he is being rather
optimistic. But let's follow suit and talk about homeomorphism classes.
Folks consider two cases, depending on whether the intersection form on
the second cohomology is even or odd. If our 4-manifold has an odd
intersection form, Donaldson showed that it is an connected sum of
copies of \(\mathbb{C}\mathrm{P}^2\) and \(-\mathbb{C}\mathrm{P}^2\). If its intersection
form is even, we don't know yet, but if the ``11/8 conjecture'' is true,
it must be a connected sum of \(\mathrm{K}3\)'s and \(S^2 \times S^2\)'s. Here I
cannot resist adding that \(\mathrm{K}3\) is a 4-manifold whose intersection form is
\(\mathrm{E}_8 \oplus \mathrm{E}_8 \oplus H \oplus H \oplus H\), where
\(H\) is the \(2\times2\) matrix \[
  \left(
    \begin{array}{cc}
      0&1\\1 &0
    \end{array}
  \right)
\] and \(\mathrm{E}_8\) is the nondegenerate symmetric \(8\times8\)
matrix describing the inner products of vectors in the wonderful
lattice, also called \(\mathrm{E}_8\), which I discussed in
\protect\hyperlink{week65}{``Week 65''}. So \(\mathrm{E}_8\) raises its
ugly head yet again.... By the way, \(H\) is just the intersection
form of \(S^2 \times S^2\), while the intersection form of
\(\mathbb{C}\mathrm{P}^2\) is just the \(1\times1\) matrix \((1)\).

Even if the 11/8 conjecture is not true, we could if necessary resort to
Wall's theorem, which implies that any 4-manifold becomes homeomorphic
--- even diffeomorphic --- to a connected sum of the 5 basic types of
``bubbles'' after one takes its connected sum with sufficiently many
copies of \(S^2 \times S^2\). This suggests that if Euclidean path
integral is dominated by the case where there are lots of virtual black
holes around, Hawking's arguments could be correct at the level of
diffeomorphism, rather than merely homeomorphism. Indeed, he says that
``in the wormhole picture, one considered metrics that were multiply
connected by wormholes. Thus one concentrated on metrics {[}I'd say
topologies!{]} with large values of the first Betti number{[}....{]}
However, in the quantum bubbles picture, one concentrates on spaces with
large values of the second Betti number.''

\begin{enumerate}
\def\labelenumi{\arabic{enumi})}
\setcounter{enumi}{3}
\tightlist
\item
  Ted Jacobson, ``Thermodynamics of spacetime: the Einstein equation of
  state'', \emph{Phys.\ Rev.\ Lett.} \textbf{75} (1995), 1260--1263.  Also available as
  \href{https://arxiv.org/abs/gr-qc/9504004}{\texttt{gr-qc/9504004}}.
\end{enumerate}

\noindent
Well, here's another paper on quantum gravity, also very good, which
seems at first to directly contradict Hawking's paper. Actually,
however, I think it's another piece in the puzzle. The idea here is to
derive Einstein's equation from thermodynamics! More precisely, ``The
key idea is to demand that this relation hold for all the local Rindler
causal horizons through each spacetime point, with {[}the change in
heat{]} and {[}the temperature{]} interpreted as the energy flux and
Unruh temperature seen by an accelerated observer just inside the
horizon. This requires that gravitational lensing by matter energy
distorts the causal structure of spacetime in just such a way that the
Einstein equation holds''. It's a very clever mix of classical and
quantum (or semiclassical) arguments. It suggests that all sorts of
quantum theories on the microscale could wind up yielding Einstein's
equation on the macroscale.

\begin{enumerate}
\def\labelenumi{\arabic{enumi})}
\setcounter{enumi}{4}
\tightlist
\item
  Lee Smolin, ``The Bekenstein bound, topological quantum field theory
  and pluralistic quantum field theory'', available as
  \href{https://arxiv.org/abs/gr-qc/9508064}{\texttt{gr-qc/9508064}}.
\end{enumerate}

\noindent
This is a continued exploration of the themes of Smolin's earlier paper,
reviewed in \protect\hyperlink{week56}{``Week 56''} and
\protect\hyperlink{week57}{``Week 57''}. Particularly interesting is the
general notion of ``pluralistic quantum field theory'', in which
different observers have different Hilbert spaces. This falls out
naturally in the \(n\)-categorical approach pursued by Crane (see
\protect\hyperlink{week56}{``Week 56''}), which I am also busily
studying.

\begin{enumerate}
\def\labelenumi{\arabic{enumi})}
\setcounter{enumi}{5}
\tightlist
\item
  Daniel Armand-Ugon, Rodolfo Gambini, Octavio Obregon and Jorge Pullin, ``Towards a loop
  representation for quantum canonical supergravity'', \emph{Nucl.\ Phys.\ B}
  \textbf{460} (1996), 615--631.  Also available as
  \href{https://arxiv.org/abs/hep-th/9508036}{\texttt{hep-th/9508036}}.
\end{enumerate}

\noindent
Some knot theorists and quantum group theorists had better take a look
at this! This paper considers the analog of \(\mathrm{SU}(2)\)
Chern--Simons theory where you use the supergroup \(\mathrm{GSU}(2)\),
and perturbatively work out the skein relations of the associated link
invariant (up to a certain low order in \(\hbar\)). If someone
understood the quantum supergroup ``quantum \(\mathrm{GSU}(2)\)'', they
could do this stuff nonperturbatively, and maybe get an interesting
invariant of links and 3-manifolds, and make some physicists happy in
the process.

\begin{enumerate}
\def\labelenumi{\arabic{enumi})}
\setcounter{enumi}{6}
\tightlist
\item
  Roh Suan Tung and Ted Jacobson, ``Spinor one-forms as gravitational
  potentials'', \emph{Class.\ Quant.\ Grav.} \textbf{12} (1995), L51--L55.  Also available as
  \href{https://arxiv.org/abs/gr-qc/9502037}{\texttt{gr-qc/9502037}}.
\end{enumerate}

\noindent
This paper writes out a new Lagrangian for general relativity, closely
related to the action that gives general relativity in the Ashtekar
variables. It's incredibly simple and beautiful! I am hoping that if I
work on it someday, it will explain to me the mysterious relation
between quantum gravity and spinor fields (see the end of
\protect\hyperlink{week60}{``Week 60''}).

\begin{enumerate}
\def\labelenumi{\arabic{enumi})}
\setcounter{enumi}{7}
\tightlist
\item
  Joseph Polchinski and Edward Witten, ``Evidence for heterotic--type
  I string duality'', \emph{Nucl.\ Phys.\ B} \textbf{460} (1996), 525--540.  Also available as
  \href{https://arxiv.org/abs/hep-th/9510169}{\texttt{hep-th/9510169}}.
\end{enumerate}

\noindent
I'm no string theorist, so I've been lagging vastly behind the new work
on ``dualities'' that has revived interest in the subject. Roughly
speaking, though, it seems folks have discovered a host of secret
symmetries relating superficially different string theories...
making them, in some deeper sense, the same. The heterotic and type I
strings are two of the most famous string theories, so if they were
really equivalent as this paper suggests, it would be very interesting.

\hypertarget{week68}{%
\section{October 29, 1995}\label{week68}}

Okay, now the time has come to speak of many things: of topoi,
glueballs, communication between branches in the many-worlds
interpretation of quantum theory, knots, and quantum gravity.

\begin{enumerate}
\def\labelenumi{\arabic{enumi})}
\tightlist
\item
  Robert Goldblatt, \emph{Topoi: the Categorial Analysis of Logic},
  Dover, Mineola, New York, 2006.  Also available as 
  \href{https://projecteuclid.org/ebooks/books-by-independent-authors/Topoi-The-Categorial-Analysis-of-Logic/toc/bia/1403013939}{\texttt{https://projecteuclid.org/ebooks/books-by-}}
   \href{https://projecteuclid.org/ebooks/books-by-independent-authors/Topoi-The-Categorial-Analysis-of-Logic/toc/bia/1403013939}{\texttt{independent-authors/Topoi-The-Categorial-Analysis-of-Logic/toc/bia/}} \break \href{https://projecteuclid.org/ebooks/books-by-independent-authors/Topoi-The-Categorial-Analysis-of-Logic/toc/bia/1403013939}{\texttt{1403013939}}.
\end{enumerate}

\noindent
If you've ever been interested in logic, you've got to read this book.
Unless you learn a bit about topoi, you are really missing lots of the
fun. The basic idea is simple and profound: abstract the basic concepts
of set theory, so as to define the notion of a ``topos'', a kind of
universe like the world of classical logic and set theory, but far more
general!

For example, there are ``intuitionistic'' topoi in which Brouwer reigns
supreme --- that is, you can't do proof by contradiction, you can't use
the axiom of choice, etc.  There is also the ``effective topos'' of
Hyland in which Turing reigns supreme --- for example, the only
functions are the effectively computable ones. There is also a
``finitary'' topos in which all sets are finite. So there are topoi to
satisfy various sorts of ascetic mathematicians who want a
stripped-down, minimal form of mathematics.

However, there are also topoi for the folks who want a mathematical
universe with lots of horsepower and all the options! There are topoi in
which everything is a function of time: the membership of sets, the
truth-values of propositions, and so on all depend on time. There are
topoi in which everything has a particular group of symmetries. Then
there are \emph{really} high-powered things like topoi of sheaves on a
category equipped with a Grothendieck topology....

And so on: not an attempt to pick out ``the'' universe of logic and
mathematics, but instead, an effort to systematically examine a bunch of
them and how they relate to each other. The details can be intimidating,
but Goldblatt explains them very nicely. A glance at the subject
headings reveal some of the delights in store: ``elementary truth'',
``local truth'', ``geometric logic'', etc.

What is a topos, precisely? Well, most people would need to limber up a
little bit before getting the precise definition... so let me just
start you off with some mental stretching exercises. In classical logic
we are used to working with two truth values, \(\mathrm{True}\) and
\(\mathrm{False}\). Let's call the set of truth values \(\Omega\), just
to make it sound impressive --- and because it's traditional in topos
theory. So, we are used to doing all our logic with
\[\Omega = \{\mathrm{True}, \mathrm{False}\}.\] In set theory, one of
the things we do with \(\Omega\) is describe subsets of a given set
\(X\). In other words, to describe a subset \(Y\) of \(X\), we can say
for each member of \(X\), whether it is \(\mathrm{True}\) or
\(\mathrm{False}\) that it is a member of \(Y\). Thus we can describe
the subset \(Y\) by giving a function \[f\colon X \to \Omega.\] We say
\(x\) is in \(Y\) if \(f(x) = \mathrm{True}\), but \(x\) is not in \(Y\)
if \(f(x) = \mathrm{False}\).

Now say we wanted to describe a topos of ``time-dependent sets''. But
instead of ``time-dependent sets'', let's act like topos theorists and
call them simply ``objects'', and instead of talking about one being a
``subset'' of another, let's say one is a ``subobject'' of another. To
keep life simple, let's consider only two times: today and tomorrow. So
we can think of an object \(X\) in this topos as a pair \((X_1, X_2)\)
of sets: one set \(X_1\) today, and another set \(X_2\) tomorrow. We say
that an object \(Y\) is a ``subobject'' of \(X\) if \(Y_1\) is a subset
of \(X_1\) and \(Y_2\) is a subset of \(X_2\). The idea is that we want
\(Y\) to be contained in \(X\) both today and tomorrow.

Now, to describe a subobject \(Y\) of \(X\), we can say what's in \(Y\)
today, and also what's in \(Y\) tomorrow. We would like to do so using
some kind of function, or what topos theorists call a ``morphism'',
\[f \colon X \to \Omega.\] Clearly we can't do this with our old truth values
set \(\{\mathrm{True},\mathrm{False}\}\). Instead, we should use a truth
values \emph{object} \(\Omega\) that keeps track of what's true or false
today and what's true or false tomorrow. In other words, \(\Omega\)
should now be the pair of sets
\[(\{\mathrm{True}, \mathrm{False}\}, \{\mathrm{True}, \mathrm{False}\})\]
Now what is that ``morphism'' \(f\) exactly? Well, it's like one
function today and another function tomorrow, or in other words, a pair
of functions! In general, a morphism \(f\colon S \to T\) between objects
in this topos is a pair of functions \((f_1,f_2)\), with
\(f_1\colon S_1 \to T_1\) and \(f_2\colon S_2 \to T_2\). Then in our
particular case, the morphism \(f\colon X \to \Omega\) will say which
elements of \(X_1\) are in \(Y_1\), and which elements of \(X_2\) are in
\(Y_2\).

This is a pretty simple example of what the objects and morphisms in a
topos can be like. They can be a lot weirder. The key thing is that you
can do a lot of the same things with them that you can do with sets and
functions. Also, you can do a lot of the same things with \(\Omega\)
that you can with \(\{\mathrm{True}, \mathrm{False}\}\). Note that in
our example, like in the classical example where
\(\Omega = \{\mathrm{True}, \mathrm{False}\}\), the subobject classifier
has a bunch of logical operations on it: morphisms like \[
  \begin{aligned}
    \mathrm{Not}&\colon \Omega \to \Omega \\
    \mathrm{And}&\colon \Omega \times \Omega \to \Omega \\
    \mathrm{Or}&\colon \Omega \times \Omega \to \Omega
  \end{aligned}
\] and so on. In our example, and in the classical example, these make
\(\Omega\) into what folks call a boolean algebra. Basically, all the
usual rule of logic apply. In general, though, \(\Omega\) only needs to
be a Heyting algebra. This is more general than a boolean algebra, and
it can be sort of intuitionistic in flavor; for example,
\(\mathrm{Not} \circ \mathrm{Not}\) doesn't need to equal the identity
morphism \(1\colon \Omega \to \Omega\), so proof by contradiction
doesn't necessarily work. A typical example of a Heyting algebra
\(\Omega\) is the set of open sets in a topological space, with
\(\mathrm{And}\) and \(\mathrm{Or}\) being intersection and union, and
with \(\mathrm{Not}\) being the \emph{interior} of the complement. This
gives a little hint as to what topoi have to do with topology.

After studying this sort of thing a while, it's rather hard to go on
pretending that the Zermelo--Fraenkel axioms of set theory, which were
cooked up in the early 20th century to escape the logical paradoxes of
Russell and others, are the last word on ``foundations''. One can
develop topos theory within set theory if one wishes, but one can also
set up topos theory from scratch, as a kind of pluralistic foundation of
mathematics.

For a deeper but still friendly and expository introduction to topoi,
try

\begin{enumerate}
\def\labelenumi{\arabic{enumi})}
\setcounter{enumi}{1}
\tightlist
\item
  Saunders Mac Lane and Ieke Moerdijk, \emph{Sheaves in Geometry and
  Logic: A First Introduction to Topos Theory}, Springer, New
  York, 1992.
\end{enumerate}
\noindent
Here you can learn about ``Brouwer's theorem: all functions are
continuous'' (in a suitably intuitionistic topos, of course). You can
also learn topos-theoretic versions of Cohen's proofs of the
independence of the continuum hypothesis and the axiom of choice.

Goldblatt's book teaches you all the category theory you need to learn
about topoi... but for people who already know some category
theory, let me give the precise definition of a topos (or more
precisely, an elementary topos, to distinguish it from a ``Grothendieck
topos''): it's a category with finite limits and power objects. This
automatically implies a lot of things, such as the existence of the
subobject classifier \(\Omega\) that I was talking about.

To get deeper into topos theory, try:

\begin{enumerate}
\def\labelenumi{\arabic{enumi})}
\setcounter{enumi}{2}
\tightlist
\item
  Michael Barr and Charles Wells, \emph{Toposes, Triples and Theories},
  Springer, Berlin, 1983.  Also available as \href{http://www.tac.mta.ca/tac/reprints/articles/12/tr12abs.html}{\texttt{http://www.tac.mta.ca/tac/reprints/articles/12/tr12abs.html}}.
\end{enumerate}

Now let me catch up on some things more directly related to physics:

\begin{enumerate}
\def\labelenumi{\arabic{enumi})}
\setcounter{enumi}{3}
\tightlist
\item
  Frank Close, ``Are glueballs and hybrids found?'', in \emph{Hadron'95: 
  The 6th International Conference on Hadron Spectroscopy, the University 
  of Manchester, Manchester, UK, 10th-14th July 1995}, World Scientific, Singapore,
  1996.  Also available as
  \href{https://arxiv.org/abs/hep-ph/9509245}{\texttt{hep-ph/9509245}}.
\item
J.\ Sexton, A.\ Vaccarino and D.\ Weingarten, ``Numerical evidence for the
observation of a scalar glueball'', \emph{Phys.\ Rev.\ Lett.} \textbf{75} (1995), 
4563-4566.  Also available as
\href{https://arxiv.org/abs/hep-lat/9510022}{\texttt{hep-lat/9510022}}.
\end{enumerate}
\noindent
Thanks go to Greg Kilcup for bringing these to my attention. Have they
found a glueball??? That would be really exciting. What's a glueball,
you ask? Well, quantum chromodynamics, our best theory of the strong
force, says that that the strong force is carried by particles called
``gluons''. Like electromagnetism, the strong force is a gauge field,
but it's a nonabelian gauge field, so the gluons themselves have charge,
or ``color''. Thus they interact in a nonlinear way. This is what lets
them bind together quarks in such a tight way. But perhaps, in addition
to pairs of quarks and antiquarks held together by gluons --- i.e.,
mesons --- and triples of quarks held together by gluons --- i.e.,
baryons --- there could be short-lived assemblages consisting entirely
of gluons, held together by their self-interactions. These are called
glueballs, but we don't know if these exist.

However, to my surprise, it turns out that there are now some candidates
out there! The first paper suggests that the \(f_0(1500)\), a neutral
spin-zero particle with mass around 1500 MeV, is a glueball. The second
paper argues instead that this is basically a quark-antiquark pair (made
of a strange quark and a strange antiquark... where ``strange'' is
the technical name for one of the 6 quarks!). It presents evidence from
a numerical simulation and argues that the ``\(\theta\)'' or
\(f_J(1710)\), a neutral particle with even spin and mass 1710 MeV, is a
glueball. Part of the subtlety here is that, thanks to the superposition
principle, there is not a perfectly sharp distinction between a glueball
with some virtual quark-antiquark pairs in it, and a quark-antiquark
pair with a bunch of virtual gluons in it. There can be ``hybrids'' that
are a bit of both a linear combination of a meson and a glueball! (This
phenomenon of ``hybridization'' is also familiar in chemistry.)

It's tough to do nonperturbative computations in nonlinear gauge field
theories --- basically one needs to approximately compute a path
integral, using Monte Carlo technique, approximating spacetime by a
lattice (in this case, a \(16 \times 16 \times 16 \times 24\) lattice).
Computing the properties of a glueball and matching it with an
experimentally observed particle would be a marvelous confirmation of
quantum chromodynamics. In addition, I find there to be something
charming about the idea that in a nonabelian gauge theory we could have
a particle made simply of the gauge field itself.

\begin{enumerate}
\def\labelenumi{\arabic{enumi})}
\setcounter{enumi}{4}
\tightlist
\item
  R.\ Plaga, ``Proposal for an experimental test of the many-worlds
  interpretation of quantum mechanics'', available as
  \href{https://arxiv.org/abs/quant-ph/9510007}{\texttt{quant-ph/9510007}}.
\end{enumerate}

\noindent
John Gribbin brought this one to my attention and asked me what I
thought about it. Basically, the idea here is to isolate an ion from its
environment in an ``ion trap'', and then perform a measurement on it with
two possible outcomes on another quantum system, and to excite the ion
only if the first outcome occurs, before the ion has had time to
``decohere'' or get ``entangled'' with the environment. Then one checks
to see if the ion is excited. The idea is that even if we didn't see the
outcome that made us excite the ion, we might see the ion excited,
because it was excited in the other ``world'' or ``branch'' --- the one
in which we \emph{did} see the outcome that made us excite the ion. The
author gets fairly excited himself, suggesting that ``outside physics,
interworld communication would lead to truly mind-boggling
possibilities''.

Does this idea really make sense? First of all, I don't think of this
sort of thing as a test of the many-worlds interpretation; I think that
all sufficiently sensible interpretations of quantum mechanics (not
necessarily \emph{very} sensible, either!) give the same concrete
predictions for all experiments, when it comes to what we actually
observe. They may make us tell very different stories about what is
happening behind the scenes, but not of any testable sort. As soon as
one comes up with something that makes different predictions, I think it
is (more or less by definition) not a new ``interpretation'' of quantum
theory but an actual new theory. And I don't think the many-worlds
interpretation is that.

So the question as I see it is simply, will this experiment work? Will
we sometimes see the ion excited even when we didn't excite it? It seems
hard; usually the decoherence between the two ``branches'' prevents this
kind of ``inter-world communication'' (not that I'm particularly fond of
this way of talking about it). What exactly is supposed to make this
case different? The problem is that the paper is quite sketchy at the
crucial point... just when the rabbit being pulled from the hat, as
it were. I haven't put much time into analyzing it, but some people
interested in this sort of thing might enjoy having a go at it.

\begin{enumerate}
\def\labelenumi{\arabic{enumi})}
\setcounter{enumi}{5}
\tightlist
\item
  Nicholas Landsman, ``Against the Wheeler--DeWitt equation'',
  available as
  \href{https://arxiv.org/abs/gr-qc/9510033}{\texttt{gr-qc/9510033}}.
\end{enumerate}
\noindent
I haven't read this one yet, but I had some nice talks with Landsman
about his ideas on quantization of constrained systems (see
\protect\hyperlink{week60}{``Week 60''}) back when I was in Cambridge,
England. Quantizing constrained systems is the main problem with the
so-called ``canonical'' approach to quantum gravity (see
\protect\hyperlink{week43}{``Week 43''}), so I was eager to see it
applied to gravity, and I guess that's what he's done. The title of the
paper is deliberately provocative... hmmm, I guess I'd better read
it soon! Here's the abstract:

\begin{quote}
The ADM approach to canonical general relativity combined with Dirac's
method of quantizing constrained systems leads to the Wheeler--DeWitt
equation. A number of mathematical as well as physical difficulties that
arise in connection with this equation may be circumvented if one
employs a covariant Hamiltonian method in conjunction with a recently
developed, mathematically rigorous technique to quantize constrained
systems using Rieffel induction. The classical constraints are cleanly
separated into four components of a covariant momentum map coming from
the diffeomorphism group of spacetime, each of which is linear in the
canonical momenta, plus a single finite-dimensional quadratic constraint
that arises in any theory, parametrized or not. The new quantization
method is carried through in a minisuperspace example, and is found to
produce a ``wavefunction of the universe''. This differs from the
proposals of both Vilenkin and Hartle-Hawking for a closed FRW universe,
but happens to coincide with the latter in the open case.
\end{quote}

\begin{enumerate}
\def\labelenumi{\arabic{enumi})}
\setcounter{enumi}{6}
\item
  Pavel Etingof and David Kazhdan, ``Quantization of Lie bialgebras,
  I'', available as
  \href{https://arxiv.org/abs/q-alg/9506005}{\texttt{q-alg/9506005}}.

  ``Quantization of Poisson algebraic groups and Poisson homogeneous
  spaces'', available as
  \href{https://arxiv.org/abs/q-alg/9510020}{\texttt{q-alg/9510020}}.
\end{enumerate}
\noindent
It sounds like Etinghof and Kazhdan are making serious progress on some
questions of Drinfeld about when you can quantize Lie bialgebras and
their kin. More stuff I need to read! I need to invent a time machine
and keep running it backwards to make my weekends longer and read this
stuff!

\begin{enumerate}
\def\labelenumi{\arabic{enumi})}
\setcounter{enumi}{7}
\item
  Steve Carlip, ``Statistical mechanics and black hole entropy'',
   available as
  \href{https://arxiv.org/abs/gr-qc/9509024}{\texttt{gr-qc/9509024}}.

  Claudio Teitelboim, ``Statistical thermodynamics of a black hole in
  terms of surface fields'', available as
  \href{https://arxiv.org/abs/hep-th/9510180}{\texttt{hep-th/9510180}}.
\end{enumerate}
\noindent
Steve Carlip's paper is a nice introduction to recent ideas, many of
them his, on deriving black hole area/entropy relations by thinking of
the entropy as associated to degrees of freedom of a field living on the
event horizon. I haven't read Teitelboim's paper, but it sounds related.

\begin{enumerate}
\def\labelenumi{\arabic{enumi})}
\setcounter{enumi}{8}
\item
  Jorge Griego, ``Is the third coefficient of the Jones knot polynomial
  a quantum state of gravity?'', available as
  \href{https://arxiv.org/abs/gr-qc/9510051}{\texttt{gr-qc/9510051}}.

  ``The Kauffman bracket and the Jones polynomial in quantum gravity'',
  available as
  \href{https://arxiv.org/abs/gr-qc/9510050}{\texttt{gr-qc/9510050}}.
\end{enumerate}
\noindent
In the loop representation of quantum gravity, states of quantum gravity
give rise to link invariants. Which link invariants, though? The
Kauffman bracket comes from a state of quantum gravity with cosmological
constant... that is something I understand pretty well by now. But
Gambini and Pullin also have an argument suggesting that the second
coefficient of the Jones polynomial (also known as the Arf invariant) is
a state of quantum gravity without cosmological constant. I've tried to
make this argument more rigorous and never succeeded. They also floated
a conjecture that \emph{all} the coefficients of the Jones polynomial
are states of quantum gravity. This confuses me a lot, because the Jones
polynomial depends on the orientations of the components of a link,
while states of quantum gravity should give link invariants that are
independent of orientations. I guess all the odd coefficients of the
Jones polynomial are orientation dependent. Thus I'm not shocked that
Griego has done calculations indicating that the third coefficient does
\emph{not} come from a state of quantum gravity.

\hypertarget{week69}{%
\section{November 11, 1995}\label{week69}}

One of the great things about starting to work on quantum gravity was
getting to know some of the people in the field. Ever since the
development of string theory and the loop representation of quantum
gravity, there has been a fair amount of interest in understanding how
quantum theory and gravity fit together. Indeed, now that the Standard
Model seems to be giving a spectacularly accurate description of all the
forces \emph{except} gravity, quantum gravity is one of the few really
big mysteries left when it comes to working out the basic laws of
physics --- or at least, one of the few \emph{obvious} big mysteries.
(As soon as one mystery starts becoming less mysterious, new mysteries
tend to become more visible.) But back when particle physics was big
business, only a few rather special sorts of people were seriously
devoted to quantum gravity. These people seem to be often more than
averagely interested in philosophy, often more interested in mathematics
(which is one of the few solid handholds in this slippery subject), and
always more resigned to the fact that Nature does not reveal all her
secrets very readily.

One of these folks is Chris Isham, whom I first saw at a conference in
Seattle in 1991. The conference was on classical field theory but
somehow he, Abhay Ashtekar, and Renate Loll sneaked in and gave some
talks on the loop representation of quantum gravity. This is when I
first became really interested in this subject, which I was later to
work on quite a bit. I remember Isham saying how he had been working on
quantum gravity for many years, and that he'd gotten used to the fact
that nothing ever worked, but that \emph{this} approach \emph{seemed} to
be working so far. He went on to talk about work he'd done with Ashtekar
on making the loop representation rigorous, which was based on
Gelfand--Naimark spectral theory. He said that as a student, when he'd
learned about this theory, he was really excited, because it completely
depends on the fact that if we have a space \(X\), we can think of any
point \(x\) in \(X\) as a functional on the space of functions on \(X\),
basically defining by defining \(x(f)\) to be \(f(x)\). He said this
with a laugh, but I knew what he meant, because I too had found this
idea tremendously exciting when I first learned the Gelfand--Naimark
theory. I guess it's something about how what seems at first like some
sort of bizarre joke can turn out to be very useful....

Anyway, later, when I decided to work on this sort of thing and was
trying to learn more about quantum gravity, I found his review article
on the problem of time (see \protect\hyperlink{week9}{``Week 9''})
tremendously helpful, and I constantly recommend it to everyone who is
trying to get their teeth into this somewhat elusive issue. So it's not
surprising that Isham figures prominently in the following nice popular
article on the problem of time:

\begin{enumerate}
\def\labelenumi{\arabic{enumi})}
\tightlist
\item
  Marcia Bartusiak, ``When the universe began, what time was it?'',
  \emph{Technology Review} (edited at the Massachusetts Institute of
  Technology), November/December 1995, pp.~54-63.
\end{enumerate}
\noindent
If you can find this, read it: it also features Karel Kuchar and Carlo
Rovelli.

This spring, I visited Isham at Imperial College in London and found him
to be just as interesting in person as in print, and not at all
scary... a bit of an cynic about all existing approaches to quantum
gravity (probably because he sees so clearly how flawed they all are),
but thoroughly good-humored about it and perfectly open to all sorts of
ideas, even my own nutty ideas about \(n\)-categories and physics.

Anyway, Isham has recently written a review article on quantum gravity
that gives a nice overview of the basic issues of the field:

\begin{enumerate}
\def\labelenumi{\arabic{enumi})}
\setcounter{enumi}{1}
\tightlist
\item
  C. J. Isham, ``Structural issues in quantum gravity'', plenary session
  lecture given at the GR14 conference, Florence, August 1995, 
  available as
  \href{https://arxiv.org/abs/gr-qc/9510063}{\texttt{gr-qc/9510063}}.
\end{enumerate}
\noindent
One interesting thing about it is the emphasis on the question of
whether spacetime is really a manifold the way we all usually think, or
perhaps something that just looks like a manifold at sufficiently large
distance scales. This is one of those fundamental issues that is rather
hard to make direct progress on; one has to sort of sneak up on it, but
it's nice to see someone boldly holding the problem up for examination.
Often the most important issues are the ones everyone is scared to talk
about, because they are so intractable.

Much of Abhay Ashtekar's early work dealt with asymptotically flat
solutions of Einstein's equation, but in about 1986 he somehow invented
a new formulation of general relativity, which everyone now calls the
``new variables'' or ``Ashtekar variables''. In terms of these new
variables general relativity looks a whole lot more like Yang--Mills
theory (the theory of all the forces \emph{except} gravity), and this
let Rovelli and Smolin formulate a radical new approach to quantum
gravity, the ``loop representation''. (For a fun, nontechnical
introduction to this, try the article by Bartusiak reviewed in
\protect\hyperlink{week10}{``Week 10''}.)

Nowadays, Ashtekar is the main person behind the drive to make the loop
representation of quantum gravity into a mathematically rigorous theory.
Thus it's natural that after that first time in Seattle I would wind up
seeing him pretty often... first at Syracuse University and then at
the Center for Gravitational Physics and Geometry which he started at
Penn State. It's really impressive how he has organized people into an
effective team there, and how he is systematically converting
people's hopes and dreams concerning the loop representation into a
beautiful set of rigorous \emph{theorems}. For a good mathematical
introduction to his program, see his paper reviewed in
\protect\hyperlink{week7}{``Week 7''}. A less mathematical introduction
is:

\begin{enumerate}
\def\labelenumi{\arabic{enumi})}
\setcounter{enumi}{2}
\tightlist
\item
  Abhay Ashtekar, ``Polymer geometry at Planck scale and quantum
  Einstein equations'', \emph{Int.\ J.\ Mod.\ Phys.\ D} \textbf{5} (1996),
   629--648.   Also available as \href{https://arxiv.org/abs/hep-th/9601054}{\texttt{hep-th/9601054}}.
\end{enumerate}

I have also seen Renate Loll fairly often in the years since that
Seattle conference. She is younger than Ashtekar and Isham (in fact, she
was a postdoc with Isham at one point), hence less intimidating to me,
which meant that I really enjoyed pestering her with stupid questions
when I was just starting to learn about this loop representation stuff.
One of her specialities is lattice gauge theory, and recently she has
developed a lattice version of quantum gravity that is eminently
suitable for computer calculations. The last time I saw her was at a
conference in Warsaw this spring (as reported in
\protect\hyperlink{week55}{``Week 55''} and
\protect\hyperlink{week56}{``Week 56''}). In the process of working on
her lattice approach, she gave Rovelli and Smolin a big shock by turning
up an error in their computation of the volume operator in quantum
gravity. A state of quantum gravity can be visualized roughly as a graph
embedded in space, with edges labelled by spins. Rovelli and Smolin had
thought there were states of nonzero volume corresponding to graphs with
only trivalent vertices (3 edges meeting a vertex, that is). As it turns
out, they'd made a sign error, and these states have zero volume; you
need a quadrivalent vertex to get some volume. She has just written a
paper on this topic:

\begin{enumerate}
\def\labelenumi{\arabic{enumi})}
\setcounter{enumi}{3}
\tightlist
\item
  Renate Loll, ``Spectrum of the volume operator in quantum gravity'',
  14 pages in plain tex, with 4 figures (postscript, compressed and
  uu-encoded), available as
  \href{https://arxiv.org/abs/gr-qc/9511030}{\texttt{gr-qc/9511030}}.
\end{enumerate}

The abstract reads as follows:

\begin{quote}
The volume operator is an important kinematical quantity in the
non-perturbative approach to four-dimensional quantum gravity in the
connection formulation. We give a general algorithm for computing its
spectrum when acting on four-valent spin network states, evaluate some
of the eigenvalue formulae explicitly, and discuss the role played by
the Mandelstam constraints.
\end{quote}

\begin{center}\rule{0.5\linewidth}{0.5pt}\end{center}

\begin{quote}
\emph{``Nothing is too wonderful to be true, if it be consistent with
the laws of nature, and in such things as these, experiment is the best
test of such consistency.''}

--- Faraday, laboratory diaries, entry 10,040, March 19, 1849.
\end{quote}

\hypertarget{week70}{%
\section{November 26, 1995}\label{week70}}

Probably many of the mathematicians reading this know about the Newton
Institute in Cambridge, a mathematics institute run by Sir Michael
Atiyah. It's a cozy little building, in a quiet neighborhood a certain
distance from the center of town, which one can reach by taking a nice
walk or bike ride over the bridge near Trinity College, across Grange
Road, and down Clarkson Road. Inside it's one big space, with stairways
slightly reminiscent of a certain picture by Escher, with a nice little
library on the first floor, tea and coffee on the 3rd floor, blackboards
in the bathrooms... everything a mathematician could want. This is
where Wiles first announced his proof of Fermat's last theorem, and they
sell T-shirts there commemorating that fact, which are unfortunately too
small to contain the proof itself... as they do not refrain from
pointing out.

I just got back from a conference there on New Connections between
Mathematics and Computer Science. It was organized by Jeremy
Gunawardena, who was eager to expose computer scientists and
mathematicians to a wide gamut of new interactions between the two
subjects. I spoke about \(n\)-categories in logic, topology and physics.
Since I don't know anything about computer science, when I first got the
invitation I thought it was a mistake: a wrong email address or
something! But Gunawardena assured me otherwise. I assumed the idea was
that \(n\)-categories, being so abstract, must have \emph{some}
application to just about \emph{everything}, even computer science.
Luckily, some other speakers at the conference gave some very nice
applications of \(n\)-category theory to computer science, so now I know
they really exist.

Unfortunately I had to miss the beginning of the conference, and
therefore missed some interesting talks of a geometrical nature by
Smale, Gromov, Shub and others. Let me say a bit about some of the talks
I did catch. You can find a list of all the speakers and abstracts of
their talks here:

\begin{enumerate}
\def\labelenumi{\arabic{enumi})}
\tightlist
\item
  Jeremy Gunawardena, ``New
  connections between mathematics and computer science: reports, abstracts
  and bibliography of a workshop'', January 1996.  Available at \href{https://www.hpl.hp.com/techreports/96/HPL-BRIMS-96-02.pdf}{\texttt{https://www.hpl.hp.com/techreports/96/HPL-BRIMS-96-02.pdf}}
\end{enumerate}
\noindent
Richard Jozsa gave an interesting talk on quantum computers, in part
outlining Peter Shor's work (see \protect\hyperlink{week34}{``Week
34''}) on efficient factoring via quantum computation, but also
presenting some new results on ``counterfactual quantum computation''.
It turns out that --- in principle --- in some cases you can get a
quantum computer to help you answer a question, even without running it,
just as long as you \emph{could have} run it! (I should add that in practice a
lot of things make this quite impractical.) This is a new twist on the
Elitzur-Vaidman bomb-testing paradox about how if you have a bunch of
bombs and half of them are duds, and the only way you can test a bomb is
by lighting the fuse and seeing if it goes off, you can still get a bomb
you're sure will work, if you use quantum mechanics. The trick involves
getting a fuse that's so sensitive that even one photon will make the
bomb go off, and then setting up a beam-splitter, and using the bomb to
measure which path the photon followed, before recombining the beams.
Check out:

\begin{enumerate}
\def\labelenumi{\arabic{enumi})}
\setcounter{enumi}{1}
\item
  Avshalom C. Elitzur and Lev Vaidman, ``Quantum mechanical interaction-free
  measurements'', \emph{Foundations of Phys.} \textbf{23} (1993),
  987--997.  Also available as \href{https://arxiv.org/abs/hep-th/9305002}{\texttt{hep-th/9305002}}.

  Graeme Mitchison and Richard Jozsa, Counterfactual quantum
  computation, \emph{Proc. Roy. Soc. Lond.} \textbf{A457} (2001)
  1175--1194. Also available as
  \href{http://arxiv.org/abs/quant-ph/9907007}{\texttt{quant-ph/9907007}}.
\end{enumerate}

Jean-Yves Girard gave an overview of linear logic. Linear logic is a new
version of logic that he invented, which has some new operations besides
the good old ones like ``and'', ``or'', and ``not''. For example, there
are things like ``par'' (written as an upside-down ampersand), ``!''
(usually pronounced ``bang'') and ``?''. Ever since I started going to
conferences on category theory and computer science I have been hearing
a lot about it, and I keep trying to get people to explain these weird
new logical operations to me. Unfortunately, I keep getting very
different answers, so it has remained rather mysterious to me, even
though it seems like a lot of fun (see \protect\hyperlink{week40}{``Week
40''}). Thus I was eager to hear it from the horse's mouth.

Indeed, Girard gave a fascinating talk on it which almost made me feel I
understood it. I think the big thing I'd been missing was a good
appreciation of topics in proof theory like ``cut elimination''. He
noted that this subject usually appears to be all about the precise
manipulation of formulas according to purely syntactic rules: ``Very
bureaucratic'' he joked, ``one parenthesis missing and you've had it!''
(For full effect, one must imagine this being said in a French accent by
someone stylishly dressed entirely in black.) He wanted to get a more
\emph{geometrical} way to think about proofs, but to do this it turned
out to be important to refine ordinary logic in certain ways\ldots.
leading to linear logic. However, I still don't feel up to explaining
it, so let me turn you to:

\begin{enumerate}
\def\labelenumi{\arabic{enumi})}
\setcounter{enumi}{2}
\item
  Jean-Yves Girard, ``Linear logic'', \emph{Theoretical Computer
  Science} \textbf{50}, 1--102, 1987.

  Jean-Yves Girard, Y. Lafont and P. Taylor, \emph{Proofs and Types},
 Cambridge U. Press, Cambridge, 1989. Also available at
  \href{http://www.paultaylor.eu/stable/Proofs+Types.html}{\texttt{http://www.paultaylor.eu/stable/Proofs+Types.html}}.
\end{enumerate}

Eric Goubault and Vaughan Pratt talked, in somewhat different ways,
about a formalism for treating concurrency using ``higher-dimensional
automata''. The basic idea is simple: say we have two jobs to do, one of
which gets us from some starting-point \(A\) to some result \(B\), and
the other of which gets us from \(A'\) to \(B'\). We can represent each
task by an arrow, as follows: \[
  \begin{aligned}
    A&\longrightarrow B
  \\A'&\longrightarrow B'
  \end{aligned}
\] We can think of this arrow as a ``morphism'', that is, a completely
abstract sort of operation taking \(A\) to \(B\). Or, we can think of it
more concretely as an interval of time, where our computer is doing
something at each moment. Alternatively, we can think of it more
discretely as a sequence of steps, starting with \(A\) and winding up
with \(B\).

If we now consider doing both these tasks concurrently, we can represent
the situation by a square: \[
  \begin{tikzcd}
    AA' \rar\dar & BA' \dar
  \\AB' \rar & BB'
  \end{tikzcd}
\] Going first across and then down corresponds to completing one task
before starting the other, while going first down and then across
corresponds to doing the other one first. However, we can also imagine
various roughly diagonal paths through the square, corresponding to
doing both tasks at the same time. We might go horizontally for a while,
then vertically, then diagonally, and so on. Of course, if the two tasks
were not completely independent --- for example, if some steps of one
could only occur after some steps of the other were finished --- we
would have some constraints on what paths from \(AA'\) to \(BB'\) were
allowed. The idea is then to model these constaints as ``holes'' in the
square, forbidden regions where the path cannot go. There may then be
several ``essentially distinct'' ways of getting from \(AA'\) to
\(BB'\), that is, classes of paths that cannot be deformed into each
other.

To anyone who knows homotopy theory, this will seem very familiar,
homotopy theory being all about spaces with holes in them, and how those
holes prevent you from continuously deforming one path into another.
Goubault's title, ``Scheduling problems and homotopy theory'',
emphasized the relationships. But there are also some big differences.
Unlike homotopy theory, here the paths are typically required to be
``monotonic'': they can't double back and go backwards in time. And, as
I mentioned, the tasks can be thought of more abstractly than as paths
in some space. So we are really talking about \(2\)-categories here:
they give a general framework for studying situations with ``dots'' or
``objects'', ``arrows between dots'' or ``morphisms'', and ``arrows
between arrows between dots'' or ``2-morphisms''. Similarly, when we
study concurrency with more than 2 tasks at a time we can think of it in
terms of \(n\)-categories.

By the way, since I don't know much about parallel processing, I'm not
sure how much the above formalism actually helps the ``working man''.
Probably not much, yet. I get the impression, however, that parallel
processing is a complicated problem, and that people are busily dreaming
up new formalisms for talking about it, hoping they will eventually be
useful for inventing and analyzing parallel programming languages.

Some references for this are:

\begin{enumerate}
\def\labelenumi{\arabic{enumi})}
\setcounter{enumi}{3}
\item
  Eric Goubault, ``Schedulers as abstract interpretations of
  higher-dimensional automata'', in \emph{Proc. PEPM '95 (La Jolla)}, ACM
  Press, 1995. Also available at
  \href{https://dl.acm.org/doi/abs/10.1145/215465.215577}{\texttt{https://dl.acm.org/doi/abs/10.1145/215465.215577}}.

  Eric Goubault and Thomas Jensen, ``Homology of higher-dimensional
  automata'', in \emph{Proc. CONCUR '92 (New York)}, Lecture Notes in
  Computer Science \textbf{630}, Springer, 1992. Also available at
  \href{http://www.lix.polytechnique.fr/~goubault/papers/Homology.pdf}{\texttt{http://www.lix.polytechnique.fr/\(\sim\)goubault/papers/Homology.pdf}}.
\item
  Vaughan Pratt, ``Time and information in sequential and concurrent
  computation'', in \emph{Proc. Theory and Practice of Parallel
  Programming}, Sendai, Japan, 1994.  Also available at \href{http://boole.stanford.edu/pub/tppp.pdf}{\texttt{http://boole.stanford.edu/pub/tppp.pdf}}.
\end{enumerate}

Yves Lafont also gave a talk with strong connections to \(n\)-category
theory. Recall that a monoid is a set with an associative product having
a unit element. One way to describe a monoid is by giving a presentation
with ``generators'', say \[a, b, c, d,\] and ``relations'', say
\[ab = a,\quad da = ac.\] We get a monoid out of this in an obvious sort
of way, namely by taking all strings built from the generators
\(a\),\(b\),\(c\), and \(d\), and then identifying two strings if you
can get from one to the other by repeated use of the relations. In math
jargon, we form the free monoid on the generators and then mod out by
the relations.

Suppose our monoid is finitely presented, that is, there are finitely
many generators and finitely many relations. How can we tell whether two
elements of it are equal? For example, does \[dacb = acc\] in the above
monoid? Well, if the two are equal, we will always eventually find that
out by an exhaustive search, applying the relations mechanicallly in all
possible ways. But if they are not, we may never find out! (For the
above example, the answer appears at the end of this article in case
anyone wants to puzzle over it. Personally, I can't stand this sort of
puzzle.) In fact, there is no general algorithm for solving this ``word
problem for monoids'', and in fact one can even write down a
\emph{specific} finitely presented monoid for which no algorithm works.

However, sometimes things are nice. Suppose you write the relations as
``rewrite rules'', that go only one way: \[
  \begin{aligned}
    ab &\to a
  \\da &\to ac
  \end{aligned}
\] Then if you have an equation you are trying to check, you can try to
repeatedly apply the rewrite rules to each side, reducing it to ``normal
form'', and see if the normal forms are equal. This will only work,
however, if some good things happen! First of all, your rewrite rules
had better terminate: it had better be that you can only apply them
finitely many times to a given string. This happens to be true for the
above pair of rewrite rules, because both rules decrease the number of
\(b\)'s and \(c\)'s. Second of all, your rewrite rules had better be
confluent: it had better be that if I use the rules one way until I
can't go any further, and you use them some other way, that we both wind
up with the same thing! If both these hold, then we can reduce any
string to a unique normal form by applying the rules until we can't do
it any more.

Unfortunately, the rules above aren't confluent; if we start with the
word \(dab\), you can apply the rules like this \[dab \to acb\] while I
apply them like this \[dab \to da \to ac\] and we both terminate, but at
different answers. We could try to cure this by adding a new rule to our
list, \[acb \to ac.\] This is certainly a valid rule, which cures the
problem at hand... but if we foolishly keep adding new rules to our
list this way we may only succeed in getting confluence and termination
when we have an \emph{infinite} list of rules: \[
  \begin{aligned}
    ab &\to a
  \\da &\to ac
  \\acb &\to ac
  \\accb &\to acc
  \\acccb &\to accc
  \\accccb &\to acccc
  \\\vdots & \vdots
  \end{aligned}
\] and so on. I leave you to check that this is really terminating and
confluent. Because it is, and because it's a very predictable list of
rules, we can use it to write a computer program in this case to solve
the word problem for the monoid at hand. But in fact, if we had been
cleverer, we could have invented a \emph{finite} list of rules that was
terminating and confluent: \[
  \begin{aligned}
    ab &\to a
  \\ac &\to da
  \end{aligned}
\] Lafont went on to describe some work by Squier:

\begin{enumerate}
\def\labelenumi{\arabic{enumi})}
\setcounter{enumi}{5}
\item
  Craig C. Squier, ``Word problems and a homological finiteness
  condition for monoids'', \emph{Jour. Pure Appl. Algebra} \textbf{49}
  (1987), 201--217.   Also available at \href{https://link.springer.com/content/pdf/10.1007\%2F3-540-17220-3_7.pdf}{\texttt{https://link.springer.com/content/pdf/10.1007\%2F3-540-17220-3\(\underline{\;}\)7.pdf}}

  Craig C. Squier, ``A finiteness condition for rewriting systems'',
  revision by F. Otto and Y. Kobayashi,  \emph{Theoretical Computer Science}
  \textbf{131} (1994), 271--294.

  Craig C. Squier and F. Otto, ``The word problem for finitely presented
  monoids and finite canonical rewriting systems'', in \emph{Rewriting
  Techniques and Applications}, ed.~J. P. Jouannuad, Lecture Notes in
  Computer Science \textbf{256}, Springer, Berlin, 1987, pp.\ 74--82.
  Also available at \href{https://link.springer.com/content/pdf/10.1007\%2F3-540-17220-3_7.pdf}{\texttt{https://link.springer.com/content/pdf/10.1007\%2F3-540-17220-3\(\underline{\;}\)7.pdf}}
\end{enumerate}
\noindent
which gave general conditions which must hold for there to be a finite
terminating and confluent set of rewrite rules for a monoid. The nice
thing is that this relies heavily on ideas from \(n\)-category theory.
Note: we started with a monoid in which the relations are
\emph{equations}, but we then started thinking of the relations as
rewrite rules or \emph{morphisms}, so what we really have is a monoidal
category. We then started worrying about ``confluences'', or equations
between these morphisms. This is typical of ``categorification'', in
which equations are replaced by morphisms, which we then want to satisfy
new equations (see \protect\hyperlink{week38}{``Week 38''}).

For the experts, let me say exactly how it all goes. Given any monoid
\(M\), we can cook up a topological space called its ``classifying
space'' \(KM\), as follows. We can think of \(KM\) as a simplicial
complex. We start by sticking in one 0-simplex, which we can visualize
as a dot like this: \[\bullet\] Then we stick in one \(1\)-simplex for
each element of the monoid, which we can visualize as an arrow going
from the dot to itself. Unrolled a bit, it looks like this: \[
  \begin{tikzpicture}
    \draw[thick] (0,0) node{$\bullet$} to node[fill=white]{$a$} (1,0) node{$\bullet$};
  \end{tikzpicture}
\] Really we should draw an arrow going from left to right, but soon
things will get too messy if I do that, so I won't. Then, whenever we
have \(ab = c\) in the monoid, we stick in a \(2\)-simplex, which we can
visualize as a triangle like this: \[
  \begin{tikzpicture}
    \draw[thick] (0,0) node{$\bullet$} to node[fill=white]{$c$} (1.5,0) node{$\bullet$} to node[fill=white]{$b$} (0.75,1.3) node{$\bullet$} to node[fill=white]{$a$} cycle;
  \end{tikzpicture}
\] Then, whenever we have \(abc = d\) in the monoid, we stick in a
\(3\)-simplex, which we can visualize as a tetrahedron like this \[
  \begin{tikzpicture}
    \draw[thick] (0,0) node{$\bullet$} to node[fill=white]{$d$} (3,0) node{$\bullet$} to node[fill=white]{$bc$} (1.5,2.6) node{$\bullet$} to node[fill=white]{$a$} cycle;
    \draw[thick] (0,0) to node[fill=white]{$ab$} (1.5,1) node{$\bullet$};
    \draw[thick] (1.5,2.6) to node[fill=white]{$b$} (1.5,1);
    \draw[thick] (3,0) to node[fill=white]{$c$} (1.5,1);
  \end{tikzpicture}
\] And so on.... This is a wonderful space whose homology groups
depend only on the monoid, so we can call them \(H_k(M)\). If we have a
presentation of \(M\) with only finitely many generators, we can build
\(KM\) using \(1\)-simplices only for those generators, and it follows
that \(H_1(M)\) is finitely generated. (More precisely, we can build a
space with the same homotopy type as \(KM\) using only the generators in
our presentation.) Similarly, if we have a presentation with only
finitely many relations, we can build \(KM\) using only finitely many
\(2\)-simplices, so \(H_2(M)\) is finitely generated. What Squier showed
is that if we can find a finite list of rewrite rules for M which is
terminating and confluent, then we can build \(KM\) using only finitely
many \(3\)-simplices, so \(H_3(M)\) is finitely generated! What's nice
about this is that homological algebra gives an easy way to compute
\(H_k(M)\) given a presentation of \(M\), so in some cases we can
\emph{prove} that a monoid has no finite list of rewrite rules for \(M\)
which is terminating and confluent, just by showing that \(H_3(M)\) is
too big. Examples of this, and many further details, appear in Lafont's
work:

\begin{enumerate}
\def\labelenumi{\arabic{enumi})}
\setcounter{enumi}{6}
\item
  Yves Lafont and Alain Proute, ``Church-Rosser property and homology of
  monoids'', \emph{Mathematical Structures in Computer Science}
  \textbf{1} (1991), 297--326. Also available at
  \href{http://iml.univ-mrs.fr/\~lafont/publications.html}{\texttt{http://iml.univ-mrs.fr/\(\sim\)lafont/publications.html}}
 
  Yves Lafont, ``A new finiteness condition for monoids presented by
  complete rewriting systems (after Craig C. Squier)'', \emph{Journal of
  Pure and Applied Algebra} \textbf{98} (1995), 229--244. Also available
  at \href{http://iml.univ-mrs.fr/\~lafont/publications.html}
  {\texttt{http://iml.univ-mrs.fr/\(\sim\)lafont/publications.html}}
\end{enumerate}

There were many other interesting talks, but I think I will quit here.
Next time I want to talk a bit about topological quantum field theory.
(Of course, folks who read \protect\hyperlink{week38}{``Week 38''} will
know that Lafont's work is deeply related to topological quantum field
theory... but I won't go into that now.)

\begin{center}\rule{0.5\linewidth}{0.5pt}\end{center}

(Answer: \(dacb = ddab = dda = dac = acc\).)

\hypertarget{week71}{%
\section{December 3, 1995}\label{week71}}

This week I will get back to mathematical physics... but before I
do, I can't resist adding that in my talk in that conference I announced
that James Dolan and I had come up with a definition of weak
\(n\)-categories. (For what these are supposed to be, and what they have
to do with physics, see \protect\hyperlink{week38}{``Week 38''},
\protect\hyperlink{week49}{``Week 49''}, and
\protect\hyperlink{week53}{``Week 53''}.) John Power was at the talk,
and before long his collaborator Ross Street sent me some email from
Australia asking about the definition. Gordon, Power, and Street have
done a lot of work on \(n\)-categories --- see
\protect\hyperlink{week29}{``Week 29''}. Now, Dolan and I have been
struggling for several months to put this definition onto paper in a
reasonably elegant and comprehensible form, so Street's request was a
good excuse to get something done quickly... sacrificing
comprehensibility for terseness. The result can be found in

\begin{enumerate}
\def\labelenumi{\arabic{enumi})}
\tightlist
\item
  John Baez and James Dolan, ``\(n\)-Categories, sketch of a
  definition'', available at \href{http://math.ucr.edu/home/baez/ncat.def.html}{\texttt{http://math.ucr.edu/home/baez/ncat.def.html}}
\end{enumerate}
\noindent
A more readable version will appear as a paper fairly soon. I should add
that this will eventually be part of Dolan's thesis, and he has done
most of the hard work on it; my role has largely been to push him along
and get him to explain everything to me.

On to some physics....

First, it's amusing to note that Maxwell's equations are back in
fashion! In the following papers you will see how the duality symmetry
of Maxwell's equations (the symmetry between electric and magnetic
fields) plays a new role in modern work on \(4\)-dimensional gauge
theory. Also, there is some good stuff in Thompson's review article
about the Seiberg--Witten theory, which is basically just a
\(\mathrm{U}(1)\) gauge theory like Maxwell's equations... but with
some important extra twists!

\begin{enumerate}
\def\labelenumi{\arabic{enumi})}
\setcounter{enumi}{1}
\item
  Erik Verlinde, ``Global aspects of electric-magnetic duality'',
  \emph{Nucl. Phys. B} \textbf{455} (1995), 211--225.  Also available as
  \href{http://arxiv.org/abs/hep-th/9506011}{\texttt{hep-th/9506011}}.

  George Thompson, ``New results in topological field theory and abelian
  gauge theory'', available as
  \href{http://arxiv.org/abs/hep-th/9511038}{\texttt{hep-th/9511038}}.
\end{enumerate}

Next, it's nice to see that work on the loop representation of quantum
gravity continues apace:

\begin{enumerate}
\def\labelenumi{\arabic{enumi})}
\setcounter{enumi}{2}
\item
  Thomas Thiemann, ``An account of transforms on
  \(\overline{\mathcal{A}/\mathcal{G}}\)'', \emph{Acta Cosmologica}
  \textbf{21} (1996), 145--167.  Also available as
  \href{http://arxiv.org/abs/gr-qc/9511049}{\texttt{gr-qc/9511049}}.

  Thomas Thiemann, ``Reality conditions inducing transforms for quantum
  gauge field theory and quantum gravity'', \emph{Class.\ Quant.\ Grav.} \textbf{13}
   (1996) 1383-1404.   Also available as
  \href{http://arxiv.org/abs/gr-qc/9511057}{\texttt{gr-qc/9511057}}.

  Abhay Ashtekar, ``A generalized Wick transform for gravity'', \emph{Phys.\ Rev.\ D}
  \textbf{53} (1996)   2865-2869.  Also available as
  \href{http://arxiv.org/abs/gr-qc/9511083}{\texttt{gr-qc/9511083}}.

  Renate Loll, ``Making quantum gravity calculable'', \emph{Acta Cosmologica}      
  \textbf{21} (1995), 131--144.  Also available
  as \href{http://arxiv.org/abs/gr-qc/9511080}{\texttt{gr-qc/9511080}}.

  Rodolfo Gambini and Jorge Pullin, ``A rigorous solution of the quantum
  Einstein equations'', \emph{Phys.\ Rev.\ D} \textbf{54} (1996), 5935--5938.  
  Also available as
  \href{http://arxiv.org/abs/gr-qc/9511042}{\texttt{gr-qc/9511042}}.
\end{enumerate}

The first three papers here discuss some new work tackling the ``reality
conditions'' and ``Hamiltonian constraint'', two of the toughest issues
in the loop representation of quantum gravity. First, the Hamiltonian
constraint is another name for the Wheeler--DeWitt equation
\[H \psi = 0\] that every physical state of quantum gravity must satisfy
(see \protect\hyperlink{week11}{``Week 11''} for why). The ``reality
conditions'' have to do with the fact that this constraint looks
different depending on whether we are working with Riemannian or
Lorentzian quantum gravity. The constraint is simpler when we work with
Riemannian quantum gravity. Classically, in \emph{Riemannian} gravity
the metric on spacetime looks like \[dt^2 + dx^2 + dy^2 + dz^2\] at any
point, if we use suitable local coordinates. In this Riemannian world,
time is no different from space! In the real world, the world of
\emph{Lorentzian} gravity, the metric looks like
\[-dt^2 + dx^2 + dy^2 + dz^2\] at any point, in suitable coordinates.
Folks often call the Riemannian world the world of ``imaginary time'',
since in some vague sense we can get from the Lorentzian world to the
Riemannian world by making the transformation \[t \mapsto it,\] called a
``Wick transform''. It looks simple the way I have just written it, in
local coordinates. But a priori it's far from clear that this Wick
transform makes any sense globally. Apparently, however, there is
something we can do along these lines, which transforms the Hamiltonian
for Lorentzian quantum gravity to the better-understood one of
Riemannian quantum gravity! Alas, I have been too distracted by
\(n\)-categories to keep up with the latest work on this, but I'll catch
up in a bit. Maybe over Christmas I can relax a bit, lounge in front of
a nice warm fire, and read these papers.

The goal of all these machinations, of course, is to find some equations
that make mathematical sense, have a good right to be called a
``quantized version of Einstein's equation'', and let one compute
answers to some physics problems. We don't expect that quantum gravity
is enough to describe what's really going on in interesting
problems... there are lots of other forces and particles out there.
Indeed, string theory is founded on the premise that quantum gravity is
completely inseparable from the quantum theories of everything else. But
here we are taking a different tack, treating quantum gravity by itself
in as simple a way as possible, expecting that the predictions of theory
will be \emph{qualitatively} of great interest even if they are
quantitatively inaccurate.

As described in earlier Finds (\protect\hyperlink{week55}{``Week 55''},
\protect\hyperlink{week68}{``Week 68''}), Loll has been working to make
quantum gravity ``calculable'', by working on a discretized version of
the theory called lattice quantum gravity. If one does it carefully,
it's not too bad to treat space as a lattice in the loop representation
of quantum gravity, because even in the continuum the theory is discrete
in a certain sense, since the states are described by ``spin networks'',
certain graphs embedded in space. (See \protect\hyperlink{week43}{``Week
43''} for more on these.) Her latest paper is an introduction to some of
these issues.

In a somewhat different vein, Gambini and Pullin have been working on
relating the loop representation to knot theory. One of their most
intriguing results is that the second coefficient of the
Alexander--Conway knot polynomial is a solution of the Hamiltonian
constraint. Here they argue for this result using a lattice
regularization of the theory.

Now let me turn to a variety of other matters...

\begin{enumerate}
\def\labelenumi{\arabic{enumi})}
\setcounter{enumi}{3}
\item
  Matt Greenwood and Xiao-Song Lin, ``On Vassiliev knot invariants
  induced from finite type'', available as
  \href{http://arxiv.org/abs/q-alg/9506001}{\texttt{q-alg/9506001}}.

  Lev Rozansky, ``On finite type invariants of links and rational
  homology spheres derived from the Jones polynomial and
  Witten--Reshetikhin--Turaev invariant'', in \emph{Geometry and Physics},
  CRC Press, Boca Raton, 2021, pp.\ 379--397. Also available as
  \href{http://arxiv.org/abs/q-alg/9511025}{\texttt{q-alg/9511025}}.

  Scott Axelrod, ``Overview and warmup example for perturbation theory
  with instantons'',  in \emph{Geometry and Physics},
  CRC Press, Boca Raton, 2021, pp.\ 321--338.   Also available as
  \href{http://arxiv.org/abs/hep-th/9511196}{\texttt{hep-th/9511196}}.
\end{enumerate}
\noindent
These papers all deal with perturbative Chern--Simons theory and its
spinoffs. The first two consider homology 3-spheres. In Chern--Simons
theory, this makes the moduli space of flat \(\mathrm{SU}(2)\)
connections trivial, thus eliminating some subtleties in the
perturbation theory. A homology 3-sphere is a 3-manifold whose homology
is the same as that of the 3-sphere... the first one was discovered
by Poincar\'e when he was studying his conjecture that every 3-manifold
with the homology of a 3-sphere is a 3-sphere. It turns out that you can
get a counterexample if you just take an ordinary 3-sphere, cut out a
solid torus embedded in the shape of a trefoil knot, and stick it back
in with the meridian and longitude (the short way around, and the long
way around) switched --- making sure they wind up pointing in the
correct directions. This is called ``doing Dehn surgery on the
trefoil''. It gives something with the homology of a 3-sphere that's not
a 3-sphere. So Poincar\'e had to revise his conjecture to say that every
3-manifold \emph{homotopic} to a 3-sphere is (homeomorphic to) a
3-sphere. This improved ``Poincar\'e conjecture'' remains unsolved...
its analog is known to be true in every dimension other than 3! Since
every possible counterexample to the Poincar\'e conjecture is a homology
3-sphere, it makes some sense to ponder these manifolds carefully.

Now, just as perturbative Chern--Simons theory gives certain special
invariants of links, said to be of ``finite type'', the same is true for
homology 3-spheres. When we say a link invariants is of finite type, all
we mean is that it satisfies a simple property described in
\protect\hyperlink{week3}{``Week 3''}. There is a similar but subtler
definition for an invariant of homology 3-spheres to be of finite type;
they need to transform in a nice way under Dehn surgery. (See
\protect\hyperlink{week60}{``Week 60''} for more references.)

The second paper concentrates precisely on those subtleties due to the
moduli space of flat connections, developing perturbation theory in the
presence of ``instantons'' (here, nontrivial flat connections).

\begin{enumerate}
\def\labelenumi{\arabic{enumi})}
\setcounter{enumi}{4}
\item
  Alan Carey, Jouko Mickelsson, and Michael Murray, ``Index theory,
  gerbes, and Hamiltonian quantization'', available as
  \href{http://arxiv.org/abs/hep-th/9511151}{\texttt{hep-th/9511151}}.

  Alan Carey, M. K. Murray and B. L. Wang, ``Higher bundle gerbes and
  cohomology classes in gauge theories'', available as
  \href{http://arxiv.org/abs/hep-th/9511169}{\texttt{hep-th/9511169}}
\end{enumerate}

Higher-dimensional algebra is sneaking into physics in yet another way:
gerbes! What's a gerbe? Roughly speaking, it's a sheaf of groupoids, but
there are some other ways of thinking of them that come in handy in
physics. See \protect\hyperlink{week25}{``Week 25''} for a review of
Brylinski's excellent book on gerbes, and also:

\begin{enumerate}
\def\labelenumi{\arabic{enumi})}
\setcounter{enumi}{5}
\item
  Jean-Luc Brylinski, ``Holomorphic gerbes and the Beilinson
  regulator'', \emph{Ast\'erisque} \textbf{226} (1994), 145--174.  Also 
  available at 
  \url{http://www.numdam.org/item/AST_1994__226__145_0/}

  Jean-Luc Brylinski, ``The geometry of degree-four characteristic
  classes and of line bundles on loop spaces I'', Duke Math. Jour. 75
  (1994), 603--638.

  Jean-Luc Brylinski and D. A. McLaughlin, ``\v Cech cocycles for characteristic classes'',
  \emph{Commun. Math. Phys.} \textbf{178} (1996), 225--236.  
  
\item
  Joe Polchinski, ``Recent results in string duality'', \emph{Prog.\ Theor.\ Phys.\ 
  Suppl.} \textbf{123} (1996), 9--18. Also available as
  \href{http://arxiv.org/abs/hep-th/9511157}{\texttt{hep-th/9511157}}.
\end{enumerate}

This should help folks keep up with the ongoing burst of work on
dualities relating superficially different string theories.

\begin{enumerate}
\def\labelenumi{\arabic{enumi})}
\setcounter{enumi}{7}
\tightlist
\item
  Leonard Susskind and John Uglum, ``String physics and black holes'',
  \emph{Nucl.\ Phys.\ Proc.\ Suppl.} \textbf{45} (1996), 115--134.  Also
  available as
  \href{http://arxiv.org/abs/hep-th/9511227}{\texttt{hep-th/9511227}}.
\end{enumerate}

Among other things, this review discusses the ``holographic hypothesis''
mentioned in \protect\hyperlink{week57}{``Week 57''}:

\begin{enumerate}
\def\labelenumi{\arabic{enumi})}
\setcounter{enumi}{8}
\tightlist
\item
  Boguslaw Broda, ``A gauge-field approach to 3- and 4-manifold
  invariants'', available as
  \href{http://arxiv.org/abs/q-alg/9511010}{\texttt{q-alg/9511010}}.
\end{enumerate}
\noindent
This summarizes the Reshetikhin--Turaev construction of 3d topological
quantum field theories from quantum groups, and Broda's own closely
related approach to 4d topological quantum field theories.

\begin{enumerate}
\def\labelenumi{\arabic{enumi})}
\setcounter{enumi}{9}
\tightlist
\item
  John Baez and Martin Neuchl, ``Higher-dimensional algebra I: braided
  monoidal \(2\)-categories'', \emph{Adv. Math.} \textbf{121} (1996), 
  196--244.  Also available as
  \href{http://arxiv.org/abs/q-alg/9511013}{\texttt{q-alg/9511013}}.
\end{enumerate}
\noindent
In this paper, we begin with a brief sketch of what is known and
conjectured concerning braided monoidal \(2\)-categories and their
applications to 4d topological quantum field theories and 2-tangles
(surfaces embedded in \(4\)-dimensional space). Then we give concise
definitions of semistrict monoidal \(2\)-categories and braided monoidal
\(2\)-categories, and show how these may be unpacked to give long
explicit definitions similar to, but not quite the same as, those given
by Kapranov and Voevodsky. Finally, we describe how to construct a
semistrict braided monoidal \(2\)-category \(Z(\mathcal{C})\) as the
`center' of a semistrict monoidal category \(\mathcal{C}\). This is
analogous to the construction of a braided monoidal category as the
center, or `quantum double', of a monoidal category. The idea is to
develop algebra that will do for 4-dimensional topology what quantum
groups and braided monoidal categories did for 3d topology. As a
corollary of the center construction, we prove a strictification theorem
for braided monoidal \(2\)-categories.

\hypertarget{week72}{%
\section{February 1, 1996}\label{week72}}

It's been a while since I've written an issue of This Week's
Finds... due to holiday distractions and a bunch of papers that
need writing up. But tonight I just can't seem to get any work done, so
let me do a bit of catching up.

I'm no string theorist, but I still can't help hearing all the rumbling
noises over in that direction: first about all the dualities relating
seemingly different string theories, and then about the mysterious
``M-theory'' in 11 dimensions which seems to underlie all these
developments. Let me try to explain a bit of this stuff... in the
hopes that I prompt some string theorists to correct me and explain it
better! I will simplify everything a lot to keep people from getting
scared of the math involved. But I may also make some mistakes, so the
experts should be kind to me and try to distinguish between the
simplifications and the mistakes.

Recall that it's hard to get a consistent string theory --- one that's
not plagued by infinite answers to interesting questions. But this
difficulty is generally regarded as a good thing, because it drastically
limits the number of different versions of string theory one needs to
think about. It's often said that there are only 5 consistent string
theories: the type I theory, the type IIA and IIB theory, and the two
kinds of heterotic string theory. I'm not sure exactly what this
statement means, but certainly it's only meant to cover supersymmetric
string theories, which can handle fermions (like the electron and
neutrino) in addition to bosons (like the photon).

Type I strings are ``open strings'' --- not closed loops --- and they
live in 10 dimensional spacetime, meaning that you need the dimension to
be 10 to make certain nasty infinities cancel out. Type II strings also
live in 10 dimensions, but they are ``closed strings''. That means that
they look like a circle, so there are vibrational modes that march
around clockwise and other modes that march around counterclockwise, and
these are supposed to correspond to different particles that we see. We
can think of these vibrational modes as moving around the circle at the
speed of light; they are called ``left-movers'' and ``right-movers''.
Now fermions which move at the speed of light are able to be rather
asymmetric and only spin one way (when viewed head-on). We say they have
a ``chirality'' or handedness. Ordinary neutrinos, for example, are
left-handed. This asymmetry of nature shocked everyone when first
discovered, but it appears to be a fact of life, and it's certainly a
fact of mathematics. In the type IIA string theory, the left-moving and
right-moving fermionic vibrational modes have opposite chiralities,
while in the IIB theory, they have the same chirality. When I last
checked, the type IIA theory seemed to fit our universe a bit better
than the IIB theory.

But lots of people say the heterotic theory matches our universe even
better. The name ``heterotic'' refers to the fact that this theory is
supposed to have ``hybrid vigor''. It's quite bizarre: the left-movers
are purely bosonic --- no fermions --- and live in \(26\)-dimensional
spacetime, the way non-supersymmetric string theories do. The
right-movers, on the other hand, are supersymmetric and live in 
\(10\)-dimensional spacetime. It sounds not merely heterotic, but downright
schizophrenic! But in fact, the \(26\)-dimensional spacetime can also be
thought of as being \(10\)-dimensional, with \(16\) extra ``curled-up
dimensions'' in the shape of a torus. This torus has two possible
shapes: \(\mathbb{R}^{16}\) modulo the
\(\mathrm{E}_8 \times \mathrm{E}_8\) lattice or the \(D_{16}^*\)
lattice. (For some of the wonders of \(\mathrm{E}_8\) and other
lattices, check out \protect\hyperlink{week64}{``Week 64''} and
\protect\hyperlink{week65}{``Week 65''}. The \(D_{16}^*\) lattice is
related to the \(D_{16}\) lattice described in those Weeks, but not
quite the same.)

Now there is still lots of room for toying with these theories depending
on how you ``compactify'': how you think of \(10\)-dimensional spacetime
as 4-dimensional spacetime plus 6 curled-up dimensions. That's because
there are lots of \(6\)-dimensional manifolds that will do the job (the
so-called ``Calabi-Yau'' manifolds). Different choices give different
physics, and there is a lot of work to be done to pick the right one.

However, recently it's beginning to seem that all five of the basic
sorts of string theory are beginning to look like different
manifestations of the same theory in 11 dimensions... some
monstrous thing called M-theory! Let me quote the following paper:

\begin{enumerate}
\def\labelenumi{\arabic{enumi})}
\tightlist
\item
  Kelly Jay Davis, ``M-Theory and string-string duality'', 
  available as
  \href{https://arxiv.org/abs/hep-th/9601102}{\texttt{hep-th/9601102}}.
\end{enumerate}

The idea seems to be roughly that depending on how one compactifies the
11th dimension, one gets different \(10\)-dimensional theories from
M-theory:

\begin{quote}
``In the past year much has happened in the field of string theory. Old
results relating the two Type II string theories and the two Heterotic
string theories have been combined with newer results relating the Type
II theory and the Heterotic theory, as well as the Type I theory and the
Heterotic theory, to obtain a single ``String Theory." In addition, there
has been much recent progress in interpreting some, if not all,
properties of String Theory in terms of an eleven-dimensional M-Theory.
In this paper we will perform a self-consistency check on the various
relations between M-Theory and String Theory. In particular, we will
examine the relation between String Theory and M-Theory by examining its
consistency with the string-string duality conjecture of six-dimensional
String Theory. So, let us now take a quick look at the relations between
M-Theory and String Theory some of which we will be employing in this
article.

In Witten's paper he established that the strong coupling limit of Type
IIA string theory in ten dimensions is equivalent to eleven-dimensional
supergravity on a ``large'' \(S^1\). {[}Note: \(S^1\) just means the
circle --- jb.{]} As the low energy limit of M-theory is
eleven-dimensional supergravity, this relation states that the strong
coupling limit of Type IIA string theory in ten-dimensions is equivalent
to the low-energy limit of M-Theory on a ``large'' \(S^1\). In the paper
of Witten and Horava, they establish that the strong coupling limit of
the ten-dimensional \(\mathrm{E}_8 \times \mathrm{E}_8\) Heterotic
string theory is equivalent to M-Theory on a ``large''
\(S^1/\mathbb{Z}_2\).

Recently, Witten, motivated by Dasgupta and Mukhi, examined M-Theory on
a \(\mathbb{Z}_2\) orbifold of the five-torus and established a relation
between M-Theory on this orbifold and Type IIB string theory on \(\mathrm{K}3\).
{[}Note: most of these undefined terms refer to various spaces; for
example, the five-torus is the \(5\)-dimensional version of a doughnut,
while \(\mathrm{K}3\) is a certain \(4\)-dimensional manifold --- jb.{]} Also,
Schwarz very recently looked at M-Theory and its relation to T-Duality.

As stated above, M-Theory on a ``large'' \(S^1\) is equivalent to a
strongly coupled Type IIA string theory in ten-dimensions. Also,
M-theory on a ``large'' \(S^1/\mathbb{Z}_2\) is equivalent to a strongly
coupled \(\mathrm{E}_8 \times \mathrm{E}_8\) Heterotic string theory in
ten dimensions. However, the string-string duality conjecture in six
dimensions states that the strongly coupled limit of a Heterotic string
theory in six dimensions on a four-torus is equivalent to a weakly
coupled Type II string theory in six dimensions on \(\mathrm{K}3\). Similarly, it
states that the strongly coupled limit of a Type II theory in six
dimensions on \(\mathrm{K}3\) is equivalent to a weakly coupled Heterotic string
theory in six-dimensions on a four-torus. Now, as a strongly coupled
Type IIA string theory in ten-dimensions is equivalent to the low energy
limit of M-Theory on a ``large'' \(S^1\), the low energy limit of
M-Theory on \(S^1 \times \mathrm{K}3\) should be equivalent to a weakly coupled
Heterotic string theory on a four-torus by way of six-dimensional
string-string duality. Similarly, as a strongly coupled
\(\mathrm{E}_8 \times \mathrm{E}_8\) Heterotic string theory in
ten-dimensions is equivalent to the low energy limit of M-Theory on a
``large'' \(S^1/\mathbb{Z}_2\), the low energy limit of M-Theory on
\(S^1/\mathbb{Z}_2 \times T^4\) should be equivalent to a weakly coupled
Type II string theory on \(\mathrm{K}3\). The first of the above two consistency
checks on the relation between M-Theory and String Theory will be the
subject of this article. However, we will comment on the second
consistency check in our conclusion."
\end{quote}

So, as you can see, there is a veritable jungle of relationships out
there. But you must be wondering by now: \emph{what's M-theory?}
According to

\begin{enumerate}
\def\labelenumi{\arabic{enumi})}
\setcounter{enumi}{1}
\tightlist
\item
  Edward Witten, ``Five-branes and M-theory on an orbifold'', available
  as
  \href{https://arxiv.org/abs/hep-th/9512219}{\texttt{hep-th/9512219}}.
\end{enumerate}
\noindent
the M stands for ``magic'', ``mystery'', or ``membrane'', according to
taste. From a mathematical viewpoint a better term might be ``murky'',
since apparently everything known about M-theory is indirect and
circumstantial, except for the classical limit, in which it seems to act
as a theory of \(2\)-branes and \(5\)-branes, where an ``\(p\)-brane'' is a
\(p\)-dimensional analog of a membrane or surface.

Well, here I must leave off, for reasons of ignorance. I don't really
understand the evidence for the existence of the M-theory... I can
only await the day when the murk clears and it becomes possible to learn
about this stuff a bit more easily. It has been suggested that string
theory is a bit of 21st-century mathematics that accidentally fell into
the 20th century. I think this is right, and that eventually much of
this stuff will be seen as much simpler than it seems now.

Now let me briefly describe some papers I actually sort of understand.

\begin{enumerate}
\def\labelenumi{\arabic{enumi})}
\setcounter{enumi}{2}
\item
  Abhay Ashtekar, ``Polymer geometry at Planck scale and quantum
  Einstein equations'', available as
  \href{https://arxiv.org/abs/hep-th/9601054}{\texttt{hep-th/9601054}}.

  Roumen Borissov, Seth Major and Lee Smolin, ``The geometry of quantum
  spin networks'', available as
  \href{https://arxiv.org/abs/gr-qc/9512043}{\texttt{gr-qc/9512043}},

  Bernd Br\"ugmann, ``On the constraint algebra of quantum gravity in the
  loop representation'', available as
  \href{https://arxiv.org/abs/gr-qc/9512036}{\texttt{gr-qc/9512036}}.

  Kiyoshi Ezawa, ``Nonperturbative solutions for canonical quantum
  gravity: an overview'', available as
  \href{https://arxiv.org/abs/gr-qc/9601050}{\texttt{gr-qc/9601050}}

  Kiyoshi Ezawa, ``A semiclassical interpretation of the topological
  solutions for canonical quantum gravity'', available as
  \href{https://arxiv.org/abs/gr-qc/9512017}{\texttt{gr-qc/9512017}}.

  Jorge Griego, ``Extended knots and the space of states of quantum
  gravity'', available as
  \href{https://arxiv.org/abs/gr-qc/9601007}{\texttt{gr-qc/9601007}}.

  Seth Major and Lee Smolin, ``Quantum deformation of quantum gravity'',
  available as
  \href{https://arxiv.org/abs/gr-qc/9512020}{\texttt{gr-qc/9512020}}.
\end{enumerate}

Work on the loop representation of quantum gravity proceeds apace. The
paper by Ashtekar and the first one by Ezawa review various recent
developments and might be good to look at if one is just getting
interested in this subject. Smolin has been pushing the idea of
combining ideas about the quantum group \(SU_q(2)\) with the loop
representation, and his papers with Borissov and Major are about that.
This seems rather interesting but still a bit mysterious to me. I
suspect that what it amounts to is thinking of loops as excitations not
of the Ashtekar--Lewandowksi vacuum state but the Chern--Simons state. I'd
love to see this clarified, since these two states are two very
important exact solutions of quantum gravity, and the latter has the
former as a limit as the cosmological constant goes to zero. In the
second paper listed, Ezawa gives semiclassical interpretations of these
and other exact solutions of quantum gravity.

\begin{enumerate}
\def\labelenumi{\arabic{enumi})}
\setcounter{enumi}{3}
\tightlist
\item
  Thomas Kerler, ``Genealogy of nonperturbative quantum-invariants of
  3-Manifolds: the surgical family'', available as
  \href{https://arxiv.org/abs/q-alg/9601021}{\texttt{q-alg/9601021}}.
\end{enumerate}

Kerler brings a bit more order to the study of quantum invariants of
3-manifolds, in particular, the Reshetikhin--Turaev,
Hennings-Kauffman-Radford, and Lyubashenko invariants. All of these are
constructed using certain braided monoidal categories, like the category
of (nice) representations of a quantum group. He describes how
Lyubashenko's invariant specializes to the Reshetikhin--Turaev invariant
for semisimple categories and to the Hennings-Kauffman-Radford invariant
for Tannakian categories. People interested in extended TQFTs and
\(2\)-categories will find his work especially interesting, because he
works with these invariants using these techniques. James Dolan and I
have argued that it's only this way that one will really understand
these TQFTs (see \protect\hyperlink{week49}{``Week 49''}).

In future editions of This Week's Finds I will say more about
\(n\)-categories and topological quantum field theory. I have a feeling
that while I've discussed these a lot, I have never really explained the
basic ideas very well. As I gradually understand the basic ideas better,
they seem simpler and simpler to me, so I think I should try to explain
them.

\hypertarget{week73}{%
\section{February 24, 1996}\label{week73}}

In this and future issues of This Week's Finds, I'd like to talk a bit
more about higher-dimensional algebra, and how it should lead to many
exciting developments in mathematics and physics in the 21st century.
I've talked quite a bit about this already, but I hear from some people
that the ``big picture'' remained rather obscure. The main reason, I
suppose, is that I was just barely beginning to see the big picture
myself! As Louis Crane noted, in this subject it often feels that we are
unearthing the fossilized remains of some enormous prehistoric beast,
still unsure of its extent or how it all fits together. Of course that's
what makes it so exciting, but I'll try to make sense what we've found
so far, and where it may lead. In the Weeks to come, I'll start out
describing some basic stuff, and work my way up to some very new ideas.

However, before I get into that, I'd like to say a bit about something
completely different: biology.

\begin{enumerate}
\def\labelenumi{\arabic{enumi})}
\item
  \emph{Biological Asymmetry and Handedness}, Ciba Foundation Symposium
  \textbf{162}, John Wiley and Sons, 1991.

  D. K. Kondepudi and D. K. Nelson, ``Weak neutral currents and the
  origins of molecular chirality'', \emph{Nature} \textbf{314},
  pp.~438--441.
\end{enumerate}

It's always puzzled me how humans and other animals could be
consistently asymmetric. A 50-50 mix of two mirror-image forms could
easily be explained by ``spontaneously broken symmetry'', but in fact
there are many instances of populations with a uniform handedness. Many
examples appear in Weyl's book ``Symmetry'' (see
\protect\hyperlink{week63}{``Week 63''}). To take an example close to
home, the human brain appears to be lateralized in a fairly consistent
manner; for example, most people have the speech functions concentrated
in the left hemisphere of their cerebrum --- even most, though not all,
left-handers.

One might find this unsurprising: it just means that the asymmetry is
encoded in the genes. But think about it: how are the genes supposed to
tell the embryo to develop in an asymmetric way? How do they explain the
difference between right and left? That's what intrigues me.

Of course, genes code for proteins, and most proteins are themselves
asymmetric. Presumably the answer lurks somewhere around here. Indeed,
even the amino acids of which the proteins are composed are asymmetric,
as are many sugars and for that matter, the DNA itself, which is
composed of two spirals, each of which has an intrinsic directionality
and hence a handedness. The handedness of many of these basic
biomolecules is uniform for all life on the globe, as far as I know.

In the conference proceedings on biological asymmetry, there is an
interesting article on the development of asymmetry in \emph{C.
elegans}. Ever since the 1960s, this little nematode has been a favorite
among biologists because of its simplicity, and because of the
advantages understanding one organism thoroughly rather than many
organisms in a sketchy way. I'm sure most of you know about the fondness
geneticists have for the fruit fly, but Caenorhabditis elegans is a far
simpler critter: it only has 959 cells, all of which have been
individually named and studied! There are over 1000 people studying it
by now, there is a journal devoted to it --- The Worm Breeder's Gazette
--- and it has its own world-wide web server. Moreover, folks are busily
sequencing not only the complete human genome but also all 100 million
bases of the DNA of \emph{C. elegans}.

But I digress! The point here is that \emph{C. elegans} is asymmetric,
and exhibits a consistent handedness. And the cool thing is that in the
conference proceedings, Wood and Kevshan report on experiments where
they artificially changed the handedness of \emph{C. elegans} embryos
when they consisted of only 6 cells! The embryos look symmetric when
they have 4 cells; by the time they have 8 cells the asymmetry is
marked. By moving some cells around at the 6-cell stage, Wood and
Kevshan were able to create fully functional \emph{C. elegans} having
opposite the usual handedness.

The question of exactly how the embryo's asymmetry originates from some
asymmetry at the molecular still seems shrouded in mystery. And there is
another puzzle: how did the biomolecules choose their handedness in the
first place? Here spontaneous symmetry breaking --- an essentially
random choice later amplified by selection --- seems a natural
hypothesis. But physicists should be interested to note that another
alternative has been seriously proposed. Weak interactions violate
parity and thus endow the laws of nature with an intrinsic handedness.
This means there is a slight difference in energies between any
biomolecule and its enantiomer, or mirror-image version. According to S.
F. Mason's article in the conference proceedings, this difference indeed
favors the observed forms of amino acids and sugars --- the left-handed
or ``L'' amino acids and the right-handed or ``D'' sugars. But the
difference is is incredibly puny --- typically it amounts to
\(10^{-14}\) joules per mole! How could such a small difference matter?
Well, Kondepudi and Nelson have done calculations suggesting that in
certain situations where there is both autocatalysis of both L and D
forms of these molecules, and also competition between them, random
fluctuations can be averaged out, while small energy level differences
can make a big difference.

That would be rather satisfying to me: knowing that my heart is where it
is for the same reason that neutrinos are left-handed. But in fact this
theory is very controversial\ldots. I mention it only because of its
charm.

\hypertarget{week73_tale}{If we think of the universe as passing through the course of history
from simplicity to complexity, from neutrinos to nematodes to humans,
it's natural to wonder what's at the bottom, where things get very
simple, where physics blurs into pure logic\ldots. far from the ``spires
of form''. Ironically, even the simplest things may be hard to
understand, because they are so abstract.   But let's give it a try ---
here comes the Tale of \(n\)-Categories.}

Let's begin with the world of sets. In a certain sense, there is nothing
much to a set except its cardinality, the number of elements it has. Of
course, set theorists work hard to build up the universe of sets from
the empty set, each set being a set of sets, with its own distinctive
personality:
\[\{\} ,\, \{\{\}\} ,\, \{\{\{\}\}\} ,\, \{\{\},\{\{\}\}\} ,\, \{\{\},\{\{\{\}\}\}\} ,\, \{\{\},\{\{\}\},\{\{\},\{\{\}\}\}\}\]
and the like. But for many purposes, a one-to-one and onto function
between two sets allows us to treat them as the same. So if necessary,
we could actually get by with just one set of each cardinality. For
example
\[\{\} ,\, \{\{\}\} ,\, \{\{\},\{\{\}\}\} ,\, \{\{\},\{\{\}\},\{\{\},\{\{\}\}\}\}\]
and so on. For short, people like to call these \[0 ,\, 1 ,\, 2 ,\, 3\]
and so on. We could wonder what comes after all these finite cardinals,
and what comes after that, and so on, but let's not. Instead, let's
ponder what we've done so far. We started with the universe of sets ---
not exactly the set of all sets, but pretty close --- but very soon we
started playing with functions between sets. This is what allowed us to
speak of two sets with the same cardinality as being isomorphic.

In short, we are really working with the \emph{category} of sets. A
category is something just as abstract as a set, but a bit more
structured. It's not a mere collection of objects; there are also
morphisms between objects, in this case the functions between sets.

Some of you might not know the precise definition of a category; let me
state it just for completeness. A category consists of a collection of
``objects'' and a collection of ``morphisms''. Every morphism \(f\) has
a ``source'' object and a ``target'' object. If the source of \(f\) is
\(X\) and its target is \(Y\), we write \(f\colon X \to Y\). In
addition, we have:

\begin{enumerate}
\def\labelenumi{\arabic{enumi})}
\item
  Given a morphism \(f\colon X \to Y\) and a morphism
  \(g\colon Y \to Z\), there is a morphism \(fg\colon X \to Z\), which
  we call the ``composite'' of \(f\) and \(g\).
\item
  Composition is associative: \((fg)h = f(gh)\).
\item
  For each object \(X\) there is a morphism \(1_X\colon X \to X\),
  called the ``identity'' of \(X\). For any \(f\colon X \to Y\) we have
  \(1_X f = f 1_Y = f\).
\end{enumerate}

That's it.

(Note that we are writing the composite of \(f\colon X \to Y\) and
\(g\colon Y \to Z\) as \(fg\), which is backwards from the usual order.
This will make life easier in the long run, though, since \(fg\) will
mean ``first do \(f\), then \(g\)''.)

Now, there are lots of things one can do with sets, which lead to all
sorts of interesting examples of categories, but in a sense the
primordial category is \(\mathsf{Set}\), the category of sets and
functions. (One might try to make this precise, by trying to prove that
every category is a subcategory of \(\mathsf{Set}\), or something like
that. Actually the right way to say how \(\mathsf{Set}\) is primordial
is called the ``Yoneda lemma''. But to understand this lemma, one needs
to understand categories a little bit.)

When we get to thinking about categories a lot, it's natural to think
about the ``category of all categories''. Now just as it's a bit bad to
speak of the set of all sets, it's bad to speak of the category of all
categories. This is true, not only because Russell's paradox tends to
ruin attempts at a consistent theory of the ``thing of all things'', but
because, just as what really counts is the \emph{category} of all sets,
what really counts is the \emph{\(2\)-category} of all categories.

To understand this, note that there is a very sensible notion of a
morphism between categories. It's called a ``functor'', and a functor
\(F\colon \mathcal{C} \to \mathcal{D}\) from a category \(\mathcal{C}\)
to a category \(\mathcal{D}\) is just something that assigns to each
object \(x\) of \(\mathcal{C}\) an object \(F(x)\) of \(\mathcal{D}\),
and to each morphism \(f\) of \(\mathcal{C}\) a morphism \(F(f)\) of
\(\mathcal{D}\), in such a way that ``all structure in sight is
preserved''. More precisely, we want:

\begin{enumerate}
\def\labelenumi{\arabic{enumi})}
\item
  If \(f\colon x \to y\), then \(F(f)\colon F(x) \to F(y)\).
\item
  If \(fg = h\), then \(F(f)F(g) = F(h)\).
\item
  If \(1_x\) is the identity morphism of \(x\), then \(F(1_x)\) is the
  identity morphism of \(F(x)\).
\end{enumerate}

It's good to think of a category as a bunch of dots --- objects --- and
arrows going between them --- morphisms. I would draw one for you if I
could here. Category theorists love drawing these pictures. In these
terms, we can think of the functor
\(F\colon \mathcal{C} \to \mathcal{D}\) as putting a little picture of
the category \(\mathcal{C}\) inside the category \(\mathcal{D}\). Each
dot of \(\mathcal{C}\) gets drawn as a particular dot in
\(\mathcal{D}\), and each arrow in \(\mathcal{C}\) gets drawn as a
particular arrow in \(\mathcal{D}\). (Two dots or arrows in
\(\mathcal{C}\) can get drawn as the same dot or arrow in
\(\mathcal{D}\), though.)

In addition, however, there is a very sensible notion of a
``2-morphism'', that is, a morphism between morphisms between
categories! It's called a ``natural transformation''. The idea is this.
Suppose we have two functors \(F\colon \mathcal{C} \to \mathcal{D}\) and
\(G\colon \mathcal{C} \to \mathcal{D}\). We can think of these as giving
two pictures of \(\mathcal{C}\) inside \(\mathcal{D}\). So for example,
if we have any object \(x\) in \(\mathcal{C}\), we get two objects in
\(\mathcal{D}\), \(F(x)\) and \(G(x)\). A ``natural transformation'' is
then a gadget that draws an arrow from each dot like \(F(x)\) to the dot
like \(G(x)\). In other words, for each \(x\), the natural
transformation \(T\) gives a morphism \(T_x\colon F(x) \to G(x)\). But
we want a kind of compatibility to occur: if we have a morphism
\(f\colon x \to y\) in \(\mathcal{C}\), we want \[
  \begin{tikzcd}
    F(x) \rar["F(f)"] \dar[swap,"T_x"]
    & F(y) \dar["T_y"]
  \\G(x) \rar[swap,"G(f)"]
    & G(y)
  \end{tikzcd}
\] to commute; in other words, we want \(T_x G(f) = F(f) T_y\).

This must seem very boring to the people who understand it and very
mystifying to those who don't. I'll need to explain it more later. For
now, let me just say a bit about what's going on. Sets are
``zero-dimensional'' in that they only consist of objects, or ``dots''.
There is no way to ``go from one dot to another'' within a set.
Nonetheless, we can go from one set to another using a function. So the
category of all sets is ``one-dimensional'': it has both objects or
``dots'' and morphisms or ``arrows between dots''. In general,
categories are ``one-dimensional'' in this sense. But this in turn makes
the collection of all categories into a ``two-dimensional'' structure, a
\(2\)-category having objects, morphisms between objects, and
\(2\)-morphisms between morphisms.

This process never stops. The collection of all \(n\)-categories is an
\((n+1)\)-category, a thing with objects, morphisms, \(2\)-morphisms,
and so on all the way up to \(n\)-morphisms. To study sets carefully we
need categories, to study categories well we need \(2\)-categories, to
study \(2\)-categories well we need \(3\)-categories, and so on...
so ``higher- dimensional algebra'', as this subject is called, is
automatically generated in a recursive process starting with a careful
study of set theory.

If you want to show off, you can call the \(2\)-category of all
categories \(\mathsf{Cat}\), and more generally, you can call the
\((n+1)\)-category of all \(n\)-categories \(n\mathsf{Cat}\).
\(n\mathsf{Cat}\) is the primordial example of an \((n+1)\)-category!

Now, just as you might wonder what comes after \(0,1,2,3,\ldots\), you
might wonder what comes after all these \(n\)-categories. The answer is
``\(\omega\)-categories''.

What comes after these? Well, let us leave that for another time. I'd
rather conclude by mentioning the part that's the most fascinating to me
as a mathematical physicist. Namely, the various dimensions of category
turn out to correspond in a very beautiful --- but still incompletely
understood --- way to the various dimensions of spacetime. In other
words, the study of physics in imaginary \(2\)-dimensional spacetimes
uses lots of \(2\)-categories, the study of physics in a 3d spacetimes
uses 3-categories, the study of physics in 4d spacetimes appears to use
4-categories, and so on. It's very surprising at first that something so
simple and abstract as the process of starting with sets and recursively
being led to study the \((n+1)\)-category of all \(n\)-categories could
be related to the dimensionality of spacetime. In particular, what could
possibly be special about 4 dimensions?

Well, it turns out that there \emph{are} some special things about 4
dimensions. But more on that later.  To continue reading the ``Tale of \(n\)-Categories'', see \hyperlink{week74_tale}{``Week 74''}.

\begin{center}\rule{0.5\linewidth}{0.5pt}\end{center}

\textbf{Addendum}: Long after writing the above, I just saw an
interesting article on chirality in biology:

\begin{enumerate}
\def\labelenumi{\arabic{enumi})}
\setcounter{enumi}{1}
\tightlist
\item
  N. Hirokawa, Y. Tanaka, Y. Okada and S. Takeda, ``Nodal flow and the
  generation of left-right asymmetry'', \emph{Cell} \textbf{125} 1
  (2006), 33--45.
\end{enumerate}
\noindent
It reports on detailed studies of how left-right asymmetry first shows
in the development of animal embryos. It turns out this asymmetry is
linked to certain genes with names like \emph{Lefty-1}, \emph{Lefty-2},
\emph{Nodal} and \emph{Pitx2}. About half of the people with a genetic
disorder called Kartagener's Syndrome have their organs in the reversed
orientation. These people also have immotile sperm and defective cilia
in their airway. This suggests that the genes controlling left-right
asymmetry also affect the development of cilia! And the link has
recently been understood....

The first visible sign of left-right asymmetry in mammal embryos is the
formation of a structure called the ``ventral node'' after the
front-back (dorsal-ventral) and top-bottom (anterior-posterior)
symmetries have been broken. This node is a small bump on the front of
the embryo.

It has recently been found that cilia on this bump wiggle in a way that
makes the fluid the embryo is floating in flow towards the \emph{left}.
It seems to be this leftward flow that generates many of the more fancy
left-right asymmetries that come later.

How do these cilia generate a leftward flow? It seems they spin around
\emph{clockwise}, and are tilted in such a way that they make a leftward
swing when they are near the surface of the embryo, and a rightward
swing when they are far away. This manages to do the job... the
article discusses the hydrodynamics involved.

I guess now the question becomes: why do these cilia spin clockwise?

\hypertarget{week74}{%
\section{March 5, 1996}\label{week74}}

Before continuing my story about higher-dimensional algebra, let me say
a bit about gravity. Probably far fewer people study general relativity
than quantum mechanics, which is partially because quantum mechanics is
more practical, but also because general relativity is mathematically
more sophisticated. This is a pity, because general relativity is so
beautiful!

Recently, I have been spending time on \texttt{sci.physics} leading an
informal (nay, chaotic) ``general relativity tutorial''. The goal is to
explain the subject with a minimum of complicated equations, while still
getting to the mathematical heart of the subject. For example, what does
Einstein's equation REALLY MEAN? It's been a lot of fun and I've learned
a lot! Now I've gathered up some of the posts and put them on a web
site:

\begin{enumerate}
\def\labelenumi{\arabic{enumi})}
\tightlist
\item
  John Baez \emph{et al}, ``General relativity tutorial'',
  \href{http://math.ucr.edu/home/baez/gr/gr.html}{\texttt{gr/gr.html}}
\end{enumerate}
\noindent
I hope to improve this as time goes by, but it should already be fun to
look at.

Let me also list a couple new papers on the loop representation of
quantum gravity, dealing with ways to make volume and area into
observables in quantum gravity:

\begin{enumerate}
\def\labelenumi{\arabic{enumi})}
\setcounter{enumi}{1}
\item
  Abhay Ashtekar and Jerzy Lewandowski, ``Quantum theory of geometry I:
  area operators'', to appear in Classical and
  Quantum Gravity, available as
  \href{https://arxiv.org/abs/gr-qc/9602046}{\texttt{gr-qc/9602046}}.

  Jerzy Lewandowski, ``Volume and quantizations'', available as
  \href{https://arxiv.org/abs/gr-qc/9602035}{\texttt{gr-qc/9602035}}.

  Roberto De Pietri and Carlo Rovelli, ``Geometry eigenvalues and scalar
  product from recoupling theory in loop quantum gravity'', 
  available as
  \href{https://arxiv.org/abs/gr-qc/9602023}{\texttt{gr-qc/9602023}}.
\end{enumerate}

I won't say anything about these now, but see
\protect\hyperlink{week55}{``Week 55''} for some information on area
operators.

\begin{center}\rule{0.5\linewidth}{0.5pt}\end{center}

\hypertarget{week74_tale}{Okay, where were we?}
 We had started messing around with sets, and we
noted that sets and functions between sets form a category, called Set.
Then we started messing around with categories, and we noted that not
only are there ``functors'' between categories, there are things that
ply their trade between functors, called ``natural transformations''. I
then said that categories, functors, and natural transformations form a
\(2\)-category. I didn't really say what a \(2\)-category is, except to
say that it has objects, morphisms between objects, and \(2\)-morphisms
between morphisms. Finally, I said that this pattern continues:
\(n\mathsf{Cat}\) forms an \((n+1)\)-category.

By the way, I said last time that \(\mathsf{Set}\) was ``the primordial
category''. Keith Ramsay reminded me by email that this can be
misleading. There are other categories that act a whole lot like
\(\mathsf{Set}\) and can serve equally well as ``the primordial
category''. These are called topoi. Poetically speaking, we can think of
these as alternate universes in which to do mathematics. For more on
topoi, see \protect\hyperlink{week68}{``Week 68''}. All I meant by
saying that \(\mathsf{Set}\) was ``the primordial category'' is that, if
we start from \(\mathsf{Set}\) and various categories of structures
built using sets --- groups, rings, vector spaces, topological spaces,
manifolds, and so on --- we can then abstract the notion of
``category'', and thus obtain \(\mathsf{Cat}\). In the same sense,
\(\mathsf{Cat}\) is the primordial \(2\)-category, and so on.

I mention this because it is part of a very important broad pattern in
higher-dimensional algebra. For example, we will see that the complex
numbers are the primordial Hilbert space, and that the category of
Hilbert spaces is the primordial ``2-Hilbert space'', and that the
\(2\)-category of 2-Hilbert spaces is the primordial ``3-Hilbert
space'', and so on. This leads to a quantum-theoretic analog of the
hierarchy of \(n\)-categories, which plays an important role in
mathematical physics. But I'm getting ahead of myself!

Let's start by considering a few examples of categories. I want to pick
some examples that will lead us naturally to the main themes of
higher-dimensional algebra. Beware: it will take us a while to get
rolling. For a while --- maybe a few issues of This Week's Finds ---
everything may seem somewhat dry, pointless and abstract, except for
those of you who are already clued in. It has the flavor of
``foundations of mathematics,'' but eventually we'll see these new
foundations reveal topology, representation theory, logic, and quantum
theory to be much more tightly interknit than we might have thought. So
hang in there.

For starters, let's keep the idea of ``symmetry'' in mind. The typical
way to think about symmetry is with the concept of a ``group''. But to
get a concept of symmetry that's really up to the demands put on it by
modern mathematics and physics, we need --- at the very least --- to
work with a \emph{category} of symmetries, rather than a group of
symmetries.

To see this, first ask: what is a category with one object? It is a
``monoid''. The \emph{usual} definition of a monoid is this: a set \(M\)
with an associative binary product and a unit element \(1\) such that
\(a1 = 1a = a\) for all \(a\) in \(M\). Monoids abound in mathematics;
they are in a sense the most primitive interesting algebraic structures.

To check that a category with one object is ``essentially just a
monoid'', note that if our category \(\mathcal{C}\) has one object
\(x\), the set \(\operatorname{Hom}(x,x)\) of all morphisms from \(x\)
to \(x\) is indeed a set with an associative binary product, namely
composition, and a unit element, namely \(1_x\). (Actually, in an
arbitrary category \(\operatorname{Hom}(x,y)\) could be a class rather
than a set. But let's not worry about such nuances.) Conversely, if you
hand me a monoid \(M\) in the traditional sense, I can easily cook up a
category with one object \(x\) and \(\operatorname{Hom}(x,x) = M\).

How about categories in which every morphism is invertible? We say a
morphism \(f\colon x\to y\) in a category has inverse \(g\colon y\to x\)
if \(fg = 1_x\) and \(gf = 1_y\). Well, a category in which every
morphism is invertible is called a ``groupoid''.

Finally, a group is a category with one object in which every morphism
is invertible. It's both a monoid and a groupoid!

When we use groups in physics to describe symmetry, we think of each
element \(g\) of the group \(G\) as a ``process''. The element \(1\)
corresponds to the ``process of doing nothing at all''. We can compose
processes \(g\) and \(h\) --- do \(h\) and then \(g\) --- and get the
product \(gh\). Crucially, every process \(g\) can be ``undone'' using
its inverse \(g^{-1}\).

We tend to think of this ability to ``undo'' any process as a key aspect
of symmetry. I.e., if we rotate a beer bottle, we can rotate it back so
it was just as it was before. We don't tend to think of \emph{smashing} the
beer bottle as a symmetry, because it can't be undone. But while
processes that can be undone are especially interesting, it's also nice
to consider other ones... so for a full understanding of symmetry
we should really study monoids as well as groups.

But we also should be interested in ``partially defined'' processes,
processes that can be done only if the initial conditions are right.
This is where categories come in! Suppose that we have a bunch of boxes,
and a bunch of processes we can do to a bottle in one box to turn it
into a bottle in another box: for example, ``take the bottle out of box
\(x\), rotate it 90 degrees clockwise, and put it in box \(y\)''. We can
then think of the boxes as objects and the processes as morphisms: a
process that turns a bottle in box \(x\) to a bottle in box \(y\) is a
morphism \(f\colon x\to y\). We can only do a morphism
\(f\colon x\to y\) to a bottle in box \(x\), not to a bottle in any
other box, so \(f\) is a ``partially defined'' process. This implies we
can only compose \(f\colon x\to y\) and \(g\colon u \to v\) to get
\(fg\colon x \to v\) if \(y = u\).

So: a monoid is like a group, but the ``symmetries'' no longer need be
invertible; a category is like a monoid, but the ``symmetries'' no
longer need to be composable!

Note for physicists: the operation of ``evolving initial data from one
spacelike slice to another'' is a good example of a ``partially
defined'' process: it only applies to initial data on that particular
spacelike slice. So dynamics in special relativity is most naturally
described using groupoids. Only after pretending that all the spacelike
slices are the same can we pretend we are using a group. It is very
common to pretend that groupoids are groups, since groups are more
familiar, but often insight is lost in the process. Also, one can only
pretend a groupoid is a group if all its objects are isomorphic.
Groupoids really are more general.

Physicists wanting to learn more about groupoids might try:

\begin{enumerate}
\def\labelenumi{\arabic{enumi})}
\setcounter{enumi}{2}
\tightlist
\item
  Alan Weinstein, ``Groupoids: unifying internal and external
  symmetry'', available at
   \href{http://www.ams.org/notices/199607/weinstein.pdf}{
  \texttt{http://www.ams.org/notices/199607/weinstein.pdf}}
\end{enumerate}

So: in contrast to a set, which consists of a static collection of
``things'', a category consists not only of objects or ``things'' but
also morphisms which can viewed as ``processes'' transforming one thing
into another. Similarly, in a \(2\)-category, the \(2\)-morphisms can be
regarded as ``processes between processes'', and so on. The eventual
goal of basing mathematics upon \(\omega\)-categories is thus to allow
us the freedom to think of any process as the sort of thing higher-level
processes can go between. By the way, it should also be very interesting
to consider ``\(\mathbb{Z}\)-categories'' (where \(\mathbb{Z}\) denotes
the integers), having \(j\)-morphisms not only for \(j = 0,1,2,\ldots\)
but also for negative \(j\). Then we may also think of any thing as a
kind of process.

How do the above remarks about groups, monoids, groupoids and categories
generalize to the \(n\)-categorical context? Well, all we did was start
with the notion of category and consider two sorts of requirement: that
the category have just one object, or that all morphisms be invertible.

A category with just one object --- a monoid --- could also be seen as a
set with extra algebraic structure, namely a product and unit. Suppose
we look at an \(n\)-category with just one object? Well, it's very
similar: then we get a special sort of \((n-1)\)-category, one with a
product and unit! We call this a ``monoidal \((n-1)\)-category''. I will
explain this more thoroughly later, but let me just note that we can
keep playing this game, and consider a monoidal \((n-1)\)-category with
just one object, which is a special sort of \((n-2)\)-category, which we
could call a ``doubly monoidal \((n-2)\)-category'', and so on. This
game must seem very abstract and mysterious when one first hears of it.
But it turns out to yield a remarkable set of concepts, some already
very familiar in mathematics, and it turns out to greatly deepen our
notion of ``commutativity''. For now, let me simply display a chart of
``\(k\)-tuply monoidal \(n\)-categories'' for certain low values of
\(n\) and \(k\):

\begin{longtable}[]{@{}llll@{}}
\caption{\(k\)-tuply monoidal \(n\)-categories}\tabularnewline
\toprule
\begin{minipage}[b]{0.26\columnwidth}\raggedright
\strut
\end{minipage} & \begin{minipage}[b]{0.21\columnwidth}\raggedright
\(n=0\)\strut
\end{minipage} & \begin{minipage}[b]{0.21\columnwidth}\raggedright
\(n=1\)\strut
\end{minipage} & \begin{minipage}[b]{0.21\columnwidth}\raggedright
\(n=2\)\strut
\end{minipage}\tabularnewline
\midrule
\endfirsthead
\toprule
\begin{minipage}[b]{0.26\columnwidth}\raggedright
\strut
\end{minipage} & \begin{minipage}[b]{0.21\columnwidth}\raggedright
\(n=0\)\strut
\end{minipage} & \begin{minipage}[b]{0.21\columnwidth}\raggedright
\(n=1\)\strut
\end{minipage} & \begin{minipage}[b]{0.21\columnwidth}\raggedright
\(n=2\)\strut
\end{minipage}\tabularnewline
\midrule
\endhead
\begin{minipage}[t]{0.26\columnwidth}\raggedright
\(k=0\)\strut
\end{minipage} & \begin{minipage}[t]{0.21\columnwidth}\raggedright
sets\strut
\end{minipage} & \begin{minipage}[t]{0.21\columnwidth}\raggedright
categories\strut
\end{minipage} & \begin{minipage}[t]{0.21\columnwidth}\raggedright
\(2\)-categories\strut
\end{minipage}\tabularnewline
\begin{minipage}[t]{0.26\columnwidth}\raggedright
\strut
\end{minipage} & \begin{minipage}[t]{0.21\columnwidth}\raggedright
\strut
\end{minipage} & \begin{minipage}[t]{0.21\columnwidth}\raggedright
\strut
\end{minipage} & \begin{minipage}[t]{0.21\columnwidth}\raggedright
\strut
\end{minipage}\tabularnewline
\begin{minipage}[t]{0.26\columnwidth}\raggedright
\(k=1\)\strut
\end{minipage} & \begin{minipage}[t]{0.21\columnwidth}\raggedright
monoids\strut
\end{minipage} & \begin{minipage}[t]{0.21\columnwidth}\raggedright
monoidal categories\strut
\end{minipage} & \begin{minipage}[t]{0.21\columnwidth}\raggedright
monoidal \(2\)-categories\strut
\end{minipage}\tabularnewline
\begin{minipage}[t]{0.26\columnwidth}\raggedright
\strut
\end{minipage} & \begin{minipage}[t]{0.21\columnwidth}\raggedright
\strut
\end{minipage} & \begin{minipage}[t]{0.21\columnwidth}\raggedright
\strut
\end{minipage} & \begin{minipage}[t]{0.21\columnwidth}\raggedright
\strut
\end{minipage}\tabularnewline
\begin{minipage}[t]{0.26\columnwidth}\raggedright
\(k=2\)\strut
\end{minipage} & \begin{minipage}[t]{0.21\columnwidth}\raggedright
commutative monoids\strut
\end{minipage} & \begin{minipage}[t]{0.21\columnwidth}\raggedright
braided monoidal categories\strut
\end{minipage} & \begin{minipage}[t]{0.21\columnwidth}\raggedright
braided monoidal \(2\)-categories\strut
\end{minipage}\tabularnewline
\begin{minipage}[t]{0.26\columnwidth}\raggedright
\strut
\end{minipage} & \begin{minipage}[t]{0.21\columnwidth}\raggedright
\strut
\end{minipage} & \begin{minipage}[t]{0.21\columnwidth}\raggedright
\strut
\end{minipage} & \begin{minipage}[t]{0.21\columnwidth}\raggedright
\strut
\end{minipage}\tabularnewline
\begin{minipage}[t]{0.26\columnwidth}\raggedright
\(k=3\)\strut
\end{minipage} & \begin{minipage}[t]{0.21\columnwidth}\raggedright
" "\strut
\end{minipage} & \begin{minipage}[t]{0.21\columnwidth}\raggedright
symmetric monoidal categories\strut
\end{minipage} & \begin{minipage}[t]{0.21\columnwidth}\raggedright
weakly involutory monoidal \(2\)-categories\strut
\end{minipage}\tabularnewline
\begin{minipage}[t]{0.26\columnwidth}\raggedright
\strut
\end{minipage} & \begin{minipage}[t]{0.21\columnwidth}\raggedright
\strut
\end{minipage} & \begin{minipage}[t]{0.21\columnwidth}\raggedright
\strut
\end{minipage} & \begin{minipage}[t]{0.21\columnwidth}\raggedright
\strut
\end{minipage}\tabularnewline
\begin{minipage}[t]{0.26\columnwidth}\raggedright
\(k=4\)\strut
\end{minipage} & \begin{minipage}[t]{0.21\columnwidth}\raggedright
" "\strut
\end{minipage} & \begin{minipage}[t]{0.21\columnwidth}\raggedright
" "\strut
\end{minipage} & \begin{minipage}[t]{0.21\columnwidth}\raggedright
strongly involutory monoidal \(2\)-categories\strut
\end{minipage}\tabularnewline
\begin{minipage}[t]{0.26\columnwidth}\raggedright
\strut
\end{minipage} & \begin{minipage}[t]{0.21\columnwidth}\raggedright
\strut
\end{minipage} & \begin{minipage}[t]{0.21\columnwidth}\raggedright
\strut
\end{minipage} & \begin{minipage}[t]{0.21\columnwidth}\raggedright
\strut
\end{minipage}\tabularnewline
\begin{minipage}[t]{0.26\columnwidth}\raggedright
\(k=5\)\strut
\end{minipage} & \begin{minipage}[t]{0.21\columnwidth}\raggedright
" "\strut
\end{minipage} & \begin{minipage}[t]{0.21\columnwidth}\raggedright
" "\strut
\end{minipage} & \begin{minipage}[t]{0.21\columnwidth}\raggedright
" "\strut
\end{minipage}\tabularnewline
\bottomrule
\end{longtable}

The quotes indicate that each column ``stabilizes'' past a certain
point. If you can't wait to read more about this, you might try
\protect\hyperlink{week49}{``Week 49''} for more, but I will explain it
all in more detail in future issues.

What if we take an \(n\)-category and demand that all \(j\)-morphisms
(\(j > 0\)) be invertible? Well, then we get something we could call an
``\(n\)-groupoid''. However, there are some important subtle issues
about the precise sense in which we might want all \(j\)-morphisms to be
invertible. I will have to explain that, too.

Let me conclude, though, by mentioning something the experts should
enjoy. If we define \(n\)-groupoids correctly, and then figure out how
to define \(\omega\)-groupoids correctly, the homotopy category of
\(\omega\)-groupoids turns out to be equivalent to the homotopy category
of topological spaces. The latter category is something algebraic
topologists have spent decades studying. This is one of the main ways
\(n\)-categories are important in topology. Using this correspondence
between \(n\)-groupoid theory and homotopy theory, the ``stabilization''
property described above is then related to a subject called ``stable
homotopy theory'', and ``\(\mathbb{Z}\)-groupoids'' are a way of talking
about ``spectra'' --- another important tool in homotopy theory.

The above paragraph is overly erudite and obscure, so let me explain the
gist: there is a way to think of a topological space as giving us an
\(\omega\)-groupoid, and the \(\omega\)-groupoid then captures all the
information about its topology that homotopy theorists find interesting.
(I will explain in more detail how this works later.) If this is
\emph{all} \(n\)-category theory did, it would simply be an interesting
language for doing topology. But as we shall see, it does a lot more.
One reason is that, not only can we use \(n\)-categories to think about
spaces, we can also use them to think about symmetries, as described
above. Of course, physicists are very interested in space and also
symmetry. So from the viewpoint of a mathematical physicist, one
interesting thing about \(n\)-categories is that they \emph{unify} the
study of space (or spacetime) with the study of symmetry.

I will continue along these lines next time and try to fill in some of
the big gaping holes.

To continue reading the ``Tale of \(n\)-Categories'', see
\protect\hyperlink{week75_tale}{``Week 75''}.

\hypertarget{week75}{%
\section{March 6, 1996}\label{week75}}

\hypertarget{week75_tale}{If you've been following my recent introduction to \(n\)-categories,
you'll note that I haven't actually given the definition of
\(n\)-categories! I've only defined categories, and \emph{hinted} at the
definition of \(2\)-categories by giving an example, namely
\(\mathsf{Cat}\).} This is the \(2\)-category whose objects are
categories, whose morphisms are functors, and whose \(2\)-morphisms are
natural transformations.

The definition of \(n\)-categories --- or maybe I should say the problem
of defining \(n\)-categories --- is actually surprisingly subtle. Since
I want to proceed at a gentle pace here, I think I should first get
everyone used to the \(2\)-category \(\mathsf{Cat}\), then define
\(2\)-categories in general, then play around with those a bit, and then
move on to \(n\)-categories for higher \(n\).

So recall what the objects, morphisms and \(2\)-morphisms in
\(\mathsf{Cat}\) are: they are categories, functors and natural
transformations. A functor \(F\colon \mathcal{C}\to\mathcal{D}\) assigns
to each object \(x\) of \(\mathcal{C}\) an object \(F(x)\) of
\(\mathcal{D}\), and to each morphism \(f\) of \(\mathcal{C}\) a
morphism \(F(f)\) of \(\mathcal{D}\), and has

\begin{enumerate}
\def\labelenumi{\arabic{enumi}.}
\tightlist
\item
  If \(f\colon x \to y\), then \(F(f)\colon F(x) \to F(y)\).
\item
  If \(fg = h\), then \(F(f)F(g) = F(h)\).
\item
  If \(1_x\) is the identity morphism of \(x\), then \(F(1_x)\) is the
  identity morphism of \(F(x)\).
\end{enumerate}

Given two functors \(F\colon\mathcal{C}\to\mathcal{D}\) and
\(G\colon\mathcal{C}\to\mathcal{D}\), a ``natural transformation''
\(T\colon F\Rightarrow G\) assigns to each object \(x\) of \(\mathcal{C}\) a
morphism \(T_x\colon F(x)\to G(x)\), such that for any morphism
\(f\colon x\to y\) in \(C\), the diagram 
\[
  \begin{tikzcd}
    F(x) \rar["F(f)"] \dar[swap,"T_x"]
    & F(y) \dar["T_y"]
  \\G(x) \rar[swap,"G(f)"]
    & G(y)
  \end{tikzcd}
\] commutes.

Let me give a nice example. Let \(\mathsf{Top}\) be the category with
topological spaces and continuous functions between them as morphisms.
Let \(\mathsf{Gpd}\) be the category with groupoids as objects and
functors between them as morphisms. (Remember from
\protect\hyperlink{week74_tale}{``Week 74''} that a groupoid is a category
with all morphisms invertible.) Then there is a functor
\[\Pi_1\colon\mathsf{Top}\to\mathsf{Gpd}\] called the ``fundamental
groupoid'' functor, which plays a very basic role in algebraic topology.

Here's how we get from any space \(X\) its ``fundamental groupoid''
\(\Pi_1(X)\). (If perchance you already know about the ``fundamental
group'', fear not, what we're doing now is very similar.) To say what
the groupoid \(\Pi_1(X)\) is, we need to say what its objects and
morphisms are. The objects in \(\Pi_1(X)\) are just the \emph{points} of \(X\)
and the morphisms are just certain equivalence classes of \emph{paths} in \(X\).
More precisely, a morphism \(f\colon x\to y\) in \(\Pi_1(X)\) is just an
equivalence class of continuous paths from \(x\) to \(y\), where two
paths from \(x\) to \(y\) are decreed equivalent if one can be
continuously deformed to the other while not moving the endpoints. (If
this equivalence relation holds we say the two paths are ``homotopic'',
and we call the equivalence classes ``homotopy classes of paths.'')

This is a truly wonderful thing! Recall the idea behind categories. A
morphism \(f\colon x\to y\) is supposed to be some abstract sort of
``process going from \(x\) to \(y\).'' The human mind often works by
visual metaphors, and our visual image of a ``process'' from \(x\) to
\(y\) is some sort of ``arrow'' from \(x\) to \(y\):
\[x\xrightarrow{f}y.\] That's why we write \(f\colon x\to y\), of
course. But now what we are doing is taking this visual metaphor
literally! We have a space \(X\), like the piece of the computer screen
on which you are actually reading this text. The objects in \(\Pi_1(X)\)
are then points in \(X\), and a morphism is basically just a path from
\(x\) to \(y\): \[x\xrightarrow{f}y.\] Well, not quite; it's a homotopy
class of paths. But still, something very interesting is going on here:
we are turning something ``concrete'', namely the space \(X\), into a
corresponding ``abstract'' thing, namely the groupoid \(\Pi_1(X)\), by
thinking of ``concrete'' paths as ``abstract'' morphisms. Here I am
thinking of geometrical concepts as ``concrete'', and algebraic ones as
``abstract''.

You may wonder why we use homotopy classes of paths, rather than paths.
One reason is that topologists want to use \(\Pi_1\) to distill a very
abstract ``essence'' of the topological space \(X\), an essence that
conveys only information that's invariant under ``homotopy
equivalence''. Another reason is that we can easily get homotopy classes
of paths to compose associatively, as they must if they are to be
morphisms in a category. We simply glom them end to end:
\[x\xrightarrow{f}y\xrightarrow{g}z\] and there is no problem with
associativity, since we can easily reparametrize the paths to make
\((fg)h = f(gh)\). (If you haven't thought about this, please do!) If we
do not work with homotopy classes, it's a pain to define composition in
such a way that \((fg)h = f(gh)\). Unless you are sneaky, you only get
that \((fg)h\) is \emph{homotopic} to \(f(gh)\); there are ways to get
composition to come out associative, but they are all somewhat immoral.
As we shall see, if we want to preserve the ``concreteness'' of \(X\) as
much as possible, and work with morphisms that are actual paths in \(X\)
rather than homotopy equivalence classes, the best thing is to work with
\(n\)-categories. More on that later.

Let's see; I claimed there is a functor
\(\Pi_1\colon\mathsf{Top}\to\mathsf{Gpd}\), but so far I have only
defined \(\Pi_1\) on the objects of \(\mathsf{Top}\); we also need to
define it for morphisms. That's easy. A morphism \(F\colon X\to Y\) in
\(\mathsf{Top}\) is a continuous map from the space \(X\) to the space
\(Y\); this is just what we need to take points in \(X\) to points in
\(Y\), and homotopy classes of paths in \(X\) to homotopy classes of
paths in \(Y\). So it gives us a morphism in \(\mathsf{Gpd}\) from the
fundamental groupoid \(\Pi_1(X)\) to the fundamental groupoid
\(\Pi_1(Y)\). There are various things to check here, but it's not hard.
We call this morphism \(\Pi_1(F)\colon\Pi_1(X)\to\Pi_1(Y)\). With a
little extra work, we can check that
\(\Pi_1\colon\mathsf{Top}\to\mathsf{Gpd}\), defined this way, is really
a functor.

Perhaps you are finding this confusing. If so, it could easily be
because there are several levels of categories and such going on here.
For example, \(\Pi_1(X)\) is a groupoid, and thus a category, so there
are morphisms like \(f\colon x\to y\) in it; but on the other hand
\(\mathsf{Gpd}\) itself is a category, and there are morphisms like
\(\Pi_1(F)\colon\Pi_1(X)\to\Pi_1(Y)\) in it, which are functors! If you
find this confusing, take heart. Getting confused this way is crucial to
learning \(n\)-category theory! After all, \(n\)-category theory is all
about how every ``process'' is also a ``thing'' which can undergo
higher-level ``processes''. Complex, interesting structures emerge from
very simple ones by the interplay of these different levels. It takes
work to mentally hop up and down these levels, and to weather the
inevitable ``level slips'' one makes when one screws up. If you expect
it to be easy and are annoyed when you mess up, you will hate this
subject. When approached in the right spirit, it is very fun; it teaches
one a special sort of agility. (We are, by the way, nowhere near the
really tricky stuff yet.)

Okay, so we have seen an interesting example of a functor
\[\Pi_1\colon\mathsf{Top}\to\mathsf{Gpd}.\] As I said, we can think of
this as going from the concrete world of spaces to the abstract world of
groupoids, turning concrete paths into abstract ``morphisms''.
Wonderfully, there is a kind of ``reverse'' functor,
\[K\colon\mathsf{Gpd}\to\mathsf{Top}\] 
sending any groupoid \(G\) into a topological space \(KG\) called the 
``classifying space'' of \(G\).   This turns the abstract into the
concrete, by making abstract morphisms into concrete paths! Basically,
it goes like this. Say we have a groupoid \(G\). We will build the space
\(K(G)\) out of simplices as follows. Start with one 0-simplex for each
object in \(G\). A 0-simplex is simply a point. We can draw the
0-simplex for the object \(x\) as follows: \[x\] Then put in one
\(1\)-simplex for each morphism in \(G\). A \(1\)-simplex is just a line
segment. We can draw the \(1\)-simplex for the morphism
\(f\colon x\to y\) as follows: \[
  \begin{tikzpicture}
    \node (x) at (0,0) {$x$};
    \node (y) at (1.5,0) {$y$};
    \draw[thick] (x) to node[fill=white]{$f$} (y);
  \end{tikzpicture}
\] Really we should draw an arrow going from left to right, but soon
things will get too messy if I do that, so I won't. Then, whenever we
have \(fg=h\) in the groupoid, we stick in a \(2\)-simplex. A
\(2\)-simplex is just a triangle and we visualize it as follows: \[
  \begin{tikzpicture}
    \node (x) at (0,0) {$x$};
    \node (y) at (1,1.7) {$y$};
    \node (z) at (2,0) {$z$};
    \draw[thick] (x) to node[fill=white]{$f$} (y);
    \draw[thick] (x) to node[fill=white]{$h$} (z);
    \draw[thick] (y) to node[fill=white]{$g$} (z);
    \node at (4,0.8) {$
      \begin{aligned}
        f&\colon x\to y
      \\g&\colon x\to z
      \\h&\colon y\to z
      \end{aligned}
    $};
  \end{tikzpicture}
\] Then, whenever we have \(fgh = k\) in the groupoid, we stick in a
\(3\)-simplex, which we can visualize as a tetrahedron like this \[
  \begin{tikzpicture}
    \node (w) at (0,0) {$x$};
    \node (x) at (1.5,2.6,0) {$x$};
    \node (y) at (1.5,1) {$y$};
    \node (z) at (3,0) {$z$};
    \draw[thick] (w) to node[fill=white]{$f$} (x);
    \draw[thick] (x) to node[fill=white]{$g$} (y);
    \draw[thick] (y) to node[fill=white]{$h$} (z);
    \draw[thick] (w) to node[fill=white]{$k$} (z);
    \draw[thick] (w) to node[fill=white]{$fg$} (y);
    \draw[thick] (x) to node[fill=white]{$gh$} (z);
    \node at (5.5,1.2) {$
      \begin{aligned}
        f&\colon w\to x
      \\g&\colon x\to y
      \\h&\colon y\to z
      \\k&\colon w\to z
      \end{aligned}
    $};
  \end{tikzpicture}
\] And so on... we do this forever and get a big ``simplicial
set,'' which we can think of as the topological space \(KG\). The
reader might want to compare \protect\hyperlink{week70}{``Week 70''},
where we do the same thing for a monoid instead of a groupoid. Really, one
can do it for any category.

That's how we define \(K\) on objects; it's not hard to define \(K\) on
morphisms too, so we get \[K\colon\mathsf{Gpd}\to\mathsf{Top}\] In case
you study this in more detail at some point, I should add that folks
often think of simplicial set as somewhat abstract combinatorial
objects in their own right, and then they break down \(K\) into two steps:
first they take the ``nerve'' of a groupoid and get a simplicial
set, and then they take the ``geometrical realization'' of the
simplicial set to get a topological space. For more on simplicial
sets and the like, try:

\begin{enumerate}
\def\labelenumi{\arabic{enumi})}
\tightlist
\item
  J. P. May, \emph{Simplicial Objects in Algebraic Topology}, Van
  Nostrand, Princeton, 1968.
\end{enumerate}

Now, in what sense is the functor \(K\colon\mathsf{Gpd}\to\mathsf{Top}\)
the ``reverse'' of the functor
\(\Pi_1\colon\mathsf{Top}\to\mathsf{Gpd}\)? Is it just the ``inverse''
in the traditional sense? No! It's something more subtle. As we shall
see, the fact that \(\mathsf{Cat}\) is a \(2\)-category means that a
functor can have a more subtle and interesting sorts of ``reverse'' than
one might expect if one were used to the simple ``inverse'' of a
function. This is something I alluded to in
\protect\hyperlink{week74_tale}{``Week 74''}: inverses become subtler as we
march up the \(n\)-categorical hierarchy.

I'll talk about this more later. But let me just drop a teaser...
Quantum mechanics is all about Hilbert spaces and linear operators
between them. Given any (bounded) linear operator \(F\colon H\to H'\)
from one Hilbert space to another, there is a subtle kind of ``reverse''
operator, called the ``adjoint'' of \(F\) and written
\(F^*\colon H'\to H\), defined by
\[\langle x,F^*y \rangle = \langle Fx,y \rangle\] for all \(x\) in \(H\)
and \(y\) in \(H'\). This is not the same as the ``inverse'' of \(F\);
indeed, it exists even if \(F\) is not invertible. This sort of
``reverse'' operator is deeply related to the ``reverse'' functors I am
hinting at above, and for this reason those ``reverse'' functors are
also called ``adjoints''. This is part of a profound and somewhat
mysterious relationship between quantum theory and category
theory... which I eventually need to describe.

To continue reading the ``Tale of \(n\)-Categories'', see
\protect\hyperlink{week76_tale}{``Week 76''}.

\hypertarget{week76}{%
\section{March 9, 1996}\label{week76}}

Yesterday I went to the oral exam of Hong Xiang, a student of Richard
Seto who is looking for evidence of quark-gluon plasma at Brookhaven.
The basic particles interacting via the strong force are quarks and
gluons; these have an associated kind of ``charge'' known as color.
Under normal conditions, quarks and gluons are confined to lie within
particles with zero total color, such as protons and neutrons, and more
generally the baryons and mesons seen in particle acccelerators --- and
possibly glueballs, as well. (See \protect\hyperlink{week68}{``Week
68''} for more on glueballs.)

However, the current theory of the strong force --- quantum
chromodynamics --- predicts that at sufficiently high densities and/or
pressures, a plasma of protons and neutrons should undergo a phase
transition called ``deconfinement'', past which the quarks and gluons
will roam freely. At low densities, this is expected to happen at a
temperature corresponding to about 200 MeV per nucleon (i.e., proton or
neutron). If my calculation is right, this is about 2 trillion Kelvin!
At low temperatures, it's expected to happen at about 5 to 20 times the
density of an atomic nucleus. (Normal nuclear matter has about 0.16
nucleons per femtometer cubed.) For more on this, check out these

\begin{enumerate}
\def\labelenumi{\arabic{enumi})}
\item
  Wikipedia, Quark-gluon plasma, \hfill \break 
  \href{https://en.wikipedia.org/wiki/Quark-gluon_plasma}{\texttt{https://en.wikipedia.org/wiki/Quark-gluon\(\underline{\;}\)plasma}}.
  
  Relativistic Heavy Ion Collider homepage,
  \href{http://www.bnl.gov/RHIC/}{\texttt{http://www.bnl.gov/RHIC/}}.
\end{enumerate}

The folks at Brookhaven are attempting to get high densities \emph{and}
temperatures by slamming two gold nuclei together. They are getting
densities of about 9 times that of a nucleus... and I forget what
sort of temperature, but there is reason to hope that the combined high
density and pressure might be enough to cause deconfinement and create a
``quark-gluon plasma''. Colliding gold on gold at high energies produces
a enormous spray of particles, but amidst this they are looking for a
particular signal of deconfinement. They are looking for
\(\varphi\)-mesons and looking to see if their lifetime is modified. A
\(\varphi\)-meson is a spin-\(1\) particle made of a strange quark /
strange antiquark pair; strange quarks and antiquarks are supposed to be
common in the quark-gluon plasma formed by the collision. Folks think
the lifetime of a \(\varphi\)-meson will be affected by the medium it
finds itself in, and that this should serve as a signature of
deconfinement. In fact, they may have already seen this!

People might also enjoy looking at this review article:

\begin{enumerate}
\def\labelenumi{\arabic{enumi})}
\setcounter{enumi}{1}
\tightlist
\item
  Adriano Di Giacomo, ``Mechanisms of colour confinement'',
  available as
  \href{https://arxiv.org/abs/hep-th/9603029}{\texttt{hep-th/9603029}}.
\end{enumerate}

\hypertarget{week76_tale}{
Okay, let me continue the ``Tale of \(n\)-Categories''.}  I want to lead up to
their role in physics, but to do it well, there are quite a few things I
need to explain first. One of the important things about \(n\)-category
theory is that they allow a much more fine-grained approach to the
notion of ``sameness'' than we would otherwise be able to achieve.

In a bare set, two elements \(x\) and \(y\) are either equal or not
equal; there is nothing much more to say.

In a category, two objects \(x\) and \(y\) can be equal or not equal,
but more interestingly, they can be \emph{isomorphic} or not, and if
they are, they can be isomorphic in many different ways. An isomorphism
between \(x\) and \(y\) is simply a morphism \(f\colon x\to y\) which
has an inverse \(g\colon y\to x\). (For a discussion of inverse
morphisms, see \protect\hyperlink{week74_tale}{``Week 74''}.)

For example, in the category Set an isomorphism is just a one-to-one and
onto function. If we know two sets \(x\) and \(y\) are isomorphic we
know that they are ``the same in a way'', even if they are not equal.
But specifying an isomorphism \(f\colon x\to y\) does more than say
\(x\) and \(y\) are the same in a way; it specifies a \emph{particular
way} to regard \(x\) and \(y\) as the same.

In short, while equality is a yes-or-no matter, a mere \emph{property},
an isomorphism is a \emph{structure}. It is quite typical, as we climb
the categorical ladder (here from elements of a set to objects of a
category) for properties to be reinterpreted as structures, or sometimes
vice-versa.

What about in a \(2\)-category? Here the notion of equality sprouts
still further nuances. Since I haven't defined \(2\)-categories in
general, let me work with an example, Cat. This has categories as its
objects, functors as its morphisms, and natural transformations as its
2-morphisms.

So... we can certainly speak, as before, of the \emph{equality} of
categories. We can also speak of the \emph{isomorphism} of categories:
an isomorphism between \(\mathcal{C}\) and \(\mathcal{D}\) is a functor
\(F\colon\mathcal{C}\to\mathcal{D}\) for which there is an inverse
functor \(G\colon\mathcal{D}\to\mathcal{C}\). I.e., \(FG\) is the
identity functor on \(\mathcal{C}\) and \(GF\) is the identity on
\(\mathcal{D}\), where we define the composition of functors in the
obvious way. But because we also have natural transformations, we can
also define a subtler notion, the \emph{equivalence} of categories. An
equivalence is a functor \(F\colon\mathcal{C}\to\mathcal{D}\) together
with a functor \(G\colon\mathcal{D}\to\mathcal{C}\) and natural
isomorphisms \(a\colon FG\to 1_C\) and \(b\colon GF \to 1_D\). A
``natural isomorphism'' is a natural transformation which has an
inverse.

Abstractly, I hope you can see the pattern here: just as we can
``relax'' the notion of equality to the notion of isomorphism when we
pass from sets to categories, we can relax the condition that \(FG\) and
\(GF\) equal identity functors to the condition that they be isomorphic
to identity functors when we pass from categories to the \(2\)-category
\(\mathsf{Cat}\). We need to have the natural transformations to be able
to speak of functors being isomorphic, just as we needed functions to be
able to speak of sets being isomorphic. In fact, with each extra level
in the theory of \(n\)-categories, we will be able to come up with a
still more refined notion of ``\(n\)-equivalence'' in this way. That's
what ``processes between processes between processes...'' allow us
to do.

But let me attempt to bring this notion of equivalence of categories
down to earth with some examples. Consider first a little category
\(\mathcal{C}\) with only one object \(x\) and one morphism, the
identity morphism \(1_x\colon x\to x\). We can draw \(\mathcal{C}\) as
follows: \[x\] where we don't bother drawing the identity morphism
\(1_x\). This category, by the way, is called the ``terminal category''.
Next consider a little category \(\mathcal{D}\) with two objects \(y\)
and \(z\) and only four morphisms: the identities \(1_y\) and \(1_z\),
and two morphisms \(f\colon y\to z\) and \(g\colon z\to y\) which are
inverse to each other. We can draw \(\mathcal{D}\) as follows: \[
  \begin{tikzcd}
    y \rar[bend right=40,swap,"f"] & z \lar[bend right=40,swap,"g"]
  \end{tikzcd}
\] where again we don't draw identities.

So: \(\mathcal{C}\) is a little world with only one object, while D is a
little world with only two isomorphic objects... that are
isomorphic in precisely one way! \(\mathcal{C}\) and \(\mathcal{D}\) are
clearly not isomorphic, because for a functor
\(F\colon\mathcal{C}\to\mathcal{D}\) to be invertible it would need to
be one-to-one and onto on objects, and also on morphisms.

However, \(\mathcal{C}\) and \(\mathcal{D}\) are equivalent. For
example, we can let \(F\colon\mathcal{C}\to\mathcal{D}\) be the unique
functor with \(F(x) = y\), and let \(G\colon\mathcal{D}\to\mathcal{C}\)
be the unique functor from \(\mathcal{D}\) to \(\mathcal{C}\). (There is
only one functor from any category to \(\mathcal{C}\), since
\(\mathcal{C}\) has only one object and one morphism; this is why we
call \(\mathcal{C}\) the terminal category.) Now, if we look at the
functor \(FG\colon\mathcal{C}\to\mathcal{C}\), it's not hard to see that
this is the identity functor on \(\mathcal{C}\). But the composite
\(GF\colon\mathcal{D}\to\mathcal{D}\) is not the identity functor on
\(\mathcal{D}\). Instead, it sends both \(y\) and \(z\) to \(y\), and
sends all the morphisms in \(\mathcal{D}\) to \(1_y\). But while not
\emph{equal} to the identity functor on \(\mathcal{D}\), the functor
\(GF\) is \emph{naturally isomorphic} to it. We can define a natural
transformation \(b\colon GF\to 1_D\) by setting \(b_y = 1_y\) and
\(b_z = f\). Here some folks may want to refresh themselves on the
definition of natural transformation, given in
\protect\hyperlink{week75_tale}{``Week 75''}, and check that \(b\) is really
one of these, and that \(b\) is a natural isomorphism because it has an
inverse.

The point is, basically, that having two uniquely isomorphic things with
no morphisms other than the isomorphisms between them and the identity
morphisms isn't really all that different from having one thing with
only the identity morphism. Category theorists generally regard
equivalent categories as ``the same for all practical purposes''. For
example, given any category we can find an equivalent category in which
any two isomorphic objects are equal. We call these ``skeletal''
categories because all the fat is gone and just the essential bones are
left. For example, the category \(\mathsf{FinSet}\) of finite sets, with
functions between them as morphisms, is equivalent to the category with
just the sets 
\[
  \begin{aligned}
    0 &= \{\}
  \\1 &= \{0\}
  \\2 &= \{0,1\}
  \\3 &= \{0,1,2\}
  \end{aligned}
\] 
etc., and functions between them as morphisms (see
\protect\hyperlink{week73_tale}{``Week 73''}). Essentially all the
mathematics that can be done in \(\mathsf{FinSet}\) can be done in this
skeletal category. This may seem shocking, but it's true\ldots.
Similarly, the category \(\mathsf{Set}\) is equivalent to the category
\(\mathsf{Card}\) having one set of each cardinality. Also, the category
\(\mathsf{Vect}\) of complex finite--dimensional vector spaces, with
linear functions between them as morphisms, is equivalent to a skeletal
category where the only objects are those of the form \(\mathbb{C}^n\).
\emph{This} example should not seem shocking; it's this fact which
allows unsophisticated people to do linear algebra under the impression
that all finite-dimensional vector spaces are of the form
\(\mathbb{C}^n\), and still manage to do all the practical computations
that more sophisticated people can do, who know the abstract definition
of vector space and thus know of lots more finite-dimensional vector
spaces.

However, there is another thing we can do in \(\mathsf{Cat}\), another
refinement of the notion of isomorphism, which I alluded to in
\protect\hyperlink{week75_tale}{``Week 75''}. This is the notion of ``adjoint
functor''. Let me mention a few examples (in addition to the example
given in \protect\hyperlink{week75_tale}{``Week 75''}) and let the reader
ponder them before giving the definition. Let \(\mathsf{Grp}\) denote
the category with groups as objects and homomorphisms as morphisms, a
homomorphism \(f\colon G\to H\) between groups being a function with
\(f(1) = 1\) and \(f(gh) = f(g)f(h)\) for all \(g, h\) in \(G\). Then
there is a nice functor \[L\colon\mathsf{Set}\to\mathsf{Grp}\] which
takes any set \(S\) to the free group on \(S\): this is the group
\(L(S)\) formed by all formal products of elements in \(S\) and inverses
thereof, with no relations other than those in the definition of a
group. For example, a typical element of the free group on \(\{x,y,z\}\)
is \(xyzy^{-1}xxy\).

(It's easy to see that \(f\colon S\to T\) is a function between sets,
there is a unique homomorphism \(L(f)\colon L(S)\to L(T)\) with
\(L(f)(x) = f(x)\) for all \(x\) in \(S\), and that this makes \(L\)
into a functor.)

There is also a nice functor \[R\colon\mathsf{Grp}\to\mathsf{Set}\]
taking any group to its underlying set, and taking any homomorphism to
its underlying function. We call this a ``forgetful'' functor since it
simply amounts to forgetting that we are working with groups and just
thinking of them as sets.

Now there is a sense in which \(L\) and \(R\) are reverse processes, but
it is delicate. They certainly aren't inverses, and they aren't even
part of an equivalence between \(\mathsf{Set}\) and \(\mathsf{Grp}\).
Nonetheless they are ``adjoints''. If the reader hasn't thought about
this, she may enjoy figuring out what this might mean... perhaps
keeping the adjoint operators mentioned in
\protect\hyperlink{week75_tale}{``Week 75''} in mind.

To continue reading the ``Tale of \(n\)-Categories'', see
\protect\hyperlink{week77_tale}{``Week 77''}

\hypertarget{week77}{%
\section{March 23, 1996}\label{week77}}

I spent last week at Penn State visiting the CGPG --- the Center for
Gravitational Physics and Geometry. I like to visit this place whenever
I can, because I've never found anywhere else that's as good for talking
about quantum gravity.

The CGPG is run by Abhay Ashtekar, who introduced the ``new variables''
for general relativity (see \protect\hyperlink{week7}{``Week 7''}). This
formulation of general relativity allowed Carlo Rovelli and Lee Smolin
to develop a new approach to quantum gravity, called the ``loop
representation''. Smolin is at the CGPG, while Rovelli teaches at
Pittsburgh, only a brief plane ride away: he was heading back just when
I showed up. Jorge Pullin, who has done a lot of work on knot theory and
quantum gravity, is also at the CGPG. Roger Penrose visits it regularly,
and happened to be there last week. There is always a peppy bunch of
grad students and postdocs wandering about the place, and some
interesting mathematicians across the street. I have a particular
interest in the work of Jean-Luc Brylinski, since he has thought a lot
about the relationships between conformal field theory and category
theory (see \protect\hyperlink{week25}{``Week 25''}).

You can find out more about the new variables at the
following web sites:

\begin{enumerate}
\def\labelenumi{\arabic{enumi})}
\item
  Center for Gravitational Physics and Geometry (CGPG) home page,
  \href{https://web.archive.org/web/20010302061610if_/http://vishnu.nirvana.phys.psu.edu:80/}{\texttt{https://web.}}  \href{https://web.archive.org/web/20010302061610if_/http://vishnu.nirvana.phys.psu.edu:80/}{\texttt{archive.org/web/20010302061610if\(\underline{\;}\)/}}\href{https://web.archive.org/web/20010302061610if_/http://vishnu.nirvana.phys.psu.edu:80/}{\texttt{http://vishnu.nirvana.phys.psu.edu:80/}} 

  Reading list on the new variables:
  \href{https://web.archive.org/web/20010210203335/http://vishnu.nirvana.phys.psu.edu/readinglist/readinglist.html}{\texttt{https://web.archive.org/web/20010210203335/}} \href{https://web.archive.org/web/20010210203335/http://vishnu.nirvana.phys.psu.edu/readinglist/readinglist.html}{\texttt{http://vishnu.nirvana.phys.psu.edu/readinglist/readinglist.html}}.
\end{enumerate}

I had two goals at the CGPG. One was to get people interested in the
uses of higher-dimensional algebra in physics, and the other was to find
out where folks were heading in quantum gravity. I made decent headway
on the first front, but let me talk about the second one.

In the last few years, Abhay Ashtekar has been working hard with a bunch
of collaborators on getting the loop representation set up on a
mathematically rigorous basis, and making good progress. There is a
natural order in which to set things up, and the next problem to deal
with is the so-called Hamiltonian constraint (see
\protect\hyperlink{week43}{``Week 43''}). I have always been very
worried about this, and I'm not alone, since this all the dynamics of
quantum gravity is in this operator. Ashtekar and Lewandowski have a
paper partially written in which they rigorously define an operator
along these lines, using earlier ideas of Rovelli and Smolin. I have
been hoping that this answer could be tested somehow... for
example, checking out its commutation relations with the other
constraints. It turns out that they have already done this to extent
that seems possible. So then the question is, what next? March on, or
continue trying to make sure the Hamiltonian constraint is right?

I should add that Pullin and Gambini have another proposal regarding the
Hamiltonian constraint:

\begin{enumerate}
\def\labelenumi{\arabic{enumi})}
\setcounter{enumi}{1}
\tightlist
\item
  Rodolfo Gambini and Jorge Pullin, ``The general solution of the
  quantum Einstein equations?'', \href{Class. Quant. Grav. } \textbf{13} (1996),
  L125--L128.  Also available as
  \href{https://arxiv.org/abs/gr-qc/9603019}{\texttt{gr-qc/9603019}}.
\end{enumerate}

This is not as fully worked out, but it has a certain mathematical charm
to it so far. Thus we may eventually be in a situation where there are
various competing quantizations of gravity using the loop
representation, differing mainly in their choice of Hamiltonian
constraint. This suggests that we need further tests for what counts as
the ``right'' Hamiltonian constraint.

When we spoke this time, Ashtekar was in favor of testing Hamiltonian
constraints by seeing whether they implied the ``Bekenstein bound''.
This bound says that the maximal entropy of a physical system is
proportional to its surface area when we take quantum gravity into
account. There are a number of heuristic derivations of this bound, so
lots of people hope it would follow from any good theory of quantum
gravity. Since the ``physical states'' of quantum gravity must be
annihilated by the Hamiltonian constraint, and the maximal entropy of a
system is just the logarithm of the number of physical states, the
Hamiltonian constraint must have some interesting properties to get the
Bekenstein bound to work out. So we can expect some work along these
lines in the near future.

I also talked to Lee Smolin. He has been very interested in the relation
between the loop representation and certain simplified versions of
quantum gravity called topological quantum field theories (TQFTs). He
has ideas on how to derive the Bekenstein bound using this relationship
--- see \protect\hyperlink{week56}{``Week 56''} and
\protect\hyperlink{week57}{``Week 57''} for a description.

The funny thing is, some of the mathematics connecting TQFTs to the loop
representation of quantum gravity also connects TQFTs to another
well-known approach to quantum gravity --- string theory! Smolin has
been boning up on string theory lately, in part by giving a course on
the subject, and presently he is eager to bring string theory and the
loop representation closer together. So we can also expect to see more
work on attempts to unify string and loops. (See
\protect\hyperlink{week18}{``Week 18''} for a bit more on strings and
loops.)

So I left feeling reinvigorated and eager to continue my own work on
higher-dimensional algebra and physics... which is what I have
talking about here ever since \protect\hyperlink{week73_tale}{``Week 73''}.
In fact, I have been engaging in a lengthy warmup, a minicourse in
category theory, with an eye to the basic themes of \(n\)-category
theory. That way, when I get around to the really cool stuff, everyone
out there will know what the heck I'm talking about. In theory, anyway.
You gotta work a bit to wrap your mind around these concepts!

\begin{center}\rule{0.5\linewidth}{0.5pt}\end{center}

\hypertarget{week77_tale}{So, let's recall where we are in our tale of \(n\)-categories.}
We were
studying increasingly subtle variations on the theme of identity and
difference. Given two categories \(\mathcal{C}\) and \(\mathcal{D}\), we
can ask if they are \emph{equal} or not. We can also discuss
\emph{isomorphisms} between \(\mathcal{C}\) and \(\mathcal{D}\). An
isomorphism is a functor \(F\colon\mathcal{C}\to\mathcal{D}\) having an
inverse: a functor \(G\colon\mathcal{D}\to\mathcal{C}\) such that \(FG\)
is equal to the identity functor on \(\mathcal{D}\) and \(GF\) is equal
to the identity on \(\mathcal{C}\).

We can also discuss \emph{equivalences} between \(\mathcal{C}\) and
\(\mathcal{D}\). An equivalence is a functor
\(F\colon\mathcal{C}\to\mathcal{D}\) together with a functor
\(G\colon\mathcal{D}\to\mathcal{C}\) such that \(FG\) is naturally
isomorphic to the identity functor on \(\mathcal{D}\), and \(GF\) is
naturally isomorphic to the identity functor on \(\mathcal{C}\).
Remember, two functors from one category to another are ``naturally
isomorphic'' if there is a natural transformation from the first to the
second, and that natural transformation has an inverse.

In math jargon we say it this way: two categories are equivalent if
there is a functor from one to the other which is invertible ``up to a
natural isomorphism''. The most useful notion of categories being 
``the same'' turns out to be not equality, or isomorphism, but this
 more supple notion of ``equivalence''!

(As we shall see later, this is because \(\mathsf{Cat}\) is a
\(2\)-category. Remember, an \(n\)-category is some sort of thing with
objects, morphisms, \(2\)-morphisms, and so on up to \(n\)-morphisms. One of
the of the main themes of \(n\)-category theory is that we may regard
two things are ``the same'', or ``equivalent'', if there is some sort of
process to get from one to the other, and this process is
invertible... up to equivalence! More precisely, we say an
\(n\)-morphism is an equivalence if it's invertible, and then we work
our way down, inductively defining a \((j-1)\)-morphism to be an
equivalence if it's invertible up to an equivalence. This downwards
induction leaves off when we define equivalence for ``\(0\)-morphisms'',
meaning objects.)

We have also begun talking about a curious situation where the
categories \(\mathcal{C}\) and \(\mathcal{D}\) are not at all ``the
same,'' but there are ``adjoint'' functors
\(L\colon\mathcal{C}\to\mathcal{D}\) and \(R\colon\mathcal{D}\to C\).
Let me list some examples before defining the concept of adjoint functor
and talking about it:

\begin{enumerate}
\def\labelenumi{\arabic{enumi}.}
\tightlist
\item
  First for the one we discussed in \protect\hyperlink{week76}{``Week
  76''}. Let \(\mathsf{Set}\) be the category of sets, and
  \(\mathsf{Grp}\) the category of groups. Let
  \(L\colon\mathsf{Set}\to\mathsf{Grp}\) be the functor taking each set
  \(S\) to the free group on \(S\), and doing the obvious thing to
  morphisms. Let \(R\colon\mathsf{Grp}\to\mathsf{Set}\) be the functor
  taking each group to its underlying set.
\item
  Let \(\mathsf{Ab}\) be the category of abelian (i.e., commutative)
  groups. Let \(L\colon\mathsf{Set}\to\mathsf{Ab}\) be the functor
  taking each set \(S\) to the free abelian group on \(S\). The ``free
  abelian group'' on \(S\) is what we get by taking the free group on
  \(S\) and imposing commutativity relations like \(xy = yx\) for all
  elements \(x,y\) in \(S\). Let \(R\colon\mathsf{Ab}\to\mathsf{Set}\)
  be the functor taking each abelian group to its underlying set.
\item
  Let \(L\colon\mathsf{Grp}\to\mathsf{Ab}\) be the functor taking each
  group \(G\) to its ``abelianization''. The abelianization of \(G\) is
  what we get when we impose the extra relations \(xy = yx\) for all
  elements \(x,y\) in \(G\). Let \(R\colon\mathsf{Ab}\to\mathsf{Grp}\)
  be the functor taking each abelian group to its underlying group.
\item
  Let \(\mathsf{Mon}\) be the category of monoids, where the objects are
  monoids and the morphisms are monoid homomorphisms. (Remember that a
  monoid is a set with an associative product and a unit; a monoid
  morphism \(f\colon M\to N\) is a function between monoids such that
  \(f(xy) = f(x)f(y)\) and \(f(1) = 1\).) Let
  \(L\colon\mathsf{Set}\to\mathsf{Mon}\) be the functor taking each set
  \(S\) to the free monoid on \(S\). This is simply the set of words
  whose letters are elements of \(S\), with the product given by
  concatenation of words, and the unit being the empty word. Let
  \(R\colon\mathsf{Mon}\to\mathsf{Set}\) be the functor taking each
  monoid to its underlying set.
\item
  Let \(L\colon\mathsf{Mon}\to\mathsf{Grp}\) be the functor taking each
  monoid \(M\) to the group obtained by throwing in formal inverses for
  every element of \(M\). The famous example of this is when
  \(\mathbb{N} = \{0,1,2,...\}\), which is a monoid whose ``product'' is
  addition. Then \(L(\mathbb{N}) = \mathbb{Z}\), the integers, since we
  have thrown in the negative integers. Let
  \(R\colon\mathsf{Grp}\to\mathsf{Mon}\) be the functor taking each
  group to its underlying monoid. I.e., \(R\) simply forgets that our
  group has inverses and thinks of it as a monoid.
\end{enumerate}

Note the common aspects of these examples! In most of them,
\(L\colon\mathcal{C} \to \mathcal{D}\) gives us a ``free'' object of
\(\mathcal{D}\) for every object of \(\mathcal{C}\), while
\(R\colon\mathcal{D}\to\mathcal{C}\) gives us an ``underlying'' object
of \(\mathcal{C}\) for every object of \(\mathcal{D}\). This is an
especially good way to think about it when the objects of
\(\mathcal{D}\) are objects of \(\mathcal{C}\) equipped with extra
structure, as in examples 1, 2, 4, and 5. For example, a group is a set
equipped with some extra structure, the group operations. In this case,
the functor \(L\colon\mathcal{C}\to\mathcal{D}\) turns an object of
\(\mathcal{C}\) into an object of \(\mathcal{D}\) by ``freely throwing
in whatever extra stuff is necessary, without putting in any relations
other than those needed to get an object of \(\mathcal{D}\)''.

It's not quite the same when the objects of \(\mathcal{D}\) are objects
of \(\mathcal{C}\) with extra \emph{properties}, as in example 3. In
this case, the functor \(L\colon\mathcal{C}\to\mathcal{D}\) forces an
object of \(\mathcal{C}\) to have the properties needed to be an object
of \(\mathcal{D}\). It does so in as nonviolent a manner as possible.

In either of these situations, \(R\colon\mathcal{D}\to\mathcal{C}\) has
the flavor of what we call a ``forgetful'' functor. This is not a
precisely defined term, but folks use it whenever we can simply
``forget'' something about an object of \(\mathcal{D}\) and think of it
as an object of \(\mathcal{C}\). For example, we can take a group, and
forget about the group operations, thinking of it as merely a set. Here
we are forgetting extra structure; we can also forget extra properties.

The crucial thing here is that unlike in an equivalence, there is a
built-in asymmetry here: \(L\) and \(R\) have very different flavors,
and serve different mathematical purposes. We call \(L\) the ``left
adjoint'' of \(R\), and we call \(R\) the ``right adjoint'' of \(L\).

There are situations where adjoint functors \(L\) and \(R\) aren't so
immediately reminiscent of the concepts ``free'' and ``underlying''. But
it's good to keep these ideas in mind when learning about adjoint
functors. I used to have trouble remembering which was supposed to be
the left adjoint and which was the right. The honest way to do this is
to remember the definition (coming up soon), but for a cheap mnemonic,
you can think of the L in a left adjoint as standing for ``liberty'' ---
that is, freedom!

So what's the definition of ``adjoint''? Roughly speaking, it's that for
any object \(c\) of \(\mathcal{C}\) and any object \(d\) of
\(\mathcal{D}\), we have
\[\operatorname{Hom}(Lc,d) = \operatorname{Hom}(c,Rd).\] Actually this
is a slight exaggeration: we don't want these to be equal. The guy on
the left is the set of morphisms from \(Lc\) to \(d\) in the category
\(\mathcal{D}\). The guy on the right is the set of morphisms from \(c\)
to \(Rd\) in the category \(\mathcal{C}\). So it's evil to want them to
be \emph{equal}. As you might guess, it's enough for them to be
naturally isomorphic in some sense. Let's not worry about that too much
yet, though. Let's get the basic idea here!

Consider example 1. Say \(S\) is a set and \(G\) is a group. Why is
\[\operatorname{Hom}(LS,G)\] naturally isomorphic to
\[\operatorname{Hom}(S,RG) \,\text{?}\]

In other words, why is the set of homomorphisms from the free group on
\(S\) to \(G\) naturally isomorphic to the set of functions from \(S\)
to the underlying set of \(G\)?

Well, say we have a homomorphism \(f\colon LS \to G\). Since \(LS\) is a
free group, we know \(f\) if we know what it does to each element of
\(S\)... and it can do whatever it wants to these elements! So we
can think of it as just a function from \(S\) to the underlying set of
\(G\). In other words, we can think of it as a function
\(f'\colon S \to RG\). Conversely, any function \(f'\colon S \to RG\)
gives us a homomorphism \(f\colon LS \to G\).

So this is the idea. Say we have an object \(c\) of \(\mathcal{C}\) and
an object \(d\) of \(\mathcal{D}\). Then:

\begin{quote}
``The set of morphisms from the free \(\mathcal{D}\)-object on \(c\) to
\(d\) is naturally isomorphic to the set of morphisms from \(c\) to the
underlying \(\mathcal{C}\)-object of \(d\).''
\end{quote}

Next time I will finish off the definition of adjoint functors, by
making this ``naturally isomorphic'' stuff precise. I will also begin to
explain what adjoint functors have to do with adjoint operators in
quantum mechanics. Remember that an ``observable'' in quantum theory is
an operator on a Hilbert space which is its own adjoint, while a
``symmetry'' in quantum theory is an operator whose adjoint is its
inverse. I eventually hope to show that this, and many other shocking
aspects of quantum theory, become less shocking when we think of the
world in terms of categories (or \(n\)-categories) rather than sets. The
way I think of it these days, the mysterious way quantum theory slammed
into physics in the early 20th century was just nature's way of telling
us we'd better learn \(n\)-category theory.

I'll also explain what adjoint functors have to do with the following
topological equations: \[
  \begin{tikzpicture}
    \begin{knot}
      \strand[thick] (0,0)
      to (0,1)
      to [out=up,in=up,looseness=2] (1,1)
      to [out=down,in=down,looseness=2] (2,1)
      to (2,2);
    \end{knot}
    \node at (3,1) {$=$};
    \begin{scope}[shift={(4,0)}]
      \begin{knot}
        \strand[thick] (0,0) to (0,2);
      \end{knot}
    \end{scope}
  \end{tikzpicture}
\] \[
  \begin{tikzpicture}
    \begin{scope}[xscale=-1,shift={(-2,0)}]
      \begin{knot}
        \strand[thick] (0,0)
        to (0,1)
        to [out=up,in=up,looseness=2] (1,1)
        to [out=down,in=down,looseness=2] (2,1)
        to (2,2);
      \end{knot}
    \end{scope}
    \node at (3,1) {$=$};
    \begin{scope}[shift={(4,0)}]
      \begin{knot}
        \strand[thick] (0,0) to (0,2);
      \end{knot}
    \end{scope}
  \end{tikzpicture}
\]

To continue reading the ``Tale of \(n\)-Categories'', see
\protect\hyperlink{week78_tale}{``Week 78''}.

\hypertarget{week78}{%
\section{March 28, 1996}\label{week78}}

\hypertarget{week78_tale}{
Last Week I began explaining the concept of ``adjoint functor''.} This
Week I want to finish explaining it --- or at least finish one round of
explanation! Then we'll begin to be able to see the unity of category
theory, topology, and quantum theory. These may seem rather distinct
subjects, but they have an interesting tendency to blur together when
one is doing topological quantum field theory. Part of the point of
higher-dimensional algebra is to explain this.

So, remember the idea of adjoint functors. Say we have categories
\(\mathcal{C}\) and \(\mathcal{D}\) and functors
\(L\colon\mathcal{C}\to\mathcal{D}\) and
\(R\colon\mathcal{D}\to\mathcal{C}\). Then we say \(L\) is the ``left
adjoint'' of \(R\), or that \(R\) is the ``right adjoint'' of \(L\), if
for any object \(c\) of \(\mathcal{C}\) and object \(d\) of
\(\mathcal{D}\), there is a natural one-to-one correspondence between
\(\operatorname{Hom}(Lc,d)\) and \(\operatorname{Hom}(c,Rd)\). An
example to keep in mind is when \(\mathcal{C}\) is the category of sets
and \(\mathcal{D}\) is the category of groups. Then \(L\) turns any set
into the free group on that set, while \(R\) turns any group into the
underlying set of that group. All sorts of other ``free'' and
``underlying'' constructions are also left and right adjoints,
respectively.

Now the only thing I didn't make very precise is what I mean by
``natural'' in the above paragraph. Informally, the idea of a
``natural'' one-to-one correspondence is that doesn't depend on any
arbitrary choices. The famous example is that if we have a
finite-dimensional vector space \(V\), it's always isomorphic to its
dual \(V^*\), but not naturally so: to set up an isomorphism we need to
pick a basis \(e_i\) of \(V\), and this gives a dual basis \(f^i\) of
\(V^*\), and then we get an isomorphism sending \(e_i\) to \(f^i\), but
this isomorphism depends on our choice of basis. But \(V\) is
\emph{naturally} isomorphic to its double dual \((V^*)^*\).

Now, it's hard to formalize the idea of ``not depending on any arbitrary
choices'' directly, so one needs to reflect on why it's bad for an
isomorphism to depend on arbitrary choices. The main reason is that the
arbitrariness may break a useful symmetry. In fact, Eilenberg and Mac
Lane invented category theory in order to formalize this idea of
``naturality as absence of symmetry-breaking''. Of course, they needed
the notion of category to get a sufficiently general concept of
``symmetry''. They realized that a nice way to turn things in the
category \(\mathcal{C}\) into things in the category \(\mathcal{D}\) is
typically a functor \(F\colon\mathcal{C}\to\mathcal{D}\). And then, if
we have two functors \(F,G\colon\mathcal{C}\to\mathcal{D}\), they
defined a ``natural transformation'' from \(F\) to \(G\) to be a bunch
of morphisms \[T_c\colon F(c)\to G(c),\] one for each object \(c\) of
\(\mathcal{C}\), such that the following diagram commutes for every
morphism \(f\colon c\to c'\) in \(\mathcal{C}\): \[
  \begin{tikzcd}
    F(c) \rar["F(f)"] \dar[swap,"T_c"]
    & F(c') \dar["T_{c'}"]
  \\G(c) \rar[swap,"G(f)"]
    & G(c')
  \end{tikzcd}
\] 
This condition says that the transformation \(T\) gets along with all
the ``symmetries'', or more precisely morphisms \(f\), in the category
\(\mathcal{D}\). We can do it either before or after applying one of
these symmetries, and we get the same result. For example, a vector
space construction which depends crucially on a choice of basis will
fail this condition if we take \(f\) to be a linear transformation
corresponding to a change of basis.

A ``natural isomorphism'' is then just a natural transformation that's
invertible, or in other words, one for which all the morphisms
T\textsubscript{c} are isomorphisms.

Okay. Hopefully that explains the idea of ``naturality'' a bit better.
But right now we are trying to figure out what we mean by saying that
\(\operatorname{Hom}(Lc,d)\) and \(\operatorname{Hom}(c,Rd)\) are
naturally isomorphic. To do this, we need to introduce a couple more
ideas: the product of categories, and the opposite of a category.

First, just as you can take the cartesian product of two sets, you can
take the cartesian product of two categories, say \(\mathcal{C}\) and
\(\mathcal{D}\). It's not hard. An object of
\(\mathcal{C}\times\mathcal{D}\) is just a pair of objects, one from
\(\mathcal{C}\) and one from \(\mathcal{D}\). A morphism in
\(\mathcal{C}\times\mathcal{D}\) is just a pair of morphisms, one from
\(\mathcal{C}\) and one from \(\mathcal{D}\). And everything works sort
of the way you'd expect.

Second, if you have a category \(\mathcal{C}\), you can form a new
category \(\mathcal{C}^\mathrm{op}\), the opposite of \(\mathcal{C}\),
which has the same objects as \(\mathcal{C}\), and has the arrows in
\(\mathcal{C}\) turned around backwards. In other words, a morphism
\(f\colon x\to y\) in \(\mathcal{C}^\mathrm{op}\) is defined to be a
morphism \(f\colon y\to x\) in \(\mathcal{C}\), and the composite \(fg\) of
morphisms in \(\mathcal{C}^\mathrm{op}\) is defined to be their
composite \(gf\) in \(\mathcal{C}\). So \(\mathcal{C}^\mathrm{op}\) is
like a through-the-looking-glass version of \(\mathcal{C}\) where they
do everything backwards. A functor
\(F\colon\mathcal{C}^\mathrm{op}\to\mathcal{D}\) is also called a
``contravariant'' functor from \(\mathcal{C}\) to \(\mathcal{D}\), since
we can think of it as a screwy functor from \(\mathcal{C}\) to
\(\mathcal{D}\) with \(F(fg) = F(g)F(f)\) instead of the usual
\(F(fg) = F(f)F(g)\). Whenever you see this perverse contravariant
behavior going on, you should suspect an opposite category is to blame.

Now, it turns out that we can think of the ``\(\operatorname{Hom}\)'' in
a category \(\mathcal{C}\) as a functor
\[\operatorname{Hom}(-,-)\colon\mathcal{C}^\mathrm{op}\times\mathcal{C}\to\mathsf{Set}\]
Here the \(-\)'s denote blanks to be filled in. Obviously, for any
object \((x,x')\) in \(\mathcal{C}^\mathrm{op}\times\mathcal{C}\), there
is a nice juicy set \(\operatorname{Hom}(x,x')\), the set of all
morphisms from \(x\) to \(x'\). But what if we have a morphism
\[(f,f')\colon(x,x')\to(y,y')\] in
\(\mathcal{C}^\mathrm{op}\times\mathcal{C}\)? For
\(\operatorname{Hom}(-,-)\) to be a functor, we should get a nice juicy
function
\[\operatorname{Hom}(f,f')\colon\operatorname{Hom}(x,x')\to\operatorname{Hom}(y,y').\]
How does this work? Well, remember that a morphism \((f,f')\) as above
is really just a pair consisting of a morphism \(f\colon x\to y\) in
\(\mathcal{C}^\mathrm{op}\) and a morphism \(f'\colon x'\to y'\) in
\(\mathcal{D}\). A morphism \(f\colon x\to y\) in
\(\mathcal{C}^\mathrm{op}\) is just a morphism \(f\colon y\to x\) in
\(\mathcal{D}\). Now say we have an unsuspecting element \(g\) of
\(\operatorname{Hom}(x,x')\) and we want to hit it with
\(\operatorname{Hom}(f,f')\) to get something in
\(\operatorname{Hom}(y,y')\). Here's how we do it:
\[\operatorname{Hom}(f,f')\colon g \mapsto f'gf\] We compose it with \(f'\)
on the left and f on the right! Composing on the left is a nice
covariant thing to do, but composing on the right is contravariant,
which is why we needed the opposite category
\(\mathcal{C}^\mathrm{op}\).

Okay, now back to our adjoint functors
\(L\colon\mathcal{C}\to\mathcal{D}\) and
\(R\colon\mathcal{D}\to\mathcal{C}\). Now we are ready to say what we
mean by \(\operatorname{Hom}(Lc,d)\) and \(\operatorname{Hom}(c,Rd)\)
being naturally isomorphic. Using the stuff we have set up, we can
define two functors
\[\operatorname{Hom}(L-,-)\colon\mathcal{C}^\mathrm{op}\times\mathcal{D}\to\mathsf{Set}\]
and
\[\operatorname{Hom}(-,R-)\colon\mathcal{C}^\mathrm{op}\times \mathcal{D} \to\mathsf{Set}\]
and we are simply saying that for \(L\) and \(R\) to be adjoints, we
demand the existence of a natural isomorphism between these functors!

Of course, this seems abstract, but if you work it out in some of the
examples of adjoint functors given in \protect\hyperlink{week76_tale}{``Week
76''} you'll see it all makes good sense.

Now let me start explaining what this all has to do with quantum theory.
(I'll put off the topology until next Week.) First of all, the
``\(\operatorname{Hom}\) functor'' we introduced,
\[\operatorname{Hom}(-,-)\colon\mathcal{C}^\mathrm{op}\times\mathcal{C}\to\mathsf{Set}\]
should remind you a whole lot of the inner product on a Hilbert space
\(H\). The inner product is linear in one slot and conjugate-linear in
the other, just like \(\operatorname{Hom}\) is covariant in one slot and
contravariant in the other. In fact, the inner product can be thought of
as a bilinear map
\[\langle -,- \rangle\colon H^\mathrm{op}\times H \to\mathbb{C}\] where
\(H^\mathrm{op}\), the ``opposite'' Hilbert space, is like \(H\) but
with a complex conjugate thrown into the definition of scalar
multiplication, and here \(\mathbb{C}\) denotes the complex numbers!

Second of all, the definition of adjoint functor, with
\(\operatorname{Hom}(Lc,d)\) and \(\operatorname{Hom}(c,Rd)\) being
naturally isomorphic, should remind you of adjoint linear operators on
Hilbert spaces. If we have a linear operator \(L\colon H\to K\) from a
Hilbert space \(H\) to a Hilbert space \(K\), its adjoint
\(R\colon K \to H\) is given by
\[\langle Lh,k \rangle = \langle h,Rk \rangle\] for all \(h\) in \(H\)
and \(k\) in \(K\).

In fact, the whole situation with adjoint functors is a kind of
``categorified'' version of the situation with adjoint linear operators.
Everything has been boosted up one notch on the \(n\)-categorical
ladder. What I mean is this: the Hilbert spaces \(H\) and \(K\) above
are \emph{sets}, with \emph{elements} \(h\) and \(k\), while the
categories \(\mathcal{C}\) and \(\mathcal{D}\) are \emph{categories},
with \emph{objects} \(c\) and \(d\). The inner product of two elements
of a Hilbert space is a \emph{number}, while the hom of two objects in a
category is a \emph{set}. Most interesting, the definition of adjoint
operators requires that \(\langle Lh,k \rangle\) and
\(\langle h,Rk \rangle\) be \emph{equal}, while the definition of
adjoint functors requires only that \(\langle Lc,d \rangle\) and
\(\langle c,Rd \rangle\) be \emph{naturally isomorphic}.

So we can think of adjoints in category theory as a boosted-up versionf
of the adjoints in quantum theory. But these days, I prefer to think of
the adjoints in quantum theory as a watered-down or ``decategorified''
version of the adjoints in category theory. The reason is that
categorification --- as noted by Louis Crane, who I believe invented the
term --- is a risky, hit-or-miss business, while decategorification is
much more systematic. Decategorification is simply the process of
neglecting the difference between isomorphism and equality. If we start
with an \(n\)-category and then get lazy and decide to think of
invertible \(n\)-morphisms as \emph{equations} between the
\((n-1)\)-morphisms, we get an \((n-1)\)-category. If we keep slacking
off like this, before you know it we're doing set theory! The final
stage of decategorification is when we get sloppy and instead of keeping
track of \emph{set}, we merely record the \emph{number} of its elements.

It's amusing to imagine this process of decategorification as one of
those elaborate Gnostic myths about the Fall. We start in the paradise
of \(\omega\)-categories (or perhaps even higher up), but by the
repeated sin of confusing equality with isomorphism we fall all the way
down the \(n\)-categorical ladder to the crude world of sets, or worse,
simply numbers. But all this happened a long time ago: now we need to
work our butt off to climb back up! In other words, historically our
early ancestors dealt with finite sets by replacing them with something
cruder: their numbers of elements. Counting is actually very handy, of
course, but it can only tell if the cardinalities of two sets are
\emph{equal}; it doesn't address the problem of specific
\emph{isomorphisms} between sets. To climb back up the \(n\)-categorical
ladder, we needed to start with the set \(\mathbb{N}\) of natural
numbers \[0, 1, 2, 3, \ldots\] and by dint of strenous mental effort
realize that this is just the decategorification of the category
\(\mathsf{FinSet}\) of finite sets. (In fact, category-theorists
routinely use \(2\) to stand for the 2-element set in the skeletal
category equivalent to \(\mathsf{FinSet}\), and so on --- see
\protect\hyperlink{week76_tale}{``Week 76''}.)

Now, you are certainly entitled to wonder if this elaborate
mathematical-theological fantasy is just a joke or if it has some
practical spinoffs. For example, is there anything we can \emph{do} with
the analogy between adjoint operators and adjoint functors? As it turns
out, there is. The point is that the analogy is not quite precise. For
example, every linear operator has an adjoint, but not every functor has
an adjoint --- nor need it be ``linear'' in any sense. If we endeavor to
make the analogy precise, we will invent a special sort of category
called a ``2-Hilbert space'' which is the precise categorified analog of
a Hilbert space. And we will invent a nice sort of ``linear'' functor
between these, and all such functors will have adjoints. Furthermore, in
this situation all left adjoints will also be right adjoints...
fixing another funny discrepancy. And these 2-Hilbert spaces turn out to
be closely related to \(2\)-dimensional topological quantum field
theories (in general, \(n\)-Hilbert spaces appear to be related to
\(n\)-dimensional TQFTs), as well as some interesting aspects of group
representation theory.

I'm busily writing a paper on exactly this stuff, but I have not
explained enough category theory here to describe it in detail yet. For
now, let me just make the connection between the
\(\operatorname{Hom}(-,-)\) of category theory and the
\(\langle -,-\rangle\) of quantum theory more clear, and hopefully more
plausible. If we have states \(h\) and \(h'\) in a Hilbert space,
\(\langle h,h'\rangle\) keeps track of the \emph{amplitude} of getting
from \(h\) to \(h'\). (Often people will say ``from \(h'\) to \(h\)'',
but here I think I really want to go the other way.) This is a mere
\emph{number}. If we have objects \(c\) and \(c'\) in a category,
\(\operatorname{Hom}(c,c')\) is the actual \emph{set} of ways to get
from \(c\) to \(c'\), that is, the set of morphisms from \(c\) to
\(c'\).

When one computes transition amplitudes by summing over paths, as in
Feynman path integrals, one is in a sense decategorifying, that is,
turning a set of ways to get from here to there into a number, the
transition amplitude. However, one is not just counting the ways, one is
counting them ``with phase''\ldots. and I must admit that the role of
the \emph{complex numbers} in quantum theory is still puzzling from this
viewpoint. For some food for thought, you might want to check out Dan
Freed's work on torsors, which are a categorified version of phases:

\begin{enumerate}
\def\labelenumi{\arabic{enumi})}
\tightlist
\item
  ``Higher algebraic structures and quantization'', by Daniel Freed,
  \emph{Commun. Math. Phys.} \textbf{159} (1994), 343--398, also
  available as
  \href{https://arxiv.org/abs/hep-th/9212115}{\texttt{hep-th/9212115}}.
\end{enumerate}

To continue reading the ``Tale of \(n\)-Categories'', see
\protect\hyperlink{week79_tale}{``Week 79''}.

\hypertarget{week79}{%
\section{April 1, 1996}\label{week79}}

Before I continue my tale of adjoint functors I want to say a little bit
about icosahedra, buckyballs, and last letter Galois wrote before his
famous duel\ldots. all of which is taken from the following marvelous
article:

\begin{enumerate}
\def\labelenumi{\arabic{enumi})}
\tightlist
\item
  Bertram Kostant, ``The graph of the truncated icosahedron and the last
  letter of Galois'', \emph{Notices Amer. Math. Soc.} \textbf{42} (September
  1995), 959--968. Also available at
  \href{http://www.ams.org/notices/199509/199509-toc.html}{\texttt{http://www.ams.org/notices/199509/199509-toc.html}}.
\end{enumerate}

When I was a graduate student at MIT I realized that Kostant (who
teaches there) was deeply in love with symmetry, and deeply
knowledgeable about some of its more mysterious byways. Unfortunately I
didn't dig too deeply into group theory at the time, and now I am
struggling to catch up.

Let's start with the Platonic solids. Note that the cube and the
octahedron are dual --- putting a vertex in the center of each of the
cube's faces gives you an octahedron, and vice versa. So every
rotational symmetry of the cube can be reinterpreted as a symmetry of
the octahedron, and vice versa. Similarly, the dodecahedron and the
icosahedron are dual, while the tetrahedron is self-dual. So while there
are 5 Platonic solids, there are really only 3 different symmetry groups
here.

These 3 ``Platonic groups'' are very interesting. The symmetry group of
the tetrahedron is the group \(A_4\) of all \emph{even} permutations of
4 things, since by rotating the tetrahedron we can achieve any even
permutation of its 4 vertices. The symmetry group of the cube is
\(S_4\), the group of \emph{all} permutations of 4 things. What are the
4 things here? Well, we can draw 4 line segments connecting opposite
vertices of the cube; these are the 4 things! The symmetry group of the
icosahedron is \(A_5\), the group of even permutations of 5 things. What
are the 5 things? It we take all the line segments connecting opposite
vertices we get 6 things, not 5, but we can't get all even permutations
of those by rotating the icosahedron. To find the \emph{5} things is a
bit trickier; I leave it as a puzzle here. See

\begin{enumerate}
\def\labelenumi{\arabic{enumi})}
\setcounter{enumi}{1}
\tightlist
\item
  John Baez, ``Some thoughts on the number 6'', available at \hfill \break
  \href{http://math.ucr.edu/home/baez/six.htm}{\texttt{http://math.ucr.edu/home/baez/six.html}}.
\end{enumerate}
\noindent
for an answer.

Once we convince ourselves that the rotational symmetry group of the
icosahedron is \(A_5\), it follows that it has \(5!/2 = 60\) elements.
But there is another nice way to see this. Take an icosahedron and chop
off all 12 corners, getting a truncated icosahedron with 12 regular
pentagonal faces and 20 regular hexagonal faces, with all edges the same
length. It looks just like a soccer ball. It's called an Archimedean
solid because, while not quite Platonic in its beauty, every face is a
regular polygon and every vertex looks alike: two pentagons abutting one
hexagon.

\[\href{http://en.wikipedia.org/wiki/Truncated_icosahedron}{\includegraphics[scale=0.6]{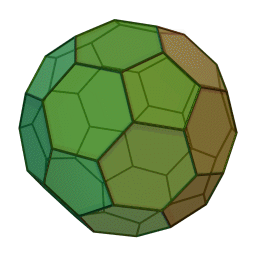}}\]

The truncated icosahedron has \(5 \times 12 = 60\) vertices. Every
symmetry of the icosahedron is a symmetry of the truncated icosahedron,
so \(A_5\) acts to permute these 60 vertices. Moreover, we can find an
element of \(A_5\) that moves a given vertex of the truncated
icosahedron to any other one, since ``every vertex looks alike''. Also,
there is a \emph{unique} element of \(A_5\) that does the job. So there
must be precisely as many elements of \(A_5\) as there are vertices of
the truncated icosahedron, namely 60.

There is a lot of interest in the truncated icosahedron recently,
because chemists had speculated for some time that carbon might form
\(C_{60}\) molecules with the atoms at the vertices of this solid, and a
while ago they found this was true. In fact, while \(C_{60}\) in this
shape took a bit of work to get ahold of at first, it turns out that
lowly soot contains lots of this stuff!
\[\href{http://en.wikipedia.org/wiki/Fullerene#Buckminsterfullerene}{\includegraphics[scale=0.8]{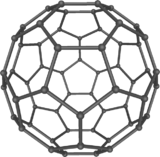}}\]
Since Buckminster Fuller was fond of using truncated icosahedra in his
geodesic domes, \(C_{60}\) and its relatives are called fullerenes, and
the shape is affectionately called a buckyball. For more about this
stuff, try:

\begin{enumerate}
\def\labelenumi{\arabic{enumi})}
\setcounter{enumi}{2}
\item
  P. W. Fowler and D. E. Manolpoulos, \emph{An Atlas of Fullerenes},
  Oxford University Press, Oxford, 1995.

  M. S. Dresselhaus, G. Dresselhaus, and P. C. Eklund, \emph{Science of
  Fullerenes and Carbon Nanotubules}, Academic Press, New York, 1994.

  G. Chung, B. Kostant and S. Sternberg, ``Groups and the buckyball'',
  in \emph{Lie Theory and Geometry}, eds.~J.-L. Brylinski, R. Brylinski,
  V. Guillemin and V. Kac, Birkh\"auser, Basel, 1994.
\end{enumerate}
\noindent
In fact, for the person who has everything, you can now buy 99.95\% pure
\(C_{60}\) online.

But I digress. Coming back to the 3 Platonic groups... there is
much more that's special about them. Most of it requires a little
knowledge of group theory to understand. For example, they are the 3
different finite subgroups of \(\mathrm{SO}(3)\) having irreducible
representations on \(\mathbb{R}^3\). And they are nice examples of
finite reflection groups. For more about them from this viewpoint, try
\protect\hyperlink{week62}{``Week 62''} and
\protect\hyperlink{week63}{``Week 63''}. Also, via the McKay
correspondence they correspond to the exceptional Lie groups
\(\mathrm{E}_6\), \(\mathrm{E}_7\), and \(\mathrm{E}_8\) --- see
\protect\hyperlink{week65}{``Week 65''} for an explanation of this!

Yet another interesting fact about these groups is buried in Galois'
last letter, written to the mathematician Chevalier on the night before
Galois' fatal duel. He was thinking about some groups we'd now call
\(\mathrm{PSL}(2,F)\). Here \(F\) is a field (for example, the real
numbers, the complex numbers, or \(\mathbb{Z}_p\), the integers
mod \(p\) where \(p\) is prime). \(\mathrm{PSL}(2,F)\) is a
``projective special linear group over \(F\).'' What does that mean?
Well, first of all, \(\mathrm{SL}(2,F)\) is the \(2\times2\) matrices
with entries in \(F\) having determinant equal to \(1\). These form a
group under good old matrix multiplication. The matrices in
\(\mathrm{SL}(2,F)\) that are scalar multiples of the identity matrix
form the ``center'' \(Z\) of \(\mathrm{SL}(2,F)\) --- the group of guys
who commute with everyone else. We can form the quotient group
\(\mathrm{SL}(2,F)/Z\), and get a new group called
\(\mathrm{PSL}(2,F)\).

Now, Galois was thinking about \(\mathrm{PSL}(2,\mathbb{Z}_p)\) where
\(p\) is prime. There's an obvious way to get this group to act as
permutations of \(p+1\) things. Here's how! For any field \(F\), the
group \(\mathrm{SL}(2,F)\) acts as linear transformations of the
\(2\)-dimensional vector space over \(F\), and it thus acts on the set
of lines through the origin in this vector space... which is called
the ``projective line'' over \(F\). But anything in \(\mathrm{SL}(2,F)\)
that's a scalar multiple of the identity doesn't move lines around, so
we can mod out by the center and think of the quotient group
\(\mathrm{PSL}(2,F)\) as acting on projective line. (By the way, this
explains the point of working with \(\mathrm{PSL}\) instead of plain old 
\(\mathrm{SL}\).)

Now, an element of the projective line is just a line through the origin
in \(F^2\). We can specify such a line by taking any nonzero vector
\((x,y)\) in \(F^2\) and drawing the line through the origin and this
vector. However, \((x',y')\) and \((x,y)\) determine the same line if
\((x',y')\) is a scalar multiple of \((x,y)\). Thus lines are in 1-1
correspondence with vectors of the form \((1,y)\) or \((x,1)\). When our
field \(F\) is \(\mathbb{Z}_p\), there are just \(p+1\) of these. So
\(\mathrm{PSL}(2,\mathbb{Z}_p)\) acts naturally on a set of \(p+1\)
things.

What Galois told Chevalier is that \(\mathrm{PSL}(2,\mathbb{Z}_p)\)
doesn't act nontrivially as permutation of any set with fewer than
\(p+1\) elements if \(p > 11\). This presumably means he knew that
\(\mathrm{PSL}(2,\mathbb{Z}_p)\) \emph{does} act nontrivially on a set
with only \(p\) elements if \(p = 5\), \(7\), or \(11\). For example,
\(\mathrm{PSL}(2,5)\) turns out to be isomorphic to \(A_5\), which acts
on a set of 5 elements in an obvious way. \(\mathrm{PSL}(2,7)\) and
\(\mathrm{PSL}(2,11)\) act on a 7-element set and an 11-element set,
respectively, in sneaky ways which Kostant describes.

These cases, \(p = 5\), \(7\) and \(11\), are the the only cases where
this happens and \(\mathrm{PSL}(2,\mathbb{Z}_p)\) is simple. (See
\protect\hyperlink{week66}{``Week 66''} if you don't know what
``simple'' means.) In each case it is very amusing to look at how
\(\mathrm{PSL}(2,\mathbb{Z}_p)\) acts nontrivially on a set with \(p\)
elements and consider the subgroup that doesn't move a particular
element of this set. For example, when \(p = 5\) we have
\(\mathrm{PSL}(2,5) = A_5\), and if we look at the subgroup of even
permutations of 5 things that leaves a particular thing alone, we get
\(A_4\). Kostant explains how if we play this game with
\(\mathrm{PSL}(2,7)\) we get \(S_4\), and if we play this game with
\(\mathrm{PSL}(2,11)\) we get \(A_5\). These are the 3 Platonic groups
again!!

But notice an extra curious coincidence. \(A_5\) is both
\(\mathrm{PSL}(2,5)\) and the subgroup of \(\mathrm{PSL}(2,11)\) that
fixes a point of an 11-element set. This gives a lot of relationships
between \(A_5\), \(\mathrm{PSL}(2,5)\), and \(\mathrm{PSL}(2,11)\). What
Kostant does is take this and milk it for all it's worth! In particular,
it turns out that one can think of \(A_5\) as the vertices of the
buckyball, and describe which vertices are connected by an edge using
the embedding of \(A_5\) in \(\mathrm{PSL}(2,11)\). I won't say how this
goes... read his paper!

This may even have some applications for fullerene spectroscopy, since
one can use symmetry to help understand spectra of compounds. (Indeed,
this is one way group theory entered chemistry in the first place.)

\begin{center}\rule{0.5\linewidth}{0.5pt}\end{center}

\hypertarget{week79_tale}{
Now let me return to talking about adjoint functors!} I have been stressing
the fact that two functors \(L\colon\mathcal{C}\to\mathcal{D}\) and
\(R\colon\mathcal{D}\to\mathcal{C}\) are adjoint if there is a natural
isomorphism between \(\operatorname{Hom}(Lc,d)\) and
\(\operatorname{Hom}(c,Rd)\). We can say that an ``adjunction'' is a
pair of functors \(L\colon\mathcal{C}\to\mathcal{D}\) and
\(R\colon\mathcal{D}\to\mathcal{C}\) together with a natural isomorphism
between \(\operatorname{Hom}(Lc,d)\) and \(\operatorname{Hom}(c,Rd)\).
But there is another way to think about adjunctions which is also good.

In \protect\hyperlink{week76_tale}{``Week 76''} we talked about an
``equivalence'' of categories. We can summarize it this way: an
``equivalence'' of the categories \(\mathcal{C}\) and \(\mathcal{D}\) is
a pair of functors \(F\colon\mathcal{C}\to\mathcal{D}\) and
\(G\colon\mathcal{D}\to\mathcal{C}\) together with natural
transformations \(e: FG \Rightarrow 1_\mathcal{D}\) and
\(i\colon 1_\mathcal{C} \Rightarrow GF\) that are themselves invertible.
(Note that we are now writing products of functors in the order that
ordinary mortals typically do, instead of the backwards way we
introduced in \protect\hyperlink{week73_tale}{``Week 73''}. Sorry! It just
happens to be better to write it this way now.) Now, the concept of
``adjunction'' is a cousin of the concept of ``equivalence'', and it's
nice to have a definition of adjunction that brings out this
relationship.

First, consider what happens in the definition of adjunction if we take
\(c = Rd\). Then we have a natural isomorphism between
\(\operatorname{Hom}(LRd,d)\) and \(\operatorname{Hom}(Rd,Rd)\). Now
there is a special element of \(\operatorname{Hom}(Rd,Rd)\), namely the
identity \(1_{Rd}\). This gives us a special element of
\(\operatorname{Hom}(LRd,d)\). Let's call this \[e_d\colon LRd \to d.\]
What is this morphism like in an example? Say
\(L\colon\mathsf{Set}\to\mathsf{Grp}\) takes each set to the free group
on that set, and \(R\colon\mathsf{Grp}\to\mathsf{Set}\) takes each group
to its underlying set. Then if \(d\) is a group, \(LRd\) is the free
group on the underlying set of \(d\). There's an obvious homomorphism
from \(LRd\) to \(d\), taking each word of elements in \(d\) and their
inverses to their product in \(d\). That's \(e_d\). It goes from the
free thing on the underlying thing of \(d\) to the thing \(d\) itself!

In fact, since everything in sight is natural, whenever we have an
adjunction the morphisms \(e_d\) define a natural transformation
\[e\colon LR \Rightarrow 1_\mathcal{D}\] Next, consider what happens in
the definition of adjunction if we take \(d = Lc\). Then we have a
natural isomorphism between \(\operatorname{Hom}(c,RLc)\) and
\(\operatorname{Hom}(Lc,Lc)\). Now there is a special element in
\(\operatorname{Hom}(Lc,Lc)\), namely the identity \(1_{Lc}\). This
gives us a special element in \(\operatorname{Hom}(c,RLc)\). Let's call
this \[i_c\colon c \to RLc.\] Again, it's good to consider the example
of sets and groups: if \(c\) is a set, \(RLc\) is the underlying set of
the free group on \(c\). There is an obvious way to map \(c\) into
\(RLc\). That's \(i_c\). It goes from the thing \(c\) to the underlying
thing of the free thing on \(c\).

As before, we get a natural transformation
\[i: 1_\mathcal{C} \Rightarrow RL\] So, as in an equivalence, when we
have an adjunction we have natural transformations
\(e: LR \Rightarrow 1_\mathcal{D}\) and
\(i: 1_\mathcal{C} \Rightarrow RL\). Unlike in an equivalence, they
needn't be natural \emph{isomorphisms}, as the example of sets and
groups shows. But they do have some cool properties, which are nice to
draw using pictures.

First, we draw \(e\) as a U-shaped thing: \[
  \begin{tikzpicture}
    \begin{knot}
      \strand[thick] (0,0)
      to [out=down,in=down,looseness=2] (1,0);
    \end{knot}
    \node[label=above:{$L$}] at (0,0) {};
    \node[label=above:{$R$}] at (1,0) {};
  \end{tikzpicture}
\] The idea here is that \(e\) goes from \(LR\) down to the identity
\(1_\mathcal{D}\), which we draw as ``nothing''. We can think of \(L\)
and \(R\) as processes, and the U-shaped thing as the meta-process of
\(L\) and \(R\) ``colliding into each other and cancelling out'', like a
particle and antiparticle. (Lest you think that's just purple prose,
wait and see! Eventually I'll explain what all this has to do with
antiparticles!) Similarly, we draw \(i\) as an upside-down-U-shaped
thing: \[
  \begin{tikzpicture}
    \begin{knot}
      \strand[thick] (0,0)
      to [out=up,in=up,looseness=2] (1,0);
    \end{knot}
    \node[label=below:{$R$}] at (0,0) {};
    \node[label=below:{$L$}] at (1,0) {};
  \end{tikzpicture}
\] In other words, \(i\) goes from the identity \(1_\mathcal{C}\) to
\(RL\).

We can also use this sort of notation to talk about identity natural
transformations. For example, if we have any old functor \(F\), there is
an identity natural transformation \(1_F\colon F\Rightarrow F\), which
we can draw as follows: \[
  \begin{tikzpicture}
    \begin{knot}
      \strand[thick] (0,0) to (0,2);
    \end{knot}
    \node[label=below:{$F$}] at (0,0) {};
    \node[label=above:{$F$}] at (0,2) {};
  \end{tikzpicture}
\] We draw it as a boring vertical line because ``nothing is happening''
as we go from \(F\) to \(F\).

Now, I haven't talked much about the ways one can compose natural
transformations like \(i\) and \(e\), but remember that they are
\(2\)-morphisms, or morphisms-between-morphisms, in \(\mathsf{Cat}\)
(the \(2\)-category of all categories). This means that they are
inherently \(2\)-dimensional, and in particular, one can compose them
both ``horizontally'' and ``vertically''. I'll explain this more next
time, but for now please take my word for it! Using these composition
operations, one can make sense of the following equations involving
\(i\) and \(e\): \[
  \begin{tikzpicture}
    \begin{knot}
      \strand[thick] (0,0)
      to (0,1)
      to [out=up,in=up,looseness=2] (1,1)
      to [out=down,in=down,looseness=2] (2,1)
      to (2,2);
    \end{knot}
    \node[label=below:{$R$}] at (0,0) {};
    \node[label=above:{$R$}] at (2,2) {};
    \node at (3,1) {$=$};
    \begin{scope}[shift={(4,0)}]
      \begin{knot}
        \strand[thick] (0,0) to (0,2);
      \end{knot}
      \node[label=below:{$R$}] at (0,0) {};
      \node[label=above:{$R$}] at (0,2) {};
    \end{scope}
  \end{tikzpicture}
\] and \[
  \begin{tikzpicture}
    \begin{scope}[xscale=-1,shift={(-2,0)}]
      \begin{knot}
        \strand[thick] (0,0)
        to (0,1)
        to [out=up,in=up,looseness=2] (1,1)
        to [out=down,in=down,looseness=2] (2,1)
        to (2,2);
      \end{knot}
      \node[label=below:{$L$}] at (0,0) {};
      \node[label=above:{$L$}] at (2,2) {};
    \end{scope}
    \node at (3,1) {$=$};
    \begin{scope}[shift={(4,0)}]
      \begin{knot}
        \strand[thick] (0,0) to (0,2);
      \end{knot}
      \node[label=below:{$L$}] at (0,0) {};
      \node[label=above:{$L$}] at (0,2) {};
    \end{scope}
  \end{tikzpicture}
\] In the first equation we are asserting that a certain way of sticking
together \(i\) and \(e\) and some identity natural transformations gives
\(1_R\colon R\Rightarrow R\). In the second we are asserting that some
other way gives \(1_L\colon L\Rightarrow L\).

I will explain these more carefully next time, but for now I mainly want
to state that we can also \emph{define} an adjunction to be a pair of
functors \(L\colon\mathcal{C}\to\mathcal{D}\) and
\(R\colon\mathcal{D}\to\mathcal{C}\) together with natural
transformations \(e\colon LR\Rightarrow 1_\mathcal{D}\) and
\(i\colon1_\mathcal{C}\Rightarrow RL\) making the above 2 equations
hold! This is the definition of ``adjunction'' that is the most similar
to the definition of ``equivalence''.

Now, topologically, these 2 equations simply say that if you have a
wiggly curve like \[
  \begin{tikzpicture}
    \begin{knot}
      \strand[thick] (0,0)
      to (0,1)
      to [out=up,in=up,looseness=2] (1,1)
      to [out=down,in=down,looseness=2] (2,1)
      to (2,2);
    \end{knot}
  \end{tikzpicture}
\] or \[
  \begin{tikzpicture}
    \begin{scope}[xscale=-1,shift={(-2,0)}]
      \begin{knot}
        \strand[thick] (0,0)
        to (0,1)
        to [out=up,in=up,looseness=2] (1,1)
        to [out=down,in=down,looseness=2] (2,1)
        to (2,2);
      \end{knot}
    \end{scope}
  \end{tikzpicture}
\] you can pull it tight to get \[
  \begin{tikzpicture}
    \begin{knot}
      \strand[thick] (0,0) to (0,2);
    \end{knot}
  \end{tikzpicture}
\] Thus, while \[
  \begin{tikzpicture}
    \begin{knot}
      \strand[thick] (0,0)
      to [out=down,in=down,looseness=2] (1,0);
    \end{knot}
  \end{tikzpicture}
\] and \[
  \begin{tikzpicture}
    \begin{knot}
      \strand[thick] (0,0)
      to [out=up,in=up,looseness=2] (1,0);
    \end{knot}
  \end{tikzpicture}
\] are not exactly ``inverses'', there is some subtler sense in which
they ``cancel out''. This corresponds to the notion that while adjoint
functors are not inverses, not even up to a natural isomorphism, they
still are ``like inverses'' in a subtler sense.

Now this may seem like a silly game, drawing natural transformations as
``string diagrams'' and interpreting adjoint functors as wiggles in the
string. But in fact this is part of a very big, very important, and very
fun game: the relation between \(n\)-category theory and the topology of
submanifolds of \(\mathbb{R}^n\). Right now we are dealing with
\(\mathsf{Cat}\), which is a \(2\)-category, so we are getting into
\(2\)-dimensional pictures. But when we get into \(3\)-categories we
will get into \(3\)-dimensional pictures, and knot theory... and
what got me into this whole business in the first place: the relation
between knots and physics. In higher dimensions it gets even fancier.

So I will continue next time and explain the recipes for composing
natural transformations, and the associated string diagrams, more
carefully. To continue reading the ``Tale of \(n\)-Categories'', see
\protect\hyperlink{week80_tale}{``Week 80''}.

\hypertarget{week80}{%
\section{April 20, 1996}\label{week80}}

There are a number of interesting books I want to mention.

Huw Price's book on the arrow of time is finally out! It's good to see a
philosopher of science who not only understands what modern physicists
are up to, but can occasionally beat them at their own game.

Why is the future different from the past? This has been vexing people
for a long time, and the stakes went up considerably when Boltzmann
proved his ``H-theorem'', which seems at first to show that the entropy
of a gas always increases, despite the time-reversibility of the laws of
classical mechanics. However, to prove the H-theorem he needed an
assumption, the ``assumption of molecular chaos''. It says roughly that
the positions and velocities of the molecules in a gas are uncorrelated
before they collide. This seems so plausible that one can easily
overlook that it has a time-asymmetry built into it --- visible in the
word ``before''. In fact, we aren't getting something for nothing in the
H-theorem; we are making a time-asymmetric assumption in order to
conclude that entropy increases with time!

The ``independence of incoming causes'' is very intuitive: if we do an
experiment on an electron, we almost always assume our choice of how to
set the dials is not correlated to the state of the electron. If we drop
this time-asymmetric assumption, the world looks rather
different... but I'll let Price explain that to you.

Anyway, Price is an expert at spotting covertly time-asymmetric
assumptions. you may remember from \protect\hyperlink{week26}{``Week
26''} that he even got into a nice argument with Stephen Hawking about
the arrow of time, thanks to this habit of his. You can read more about
it in:

\begin{enumerate}
\def\labelenumi{\arabic{enumi})}
\tightlist
\item
  Huw Price, \emph{Time's Arrow and Archimedes' Point: New Directions
  for a Physics of Time}, Oxford U.\ Press,  Oxford, 1996.
\end{enumerate}

Also, there is a new book out by Hawking and Roger Penrose on quantum
gravity. First they each present their own ideas, and then they duke it
out in a debate in the final chapter. This book is an excellent place to
get an overview of some of the main ideas in quantum gravity. It helps
if you have a little familiarity with general relativity, or
differential geometry, or are willing to fake it.

There is even some stuff here about the arrow of time! Hawking has a
theory of how it arose, starting from his marvelous ``no-boundary
boundary conditions'', which say that the wavefunction of the universe
is full of quantum fluctuations corresponding to big bangs which erupt
and then recollapse in big crunches. The wavefunction itself has no
obvious ``time-asymmetry'', indeed, time as we know it only makes sense
\emph{within} any one of the quantum fluctuations, one of which is
presumably the world we know! But Hawking thinks that each of these
quantum fluctuations, or at least most of them, should have an arrow of
time. This is what Price raised some objections to. Hawking seems to
argue that each quantum fluctuation should ``start out'' rather smooth
near its big bang and develop more inhomogeneities as time passes,
``winding up'' quite wrinkly near its big crunch. But it's not at all
clear what this ``starting out'' and ``winding up'' means. Possibly he
is simply speaking vaguely, and all or most of the quantum fluctuations
can be shown to have one smooth end and wrinkly at the other. That would
be an adequate resolution to the arrow of time problem. But it's not
clear, at least not to me, that Hawking really showed this.

Penrose, on the other hand, has some closely related ideas. His ``Weyl
curvature hypothesis'' says that the Weyl curvature of spacetime goes to
zero at initial singularities (e.g.~the big bang) and infinity at final
ones (e.g.~black holes). The Weyl curvature can be regarded as a measure
of the presence of inhomogeneity --- the ``wrinkliness'' I alluded to
above. The Weyl curvature hypothesis can be regarded as a
time-asymmetric law built into physics from the very start.

To see them argue it out, read

\begin{enumerate}
\def\labelenumi{\arabic{enumi})}
\setcounter{enumi}{1}
\tightlist
\item
  Stephen Hawking and Roger Penrose, \emph{The Nature of Space and
  Time}, Princeton U.\ Press, Princeton, 1996.
\end{enumerate}

There are also a couple of more technical books on general relativity
that I'd been meaning to get ahold of for a long time. They both feature
authors of that famous book,

\begin{enumerate}
\def\labelenumi{\arabic{enumi})}
\setcounter{enumi}{2}
\tightlist
\item
  Charles Misner, Kip Thorne and John Wheeler, \emph{Gravitation},
  Freeman Press, 1973,
\end{enumerate}

which was actually the book that made me decide to work on quantum
gravity, back at the end of my undergraduate days. They are:

\begin{enumerate}
\def\labelenumi{\arabic{enumi})}
\setcounter{enumi}{3}
\tightlist
\item
  Ignazio Ciufolini and John Archibald Wheeler, \emph{Gravitation and
  Inertia}, Princeton U.\ Press, Princeton, 1995.
\end{enumerate}

and

\begin{enumerate}
\def\labelenumi{\arabic{enumi})}
\setcounter{enumi}{4}
\tightlist
\item
  Kip Thorne, Richard Price and Douglas Macdonald, eds., \emph{Black
  Holes: The Membrane Paradigm}, Yale U.\ Press, New Haven, 1986.
\end{enumerate}

The book by Ciufolini and Wheeler is full of interesting stuff, but it
concentrates on ``gravitomagnetism'': the tendency, predicted by general
relativity, for a massive spinning body to apply a torque to nearby
objects. This is related to Mach's old idea that just as spinning a
bucket pulls the water in it up to the edges, thanks to the centrifugal
force, the same thing should happen if instead we make lots of stars
rotate around the bucket! Einstein's theory of general relativity was
inspired by Mach, but there has been a long-running debate over whether
general relativity is ``truly Machian'' --- in part because nobody knows
what ``truly Machian'' means. In any event, Ciufolini and Wheeler argue
that gravitomagnetism exhibits the Machian nature of general relativity,
and they give a very nice tour of gravitomagnetic effects.

That is fine in theory. However, the gravitomagnetic effect has never
yet been observed! It was supposed to be tested by Gravity Probe B, a
satellite flying at an altitude of about 650 kilometers, containing a
superconducting gyroscope that should precess at a rate of 42
milliarcseconds per year thanks to gravitomagnetism. I don't know what
ever happened with this, though: the following web page says ``Gravity
Probe B is expected to fly in 1995'', but now it's 1996, right? Maybe
someone can clue me in to the latest news\ldots. I seem to remember some
arguments about funding the program.

\begin{enumerate}
\def\labelenumi{\arabic{enumi})}
\setcounter{enumi}{5}
\tightlist
\item
  Gravity Probe B, \href{https://web.archive.org/web/19970506010314/http://stugyro.stanford.edu/RELATIVITY/GPB/GPB.html}{\texttt{https://web.archive.org/web/19970506010314/}}\hfill \break  \href{https://web.archive.org/web/19970506010314/http://stugyro.stanford.edu/RELATIVITY/GPB/GPB.html}{\texttt{http://stugyro.stanford.edu/RELATIVITY/GPB/GPB.html}}.
\end{enumerate}

(Note added in 2002: now this webpage is gone; see
\url{http://einstein.stanford.edu/} for the latest story.)

Kip Thorne's name comes up a lot in conjuction with black holes and the
LIGO --- or Laser-Interferometer Gravitational-Wave Observatory ---
project. As pairs of black holes or neutron stars emit
gravitational radiation, they should spiral in towards each other. In
their final moments, as they merge, they should emit a ``chirp'' of
gravitational radiation, increasing in frequency and amplitude until
their ecstatic union is complete. The LIGO project aims to observe these
chirps, and any other sufficiently strong gravitational radiation that
happens to be passing by our way. LIGO aims to do this by using laser
interferometry to measure the distance between two points about 4
kilometers apart to an accuracy of about \(10^{-18}\) meters, thus
detecting tiny ripples in the spacetime metric. For more on LIGO, try:

\begin{enumerate}
\def\labelenumi{\arabic{enumi})}
\setcounter{enumi}{6}
\tightlist
\item
  LIGO project home page, \href{http://www.ligo.caltech.edu/}{\texttt{http://www.ligo.caltech.edu/}}.
\end{enumerate}

Thorne helped develop a nice way to think of black holes by envisioning
their event horizon as a kind of ``membrane'' with well-defined
mechanical, electrical and magnetic properties. This is called the
membrane paradigm, and is useful for calculations and understanding what
black holes are really like. The book ``Black Holes: The Membrane
Paradigm'' is a good place to learn about this.

\begin{center}\rule{0.5\linewidth}{0.5pt}\end{center}

\hypertarget{week80_tale}{
Now let me return to the ``Tale of \(n\)-Categories'' --- and in particular,
\(2\)-categories.} So far I've said only
that a \(2\)-category is some sort of structure with objects, morphisms
between objects, and \(2\)-morphisms between morphisms. But I have been
attempting to develop your intuition for \(\mathsf{Cat}\), the
primordial example of a \(2\)-category. Remember, \(\mathsf{Cat}\) is
the \(2\)-category of all categories! Its objects are categories, its
morphisms are functors, and its \(2\)-morphisms are natural
transformations --- these being defined in
\protect\hyperlink{week73_tale}{``Week 73''} and again in
\protect\hyperlink{week75_tale}{``Week 75''}.

How can you learn more about \(2\)-categories? Well, a really good place
is the following article by Ross Street, who is one of the great gurus
of \(n\)-category theory. For example, he was the one who invented
\(\omega\)-categories!

\begin{enumerate}
\def\labelenumi{\arabic{enumi})}
\setcounter{enumi}{7}
\tightlist
\item
  Ross Street, ``Categorical structures'', in \emph{Handbook of
  Algebra}, vol.~\textbf{1}, ed.~M. Hazewinkel, Elsevier, Amsterdam, 1996.
\end{enumerate}
\noindent
Physicists should note his explanation of the Yang--Baxter and
Zamolodchikov equations in terms of category theory. If you have trouble
finding this, you might try

\begin{enumerate}
\def\labelenumi{\arabic{enumi})}
\setcounter{enumi}{8}
\tightlist
\item
  G.\ Maxwell Kelly and Ross Street, ``Review of the elements of
  \(2\)-categories'', in Lecture Notes in Mathematics \textbf{420}, Springer,
  Berlin, 1974, pp.~75--103.
\end{enumerate}

I can't really compete with these for thoroughness, but at least let me
give the definition of a \(2\)-category. I'll give a pretty
nuts-and-bolts definition; later I'll give a more elegant and abstract
one. Readers who are familiar with \(\mathsf{Cat}\) should keep this
example in mind at all times!

This definition is sort of long, so if you get tired of it, concentrate
on the pictures! They convey the basic idea. Also, keep in mind is that
this is going to be sort of like the definition of a category, but with
an extra level on top, the \(2\)-morphisms.

So: first of all, a \(2\)-category consists of a collection of
``objects'' and a collection of ``morphisms''. Every morphism \(f\) has
a ``source'' object and a ``target'' object. If the source of \(f\) is
\(x\) and its target is y, we write \(f\colon x\to y\). In addition, we
have:

\begin{enumerate}
\def\labelenumi{\arabic{enumi})}
\item
  Given a morphism \(f\colon x\to y\) and a morphism \(g\colon y\to Z\),
  there is a morphism \(fg\colon x\to Z\), which we call the
  ``composite'' of \(f\) and \(g\).
\item
  Composition is associative: \((fg)h = f(gh)\).
\item
  For each object \(x\) there is a morphism \(1_x\colon x \to x\),
  called the ``identity'' of x. For any \(f\colon x\to y\) we have
  \(1_x f = f 1_y = f\).
\end{enumerate}

You should visualize the composite of \(f\colon x\to y\) and
\(g\colon y\to Z\) as follows: \[x\xrightarrow{f}y\xrightarrow{g}Z.\] So
far this is exactly the definition of a category! But a \(2\)-category
\emph{also} consists of a collection of ``2-morphisms''. Every \(2\)-morphism
\(T\) has a ``source'' morphism \(f\) and a target morphism \(g\). If
the source of \(T\) is \(f\) and its target is \(g\), we write
\(T\colon f\Rightarrow g\). If \(T\colon f\Rightarrow g\), we require
that \(f\) and \(g\) have the same source and the same target; for
example, \(f\colon x\to y\) and \(g\colon x\to y\).  You should visualize
\(T\) as follows:
\[\includegraphics[scale=0.3]{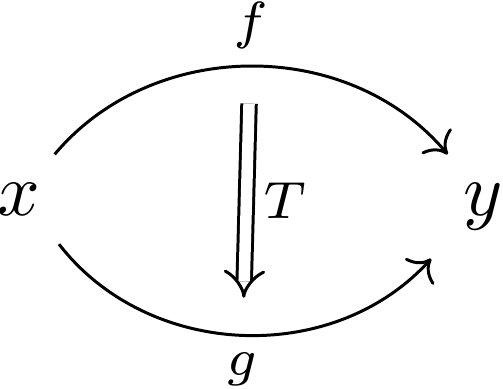}\] In addition, we
have:
\vskip 0.5em

1\({}^\prime\)) Given a \(2\)-morphism \(S\colon f\Rightarrow g\) and a
\(2\)-morphism \(T\colon g\Rightarrow h\), there is a \(2\)-morphism
\(ST\colon f\Rightarrow h\), which we call the ``vertical composite'' of
\(S\) and \(T\).

\vskip 0.5em
2\({}^\prime\)) Vertical composition is associative: \((ST)U = S(TU)\).

\vskip 0.5em
3\({}^\prime\)) For each morphism \(f\) there is a \(2\)-morphism
\(1_f\colon f\Rightarrow f\), called the ``identity'' of \(f\). For any
\(T\colon f\Rightarrow g\) we have \(1_f T = T 1_g = T\).

\vskip 0.5em
Note that these are just like the previous 3 rules. We draw the vertical
composite of \(S\colon f\Rightarrow g\) and \(T\colon g\Rightarrow h\)
like this: \[\includegraphics[scale=0.3]{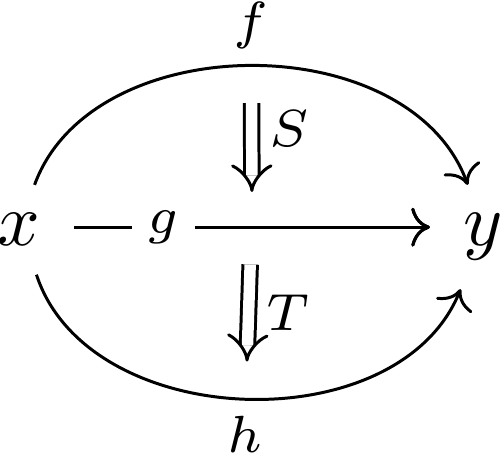}\] Now
for a twist. We also require that we can ``horizontally'' compose
2-morphisms as follows:
\[\includegraphics[scale=0.3]{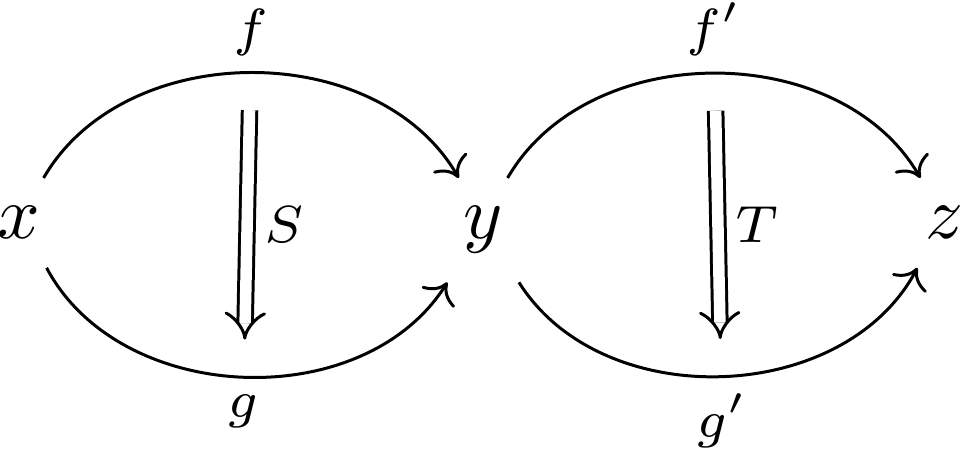}\] So we also
demand:

\vskip 0.5em
1\({}^{\prime\prime}\)) Given morphisms \(f,g\colon x\to y\) and \(f',g'\colon y\to z\),
and \(2\)-morphisms \(S\colon f\Rightarrow g\) and
\(T\colon f'\Rightarrow g'\), there is a \(2\)-morphism
\(S\cdot T\colon ff' \Rightarrow gg'\), which we call the ``horizontal
composite'' of \(S\) and \(T\).

\vskip 0.5em
2\({}^{\prime\prime}\))  Horizontal composition is associative:
\((S\cdot T)\cdot U = S\cdot (T\cdot U)\).

\vskip 0.5em
3\({}^{\prime\prime}\))  The identities for vertical composition are also the identities
for horizontal composition. That is, given \(f,g\colon x\to y\) and
\(T\colon f\Rightarrow g\) we have
\(1_{1_x}\cdot T = T\cdot 1_{1_y} = T\).

\vskip 0.5em \noindent
Finally, we demand the ``exchange law'' relating horizontal and vertical
composition: \[(ST)\cdot (S'T') = (S\cdot S')(T\cdot T')\] This makes
the following \(2\)-morphism unambiguous:
\[\includegraphics[scale=0.3]{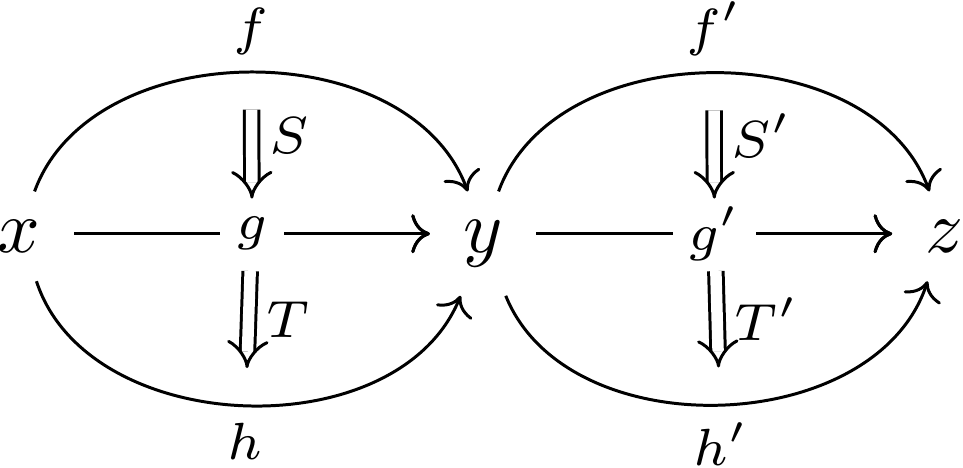}\] We can think of it
either as the result of first doing two vertical composites, and then
one horizontal composite, or as the result of first doing two horizontal
composites, and then one vertical composite!

Here we can really see why higher-dimensional algebra deserves its name.
Unlike category theory, where we can visualize morphisms as
1-dimensional arrows, here we have \(2\)-morphisms which are
intrinsically 2-dimensional, and can be composed both vertically and
horizontally.

Now if you are familiar with \(\mathsf{Cat}\), you may be wondering how
we vertically and horizontally compose natural transformations, which
are the 2-morphisms in \(\mathsf{Cat}\). Let me leave this as an
exercise for now... there's a nice way to do it that makes
\(\mathsf{Cat}\) into a \(2\)-category. This exercise is a good one to
build up your higher-dimensional algebra muscles.

In fact, we could have invented the above definition of \(2\)-category
simply by thinking a lot about \(\mathsf{Cat}\) and what you can do with
categories, functors, and natural transformations. I'm pretty sure
that's more or less what happened, historically! Thinking hard enough
about \(n\mathsf{Cat}\) leads us on to the definition of \((n+1)\)-categories\ldots.

But that's enough for now. Typing those diagrams is hard work.  To continue reading the ``Tale of \(n\)-Categories'', see
\protect\hyperlink{week83_tale}{``Week 83''}.

\begin{center}\rule{0.5\linewidth}{0.5pt}\end{center}

I thank Keith Harbaugh for catching lots of typos and other mistakes in
\protect\hyperlink{week73}{``Week 73''} --
\protect\hyperlink{week80}{``Week 80''}.

\hypertarget{week81}{%
\section{May 12, 1996}\label{week81}}

I think I'll take a little break on the continuing saga of \(n\)-categories.
Instead I'll talk about something which is secretly the very same
subject, namely loop groups and their central extensions. This is
important in string theory. But first I want to say a bit about some
very high-energy physics.

\begin{enumerate}
\def\labelenumi{\arabic{enumi})}
\tightlist
\item
  D.\ J.\ Bird \emph{et al}, ``Detection of a cosmic ray with measured energy
  well beyond the expected spectral cutoff due to cosmic microwave
  radiation'', \emph{The Astrophysical Journal} \textbf{441} (1995), 144--150.  
  Also available as
  \href{https://arxiv.org/abs/astro-ph/9410067}{\texttt{astro-ph/9410067}}

  P.\ Bhattacharjee and G.\ Sigl, ``Monopole annihilation and highest energy
  cosmic rays'', \emph{Phys. Rev. D} \textbf{51} (1995), 4079.  Also available as
  \href{https://arxiv.org/abs/astro-ph/9412053}{\texttt{astro-ph/9412053}}.

   R. J. Protheroe and P. A. Johnson, ``Are topological defects responsible
  for the 300 EeV cosmic rays?'', \emph{Nucl. Phys. B - Proc. Suppl. } \textbf{48} 
  (1996), 485--487. Also available as
\href{https://arxiv.org/abs/astro-ph/9605006}{\texttt{astro-ph/9605006}}.
\end{enumerate}

In 1994, folks at the Fly's Eye air shower detector in Utah observed a
cosmic ray whose energy was about 320 EeV. In case you forget what an
EeV is, it's a unit of energy equal to a billion GeV, and a Gev is equal
to a billion ev (electron volts). Current particle accelerators
routinely create particles with energies about a hundred GeV, but a few hundred
\emph{EeV} is a whole different story! That's about 50 joules, the
energy of a one-kilogram mass moving at 10 meters/second, all packed in
one particle!

Nobody knows what would produce cosmic rays of this energy. To make the
puzzle more mysterious, this energy is above the Greisen--Zatsepin--Kuz'min 
(or GZK) cutoff for cosmic rays produced at moderate
extragalactic distances (30 megaparsecs). The idea of the GZK cutoff is
that particles of extremely high energies whizzing through space would
interact significantly with the cosmic microwave background radiation,
losing energy to produce pions.

So it seems that something is producing really high energy particles,
and this something is not too far away, by cosmic standards. Established
mechanisms don't get energies that high. A possibility studied by
various authors including P. Bhattacharjee and G. Sigl is that these
super-energetic cosmic rays are produced by the decay of ``topological
defects''. Various grand unified theories, or GUTs, predict that the
strong, weak, and electromagnetic forces all act the same at really
high temperatures, while at low temperatures (like any sort of
temperature you'd find around here) a ``spontaneous symmetry breaking''
occurs which makes them split up into their observed distinct
personalities.

Mathematically this is a bit like how a magnet at low temperatures
randomly picks out a certain axis of magnetization, breaking the
rotational symmetry it possesses at high temperatures. And like in the
case of a magnet, one would expect the possibility of ``topological
defects'' where different regions of space pick different ways to break
the symmetry, leaving nasty spots like lumps in the carpet that can't be
straightened out. Ordinary magnets typically exhibit \(2\)-dimensional
``domain walls'' between domains having different axes of magnetization.
But in various GUTs folks have considered, one can also get
1-dimensional ``cosmic strings'' and 0-dimensional ``topological
solitons'' including magnetic monopoles --- particles with magnetic
charge. None of these topological defects have ever been observed,
despite a fair amount of searching. Could super-energetic cosmic rays be
the result of a monopole-antimonopole collision?

Alas, Protheroe and Johnson's paper argues that in such decays lots of
the energy would go into the production of high-energy \(\gamma\)
rays... more than has been observed in the super-energetic cosmic
ray showers. So maybe the puzzle has some other answer.

The weekend before last I went to the 11th Geometry Festival, which was
held at the University of Maryland. Since I work on quantum gravity, I
could be said to be a geometer of sorts --- perhaps a quantum geometer.
But geometry means a lot of different things to different people, and
this conference concentrated on some aspects of geometry that I don't
know much about. In particular, there were talks by Schmuel Weinberger,
Bruce Kleiner and G. Wei on the implications of positive and negative
curvature for Riemannian geometry.

A talk that was right up my alley was given by Jean-Luc Brylinski. It
dealt with themes from his papers with McLaughlin:

\begin{enumerate}
\def\labelenumi{\arabic{enumi})}
\setcounter{enumi}{1}
\item
  Jean-Luc Brylinski and Dennis A. McLaughlin, ``The geometry of degree
  four characteristic classes and of line bundles on loop spaces, I'',
  \emph{Duke Math.\ Journal} \textbf{75} (1994), 603--638. ``II'',
  \emph{Duke Math. Journal} \textbf{83} (1996), 105--139.

  Jean-Luc Brylinski, ``Central extensions and reciprocity laws'',
  \emph{Cahiers de Topologie et G\'eom\'etrie Diff\'erentielle Cat\'egoriques} 
  \textbf{38} (1997), 193--215. 

  Jean-Luc Brylinski, ``Coadjoint orbits of central extensions of gauge
  groups'',  \emph{Commun. Math. Phys.} \textbf{188} (1997), 351--365.

  Jean-Luc Brylinski and Dennis A. McLaughlin, ``The geometry of two
  dimensional symbols'', \emph{K-Theory} \textbf{10} (1996), 215--237.
\end{enumerate}

Let me say a bit about the math underlying these papers, the basic stuff
that they build on. One hot topic in mathematical physics in the last
decade has been the study of ``loop groups''. Say you take any Lie group
\(G\). Then the ``loop group'' \(LG\) is the set of smooth functions
from the circle to \(G\). This becomes a group with pointwise
multiplication as the group operation. This sort of group shows up in
\(2\)-dimensional quantum field theory, where spacetime could be the
cylinder. Then ``space'' is the circle, and if we are studying gauge
theory with gauge group \(G\), the group of gauge transformations over
space would be the loop group \(LG\). One main reason for being
interested in \(2\)-dimensional quantum field theory is string theory:
here we think of the \(2\)-dimensional worldsheet of the string as a
spacetime in its own right, and we are often interested in doing gauge
theory over this spacetime. For this reason, folks in string theory need
to understand all they can about unitary representations of loop groups.

Actually they are interested in \emph{projective} representations of
loop groups. Remember, in quantum mechanics two vectors in a Hilbert
space give the same expectation values for any observable if they differ
only by a phase. So it is perfectly fine for a group representation
\(R\) to satisfy the usual law \[R(g)R(h) = R(gh)\] where \(g\), \(h\)
are group elements, \emph{only up to a phase}. So in the definition of a
projective representation we weaken the above requirement to
\[R(g) R(h) = c(g,h) R(gh)\] where \(c(g,h)\) is a phase depending on
\(g\) and \(h\). Folks call \(c(g,h)\) the ``cocycle'' of the projective
representation.

A projective unitary representation of a group \(H\) can also be thought
of as a representation of a bigger group \(\widetilde{H}\) called a
``central extension'' of \(H\). The idea is that this bigger group has a
bunch of phases built into it to absorb the phase ambiguities in the
projective representation of \(H\). Let \(\mathrm{U}(1)\) be the unit
circle in the complex plane, a group under multiplication. This is the
group of phases. We can think of \(\widetilde{H}\) as
\(H \times \mathrm{U}(1)\) given a sneaky product designed to soak up
the phases produced by the cocycle: \[(g, a)(h, b) = (gh, ab c(g,h)).\]
We can define a unitary representation \(S\) of \(\widetilde{H}\) as
follows: \[S(g, a) = R(g)a.\] It's then obvious that
\[S(g, a) S(h, b) = S((g, a)(h, b))\] so \(S\) is really a
representation. For this reason, if we are doing quantum theory and we
don't like projective representations, it's okay as long as we
understand the central extensions of our group of symmetries.

So, instead of thinking about projective representations of loop groups,
we can think about central extensions of loop groups. How does one get
ahold of these? There is a nice trick which Brylinski described in his
talk. To get this trick, we need to assume a bit about the group \(G\).
Let's assume it's a connected and simply-connected simple Lie group.
I'll explain that in a minute, but some good examples to keep in mind
are \(\mathrm{SU}(n)\) and \(\mathrm{Spin}(n)\); see
\protect\hyperlink{week61}{``Week 61''} for the definitions and a bit of
information about these groups.

Now remember that \(S^k\) stands for the \(k\)-dimensional sphere, and
\(\pi_k(X)\) of a topological space \(X\) stands for the set of
continuous maps from \(S^k\) to \(X\), modulo homotopy. In other words,
two continuous maps from \(S^k\) to \(X\) define the same element of
\(\pi_k(X)\) if one can be continuously deformed to the other.

Saying that \(G\) is connected means that \(\pi_0(G) = 0\). To
understand this you need to realize that \(S^0\) consists of two points.
So \(\pi_0(G) = 0\) means that \(G\) consists of one piece, any two
points of which can be connected by a continuous path.

Saying that \(G\) is simply connected means that \(\pi_1(G) = 0\). In
other words, all loops in \(G\) can be ``pulled tight''. A good example
of a group that's \emph{not} simply connected is the group \(\mathrm{SO}(n)\)
of rotations in \(n\) dimensional space. This flaw with
\(\mathrm{SO}(n)\) is why they needed to invent \(\mathrm{Spin}(n)\);
see \protect\hyperlink{week61}{``Week 61''}.

As it turns out, every Lie group has \(\pi_2(G) = 0\). So all 2-spheres
in \(G\) can be pulled tight too. Imagine taking a balloon and sticking
it in \(G\); then you can always shrink it down to a point in a
continuous way without it getting stuck around a hole in \(G\).

Saying that \(G\) is simple is an algebraic rather than topological
condition, and I explained this condition in
\protect\hyperlink{week63}{``Week 63''}. But it has topological
ramifications. It implies, for example, that \(\pi_3(G) = \mathbb{Z}\),
the group of integers. So to each way of sticking a 3-sphere in \(G\) we
can associate an integer. One way to compute this integer involves the
Killing form on the Lie algebra of \(G\). This is a special inner
product on the Lie algebra of \(G\). Using this inner product and the
bracket in the Lie algebra we can convert 3 vectors \(u\), \(v\), and
\(w\) in the Lie algebra into a number as follows:
\[W(u,v,w) = k \langle [u,v],w \rangle\] Here \(k\) is a constant that
will make life simpler later. The special property of the Killing form
implies that \(W\) is totally antisymmetric, and we can use left
multiplication to translate \(W\) all over the group \(G\) obtaining a
closed \(3\)-form on \(G\). Call this \(3\)-form \(W\), too. Then, given
any smooth function from \(S^3\) into \(G\) we can pull back \(W\) to
\(S^3\) and integrate it over \(S^3\). If we choose the constant \(k\)
right, the result will be an integer --- the integer we were looking
for.

Hmm, this is getting technical. Well, it'll keep getting more technical.
Just stop reading when it becomes unpleasant.

Okay, these topological facts about the group \(G\) have a bunch of cool
consequences. One trick usually goes by the name of the ``WZW action'',
which refers to Wess, Zumino, and Witten. Say we have smooth function
\(f\) from \(S^2\) to \(G\). Since \(\pi_2(G) = 0\) we can extend \(f\)
to a smooth function \(F\) from the \(3\)-dimensional ball, \(D^3\), to
\(G\). (This is just another way of ``pulling the balloon tight'' as
mentioned above.) Now we can use \(F\) to pull back the magic \(3\)-form
\(W\) to \(D^3\), and then we can integrate the resulting \(3\)-form
over \(D^3\) to get a number \(S(f)\) called the Wess--Zumino--Witten
action.

Unfortunately, this number depends on the choice of the function \(F\)
extending \(f\) to the ball. Fortunately, it doesn't depend too much on
\(F\). Say we extended \(f\) to some other function \(F'\) from the ball
to \(G\). Then \(F\) together with \(F'\) define a map from \(S^3\) to
\(G\), and we know from the previous stuff that the integral of the
pullback of \(W\) over this \(S^3\) is an integer. So changing our
choice of an extension of \(f\) only changes \(S(f)\) by an integer.
This means that the exponential of the WZW action:
\[\exp(2 \pi i S(f))\] is completely well-defined. This is nice in
quantum physics, where the exponential of the action is what really
matters. Note also that this exponential is just a phase! So we are
getting a function \[A\colon\mathrm{Maps}(S^2,G)\to \mathrm{U}(1)\]
assigning a phase to any map \(f\) from \(S^2\) to \(G\).

Now \(\mathrm{Maps}(S^2,G)\) is sort of like the loop group, since the
loop group is just \(\mathrm{Maps}(S^1,G)\). In particular, it too
becomes a group by pointwise multiplication. A bit of calculation shows
that \(A\) above is a group homomorphism: \[A(f) A(g) = A(fg).\] This
homomorphism is the key to finding the central extension of the loop
group. Here's how we do it. By now everyone but the experts has probably
fallen asleep at the screen, so I can pull out all the stops.

Here's a useful way to think of a central extensions: a central
extension \(\widetilde{H}\) of the group \(H\) by the group
\(\mathrm{U}(1)\) is a special sort of short exact sequence of groups:
\[1 \to \mathrm{U}(1) \to \widetilde{H} \to H \to 1\] By ``short exact
sequence of groups'' I simply mean that \(\mathrm{U}(1)\) is a subgroup
of \(\widetilde{H}\) and that \(\widetilde{H}\) modulo \(\mathrm{U}(1)\)
is \(H\). What's special about central extensions is that
\(\mathrm{U}(1)\) is in the \emph{center} of \(\widetilde{H}\). You can
check that this definition of central extension matches up with our
earlier more lowbrow definition.

Now how do we get this short exact sequence? Well, it comes from a short
exact sequence of spaces: \[\{*\} \to S^1 \to D^2 \to S^2 \to \{*\}\]

This diagram means simply that we can think of the circle as a subspace
of the \(2\)-dimensional disc \(D^2\) in an obvious way, and then if we
collapse this circle to a point the disc gets squashed down to a
2-sphere. Now, from this short exact sequence we get a short exact
sequence of groups
\[1 \to \mathrm{Maps}(S^2,G) \to \mathrm{Maps}(D^2,G) \to \mathrm{Maps}(S^1,G) \to 1\]

In other words, \(\mathrm{Maps}(S^2,G)\) is a normal subgroup of
\(\mathrm{Maps}(D^2,G)\), and if we mod out by this subgroup we get
\(\mathrm{Maps}(S^1,G)\). Now we can use the homomorphism
\(A\colon\mathrm{Maps}(S^2,G)\to \mathrm{U}(1)\) to get ourselves
another exact sequence like this: \[
  \begin{tikzcd}
    1 \rar
    & \mathrm{Maps}(S^2,G)
      \rar["i"] \dar["A"]
    & \mathrm{Maps}(D^2,G)
      \rar["j"] \dar
    & \mathrm{Maps}(S^1,G)
      \rar \dar["1"]
    & 1
  \\1 \rar
    & \mathrm{U}(1)
      \rar["i"]
    & \widetilde{L}
      \rar["j"]
    & \mathrm{Maps}(S^1,G)
      \rar
    & 1
  \end{tikzcd}
\] Remembering that \(\mathrm{Maps}(S^1,G)\) is the loop group,
\(\widetilde{L}\) turns out to be the desired central extension!
Concretely we can think of \(\widetilde{L}\) as a quotient group of
\(\mathrm{Maps}(D^2,G)\times \mathrm{U}(1)\) by the subgroup of pairs of
the form \((i(f),A(f))\) with \(f\) in \(\mathrm{Maps}(S^2,G)\).

There is something fascinating about how spheres of different dimensions
--- \(S^0\), \(S^1\), \(S^2\), and \(S^3\) --- conspire together with
the topology of the group \(G\) to yield the central extension of the
loop group \(LG\). It appears that what we are really studying are the
closely related cohomology groups:

\begin{itemize}
\tightlist
\item
  \(H^0(\mathrm{Maps}(S^3,G))\) which is just another way of saying
  \(\pi_3(G)\)
\item
  \(H^1(\mathrm{Maps}(S^2,G))\) which describes homomorphisms from
  \(\mathrm{Maps}(S^2,G)\) to \(\mathrm{U}(1)\)
\item
  \(H^2(\mathrm{Maps}(S^1,G))\) which describes central extensions of
  \(\mathrm{Maps}(S^1,G)\)
\item
  \(H^3(\mathrm{Maps}(S^0,G))\) which is just another way of saying
  \(H^3(G)\), where \(W\) lives.
\end{itemize}

There is a fourth term in this series which I didn't get around to
talking about; it's

\begin{itemize}
\tightlist
\item
  \(H^4(\mathcal{B}G)\) where the degree 4 characteristic class for
  \(G\)-bundles, e.g.~the 2nd Chern class for \(\mathrm{SU}(n)\), lives.
\end{itemize}

Here \(\mathcal{B}G\) is the ``classifying space'' of \(G\). I would
like to understand more deeply what's going on with this series, because
the different terms have a lot to do with physics in different
dimensions --- dimensions 0 to 4, the ``low dimensions'' which are so
specially interesting.

I should conclude by noting that while a lot of this appeared in
Brylinski's talk, and a lot of it is probably common knowledge among
topologists, it was in some conversations with James Dolan that we
worked out some of the patterns I mention here.

\hypertarget{week82}{%
\section{May 17, 1996}\label{week82}}

I will continue to take a break from the tale of \(n\)-categories. As
the academic year winds to an end, an enormous pile of articles and
books is building up on my desk. I can kill two birds with one stone if
I list some of them while filing them. Here is a sampling:

\begin{enumerate}
\def\labelenumi{\arabic{enumi})}
\tightlist
\item
  \emph{Advances in Applied Clifford Algebras}, ed.~Jaime Keller.
\end{enumerate}
\noindent
This is a homegrown journal for fans of Clifford algebras. What are
Clifford algebras? Well, let's start at the beginning, with the
quaternions\ldots.

As J. Lambek has pointed out, not many mathematicians can claim to have
introduced a new kind of number. One of them was the Sir William Rowan
Hamilton. He knew about the real numbers \(\mathbb{R}\), of course, and
also the complex numbers \(\mathbb{C}\), which are the reals with a
square root of \(-1\), usually called \(i\), thrown in. Why not try
putting in another square root of \(-1\)? This might give a
\(3\)-dimensional algebra that'd help with \(3\)-dimensional space as
much as the complex numbers help with 2 dimensions. He tried this but
couldn't get division to work out well. He struggled this for a long
time. On the 16th of October, 1843, he was walking along the Royal Canal
with his wife to a meeting of the Royal Irish Academy when he had a good
idea: ``\ldots there dawned on me the notion that we must admit, in some
sense, a fourth dimension of space for the purpose of calculating with
triples... An electric circuit seemed to close, and a spark
flashed forth.'' He carved the decisive relations
\[i^2 = j^2 = k^2 = ijk = -1\] in the stone of Brougham Bridge as he
passed it. This was bold: a \emph{noncommutative} algebra, since
\(ij = -ji\), \(jk = -kj\), and \(ik = -ki\) follow from the above
equations. These are the quaternions, which now we call \(\mathbb{H}\)
after Hamilton.

Hamilton wound up spending much of his time on quaternions. The lawyer
and mathematician Arthur Cayley heard Hamilton lecture on quaternions
and --- I imagine --- was influenced by this to invent his
``octonions'', an 8-dimensional nonassociative algebra in which division
still works nicely. For more on quaternions, octonions, and the general
subject of division algebras, try \protect\hyperlink{week59}{``Week
59''} and \protect\hyperlink{week61}{``Week 61''}.

In 1845, two years after the birth of the quaternions, the visionary
William Clifford was born in Exeter, England. He only lived to the age
of 37: despite suffering from lung disease, he worked with incredible
intensity, and his closest friend wrote that ``He could not be induced,
or only with the utmost difficulty, to pay even moderate attention to
the cautions and observances which are commonly and aptly described as
`taking care of one's self'\,''. But in his short life, he pushed quite
far into the mathematics that would become the physics of the 20th
century. He studied the geometry of Riemann and prophetically envisioned
general relativity in 1876, in the following famous remarks:

\begin{quote}
"Riemann has shown that as there are different kinds of lines and
surfaces, so there are different kinds of space of three dimensions; and
that we can only find out by experience to which of these kinds the
space in which we live belongs. I hold in fact

\begin{enumerate}
\def\labelenumi{(\arabic{enumi})}
\item
  That small portions of space \emph{are} in fact of a nature analogous
  to little hills on a surface which is on the average flat; namely,
  that the ordinary laws of geometry are not valid for them.
\item
  That this property of being curved or distorted is continually being
  passed on from one portion of space to another after the manner of a
  wave.
\item
  That this variation of the curvature of space is what really happens
  in that phenomenon which we call the \emph{motion of matter}, whether
  ponderable or etherial.
\item
  That in the physical world nothing else takes place but this
  variation, subject (possibly) to the law of continuity.
\end{enumerate}
\end{quote}

He also substantially generalized Hamilton's quaternions, dropping the
condition that one have a division algebra, and focusing on the aspects
crucial to \(n\)-dimensional geometry. He obtained what we call the Clifford
algebras.

What's a Clifford algebra? Well, there are various flavors. But one of
the nicest --- let's call it \(\mathrm{C}_n\) --- is just the
associative algebra over the real numbers generated by \(n\)
anticommuting square roots of \(-1\). That is, we start with \(n\)
fellows called \[e_1, \ldots , e_n\] and form all formal products of
them, including the empty product, which we call \(1\). Then we form all
real linear combinations of these products, and then we impose the
relations \[
  \begin{aligned}
    e_i^2 &= -1
  \\e_ie_j &= -e_je_i.
  \end{aligned}
\] What are these algebras like? Well, \(C_0\) is just the real numbers,
since none of these \(e_i\)'s are thrown into the stew. \(C_1\) has one
square root of \(-1\), so it is just the complex numbers. \(C_2\) has
two square roots of \(-1\), \(e_1\) and \(e_2\), with
\[e_1 e_2 = - e_2 e_1.\] Thus \(C_2\) is just the quaternions, with
\(e_1\), \(e_2\), and \(e_1 e_2\) corresponding to Hamilton's \(i\),
\(j\), and \(k\).

How about the \(\mathrm{C}_n\) for larger values of \(n\)? Well, here is
a little table up to \(n = 8\):

\begin{longtable}[]{@{}rl@{}}
\toprule
\endhead
\(C_0\) & \(\mathbb{R}\)\tabularnewline
\(C_1\) & \(\mathbb{C}\)\tabularnewline
\(C_2\) & \(\mathbb{H}\)\tabularnewline
\(C_3\) & \(\mathbb{H}+\mathbb{H}\)\tabularnewline
\(C_4\) & \(\mathbb{H}(2)\)\tabularnewline
\(C_5\) & \(\mathbb{C}(4)\)\tabularnewline
\(C_6\) & \(\mathbb{R}(8)\)\tabularnewline
\(C_7\) & \(\mathbb{R}(8)+\mathbb{R}(8)\)\tabularnewline
\(C_8\) & \(\mathbb{R}(16)\)\tabularnewline
\bottomrule
\end{longtable}

What do these entries mean? Well, \(\mathbb{R}(n)\) means the
\(n\times n\) matrices with real entries. Similarly, \(\mathbb{C}(n)\)
means the \(n\times n\) complex matrices, and \(\mathbb{H}(n)\) means
the \(n\times n\) quaternionic matrices. All these become algebras with
the usual matrix addition and matrix multiplication. Finally, if \(A\)
is an algebra, \(A + A\) means the algebra consisting of pairs of guys
in \(A\), with the obvious rules for addition and multiplication: \[
  \begin{aligned}
    (a, a') + (b, b') &= (a + b, a' + b')
  \\(a, a') (b, b') &= (ab, a'b')
  \end{aligned}
\]

You might enjoy checking some of these descriptions of the Clifford
algebras \(\mathrm{C}_n\) for \(n\) up to 8. You have to find the
``isomorphism'' --- the correspondence between the Clifford algebra and
the algebra I claim is really the same. This gets pretty tricky when
\(n\) gets big.

How about \(n\) larger than 8? Well, here a remarkable fact comes into
play. Clifford algebras display a certain sort of ``period 8''
phenomenon. Namely, \(C_{n+8}\) consists of \(16\times 16\) matrices
with entries in \(\mathrm{C}_n\)! For a proof you might try

\begin{enumerate}
\def\labelenumi{\arabic{enumi})}
\setcounter{enumi}{1}
\tightlist
\item
  H.\ Blaine Lawson, Jr.~and Marie-Louise Michelson, \emph{Spin
  Geometry}, Princeton U. Press, Princeton, 1989.
\end{enumerate}

or

\begin{enumerate}
\def\labelenumi{\arabic{enumi})}
\setcounter{enumi}{2}
\tightlist
\item
  Dale Husemoller, \emph{Fibre Bundles}, Springer, Berlin, 1994.
\end{enumerate}

These books also describe some of the amazing consequences of this
periodicity phenomenon. The topology of \(n\)-dimensional manifolds is
very similar to the topology of \((n+8)\)-dimensional manifolds in some
subtle but important ways! I should describe this ``Bott periodicity''
sometime, but for now let me leave it as a tantalizing mystery.

I will also have to take a rain check when it comes to describing the
importance of Clifford algebras in physics... let me simply note
that they are crucial for understanding spin-\(1/2\) particles. I talked
a bit about this in \protect\hyperlink{week61}{``Week 61''}.

The ``Spin Geometry'' book goes into a lot of detail on Clifford
algebras, spinors, the Dirac equation and more. The ``Fibre Bundles''
book concentrates on the branch of topology called K-theory, and uses
this together with Clifford algebras to tackle various subtle questions,
such as how many linearly independent vector fields you can find on a
sphere.

\begin{enumerate}
\def\labelenumi{\arabic{enumi})}
\setcounter{enumi}{3}
\tightlist
\item
  Ralph L.\ Cohen, John D.\ S.\ Jones, and Graeme Segal, ``Morse theory
  and classifying spaces'', preprint as of Sept.~13, 1991.
\end{enumerate}

This is a nice way to think about what's really going on with Morse
theory. In Morse theory we study the topology of a compact Riemannian
manifold by putting a ``Morse function'' on it: a real-valued smooth
function with only nondegenerate critical points. The gradient of this
function defines a vector field and we use the way points flow along
this vector field to chop the manifold up into convenient pieces or
``cells''. A while back, Witten discovered, or rediscovered, a very cute
way to compute a topological invariant called the ``homology'' of the
manifold using Morse theory. (I've heard that this was previously known
and then largely forgotten.)

Here the authors refine this construction. They cook up a category
\(\mathcal{C}\) from the Morse function: the objects of \(\mathcal{C}\)
are critical points of the Morse function, and the morphisms are
piecewise gradient flow lines. This is a topological category, meaning
that for any pair of objects \(x\) and \(y\) the morphisms in
\(\operatorname{Hom}(x,y)\) form a topological space, and composition is
a continuous map. There is a standard recipe to construct the
``classifying space'' of any topological category, invented by Segal in
the following paper:

\begin{enumerate}
\def\labelenumi{\arabic{enumi})}
\setcounter{enumi}{4}
\tightlist
\item
  Graeme B. Segal, ``Classifying spaces and spectral sequences'',
  \emph{Pub. IHES} \textbf{34} (1968), 105--112.
\end{enumerate}

\noindent
I described classifying spaces for discrete groups in
\protect\hyperlink{week70}{``Week 70''}, and the more general case of
discrete groupoids in \protect\hyperlink{week75}{``Week 75''}. The
construction for topological categories is similar: we make a big space
by sticking in one point for each object, one edge for each morphism,
one triangle for each composable pair of morphisms: \[
  \begin{tikzpicture}
    \node (x) at (0,0) {$x$};
    \node (y) at (1,1.7) {$y$};
    \node (z) at (2,0) {$z$};
    \draw[thick] (x) to node[fill=white]{$f$} (y);
    \draw[thick] (x) to node[fill=white]{$gf$} (z);
    \draw[thick] (y) to node[fill=white]{$g$} (z);
    \node at (4,0.8) {$
      \begin{aligned}
        f&\colon x\to y
      \\g&\colon y\to z
      \\gf&\colon x\to z
      \end{aligned}
    $};
  \end{tikzpicture}
\] and so on. The only new trick is to make sure this space gets a
topology in the right way using the topologies on the spaces
\(\operatorname{Hom}(x,y)\).

Anyway, if we form this classifying space from the topological category
\(\mathcal{C}\) coming from the Morse function on our manifold \(M\), we
get a space homotopic to \(M\)! In other words, for many topological
purposes the category \(\mathcal{C}\) is just as good as the manifold
\(M\) itself.

\begin{enumerate}
\def\labelenumi{\arabic{enumi})}
\setcounter{enumi}{5}
\item
  Ross Street, ``Descent theory'', preprint of talks given at
  Oberwolfach, Sept.~17--23, 1995.

  Ross Street, ``Fusion operators and cocycloids in monoidal
  categories'', \emph{Appl. Cat. Str.} \textbf{6} (1998), 177--191.
\end{enumerate}
\noindent
Street is one of the gurus of \(n\)-category theory, which he notes
``might be called post-modern algebra (or even `post-modern mathematics'
since geometry and algebra are handled equally well by higher
categories).'' His paper ``Descent theory'' serves as a rapid
introduction to \(n\)-categories. But the real point of the paper is to
explain the role \(n\)-categories play in cohomology theory: in particular,
how to do cohomology with coefficients in an \(\omega\)-category!

\begin{enumerate}
\def\labelenumi{\arabic{enumi})}
\setcounter{enumi}{6}
\item
  Viqar Husain, ``Intersecting-loop solutions of the hamiltonian
  constraint of quantum general relativity'', \emph{Nucl. Phys. B}
  \textbf{313} (1989), 711--724.

  Viqar Husain and Karel V. Kuchar, ``General covariance, new variables,
  and dynamics without dynamics'', \emph{Phys. Rev.~D} \textbf{42}
  (1990), 4070--4077.

  Viqar Husain, ``Einstein's equations and the chiral model'', 
  \emph{Phys.\ Rev.\ D} \textbf{53} (1996), 4327.  Also available as
  \href{https://arxiv.org/abs/gr-qc/9602050}{\texttt{gr-qc/9602050}}.
\end{enumerate}
\noindent
Viqar is one of the excellent younger folks at the Center for
Gravitational Physics and Geometry at Penn State; I only had a bit of
time to speak with him during my last visit there, but I got some of his
papers. The first paper is from the good old days when folks were just
beginning to find explicit solutions of the constraints of quantum
gravity using the loop representation --- it's still worth reading! The
second introduced a field theory now called the Husain--Kuchar model,
which has the curious property of resembling gravity without the
dynamics. The third studies \(4\)-dimensional general relativity
assuming there are two commuting spacelike Killing vector fields. Here
he reduces the theory to a \(2\)-dimensional theory which appears to be
completely integrable --- though it has not been proved to be so in the
sense of admitting a complete set of Poisson-commuting conserved
quantities.

\begin{enumerate}
\def\labelenumi{\arabic{enumi})}
\setcounter{enumi}{7}
\tightlist
\item
  \emph{The Interface of Knots and Physics}, ed.~Louis H. Kauffman,
  Proc. Symp. Appl. Math. \textbf{51}, AMS,
  Providence, 1996.
\end{enumerate}
\noindent
This slim volume contains the proceedings of an AMS ``short course'' on
knots and physics held in San Francisco in January 1995, namely:

\begin{itemize}
\tightlist
\item
  Louis H. Kauffman, ``Knots and statistical mechanics''
\item
  Ruth J. Lawrence, ``An introduction to topological field theory''
\item
  Dror Bar-Natan, ``Vassiliev and quantum invariants of braids''
\item
  Samuel J. Lomonaco, ``The modern legacies of Thomson's atomic vortex
  theory in classical electrodynamics''
\item
  John C. Baez, ``Spin networks in nonperturbative quantum gravity''
\end{itemize}

\begin{center}\rule{0.5\linewidth}{0.5pt}\end{center}

\begin{quote}
\emph{William Kingon Clifford \\
Born May 4th, 1845 \\
Died March 3rd, 1879 \\
\\
I was not, and was conceived \\
I loved, and did a little work \\
I am not, and grieve not. \\
\\
And \\
\\
Lucy, his wife \\
Died April 21st, 1929 \\
\\
Oh, two such silver currents when they join \\
Do glorify the banks that bound them in.}
\vskip 1em
--- William Clifford's tomb
\end{quote}

\hypertarget{week83}{%
\section{June 10, 1996}\label{week83}}

I'll return to the Tale of \(n\)-Categories this week, and continue to
explain the mysteries of duals and inverses. But first let me describe
two new papers by Connes.

\begin{enumerate}
\def\labelenumi{\arabic{enumi})}
\item
  Alain Connes, ``Gravity coupled with matter and the foundation of
  non-commutative geometry'', available as
  \href{https://arxiv.org/abs/hep-th/9603053}{\texttt{hep-th/9603053}}.

  Ali H.\ Chamseddine and Alain Connes, ``The spectral action
  principle'', available as
  \href{https://arxiv.org/abs/hep-th/9606001}{\texttt{hep-th/9606001}}.
\end{enumerate}
\noindent
The second paper here fills in details that are missing from the first.
Hopefully lots of you know that Connes is the wizard of operator theory
who turned to inventing a new branch of geometry, ``noncommutative
geometry''. The idea of algebraic geometry is that we can study a space
by studying the functions on that space --- which typically form some
kind of commutative algebra. If we let the algebra become
noncommutative, it is no longer functions on some space, but we can
pretend it is nonetheless, and do geometry by analogy with the
commutative case. This is very much based on the philosophy of quantum
mechanics, where the observables form a noncommutative algebra, yet are
analogous to the commutative algebras of observables of classical
mechanics, these commutative algebras consisting simply of functions on
the classical space states.

In quantum mechanics, the failure of two observables to commute implies
that they cannot always be simultaneously measured with arbitrary
accuracy; there is a very precise mathematical statement of Heisenberg's
uncertainty principle that makes this quantitative. We can thus think of
noncommutative geometry as ``quantum geometry'', geometry where the
uncertainty principle of quantum mechanics has infected the very notion
of space itself! In noncommutative geometry it impossible to
simultaneously measure all the coordinates of a point with arbitrary
accuracy, because they do not commute!

For the definitive introduction to noncommutative geometry, see Connes'
book ``Noncommutative Geometry'', reviewed in
\protect\hyperlink{week39}{``Week 39''}. Already in this book Connes,
working with Lott, was beginning to explore the idea that the geometry
of our physical universe is noncommutative. Actually, they used ideas
from noncommutative geometry to study a weird kind of commutative
geometry in which spacetime is ``two-sheeted''---two copies of standard
\(4\)-dimensional spacetime, very close together. In normal geometry it
doesn't even make sense to speak of two separate copies of spacetime
being ``close together'', since there is no way to get from one to the
other! Tricks from noncommutative geometry allow it to make sense. They
found something amazing: if you do
\(\mathrm{U}(1)\times \mathrm{SU}(2)\) Yang--Mills theory on this
spacetime, you get the Higgs particle for free!

Sorry for the jargon. What it means is this: in the Standard Model of
particle physics we describe the electromagnetic force and the weak
nuclear force in a unified way using a theory called
``\(\mathrm{U}(1)\times \mathrm{SU}(2)\) Yang--Mills theory'', but then
we postulate an extra particle, the Higgs particle, which has the effect
of making the electromagnetic force work quite differently from the weak
force. We say it ``breaks the symmetry'' between the two forces. It has
not yet been observed, though particle physicists hope to see it (or
not!) in experiments coming up fairly soon. It is a rather puzzling, ad
hoc element of the Standard Model. The amazing thing about the
Connes-Lott model is that it arises in a natural way from the fact that
spacetime has two sheets.

Connes and Lott also studied the strong force, but now Connes has
introduced gravity into his model. I haven't had time to absorb this new
work yet, so let me simply say what his current model of spacetime is,
and list some of the concrete predictions the new theory makes. His
spacetime is the noncommutative algebra consisting of smooth functions
on good old \(4\)-dimensional Minkowski spacetime, taking values in the
algebra \(A\) given by the direct sum
\[A = \mathbb{C} + \mathbb{H} + M_3(\mathbb{C})\] where \(\mathbb{C}\)
is the complex numbers, \(\mathbb{H}\) is the quaternions, and
\(M_3(\mathbb{C})\) is the \(3\times3\) complex matrices. (Exercise:
redo Connes' model, replacing \(M_3(\mathbb{C})\) with the octonions.
Hint: develop nonassociative geometry and use Geoffrey Dixon's theory
relating the electromagnetic, weak, and strong forces to the complex
numbers, quaternions, and octonions, respectively. See
\protect\hyperlink{week59}{``Week 59''} for references to Dixon's work,
and an explanation of quaternions and octonions.)

The Chamseddine-Connes model predicts that the sine squared of the
Weinberg angle --- an important constant in the theory of the
electroweak force --- is between \(.206\) and \(.210\). Unfortunately
this disagrees with the experimental value of \(.2325\), but it's sort
of surprising that they can derive something this close, since in the
Standard Model the Weinberg angle is just an arbitrary parameter. They also
derive a Higgs mass of 160--180 GeV, and expect accuracy comparable to
their prediction of the Weinberg angle (about 10\%).

Well worth pondering!

\begin{center}\rule{0.5\linewidth}{0.5pt}\end{center}

\hypertarget{week83_tale}{
There is an interesting analogy between the dual of a vector space and
the inverse of a number which I would like to explain now.} I assume you
know that multiplying numbers is a lot like tensoring vector spaces. For
example, just as multiplication distributes over addition, tensoring
distributes over direct sums. Also, just as there is a number called
\(1\) which is the unit for multiplication, there is a \(1\)-dimensional
vector space, the ground field itself, which is the unit for tensoring.
Let me take the unusual liberty of writing tensor products by
juxtaposition, so that \(xy\) is the tensor product of the vector space
\(x\) and the vector space \(y\), and let me call the \(1\)-dimensional
vector space that's the unit for tensoring simply ``\(1\)''.

Now, if a number \(x\) has an inverse \(y\), we have \[yx = 1\] and
\[1 = xy.\] Similarly, if a vector space \(x\) has a dual \(y\), we have
linear maps \[e\colon yx\to 1\] and \[i\colon 1\to xy\] What are these
linear maps? Well, the whole point of the dual vector space y is that a
vector in \(y\) is a linear functional from \(x\) to \(1\). This ``dual
pairing'' between vectors in \(y\) and those in \(x\) defines a linear
map \(e\colon yx\to 1\), often called the ``counit''. On the other hand,
the space \(xy\) can be thought of as the space of linear
transformations of \(x\). The linear map \(i\colon 1\to xy\) sends any
scalar (i.e., any vector in \(1\)) to the corresponding scalar multiple
of the identity transformation of \(x\).

So we see that dual vector spaces are a bit like inverse numbers, except
that we don't have \(yx = 1\) and \(1 = xy\), and we don't even have
that \(yx\) is \emph{isomorphic} to \(1\) and \(1\) is \emph{isomorphic}
to \(xy\). We just have some maps going from \(yx\) to \(1\), and from
\(1\) to \(xy\).

These maps satisfy two equations, though. Here's the first. We start
with \(x\), use the obvious isomorphism to map to \(1x\), then use
\(i\colon 1\to xy\) to map this to \(xyx\), then use \(e\colon yx\to 1\)
to map this to \(x1\), and then use the other obvious isomorphism to map
back to \(x\). This composite of maps should be the identity on \(x\).
What this says is that the identity linear transformation of \(x\)
really acts as the identity!

Stealing a trick from \protect\hyperlink{week79}{``Week 79''}, we can
draw this as follows. Draw the counit \(e\colon yx\to 1\) as follows: \[
  \begin{tikzpicture}
    \begin{knot}
      \strand[thick] (0,0.5)
        to (0,0)
        to [out=down,in=down,looseness=2] (1,0)
        to (1,0.5);
    \end{knot}
    \node[fill=white] at (0,0) {$y$};
    \node[fill=white] at (1,0) {$x$};
    \node[label=below:{$e$}] at (0.5,-0.6) {$\bullet$};
  \end{tikzpicture}
\] and draw the unit \(i\colon1\to xy\) as follows: \[
  \begin{tikzpicture}
    \begin{knot}
      \strand[thick] (0,-0.5)
        to (0,0)
        to [out=up,in=up,looseness=2] (1,0)
        to (1,-0.5);
    \end{knot}
    \node[fill=white] at (0,0) {$x$};
    \node[fill=white] at (1,0) {$y$};
    \node[label=above:{$i$}] at (0.5,0.57) {$\bullet$};
  \end{tikzpicture}
\] Then the above equation says that \[
  \begin{tikzpicture}
    \begin{knot}
      \strand[thick] (0,0)
      to (0,1)
      to [out=up,in=up,looseness=2] (1,1)
      to [out=down,in=down,looseness=2] (2,1)
      to (2,2);
    \end{knot}
    \node[fill=white] at (0,0.25) {$x$};
    \node[fill=white] at (2,1.75) {$x$};
    \node[fill=white] at (0,1) {$x$};
    \node[fill=white] at (1,1) {$y$};
    \node[fill=white] at (2,1) {$x$};
    \node at (3,1) {$=$};
    \begin{scope}[shift={(4,0)}]
      \begin{knot}
        \strand[thick] (0,0) to (0,2);
      \end{knot}
      \node[fill=white] at (0,0.25) {$x$};
      \node[fill=white] at (0,1.75) {$x$};
    \end{scope}
  \end{tikzpicture}
\] Here the left side, which we read from top to bottom, corresponds to
the composite \(x\to 1x\to xyx\to x1\to x\). (The factors of \(1\) are
invisible in the picture, since they don't do much.) The left side
corresponds to the identity map \(x\to x\).

The second equation goes like this. We start with \(y\), use the obvious
isomorphism to map to \(y1\), then use the unit to map this to \(yxy\),
then use the counit to map this to \(1y\), and then use the other
obvious isomorphism to map back to \(y\). This composite should be the
identity on \(y\). What this says is that the identity linear
transformation of \(x\) also acts dually as the identity on \(y\)! We
can draw this as follows: \[
  \begin{tikzpicture}
    \begin{scope}[xscale=-1,shift={(-2,0)}]
      \begin{knot}
        \strand[thick] (0,0)
        to (0,1)
        to [out=up,in=up,looseness=2] (1,1)
        to [out=down,in=down,looseness=2] (2,1)
        to (2,2);
      \end{knot}
      \node[fill=white] at (0,0.5) {$y$};
      \node[fill=white] at (2,1.5) {$y$};
      \node[fill=white] at (2,1) {$y$};
      \node[fill=white] at (1,1) {$x$};
      \node[fill=white] at (0,1) {$y$};
    \end{scope}
    \node at (3,1) {$=$};
    \begin{scope}[shift={(4,0)}]
      \begin{knot}
        \strand[thick] (0,0) to (0,2);
      \end{knot}
      \node[fill=white] at (0,0.5) {$y$};
      \node[fill=white] at (0,1.7) {$y$};
    \end{scope}
  \end{tikzpicture}
\] If you now steal a peek at \protect\hyperlink{week79_tale}{``Week 79''},
you'll see that these two equations are just the same equations used to
define adjoint functors in category theory! What's going on? Well, dual
vector spaces are analogous to adjoint functors, clearly. But more
deeply, what we have is an analogy between duals in any category with
tensor products --- or ``monoidal category'' --- and adjoints in any
\(2\)-category.

What's a monoidal category, exactly? Roughly it's a category with some
sort of ``tensor product'' and ``unit object''. But we can precisely
define the so-called ``strict'' monoidal categories as follows: they are
simply \(2\)-categories with one object. (Turn to
\protect\hyperlink{week80}{``Week 80''} for a definition of
\(2\)-categories.) A \(2\)-category has objects, morphisms, and
\(2\)-morphisms, but if there is only one object, we can do the
following relabelling trick: \[
  \begin{aligned}
    \text{2-morphisms} &\mapsto \text{morphisms}
  \\\text{morphisms} &\mapsto \text{objects}
  \\\text{objects} &\mapsto 
  \end{aligned}
\] Namely, we can forget about the object, call the morphisms
``objects'', and call the \(2\)-morphisms ``morphisms''. But since all
the new ``objects'' were really morphisms from the original single
object to itself, they can all be composed, or ``tensored''. That's why
we get a category with ``tensor product'', and similarly, a ``unit
object''.

So, just as a category with one object is just a monoid, a
\(2\)-category with one object is a monoidal category! This is one
instance of a trick that I sketched many more cases of in
\protect\hyperlink{week74_tale}{``Week 74''}.

Now, in \protect\hyperlink{week79_tale}{``Week 79''} I defined left and right
adjoints of functors between categories. Here the only thing I really
needed about category theory was that \(\mathsf{Cat}\) is a
\(2\)-category with categories as its objects, functors as its
morphisms, and natural transformations as its 2-morphisms. So we can
define left and right adjoints of morphisms in any \(2\)-category by
analogy as follows:

Suppose \(a\) and \(b\) are objects in a \(2\)-category. Then we say
that the morphism \[L\colon a\to b\] is a ``left adjoint'' of the
morphism \[R\colon b\to a\] (and \(R\) is a ``right adjoint'' of \(L\))
if there are \(2\)-morphisms \[
  \begin{aligned}
    e&\colon RL\Rightarrow 1_b
  \\i&\colon 1_a\Rightarrow LR
  \end{aligned}
\] satisfying two magic equations. If we draw \(e\) and \(i\) as we did
above, \[
  \begin{tikzpicture}
    \begin{knot}
      \strand[thick] (0,0.5)
        to (0,0)
        to [out=down,in=down,looseness=2] (1,0)
        to (1,0.5);
    \end{knot}
    \node[fill=white] at (0,0) {$y$};
    \node[fill=white] at (1,0) {$x$};
    \node[label=below:{$e$}] at (0.5,-0.6) {$\bullet$};
  \end{tikzpicture}
  \qquad
  \begin{tikzpicture}
    \begin{knot}
      \strand[thick] (0,-0.5)
        to (0,0)
        to [out=up,in=up,looseness=2] (1,0)
        to (1,-0.5);
    \end{knot}
    \node[fill=white] at (0,0) {$x$};
    \node[fill=white] at (1,0) {$y$};
    \node[label=above:{$i$}] at (0.5,0.57) {$\bullet$};
  \end{tikzpicture}
\] then the two magic equations are \[
  \begin{tikzpicture}
    \begin{knot}
      \strand[thick] (0,0)
      to (0,1)
      to [out=up,in=up,looseness=2] (1,1)
      to [out=down,in=down,looseness=2] (2,1)
      to (2,2);
    \end{knot}
    \node[fill=white] at (0,0.25) {$L$};
    \node[fill=white] at (2,1.75) {$L$};
    \node[fill=white] at (0,1) {$L$};
    \node[fill=white] at (1,1) {$R$};
    \node[fill=white] at (2,1) {$L$};
    \node at (3,1) {$=$};
    \begin{scope}[shift={(4,0)}]
      \begin{knot}
        \strand[thick] (0,0) to (0,2);
      \end{knot}
      \node[fill=white] at (0,0.25) {$L$};
      \node[fill=white] at (0,1.75) {$L$};
    \end{scope}
  \end{tikzpicture}
\] and \[
  \begin{tikzpicture}
    \begin{scope}[xscale=-1,shift={(-2,0)}]
      \begin{knot}
        \strand[thick] (0,0)
        to (0,1)
        to [out=up,in=up,looseness=2] (1,1)
        to [out=down,in=down,looseness=2] (2,1)
        to (2,2);
      \end{knot}
      \node[fill=white] at (0,0.5) {$R$};
      \node[fill=white] at (2,1.5) {$R$};
      \node[fill=white] at (2,1) {$R$};
      \node[fill=white] at (1,1) {$L$};
      \node[fill=white] at (0,1) {$R$};
    \end{scope}
    \node at (3,1) {$=$};
    \begin{scope}[shift={(4,0)}]
      \begin{knot}
        \strand[thick] (0,0) to (0,2);
      \end{knot}
      \node[fill=white] at (0,0.5) {$R$};
      \node[fill=white] at (0,1.7) {$R$};
    \end{scope}
  \end{tikzpicture}
\]

Alternatively, we can state these equations using the \(2\)-categorical
notation described in \protect\hyperlink{week80_tale}{``Week 80''}, by saying
that the following vertical composites of \(2\)-morphisms are identity
morphisms:
\[L = 1_aL\xRightarrow{i\cdot1_L}LRL\xRightarrow{1_L\cdot e}L1_a = L\]
and
\[R = R1_a\xRightarrow{1_R\cdot i}RLR\xRightarrow{e\cdot1_R}1_bR = R\]
where \(\cdot\) denotes the horizontal composite. If you look at these,
and compare them to the graphical notation above, you'll see they are
really saying the same thing.

The punchline is, when our \(2\)-category has one object, we can think
of it as a monoidal category, and then these equations are the
definition of ``duals'' --- one example being our earlier definition of
dual vector spaces in the monoidal category Vect of vector spaces!

So adjoint functors and dual vector spaces are both instances of the
general notion of adjoint \(1\)-morphisms in a \(2\)-category.
Adjointness is a very basic concept.

I hope all that made some sense.

If this category theory stuff seems confusing, maybe you should read a
3-volume book about it! I can see you smiling, but seriously, I find the
following reference very useful (despite a certain number of annoying
errors). You can find a lot of good stuff about adjoint functors,
monoidal categories, and much much more in here:

\begin{enumerate}
\def\labelenumi{\arabic{enumi})}
\setcounter{enumi}{1}
\tightlist
\item
  Francis Borceux, \emph{Handbook of Categorical Algebra}, Cambridge U.
  Press 1994. \emph{Volume 1: Basic Category Theory}. \emph{Volume 2:
  Categories and Structure}. \emph{Volume 3: Categories of Sheaves}.
\end{enumerate}

To continue reading the ``Tale of \(n\)-Categories'', see
\protect\hyperlink{week84_tale}{``Week 84''}.

\hypertarget{week84}{%
\section{June 27, 1996}\label{week84}}

While I try to limit myself to mathematical physics in This Week's
Finds, I can't always keep it from spilling over into other
subjects... so if you're not interested in computers, just skip
down to reference 8 below. A while back on \texttt{sci.physics} Matt
McIrvin pointed out that the closest thing we have to the computer of
old science fiction --- the underground behemoth attended by technicians
in white lab coats that can answer any question you type in --- is
AltaVista. I agree wholeheartedly.

In case you are a few months or years behind on the technological front,
let me explain: these days there is a vast amount of material available
on the World-Wide Web, so that the problem has become one of locating
what you are interested in. You can do this with programs known as
search engines. There are lots of search engines, but my favorite these
days is AltaVista, which is run by DEC, and seems to be especially
comprehensive. So these days if you want to know, say, the meaning of
life, you can just go to

\begin{enumerate}
\def\labelenumi{\arabic{enumi})}
\tightlist
\item
  AltaVista, \texttt{http://www.altavista.digital.com/}
\end{enumerate}

type in ``meaning of life'', and see what everyone has written about it.
You'll be none the wiser, of course, but that's how it always worked in
those old science fiction stories, too.

The intelligence of AltaVista is of course far less than that of a fruit
fly. But Matt's comment made me think about how now, as soon as we
develop even a rudimentary form of artificial intelligence, it will
immediately have access to vast reams of information on the Web...
and may start doing some surprising things.

An example of what I'm talking about is the CYC project, Doug Lenat's
\$35 million project, begun in 1984, to write a program with common
sense. In fact, the project plans to set CYC loose on the web once it
knows enough to learn something from it.

\begin{enumerate}
\def\labelenumi{\arabic{enumi})}
\setcounter{enumi}{1}
\tightlist
\item
  CYC project homepage, \texttt{http://www.cyc.com/}
\end{enumerate}

The idea behind CYC is to encode ``common sense'' as about half a
million rules of thumb, declarative sentences which CYC can use to
generate inferences. To have any chance of success, these rules of thumb
must be organized and manipulated very carefully. One key aspect of this
is CYC's ontology --- the framework that lets it know, for example, that
you can eat 4 sandwiches, but not 4 colors or the number 4. Most of the
CYC code is proprietary, but the ontology will be made public in July of
this year, they say. One can already read about aspects of it in

\begin{enumerate}
\def\labelenumi{\arabic{enumi})}
\setcounter{enumi}{2}
\tightlist
\item
  Douglas B. Lenat and R.V. Guha, \emph{Building Large Knowledge-Based
  Systems: Representation and Inference in the Cyc Project},
  Addison-Wesley, Reading, Mass., 1990.
\end{enumerate}

My network of spies informs me that many hackers are rather suspicious
of CYC. For an interesting and somewhat critical account of CYC at one
stage of its development, see

\begin{enumerate}
\def\labelenumi{\arabic{enumi})}
\setcounter{enumi}{3}
\tightlist
\item
  Vaughan Pratt, ``CYC Report'',
  \texttt{http://boole.stanford.edu/pub/cyc.report}
\end{enumerate}

Turning to something that's already very practical, I was very pleased
when I found one could use AltaVista to do ``backlinks''. Certainly the
World-Wide Web is in part inspired by Ted Nelson's visionary but
ill-starred Xanadu project.

\begin{enumerate}
\def\labelenumi{\arabic{enumi})}
\setcounter{enumi}{4}
\tightlist
\item
  Project Xanadu, \texttt{http://xanadu.net/the.project}
\end{enumerate}

Backlinking is one of the most tricky parts of Ted Nelson's vision, one
often declared infeasible, but one upon which he has always insisted.
Basically, the idea is that you should always be able to find all the
documents pointing \emph{to} a given document, as well as those to which
it points. This allows \textbf{commentary} or \textbf{annotation}: if
you read something, you can read what other people have written about
it. My spies inform me that the World-Wide Web Committee is moving in
this direction, but it is exciting that one can already do ``backlinks
browsing'' with the help of a program written by Ted Kaehler:

\begin{enumerate}
\def\labelenumi{\arabic{enumi})}
\setcounter{enumi}{5}
\tightlist
\item
  Ted Kaehler's backlinks browser,
  \texttt{http://www.foresight.org/backlinks1.3.1.html}
\end{enumerate}

Go to this page at the start of your browsing session. Follow the
directions and let it create a new Netscape window for you to browse in.
Whenever you want backlinks, click in the original page, and click
``Links to Other Page''. This launches an AltaVista search for links to
the page you were just looking at.

It works quite nicely. I hope you try it, because with backlinking the
Web will become a much more interesting and useful place, and the more
people who know about it, the sooner it will catch on. For more
discussion of backlinking, see

\begin{enumerate}
\def\labelenumi{\arabic{enumi})}
\setcounter{enumi}{6}
\item
  Backlinking news at the Foresight Institute,
  \texttt{http://www.foresight.org/backlinks.news.html}

  Robin Hanson's ideas on backlinking,
  \texttt{http://www.hss.caltech.edu/\textasciitilde{}hanson/findcritics.html}
\end{enumerate}

I thank my best buddy Bruce Smith for telling me about CYC, Project
Xanadu, and Ted Kaehler's backlinks browser.

Now let me turn to some mathematics and physics.

\begin{enumerate}
\def\labelenumi{\arabic{enumi})}
\setcounter{enumi}{7}
\item
  Francesco Fucito, Maurizio Martellini and Mauro Zeni, ``The BF
  formalism for QCD and quark confinement'', available as
  \href{https://arxiv.org/abs/hep-th/9605018}{\texttt{hep-th/9605018}}.
\item
  Ioannis Tsohantjis, Alex C Kalloniatis, and Peter D. Jarvis, ``Chord
  diagrams and BPHZ subtractions'', available as
  \href{https://arxiv.org/abs/hep-th/9604191}{\texttt{hep-th/9604191}}.
\end{enumerate}

These two papers both treat interesting relationships between topology
and quantum field theory --- not the ``topological quantum field
theory'' beloved of effete mathematicians such as myself, but the
pedestrian sort of quantum field theory that ordinary working physicists
use to study particle physics. So we are seeing an interesting
cross-fertilization here: first quantum field theory got applied to
topology, and now the resulting ideas are getting applied back to
particle physics.

Why don't we see quarks roaming the streets freely at night? Because
they are confined! Confined to the hadrons in which they reside, that
is. We mainly see two sorts of hadrons: baryons made of three quarks,
like the proton and neutron, and mesons made of a quark and an
antiquark, like the pion or kaon. Why are the quarks confined in
hadrons? Well, roughly it's because if you grab a quark inside a hadron
and try to pull it out, the force needed to pull it doesn't decrease as
you pull it farther out; instead, it remains about constant. Thus the
energy grows linearly with the distance, and eventually you have put
enough energy into the hadron to create another quark-antiquark pair,
and \emph{pop} --- you find you are holding not a single quark but a
quark together with a newly born antiquark, that is, a meson! What's
left is a hadron with a newly born quark as the replacement for the one
you tried to pull out!

That's a pretty heuristic description. In fact, particle physicists do
not usually grab hadrons and try to wrest quarks from them with their
bare hands. Instead they smash hadrons and other particles at each other
and study the debris. But as a rough sketch of the theory of quark
confinement, the above description is not \emph{completely} silly.

There are various interesting things left to do, though. One is to try
to get those quarks out by means of sneaky tricks. The only way known is
to \emph{heat} a bunch of hadrons to ridiculously high temperatures,
preferably at ridiculously high pressures. I'm talking temperatures like
2 trillion degrees, and densities comparable to that of nuclear matter!
This should yield a ``quark-gluon plasma'' in which quarks can zip
around freely at enormous energies. That's what the folks at the
Relativistic Heavy Ion Collider are doing --- see
\protect\hyperlink{week76_tale}{``Week 76''} for more.

This should certainly keep the experimentalists entertained. On the
other hand, theorists can have lots of fun trying to understand more
deeply why quarks are confined. We'd like best to derive confinement in
some fairly clear and fairly rigorous way from quantum chromodynamics,
or QCD --- our current theory of the strong force, the force that binds
the quarks together. Unfortunately, mathematical physicists are still
struggling to formulate QCD in a rigorous way, so they can't yet turn to
the exciting challenge of proving that confinement follows from QCD. And
we certainly don't expect any simple way to ``exactly solve'' QCD, since
it is complicated and highly nonlinear. So what some people do instead
is computer simulations of QCD, in which they approximate spacetime by a
lattice and do a lot of number-crunching. They usually use a fairly
measly-sounding grid of something like \(16 \times 16 \times 16 \times 16\) 
sites or so, since currently calculations take too long when the lattice 
gets much bigger than that.

Numerical calculations like these have a lot of potential. In
\protect\hyperlink{week68}{``Week 68''}, for example, I talked about how
people found numerical evidence for the existence of a ``glueball'' ---
a hadron made of no quarks, just gluons, the gluon being the particle
that carries the strong force. This glueball candidate seems to match
the features of an observed particle! Also, people have put a lot of
work into computing the masses of more familiar hadrons. So far I
believe they have concentrated on mesons, which are simpler. Eventually
we should in principle be able to calculate things like the mass of the
proton and neutron --- which would be really thrilling, I think.
Numerical calculations have also yielded a lot of numerical evidence
that QCD predicts confinement.

Still, one would very much like some conceptual explanation for
confinement, even if it's not quite rigorous. One way people try to
understand it is in terms of ``dual superconductivity''. In certain
superconductors, magnetic fields can only penetrate as long narrow tubes
of magnetic flux. (For example, this happens in neutron stars --- see
\protect\hyperlink{week37}{``Week 37''}.) Now, just as electromagnetism
consists of an ``electric'' part and a ``magnetic'' part, so does the
strong force. But the idea is that confinement is due to the
\emph{electric} part of the strong force only being able to penetrate
the vacuum in the form of long narrow tubes of field lines. The electric
and magnetic fields are ``dual'' to each other in a precise mathematical
sense, so this is referred to as dual superconductivity. Quarks have the
strong force version of electric charge --- called ``color'' --- and
when we try to pull two quarks apart, the strong electric field gets
pulled into a long tube. This is why the force remains constant rather
than diminishing as the distance between the quarks increases.

A while back, 't Hooft proposed an idea for studying confinement in
terms of dual superconductivity and certain ``order'' and ``disorder''
observables. It seems this idea has languished to some extent due to a
lack of necessary mathematical infrastructure. For the last couple of
years, Martellini has been suggesting to use ideas from topological
quantum field theory to serve this role. In particular, he suggested
treating Yang--Mills theory as a perturbation of \(BF\) theory, and
applying some of the ideas of Witten and Seiberg (who related
confinement in the supersymmetric generalization of Yang--Mills theory to
Donaldson theory). In the paper with Fucito and Zeni, they make some of
these ideas precise. There are still some potentially serious loose
ends, so I am very interested to hear the reaction of others working on
confinement.

I have not studied the paper of Tsohantjis, Kalloniatis, and Jarvis in
any detail, but people studying Vassiliev invariants might find it
interesting, since it claims to relate the renormalization theory of
\(\varphi^3\) theory to the mathematics of chord diagrams.

\begin{enumerate}
\def\labelenumi{\arabic{enumi})}
\setcounter{enumi}{9}
\tightlist
\item
  Masaki Kashiwara and Yoshihisa Saito, ``Geometric construction of
  crystal bases'',
  \href{https://arxiv.org/abs/q-alg/9606009}{\texttt{q-alg/9606009}}.
\end{enumerate}
\noindent
The ``canonical'' and ``crystal'' bases associated to quantum groups,
studied by Kashiwara, Lusztig, and others, are exciting to me because
they indicate that the quantum groups are just the tip of a still richer
structure. Whenever you see an algebraic structure with a basis in which
the structure constants are nonnegative integers, you should suspect
that you are really working with a category of some sort, but in
boiled-down or ``decategorified'' form.

Consider for example the representation ring \(R(G)\) of a group \(G\).
This is a ring whose elements are just the isomorphism classes of
finite-dimensional representations of \(G\). Addition in \(R(G)\)
corresponds to taking the direct sum of representations, while
multiplication corresponds to taking the tensor product. Thus for
example if \(x\) and \(y\) are irreducible representations of \(G\) ---
or ``irreps'' for short --- and \([x]\) and \([y]\) are the
corresponding basis elements of \(R(G)\), the product \([x][y]\) is equal to
a linear combination of the irreps appearing in \(x\otimes y\), with the
coefficients in the linear combination being the \emph{multiplicities}
with which the various irreps appear in \(x\otimes y\). These
coefficients are therefore nonnegative integers. They are an example of
what I'm calling ``structure constants''.

What's happening here is that the ring \(R(G)\) is serving as a
``decategorified'' version of the category \(\mathsf{Rep}(G)\) of
representations of the group \(G\). Almost everything about \(R(G)\) is
just a decategorified version of something about \(\mathsf{Rep}(G)\).
For example, \(R(G)\) is a monoid under multiplication, while
\(\mathsf{Rep}(G)\) is a monoidal category under tensor product.
\(R(G)\) is actually a commutative monoid, while \(\mathsf{Rep}(G)\) is
a symmetric monoidal category --- this being jargon for how the tensor
product is ``commutative'' up to a nice sort of isomorphism. In \(R(G)\)
we have addition, while in \(\mathsf{Rep}(G)\) we have direct sums,
which category theorists would call ``biproducts''. And so on. The
representation ring is a common tool in group theory, but a lot of the
reason for working with \(R(G)\) is simply that we don't know enough
about category theory and are too scared to work directly with
\(\mathsf{Rep}(G)\). There are also \emph{good} reasons for working with
\(R(G)\), basically \emph{because} it is simpler and contains less
information than \(\mathsf{Rep}(G)\).

We can imagine that if someone handed us a representation ring \(R(G)\)
we might eventually notice that it had a nice basis in which the
structure constants were nonnegative integers. And we might then realize
that lurking behind it was a category, \(\mathsf{Rep}(G)\). And then all
sorts of things about it would become clearer\ldots.

Something similar like this seems to be happening with quantum groups!
Ignoring a lot of important technical details, let me just say that
quantum groups turn out have nice bases in which the structure constants
are nonnegative integers, and the reason is that lurking behind the
quantum groups are certain categories. What sort of categories?
Categories of ``Lagrangian subvarieties of the cotangent bundles of
quiver varieties''. Yikes! I don't think I'll explain \emph{that}
mouthful! Let me just note that a quiver is itself a cute little
category that you cook up by taking a graph and thinking of the vertices
as objects and the edges as morphisms, like this:
\[\bullet\to\bullet\to\bullet\to\bullet\to\bullet\] If you do this to a
graph that's the Dynkin diagram of a Lie group --- see
\protect\hyperlink{week62}{``Week 62''} and the weeks following that ---
then the fun starts! Dynkin diagrams give Lie groups, but also quantum
groups, and now it turns out that they also give rise to certain
categories of which the quantum groups are decategoried, boiled-down
versions\ldots. I don't understand all this, but I certainly intend to,
because it's simply amazing how a world of complex symmetry emerges from
these Dynkin diagrams.

For more on this stuff try the paper by Crane and Frenkel referred to in
\protect\hyperlink{week38}{``Week 38''} and
\protect\hyperlink{week50}{``Week 50''}. It suggests some amazing
relationships between this stuff and \(4\)-dimensional topology\ldots.

\begin{center}\rule{0.5\linewidth}{0.5pt}\end{center}

\hypertarget{week84_tale}{
Let me conclude by reminding you where I am in the ``Tale of
\(n\)-Categories'', and where I want to go next.} So far I have spoken
mainly of 0-categories, \(1\)-categories, and \(2\)-categories, with
lots of vague allusions as to how various patterns generalize to higher
\(n\). Also, I have concentrated mainly on the related notions of
equality, isomorphism, equivalence, and adjointness. Equality,
isomorphism and equivalence are the most natural notions of ``sameness''
when working in 0-categories, \(1\)-categories, and \(2\)-categories,
respectively. Adjointness is a closely related but more subtle and
exciting concept that you can only start talking about once you get to
the level of \(2\)-categories. People got tremendously excited by it
when they started working with the \(2\)-category \(\mathsf{Cat}\) of
all small categories, because it explained a vast number of situations
where you have a way to go back and forth between two categories,
without the categories being ``the same'' (or equivalent). Another
exciting thing about adjointness is that it really highlights the
relation between \(2\)-categories and \(2\)-dimensional topology ---
thus pointing the way to a more general relation between
\(n\)-categories and \(n\)-dimensional topology. From this point of
view, adjointness is all about ``folds'': \[
  \begin{tikzpicture}
    \begin{knot}
      \strand[thick] (0,0)
      to [out=down,in=down,looseness=2] (1,0);
    \end{knot}
  \end{tikzpicture}
  \qquad
  \begin{tikzpicture}
    \begin{knot}
      \strand[thick] (0,0)
      to [out=up,in=up,looseness=2] (1,0);
    \end{knot}
  \end{tikzpicture}
\] and their ability to cancel: \[
  \begin{tikzpicture}
    \begin{knot}
      \strand[thick] (0,0)
      to (0,1)
      to [out=up,in=up,looseness=2] (1,1)
      to [out=down,in=down,looseness=2] (2,1)
      to (2,2);
    \end{knot}
    \node at (3,1) {$=$};
    \begin{scope}[shift={(4,0)}]
      \begin{knot}
        \strand[thick] (0,0) to (0,2);
      \end{knot}
    \end{scope}
  \end{tikzpicture}
\] \[
  \begin{tikzpicture}
    \begin{scope}[xscale=-1,shift={(-2,0)}]
      \begin{knot}
        \strand[thick] (0,0)
        to (0,1)
        to [out=up,in=up,looseness=2] (1,1)
        to [out=down,in=down,looseness=2] (2,1)
        to (2,2);
      \end{knot}
    \end{scope}
    \node at (3,1) {$=$};
    \begin{scope}[shift={(4,0)}]
      \begin{knot}
        \strand[thick] (0,0) to (0,2);
      \end{knot}
    \end{scope}
  \end{tikzpicture}
\] This concept of ``doubling back'' or ``backtracking'' is a very
simple and powerful one, so it's not surprising that it is prevalent
throughout mathematics and physics. It is an essentially
\(2\)-dimensional phenomenon (though it occurs in higher dimensions as
well), so it can be understood most generally in the language of
\(2\)-categories.

(In physics, ``doubling back'' is related to the notion of antiparticles
as particle moving backwards in time, and appears in the Feynman
diagrams for annihilation and creation of particle/antiparticle pairs.
For those familiar with the category-theoretic approach to Feynman
diagrams, the stuff in \protect\hyperlink{week83_tale}{``Week 83''} about
dual vector spaces should suffice to make this connection precise.)

Next I will talk about another \(2\)-dimensional concept, the concept of
``joining'' or ``merging'': \[
  \begin{tikzpicture}
    \begin{knot}
      \strand[thick] (0,0)
        to [out=down,in=up] (0.5,-1)
        to (0.5,-1.5);
      \strand[thick] (1,0)
        to [out=down,in=up] (0.5,-1);
    \end{knot}
  \end{tikzpicture}
\] This is probably even more powerful than the concept of ``folding'':
it shows up whenever we add numbers, multiply numbers, or in many other
ways combine things. The \(2\)-categorical way to understand this is as
follows. Suppose we have an object \(x\) in a \(2\)-category, and a
morphism \(f\colon x \to x\). Then we can ask for a \(2\)-morphism
\[M\colon f^2 \Rightarrow f.\] If we have such a thing, we can draw it
as a traditional \(2\)-categorical diagram: \[
  \begin{tikzpicture}
    \node (xl) at (0,0) {$x$};
    \node (xt) at (1.25,2) {$x$};
    \node (xr) at (2.5,0) {$x$};
    \draw[thick] (xl) to node[fill=white]{$f$} (xt);
    \draw[thick] (xt) to node[fill=white]{$f$} (xr);
    \draw[thick] (xl) to node[fill=white]{$f$} (xr);
    \draw[-implies,double equal sign distance] (xt) to (1.25,0.2);
    \node at (1,0.7) {$M$};
  \end{tikzpicture}
\] or dually as a ``string diagram'' \[
  \begin{tikzpicture}
    \begin{knot}
      \strand[thick] (0,0.5)
        to (0,0)
        to [out=down,in=up] (0.5,-1)
        to (0.5,-2);
      \strand[thick] (1,0.5)
        to (1,0)
        to [out=down,in=up] (0.5,-1);
    \end{knot}
    \node[fill=white] at (0,0) {$f$};
    \node[fill=white] at (1,0) {$f$};
    \node[label=left:{$M$}] at (0.5,-0.95) {$\bullet$};
    \node[fill=white] at (0.5,-1.5) {$f$};
  \end{tikzpicture}
\] Regardless of how you draw it, the \(2\)-morphism
\(M\colon f^2 \Rightarrow f\) represents a process turning two copies of
\(f\) into one. And as we'll see, all sorts of fancy ways mathematicians
have of studying this sort of process --- ``monoids'', ``monoidal
categories'', and ``monads'' --- are special cases of this sort of
situation.

To continue reading the ``Tale of \(n\)-Categories'', see
\protect\hyperlink{week89_tale}{``Week 89''}.

\hypertarget{week85}{%
\section{July 14, 1996}\label{week85}}

I'm spending this month at the Erwin Schr\"odinger Institute in Vienna,
where Abhay Ashtekar and Peter Aichelburg are running a workshop called
Mathematical Problems of Quantum Gravity.

Ashtekar is one of the founders of an approach to quantizing gravity
called the loop representation. I've explained this approach in
\protect\hyperlink{week7}{``Week 7''}, \protect\hyperlink{week43}{``Week
43''}, and other places, but let me just remind you of the basic idea.
In the traditional approach to reconciling general relativity with
quantum theory, excitations of the gravitational field were described by
small ripples in the geometry of flat spacetime, or ``gravitons''. In
the loop representation, they are instead described by collections of
loops, which we can think of as ``flux tubes of area'' floating in an
otherwise utterly featureless void. More recently, the loop approach has
been supplemented by a technical device known as ``spin networks'':
roughly speaking, a spin network is a graph whose edges are labelled by
spins \(0,1/2,1,3/2,\ldots\) with an edge of spin \(j\) corresponding to
a flux tube carrying area equal to \(\sqrt{j(j+1)}\) times the square of
the Planck length --- the fundamental length scale in quantum gravity,
about \(10^{-35}\) meters. For more on spin networks, try
\protect\hyperlink{week55}{``Week 55''}.

Quantum gravity has always been a tough subject. After a lot of work, a
lot of people concluded that the traditional approach to quantum gravity
didn't make sense, mathematically. This led to string theory, an attempt
to quantize gravity together with all the other forces and particles.
But in the late 1980s, Rovelli and Smolin revived hopes of quantizing
gravity alone by introducing the loop representation.

One doesn't expect the loop representation to describe much real physics
until one introduces other forces and particles. Pure gravity is just a
warm-up exercise --- but it's not at all easy! When the loop
representation was born, it was rather sketchy at many points. A lot of
mathematical problems had to be overcome to make it as precise as it is
now\ldots. and there are a lot of formidable difficulties left, any one
of which could spell doom for the theory. Luckily, progress has been
rapid. Many of the problems which seemed hopelessly beyond our reach a
few years ago can now be formulated precisely, and maybe even solved.
The idea of this workshop is to start tackling these problems.

A lot has been going on! People give talks at 11 in the morning, while
afternoons are devoted to more informal discussions in small groups.
There are general introductory talks on Tuesdays, more technical talks
on Thursdays, and short talks on research in progress on some other
days.

To give a bit of the flavor of the workshop, let me describe things day
by day. I'll need to describe some days very sketchily, though, or I'll
never finish writing this!

\begin{itemize}
\item
  \textbf{Wednesday, July 3} --- Rodolfo Gambini spoke on
  gauge-invariance in the extended loop representation. The idea of the
  loop representation is to study the gravitational vector potential by
  studying certain integrals of it around loops. Mathematicians call
  this the trace of the holonomy, and physicists call it a Wilson loop
  or the trace of a path-ordered exponential. In the loop
  representation, states of quantum gravity are described by certain
  functions that eat Wilson loops and spit out complex numbers...
  i.e., that assign an ``amplitude'' to each Wilson loop.

  In quantum field theory you often need to average a quantum field over
  some \(3\)-dimensional region of space or \(4\)-dimensional region of
  spacetime to get a well-defined operator. Wilson loops are rather
  singular because a loop is a one-dimensional object. On the other
  hand, they are nice because they are gauge-invariant: they don't
  change when we do a gauge transformation to the vector potential.

  In the ``extended'' loop representation one tries to make the integral
  less singular by not dealing with actual loops, but certain analogous
  integrals over all \(3\)-dimensional space. Heuristic calculations
  suggest that they are gauge-invariant, but Troy Schilling noticed a
  while ago that they aren't always \emph{really} gauge-invariant ---
  basically because the the path-ordered exponential is given by a
  certain Taylor series, and nasty things can happen when you manipulate
  infinite series without checking if your manipulations are legitimate!
  See:

  \begin{enumerate}
  \def\labelenumi{\arabic{enumi})}
  \tightlist
  \item
    Troy Schilling, ``Non-covariance of the generalized holonomies:
    Examples'', available as
    \href{https://arxiv.org/abs/gr-qc/9503064}{\texttt{gr-qc/9503064}}.
  \end{enumerate}

  There has been a certain amount of competition between the extended
  loop representation, developed by Gambini and various coauthors, and
  Ashtekar's approach. Thus Schilling's result was seen as a blow
  against the extended loop representation. In Gambini's talk, he argued
  that gauge-invariance is rigorously maintained by certain extended
  loops, e.g.~those for which the power series has finitely many terms.
  The most famous examples of functions of extended loops with only
  finitely many terms are the Vassiliev invariants, which come up in
  knot theory (see \protect\hyperlink{week3}{``Week 3''}). Gambini and
  Pullin have claimed that certain Vassiliev invariants are states of
  quantum gravity, so these are of special interest.

  The feeling was that we needed to compare these different loop
  representations more carefully because they both have advantages.

  Also, Renate Loll spoke about ``Lattice Gravity''. See
  \protect\hyperlink{week55}{``Week 55''} for a bit more on this. Her
  talk led to an interesting discussion of the meaning of the limit, as
  the lattice spacing goes to zero, of quantum gravity as done on a
  lattice. Does it make sense? One needs, apparently, to look at ones
  formula for the Hamiltonian constraint on the lattice, and see if it
  depends on the Planck length in a manner \emph{other than} having the
  Planck length as an overall prefactor. Various people tried to do the
  calculation on the spot, and got mixed up.
\item
  \textbf{Thursday, July 4} --- Thomas Thiemann spoke on ``The
  Hamiltonian Constraint for Lorentzian Canonical Quantum Gravity''.
  This was a big bombshell. The Hamiltonian constraint in quantum
  gravity is one of the biggest, baddest problems we are facing. It's
  the analog of Schrodinger's equation in quantum mechanics, but it's a
  constraint: \[H\psi=0.\] All the dynamics of the theory is contained
  in this equation, yet we only roughly understand how to define it in a
  rigorous way. Thiemann, a student of Ashtekar who is now a postdoc at
  Harvard, had put the following 5 papers on the general relativity
  preprint server right before the workshop. The first one gives a
  rigorous definition of the Hamiltonian constraint!

  \begin{enumerate}
  \def\labelenumi{\arabic{enumi})}
  \setcounter{enumi}{1}
  \item
    Thomas Thiemann, ``Quantum Spin Dynamics (QSD)'', available
    as
    \href{https://arxiv.org/abs/gr-qc/9606089}{\texttt{gr-qc/9606089}}.

    Thomas Thiemann, ``Quantum Spin Dynamics (QSD) II'', 
    available as
    \href{https://arxiv.org/abs/gr-qc/9606090}{\texttt{gr-qc/9606090}}.

    Thomas Thiemann, ``Anomaly-free formulation of non-perturbative,
    four-\break dimensional Lorentzian quantum gravity'', \emph{Phys. Lett. B}
    \textbf{380} (1996) 257--264.   Also available as
    \href{https://arxiv.org/abs/gr-qc/9606088}{\texttt{gr-qc/9606088}}.

    Thomas Thiemann, ``Closed formula for the matrix elements of the
    volume operator in canonical quantum gravity'', available
    as
    \href{https://arxiv.org/abs/gr-qc/9606091}{\texttt{gr-qc/9606091}}.

    Thomas Thiemann, ``A length operator for canonical quantum
    gravity'', available as
    \href{https://arxiv.org/abs/gr-qc/9606092}{\texttt{gr-qc/9606092}}.
  \end{enumerate}

  It is interesting to compare ``Quantum Spin Dynamics'' to the paper by
  Ashtekar and Lewandowksi, so far available only in draft form to a
  select few, in which they gave a rigorous definition of the square
  root of the Hamiltonian constraint. The advantage of ``QSD'' is that
  it deals directly with the Hamiltonian constraint, rather than its
  square root, and that it does this using some ingenious formulas for
  the Hamiltonian constraint of Lorentzian gravity in terms of the
  Hamiltonian constraint for Riemannian gravity and the total volume and
  total extrinsic curvature of the universe (which we assume is
  compact).

  You see, quantum gravity comes in two flavors, Lorentzian and
  Riemannian, depending on whether we work with real time --- the
  obviously sensible thing to do --- or imaginary time --- not at all
  obviously sensible, but with a curious mathematical charm to it, which
  has won many hearts. The interplay between these two has long been a
  bugaboo of the loop representation. The Lorentzian theory is harder to
  work with, so lots of people cheat and study the Riemannian theory.
  Sometimes they do this covertly, with a guilty conscience, so in some
  papers it's left unclear which theory the author is actually working
  with! Thiemann's work, however, seems to exploit the interplay between
  the theories in a benign way --- related to earlier ideas of Ashtekar,
  but different. I would like to understand this interplay more deeply.

  Due to jetlag I woke up at 4 am on the morning of this talk, and I
  couldn't get back to sleep, so I read his paper. When I came to the
  Institute at 9 am --- a shockingly early hour for people working on
  quantum gravity --- I was sure nobody would be there yet. But as I
  entered I bumped into Carlo Rovelli. It turned out he had stayed up
  all night reading Thiemann's paper, too excited to sleep!

  After this talk everyone was busily trying to learn Thiemann's stuff,
  trying to figure out if it is physically correct, and trying to figure
  out what to do next.
\item
  \textbf{Tuesday, July 9} --- Abhay Ashtekar gave a general talk on the
  ``Quantum Theory of Geometry''. Most of it was well-known stuff to
  fans of the loop representation, but one new tidbit concerned the
  noncommutativity of area operators. Since the area of surfaces in
  space depends only on the metric on space, not on its first time
  derivative, one might expect their quantum analogs to commute, since
  the metric and its first time derivative are analogous to position and
  momentum in quantum mechanics. But they don't commute! In a later
  talk, Ashtekar explained that this is not really a strange new feature
  of quantum gravity, but one which has its classical analog.
\item
  \textbf{Wednesday, July 10} --- Kirill Krasnov gave a talk on a paper
  we started working on together just recently, ``Quantization of
  diffeomorphism invariant theories with fermions''. Kirill is a young
  Ukrainian physicist whom I first met last summer in Warsaw; he had
  written a nice paper on the loop representation of quantum gravity
  coupled to electromagnetism and fermions:

  \begin{enumerate}
  \def\labelenumi{\arabic{enumi})}
  \setcounter{enumi}{2}
  \tightlist
  \item
    Kirill Krasnov, ``Quantum loop representation for fermions coupled
    to Eistein--Maxwell field'', \emph{Phys.\ Rev.\ D} \textbf{53} (1996),
    1874; available as
    \href{https://arxiv.org/abs/gr-qc/9506029}{\texttt{gr-qc/9506029}}.
  \end{enumerate}

  When I met him again here, it turned out he was continuing this work,
  and also making it more rigorous. Now, I had for some time been
  meaning to write something with Hugo Morales-Tecotl showing that a
  slight generalization of spin network states form a basis of states
  for such theories. These states had already appeared, for example, in
  his work with Rovelli:

  \begin{enumerate}
  \def\labelenumi{\arabic{enumi})}
  \setcounter{enumi}{3}
  \item
    Carlo Rovelli and Hugo Morales-Tecotl, ``Fermions in quantum
    gravity'', \emph{Phys. Rev.~Lett.} \textbf{72} (1994), 3642--3645.

    Carlo Rovelli and Hugo Morales-Tecotl, \emph{Nucl. Phys. B}
    \textbf{451} (1995), 325.  Also available as
    \href{https://arxiv.org/abs/gr-qc/9401011}{\texttt{gr-qc/9401011}}.
  \end{enumerate}

  But we had never gotten around to it. So, I decided to team up with
  Kirill and write a paper on this stuff.
\end{itemize}

\hypertarget{week86}{%
\section{August 6, 1996}\label{week86}}

Let me continue my reportage of what happened at the Mathematical
Problems of Quantum Gravity workshop in Vienna. I will only write about
quantum gravity aspects here. I had an interesting conversation with
Bertram Kostant in which he explained to me the deep inner secrets of
the exceptional Lie group \(\mathrm{E}_8\). However, my writeup of that
has grown to the point where I will save it for some other week.

By the way, my course on \(n\)-category theory is not over! I'm merely
taking a break from it, and will return to it after this workshop.

So....

\begin{itemize}
\item
  \textbf{Wednesday, July 10th} --- Jerzy Lewandowski gave a talk on the
  ``Spectrum of the Area Operator''. As I've mentioned a few times
  before, one of the exciting things about the loop representation of
  quantum gravity is that the spectrum of the area operator associated
  to any surface is discrete. In other words, area is quantized!

  Let me describe the area operator a bit more precisely. Before I tell
  you what the area operator is, I have to tell you what it operates on.
  Remember from \protect\hyperlink{week43}{``Week 43''} that there are
  various Hilbert spaces floating around in the canonical quantization
  of gravity. First there is the ``kinematical state space''. In the
  old-fashioned metric approach to quantum gravity, known as
  geometrodynamics, this was supposed to be \(L^2(\mathrm{Met})\), where
  \(\mathrm{Met}\) is the space of Riemannian metrics on space. (We take
  as space some 3-manifold \(S\), which for simplicity we assume is
  compact.) The problem was that nobody knew how to rigorously define
  this Hilbert space \(L^2(\mathrm{Met})\). In the ``new variables''
  approach to quantum gravity, the kinematical state space is taken
  instead to be \(L^2(\mathcal{A})\), where \(\mathcal{A}\) is the space
  of connections on some trivial \(\mathrm{SU}(2)\) bundle over
  \(S\). This \emph{can} be defined rigorously.

  Here's the idea, roughly. Fix any graph \(g\), with finitely many
  edges and vertices, embedded in \(S\). Let \(\mathcal{A}_g\), the
  space of connections on that graph, be \(\mathrm{SU}(2)^n\) where
  \(n\) is the number of edges in \(e\). Thus a connection on a graph
  tells us how to parallel transport things along each edge of that
  graph --- an idea familiar from lattice gauge theory.
  \(L^2(\mathcal{A}_g)\) is well-defined because \(\mathrm{SU}(2)\) has
  a nice measure on it, namely Haar measure, so \(\mathcal{A}_g\) gets a
  nice measure on it as well.

  Now if one graph \(g\) is contained in another graph \(h\), the space
  \(L^2(\mathcal{A}_g)\) is contained in the space
  \(L^2(\mathcal{A}_h)\) in an obvious way. So we can form the union of
  all the Hilbert spaces \(L^2(\mathcal{A}_g)\) and get a big Hilbert
  space \(L^2(\mathcal{A})\). Mathematicians would say that
  \(L^2(\mathcal{A})\) is the ``projective limit'' of the Hilbert spaces
  \(L^2(\mathcal{A}_g)\), but let's not worry about that.

  So that's how we get the space of ``kinematical states'' in the loop
  representation of quantum gravity. The space of physical states is
  then obtained by imposing constraints: the Gauss law constraint (i.e.,
  gauge-invariance), the diffeomorphism constraint (i.e., invariance
  under diffeomorphisms of S) and the Hamiltonian constraint (i.e.,
  invariance under time evolution). States in the physical state space
  are supposed to only contain information that's invariant under all
  coordinate transformations and gauge transformations --- the really
  physical information.

  I explained these constraints to some extent in
  \protect\hyperlink{week43}{``Week 43''}, and I don't really want to
  worry about them here. But let me just mention that the Gauss law
  constraint is easy to impose in a mathematically rigorous way. The
  diffeomorphism constraint is harder but still possible, and the
  Hamiltonian constraint is the big thorny question plaguing quantum
  gravity --- see \protect\hyperlink{week85}{``Week 85''} for some
  recent progress on this. The area operators I'll be talking about are
  self-adjoint operators on the space of kinematical states,
  \(L^2(\mathcal{A})\), and are a preliminary version of some related
  operators one hopes eventually to get on the physical state space,
  after much struggle and sweat.

  To define an operator on \(L^2(\mathcal{A})\) it's enough to define it
  on \(L^2(\mathcal{A}_g)\) for every graph \(g\) and then check that
  these definitions fit together consistently to give an operator on the
  big space \(L^2(\mathcal{A})\). So let's take a graph \(g\) and a
  surface \(s\) in space. The area operator we're after is supposed to
  be the quantum analog of the usual classical formula for the area of
  \(s\). The usual classical area is a function of the metric on space;
  similarly, the quantum area is an operator on the space
  \(L^2(\mathcal{A})\).

  The area operator only cares about the points where the graph
  intersects the surface. We assume that there are only finitely many
  points where it does so, apart from points where the edges are tangent
  to the surface. (To make this assumption reasonable, we need to
  assume, e.g., that the space \(S\) has a real-analytic structure and
  the surface and graph are analytic --- an annoying technicality that I
  have been seeking to eliminate.)

  The area operator is built using three operators on
  \(L^2(\mathrm{SU}(2))\) called \(J_1\), \(J_2\), and \(J_3\), the
  self-adjoint operators corresponding to the 3 generators of
  \(\mathrm{SU}(2)\) --- which often show up in physics as the three
  components of angular momentum! Alternatively, we can think of all
  three together as one vector-valued operator \(J\), the ``angular
  momentum operator''. Note that \(L^2(\mathcal{A}_g)\) is just the
  tensor product of one copy of the Hilbert space
  \(L^2(\mathrm{SU}(2))\) for each edge of our graph \(g\). Thus for any
  edge \(e\) we get an angular momentum operator \(J(e)\) that acts on
  the copy of \(L^2(\mathrm{SU}(2))\) corresponding to the edge \(e\) in
  question, leaving the other copies alone.

  This, then, is how we define the operator on \(L^2(\mathcal{A}_g)\)
  corresponding to the area of the surface \(s\). Pick an orientation
  for the surface \(s\). For any point where the graph \(g\) intersects
  \(s\), let \(J(\mathrm{in})\) denote the sum of the angular momentum
  operators of all edges intersecting \(s\) at the point in question and
  pointing ``inwards'' relative to our chosen orientation. Similarly,
  let \(J(\mathrm{out})\) be the sum of the angular momentum operators
  of edges intersecting \(s\) at the point in question and pointing
  ``outwards''. (Note: edges tangent to the surface contribute neither
  to \(J(\mathrm{in})\) nor \(J(\mathrm{out})\).) Now sum up, over all
  points where the graph intersects the surface, the following quantity:
  \[\sqrt{(J(\mathrm{in})-J(\mathrm{out})) \cdot (J(\mathrm{in})-J(\mathrm{out}))}\]
  where the dot denotes the obvious sort of dot product of vector-valued
  operators. Multiply by half the Planck length squared and you've got
  the area operator!

  This very beautiful and simple formula was derived by Ashtekar and
  Lewandowski, but the first people to try to quantize the area operator
  were Rovelli and Smolin; see

  \begin{enumerate}
  \def\labelenumi{\arabic{enumi})}
  \item
   Carlo Rovelli and Lee Smolin, 
  ``Discreteness of area and volume in quantum gravity'', \emph{Nucl. Phys. B}
    \textbf{442} (1995), 593--619.   Also available as
    \href{https://arxiv.org/abs/gr-qc/9411005}{\texttt{gr-qc/9411005}}.

    Abhay Ashtekar and Jerzy Lewandowski, ``Quantum theory of geometry
    I: area operators'', \emph{Class. Quant. Grav.} \textbf{14} (1997)
    A55. Also available as
    \href{https://arxiv.org/abs/gr-qc/9602046}{\texttt{gr-qc/9602046}}.
  \end{enumerate}

  In his talk Jerzy showed how to work the spectrum of the area operator
  (which is discrete) and showed how it could depend on whether the
  surface \(s\) cuts space into two parts or not.

  Later that day, Mike Reisenberger, Matthias Blau, Carlo Rovelli and I
  talked about the relation between string theory and the loop
  representation of quantum gravity.

  Mike has been working on a very interesting ``state sum model'' for
  quantum gravity; that is, a discretized model in which spacetime is
  made of \(4\)-simplices (the 4d version of tetrahedra), fields are
  thought of as ways of labelling the faces, edges and so on by spins,
  elements of \(\mathrm{SU}(2)\) and the like, and the path integral is
  replaced by a sum over these labellings. This works out quite nicely
  for quantum gravity in 3 dimensions --- see
  \protect\hyperlink{week16}{``Week 16''} --- but it's much more
  challenging in 4 dimensions.

  One nice feature of these state sum models for quantum gravity is that
  they may be reinterpreted as sums over ``worldsheets'' --- surfaces
  mapped into spacetime. Since the spacetime is discrete, so are these
  surfaces --- they're made of lots of triangles --- but apart from
  that, having a path integral that's a sum over worldsheets is
  pleasantly reminscent of string theory. Indeed, once upon a time I
  proposed that the loop representation of quantum gravity and string
  theory were two aspects of some theory still waiting to be fully
  understood:

  \begin{enumerate}
  \def\labelenumi{\arabic{enumi})}
  \setcounter{enumi}{1}
  \tightlist
  \item
    John Baez, ``Strings, loops, knots, and gauge fields'', in
    \emph{Knots and Quantum Gravity}, ed.~J. Baez, Oxford U. Press,
    Oxford, 1994.   Available as
    \href{https://arxiv.org/abs/hep-th/9309067}{\texttt{hep-th/9309067}},
    34 pages.
  \end{enumerate}

  The problem was getting a concrete way to relate the Lagrangian for
  the string theory to the Lagrangian for gravity (or whatever gauge
  theory one started with). Iwasaki tackled this problem in
  3d quantum gravity using state sum models:

  \begin{enumerate}
  \def\labelenumi{\arabic{enumi})}
  \setcounter{enumi}{2}
  \tightlist
  \item
    Junichi Iwasaki, ``A reformulation of the Ponzano-Regge quantum
    gravity model in terms of surfaces'', University of Pittsburgh, 11
    pages in LaTeX format available as
    \href{https://arxiv.org/abs/gr-qc/9410010}{\texttt{gr-qc/9410010}}.
  \end{enumerate}

  Later, Reisenberger extended this approach to deal with certain 4d
  theories which are simpler than quantum gravity, like \(BF\) theory:

  \begin{enumerate}
  \def\labelenumi{\arabic{enumi})}
  \setcounter{enumi}{3}
  \tightlist
  \item
    Michael Reisenberger, ``Worldsheet formulations of gauge theories
    and gravity'', available as
    \href{https://arxiv.org/abs/gr-qc/9412035}{\texttt{gr-qc/9412035}}.
  \end{enumerate}

  In all of these theories, one computes the action for the worldsheets
  by summing something over places where they intersect. In other words,
  they ``interact'' at intersections.

  But the really exciting thing would be to do something like this for
  Mike's new state sum model for 4d quantum gravity. And the real
  challenge would be to relate this --- if possible! --- to conventional
  string theory. In a coffeeshop I suggested that one might do this by
  using the usual formula for the action in (bosonic) string theory.
  This is simply the \emph{area} of the string worldsheet with respect
  to some background metric. The loop representation of quantum gravity
  doesn't make reference to any background metric; the closest
  approximation to a classical metric is a ``weave'' state in which
  space is tightly packed with lots of loops or spin networks. From the
  4d point of view, we'd expect this to correspond to a spacetime packed
  with lots of worldsheets. Now, given the relation between area and
  intersection number in the loop representation (see above!), one might
  expect the area of a given worldsheet to be roughly proportional to
  the number of its intersections with the other worldsheets in this
  ``weave''. But this is what one would expect in any theory where the
  worldsheets interact at intersections. So, one could hope that Mike's
  state sum model would be approximately equivalent to a string theory
  of the sort string theorists study.

  There are lots of obvious problems with this idea, but it led to an
  interesting conversation, and I am still not convinced that it is
  crazy.
\item
  \textbf{Thursday, July 11th} --- Jorge Pullin spoke on skein relations
  and the Hamiltonian constraint in lattice quantum gravity. His idea
  was that the Hamiltonian constraint contains a ``topological factor''
  which serves as a skein relation on loop states.
\item
  \textbf{Friday, July 12th} --- Abhay Ashtekar gave a talk on
  ``Noncommutativity of Area Operators''. This explained how the rather
  shocking fact that the area operators for two intersecting surfaces
  needn't commute actually has a perfect analog in classical general
  relativity.

  Mike Reisenberger spoke on ``Euclidean Simplicial GR''. This presented
  the details of his state sum model. Since he hasn't published this
  yet, and since I am getting a bit tired out, I won't describe it here.
\item
  \textbf{Monday, July 15th} --- Renate Loll gave a talk on the volume
  and area operators in lattice gravity. I wrote a bit about her work on
  the volume operator in \protect\hyperlink{week55}{``Week 55''}, and
  more can be found in:

  \begin{enumerate}
  \def\labelenumi{\arabic{enumi})}
  \setcounter{enumi}{4}
  \item
    Renate Loll, ``The volume operator in discretized quantum gravity'',
    available as
    \href{https://arxiv.org/abs/gr-qc/9506014}{\texttt{gr-qc/9506014}}.

    Renate Loll, ``Spectrum of the volume operator in quantum gravity'',
    available as
    \href{https://arxiv.org/abs/gr-qc/9511030}{\texttt{gr-qc/9511030}}.
  \end{enumerate}

  Also, Jerzy Lewandowski spoke on his work with Ashtekar on the volume
  operator in the continuum theory:

  \begin{enumerate}
  \def\labelenumi{\arabic{enumi})}
  \setcounter{enumi}{5}
  \item
    Jerzy Lewandowski, ``Volume and quantizations'', available
    as
    \href{https://arxiv.org/abs/gr-qc/9602035}{\texttt{gr-qc/9602035}}.

    Abhay Ashtekar and Jerzy Lewandowski, ``Quantum theory of geometry
    II: volume operators'', manuscript in preparation.
  \end{enumerate}

  The volume operator is more tricky than the area operator, and various
  proposed formulas for it do not agree. This is summarized quite
  clearly in Jerzy's paper.

  In fact, I have already left Vienna by now. I was too busy there to
  keep up with This Week's Finds, but my life is a bit calmer now and I
  will try to finish these reports soon.
\end{itemize}

\hypertarget{week87}{%
\section{August 20, 1996}\label{week87}}

Let me continue summarizing what happened during July at the
Mathematical Problems of Quantum Gravity workshop in Vienna. The first
two weeks concentrated on the foundations of the loop representation of
quantum gravity; the next week was all about black holes!

\begin{itemize}
\item
  \textbf{Tuesday, July 16th} --- Ted Jacobson gave an overview of
  ``Issues of Black Hole Thermodynamics''. There is a lot to say about
  this subject and I won't try to repeat his marvelous talk here. Let me
  just mention a very interesting technical point he made. The original
  Bekenstein-Hawking formula for the entropy of a black hole is
  \[S=A/(4\hbar G)\] where \(A\) is the area of the event horizon,
  \(\hbar\) is Planck's constant, and G is Newton's constant. One way to
  try to derive this is from the partition function of a quantum field
  theory involving gravity and other fields. Jacobson sketched a
  heuristic calculation along these lines. When you do this calculation
  it's natural to worry why the other fields, representing various forms
  of matter, don't seem to contribute to the answer above. Also, when we
  do quantum field theory, there is often a difference between the
  ``bare'' coupling constants we put into the theory and the
  ``renormalized'' coupling constants that are what the theory predicts
  we'll observe experimentally. So it's natural to worry about whether
  it's the bare or renormalized Newton's constant \(G\) that enters the
  above formula --- even though quantum gravity is so unlike most other
  quantum field theories that it's unclear that this worry makes sense,
  ultimately.

  Anyway, the nice thing is that these two worries cancel each other
  out. In other words: yes, it's the renormalized Newton's constant
  \(G\) --- the physically measured one --- that enters the above
  formula. But at least to first order in \(\hbar\), the difference
  between the bare \(G\) and the renormalized \(G\) is precisely due to
  the interactions between gravity and the matter fields (including the
  self-interaction of the gravitational field). In other words, the
  matter fields really \emph{do} contribute to the black hole entropy,
  but this contribution is absorbed into the definition of the
  renormalized \(G\).

  In the most extreme case, the bare \(1/G\) is zero, and the
  renormalized \(1/G\) is entirely due to interactions between matter
  and gravity. This is Andrei Sakharov's theory of ``induced gravity''.
  According to Jacobson, in this case all of the black hole entropy is
  ``entanglement entropy'' --- this being standard jargon for the way
  that two parts of a quantum system can each have entropy due to
  correlations, even though the whole system has zero entropy.
  Unfortunately my notes do not allow me to reconstruct the wonderful
  argument whereby he showed this. (See
  \protect\hyperlink{week27}{``Week 27''} for a more detailed
  explanation of entanglement entropy.)
\item
  \textbf{Wednesday July 17th} --- There was a talk on ``Colombeau
  theory'' by a mathematician whose name I unfortunately failed to
  catch. Colombeau theory is a theory that allows you to multiply
  distributions, just like they said in school that you weren't allowed
  to do. So if for example you want to square the Dirac delta function,
  you can do it in the context of Colombeau theory. There has been a
  certain amount of debate, however, on whether Colombeau theory allows
  you to this multiplication in a \emph{useful} way. There were a lot of
  physicists at this talk who would be willing and eager to master
  Colombeau theory if it let one solve the physics problems they wanted
  to solve. However, after much discussion, it appears that they didn't
  buy it. I believe that at best Colombeau theory provides a useful
  framework for understanding the ambiguities one encounters when
  multiplying distributions.

  I say ``ambiguities'' rather than ``disasters'' because while the
  square of the Dirac delta function has no sensible interpretation as a
  distribution, there are many cases, such as when you try to multiply
  the Dirac delta function and the Heaviside function, where you can
  interpret the result as a distribution in a variety of ways. These
  ambiguous cases are the ones of greatest interest in physics. A nice
  place to see this in quantum field theory is in

  \begin{enumerate}
  \def\labelenumi{\arabic{enumi})}
  \tightlist
  \item
    G. Scharf, \emph{Finite Quantum Electrodynamics: the Causal
    Approach}, Springer, Berlin, 1995.
  \end{enumerate}

  If you want to learn about Colombeau theory you can try:

  \begin{enumerate}
  \def\labelenumi{\arabic{enumi})}
  \setcounter{enumi}{1}
  \tightlist
  \item
    J. F. Colombeau, \emph{Multiplication of Distributions: a Tool in
    Mathematics, Numerical Engineering, and Theoretical Physics},
    Lecture Notes in Mathematics \textbf{1532}, Springer, Berlin, 1992.
  \end{enumerate}

  Later that day I had nice conversation with Jerzy Lewandowski on the
  approach to the loop representation where one uses smooth, rather than
  analytic, loops. (See \protect\hyperlink{week55}{``Week 55''} for more
  on this issue.)
\item
  \textbf{Thursday, July 18th} --- Carlo Rovelli spoke on ``Black Hole
  Entropy'', reporting some work he did with Kirill Krasnov. They have a
  nice approach to computing the black hole entropy using the loop
  representation of quantum gravity. A common goal among quantum gravity
  folks is to recover the Bekenstein-Hawking formula from some
  full-fledged theory of quantum gravity --- the original derivation
  being a curious ``semiclassical'' hybrid of quantum and classical
  reasoning. In a statistical mechanical approach, entropy should be the
  logarithm of the number of microstates some system can have in a given
  macrostate. So one wants to count states somehow. Basically what
  Rovelli and Krasnov do is count the number of ways a surface can be
  pierced by a spin network so as to give it a certain area. (This uses
  the formula for the area operator I described in
  \protect\hyperlink{week86}{``Week 86''}.) They get an entropy
  proportional to the area, but not with the same constant as in the
  Bekenstein-Hawking formula.

  There were some hopes that taking matter fields into account might
  give the right constant. But since everyone had been to Ted Jacobson's
  talk, this led to much interesting wrangling over whether Rovelli and
  Krasnov were using the bare or renormalized Newton's constant \(G\),
  and whether the concept of bare and renormalized \(G\) even makes
  sense, ultimately! Also, there are some extremely important puzzles
  about what the right way to count states is, in these loop
  representation computations.

  For more, try:

  \begin{enumerate}
  \def\labelenumi{\arabic{enumi})}
  \setcounter{enumi}{2}
  \item
    Carlo Rovelli, ``Loop quantum gravity and black hole physics'',
     available as
    \href{https://arxiv.org/abs/gr-qc/9608032}{\texttt{gr-qc/9608032}}.

    Kirill Krasnov, ``The Bekenstein bound and non-perturbative quantum
    gravity'', available as
    \href{https://arxiv.org/abs/gr-qc/9603025}{\texttt{gr-qc/9603025}}.

    Kirill Krasnov, ``On statistical mechanics of gravitational
    systems'', available as
    \href{https://arxiv.org/abs/gr-qc/9605047}{\texttt{gr-qc/9605047}}.
  \end{enumerate}
\item
  \textbf{Friday, July 19th} --- Don Marolf spoke on ``Black hole
  entropy in string theory''. He attempted valiantly to describe the
  string-theoretic approach to computing black hole entropy to an
  audience only generally familiar with string theory. I will not try to
  summarize his talk, except to note that he mainly discussed the case
  of a black hole in 5 dimensions, which was really a ``black string''
  in 6 dimensions --- a solution with translational symmetry in the 6th
  dimension, but where the extra 6th dimension is so tiny that ordinary
  \(5\)-dimensional folks think they've just got a black hole. (By the
  way, even the \(6\)-dimensional approach is really just a way of
  talking about a string theory that fundamentally lives in 10
  dimensions. This stuff is not for the faint-hearted.)

  Here are a few papers on this subject by Marolf and Horowitz:

  \begin{enumerate}
  \def\labelenumi{\arabic{enumi})}
  \setcounter{enumi}{3}
  \item
    Gary Horowitz, ``The origin of black hole entropy in string
    theory'', available as
    \href{https://arxiv.org/abs/gr-qc/9604051}{\texttt{gr-qc/9604051}}.

    Gary Horowitz and Donald Marolf, ``Counting states of black
    strings with traveling waves'', available as
    \href{https://arxiv.org/abs/hep-th/9605224}{\texttt{hep-th/9605224}}.

    Gary Horowitz and Donald Marolf, ``Counting states of black
    strings with traveling waves II'', available as
    \href{https://arxiv.org/abs/hep-th/9606113}{\texttt{hep-th/9606113}}.
  \end{enumerate}
\item
  \textbf{Monday, July 22nd} --- Kirill Krasnov spoke on ``The
  Eistein--Maxwell Theory of Black Hole Entropy''. This was a report on
  attempts to see how his calculations of the black entropy in the loop
  representation changed when he took the electromagnetic field into
  account. The calculations were very tentative, for certain technical
  reasons I won't go into here, but they made even clearer the
  importance of the issue of how one counts states when computing
  entropy in this approach.

  Later, I had a nice conversation with Carlo Rovelli about my hopes for
  thinking of fermions (e.g., electrons) as the ends of wormholes in the
  loop representation of quantum gravity. We came up with a nice
  heuristic argument to get the right Fermi statistics for these
  wormhole ends. Hopefully we can make this all more precise at some
  later date.
\item
  \textbf{Tuesday, July 23rd} --- Ted Jacobson gave informal talks on
  two subjects, the first of which was ``Transplanckian puzzle: origin
  of outgoing black hole modes''. This dealt with the puzzling fact that
  in the standard computation of Hawking radiation, the rather
  low-frequency radiation which leaves the hole is the incredibly
  redshifted offspring of high-frequency modes which swung past the
  horizon shortly after the hole's formation --- modes whose wavelength
  is far smaller than the Planck length!

  What if spacetime is ``grainy'' in some way at the Planck scale?
  Jacobson studied this using an analogy introduced by Unruh. If you
  have fluid flowing down a narrowing pipe, and at some point the
  velocity of the fluid flow exceeds the speed of sound in the fluid,
  there will be a ``sonic horizon''. In other words, there is a line
  where the fluid flow exceeds the speed of sound, and no sound can work
  its way upstream across that line. Now if you quantize the theory of
  sound in a simple-minded way you get ``phonons'' --- which have indeed
  been observed in solid-state physics. Unruh showed that in the case at
  hand you would get ``Hawking radiation'' of phonons from the sonic
  horizon, going upstream and getting shifted to lower frequencies as
  they go.

  Jacobson considered what would happen if you actually took into
  account the graininess of the fluid. (He considered the theory of
  liquid helium, to be specific.) The graininess at the molecular scale
  means that the group velocity of waves drops at very high frequencies.
  So what happens instead of ``Hawking radiation'' is something rather
  different. Start with a high-frequency wave attempting to go upstream,
  starting from upstream of the sonic horizon. Its group velocity is
  very slow so it fails miserably and gets swept toward the sonic
  horizon, like a hapless poor swimmer getting pulled to the edge of a
  waterfall despite trying to swim upstream. But as it gets pulled near
  the horizon its wavelength increases, and thus group velocity
  increases, thus allowing it to shoot upstream at the last minute! (An
  analogous process is apparently familiar in plasma physics under the
  name of ``mode conversion''.) In this scenario, the Hawking radiation
  winds up resulting from incoming modes through this process of mode
  conversion --- modes that have short wavelength, but not as short as
  the intermolecular spacing (or Planck length, in the gravitational
  case.)

  Ted Jacobson's second talk was even more interesting to me, but I'll
  postpone that for next Week.

  Here, by the way, is a paper related to the talk by Pullin discussed
  in \protect\hyperlink{week86}{``Week 86''}:

  \begin{enumerate}
  \def\labelenumi{\arabic{enumi})}
  \setcounter{enumi}{4}
  \tightlist
  \item
    Hugo Fort, Rodolfo Gambini and Jorge Pullin, ``Lattice knot theory
    and quantum gravity in the loop representation'', available
    as
    \href{https://arxiv.org/abs/gr-qc/9608033}{\texttt{gr-qc/9608033}}.
  \end{enumerate}
\end{itemize}

\hypertarget{week88}{%
\section{August 26, 1996}\label{week88}}

This issue concludes my report of what happened at the Mathematical
Problems of Quantum Gravity workshop in Vienna. I left the workshop at
the end of July, so my reportage ends there, but the workshop went on
for a few more weeks after that. I'll be really bummed out if I find out
that they solved all the problems with quantum gravity after I left.

Before I launch into my day-by-day account of what happened, let me note
that I've written a little introduction to Thiemann's work on the
Hamiltonian constraint, which he presented at the workshop (see
\protect\hyperlink{week85}{``Week 85''}):

\begin{enumerate}
\def\labelenumi{\arabic{enumi})}
\tightlist
\item
  John Baez, ``The Hamiltonian constraint in the loop representation of
  quantum gravity'', available at
  \href{http://math.ucr.edu/home/baez/hamiltonian/}{\texttt{http://math.ucr.edu/home/baez/hamiltonian/}}.
\end{enumerate}
\noindent
A less technical version of this appears in Jorge Pullin's newsletter
\emph{Matters of Gravity}, issue 8, at \href{http://www.phys.lsu.edu//mog/mog8/node7.html}{\texttt{http://www.phys.lsu.edu//mog/mog8/node7.html}}.

Okay... I'll start out simple today since there is something nice
and simple to ponder:

\begin{itemize}
\item
  \textbf{Tuesday, July 23rd} --- Ted Jacobson spoke on the ``Geometry
  and Evolution of Degenerate Metrics''. One of the interesting things
  about Ashtekar's reformulation of general relativity is that it
  extends general relativity to the case of degenerate metrics, that is,
  metrics where there are vectors whose dot product with all other
  vectors is zero. However, one needs to be very careful because
  different versions of Ashtekar's formulation give \emph{different}
  ways of handling degenerate metrics.

  To see why in a simple example, remember that the usual metric on
  Minkowski spacetime is nondegenerate and in nice coordinates looks
  like \[-dt^2 + dx^2 + dy^2 + dz^2\] Here we are setting the speed of
  light equal to \(1\). In general relativity, one way people describe
  the metric is using a tensor \(g_{ab}\), where the indices \(a\) and
  \(b\) go from 0 to 3. In Minkowski space this tensor equals \[
      \left(
        \begin{array}{cccc}
          -1&0&0&0
        \\0&1&0&0
        \\0&0&1&0
        \\0&0&0&1
        \end{array}
      \right)
    \] What this tensor means is that if you have two vectors \(v\) and
  \(w\), their dot product is \(g_{ab} v^a w^b\), where as usual we
  multiply the entries of the metric tensor and the vectors \(v\) and
  \(w\) as indicated, and then sum over repeated indices. So, for
  example, the dot product of the vector \[v = (1, 1, 0, 0)\] with
  itself is \(0\), though its dot product with other vectors needn't be
  zero. There is a bunch of vectors whose dot products with themselves
  are zero, and these are called lightlike vectors, because light
  travels in these directions, moving one unit in space for each unit in
  time. There is actually a cone of lightlike vectors, called the
  lightcone.

  One can imagine a world where the metric \(g_{ab}\) is \[
      \left(
        \begin{array}{cccc}
          -1&0&0&0
        \\0&1&0&0
        \\0&0&k&0
        \\0&0&0&k
        \end{array}
      \right)
    \] for some \(k > 0\). This world isn't really so different from
  Minkowski space, because you can also think of it as Minkowski space
  described in screwy coordinates where you are measuring distances in
  the \(y\) and \(z\) directions in different units than the \(x\)
  direction. When \(k\) gets small, you can check that the lightcone
  gets stretched out in the \(y\) and \(z\) directions. Alternatively,
  when \(k\) gets big, the lightcone gets squashed in the \(y\) and
  \(z\) directions.

  Another way to formulate general relativity uses the inverse metric
  \(g^{ab}\). This is just the inverse of the matrix \(g_{ab}\), which
  is indeed invertible when the metric is nondegenerate. So for example
  in the above case the inverse metric \(g^{ab}\) is \[
      \left(
        \begin{array}{cccc}
          -1&0&0&0
        \\0&1&0&0
        \\0&0&K&0
        \\0&0&0&K
        \end{array}
      \right)
    \] where \(K = 1/k\). You can think of \(K\) as the speed of light
  in the \(y\) and \(z\) directions, which is different from the speed
  of light in the \(x\) direction.

  Now there are two different limiting cases we can consider, depending
  on whether we work with the metric or the inverse metric. If we work
  with the metric, we can let \(k = 0\). This corresponds to making the
  speed of light in the \(y\) and \(z\) directions \emph{infinite}, so
  that information can go as fast as it likes in those directions and
  the lightcone gets completely stretched out in those directions. Note
  that now the metric \(g_{ab}\) is \[
      \left(
        \begin{array}{cccc}
          -1&0&0&0
        \\0&1&0&0
        \\0&0&0&0
        \\0&0&0&0
        \end{array}
      \right)
    \] so the inverse metric doesn't even make sense --- you can't
  invert this matrix. If we extend general relativity to degenerate
  metrics, we are allowing ourselves to study weird worlds like this.
  Why we'd want to --- well, that's another matter.

  If we work with the inverse metric, we can't let \(k = 0\), but we can
  let \(K = 0\). This corresponds to making the speed of light in the
  \(y\) and \(z\) directions \emph{zero}, so that information can't go
  at all in those directions: the lightcone is squashed down onto the
  \(t\)-\(x\) plane. Now it's the inverse metric that equals \[
      \left(
        \begin{array}{cccc}
          -1&0&0&0
        \\0&1&0&0
        \\0&0&0&0
        \\0&0&0&0
        \end{array}
      \right)
    \] and the metric doesn't even make sense.

  Ted Jacobson's talk was about doing general relativity in weird worlds
  like this, where the inverse metric is degenerate. Here information
  flows only along surfaces, like the \(x\)-\(t\) plane in the example
  above, and these different surfaces don't really talk to each other
  very much. It's as if the world was split up (or in math jargon,
  foliated) into lots of different \(2\)-dimensional worlds, which
  didn't know about each other. Jacobson showed that in this case, the
  equations of general relativity (extended in a certain way to
  degenerate inverse metrics) boil down to saying that there are two
  kinds of massless spin-\(1/2\) particle living on all these
  \(2\)-dimensional worlds.

  This got me quite excited because it reminded me of string theory,
  which is all about massless particles (or in physics jargon, conformal
  fields) living on the \(2\)-dimensional string worldsheet. I am always
  hunting around for relationships between string theory and the loop
  representation of quantum gravity, and I think this is an important
  clue. The reason is that I think the loop representation can be
  thought of as a quantum version of the theory of degenerate solutions
  of general relativity where the metric is \emph{zero} most places and
  less degenerate (but still degenerate) on certain surfaces. When you
  slice one of these surfaces with the hyperplane \(t = 0\) you get a
  bunch of loops (or more generally a graph), and these are the loops of
  the loop representation. Jacobson's talk may give a way to understand
  the conformal field theory living on these surfaces, which one needs
  if one wants to think of these surfaces as the ``string worldsheets''
  of string theory fame. Anyway, I am busily thrashing this stuff out
  and trying to write a paper on it, but it may or may not hang
  together.

  Jacobson's talk is based on a short paper he'd just been editing the
  galley proofs for; so it should come out soon:

  \begin{enumerate}
  \def\labelenumi{\arabic{enumi})}
  \setcounter{enumi}{1}
  \tightlist
  \item
    Ted Jacobson, ``1+1 sector of 3+1 gravity'', \emph{Class. Quant.
    Grav.} \textbf{13} (1996), L1--L6.
  \end{enumerate}

  Now around this time the Erwin Schr\"odinger Institute, where the
  workshop was being held, moved from its comfortable old spot on
  Pasteurgasse to a more spacious location on Boltzmanngasse, near the
  physics department. (In Germany the word ``Gasse'' means ``alley'',
  and one might find it disrespectful that Pasteur and Boltzmann have
  mere alleys named after them, but in Vienna even lots of large streets
  are called ``Gasse'', when in Germany they'd be called ``Strasse''.
  But then even the word for potato is different in Austria; it's all
  part of the charm of the place.) The move disrupted the schedule of
  the talks a bit, and it also seems to have disrupted my note-taking,
  which gets more sketchy from here on out. Some of the dates below
  might be a bit off.
\item
  \textbf{Thursday, July 25th} --- I spoke on ``Topological Quantum
  Field Theory''. I am always talking about this on This Week's Finds so
  I won't bore you with the details. Basically I summarized what is
  known about \(BF\) theory (a particular topological quantum field
  theory) in dimensions 2, 3, and 4, and the discrete formulation of
  \(BF\) theory where you chop spacetime into simplices and label the
  edges and so on with spins and the like --- so-called ``state sum
  models''. You can read more about this in
  \protect\hyperlink{week38}{``Week 38''}.

  Later that day, Jerzy Lewandowski spoke on ``Degenerate Metrics''.
  Being somewhat less degenerate than Ted Jacobson, he spoke about
  extending general relativity to cases where the inverse metric looks
  like \[
      \left(
        \begin{array}{cccc}
          -1&0&0&0
        \\0&1&0&0
        \\0&0&1&0
        \\0&0&0&0
        \end{array}
      \right)
    \] In other words, where the speed of light is zero only in the
  \(z\) direction. Basically what happens is that spacetime gets
  foliated with a lot of \(3\)-dimensional slices, and on each one you
  get the equations of \(3\)-dimensional general relativity.
\item
  \textbf{Friday, July 26th} --- Thomas Strobl spoke on
  \(2\)-dimensional gravity. I don't understand his work well enough yet
  to have anything much to say, but the most interesting thing about it
  to \emph{me} is that it allows one to see how quantum groups emerge
  from the \(G/G\) gauged Wess--Zumino--Witten model (a certain
  \(2\)-dimensional topological quantum field theory), by describing
  this theory as the quantization of a Poisson \(\sigma\)-model --- a
  field theory where the fields take values in a Poisson manifold. For
  more, try:

  \begin{enumerate}
  \def\labelenumi{\arabic{enumi})}
  \setcounter{enumi}{2}
  \item
    Peter Schaller and Thomas Strobl, A brief introduction to Poisson
    \(\sigma\)-models, available as
    \href{https://arxiv.org/abs/hep-th/9507020}{hep-th/9507020}.

    Peter Schaller and Thomas Strobl, Poisson \(\sigma\)-models: a
    generalization of 2d gravity-Yang--Mills systems, available
    as \href{https://arxiv.org/abs/hep-th/9411163}{hep-th/9411163}.
  \end{enumerate}

  Later, I had a great conversation with Mike Reisenberger and Carlo
  Rovelli on reformulating the loop representation of quantum gravity in
  terms of surfaces embedded in spacetime. This again touched upon my
  interest in relating string theory and the loop representation. They
  are writing a paper on this which should be on the preprint servers
  pretty soon, so I'll wait until then to talk about it.
\item
  \textbf{Saturday, July 27th} --- Carlo Rovelli explained some things
  about the problem of time to me.
\item
  \textbf{Monday, July 30th} --- I spoke about relative states and
  entanglement entropy in two-part quantum systems (see
  \protect\hyperlink{week27}{``Week 27''}) and the applications of these
  ideas to topological quantum field theory and quantum gravity. A lot
  of this came from my attempts to understand the relation between
  quantum gravity and Chern--Simons theory, and Lee Smolin's work where
  he tries to use this relation to derive the Bekenstein bound on the
  entropy of a system in terms of its surface area (see
  \protect\hyperlink{week56}{``Week 56''}).

  An interesting little fact that I needed to use is that if you have a
  two-part quantum system in a pure state --- a state of zero entropy
  --- the two parts, regarded individually, can themselves have entropy,
  but the entropies of the two parts are equal. I worked this out using
  the symmetry of the situation but Walter Thirring, who attended the
  talk, pointed out that it can also be derived from a wonderful general
  fact: the triangle inequality! Namely, if your two-part system has
  entropy \(S\), and the two parts individually have entropies \(S_1\)
  and \(S_2\), then \(S\) can never be less than \(|S_1 - S_2|\) or
  greater than \(S_1 + S_2\). (In classical mechanics it's also true
  that \(S\) can never be less than \emph{either} \(S_1\) \emph{or}
  \(S_2\), but this fails in quantum mechanics, where for example you
  can have \(S\) be zero but \(S_1 = S_2 > 0\).)
\item
  \textbf{Wednesday, August 1st} --- Full of excitement and new ideas, I
  somewhat regretfully left the workshop and flew to London. Then I
  spent most of August working at Imperial College, thanks to a kind
  offer of office space from Chris Isham. I had some nice talks with
  Isham and his students on quantum gravity and the decoherent histories
  approach to quantum mechanics. I'll say a bit about this in a while,
  but next Week I am going to talk about triality and the secret inner
  meaning of \(\mathrm{E}_8\).
\end{itemize}

\hypertarget{week89}{%
\section{September 17, 1996}\label{week89}}

This week I want to return to the tale of \(n\)-categories, from which I
have been taking a break during summer vacation. But first, here are a
few things about quantum gravity. Last time I mentioned Jorge Pullin's
newsletter on general relativity, ``Matters of Gravity''. I am pleased
to report that it is now available on the world-wide web:

\begin{enumerate}
\def\labelenumi{\arabic{enumi})}
\tightlist
\item
  Jorge Pullin, ed., \emph{Matters of Gravity}, September 1996 issue available as 
  \href{https://arxiv.org/abs/gr-qc/9609008}{\texttt{gr-qc/9609008}}.
\end{enumerate}

Anyone who wants to keep up with the latest news on general relativity
should certainly read ``Matters of Gravity'' and MacCallum's list.
MacCallum's list? Yes, I should've mentioned it earlier: it's a mailing
list where you can find out where the general relativity conferences
are, where the postdoctoral positions are, what the latest books are,
and so on.

\begin{enumerate}
\def\labelenumi{\arabic{enumi})}
\setcounter{enumi}{1}
\tightlist
\item
  MacCallum's gravity mailing list: to subscribe send polite email to \hfill \break
  \texttt{M.A.H.MacCallum@qmw.ac.uk}
\end{enumerate}

By the way, a bunch of math and physics preprints are available from the
Schr\"odinger Institute, including a lot of new stuff on quantum gravity
that came out of that workshop I've been talking about:

\begin{enumerate}
\def\labelenumi{\arabic{enumi})}
\setcounter{enumi}{2}
\item
  Erwin Schr\"odinger Institute preprint archive, available at
  \texttt{http://www.esi.ac.at/ESI-Preprints.html}. Recent preprints
  include:

  Abhay Ashtekar and Alejandro Corichi, ``Photon inner-product and the
  Gauss linking number''.

  Abhay Ashtekar, Donald Marolf, Jose Mourao and Thomas Thiemann,
  ``\(\mathrm{SU}(N)\) quantum Yang--Mills theory in 2 dimensions: a
  complete solution''.

  Hugo Fort, Rodolfo Gambini and Jorge Pullin, ``Lattice knot theory and
  quantum gravity in the loop representation'', also available as
  \href{https://arxiv.org/abs/gr-qc/9608033}{\texttt{gr-qc/9608033}}.

  Michael Reisenberger, ``A left-handed simplicial action for Euclidean
  GR''.

  Carlo Rovelli, ``Loop quantum gravity and black hole physics''.
\end{enumerate}

\noindent
I described the ideas behind some of these papers in
\protect\hyperlink{week85}{"Week 85} --
\protect\hyperlink{week88}{``Week 88''}. I didn't mention the paper by
Ashtekar and Corichi. It gives nice formula for the inner product in the
Hilbert space for photons in terms of the Gauss linking number --- a
thing that counts how many times one knot links another.

In its simplest form, the formula goes like this: say you have two
knots, and you do a line integral of the electric field around one of
them, and of the magnetic field around the other. You get two
observables which in the \emph{quantum} theory of electromagnetism do
not commute. So the uncertainty principle says you can't measure them
both exactly at once. In fact, the uncertainty in one times the
uncertainty in the other can't be less than \(\hbar/2\) times the
absolute value of the Gauss linking number of the two knots! A nice
blend of quantum theory and topology! This winds up also being relevant
to the photon inner product, because, as the experts out there should
know, the canonical commutation relations in a free field theory always
come from the imaginary part of the inner product in the single-particle
Hilbert space.

In \protect\hyperlink{week88}{``Week 88''} I also mentioned a talk by
Jerzy Lewandowski, which has now appeared as a preprint:

\begin{enumerate}
\def\labelenumi{\arabic{enumi})}
\setcounter{enumi}{3}
\tightlist
\item
  Jerzy Lewandowski and Jacek Wilsniewski, ``2+1 sector of 3+1
  gravity'', available as
  \href{https://arxiv.org/abs/gr-qc/9609019}{\texttt{gr-qc/9609019}}.
\end{enumerate}

Also, Lee Smolin has written a paper arguing that Thiemann's work has
trouble squaring with the positivity of energy and the existence of
long-range correlations (i.e., massless gravitons) that one might expect
from semi-classical approaches to quantum gravity.

\begin{enumerate}
\def\labelenumi{\arabic{enumi})}
\setcounter{enumi}{4}
\tightlist
\item
  Lee Smolin, ``The classical limit and the form of the Hamiltonian
  constraint in nonperturbative quantum gravity'', available as
  \href{https://arxiv.org/abs/gr-qc/9609034}{\texttt{gr-qc/9609034}}.
\end{enumerate}

This paper has sparked some controversy in the loop representation
community. Its arguments are heuristic rather than mathematically
rigorous, so one can certainly imagine ways to wriggle out of the
conclusions it tries to draw. Nonetheless I think it does a good service
by focusing attention on down-to-earth physical issues. If the more
mathematically inclined quantum gravity folks are able either to prove
\emph{or} refute Smolin's ideas, we'll have made lots of progress.

Smolin has also written a paper relating the loop representation to
string theory:

\begin{enumerate}
\def\labelenumi{\arabic{enumi})}
\setcounter{enumi}{5}
\tightlist
\item
  Lee Smolin, ``Three dimensional strings as collective coordinates of
  four dimensional quantum gravity'', available as
  \href{https://arxiv.org/abs/gr-qc/9609031}{\texttt{gr-qc/9609031}}.
\end{enumerate}

This paper really freaks me out, because it attempts to relate the loop
representation of quantum gravity in \(4\)-dimensional spacetime to
string theory in \emph{3-dimensional} spacetime. That's an idea that
never would have occurred to me. Smolin suggests it might possibly be
related to how supergravity in 11 dimensions is related to string theory
in 10 dimensions, but unfortunately I don't know enough about all that
to know where to go with it. I need to learn more about this string
theory duality stuff --- see \protect\hyperlink{week72}{``Week 72''} for
my pathetic attempts so far to understand it. I haven't read this yet,
but I should:

\begin{enumerate}
\def\labelenumi{\arabic{enumi})}
\setcounter{enumi}{6}
\tightlist
\item
  Michael Dine, ``String theory dualities'', available as
  \href{https://arxiv.org/abs/hep-th/9609051}{\texttt{hep-th/9609051}}.
\end{enumerate}

It's an expository article.

\begin{center}\rule{0.5\linewidth}{0.5pt}\end{center}

\hypertarget{week89_tale}{Okay, now let's go back to the tale of 
\(n\)-categories.} As promised, I will tell you all about monads, 
monoids, monoid objects, and monoidal categories.

You may or may not remember, but in \protect\hyperlink{week80_tale}{``Week
80''} I explained the idea of a ``\(2\)-category'' pretty precisely.
This is a gadget with a bunch of objects, a bunch of morphisms going
from one object to another, and a bunch of \(2\)-morphisms going from
one morphism to another. We write \(f\colon x\to y\) to denote a
morphism \(f\) from the object \(x\) to the object \(y\), and we write
\(F\colon f\Rightarrow g\) to denote a \(2\)-morphism \(F\) from the
morphism \(f\) to the morphism \(g\).

Just as in a category, in a \(2\)-category we can compose a morphism
\(f\colon x\to y\) with a morphism \(g\colon y\to z\) to get a morphism
\(fg\colon x\to z\). (Note that I write \(fg\) instead of \(gf\); I'm
going to use this ordering most of the time, though I may occasionally
change my mind just to confuse you more.) Similarly, we can compose a
\(2\)-morphism \(F\colon f\Rightarrow g\) with a 2-morphism
\(G\colon g\Rightarrow h\) to get a \(2\)-morphism
\(FG\colon f\Rightarrow h\). This is called ``vertical composition'' of
\(2\)-morphisms. We can visualize FG like this:
\[\includegraphics[scale=0.3]{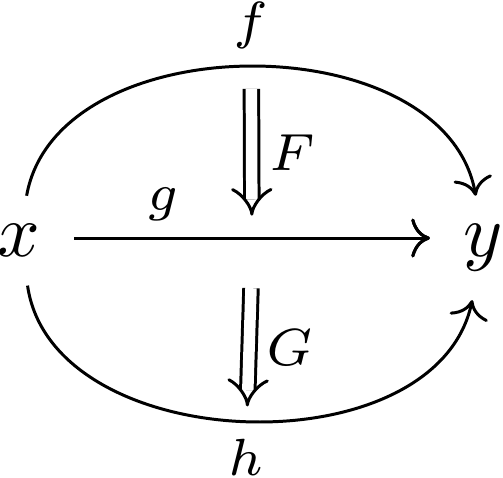}\] We stick \(F\) on
top of \(G\) to get \(FG\), which is why it's called ``vertical''
composition.

Also, if we have morphisms \(f,g\colon x\to y\) and
\(f',g'\colon y\to z\), and 2-morphisms \(F\colon f\Rightarrow g\) and
\(F'\colon f'\Rightarrow g'\), we can ``horizontally compose'' \(F\) and
\(F'\) to get \(F\cdot F'\colon ff'\Rightarrow gg'\). It looks like
this: \[\includegraphics[scale=0.3]{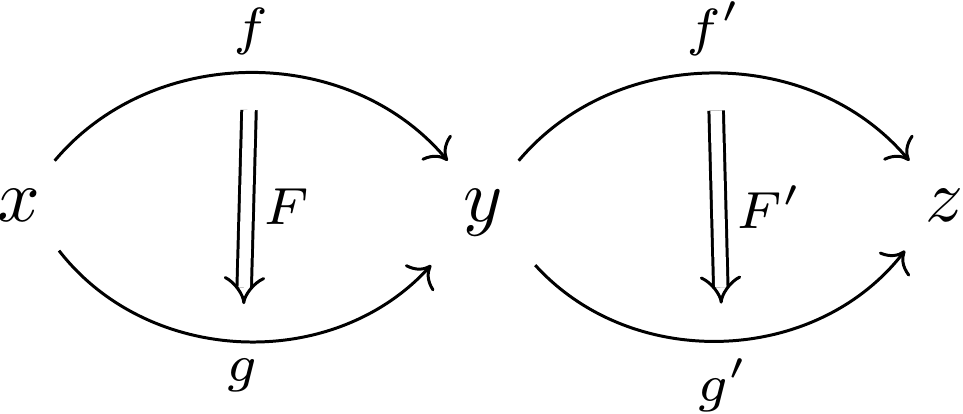}\] There are
some axioms all this stuff has to satisfy, which I described in
\protect\hyperlink{week80_tale}{``Week 80''}, but I won't repeat them here.
The main thing to keep in mind is that a \(2\)-category is like an
abstract 2-dimensional world... and the axioms for a \(2\)-category
are algebraic distillations of the rules for putting things together in
2 dimensions. In particular, you can put the \(2\)-morphisms together
side by side (horizontally) or one on top of the other (vertically), if
they fit.

Later I'll say more about what \(2\)-categories have to do with
2-dimensional physics, but right now I want to do something more
fundamental. I want to show how all sorts of concepts of
``multiplication'' or ``combination'' fit nicely into the framework of
\(2\)-categories. The basic idea is really simple: we often think of
multiplication as some sort of function \[M\colon s\times s\to s\] where
we take two elements \(a\) and \(b\) from some set \(s\), and
``multiply'' them to get a new one \(M(a,b)\). But we can visualize this
as follows: \[
  \begin{tikzpicture}
    \node (xl) at (0,0) {$\bullet$};
    \node (xt) at (1.25,2) {$\bullet$};
    \node (xr) at (2.5,0) {$\bullet$};
    \draw[thick] (xl) to node[fill=white]{$s$} (xt);
    \draw[thick] (xt) to node[fill=white]{$s$} (xr);
    \draw[thick] (xl) to node[fill=white]{$s$} (xr);
    \draw[-implies,double equal sign distance] (xt) to (1.25,0.2);
    \node at (1,0.7) {$M$};
  \end{tikzpicture}
\] I've drawn a triangular shaped gadget that takes two ``inputs'' from
the two slanted edges labelled \(s\), and spits out one ``output'' from
the horizontal edge labelled \(s\) on the bottom. It's clear from the
geometry here that \(M\) is something \(2\)-dimensional --- hence, a
\(2\)-morphism --- and that \(s\) is \(1\)-dimensional --- hence, a
morphism. Let's label the corners too: \[
  \begin{tikzpicture}
    \node (xl) at (0,0) {$x$};
    \node (xt) at (1.25,2) {$x$};
    \node (xr) at (2.5,0) {$x$};
    \draw[thick,->] (xl) to node[fill=white]{$s$} (xt);
    \draw[thick,->] (xt) to node[fill=white]{$s$} (xr);
    \draw[thick,->] (xl) to node[fill=white]{$s$} (xr);
    \draw[-implies,double equal sign distance] (xt) to (1.25,0.2);
    \node at (1,0.7) {$M$};
  \end{tikzpicture}
\] to make it clear that \(s\) is a morphism from \(x\) to itself. Here
\(x\), being 0-dimensional, is an object.

This hocus-pocus may seem mystifying, but if you bear with me and work
at it you'll see what I'm up to. I'm saying that essence of
``multiplication'' can be described very generally in a situation where
you have a \(2\)-category with an object \(x\) in it, a morphism
\(s\colon x\to x\), and a 2-morphism \(M\colon ss\Rightarrow s\). Often
we are interested in situations like this where the ``multiplication''
\(M\) is associative, meaning that the composite
\[sss\xRightarrow{M\cdot1_s}ss\xRightarrow{M}s\] equals
\[sss\xRightarrow{1_s\cdot M}ss\xRightarrow{M}s\] (Here
\(1_s\colon s\Rightarrow s\) is the identity \(2\)-morphism from \(s\)
to itself... the axioms for a \(2\)-category say that this exists.)
Also, we're often interested in situations where there is a
``multiplicative unit'', that is, a \(2\)-morphism \(I\colon 1_x\to s\)
for which \[s = 1_xs\xRightarrow{I\cdot1_s}ss\xRightarrow{M}s\] equals
\(1_s\), and so does
\[s = s1_x\xRightarrow{1_s\cdot I}ss\xRightarrow{M}s\] If we have a
\(2\)-category with stuff in it satisfying these rules, we say we have a
``monad'' in that \(2\)-category.

What is an example of a monad? Well, consider our original example where
s is a set and M is a function. We can think of this as living in a
\(2\)-category as follows. Our \(2\)-category will have only one object,
\(x\). The morphisms of this \(2\)-category are sets, and composing
morphisms corresponds to taking the cartesian product of sets. The
\(2\)-morphisms of this \(2\)-category are functions between sets.

What does a monad amount to in this case? Well, work it out! The
multiplicative unit \(1_x\) must corresponds to the one-element set;
\(s\) is some set; the \(2\)-morphism \(I\colon 1_x\Rightarrow s\) is a
function from the one-element set to \(s\), which picks out a special
\emph{element} of \(s\); the 2-morphism \(M\colon ss\Rightarrow s\) is
our multiplication operation. The axioms of a monad I gave then say that
this multiplication is associative and that the special element of \(s\)
is the multiplicative unit... that is, it serves as the left and
right identity for multiplication.

So we have a set with an associative multiplication and a unit for this
multiplication. That's what folks call a ``monoid'' --- see
\protect\hyperlink{week74_tale}{``Week 74''} for more on these. So a monoid
is a special sort of monad!

The point, however, is that there are lots of other kinds of monads, and
this \(2\)-categorical nonsense unifies the study of all of them.
Consider, for example, that trick we played of turning the category
\(\mathsf{Set}\) into a \(2\)-category with just one object \(x\). It's
a very versatile trick. In general, a \(2\)-category with just one
object is called a ``monoidal category'', because you can do this
relabelling trick: \[
  \begin{aligned}
    \text{2-morphisms} &\mapsto \text{morphisms}
  \\\text{morphisms} &\mapsto \text{objects}
  \\\text{objects} &\mapsto 
  \end{aligned}
\] You take the \(2\)-category with just one object, forget the object,
call the morphisms ``objects'' and the \(2\)-morphisms ``morphisms'',
and you've got a category! But one where you can compose or ``multiply''
or ``tensor'' objects, because they were secretly morphisms from \(x\)
to itself. For example, \(\mathsf{Set}\) is a monoidal category where we
can multiply objects (i.e., sets) with the cartesian product.

However, there are lots of other interesting monoidal categories. For
example, \(\mathsf{Vect}\) (the category of vector spaces) becomes a
monoidal category if we multiply vector spaces by tensoring them.
\(\mathsf{Top}\) (the category of topological spaces) becomes a monoidal
category if we multiply spaces by taking their cartesian product with
the usual product topology. \(\mathsf{Mon}\) (the category of monoids)
becomes a monoidal category if we multiply monoids by taking their cartesian
product. And so on\ldots.

Because a monoidal category is a \(2\)-category with one object, we can
talk about monads in any monoidal category. These are usually called
``monoid objects'', because they are like a monoid living in the
category in question. For example, a monoid object in \(\mathsf{Vect}\)
is an associative algebra. A monoid object in \(\mathsf{Top}\) is a
topological monoid.

Sometimes funny things happen: for example, a monoid object in
\(\mathsf{Mon}\) is a commutative monoid! This ``birth of
commutativity'' illustrates something called the ``Eckmann--Hilton
principle''. Some more sophisticated ramifications of this principle are
discussed in the following paper:

\begin{enumerate}
\def\labelenumi{\arabic{enumi})}
\setcounter{enumi}{7}
\tightlist
\item
  John Baez and Martin Neuchl, ``Higher-dimensional algebra I: braided
  monoidal \(2\)-categories'', \emph{Adv. Math.} \textbf{121} (1996),
  196--244. Also available as
  \href{http://arxiv.org/abs/q-alg/9511013}{\texttt{q-alg/9511013}}.
\end{enumerate}

We can get into some curious self-referential loops, too: the category
having (small) categories as objects and functors as morphisms becomes a
monoidal category with the ``cartesian product'' of categories as the
way to multiply objects... and a monoid object in this is a (small)
monoidal category! Try wrapping your brain around that! A monoid object
is something you define in a monoidal category, but a monoidal category
is itself a kind of monoid object! This illustrates something that James
Dolan and I call it the ``microcosm principle''. I should note at this
point --- I should have noted it before --- that most of this stuff
about category theory is stuff I learned from Dolan. We are writing a
paper in which we give a general definition of \(n\)-categories, and
explain this ``microcosm principle''.

Anyway, some of the most interesting monads live not in monoidal
categories but \(2\)-categories with lots of objects. The primordial
\(2\)-category is Cat, which has (small) categories as objects, functors
as morphisms and \emph{natural transformations} as \(2\)-morphisms. (A
minute ago I gave a way to think of \(\mathsf{Cat}\) as a monoidal
category. That was a bit different than this!) Monads in
\(\mathsf{Cat}\) are the first monads anyone called ``monads'', I
believe. You can read a bunch about them in the bible of category
theory:

\begin{enumerate}
\def\labelenumi{\arabic{enumi})}
\setcounter{enumi}{8}
\tightlist
\item
  \emph{Categories for the Working Mathematician}, by Saunders Mac Lane,
  Springer, Berlin, 1988.
\end{enumerate}
\noindent
Believe or not, monads in \(\mathsf{Cat}\) are nice way to think about
\emph{algebraic theories} --- a branch of logic perhaps pioneered by the
theory of ``universal algebra''. (My knowledge of the history here is
sort of fuzzy.) It would take me a while to explain this so I'll put it
off for next Week.

Let me just wrap up by saying that we can take this picture \[
  \begin{tikzpicture}
    \node (xl) at (0,0) {$x$};
    \node (xt) at (1.25,2) {$x$};
    \node (xr) at (2.5,0) {$x$};
    \draw[->] (xl) to node[fill=white]{$s$} (xt);
    \draw[->] (xt) to node[fill=white]{$s$} (xr);
    \draw[->] (xl) to node[fill=white]{$s$} (xr);
    \draw[-implies,double equal sign distance] (xt) to (1.25,0.2);
    \node at (1,0.7) {$M$};
  \end{tikzpicture}
\] and draw a ``dual'' picture like this: \[
  \begin{tikzpicture}
    \begin{knot}
      \strand[thick] (0,0.5)
        to (0,0)
        to [out=down,in=up] (0.5,-1)
        to (0.5,-2);
      \strand[thick] (1,0.5)
        to (1,0)
        to [out=down,in=up] (0.5,-1);
    \end{knot}
    \node[fill=white] at (0,0) {$s$};
    \node[fill=white] at (1,0) {$s$};
    \node[label=left:{$M$}] at (0.5,-0.95) {$\bullet$};
    \node[fill=white] at (0.5,-1.5) {$s$};
  \end{tikzpicture}
\] which illustrates perhaps more vividly how \(M\) is the process of
two copies of \(s\) getting squashed down into one copy. This sort of
picture is called a ``string diagram'', and it's literally the Poincar\'e
dual of the earlier picture, meaning that stuff that was
\(k\)-dimensional is now drawn as \((2-k)\)-dimensional. (The
0-dimensional object \(x\) is now the 2-dimensional ``background.'') For
more on string diagrams, see:

\begin{enumerate}
\def\labelenumi{\arabic{enumi})}
\setcounter{enumi}{9}
\tightlist
\item
  Ross Street, ``Categorical structures'', in \emph{Handbook of
  Algebra}, vol.~\textbf{1}, ed.~M. Hazewinkel, Elsevier, 1996.
\end{enumerate}

This diagram may also remind physicists (if any of them are still
reading this) of a Feynman diagram, in particular a 3-gluon vertex in
QCD. It's no coincidence! I'll have to say more about that later,
though.

To continue reading the ``Tale of \(n\)-Categories'', see
\protect\hyperlink{week92_tale}{``Week 92''}.

\hypertarget{week90}{%
\section{September 30, 1996}\label{week90}}

If you've been following This Week's Finds, you know that I'm in love
with symmetry. Lately I've been making up for my misspent youth by
trying to learn more about simple Lie groups. They are, roughly
speaking, the basic building blocks of the symmetry groups of physics.

In trying to learn about them, certain puzzles come up. In July I asked
Bertram Kostant about one that's been bugging me for years: ``Why does
\(\mathrm{E}_8\) exist?'' In a word, his answer was: ``Triality!'' This
was incredibly exciting to me; it completely blew my mind. But I should
start at the beginning\ldots.

In my youth, I found the classification of simple Lie groups to be
unintuitive and annoying. I still do, but over the years I've realized
that suffering through this classification theorem is the necessary
entrance fee to a whole world of symmetry. I gave a tour of this world
in \protect\hyperlink{week62}{``Week 62''} --
\protect\hyperlink{week65}{``Week 65''}, but here I want to make
everything as simple as possible, so I won't assume you've read that
stuff. Experts should jump directly to the end of this article and read
backwards until it becomes boring.

A Lie group is a group that can be given coordinates for which all the
group operations are infinitely differentiable. A good example is the
group \(\mathrm{SO}(n)\) of rotations in \(n\)-dimensional Euclidean
space. You can multiply rotations by doing first one and then the other,
or mathematically by doing matrix multiplication. Every rotation has an
inverse, given mathematically by the inverse matrix. Since matrices are
just bunches of numbers, you can coordinatize \(\mathrm{SO}(n)\), at
least locally, and in terms of these coordinates the operations of
multiplication and taking inverses are infinitely differentiable, or
``smooth'', so \(\mathrm{SO}(n)\) is a Lie group.

Using the magic of calculus, we can think of tangent vectors at the
identity element of \(\mathrm{SO}(n)\) as ``infinitesimal rotations''.
So for example, taking \(n = 3\), let's start with the rotation by the
angle \(t\) about the \(z\) axis, given by the matrix: 
\[
  \left(
    \begin{array}{ccc}
      \cos t & -\sin t & 0
    \\\sin t & \cos t & 0
    \\0 & 0 & 1
    \end{array}
  \right)
\] Then we can differentiate this and set \(t = 0\) to get an
``infinitesimal rotation about the \(z\) axis'': \[
  \left(
    \begin{array}{ccc}
      0 & -1 & 0
    \\1 & 0 & 0
    \\0 & 0 & 0
    \end{array}
  \right)
\] Let's call this \(J_z\), since it's very related to angular momentum
about the \(z\) axis. (Folks often throw in a factor of \(-i\) when they
define \(J_z\) in quantum mechanics, but let's not bother with that
here.)

Similarly we have \(J_x\) and \(J_y\). Now rotations about different
axes don't commute, so these infinitesimal rotations don't either. In
fact, we have \[
  \begin{aligned}
    J_x J_y - J_y J_x &= J_z,
  \\J_y J_z - J_z J_y &= J_x,
  \\J_z J_x - J_x J_z &= J_y.
  \end{aligned}
\]

If you have never done it, there are few things in life as rewarding at
this point as computing \(J_x\) and \(J_y\) for yourself and checking
the above ``commutation relations''.

Folks usually write the ``commutators'' on the left hand side using
brackets, like this: \[
  \begin{aligned}
    \,[J_x,J_y] &= J_z,
  \\ [J_y,J_z] &= J_x,
  \\ [J_z,J_x] &= J_y.
  \end{aligned}
\] These relations are lurking in the definition of quaternions and also
the vector cross product. Quaternions and cross products are good for
understanding rotations in \(3\)-dimensional space; they let us describe
infinitesimal rotations and their failure to commute. Here we are
calling a spade a spade and working directly with the algebra of
infinitesimal rotations, which folks call \(\mathfrak{so}(3)\). (For
related stuff, see \protect\hyperlink{week5}{``Week 5''}.)

Okay. The point is, we can do this trick for any Lie group! The space of
``infinitesimal group elements'', or more precisely tangent vectors at
the identity element of a Lie group, is called the ``Lie algebra of the
group''. It's a vector space whose dimension is the dimension of the
group, and it always has a bracket operation on it satisfying certain
axioms (listed in \protect\hyperlink{week3}{``Week 3''}).

The classification of Lie groups can be reduced to the classification of
Lie algebras, because the Lie algebra almost determines the Lie group.
More precisely, every Lie algebra is the Lie algebra of a unique Lie
group that is ``simply connected'' --- i.e., one for which every loop in
it can be continuously shrunk to a point. People understand how to get
from any Lie group to a simply connected one (called its ``universal
cover''), so if we understand simply connected Lie groups, we pretty
much understand all Lie groups. See \protect\hyperlink{week61}{``Week
61''} for an instance of this philosophy.

Now classifying Lie algebras is just a matter of heavy-duty linear
algebra. Let me explain what the ``simple'' Lie algebras are; you'll
have to take my word for it that understanding these is a big step
towards understanding all Lie algebras.

At one extreme in the world of Lie groups are the commutative, or
``abelian'' Lie groups. Here multiplication is commutative, so
\([x,y] = 0\) for all \(x\) and \(y\) in the Lie algebra of the group.
At the other extreme are the ``semisimple'' Lie groups. Here every
element in the Lie algebra is of the form \([x,y]\) for some \(x\) and
\(y\): roughly, if we bracket the whole Lie algebra with itself, we get
itself back again. The semisimple Lie algebras turn out to be incredibly
important in physics, where they are the typical ``gauge groups'' of
field theories.

The ``simple'' Lie algebras are the building blocks of the semisimple
ones: every semisimple Lie algebra can be broken down into pieces that
are simple. (Technically, we say it's a ``direct sum'' of simple Lie
algebras). We say a Lie group is simple if its Lie algebra is simple.

So: what are the simple Lie algebras? They were classified, thanks to
some heroic work by Killing and Cartan, in the early part of the 20th
century. To keep life simple (ahem) I'll only give the classification of
those simple Lie algebras whose corresponding Lie groups are
\emph{compact} --- meaning roughly that they are finite in size. (For
example, \(\mathrm{SO}(n)\) is compact.) It turns out that if we
understand the compact ones, we can understand the noncompact ones too.

So, here are the Lie algebras of the compact simple Lie groups! There
are 4 straightforward infinite families and 5 delightful and puzzling
exceptions. The 4 infinite families are easy to understand and are
called ``classical groups''. They are the workhorses of mathematics and
physics. The other 5 are called ``exceptional groups''. They have always
seemed very mysterious to me.

The 4 infinite families are:

\begin{itemize}
\tightlist
\item
  \(\mathrm{A}_n\): This is the Lie algebra of \(\mathrm{SU}(n)\), the
  group of \(n\times n\) complex matrices that preserve lengths (i.e.,
  are unitary) and have determinant \(1\). This Lie algebra is also
  called \(\mathfrak{su}(n)\).
\item
  \(\mathrm{B}_n\): This is the Lie algebra of \(\mathrm{SO}(2n+1)\),
  the group of \((2n+1)\times(2n+1)\) real matrices that preserve
  lengths (i.e., are orthogonal) and have determinant \(1\). This Lie
  algebra is also called \(\mathfrak{so}(2n+1)\).
\item
  \(\mathrm{C}_n\): This is the Lie algebra of \(\mathrm{Sp}(n)\), the
  group of \(n\times n\) quaternionic matrices that preserve lengths.
  This Lie algebra is also called \(\mathfrak{sp}(n)\).
\item
  \(\mathrm{D}_n\): This is the Lie algebra of \(\mathrm{SO}(2n)\), the
  group of \(2n\times 2n\) real matrices that preserve lengths and have
  determinant \(1\). This Lie algebra is also called
  \(\mathfrak{so}(2n)\).
\end{itemize}

You may justly wonder why the heck they are called \(\mathrm{A}_n\),
\(\mathrm{B}_n\), \(\mathrm{C}_n\), and \(\mathrm{D}_n\), and why we
separated out the even and odd cases of \(\mathrm{SO}(n)\) as we did!
This is explained in \protect\hyperlink{week64}{``Week 64''}, and I
don't want to worry about it here. Anyway, glossing over some nuances,
we see that these guys are all pretty much just groups of rotations in
real, complex, and quaternionic vector spaces.

The 5 exceptions are as follows:

\begin{itemize}
\tightlist
\item
  \(\mathrm{F}_4\): A 52-dimensional Lie algebra.
\item
  \(\mathrm{G}_2\): A 14-dimensional Lie algebra.
\item
  \(\mathrm{E}_6\): A 78-dimensional Lie algebra.
\item
  \(\mathrm{E}_7\): A 133-dimensional Lie algebra.
\item
  \(\mathrm{E}_8\): A 248-dimensional Lie algebra.
\end{itemize}

Here I am being rather reticent about what these Lie algebras --- or the
corresponding Lie groups, which go by the same names --- actually \emph{are!}
The reason is that it's not so easy to explain. One can certainly
describe the exceptional Lie groups as groups of matrices with certain
complicated properties, but often this is done in a way that leaves one
utterly puzzled as to the real reason why these simple Lie groups exist.

Of course, the answer to ``why'' a mathematical object exists is a
matter of taste. You may feel satisfied if you can easily construct it
from other objects you know and love, or you may feel satisfied once it
is so tightly woven into your overall scheme of things that you can't
imagine life without it.

In any case, I have long been asking people why the exceptional Lie
groups exist, but without much luck. Until recently I only felt happy
about one of them, the one called \(\mathrm{G}_2\): it's the group of
rotations of the octonions! The real numbers, complex numbers,
quaternions and octonions are the only ``normed division algebras'' ---
a property which makes it easy to define rotation groups --- but the
octonions are weirder than the other three because, unlike the others,
they are not associative. (See \protect\hyperlink{week59}{``Week 59''}
and \protect\hyperlink{week61}{``Week 61''} for details.) One might
expect a series of simple Lie groups coming from rotations in octonionic
vector spaces, like the other classical series... but there isn't
one! The only simple Lie group like this is the group of rotations of a
\emph{one}-dimensional octonionic vector space, \(\mathrm{G}_2\). (More
precisely, we say that \(\mathrm{G}_2\) is the group of automorphisms of
the octonions, that is, the linear transformations that preserve the
octonion product. These all preserve lengths.)

The idea that the exceptional groups are all related to octonions is
sort of pleasing, because one might easily \emph{expect} that the reals,
complexes and quaternions give nice infinite series of ``classical'' Lie
groups, while the octonions, being much more bizarre, give only 5
bizarre ``exceptional'' Lie groups. Indeed, in
\protect\hyperlink{week64}{``Week 64''} I described how \(\mathrm{F}_4\)
and \(\mathrm{E}_6\) are related to the octonions... but in a
pretty complicated way! As for \(\mathrm{E}_7\) and \(\mathrm{E}_8\),
here until recently I had always been completely in the dark. This is
all the more irksome because the biggest, most mysterious exceptional
Lie group of all, \(\mathrm{E}_8\), plays an important role in string
theory!

Luckily, on Thursday July 11th I ran into Bertram Kostant, who had been
attending the previous workshop here at the Erwin Schr\"odinger
Institute. As I described in \protect\hyperlink{week79}{``Week 79''},
Kostant is one of the expert's experts on group theory. So I got up my
nerve and asked him, ``Why does \(\mathrm{E}_8\) exist?'' And he told
me! Best of all, he explained both \(\mathrm{E}_8\) and \(\mathrm{F}_4\)
in terms of a principle that I knew was crucial for understanding
\(\mathrm{G}_2\) and the octonions: the principle of triality!

I sketched a description of triality in
\protect\hyperlink{week61}{``Week 61''}. Let me just summarize the idea
here. One of the main way to understand Lie algebras is to understand
their ``representations''. A representation of a Lie algebra is simply a
function from it to the space of \(n\times n\) matrices that preserves
the bracket operation. (The \(n\times n\) matrices form a Lie algebra
with the commutator as the bracket operation.) For example,
\(\mathfrak{so}(n)\) has a representation where we map each element to
an \(n\times n\) matrix in the most utterly obvious way: each element \emph{is}
an \(n\times n\) matrix, so don't do anything to it! This is called the
``vector'' representation, because this is how we do infinitesimal
rotations to vectors. But \(\mathfrak{so}(n)\) also has representations
called ``spinor'' representations. In physics, the vector representation
describes spin-\(1\) particles, while the spinor representations
describe spin-\(1/2\) particles.

Spinor representations work differently depending on whether the
dimension \(n\) is even or odd. (This is one reason why people
distinguish the even and odd \(n\) cases of \(\mathfrak{so}(n)\) in that
classification of simple Lie algebras above!) When \(n\) is odd there is one
spinor representation. That's why in ordinary \(3\)-dimensional space
there is just one kind of spinor to worry about, as you learn when you
learn about spin-\(1/2\) particles in undergraduate quantum mechanics.
When \(n\) is even there are two different spinor representations, called
the ``left-handed'' and ``right-handed'' spinor representations. This
shows up when you do quantum mechanics taking special relativity --- and
\(4\)-dimensional spacetime --- into account. For example, the way
neutrinos transform under rotations is described by the left-handed
spinor representation, while anti-neutrinos are described by
right-handed spinors.

When \(n\) is even, both the spinor representations of
\(\mathfrak{so}(n)\) are of dimension \(2^{n/2 - 1}\). That is, they are
functions from \(\mathfrak{so}(n)\) to the space of
\(2^{n/2 - 1} \times 2^{n/2 - 1}\) matrices. Now something marvelous
happens when \(n = 8\). Namely, \(2^{n/2 - 1} = n\), so the spinor
representations are just as big as the vector representation. This might
lead one to hope that in some sense they are ``the same'' as the vector
representation. This is actually true, but in a subtle way\ldots. they
are not ``equivalent'' representations in the standard sense of Lie
algebra theory, but something sneakier is true.

The Lie algebra \(\mathfrak{so}(8)\) has interesting symmetries! It has
a little symmetry group with 6 elements, the same as the symmetries of a
equilateral triangle, and using these 6 symmetries we can permute the
vector, left-handed spinor, and right-handed spinor representations into
each other however we please!

For example, one of these symmetries switches the left-handed and
right-handed spinor representations, but leaves the vector
representation alone. Actually, this symmetry works in any even
dimension, not just dimension 8. Its analogue in \(4\)-dimensional
spacetime is called ``parity'', a symmetry that turns left-handed
particles into right-handed ones and vice versa. The fact that there are
no right-handed neutrinos means that the laws of nature do not actually
have this symmetry... but it's still very important in math and
physics.

What's special about dimension 8 is that there are symmetries switching
the vector representation and the spinor representations. For example:
if we take an element \(x\) of \(\mathfrak{so}(8)\), apply the right
symmetry of \(\mathfrak{so}(8)\) to turn it into another element of
\(\mathfrak{so}(8)\), and then use the right-handed spinor
representation to it to turn it into a matrix, we get the same thing as
if we just used the vector representation to turn \(x\) into a matrix.

Now \(\mathfrak{so}(8)\) is the Lie algebra of the Lie group
\(\mathrm{SO}(8)\), but \(\mathrm{SO}(8)\) is not ``simply connected''
in the sense defined above. The simply connected group whose Lie algebra
is \(\mathrm{SO}(n)\) is called \(\mathrm{Spin}(n)\). I gave an
introduction to these ``spin groups'' in
\protect\hyperlink{week61}{``Week 61''}, and I don't want to say much
about them here, except for this: the triality symmetries of
\(\mathfrak{so}(8)\) do not give symmetries of \(\mathrm{SO}(8)\), but
they do give symmetries of \(\mathrm{Spin}(8)\). Experts say the group
of outer automorphisms modulo inner automorphisms of \(\mathrm{SO}(8)\)
is \(S_3\) (the group of permutations of 3 things).

Pretty sneaky, how a group of symmetries can have its own group of
symmetries, no? As we'll now see, this is what gives birth to
\(\mathrm{G}_2\), \(\mathrm{F}_4\), \(\mathrm{E}_8\), and the octonions.

To get \(\mathrm{G}_2\) is pretty simple; we look at those elements of
\(\mathrm{Spin}(8)\) that are fixed (i.e., unaffected) by all the
triality symmetries, and these form a subgroup, which is
\(\mathrm{G}_2\).

For the rest, we need one more fact: there is a way to ``multiply'' a
left-handed spinor and a right-handed spinor and get a vector. This is
true in all even dimensions, not just \(n = 8\), so in particular it is
familiar to particle theorists who live in \(4\)-dimensional spacetime.
As I noted, what happens to a neutrino when you rotate (or Lorentz
transform) it is described using left-handed spinors, while
anti-neutrinos are described by right-handed spinors. Similarly, photons
are described by vectors. So as far as \emph{rotational} properties go,
one could think of a photon as a bound state of a neutrino and an
antineutrino. This led Schr\"odinger (or someone) to propose at one point
that photons were actually neutrino-antineutrino pairs. Subsequent
experiments showed this theory has lots of problems, and nobody sane
believes it any more. Still, it's sort of cute.

Now, in 8 dimensions, it shouldn't be surprising that we can also
multiply a left-handed spinor and a vector to get a right-handed spinor,
and so on. The point is, you can just use triality to permute the three
representations whichever way you please... they are not really all
that different.

So in particular, you can multiply two \(8\)-dimensional vectors and get
another vector. And this gives us the octonions!

Now how about \(\mathrm{F}_4\) and \(\mathrm{E}_8\)? This is the cool
stuff Kostant told me about. Here I will describe the Lie algebras, not
the Lie groups.

Let's call the right-handed and left-handed spinor representations
\(S_+\) and \(S_-\), respectively. (Us left-handers are always getting
shafted, being ``sinister'' rather than ``dextrous'' and all that, so we
get \(S_-\) rather than \(S_+\).) And let's call the vector
representation \(V\). And let's be sloppy, the way people usually are,
and also use these letters to stand for the \(8\)-dimensional vector
spaces on which \(\mathfrak{so}(8)\) acts as transformations.

Now let's form the direct sum of vector spaces
\[\mathfrak{so}(8)\oplus S_+ \oplus S_- \oplus V\] 
A vector in this vector space is just a list consisting of a guy in \(\mathfrak{so}(8)\),
a guy in \(S_+\), a guy in \(S_-\), and a guy in \(V\). The dimension of
this vector space is therefore 
\[28+8+8+8=52\] 
since it takes
\(n(n-1)/2\) numbers to describe a rotation in \(n\) dimensions. Hey!
Look! 52 is the dimension of \(\mathrm{F}_4\)! So maybe this thing is
\(\mathrm{F}_4\).

Yes, it is! Here's how it works. To make this gadget into a Lie algebra
--- which turns out to be \(\mathrm{F}_4\) --- we need a way to take the
``bracket'' of any two elements in it. We already know how to take the
bracket of two guys in \(\mathfrak{so}(8)\), so that's no problem. Since
\(\mathfrak{so}(8)\) acts as transformations of \(S_+\) and \(S_-\) and
\(V\), we also know how to multiply a guy in \(\mathfrak{so}(8)\) by one
of these other guys. We also know how to multiply a guy in \(S_+\) by a
guy in \(S_-\) to get a guy in \(V\), and so on. Finally, we can
multiply two guys in \(V\) to get a guy in \(\mathfrak{so}(8)\) as
follows: two vectors determine an infinitesimal rotation which starts
rotating the first vector in the direction of the second. (More
technically, we say that \(\mathfrak{so}(8)\) is isomorphic to the
second exterior power of \(V\), so we can multiply two guys in \(V\) to
get a guy in \(\mathfrak{so}(8)\) using the wedge product.) Using
triality, we can equally well multiply two guys in \(S_+\) to get a guy
in \(\mathfrak{so}(8)\), or multiply two guys in \(S_-\) to get a guy in
\(\mathfrak{so}(8)\).

So taking all these multiplication operations together, we can cook up a
way to take the bracket of any two guys in
\(\mathfrak{so}(8)\oplus S_+\oplus S_-\oplus V\) and get another such
guy. If you do it right --- I've been pretty vague, so I leave it to you
to fill in the details --- you can get this bracket to satisfy the Lie
algebra axioms, and you get \(\mathrm{F}_4\)!

Emboldened with our success, we now look at the vector space
\[\mathfrak{so}(8)\oplus\mathfrak{so}(8)\oplus\operatorname{End}(S_+)\oplus\operatorname{End}(S_-)\oplus\operatorname{End}(V).\]
Here \(\operatorname{End}(S_+)\) is the space of all linear
transformations of the vector space \(S_+\), so if you like, it's just
the space of \(8\times8\) matrices. Similarly for
\(\operatorname{End}(S_-)\) and \(\operatorname{End}(V)\). Now the
dimension of this space is \[28+28+64+64+64=248\] Hey! This is just the
dimension of \(\mathrm{E}_8\)! Maybe this space is \(\mathrm{E}_8\)!

Yes indeed. Again, you can cook up a bracket operation on this space
using all the stuff we've got. Here's the basic idea.
\(\operatorname{End}(S_+)\), \(\operatorname{End}(S_-)\), and
\(\operatorname{End}(V)\) are already Lie algebras, where the bracket of
two guys \(x\) and \(y\) is just the commutator \([x,y]=xy-yx\), where
we multiply using matrix multiplication. Since \(\mathfrak{so}(8)\) has
a representation as linear transformations of \(V\), it has two
representations on \(\operatorname{End}(V)\), corresponding to left and
right matrix multiplication; glomming these two together we get a
representation of \(\mathfrak{so}(8)\oplus\mathfrak{so}(8)\) on
\(\operatorname{End}(V)\). Similarly we have representations of
\(\mathfrak{so}(8)\oplus\mathfrak{so}(8)\) on
\(\operatorname{End}(S_+)\) and \(\operatorname{End}(S_-)\). Putting all
this stuff together we get a Lie algebra, if we do it right --- and it's
\(\mathrm{E}_8\). At least that's what Kostant said; I haven't checked
it.

So now we see, at least roughly, how triality gives birth to the
octonions, \(\mathrm{G}_2\), \(\mathrm{F}_4\), and \(\mathrm{E}_8\).
That leaves \(\mathrm{E}_8\)'s ``little brothers'' \(\mathrm{E}_6\) and
\(\mathrm{E}_7\). These are contained in \(\mathrm{E}_8\) as Lie
subalgebras, but apart from that I don't know any especially beautiful
way to get ahold of them, except for the way to get \(\mathrm{E}_6\)
from \(3 \times 3\) matrices of octonions, which I described in
\protect\hyperlink{week64}{``Week 64''}.

For some references to this stuff, try:

\begin{enumerate}
\def\labelenumi{\arabic{enumi})}
\item
  Claude C. Chevalley, \emph{The Algebraic Theory of Spinors}, Columbia
  U.\ Press, New York, 1954.
\item
  F. Reese Harvey, \emph{Spinors and Calibrations}, Perspectives in
  Mathematics, \textbf{9}, Academic Press, Boston, 1990.
\item
  Ian R. Porteous, \emph{Topological Geometry}, Cambridge
  U.\ Press, Cambridge, 1981.
\item
  Ian R. Porteous, \emph{Clifford Algebras and the Classical Groups},
  Cambridge U.\ Press, Cambridge, 1995.
\item
  Hans Freudenthal and H. de Vries, \emph{Linear Lie Groups}, Academic
  Press, New York, 1969.
\item
  Alex J. Feingold, Igor B. Frenkel, and John F. X. Rees, \emph{Spinor
  Construction of Vertex Operator Algebras, Triality, and
  \(\mathrm{E}_8^{(1)}\)}, Contemp.\ Math.\ \textbf{121}, AMS, Providence,
  1991.
\end{enumerate}
\noindent
I haven't looked at all these books lately, and the only source I
\emph{know} contains the above construction of \(\mathrm{E}_8\) from
triality is the last one, by Feingold, Frenkel, and Rees.

Now let me allow myself to get a bit more technical.

I am still not entirely happy, by any means, because what I'd really
like would be a simple explanation of why these exceptional simple Lie
algebras arise from triality, \emph{and no others}. In other words, I'd
like a classification of the simple Lie algebras that proceeded not by
the usual exhaustive (and exhausting) case-by-case study of Dynkin
diagrams, but by some less combinatorial and more ``synthetic''
approach. For example, it would be nice to really see a good explanation
of how the reals, the complexes, the quaternions and octonions each give
rise to a family of simple Lie algebras, and one gets \emph{all} of them
this way.

On the other hand, don't think I'm knocking the Dynkin diagram stuff. As
I explained in \href{week62.html}{``Week 62''} --
\protect\hyperlink{week64}{``Week 64''}, what's really fundamental to
the Dynkin diagram approach seems to be the not the Lie algebras
themselves but their root lattices. Taking lattices as fundamental to
the study of symmetry \emph{does} seem to be a good idea, since it gets
you to not just the simple Lie algebras described above, but also the
``Kac--Moody algebras'' so important in string theory and other forms of
\(2\)-dimensional physics, as well as marvelous things like the Leech
lattice and the Monster group.

The Dynkin diagram approach also makes it clear \emph{why} triality
exists: symmetries of Dynkin diagrams always give outer automorphisms of
the corresponding Lie algebras, and as you examine the Dynkin diagrams
of \(\mathrm{D}_n\), you get \[
  \begin{tikzpicture}
    \node at (2,0) {$\mathrm{D}_2=\mathfrak{so}(4)$};
    \node at (3.25,0) {$=$};
    \draw[thick] (4,1) node {$\bullet$};
    \draw[thick] (4,-1) node {$\bullet$};
  \end{tikzpicture}
\] \[
  \begin{tikzpicture}
    \node at (1,0) {$\mathrm{D}_3=\mathfrak{so}(6)$};
    \node at (2.25,0) {$=$};
    \draw[thick] (3,0) node {$\bullet$};
    \draw[thick] (3,0) to (4,1) node {$\bullet$};
    \draw[thick] (3,0) to (4,-1) node {$\bullet$};
  \end{tikzpicture}
\] \[
  \begin{tikzpicture}
    \node at (0,0) {$\mathrm{D}_4=\mathfrak{so}(8)$};
    \node at (1.25,0) {$=$};
    \draw[thick] (2,0) node{$\bullet$} to (3,0) node {$\bullet$};
    \draw[thick] (3,0) to (4,1) node {$\bullet$};
    \draw[thick] (3,0) to (4,-1) node {$\bullet$};
  \end{tikzpicture}
\] \[
  \begin{tikzpicture}
    \node at (-1,0) {$\mathrm{D}_6=\mathfrak{so}(10)$};
    \node at (0.25,0) {$=$};
    \draw[thick] (1,0) node{$\bullet$} to (2,0) node{$\bullet$} to (3,0) node {$\bullet$};
    \draw[thick] (3,0) to (4,1) node {$\bullet$};
    \draw[thick] (3,0) to (4,-1) node {$\bullet$};
  \end{tikzpicture}
\] and you can just \emph{see} how when you get to \(\mathfrak{so}(8)\)
there is that amazing triality symmetry, flashing briefly into being
before reverting to the boring old duality symmetry which only
interchanges the left-handed and right-handed spinor representations,
corresponding to the two dots on the far right of the Dynkin diagram.
(The dot on the far left corresponds to the vector representation.)

Of course, people don't usually talk about \(\mathrm{D}_2\) or
\(\mathrm{D}_3\), because \(\mathrm{D}_2\) is two copies of
\(\mathrm{A}_1\), and \(\mathrm{D}_3\) is the same as \(\mathrm{A}_3\).
However, there is no shame in doing so, and indeed a lot of insight to
be gained: the fact that \(\mathrm{D}_2\) consists of two copies of
\(\mathrm{A}_1\) corresponds to the isomorphism
\[\mathfrak{so}(4) = \mathfrak{su}(2) \oplus \mathfrak{su}(2),\] while
the fact that \(\mathrm{D}_3\) is the same as \(\mathrm{A}_3\)
corresponds to the isomorphism \[\mathfrak{so}(6) = \mathfrak{su}(4).\]

Each of these could easily serve as the springboard for a very long and
interesting discussion. However, I will refrain. Here let me simply note
that you can always ``fold'' a Dynkin diagram using one of its
symmetries, and if you do this to \(\mathrm{D}_4\) using triality you go
from \[
  \begin{tikzpicture}
    \node at (0,0) {$\mathrm{D}_4$};
    \node at (0.75,0) {$=$};
    \draw[thick] (1.5,0) node{$\bullet$} to (3,0) node {$\bullet$};
    \draw[thick] (3,0) to (4,1) node {$\bullet$};
    \draw[thick] (3,0) to (4,-1) node {$\bullet$};
  \end{tikzpicture}
\] down to \[
  \begin{tikzpicture}
    \node at (0,0) {$\mathrm{G}_2$};
    \node at (0.75,0) {$=$};
    \draw[ thick] (1.5,0) node{$\bullet$} to node[label=above:{6}]{\textgreater} (3,0) node {$\bullet$};
  \end{tikzpicture}
\] (Here the number 6 means that the two roots are at an angle of
\(\pi/6\) from each other. People usually just draw a triple line to
indicate this. The arrow points from the long root to the shorter root.)
This corresponds to how \(\mathrm{G}_2\) is the subgroup of
\(\mathrm{Spin}(8)\) consisting of elements that are invariant under
triality. You can also go from \[
  \begin{tikzpicture}
    \node at (-1.5,0) {$\mathrm{E}_6$};
    \node at (-0.75,0) {$=$};
    \draw[thick] (0,0) node{$\bullet$} to (1,0) node{$\bullet$} to (2,0) node{$\bullet$} to (3,0) node {$\bullet$} to (4,0) node {$\bullet$};
    \draw[thick] (2,0) to (2,1) node{$\bullet$};
  \end{tikzpicture}
\] down to
\[
  \begin{tikzpicture}
  \node at (-1.5,0) {$\mathrm{F}_4$};
    \node at (-0.75,0) {$=$};
    \draw[thick] (0,0) node{$\bullet$} to (1,0) node{$\bullet$} to node[label=above:{$4$}]{\textgreater} (2,0) node{$\bullet$} to (3,0) node {$\bullet$};
  \end{tikzpicture}
\]
by folding along the reflection symmetry. And Friedrich Knop told me
a neat way to get triality symmetry \emph{from} \(\mathrm{F}_4\), if you
happen to have \(\mathrm{F}_4\) around: the long roots of
\(\mathrm{F}_4\) form a root system of type \(\mathrm{D}_4\), which
defines an embedding of \(\mathrm{Spin}(8)\) into the Lie group
\(\mathrm{F}_4\) (more precisely, the compact real form). On the other
hand, the two short simple roots define an embedding of
\(\mathrm{SU}(3)\) in \(\mathrm{F}_4\). The Weyl group of
\(\mathrm{SU}(3)\) is \(S_3\) and can be lifted to \(\mathrm{SU}(3)\),
so we have an \(S_3\) subgroup of \(\mathrm{F}_4\). This acts by
conjugation on the \(\mathrm{Spin}(8)\) subgroup, implementing the
triality symmetries!

But I digress. My main point is, the Dynkin diagram symmetries do give a
nice way to understand outer automorphisms of simple Lie groups, and
these provide some important ties between simple Lie algebras, including
triality, which links the ``classical'' world to the ``exceptional''
world. But it is also nice to try to understand these in a somewhat more
``conceptual'' way. This is one of the reasons I'm interested in
2-Hilbert spaces: they seem to help one understand this stuff
from a new angle. But more on those, later. They tie into the
\(n\)-category stuff I'm always talking about. I will return to that
tale soon, and I'll keep building up some of the tools we need, until we
are ready to launch into a description of 2-Hilbert spaces.

In writing this Week's Finds, I benefitted greatly from email
correspondence with Robt Bryant, Christopher Henrich, Geoffrey Mess,
Friedrich Knop, and others.

\hypertarget{week91}{%
\section{October 6, 1996}\label{week91}}

For a while now I've been meaning to finish talking about monads and
adjunctions, and explain what that has to do with the 4-color theorem.
But first I want to say a little bit more about ``triality'', which was
the subject of \protect\hyperlink{week90}{``Week 90''}.

Triality is a cool symmetry of the infinitesimal rotations in
8-dimensional space. It was only last night, however, that I figured out
what triality has to do with 3 dimensions. Since it's all about the
number \emph{three} obviously triality should originate in the
symmetries of \emph{three}-dimensional space, right? Well, maybe it's
not so obvious, but it does. Here's how.

Take good old three-dimensional Euclidean space with its usual basis of
unit vectors \(i\), \(j\), and \(k\). Look at the group of all
permutations of \(\{i,j,k\}\). This is a little 6-element group which
people usually call \(S_3\), the ``symmetric group on 3 letters''.

Every permutation of \(\{i,j,k\}\) defines a linear transformation of
three-dimensional Euclidean space in an obvious way. For example the
permutation \(p\) with \[p(i) = j, \quad p(j) = k, \quad p(k) = i\]
determines a linear transformation, which we'll also call \(p\), with
\[p(ai+ bj + ck) = aj + bk + ci.\] In general, the linear
transformations we get this way either preserve the cross product, or
switch its sign. If \(p\) is an even permutation we'll get
\[p(v)\times p(w) = p(v\times w)\] while if \(p\) is odd we'll get
\[p(v)\times p(w) = -p(v\times w) = p(w\times v).\] That's where
triality comes from. But now let's see what it has to do with
\emph{four}-dimensional space. We can describe four-dimensional space
using the quaternions. A typical quaternion is something like
\[a + bi + cj + dk\] where \(a\), \(b\), \(c\), \(d\) are real numbers,
and you multiply quaternions by using the usual rules together with the
rules \[
  \begin{gathered}
    i^2 = j^2 = k^2 = -1
  \\ij=k,\quad jk=i,\quad ki=j,
  \\ji=-k,\quad kj=-i,\quad ik=-j.
  \end{gathered}
\] Now, any permutation \(p\) of \(\{i,j,k\}\) also determines a linear
transformation of the quaternions, which we'll also call \(p\). For
example, the permutation \(p\) I gave above has
\[p(a + bi + cj + dk) = a + bj + ck + di.\] The quaternion product is
related to the vector cross product, and so one can check that for any
quaternions \(q\) and \(q'\) we get \[p(qq') = p(q)p(q')\] if \(p\) is
even, and \[p(q'q) = p(q')p(q)\] if \(p\) is odd. So we are getting
triality to act as some sort of symmetries of the quaternions.

Now sitting inside the quaternions there is a nice lattice called the
``Hurwitz integral quaternions''. It consists of the quaternions
\(a + bi + cj + dk\) for which either \(a\), \(b\), \(c\), \(d\) are all
integers, or all half-integers. Here I'm using physics jargon, and
referring to any number that's an integer plus \(1/2\) as a
``half-integer''. A half-integer is \emph{not} any number that's half an
integer!

You can think of this lattice as the \(4\)-dimensional version of all
the black squares on a checkerboard. One neat thing is that if you
multiply any two guys in this lattice you get another guy in this
lattice, so we have a ``subring'' of the quaternions. Another neat thing
is that if you apply any permutation of \(\{i,j,k\}\) to a guy in this
lattice, you get another guy in this lattice --- this is easy to see. So
we are getting triality to act as some sort of symmetries of this
lattice. And \emph{that} is what people \emph{usually} call triality.

Let me explain, but now let me use a lot of jargon. (Having shown it's
all very simple, I now want to relate it to the complicated stuff people
usually talk about. Skip this if you don't like jargon.) We saw how to
get \(S_3\) to act as automorphisms and antiautomorphisms of
\(\mathbb{R}^3\) with its usual vector cross product... or
alternatively, as automorphisms and antiautomorphisms of the Lie algebra
\(\mathfrak{so}(3)\). From that we got an action as automorphisms and
antiautomorphisms of the quaternions and the Hurwitz integral
quaternions. But the Hurwitz integral quaternions are just a differently
coordinatized version of the \(4\)-dimensional lattice \(D_4\)! So we
have gotten triality to act as symmetries of the \(D_4\) lattice, and
hence as automorphisms of the Lie algebra \(D_4\), or in other words
\(\mathfrak{so}(8)\), the Lie algebra of infinitesimal rotations in 8
dimensions. (For more on the \(D_4\) lattice see
\protect\hyperlink{week65}{``Week 65''}, where I describe it using
different, more traditional coordinates.)

Actually I didn't invent all this stuff, I sort of dug it out of the
literature, in particular:

\begin{enumerate}
\def\labelenumi{\arabic{enumi})}
\tightlist
\item
  John H. Conway and Neil J. A. Sloane, \emph{Sphere Packings, Lattices
  and Groups}, Grundlehren der mathematischen
  Wissenschaften \textbf{290}, Springer, Berlin, 1993.
\end{enumerate}

and

\begin{enumerate}
\def\labelenumi{\arabic{enumi})}
\setcounter{enumi}{1}
\tightlist
\item
  Frank D. (Tony) Smith, ``Sets and \(\mathbb{C}^n\); quivers and
  ADE; triality; generalized supersymmetry; and
  \(\mathrm{D}_4\)-\(\mathrm{D}_5\)-\(\mathrm{E}_6\)'', available as
  \href{https://arxiv.org/abs/hep-th/9306011}{\texttt{hep-th/9306011}}.
\end{enumerate}
\noindent
But I've never quite seen anyone come right out and admit that triality
arises from the permutations of the unit vectors \(i\), \(j\), and \(k\)
in 3d Euclidean space.

I should add that Tony Smith has a bunch of far-out stuff about
quaternions, octonions, Clifford algebras, triality, the \(D_4\) lattice
--- you name it! --- on his home page:

\begin{enumerate}
\def\labelenumi{\arabic{enumi})}
\setcounter{enumi}{2}
\tightlist
\item
  Tony Smith's home page, \url{http://valdostamuseum.org/hamsmith/}
\end{enumerate}

He engages in more free association than is normally deemed proper in
scientific literature --- you may raise your eyebrows at sentences like
``the Tarot shows the Lie algebra structure of the
\(\mathrm{D}_4\)-\(\mathrm{D}_5\)-\(\mathrm{E}_6\) model, while the 
I Ching shows its
Clifford algebra structure'' --- but don't be fooled; his mathematics is
solid. When it comes to the physics, I'm not sure I buy his theory of
everything, but that's not unusual: I don't think I buy \emph{anyone's}
theory of everything!

Let me wrap up by passing on something he told me about triality and the
exceptional groups. In \protect\hyperlink{week90}{``Week 90''} I
described how you could get the Lie groups \(\mathrm{G}_2\),
\(\mathrm{F}_4\) and \(\mathrm{E}_8\) from triality. I didn't know how
\(\mathrm{E}_6\) and \(\mathrm{E}_7\) fit into the picture. He emailed
me, saying:

\begin{quote}
"Here is a nice way: Start with \(\mathrm{D}_4 = \mathrm{Spin}(8)\):

\[28 =  28  +   0  +   0  +   0  +   0  +   0  +   0\]

Add spinors and vector to get \(\mathrm{F}_4\):

\[52 =  28  +   8  +   8  +   8  +   0  +   0  +   0\]

Now, ``complexify'' the \(8+8+8\) part of \(\mathrm{F}_4\) to get
\(\mathrm{E}_6\):

\[78 =  28  +  16  +  16  +  16  +   1  +   0  +   1\]

Then, ``quaternionify'' the \(8+8+8\) part of \(\mathrm{F}_4\) to get
\(\mathrm{E}_7\):

\[133 =  28  +  32  +  32  +  32  +   3  +   3  +   3\]

Finally, ``octonionify'' the \(8+8+8\) part of \(\mathrm{F}_4\) to get
\(\mathrm{E}_8\):

\[248 =  28  +  64  +  64  +  64  +   7  +  14  +   7\]

This way shows you that the ``second'' \(\mathrm{Spin}(8)\) in
\(\mathrm{E}_8\) breaks down as \(28 = 7 + 14 + 7\) which is globally
like two 7-spheres and a \(\mathrm{G}_2\), one \(S_7\) for left-action,
one for right-action, and a \(\mathrm{G}_2\) automorphism group of
octonions that is needed to for ``compatibility'' of the two \(S_7\)s.
The \(3+3+3\) of \(\mathrm{E}_7\), the \(1+0+1\) of \(\mathrm{E}_6\),
and the \(0+0+0\) of \(\mathrm{F}_4\) and \(D_4\) are the quaternionic,
complex, and real analogues of the \(7+14+7\)."
\end{quote}

When I asked him where he got this, he said he cooked it up himself
using the construction of \(\mathrm{E}_8\) that I learned from Kostant
together with the Freudenthal--Tits magic square. He gave some references
for the latter:

\begin{enumerate}
\def\labelenumi{\arabic{enumi})}
\setcounter{enumi}{3}
\item
  Hans Freudenthal, ``Lie groups in the foundations of geometry'', 
  \emph{Adv.\ Math.} \textbf{1} (1964) 145--190.
\item
  Jacques Tits, ``Alg\'ebres alternatives, alg\'ebres de Jordan et alg\'ebres de Lie exceptionnelles.  I. Construction'' \emph{Indag.\ Math.\ }\textbf{28} (1966), 223--237.
\item
  Kevin McCrimmon, ``Jordan algebras and their applications'',
  \emph{Bull.\ Amer.\ Math.\ Soc.} \textbf{84} (1978), 612--627, at pp.~620--621.
  Available at
   \href{https://projecteuclid.org/journals/bulletin-of-the-american-mathematical-society/volume-84/issue-4/Jordan-algebras-and-their-applications/bams/1183540925.full}{\texttt{https://projecteuclid.org/}}
\href{https://projecteuclid.org/journals/bulletin-of-the-american-mathematical-society/volume-84/issue-4/Jordan-algebras-and-their-applications/bams/1183540925.full}{\texttt{journals/bulletin-of-the-american-mathematical-society/volume-84/}}
\break
\href{https://projecteuclid.org/journals/bulletin-of-the-american-mathematical-society/volume-84/issue-4/Jordan-algebras-and-their-applications/bams/1183540925.full}{\texttt{issue-4/Jordan-algebras-and-their-applications/bams/1183540925.full}}
\end{enumerate}
\noindent
I would describe it here, but I'm running out of steam, and it's easy to
learn about it from Tony Smith's web page:

\begin{enumerate}
\def\labelenumi{\arabic{enumi})}
\setcounter{enumi}{6}
\tightlist
\item
  Tony Smith, Freudenthal--Tits magic square,
 \href{https://web.archive.org/web/20110808093537/http://valdostamuseum.org/hamsmith/FTsquare.html}{\texttt{https://web.archive.org/web/}} \break
\href{https://web.archive.org/web/20110808093537/http://valdostamuseum.org/hamsmith/FTsquare.html}{\texttt{20110808093537/http://valdostamuseum.org/hamsmith/FTsquare.html}}
\end{enumerate}

\begin{center}\rule{0.5\linewidth}{0.5pt}\end{center}

\begin{quote}
\emph{``I regret that it has been necessary for me in this lecture to
administer such a large dose of four-dimensional geometry. I do not
apologise, because I am not really responsible for the fact that nature
in its most fundamental aspect is four-dimensional.''}

--- Albert North Whitehead.
\end{quote}

\hypertarget{week92}{%
\section{October 17, 1996}\label{week92}}

\hypertarget{week92_tale}{
I'm sure most of you have lost interest in my ``Tale of \(n\)-Categories'',
because it takes a fair amount of work to keep up with all the abstract
concepts involved.} However, we are now at a point where we can have some
fun with what we've got, even if you haven't really followed all the
previous stuff. So what follows is a rambling tour through monads,
adjunctions, the 4-color theorem and the large-N limit of
\(\mathrm{SU}(N)\) gauge theory\ldots.

Okay, so in \protect\hyperlink{week89_tale}{``Week 89''} we defined a gadget
called a ``monad''. Using the string diagrams we talked about, you can
think of a monad as involving a process like this: \[
  \begin{tikzpicture}
    \begin{knot}
      \strand[thick] (0,0.5)
        to (0,0)
        to [out=down,in=up] (0.5,-1)
        to (0.5,-2);
      \strand[thick] (1,0.5)
        to (1,0)
        to [out=down,in=up] (0.5,-1);
    \end{knot}
    \node[fill=white] at (0,0) {$s$};
    \node[fill=white] at (1,0) {$s$};
    \node[label=left:{$M$}] at (0.5,-0.95) {$\bullet$};
    \node[fill=white] at (0.5,-1.5) {$s$};
  \end{tikzpicture}
\] which we read downwards as describing the ``fusion'' of two copies of
something called \(s\) into one copy of the same thing \(s\). The fusion
process itself is called \(M\).

I can hear you wonder, what exactly \emph{is} this thing s? What
\emph{is} this process \(M\)? Well, I gave the technical answer in
\protect\hyperlink{week89_tale}{``Week 89''} --- but the point is that
\(n\)-category theory is deliberately designed to be so general that it
covers pretty much anything you could want! For example, \(s\) could be
the set of real numbers and \(M\) could be multiplication of real
numbers, which is a function from \(s\times s\) to \(s\). Or we could be
doing topology in the plane, in which case the picture above stands for
exactly what it looks like: two lines merging to form one line! These
and many other situations are analogous, and the formalism allows us to
treat them all at once. Here I will not review all the rules of the
game. If you just play along and trust me everything will be all right.
If you don't trust me, go back and check the definitions.

Let me turn to the axioms for a monad. In addition to the multiplication
\(M\) we want to have a ``multiplicative identity'', \(I\), looking like
this: \[
  \begin{tikzpicture}
    \begin{knot}
      \strand[thick] (0,1)
        to (0,0);
    \end{knot}
    \node[label=above:{$I$}] at (0,1) {};
    \node[fill=white] at (0,0.5) {$s$};
  \end{tikzpicture}
\] Here nothing is coming in, and a copy of \(s\) is going out. Because
ordinary multiplication has \(1x = x\) and \(x1 = x\) for all \(x\), we
want the following axioms to hold: \[
  \begin{tikzpicture}
    \begin{knot}
      \strand[thick] (0,0)
        to [out=down,in=up] (0.5,-1)
        to (0.5,-2);
      \strand[thick] (1,1)
        to (1,0)
        to [out=down,in=up] (0.5,-1);
    \end{knot}
    \node[label=above:{$I$}] at (0,0) {};
    \node[fill=white] at (1,0.5) {$s$};
    \node[label=left:{$M$}] at (0.5,-0.95) {$\bullet$};
    \node[fill=white] at (0.5,-1.5) {$s$};
    \node at (2,-0.5) {$=$};
    \begin{knot}
      \strand[thick] (3,1) to (3,-2);
    \end{knot}
    \node[fill=white] at (3,0) {$s$};
  \end{tikzpicture}
\] and \[
  \begin{tikzpicture}
    \begin{scope}[xscale=-1,shift={(-1,0)}]
      \begin{knot}
        \strand[thick] (0,0)
          to [out=down,in=up] (0.5,-1)
          to (0.5,-2);
        \strand[thick] (1,1)
          to (1,0)
          to [out=down,in=up] (0.5,-1);
      \end{knot}
      \node[label=above:{$I$}] at (0,0) {};
      \node[fill=white] at (1,0.5) {$s$};
      \node[label=left:{$M$}] at (0.5,-0.95) {$\bullet$};
      \node[fill=white] at (0.5,-1.5) {$s$};
    \end{scope}
    \node at (2,-0.5) {$=$};
      \begin{knot}
        \strand[thick] (3,1) to (3,-2);
      \end{knot}
    \node[fill=white] at (3,0) {$s$};
  \end{tikzpicture}
\] Also, since ordinary multiplication has \((xy)z = x(yz)\), we want
the following associativity law to hold, too: \[
  \begin{tikzpicture}
    \begin{knot}
      \strand[thick] (0,0.5)
        to (0,0)
        to [out=down,in=up] (0.5,-1)
        to (0.5,-1.5)
        to [out=down,in=up] (1,-2.5)
        to (1,-3.5);
      \strand[thick] (1,0.5)
        to (1,0)
        to [out=down,in=up] (0.5,-1);
      \strand[thick] (2,0.5)
        to (2,0)
        to [out=down,in=up] (1.5,-1)
        to (1.5,-1.5)
        to [out=down,in=up] (1,-2.5);
    \end{knot}
    \node[fill=white] at (0,0) {$s$};
    \node[fill=white] at (1,0) {$s$};
    \node[fill=white] at (2,0) {$s$};
    \node[label=left:{$M$}] at (0.5,-0.95) {$\bullet$};
    \node[label=left:{$M$}] at (1,-2.45) {$\bullet$};
    \node[fill=white] at (0.5,-1.5) {$s$};
    \node[fill=white] at (1,-3) {$s$};
    \node at (3,-1.75) {$=$};
    \begin{scope}[xscale=-1,shift={(-6,0)}]
    \begin{knot}
      \strand[thick] (0,0.5)
        to (0,0)
        to [out=down,in=up] (0.5,-1)
        to (0.5,-1.5)
        to [out=down,in=up] (1,-2.5)
        to (1,-3.5);
      \strand[thick] (1,0.5)
        to (1,0)
        to [out=down,in=up] (0.5,-1);
      \strand[thick] (2,0.5)
        to (2,0)
        to [out=down,in=up] (1.5,-1)
        to (1.5,-1.5)
        to [out=down,in=up] (1,-2.5);
    \end{knot}
    \node[fill=white] at (0,0) {$s$};
    \node[fill=white] at (1,0) {$s$};
    \node[fill=white] at (2,0) {$s$};
    \node[label=left:{$M$}] at (0.5,-0.95) {$\bullet$};
    \node[label=left:{$M$}] at (1,-2.45) {$\bullet$};
    \node[fill=white] at (0.5,-1.5) {$s$};
    \node[fill=white] at (1,-3) {$s$};
    \end{scope}
  \end{tikzpicture}
\] These rules are a translation of the rules given in
\protect\hyperlink{week89_tale}{``Week 89''} into string diagram form.

If you are a physicist, you can think of these diagrams as being funny
Feynman diagrams where you've got some kind of particle \(s\) and two
processes \(M\) and \(I\). Then \(M\) is a bit like what you'd call a
``cubic self-interaction'', where two particles combine to form a third.
These interactions show up in simple textbook theories like the
``\(\varphi^3\) theory'' and, more importantly, in nonabelian gauge
field theories like quantum chromodynamics, where the gauge bosons have
cubic self-interactions. On the other hand, \(I\) is a bit like what you'd
usually call a ``source'' or an ``external potential'', some sort of
field imposed from outside that can create particles of type \(s\). You
shouldn't take the analogy with Feynman diagrams too seriously yet,
because the context we're working in is so general, and the most
interesting physics theories don't correspond to monads but to more
elaborate setups. However, we could flesh out the analogy to make it
very precise and accurate if we wanted, and this is especially important
in topological quantum field theory. More later about that.

Now in \protect\hyperlink{week83_tale}{``Week 83''} I discussed a different
sort of gadget, called an ``adjunction''. Here you have two guys \(x\)
and \(x^*\), and two processes \(U\) and \(C\) called the ``unit'' and
``counit'', which look like this: \[
  \begin{tikzpicture}
    \begin{knot}
      \strand[thick] (0,-0.5)
        to (0,0)
        to [out=up,in=up,looseness=2] (1,0)
        to (1,-0.5);
    \end{knot}
    \node[fill=white] at (0,0) {$x$};
    \node[fill=white] at (1,0) {$x^*$};
    \node[label=above:{$U$}] at (0.5,0.57) {$\bullet$};
  \end{tikzpicture}
  \qquad\raisebox{2em}{\text{and}}\qquad
  \begin{tikzpicture}
    \begin{knot}
      \strand[thick] (0,0.5)
        to (0,0)
        to [out=down,in=down,looseness=2] (1,0)
        to (1,0.5);
    \end{knot}
    \node[fill=white] at (0,0) {$x^*$};
    \node[fill=white] at (1,0) {$x$};
    \node[label=below:{$C$}] at (0.5,-0.6) {$\bullet$};
  \end{tikzpicture}
\] They satisfy the following axioms: \[
  \begin{tikzpicture}
    \begin{knot}
      \strand[thick] (0,0)
      to (0,1)
      to [out=up,in=up,looseness=2] (1,1)
      to [out=down,in=down,looseness=2] (2,1)
      to (2,2);
    \end{knot}
    \node[fill=white] at (0,0.5) {$x$};
    \node[fill=white] at (2,1.5) {$x$};
    \node[fill=white] at (1,1) {$x^*$};
    \node[label=above:{$U$}] at (0.5,1.57) {$\bullet$};
    \node[label=below:{$C$}] at (1.5,0.4) {$\bullet$};
    \node at (3,1) {$=$};
    \begin{scope}[shift={(4,0)}]
      \begin{knot}
        \strand[thick] (0,0) to (0,2);
      \end{knot}
      \node[fill=white] at (0,1.7) {$x$};
    \end{scope}
  \end{tikzpicture}
\] \[
  \begin{tikzpicture}
    \begin{scope}[xscale=-1,shift={(-2,0)}]
      \begin{knot}
        \strand[thick] (0,0)
        to (0,1)
        to [out=up,in=up,looseness=2] (1,1)
        to [out=down,in=down,looseness=2] (2,1)
        to (2,2);
      \end{knot}
      \node[fill=white] at (0,0.5) {$x^*$};
      \node[fill=white] at (2,1.5) {$x^*$};
      \node[fill=white] at (1,1) {$x$};
      \node[label=above:{$U$}] at (0.5,1.57) {$\bullet$};
      \node[label=below:{$C$}] at (1.5,0.4) {$\bullet$};
    \end{scope}
    \node at (3,1) {$=$};
    \begin{scope}[shift={(4,0)}]
      \begin{knot}
        \strand[thick] (0,0) to (0,2);
      \end{knot}
      \node[fill=white] at (0,1.7) {$x^*$};
    \end{scope}
  \end{tikzpicture}
\] Physically, we can think of \(x^*\) as the antiparticle of \(x\), and
then \(U\) is the process of creation of a particle-antiparticle pair,
while \(C\) is the process of annihilation. The axioms just say that for
a particle or antiparticle to ``double back in time'' by means of these
processes isn't really different than for it to march obediently along
forwards. Mathematically, one nice example of an adjunction involves a
vector space \(x\) and its dual vector space \(x^*\). This is really the
same example, since if the behavior of a particle under symmetry
transformations is described by some group representation, its
antiparticle is described by the dual representation. For more details
on the math, see \protect\hyperlink{week83_tale}{``Week 83''}.

Now, let's see how to get a monad from an adjunction! We need to get
\(s\), \(M\), and \(I\) from \(x\), \(x^*\), \(U\), and \(C\). To do
this, we first define \(s\) to be \(xx^*\). Then define \(M\) to be \[
  \begin{tikzpicture}
    \begin{knot}
      \strand[thick] (0,0.5)
        to (0,0)
        to [out=down,in=down,looseness=2] (1,0)
        to (1,0.5);
    \end{knot}
    \node[fill=white] at (0,0) {$x^*$};
    \node[fill=white] at (1,0) {$x$};
    \node[label=below:{$C$}] at (0.5,-0.6) {$\bullet$};
    \begin{knot}
      \strand[thick] (-0.75,0.5)
        to (-0.75,0)
        to [out=down,in=up] (0.125,-1.75)
        to (0.125,-2.5);
      \strand[thick] (1.75,0.5)
        to (1.75,0)
        to [out=down,in=up] (0.875,-1.75)
        to (0.875,-2.5);
    \end{knot}
    \node[fill=white] at (-0.75,0) {$x$};
    \node[fill=white] at (1.75,0) {$x^*$};
    \node[fill=white] at (0,-2) {$x$};
    \node[fill=white] at (1,-2) {$x^*$};
  \end{tikzpicture}
\] Again, to really understand the rules of the game you need to learn a
bit about string diagrams and \(2\)-categories, but the basic idea is
supposed to be simple: we can get two \(xx^*\)'s to turn into one
\(xx^*\) by letting an \(x^*\) and \(x\) annihilate each other!

Finally, we define \(I\) to be \[
  \begin{tikzpicture}
    \begin{knot}
      \strand[thick] (0,-0.5)
        to (0,0)
        to [out=up,in=up,looseness=2] (1,0)
        to (1,-0.5);
    \end{knot}
    \node[fill=white] at (0,0) {$x$};
    \node[fill=white] at (1,0) {$x^*$};
    \node[label=above:{$U$}] at (0.5,0.57) {$\bullet$};
  \end{tikzpicture}
\] In other words, an \(xx^*\) can be created out of nothing since it's
a ``particle/antiparticle pair''.

Now one can check that all the axioms for a monad hold. You really need
to know a bit about \(2\)-categories to do it carefully, but basically
you just let yourself deform the pictures, in part with the help of the
axioms for an adjunction, which let you straighten out curves that
``double back in time.'' So for example, we can prove the identity law
\[
  \begin{tikzpicture}
    \begin{knot}
      \strand[thick] (0,-0.5)
        to (0,0)
        to [out=up,in=up,looseness=2] (1,0)
        to (1,-0.5);
    \end{knot}
    \node[fill=white] at (0,0) {$x$};
    \node[fill=white] at (1,0) {$x^*$};
    \node[label=above:{$U$}] at (0.5,0.57) {$\bullet$};
    \begin{knot}
      \strand[thick] (0,-0.5)
        to [out=down,in=up,looseness=1.5] (1,-3)
        to (1,-3.5);
    \end{knot}
    \node[fill=white] at (1,-3) {$x$};
    \begin{knot}
      \strand[thick] (1,-0.5)
        to [out=down,in=down,looseness=2] (2,-0.5)
        to (2,1.5);
    \end{knot}
    \node[label=below:{$C$}] at (1.5,-1.1) {$\bullet$};
    \node[fill=white] at (2,1) {$x$};
    \begin{knot}
      \strand[thick] (2,-3.5)
        to (2,-3)
        to [out=up,in=down,looseness=1.5] (3,-0.5)
        to (3,1.5);
    \end{knot}
    \node[fill=white] at (2,-3) {$x^*$};
    \node[fill=white] at (3,1) {$x^*$};
    \node at (4,-1) {$=$};
    \begin{knot}
      \strand[thick] (5,1.5) to (5,-3.5);
      \strand[thick] (6,1.5) to (6,-3.5);
    \end{knot}
    \node[fill=white] at (5,1) {$x$};
    \node[fill=white] at (6,1) {$x^*$};
  \end{tikzpicture}
\] by canceling the \(U\) and the \(C\) on the left using one of the
axioms for an adjunction. Similarly, associativity holds because the
following two pictures are topologically the same: \[
  \begin{tikzpicture}
    \begin{knot}
      \strand[thick] (0,0.5)
        to (0,0)
        to [out=down,in=down,looseness=2] (1,0)
        to (1,0.5);
    \end{knot}
    \node[fill=white] at (0,0) {$x^*$};
    \node[fill=white] at (1,0) {$x$};
    \node[label=below:{$C$}] at (0.5,-0.6) {$\bullet$};
    \begin{knot}
      \strand[thick] (-0.75,0.5)
        to (-0.75,0)
        to [out=down,in=up] (0.125,-1.75)
        to (0.125,-2.5);
      \strand[thick] (1.75,0.5)
        to (1.75,0)
        to [out=down,in=up] (0.875,-1.75)
        to (0.875,-2.5);
    \end{knot}
    \node[fill=white] at (-0.75,0) {$x$};
    \node[fill=white] at (1.75,0) {$x^*$};
    \node[fill=white] at (0,-2) {$x$};
    \node[fill=white] at (1,-2) {$x^*$};
    \begin{scope}[shift={(0.875,-3)}]
      \begin{knot}
        \strand[thick] (0,0.5)
          to (0,0)
          to [out=down,in=down,looseness=2] (1,0)
          to (1,0.5);
      \end{knot}
      \node[fill=white] at (0,0) {$x^*$};
      \node[fill=white] at (1,0) {$x$};
      \node[label=below:{$C$}] at (0.5,-0.6) {$\bullet$};
      \begin{knot}
        \strand[thick] (-0.75,0.5)
          to (-0.75,0)
          to [out=down,in=up] (0.125,-1.75)
          to (0.125,-2.5);
        \strand[thick] (1.75,0.5)
          to (1.75,0)
          to [out=down,in=up] (0.875,-1.75)
          to (0.875,-2.5);
      \end{knot}
      \node[fill=white] at (-0.75,0) {$x$};
      \node[fill=white] at (1.75,0) {$x^*$};
      \node[fill=white] at (0,-2) {$x$};
      \node[fill=white] at (1,-2) {$x^*$};
    \end{scope}
    \begin{scope}[shift={(1.875,-2.5)}]
      \begin{knot}
        \strand[thick] (0,0)
          to (0,0.5)
          to [out=up,in=down,looseness=0.75] (1,2.5)
          to (1,3);
        \strand[thick] (0.75,0)
          to (0.75,0.5)
          to [out=up,in=down,looseness=0.75] (1.75,2.5)
          to (1.75,3);
      \end{knot}
      \node[fill=white] at (0,0.5) {$x$};
      \node[fill=white] at (0.75,0.5) {$x^*$};
      \node[fill=white] at (1,2.5) {$x$};
      \node[fill=white] at (1.75,2.5) {$x^*$};
    \end{scope}
    \node at (4.5,-2.5) {$=$};
    \begin{scope}[xscale=-1,shift={(-9,0)}]
      \begin{knot}
        \strand[thick] (0,0.5)
          to (0,0)
          to [out=down,in=down,looseness=2] (1,0)
          to (1,0.5);
      \end{knot}
      \node[fill=white] at (0,0) {$x$};
      \node[fill=white] at (1,0) {$x^*$};
      \node[label=below:{$C$}] at (0.5,-0.6) {$\bullet$};
      \begin{knot}
        \strand[thick] (-0.75,0.5)
          to (-0.75,0)
          to [out=down,in=up] (0.125,-1.75)
          to (0.125,-2.5);
        \strand[thick] (1.75,0.5)
          to (1.75,0)
          to [out=down,in=up] (0.875,-1.75)
          to (0.875,-2.5);
      \end{knot}
      \node[fill=white] at (-0.75,0) {$x^*$};
      \node[fill=white] at (1.75,0) {$x$};
      \node[fill=white] at (0,-2) {$x^*$};
      \node[fill=white] at (1,-2) {$x$};
      \begin{scope}[shift={(0.875,-3)}]
        \begin{knot}
          \strand[thick] (0,0.5)
            to (0,0)
            to [out=down,in=down,looseness=2] (1,0)
            to (1,0.5);
        \end{knot}
        \node[fill=white] at (0,0) {$x$};
        \node[fill=white] at (1,0) {$x^*$};
        \node[label=below:{$C$}] at (0.5,-0.6) {$\bullet$};
        \begin{knot}
          \strand[thick] (-0.75,0.5)
            to (-0.75,0)
            to [out=down,in=up] (0.125,-1.75)
            to (0.125,-2.5);
          \strand[thick] (1.75,0.5)
            to (1.75,0)
            to [out=down,in=up] (0.875,-1.75)
            to (0.875,-2.5);
        \end{knot}
        \node[fill=white] at (-0.75,0) {$x^*$};
        \node[fill=white] at (1.75,0) {$x$};
        \node[fill=white] at (0,-2) {$x^*$};
        \node[fill=white] at (1,-2) {$x$};
      \end{scope}
      \begin{scope}[shift={(1.875,-2.5)}]
        \begin{knot}
          \strand[thick] (0,0)
            to (0,0.5)
            to [out=up,in=down,looseness=0.75] (1,2.5)
            to (1,3);
          \strand[thick] (0.75,0)
            to (0.75,0.5)
            to [out=up,in=down,looseness=0.75] (1.75,2.5)
            to (1.75,3);
        \end{knot}
        \node[fill=white] at (0,0.5) {$x^*$};
        \node[fill=white] at (0.75,0.5) {$x$};
        \node[fill=white] at (1,2.5) {$x^*$};
        \node[fill=white] at (1.75,2.5) {$x$};
      \end{scope}
    \end{scope}
  \end{tikzpicture}
\] Whew! Drawing these is tough work.

Now, as I said, an example of an adjunction is a vector space \(x\) and
its dual \(x^*\). What monad do we get in this case? Well, the vector
space \(x\) tensored with \(x^*\) is just the vector space of linear
transformations of \(x\), so that's our monad in this case. In less
high-brow terms, we've proven that matrices form an algebra when you
define matrix multiplication in the usual way! In particular, the above
picture serves as a diagrammatic proof that matrix multiplication is
associative.

Of course, people didn't invent all this fancy-looking (but actually
very basic) stuff just to deal with matrix multiplication! Or did they?
Well, actually, Penrose \emph{did} invent a diagrammatic notation for
tensors which is just a slightly souped-up version of the above stuff.
You can find it in:

\begin{enumerate}
\def\labelenumi{\arabic{enumi})}
\tightlist
\item
  ``Applications of negative dimensional tensors'', by Roger Penrose, in
  \emph{Combinatorial Mathematics and its Applications}, ed.~D. J. A.
  Welsh, Academic Press, 1971.
\end{enumerate}

But most of the work on this sort of thing has been aimed at
applications of other sorts.

Now let me drift over to a related subject, the large-\(N\) limit of
\(\mathrm{SU}(N)\) gauge theory. Quantum chromodynamics, or QCD, is an
\(\mathrm{SU}(N)\) gauge theory with \(N = 3\), but it turns out that
things simplify a lot in the limit as \(N\to\infty\), and one gets some
nice qualitative insight into the strong force by considering this
simplified theory. One can even treat the number \(3\) as a small
perturbation around the number \(\infty\) and get some decent answers! A
good introduction to this appears in Coleman's delightful book,
essential reading for anyone learning particle physics:

\begin{enumerate}
\def\labelenumi{\arabic{enumi})}
\setcounter{enumi}{1}
\tightlist
\item
  Sidney Coleman, \emph{Aspects of Symmetry}, Cambridge U.\
  Press, Cambridge, 1989.
\end{enumerate}
\noindent
Check out section 8.3.1, entitled ``the double line representation and
the dominance of planar graphs''. Coleman considers Yang--Mills theories,
like QCD, but many of the same ideas apply to other gauge theories.

The idea is that if we start out studying the Feynman diagrams for a
gauge field theory with gauge group \(\mathrm{SU}(N)\), and see how much
various diagrams contribute to any process for large \(N\), the diagrams
that contribute the most are those that can be drawn on a plane without
any lines crossing. Technically, the reason is that diagrams that can
only be drawn on a surface of genus \(g\) grow like \(N^{2-2g}\) as
\(N\) increases. This number \(2-2g\) is called the Euler characteristic
and it's biggest when your surface has no handles.

Even better, in the \(N\to\infty\) limit we can think of the Feynman
diagrams using diagrams like the ones above. For example, we can think
of the cubic self-interaction in Yang--Mills theory as simply matrix
multiplication: \[
  \begin{tikzpicture}
    \draw[thick] (1,0) to node[fill=white]{$x^*$} (1.5,-1) node[label=below:{$C$}]{$\bullet$} to node[fill=white]{$x$} (2,0);
    \draw[thick] (0,0) to node[fill=white]{$x$} (1,-2) to node[fill=white]{$x$} (1,-3);
    \draw[thick] (3,0) to node[fill=white]{$x^*$} (2,-2) to node[fill=white]{$x^*$} (2,-3);
  \end{tikzpicture}
\] and the quartic self-interaction as something a wee bit fancier: \[
  \begin{tikzpicture}
    \draw[thick] (1,0) to node[fill=white]{$x^*$} (1.5,-1) node[label=below:{$C$}]{$\bullet$} to node[fill=white]{$x$} (2,0);
    \draw[thick] (0,0) to node[fill=white]{$x$} (1,-2) to node[fill=white]{$x$} (0,-4);
    \draw[thick] (3,0) to node[fill=white]{$x^*$} (2,-2) to node[fill=white]{$x^*$} (3,-4);
    \draw[thick] (1,-4) to node[fill=white]{$x^*$} (1.5,-3) node[label=above:{$U$}]{$\bullet$} to node[fill=white]{$x$} (2,-4);
  \end{tikzpicture}
\] Apparently these ideas have spawned a whole field of physics called
``matrix models''.

These ideas work not only for Yang--Mills theory but also for
Chern--Simons theory, which is a topological quantum field theory: a
theory that doesn't require any metric on spacetime to make sense. Here
they have been exploited by Dror Bar-Natan to come up with a new
formulation of the famous 4-color theorem:

\begin{enumerate}
\def\labelenumi{\arabic{enumi})}
\setcounter{enumi}{2}
\tightlist
\item
  Dror Bar-Natan, ``Lie algebras and the four color theorem'', 
  available as
  \href{https://arxiv.org/ps/q-alg/9606016}{\texttt{q-alg/9606016}}.
\end{enumerate}

As I explained in \protect\hyperlink{week8}{``Week 8''} and
\protect\hyperlink{week22}{``Week 22''}, there is a way to formulate
the 4-color theorem as a statement about trivalent graphs. In
particular, Penrose invented a little recipe that lets us calculate an
invariant of trivalent graphs, which is zero for some \emph{planar}
graph only if some corresponding map can't be 4-colored. This recipe
involves the vector cross product, or equivalently, the Lie algebra of
the group \(\mathrm{SU}(2)\). You can generalize it to work for
\(\mathrm{SU}(N)\). And if you then consider the \(N\to\infty\) limit,
you get the above stuff! (The point is that the above stuff also gives a
rule for computing a number from any trivalent graph.)

Now as I said, in the \(N\to\infty\) limit all the nonplanar Feynman
diagrams give negligible results compared to the planar ones. So another
way to state the 4-color theorem is this: if the \(\mathrm{SU}(2)\)
invariant of a trivalent graph is zero, the \(\mathrm{SU}(N)\) invariant
is negligible in the \(N\to\infty\) limit.

This doesn't yet give a new proof of the 4-color theorem. But it makes
it into sort of a \emph{physics} problem: a problem about the relation
of \(\mathrm{SU}(2)\) Chern--Simons theory and the \(N\to\infty\) limit
of Chern--Simons theory.

Now, the 4-color theorem is one of the two deep mysteries of
2-dimensional topology --- a subject too often considered trivial. The
other mystery is the Andrews--Curtis conjecture, discussed in
\protect\hyperlink{week23}{``Week 23''}. Often a problem is hard or
unsolvable until you get the right tools. Topological quantum field
theory is a new tool in topology, so one could hope it'll shed some
light on these problems. Bar-Natan's paper is a tantalizing piece of
evidence that maybe, just maybe, it will.

One can't really tell yet.

Anyway, I don't really care much about the 4-color theorem per se. If I
ever need to color a map I'll hire a cartographer. It's the connections
between seemingly disparate subjects that I find interesting.
\(2\)-categories are a very abstract formalism developed to describe
\(2\)-dimensional ways of glomming things together. Starting from the
study of \(2\)-categories, we very naturally get the notions of
``monad'' and ``adjunction''. And before we know it, this leads us to
some interesting questions about \(2\)-dimensional quantum field theory:
for really, the dominance of planar diagrams in the \(N\to\infty\) limit
of gauge theory is saying that in this limit the theory becomes
essentially a 2-dimensional field theory, in some funny sense. And then,
lo and behold, this turns out to be related to the 4-color theorem!

By the way, I guess you all know that the 4-color theorem was proved
using a computer, by breaking things down into lots of separate cases.
(See \protect\hyperlink{week22}{``Week 22''} for references.) Well,
there's a new proof out, which also uses a computer, but is supposed to
be simpler:

\begin{enumerate}
\def\labelenumi{\arabic{enumi})}
\setcounter{enumi}{3}
\tightlist
\item
  Neil Robertson, Daniel P. Sanders, Paul Seymour, and Robin Thomas, ``A
  new proof of the four-colour theorem'', \emph{Electronic Research
  Announcements of the American Mathematical Society} \textbf{2} (1996),
  17--25. Available at \href{http://www.ams.org/journals/era/1996-02-01/}
  {\texttt{http://www.ams.org/journals/}} \href{http://www.ams.org/journals/era/1996-02-01/}{\texttt{journals/era/1996-02-01/}}.
\end{enumerate}
\noindent
I'm still hoping for the 2-page ``physicist's proof'' using path
integrals!

To continue reading the ``Tale of \(n\)-Categories'', see
\protect\hyperlink{week99_tale}{``Week 99''}.  For more on adjunctions and monoid objects, try ``Week173'' and especially ``Week 174''.

\hypertarget{week93}{%
\section{October 27, 1996}\label{week93}}

Lately I've been trying to learn more about string theory. I've always
had grave doubts about string theory, but it seems worth knowing about.
As usual, when I'm trying to learn something I find it helpful to write
about it --- it helps me remember stuff, and it points out gaps in my
understanding. So I'll start trying to explain some string theory in
this and forthcoming Week's Finds.

However: watch out! This isn't going to be a systematic introduction to
the subject. First of all, I don't know enough to do that. Secondly, it
will be very quirky and idiosyncratic, because the aspects of string
theory I'm interested in now aren't necessarily the ones most string
theorists would consider central. I've been taking as my theme of
departure, ``What's so great about 10 and 26 dimensions?'' When one
reads about string theory, one often hears that it only works in 10 or
26 dimensions --- and the obvious question is, why?

This question leads one down strange roads, and one runs into lots of
surprising coincidences, and spooky things that sound like coindences
but might NOT be coincidences if we understood them better.

For example, when we have a string in 26 dimensions we can think of it
as wiggling around in the 24 directions perpendicular to the
2-dimensional surface the string traces out in spacetime (the ``string
worldsheet''). So the number 24 plays an especially important role in
26-dimensional string theory. It turns out that
\[1^2 + 2^2 + 3^2 + \cdots + 24^2 = 70^2.\] In fact, 24 is the
\emph{only} integer \(n > 1\) such that the sum of squares from \(1^2\)
to \(n^2\) is itself a perfect square. Is this a coincidence? Probably
not, as I'll eventually explain! This is just one of many eerie facts
one meets when trying to understand this stuff.

For starters I just want to explain why dimensions of the form
\(8k + 2\) are special. Notice that if we take \(k = 0\) here we get
\(2\), the dimension of the string worldsheet. For \(k = 1\) we get
\(10\), the dimension of spacetime in ``supersymmetric string theory''.
For \(k = 3\) we get \(26\), the dimension of spacetime in ``purely
bosonic string theory''. So these dimensions are important. What about
\(k = 2\) and the dimension \(18\), I hear you ask? Well, I don't know
what happens there yet... maybe someone can tell me! All I want to
do now is to explain what's good about \(8k+2\).

But I need to start by saying a bit about fermions.

Remember that in the Standard Model of particle physics --- the model
that all fancier theories are trying to outdo --- elementary particles
come in 3 basic kinds. There are the basic fermions. In general a
``fermion'' is a particle whose angular momentum comes in units of
Planck's constant \(\hbar\) times \(1/2\), \(3/2\), \(5/2\), and so on.
Fermions satisfy the Pauli exclusion principle --- you can't put two
identical fermions in the same state. That's why we have chemistry: the
electrons stack up in ``shells'' at different energy levels, instead of
all going to the lowest-energy state, because they are fermions and
satisfy the exclusion principle. In the Standard Model the fermions go
like this:

\begin{longtable}[]{@{}llll@{}}
\toprule
\textbf{Leptons} & & \textbf{Quarks} &\tabularnewline
\midrule
\endhead
electron & electron neutrino & down quark & up quark\tabularnewline
muon & muon neutrino & strange quark & charm quark\tabularnewline
tauon & tauon neutrino & bottom quark & top quark\tabularnewline
\bottomrule
\end{longtable}

There are three ``generations'' here, all rather similar to each other.

There are also particles in the Standard Model called ``bosons'' having
angular momentum in units of \(\hbar\) times \(0\),\(1\),\(2\), and so
on. Identical bosons, far from satisfying the exclusion principle, sort
of like to all get into the same state: one sees this in phenomena such
as lasers, where lots of photons occupy the same few states. Most of the
bosons in the Standard Model are called ``gauge bosons''. These carry
the different forces in the standard model, by which the particles
interact:

\begin{longtable}[]{@{}lll@{}}
\toprule
Electromagnetic force & Weak force & Strong force\tabularnewline
\midrule
\endhead
photon & W${}_+$, W${}_-$, Z & 8
gluons\tabularnewline
\bottomrule
\end{longtable}

Finally, there is also a bizarre particle in the Standard Model called
the ``Higgs boson''. This was first introduced as a rather ad hoc
hypothesis: it's supposed to interact with the forces in such a way as
to break the symmetry that would otherwise be present between the
electromagnetic force and the weak force. It has not yet been observed;
finding it would would represent a great triumph for the Standard Model,
while \emph{not} finding it might point the way to better theories.

Indeed, while the Standard Model has passed many stringent experimental
tests, and successfully predicted the existence of many particles which
were later observed (like the W, the Z, and the charm and top quarks),
it is a most puzzling sort of hodgepodge. Could nature really be this
baroque at its most fundamental level? Few people seem to think so; most
hope for some deeper, simpler theory.

It's easy to want a ``deeper, simpler theory'', but how to get it? What
are the clues? What can we do? Experimentalists certainly have their
work cut out for them. They can try to find or rule out the Higgs. They
can also try to see if neutrinos, assumed to be massless in the Standard
Model, actually have a small mass --- for while the Standard Model could
easily be patched if this were the case, it would shed interesting light
on one of the biggest mysteries in physics, namely why the fermions in
nature seem not to be symmetric under reflection, or ``parity''. Right
now, we believe that neutrinos only exist in a left-handed form,
rotating one way but not the other around the axis they move along. This
is intimately related to their apparent masslessness. In fact, for
reasons that would take a while to explain, the lack of parity symmetry
in the Standard Model forces us to assume all the observed fermions
acquire their mass only through interaction with the Higgs particle! For
more on the neutrino mass puzzle, try:

\begin{enumerate}
\def\labelenumi{\arabic{enumi})}
\tightlist
\item
  Paul Langacker, Implications of Neutrino Mass,
  \href{https://www.physics.upenn.edu/~pgl/neutrino/jhu/jhu.html}{\texttt{https://www.physics.upenn.edu/}}  \href{https://www.physics.upenn.edu/~pgl/neutrino/jhu/jhu.html}{\texttt{\(\sim\)pgl/neutrino/jhu/jhu.html}}
\end{enumerate}
\noindent
And, of course, experimentalists can continue to do what they always do
best: discover the utterly unexpected.

Theorists, on the other hand, have been spending the last couple of
decades poring over the standard model and trying to understand what
it's telling us. It's so full of suggestive patterns and partial
symmetries! First, why are there 3 forces here? Each force goes along
with a group of symmetries called a ``gauge group'', and
electromagnetism corresponds to \(\mathrm{U}(1)\), while the weak force
corresponds to \(\mathrm{SU}(2)\) and the strong force corresponds to
\(\mathrm{SU}(3)\). (Here \(\mathrm{U}(n)\) is the group of
\(n\times n\) unitary complex matrices, while \(\mathrm{SU}(n)\) is the
subgroup consisting of those with determinant equal to \(1\).) Well,
actually the Standard Model partially unifies the electromagnetic and
weak force into the ``electroweak force'', and then resorts to the Higgs
to explain why these forces are so different in practice. Various
``grand unified theories'' or ``GUTs'' try to unify the forces further
by sticking the group
\(\mathrm{SU}(3)\times\mathrm{SU}(2)\times\mathrm{U}(1)\) into a bigger
group --- but then resort to still more Higgses to break the symmetry
between them!

Then, there is the curious parallel between the leptons and quarks in
each generation. Each generation has a lepton with mass, a massless or
almost massless neutrino, and two quarks. The massive lepton has charge
\(-1\), the neutrino has charge \(0\) as its name suggests, the ``down''
type quark has charge \(-1/3\), and the ``up'' type quark has charge
\(2/3\). Funny pattern, eh? The Standard Model does not really explain
this, although it would be ruined by ``anomalies'' --- certain
nightmarish problems that can beset a quantum field theory --- if one
idly tried to mess with the generations by leaving out a quark or the
like. It's natural to try to ``unify'' the quarks and leptons, and
indeed, in grand unified theories like the \(\mathrm{SU}(5)\) theory
proposed in 1974 of Georgi and Glashow, the quarks and leptons are
treated in a unified way.

Another interesting pattern is the repetition of generations itself. Why
is there more than one? Why are there three, almost the same, but with
the masses increasing dramatically as we go up? The Standard Model makes
no attempt to explain this, although it does suggest that there had
better not be more than 17 quarks --- more, and the strong force would
not be ``asymptotically free'' (weak at high energies), which would
cause lots of problems for the theory. In fact, experiments strongly
suggest that there are no more than 3 generations. Why?

Finally, there is the grand distinction between bosons and fermions.
What does this mean? Here we understand quite a bit from basic
principles. For example, the ``spin-statistics theorem'' explains why
particles with half-integer spin should satisfy the Pauli exclusion
principle, while those with integer spin should like to hang out
together. This is a very beautiful result with a deep connection to
topology, which I try to explain here:

\begin{enumerate}
\def\labelenumi{\arabic{enumi})}
\setcounter{enumi}{1}
\tightlist
\item
  John Baez, Spin, statistics, CPT and all that jazz,
  \href{http://math.ucr.edu/home/baez/spin.stat.html}{\texttt{http://math.ucr.edu/home/}}  \href{http://math.ucr.edu/home/baez/spin.stat.html}{\texttt{baez/spin.stat.html}}
\end{enumerate}
\noindent
But many people have tried to bridge the chasm between bosons and
fermions, unifying them by a principle called ``supersymmetry''. As in
the other cases mentioned above, when they do this, they then need to
pull tricks to ``break'' the symmetry to get a theory that fits the
experimental fact that bosons and fermions are very different.
Personally, I'm suspicious of all these symmetries postulated only to be
cleverly broken; this approach was so successful in dealing with the
electroweak force --- modulo the missing Higgs! --- that it seems to have
been accepted as a universal method of having ones cake and eating it
too.

Now, string theory comes in two basic flavors. Purely bosonic string
theory lives in 26 dimensions and doesn't have any fermions in it.
Supersymmetric string theories live in 10 dimensions and have both
bosons and fermions, unified via supersymmetry. To deal with the
fermions in nature, most work in physics has focused on the
supersymmetric case. Just for completeness, I should point out that
there are 5 different supersymmetric string theories: type I, type IIA,
type IIB, \(\mathrm{E}_8\times\mathrm{E}_8\) heterotic and
\(\mathrm{SO}(32)\) heterotic. For more on these, see
\protect\hyperlink{week72}{``Week 72''}. We won't be getting into them
here. Instead, I just want to explain how fermions work in different
dimensions, and why nice things happen in dimensions of the form
\(8k + 2\). Most of what I say is in Section 3 of

\begin{enumerate}
\def\labelenumi{\arabic{enumi})}
\setcounter{enumi}{2}
\tightlist
\item
  John H. Schwarz, ``Introduction to supersymmetry'', in
  \emph{Superstrings and Supergravity, Proc. of the 28th Scottish
  Universities Summer School in Physics}, ed.~A. T. Davies and D. G.
  Sutherland, University Printing House, Oxford, 1985.
\end{enumerate}
\noindent
but mathematicians may also want to supplement this with material from
the book \emph{Spin Geometry} by Lawson and Michelson, cited in
\protect\hyperlink{week82}{``Week 82''}.

To understand fermions in different dimensions we need to understand
Clifford algebras. As far as I know, when Clifford originally invented
these algebras in the late 1800s, he was trying to generalize Hamilton's
quaternion algebra by considering algebras that had lots of different
anticommuting square roots of \(-1\). In other words, he considered an
associative algebra generated by a bunch of guys \(e_1,\ldots,e_n\),
satisfying \[e_i^2 = -1\] for all \(i\), and \[e_i e_j = - e_j e_i\]
whenever \(i\) is not equal to \(j\). I discussed these algebras in
\protect\hyperlink{week82}{``Week 82''} and I said what they all were
--- they all have nice descriptions in terms of the reals, the
complexes, and the quaternions.

These original Clifford algebras are great for studying rotations in
\(n\)-dimensional Euclidean space --- please take my word for this for
now. However, here we want to study rotations and Lorentz
transformations in \(n\)-dimensional Minkowski spacetime, so we need to
work with a slightly Different kind of Clifford algebra, which was
probably invented by Dirac. In \(n\)-dimensional Euclidean space the
metric (used for measuring distances) is
\[dx_1^2 + dx_2^2 + \cdots + dx_n^2\] while in \(n\)-dimensional
Minkowski spacetime it is \[dx_1^2 + dx_2^2 + \cdots - dx_n^2\] or if
you prefer (it's just a matter of convention), you can take it to be
\[-dx_1^2 - dx_2^2 - \cdots + dx_n^2\] So it turns out that we need to
switch some signs in the definition of the Clifford algebra when working
in Minkowski spacetime.

In general, we can define the Clifford algebra \(C_{p,q}\) to be the
algebra generated by a bunch of elements \(e_i\), with \(p\) of them
being square roots of \(-1\) and \(q\) of them being square roots of
\(1\). As before, we require that they anticommute:
\[e_i e_j = - e_j e_i\] when \(i\) and \(j\) are different. Physicists
usually call these guys ``gamma matrices''. For \(n\)-dimensional
Minkowski space we can work either with \(C_{n-1,1}\) or \(C_{1,n-1}\),
depending on our preference. As Cecile DeWitt has pointed out, it
\emph{does} make a difference which one we use.

With some work, one can check that these algebras go like this:

\begin{longtable}[]{@{}rlrl@{}}
\toprule
\endhead
\(C_{0,1}\) & \(\mathbb{R}+\mathbb{R}\) & \(C_{1,0}\) &
\(\mathbb{C}\)\tabularnewline
\(C_{1,1}\) & \(\mathbb{R}(2)\) & \(C_{1,1}\) &
\(\mathbb{R}(2)\)\tabularnewline
\(C_{2,1}\) & \(\mathbb{C}(2)\) & \(C_{1,2}\) &
\(\mathbb{R}(2)+\mathbb{R}(2)\)\tabularnewline
\(C_{3,1}\) & \(\mathbb{H}(2)\) & \(C_{1,3}\) &
\(\mathbb{R}(4)\)\tabularnewline
\(C_{4,1}\) & \(\mathbb{H}(2)+\mathbb{H}(2)\) & \(C_{1,4}\) &
\(\mathbb{C}(4)\)\tabularnewline
\(C_{5,1}\) & \(\mathbb{H}(4)\) & \(C_{1,5}\) &
\(\mathbb{H}(4)\)\tabularnewline
\(C_{6,1}\) & \(\mathbb{C}(8)\) & \(C_{1,6}\) &
\(\mathbb{H}(4)+\mathbb{H}(4)\)\tabularnewline
\(C_{7,1}\) & \(\mathbb{R}(16)\) & \(C_{1,7}\) &
\(\mathbb{H}(8)\)\tabularnewline
\bottomrule
\end{longtable}

I've only listed these up to \(8\)-dimensional Minkowski spacetime, and
the cool thing is that after that they sort of repeat --- more
precisely, \(C_{n+8,1}\) is just the same as \(16\times16\) matrices
with entries in \(C_{n,1}\), and \(C_{1,n+8}\) is just \(16\times16\)
matrices with entries in \(C_{1,n}\)! This ``period-8'' phenomenon,
sometimes called Bott periodicity, has implications for all sorts of
branches of math and physics. This is why fermions in 2 dimensions are a
bit like fermions in 10 dimensions and 18 dimensions and 26
dimensions\ldots.

In physics, we describe fermions using ``spinors'', but there are
different kinds of spinors: Dirac spinors, Weyl spinors, Majorana
spinors, and even Majorana--Weyl spinors. This is a bit technical but I
want to dig into it here, since it explains what's special about
\(8k + 2\) dimensions and especially 10 dimensions.

Before I get technical, though, let me just summarize the point for
those of you who don't want all the gory details. ``Dirac spinors'' are
what you use to describe spin-\(1/2\) particles that come in both
left-handed and right-handed forms and aren't their own antiparticle ---
like the electron. Weyl spinors have half as many components, and
describe spin-\(1/2\) particles with an intrinsic handedness that aren't
their own antiparticle --- like the neutrino. ``Weyl spinors'' are only
possible in even dimensions!

Both these sorts of spinors are ``complex'' --- they have complex-valued
components. But there are also real spinors. These are used for
describing particles that are their own antiparticle, because the
operation of turning a particle into an antiparticle is described
mathematically by complex conjugation. ``Majorana spinors'' describe
spin-\(1/2\) particles that come in both left-handed and right-handed
forms and are their own antiparticle. Finally, ``Majorana--Weyl spinors''
are used to describe spin-\(1/2\) particles with an intrinsic handedness
that are their own antiparticle.

As far as we can tell, none of the particles we've seen are Majorana or
Majorana--Weyl spinors, although if the neutrino has a mass it might be a
Majorana spinor. Majorana and Majorana--Weyl spinors only exist in
certain dimensions. In particular, Majorana--Weyl spinors are very
finicky: they only work in dimensions of the form \(8k + 2\). This is
part of what makes supersymmetric string theory work in 10 dimensions!

Now let me describe the technical details. I'm doing this mainly for my
own benefit; if I write this up, I'll be able to refer to it whenever I
forget it. For those of you who stick with me, there will be a little
reward: we'll see that a certain kind of supersymmetric gauge theory, in
which there's a symmetry between gauge bosons and fermions, only works
in dimensions 3, 4, 6, and 10. Perhaps coincidentally --- I don't
understand this stuff well enough to know --- these are also the
dimensions when supersymmetric string theory works classically. (It's
the quantum version that only works in dimension 10.)

So: part of the point of these Clifford algebras is that they give
representations of the double cover of the Lorentz group in different
dimensions. In \protect\hyperlink{week61}{``Week 61''} I explained this
double cover business, and how the group \(\mathrm{SO}(n)\) of rotations
of \(n\)-dimensional Euclidean space has a double cover called
\(\mathrm{Spin}(n)\). Similarly, the Lorentz group of \(n\)-dimensional
Minkowski space, written \(\mathrm{SO}(n-1,1)\), has a double cover we
could call \(\mathrm{Spin}(n-1,1)\). The spinors we'll discuss are all
representations of this group.

The way Clifford algebras help is that there is a nice way to embed
\(\mathrm{Spin}(n-1,1)\) in either \(C_{n-1,1}\) or \(C_{1,n-1}\), so
any representation of these Clifford algebras gives a representation of
\(\mathrm{Spin}(n-1,1)\). We have a choice of dealing with real
representations or complex representations. Any complex representation
of one of these Clifford algebras is also a representation of the
\emph{complexified} Clifford algebra. What I mean is this: above I
implicitly wanted \(C_{p,q}\) to consist of all \emph{real} linear
combinations of products of the \(e_i\), but we could have worked with
\emph{complex} linear combinations instead. Then we would have
``complexified'' \(C_{p,q}\). Since the complex numbers include a square
root of minus 1, the complexification of \(C_{p,q}\) only depends on the
dimension \(p + q\), not on how many minus signs we have.

Now, it is easy and fun and important to check that if you complexify
\(\mathbb{R}\) you get \(\mathbb{C}\), and if you complexify
\(\mathbb{C}\) you get \(\mathbb{C}+\mathbb{C}\), and if you complexify
\(\mathbb{H}\) you get \(\mathbb{C}(2)\). Thus from the above table we
get this table:

\begin{longtable}[]{@{}ll@{}}
\toprule
dimension \(n\) & complexified Clifford algebra\tabularnewline
\midrule
\endhead
1 & \(\mathbb{C}+\mathbb{C}\)\tabularnewline
2 & \(\mathbb{C}(2)\)\tabularnewline
3 & \(\mathbb{C}(2)+\mathbb{C}(2)\)\tabularnewline
4 & \(\mathbb{C}(4)\)\tabularnewline
5 & \(\mathbb{C}(4)+\mathbb{C}(4)\)\tabularnewline
6 & \(\mathbb{C}(8)\)\tabularnewline
7 & \(\mathbb{C}(8)+\mathbb{C}(8)\)\tabularnewline
8 & \(\mathbb{C}(16)\)\tabularnewline
\bottomrule
\end{longtable}
\noindent
Notice this table is a lot simpler --- complex Clifford algebras are
``period-2'' instead of period-8.

Now the smallest complex representation of the complexified Clifford
algebra in dimension \(n\) is what we call a ``Dirac spinor''. We can figure
out what this is using the above table, since the smallest complex
representation of \(\mathbb{C}(n)\) or \(\mathbb{C}(n) + \mathbb{C}(n)\)
is on the \(n\)-dimensional complex vector space \(\mathbb{C}^n\), given
by matrix multiplication. Of course, for
\(\mathbb{C}(n) + \mathbb{C}(n)\) there are \emph{two} representations
depending on which copy of \(\mathbb{C}(n)\) we use, but these give
equivalent representations of \(\mathrm{Spin}(n-1,1)\), which is what
we're really interested in, so we still speak of ``the'' Dirac spinors.

So we get:

\begin{longtable}[]{@{}ll@{}}
\toprule
dimension \(n\) & Dirac spinors\tabularnewline
\midrule
\endhead
1 & \(\mathbb{C}\)\tabularnewline
2 & \(\mathbb{C}(2)\)\tabularnewline
3 & \(\mathbb{C}(2)\)\tabularnewline
4 & \(\mathbb{C}(4)\)\tabularnewline
5 & \(\mathbb{C}(4)\)\tabularnewline
6 & \(\mathbb{C}(8)\)\tabularnewline
7 & \(\mathbb{C}(8)\)\tabularnewline
8 & \(\mathbb{C}(16)\)\tabularnewline
\bottomrule
\end{longtable}

The dimension of the Dirac spinors doubles as we go to each new even
dimension.

We can also look for the smallest real representation of \(C_{n-1,1}\)
or \(C_{1,n-1}\). This is easy to work out from our tables using the
fact that the algebra \(\mathbb{R}\) has its smallest real
representation on \(\mathbb{R}\), while for \(\mathbb{C}\) it's on
\(\mathbb{R}^2\) and for \(\mathbb{H}\) it's on \(\mathbb{R}^4\).

Sometimes this smallest real representation is secretly just the Dirac
spinors \emph{viewed as a real representation} --- we can view
\(\mathbb{C}^n\) as the real vector space \(\mathbb{R}^{2n}\). But
sometimes the Dirac spinors are the \emph{complexification} of the
smallest real representation --- for example, \(\mathbb{C}^n\) is the
complexification of \(\mathbb{R}^n\). In this case folks call the
smallest real representation ``Majorana spinors''.

When we are looking for the smallest real representations, we get
different answers for \(C_{n-1,1}\) and \(C_{1,n-1}\). Here is what we
get:

\begin{longtable}[]{@{}lllcllc@{}}
\toprule
\begin{minipage}[b]{0.04\columnwidth}\raggedright
\(n\)\strut
\end{minipage} & \begin{minipage}[b]{0.15\columnwidth}\raggedright
\(C_{n-1,1}\)\strut
\end{minipage} & \begin{minipage}[b]{0.19\columnwidth}\raggedright
smallest \(\mathbb{R}\) rep.\strut
\end{minipage} & \begin{minipage}[b]{0.04\columnwidth}\centering
M.s.?\strut
\end{minipage} & \begin{minipage}[b]{0.15\columnwidth}\raggedright
\(C_{1,n-1}\)\strut
\end{minipage} & \begin{minipage}[b]{0.19\columnwidth}\raggedright
smallest \(\mathbb{R}\) rep.\strut
\end{minipage} & \begin{minipage}[b]{0.04\columnwidth}\centering
M.s?\strut
\end{minipage}\tabularnewline
\midrule
\endhead
\begin{minipage}[t]{0.04\columnwidth}\raggedright
1\strut
\end{minipage} & \begin{minipage}[t]{0.15\columnwidth}\raggedright
\(\mathbb{R}+\mathbb{R}\)\strut
\end{minipage} & \begin{minipage}[t]{0.19\columnwidth}\raggedright
\(\mathbb{R}\)\strut
\end{minipage} & \begin{minipage}[t]{0.04\columnwidth}\centering
\checkmark\strut
\end{minipage} & \begin{minipage}[t]{0.15\columnwidth}\raggedright
\(\mathbb{C}\)\strut
\end{minipage} & \begin{minipage}[t]{0.19\columnwidth}\raggedright
\(\mathbb{R}^2\)\strut
\end{minipage} & \begin{minipage}[t]{0.04\columnwidth}\centering
\strut
\end{minipage}\tabularnewline
\begin{minipage}[t]{0.04\columnwidth}\raggedright
2\strut
\end{minipage} & \begin{minipage}[t]{0.15\columnwidth}\raggedright
\(\mathbb{R}(2)\)\strut
\end{minipage} & \begin{minipage}[t]{0.19\columnwidth}\raggedright
\(\mathbb{R}^2\)\strut
\end{minipage} & \begin{minipage}[t]{0.04\columnwidth}\centering
\checkmark\strut
\end{minipage} & \begin{minipage}[t]{0.15\columnwidth}\raggedright
\(\mathbb{R}(2)\)\strut
\end{minipage} & \begin{minipage}[t]{0.19\columnwidth}\raggedright
\(\mathbb{R}^2\)\strut
\end{minipage} & \begin{minipage}[t]{0.04\columnwidth}\centering
\checkmark\strut
\end{minipage}\tabularnewline
\begin{minipage}[t]{0.04\columnwidth}\raggedright
3\strut
\end{minipage} & \begin{minipage}[t]{0.15\columnwidth}\raggedright
\(\mathbb{C}(2)\)\strut
\end{minipage} & \begin{minipage}[t]{0.19\columnwidth}\raggedright
\(\mathbb{R}^4\)\strut
\end{minipage} & \begin{minipage}[t]{0.04\columnwidth}\centering
\strut
\end{minipage} & \begin{minipage}[t]{0.15\columnwidth}\raggedright
\(\mathbb{R}(2)+\mathbb{R}(2)\)\strut
\end{minipage} & \begin{minipage}[t]{0.19\columnwidth}\raggedright
\(\mathbb{R}^2\)\strut
\end{minipage} & \begin{minipage}[t]{0.04\columnwidth}\centering
\checkmark\strut
\end{minipage}\tabularnewline
\begin{minipage}[t]{0.04\columnwidth}\raggedright
4\strut
\end{minipage} & \begin{minipage}[t]{0.15\columnwidth}\raggedright
\(\mathbb{H}(2)\)\strut
\end{minipage} & \begin{minipage}[t]{0.19\columnwidth}\raggedright
\(\mathbb{R}^8\)\strut
\end{minipage} & \begin{minipage}[t]{0.04\columnwidth}\centering
\strut
\end{minipage} & \begin{minipage}[t]{0.15\columnwidth}\raggedright
\(\mathbb{R}(4)\)\strut
\end{minipage} & \begin{minipage}[t]{0.19\columnwidth}\raggedright
\(\mathbb{R}^4\)\strut
\end{minipage} & \begin{minipage}[t]{0.04\columnwidth}\centering
\checkmark\strut
\end{minipage}\tabularnewline
\begin{minipage}[t]{0.04\columnwidth}\raggedright
5\strut
\end{minipage} & \begin{minipage}[t]{0.15\columnwidth}\raggedright
\(\mathbb{H}(2)+\mathbb{H}(2)\)\strut
\end{minipage} & \begin{minipage}[t]{0.19\columnwidth}\raggedright
\(\mathbb{R}^8\)\strut
\end{minipage} & \begin{minipage}[t]{0.04\columnwidth}\centering
\strut
\end{minipage} & \begin{minipage}[t]{0.15\columnwidth}\raggedright
\(\mathbb{C}(4)\)\strut
\end{minipage} & \begin{minipage}[t]{0.19\columnwidth}\raggedright
\(\mathbb{R}^8\)\strut
\end{minipage} & \begin{minipage}[t]{0.04\columnwidth}\centering
\strut
\end{minipage}\tabularnewline
\begin{minipage}[t]{0.04\columnwidth}\raggedright
6\strut
\end{minipage} & \begin{minipage}[t]{0.15\columnwidth}\raggedright
\(\mathbb{H}(4)\)\strut
\end{minipage} & \begin{minipage}[t]{0.19\columnwidth}\raggedright
\(\mathbb{R}^{16}\)\strut
\end{minipage} & \begin{minipage}[t]{0.04\columnwidth}\centering
\strut
\end{minipage} & \begin{minipage}[t]{0.15\columnwidth}\raggedright
\(\mathbb{H}(4)\)\strut
\end{minipage} & \begin{minipage}[t]{0.19\columnwidth}\raggedright
\(\mathbb{R}^{16}\)\strut
\end{minipage} & \begin{minipage}[t]{0.04\columnwidth}\centering
\strut
\end{minipage}\tabularnewline
\begin{minipage}[t]{0.04\columnwidth}\raggedright
7\strut
\end{minipage} & \begin{minipage}[t]{0.15\columnwidth}\raggedright
\(\mathbb{C}(8)\)\strut
\end{minipage} & \begin{minipage}[t]{0.19\columnwidth}\raggedright
\(\mathbb{R}^{16}\)\strut
\end{minipage} & \begin{minipage}[t]{0.04\columnwidth}\centering
\strut
\end{minipage} & \begin{minipage}[t]{0.15\columnwidth}\raggedright
\(\mathbb{H}(4)+\mathbb{H}(4)\)\strut
\end{minipage} & \begin{minipage}[t]{0.19\columnwidth}\raggedright
\(\mathbb{R}^{16}\)\strut
\end{minipage} & \begin{minipage}[t]{0.04\columnwidth}\centering
\strut
\end{minipage}\tabularnewline
\begin{minipage}[t]{0.04\columnwidth}\raggedright
8\strut
\end{minipage} & \begin{minipage}[t]{0.15\columnwidth}\raggedright
\(\mathbb{R}(16)\)\strut
\end{minipage} & \begin{minipage}[t]{0.19\columnwidth}\raggedright
\(\mathbb{R}^{16}\)\strut
\end{minipage} & \begin{minipage}[t]{0.04\columnwidth}\centering
\checkmark\strut
\end{minipage} & \begin{minipage}[t]{0.15\columnwidth}\raggedright
\(\mathbb{H}(8)\)\strut
\end{minipage} & \begin{minipage}[t]{0.19\columnwidth}\raggedright
\(\mathbb{R}^{32}\)\strut
\end{minipage} & \begin{minipage}[t]{0.04\columnwidth}\centering
\strut
\end{minipage}\tabularnewline
\bottomrule
\end{longtable}
\noindent
I've noted when the representations are Majorana spinors. Everything
repeats with period 8 after this, in an obvious way.

Finally, sometimes there are ``Weyl spinors'' or ``Majorana--Weyl''
spinors. The point is that sometimes the Dirac spinors, or Majorana
spinors, are a \emph{reducible} representation of
\(\mathrm{Spin}(1,n-1)\). For Dirac spinors this happens in every even
dimension, because the Clifford algebra element
\[\Gamma = e_1 \cdots e_n\] commutes with everything in
\(\mathrm{Spin}(1,n-1)\) and \(\Gamma^2\) is \(1\) or \(-1\), so we can
break the space of Dirac spinors into the two eigenspaces of \(\Gamma\),
which will be smaller reps of \(\mathrm{Spin}(1,n-1)\) --- the ``Weyl
spinors''. Physicists usually call this \(\Gamma\) thing
``\(\gamma_5\)'', and it's an operator that represents parity
transformations. We get ``Majorana--Weyl'' spinors only when we have
Majorana spinors, \(n\) is even, and \(\Gamma^2 = 1\), since we are then
working with real numbers and \(-1\) doesn't have a square root. You can
work out \(\Gamma^2\) for either \(C_{n-1,1}\) or \(C_{1,n-1}\), and see
that we'll only get Majorana--Weyl spinors when \(n = 8k + 2\).

Whew! Let me summarize some of our results:

\begin{longtable}[]{@{}lllll@{}}
\toprule
\(n\) & Dirac & Majorana & Weyl & Majorana--Weyl\tabularnewline
\midrule
\endhead
1 & \(\mathbb{C}\) & \(\mathbb{R}\) & &\tabularnewline
2 & \(\mathbb{C}^2\) & \(\mathbb{R}^2\) & \(\mathbb{C}\) &
\(\mathbb{R}\)\tabularnewline
3 & \(\mathbb{C}^2\) & \(\mathbb{R}^2\) & &\tabularnewline
4 & \(\mathbb{C}^4\) & \(\mathbb{R}^4\) & \(\mathbb{C}^2\)
&\tabularnewline
5 & \(\mathbb{C}^4\) & & &\tabularnewline
6 & \(\mathbb{C}^8\) & & \(\mathbb{C}^4\) &\tabularnewline
7 & \(\mathbb{C}^8\) & & &\tabularnewline
8 & \(\mathbb{C}^{16}\) & \(\mathbb{R}^{16}\) & \(\mathbb{C}^8\)
&\tabularnewline
\bottomrule
\end{longtable}
\noindent
When there are blanks here, the relevant sort of spinor doesn't exist.
Here I'm not distinguishing Majorana spinors that come from
\(C_{n-1,1}\) and those that come from \(C_{1,n-1}\); you can do that
with the previous table. Again, things continue for larger \(n\) in an
obvious way.

Now, let's imagine a theory that has a supersymmetry between a gauge
bosons and a fermion. We'll assume there are as many physical degrees of
freedom for the gauge boson as there are for the fermion. Gauge bosons
have \(n - 2\) physical degrees of freedom in \(n\) dimensions: for example,
in dimension 4 the photon has 2 degrees of freedom, the left and right
polarized states. So we want to find a kind of spinor that has \(n - 2\)
physical degrees of freedom. But the number of physical degrees of
freedom of a spinor field is half the number of (real) components of the
spinor, since the Dirac equation relates the components. So we are
looking for a kind of spinor that has \(2(n - 2)\) real components. This
occurs in only 4 cases:

\begin{itemize}
\item
  \(n = 3\): then \(2(n-2) = 2\), and Majorana spinors have 2 real
  components
\item
  \(n = 4\): then \(2(n-2) = 4\), and Majorana or Weyl spinors have 4
  real components
\item
  \(n = 6\): then \(2(n-2) = 8\), and Weyl spinors have 8 real
  components
\item
  \(n = 10\): then \(2(n-2) = 16\), and Majorana--Weyl spinors have 16
  real components
\end{itemize}
\noindent
Note we count complex components as two real components. And note how
dimension 10 works: the dimension of the spinors grows pretty fast as \(n\)
increases, but the Majorana--Weyl condition reduces the dimension by a
factor of 4, so dimension 10 just squeaks by!

Here John Schwarz explains how nice things happen in the same dimensions
for superstring theory:

\begin{enumerate}
\def\labelenumi{\arabic{enumi})}
\setcounter{enumi}{3}
\tightlist
\item
  John H. Schwarz, ``Introduction to superstrings'', in
  \emph{Superstrings and Supergravity, Proc. of the 28th Scottish
  Universities Summer School in Physics}, ed.~A. T. Davies and D. G.
  Sutherland, University Printing House, Oxford, 1985.
\end{enumerate}

He also makes a tantalizing remark: perhaps these 4 cases correspond
somehow to the reals, complexes, quaternions and octonions. Note:
\(3 = 1 + 2\), \(4 = 2 + 2\), \(6 = 4 + 2\) and \(10 = 8 + 2\). You can
never tell with this stuff... everything is related.

I thank Joshua Burton for helping me overcome my fear of Majorana
spinors, and for correcting a number of embarrassing errors in the first
version of this article.

\begin{center}\rule{0.5\linewidth}{0.5pt}\end{center}

\textbf{Addendum:} In July 2001, long after the above article was
written, Lubos Motl explained where the number 18 shows up in string
theory:

\begin{quote}
Today we know that the two heterotic string theories are related by
various dualities. For example, in 17+1 dimension, the lattices
\(\Gamma_{16}\) and \(\Gamma_8+\Gamma_8\), with an added Lorentzian
\(\Gamma_{1,1}\), become isometric. There is a single even self-dual
lattice in 17+1 dimensions, \(\Gamma_{17,1}\). This is the reason why
two heterotic string theories are T-dual to each other. The
compactification on a circle adds two extra \(\mathrm{U}(1)\)s (from
Kaluza-Klein graviphoton and the B-field), and with appropriate Wilson
lines, a compactification of one heterotic string theory on radius \(R\)
is equivalent to the other on radius \(1/R\), using correct units.
\end{quote}

Also, in \protect\hyperlink{week104}{``Week 104''}, and especially in
the Addendum written by Robert Helling, we'll see that it's \emph{not} a
coincidence that super-Yang--Mills theory works nicely in dimensions that
are 2 more than the dimensions of the reals, complex numbers,
quaternions and octonions.

\begin{center}\rule{0.5\linewidth}{0.5pt}\end{center}

\begin{quote}
\emph{Since the mathematicians have grabbed ahold of the theory of
relativity, I no longer understand it.}

--- Albert Einstein
\end{quote}

\hypertarget{week94}{%
\section{November 11, 1996}\label{week94}}

Today I want to talk a bit about asymptotic freedom.

First of all, remember that in quantum field theory, studying very small
things is the same as studying things at very high energies. The reason
is that in quantum mechanics you need to collide two particles at a
large relative momentum \(p\) to make sure the distance \(x\) between them
gets small, thanks to the uncertainty principle. But in special
relativity the energy \(E\) and momentum \(p\) of a particle of mass
\(m\) are related by \[E^2 = p^2 + m^2,\] in God's units, where the
speed of light is \(1\). So small \(x\) also corresponds to large \(E\).

``Asymptotic freedom'' refers to fact that some forces become very weak
at high energies, or equivalently, at very short distances. The most
interesting example of this is the so-called ``strong force'', which
holds the quarks together in a hadron, like a proton or neutron. True to
its name, it is very strong at distances comparable to the radius of
proton, or at energies comparable to the mass of the proton (where if we
don't use God's units, we have to use \(E = mc^2\) to convert units of
mass to units of energy). But if we smash protons at each other at much
higher energies, the constituent quarks act almost as free particles,
indicating that the strong force gets weak when the quarks get really
close to each other.

Now in \protect\hyperlink{week76}{``Week 76''} and
\protect\hyperlink{week84}{``Week 84''} I talked about another
phenomenon, called ``confinement''. This simply means that at lower
energies, or larger distance scales, the strong force becomes so strong
that it is \emph{impossible} to pull a quark out of a hadron. Asymptotic
freedom and confinement are two aspects of the same thing: the
dependence of the strength of the strong force on the energy scale.
Asymptotic freedom is better understood, though, because the weaker a
force is, the better we can apply the methods of perturbation theory ---
a widely used approach where we try to calculate everything as a Taylor
series in the ``coupling constant'' measuring the strength of the force
in question. This is often successful when the coupling constant is
small, but not when it's big.

The interesting thing is that in quantum field theory the coupling
constants ``run''. This is particle physics slang for the fact that they
depend on the energy scale at which we measure them. ``Asymptotic
freedom'' happens when the coupling constant runs down to zero as we
move up to higher and higher energy scales. If you want to impress
someone about your knowledge of this, just mutter something about the
``beta function'' being negative --- this is a fancy way of saying the
coupling constant decreases as you go to higher energies. You'll sound
like a real expert.

Now, Frank Wilczek is one of the original discoverers of asymptotic
freedom. He \emph{is} a real expert. He recently won a prize for this
work, and he gave a nice talk which he made into a paper:

\begin{enumerate}
\def\labelenumi{\arabic{enumi})}
\tightlist
\item
  Frank Wilczek, ``Asymptotic freedom'', available as
  \href{https://arxiv.org/abs/hep-th/9609099}{\texttt{hep-th/9609099}}.
\end{enumerate}
\noindent
Among other things, he gives a nice summary of the work of Nielsen and
Hughes, which gave the first really easy to understand explanation of
asymptotic freedom. For the original work, try:

\begin{enumerate}
\def\labelenumi{\arabic{enumi})}
\setcounter{enumi}{1}
\item
  N. K. Nielsen, Asymptotic freedom as a spin effect, \emph{Amer. Jour. Phys.} \textbf{49} (1981), 1171--1178.  Available
  at \href{https://inpp.ohio.edu/~inpp/nuclear_lunch/archive/2008/ajp-49-1171-1981.pdf}{\texttt{https://inpp.ohio.edu/\(\sim\)inpp/nuclear\(\underline{\;}\)lunch/archive}} \href{https://inpp.ohio.edu/~inpp/nuclear_lunch/archive/2008/ajp-49-1171-1981.pdf}{\texttt{2008/ajp-49-1171-1981.pdf}}.
\item
  R. J. Hughes, More comments on asymptotic freedom, \emph{Nucl. Phys. B} \textbf{186} (1981), 376--396.
\end{enumerate}

Why would a force get weak at short distance scales? Actually it's
easier to imagine why it would get \emph{strong} --- and sometimes that
is what happens. Of course there are lots of forces that decrease with
distance like \(1/r^2\), but I'm talking about something more drastic:
I'm talking about ``screening''.

For example, say you have an electron in some water. It'll make an
electric field, but this will push all the other negatively charged
particles little bit \emph{away} from your electron and pull all the
positively charged ones a little bit \emph{towards} your electron:

\begin{verbatim}
                                   -
                                     +

                         your electron: -        +-
 
                                            +
                                              -
                      
\end{verbatim}

In other words, it will ``polarize'' all the neighboring water
molecules. But this will create a counteracting sort of electric field,
since it means that if you draw any sphere around your electron, there
will be a bit more \emph{positively} charged other stuff in that sphere
than negatively charged other stuff. The bigger the sphere is, the more
this effect occurs --- though there is a limit to how much it occurs. We
say that the further you go from your electron, the more its electric
charge is ``screened'', or hidden, behind the effect of the
polarization.

This effect is very common in materials that don't conduct electricity,
like water or plastics or glass. They're called ``dielectrics'', and the
dielectric constant, \(\varepsilon\), measures the strength of this
screening effect. Unlike in math, this \(\varepsilon\) is typically
bigger than \(1\). If you apply an electric field to a dielectric
material, the electric field inside the material is only
\(1/\varepsilon\) as big as you'd expect if this polarization wasn't
happening.

What's cool is that according to quantum field theory, screening occurs
even in the vacuum, thanks to ``vacuum polarization''. One can visualize
it rather vaguely as due to a constant buzz of virtual
particle-antiparticle pairs getting created and then annihilating ---
called ``vacuum bubbles'' in the charming language of Feynman diagrams,
because you can draw them like this: \[
  \begin{tikzpicture}
    \begin{knot}
      \strand[thick] (0,0)
        to [out=up,in=up,looseness=2] (1.5,0);
      \strand[thick] (0,0)
        to [out=down,in=down,looseness=2] (1.5,0);
    \end{knot}
    \node[fill=white] at (0.2,0.5) {e\textsuperscript{+}};
    \node[fill=white] at (1.3,0.5) {e\textsuperscript{-}};
  \end{tikzpicture}
\] Here I've drawn a positron-electron pair getting created and then
annihilating as time passes.

There is a lot I should say about virtual particles, and how despite the
fact that they aren't ``real'' they can produce very real effects like
vacuum polarization. A strong enough electric field will even ``spark
the vacuum'' and make the virtual particles \emph{become} real! But
discussing this would be too big of a digression. Suffice it to say that
you have to learn quantum field theory to see how something that starts
out as a kind of mathematical book-keeping device --- a line in a
Feynman diagram --- winds up acting a bit like a real honest particle.
It's a case of a metaphor gone berserk, but in an exceedingly useful
way.

Anyway, so much for screening. Asymptotic freedom requires something
opposite, called ``anti-screening''! That's why it's harder to
understand.

Nielsen and Hughes realized that anti-screening is easier to understand
using magnetism than electricity. In analogy to dielectrics, there are
some materials that screen magnetic fields, and these are called
``diamagnetic'' --- for example, one of the strongest diamagnets is
bismuth. But in addition, there are materials that ``anti-screen''
magnetic fields --- the magnetic field inside them is stronger than the
externally applied magnetic field --- and these are called
``paramagnetic''. For example, aluminum is paramagnetic. People keep
track of paramagnetism using a constant called the magnetic
permeability, \(\mu\). Just to confuse you, this works the opposite way
from the dielectric constant. If you apply a magnetic field to some
material, the magnetic field inside it is \(\mu\) times as big as you'd
expect if there were no magnetic effects going on.

The nice thing is that there are lots of examples of paramagnetism, and
we can sort of understand it if we think about it. It turns out that
paramagnetism in ordinary matter is due to the spin of the electrons in
it. The electrons are like little magnets --- they have a little
``magnetic moment'' pointing along the axis of their spin. Actually,
purely by convention it points in the direction opposite their spin,
since for some stupid reason Benjamin Franklin decided to decree that
electrons were \emph{negative}. But don't worry about this --- it
doesn't really matter. The point is that when you put electrons in a
magnetic field, their spins like to line up in such a way that their
magnetic field points the same way as the externally applied magnetic
field, just like a compass needle does in the Earth's magnetic field. So
they \emph{add} to the magnetic field. Ergo, paramagnetism.

Now, spin is a form of angular momentum intrinsic to the electron, but
there is another kind of angular momentum, namely orbital angular
momentum, caused by how the electron (or whatever particle) is moving
around in space. It turns out that orbital angular momentum also has
magnetic effects, but only causes diamagnetism. The idea is that when you
apply a magnetic field to some material, it can also make the electrons
in it tend to move in orbits perpendicular to the magnetic field, and
the resulting current creates a magnetic field. But this magnetic field
must \emph{oppose} the external magnetic field. Ergo, diamagnetism.

Why does orbital angular momentum work one way, while spin works the
other way? I'll say a bit more about that later. Now let me get back to
asymptotic freedom.

I've talked about screening and antiscreening for both electric and
magnetic fields now. But say the ``substance'' we're studying is the
\emph{vacuum}. Unlike most substances, the vacuum doesn't look different
when we look at it from a moving frame of reference. We say it's
``Lorentz-invariant''. But if we look at an electric field in a moving
frame of reference, we see a bit of magnetic field added on, and vice
versa. We say that the electric and magnetic fields transform into each
other... they are two aspects of single thing, the electromagnetic
field. So the amount of \emph{electric} screening or antiscreening in
the vacuum has to equal the amount of the \emph{magnetic} screening or
antiscreening. In other words, thanks to the silly way we defined
\(\varepsilon\) differently from \(\mu\), we must have
\[\varepsilon = 1/\mu\] in the vacuum.

Now the cool thing is that the Yang--Mills equations, which describe the
strong force, are very similar to Maxwell's equations. In particular,
the strong force, also known as the ``color'' force, consists of two
aspects, the ``chromoelectric'' field and ``chromomagnetic'' field.
Moreover, the same argument above applies here: the vacuum must give the
same antiscreening for the chromoelectric field as it does for the
chromomagnetic field, so \(\varepsilon = 1/\mu\) here too.

So to understand asymptotic freedom it is sufficient to see why the
vacuum acts like a paramagnet for the strong force! This depends on a
big difference between the strong force and electromagnetism. Just as
the electromagnetic field is carried by photons, which are spin-\(1\)
particles, the strong force is carried by ``gluons'', which are also
spin-\(1\) particles. But while the photon is electrically uncharged,
the gluon is charged as far as the strong force goes: we say it has
``color''.

The vacuum is bustling with virtual gluons. When we apply a
chromomagnetic field to the vacuum, we get two competing effects:
paramagnetism thanks to the \emph{spin} of the gluons, and diamagnetism
due to their \emph{orbital angular momentum}. But --- the spin effect is
stronger. The vacuum acts like a paramagnet for the strong force. So we
get asymptotic freedom!

That's the basic idea. Of course, there are some loose ends. To see why
the spin effect is stronger, you have to calculate a bit. At least I
don't know how to see it without calculating --- but Wilczek sketches
the calculation, and it doesn't look too bad. It's also true in most
metals that the spin effect wins, so they are paramagnetic.

You might also wonder why spin and orbital angular momentum work
oppositely as far as magnetism goes. Unfortunately I don't have any
really simple slick answer. One thing is that it seems any answer must
involve quantum mechanics. {[}Note: later I realized some very basic
things about this, which I append below.{]} In volume II of his
magnificent series:

\begin{enumerate}
\def\labelenumi{\arabic{enumi})}
\setcounter{enumi}{3}
\tightlist
\item
  Richard Feynman, Robert Leighton, and Matthew Sands, \emph{The Feynman
  Lectures on Physics}, Addison-Wesley, Reading, 1964.
\end{enumerate}

Feynman notes: ``It is a consequence of classical mechanics that if you
have any kind of system --- a gas with electrons, protons, and whatever
--- kept in a box so that the whole thing can't turn, there will be no
magnetic effect. {[}\ldots.{]} The theorem then says that if you turn on
a magnetic field and wait for the system to get into thermal
equilibrium, there will be no paramagnetism or diamagnetism --- there
will be no induced magnetic moment. Proof: According to statistical
mechanics, the probability that a system will have any given state of
motion is proportional to \(\exp(-U/kT)\), where \(U\) is the energy of
that motion. Now what is the energy of motion. For a particle moving in
a constant magnetic field, the energy is the ordinary potential energy
plus \(mv^2/2\), with nothing additional for the magnetic field. (You
know that the forces from electromagnetic fields are
\(q(E + v \times B)\), and that the rate of work \(F\cdot v\) is just
\(qE\cdot v\), which is not affected by the magnetic field.) So the
energy of a system, whether it is in a magnetic field or not, is always
given by the kinetic energy plus the potential energy. Since the
probability of any motion depends only on the energy --- that is, on the
velocity and position --- it is the same whether or not there is a
magnetic field. For \emph{thermal} equilibrium, therefore, the magnetic
field has no effect.''

So to understand magnetism we really need to work quantum-mechanically.
Laurence Yaffe has brought to my attention a nice path-integral argument
as to why orbital angular momentum can only yield diamagnetism; this can
be found in his charming book:

\begin{enumerate}
\def\labelenumi{\arabic{enumi})}
\setcounter{enumi}{4}
\tightlist
\item
  Barry Simon, \emph{Functional Integration and Quantum Physics},
  Academic Press, New York, 1979.
\end{enumerate}
\noindent
This argument is very simple if you know about path integrals, but I
think there should be some more lowbrow way to see it, too. I think it's
good to make all this stuff as simple as possible, because the phenomena
of asympotic freedom and confinement are very important and shouldn't
only be accessible to experts.

I'd like to thank Douglas Singleton, Matt McIrvin, Mike Kelsey, and
Laurence Yaffe for some posts on sci.physics.research that helped me
understand this stuff.

\begin{center}\rule{0.5\linewidth}{0.5pt}\end{center}

\textbf{Addendum} \emph{(November 13, 1996)}. Thanks to emails from
Yehuda Naveh and Bruce Smith I'm beginning to understand this stuff at
the 13-year-old level it deserves. If you want to jump to the punchline,
skip down to the stuff between double lines --- that's the part I should
have known ages ago!

Here's the deal. Feynman's theorem deals with classical systems made
only of a bunch of electrically charged point particles. Remember how it
goes: A magnetic field can never do work on such a system, because it
always exerts a force perpendicular to the velocity of an electrically
charged particle. So the energy of such a system is independent of the
externally applied magnetic field. Now, in statistical mechanics the
equilibrium state of a system depends only on the energy of each state,
since the probability of being in a state with energy \(E\) is
proportional to \(\exp(-E/kT)\). So an external magnetic field doesn't
affect the equilibrium state of this sort of system. So there can't be
anything like paramagnetism or diamagnetism, where the equilibrium state
is affected by an external magnetic field.

But suppose instead we allowed an extra sort of building block of our
system, in addition to electrically charged particles. Suppose we allow
little ``current loops''. We take these as ``primitives'', in the sense
that we don't ask how or why the current keeps flowing around the loop,
we just assume it does. We just \emph{define} one of these ``current
loops'' to be a little circle of stuff with a constant mass per unit
length, with a constant current that flows around it. This may or may
not be physically reasonable, but we're gonna do it anyway!

Note: If we tried to make a current loop out of classical electrically
charged point particles, the current loop would tend to fall apart! A
loop is not going to be the equilibrium state of a bunch of charged
particles. So we are going to get around this by taking current loops as
new primitives --- simply \emph{assuming} they exist and have the
properties given above.

If we build our system out of current loops and point particles,
Feynman's theorem no longer applies. Why? Well, a constant magnetic
field exerts a force perpendicular to the direction of the current, and
this applies a \emph{torque} to the current loop --- no net force, just
a torque. But since the current loop is made out of stuff that has a
constant mass per unit length, when the current loop is rotating it will
have kinetic energy. So by applying a torque to the current loop, the
magnetic field does \emph{work} on the current loop. Thus Feynman's
reasoning no longer applies to this case.

In particular, what happens is just what we expect. The torque on the
little current loops makes them want to line up with the external
magnetic field. In other words, they will have less energy when they are
lined up like this. In particular, the energy of the system \emph{does}
depend on the external magnetic field, and the equilibrium state will
tend to have more little current loops lined up with the field than not.

Now if we keep track of the magnetic field produced by these current
loops, we see it points the same way as the externally applied field. So
we get paramagnetism.

Now, even without doing a detailed quantum-mechanical treatment of this
problem, we see what's special about spin: a particle with spin is a bit
like one of our imaginary ``primitive current loops''. This is how spin
can give paramagnetism.

Great. But what had always been bugging me is this! If you put a charged
particle in a constant magnetic field, it moves in a circular or spiral
orbit. For simplicity let's say it moves in a circle. You can think of
this, if you like, as a kind of current loop --- but a very different
sort of current loop than the one we've just been considering! In
particular, if you work it out, this particle circling around will
produce a magnetic field that \emph{opposes} the external magnetic
field. On the other hand, our primitive current loops are in the state
of least energy when they're lined up to produce a magnetic field that
\emph{goes with} the external field.

What's the deal? Well, it's just something about how the vector cross
product works; you gotta work it out yourself to believe it. All you
need to know is that the force on a charged particle is
\(q v \times B\). It boils down to this:

\begin{center}\rule{0.5\linewidth}{0.5pt}\end{center}

\begin{quote}
A positively charged particle orbiting in a magnetic field pointing
along the \(z\)-axis will orbit CLOCKWISE in the \(xy\) plane.
However, a primitive current loop in a magnetic field pointing along the
\(z\)-axis will be in its state of least energy when the current runs
COUNTERCLOCKWISE in the \(xy\) plane.
\end{quote}

\begin{center}\rule{0.5\linewidth}{0.5pt}\end{center}
\noindent
I'm sure this is what was nagging at me. It's just one of those basic
funny little things. If I'm still mixed up, someone had better let me
know.

There are a couple other things perhaps worth saying about this:

\begin{enumerate}
\def\labelenumi{\arabic{enumi}.}
\tightlist
\item
  In our calculation of the energy of the system, we have been
  neglecting the energy due to the electric and magnetic fields
  \emph{produced} by our point particles and current loops. A more
  careful analysis would take these into account. In particular, the
  reason ferromagnets prefer to have lots of ``domains'' than to have
  all their little current loops lined up, is to keep the energy due to
  the magnetic field produced by these loops from getting too big.
\item
  A little current loop acts like a magnetic dipole. We'd also get
  interesting effects if we had magnetic monopoles. Here I simply assume
  that, just as an electric field exerts a force on a electrically
  charged particle equal to \(q E\), a magnetic field exerts a force on
  a magnetically charged particle equal to \(m B\), where \(m\) is the
  magnetic charge. A magnetic field would then be able to do work on a
  magnetic monopole, and again Feynman's theorem would not apply. So
  it's perhaps not so surprising that Feynman's theorem fails when we
  have magnetic dipoles as primitive constituents of our system, too
  (although these dipoles had better not be points --- they need a
  moment of inertia for a torque on them to do work).
\end{enumerate}

\hypertarget{week95}{%
\section{November 26, 1996}\label{week95}}

Last week I talked about asymptotic freedom --- how the ``strong'' force
gets weak at high energies. Basically, I was trying to describe an
aspect of ``renormalization'' without getting too technical about it. By
coincidence, I recently got my hands on a book I'd been meaning to read
for quite a while:

\begin{enumerate}
\def\labelenumi{\arabic{enumi})}
\tightlist
\item
  Laurie M. Brown, ed., \emph{Renormalization: From Lorentz to Landau
  (and Beyond)}, Springer, Berlin, 1993.
\end{enumerate}
\noindent
It's a nice survey of how attitudes to renormalization have changed over
the years. It's probably the most fun to read if you know some quantum
field theory, but it's not terribly technical, and it includes a
``Tutorial on infinities in QED'', by Robert Mills, that might serve as
an introduction to renormalization for folks who've never studied it.

Okay, on to some new stuff\ldots.

It's a bit funny how one of the most curious features of bosonic string
theory in 26 dimensions was anticipated by the number theorist Edouard
Lucas in 1875. I assume this is the same Lucas who is famous for the
Lucas numbers: 1,3,4,7,11,18,\ldots, each one being the sum of the
previous two, after starting off with 1 and 3. They are not quite as
wonderful as the Fibonacci numbers, but in a study of pine cones it was
found that while \emph{most} cones have consecutive Fibonacci numbers of
spirals going around clockwise and counterclockwise, a small minority of
deviant cones use Lucas numbers instead.

Anyway, Lucas must have liked playing around with numbers, because in
one publication he challenged his readers to prove that: ``A square
pyramid of cannon balls contains a square number of cannon balls only
when it has 24 cannon balls along its base''. In other words, the only
integer solution of \[1^2 + 2^2 + \cdots + n^2 = m^2,\]
is the solution \(n = 24\), not counting silly solutions like \(n=0\)
and \(n=1\).

It seems that Lucas didn't have a proof of this; the first proof is due
to G. N. Watson in 1918, using elliptic functions. Apparently an
elementary proof appears in the following ridiculously overpriced book:

\begin{enumerate}
\def\labelenumi{\arabic{enumi})}
\setcounter{enumi}{1}
\tightlist
\item
  W. S. Anglin, \emph{The Queen of Mathematics: An Introduction to
  Number Theory}, Kluwer, Dordrecht, 1995.
\end{enumerate}
\noindent
For more historical details, see the review in

\begin{enumerate}
\def\labelenumi{\arabic{enumi})}
\setcounter{enumi}{2}
\tightlist
\item
  Jet Wimp, ``Eight recent mathematical books'', \emph{Math.
  Intelligencer} \textbf{18} (1996), 72--79.
\end{enumerate}
\noindent
Unfortunately, I haven't seen these proofs of Lucas' claim, so I don't
know why it's true. I do know a little about its relation to string
theory, so I'll talk about that.

There are two main flavors of string theory, ``bosonic'' and
``supersymmetric''. The first is, true to its name, just the quantized,
special-relativistic theory of little loops made of some abstract string
stuff that has a certain tension --- the ``string tension''. Classically
this theory would make sense in any dimension, but quantum-mechanically,
for reasons that I want to explain \emph{someday} but not now, this
theory works best in 26 dimensions. Different modes of vibration of the
string correspond to different particles, but the theory is called
``bosonic'' because these particles are all bosons. That's no good for a
realistic theory of physics, because the real world has lots of
fermions, too. (For a bit about bosons and fermions in particle physics,
see \protect\hyperlink{week93}{``Week 93''}.)

For a more realistic theory people use ``supersymmetric'' string theory.
The idea here is to let the string be a bit more abstract: it vibrates
in ``superspace'', which has in addition to the usual coordinates some
extra ``fermionic'' coordinates. I don't want to get too technical here,
but the basic idea is that while the usual coordinates commute as usual:
\[x_i x_j = x_j x_i\] the fermionic coordinates ``anticommute''
\[y_i y_j = -y_j y_i\] while the bosonic coordinates commute with
fermionic ones: \[x_i y_j = y_j x_i.\] If you've studied bosons and
fermions this will be sort of familiar; all the differences between them
arise from the difference between commuting and anticommuting variables.
For a little glimpse of this subject try:

\begin{enumerate}
\def\labelenumi{\arabic{enumi})}
\setcounter{enumi}{3}
\tightlist
\item
  John Baez, Spin and the harmonic oscillator,
  \href{http://math.ucr.edu/home/baez/harmonic.html}{\texttt{http://math.ucr.edu/home/baez/}} \href{http://math.ucr.edu/home/baez/harmonic.html}{\texttt{harmonic.html}}.
\end{enumerate}

As it so happens, supersymmetric string theory --- often abbreviated to
``superstring theory'' --- works best in 10 dimensions. There are five
main versions of superstring theory, which I described in
\protect\hyperlink{week74}{``Week 74''}. The type I string theory
involves open strings --- little segments rather than loops. The type
IIA and type IIB theories involve closed strings, that is, loops. But
the most popular sort of superstring theories are the ``heterotic
strings''. A nice introduction to these, written by one of their
discoverers, is:

\begin{enumerate}
\def\labelenumi{\arabic{enumi})}
\setcounter{enumi}{4}
\tightlist
\item
  David J. Gross, ``The heterotic string'', in \emph{Workshop on Unified
  String Theories}, eds.~M. Green and D. Gross, World Scientific,
  Singapore, 1986, pp.~357--399.
\end{enumerate}
\noindent
These theories involve closed strings, but the odd thing about them,
which accounts for the name ``heterotic'', is that vibrations of the
string going around one way are supersymmetric and act as if they were
in 10 dimensions, while the vibrations going around the other way are
bosonic and act as if they were in 26 dimensions!

To get this string with a split personality to make sense, people
cleverly think of the 26 dimensional spacetime for the bosonic part as a
10-dimensional spacetime times a little \(16\)-dimensional curled-up
space, or ``compact manifold''. To get the theory to work, it seems that
this compact manifold needs to be flat, which means it has to be a torus
--- a 16-dimensional torus. We can think of any such torus as
\(16\)-dimensional Euclidean space ``modulo a lattice''. Remember, a
lattice in Euclidean space is something that looks sort of like this: \[
  \begin{tikzpicture}[scale=0.7]
    \draw[->] (-3,0) to (4,0) node[label=below:{$x$}]{};
    \draw[->] (0,-3) to (0,4) node[label=left:{$y$}]{};
    \foreach \m in {-1,0,1,2}
    {
      \foreach \n in {-1,0,1,2}
      {
        \node at ({\m*1.5-\n/3-0.2},{1.5*\n+\m-0.5}) {$\bullet$};
      }
    }
  \end{tikzpicture}
\] Mathematically, it's just a discrete subset \(L\) of \(\mathbb{R}^n\)
(\(n\)-dimensional Euclidean space, with its usual coordinates) with the
property that if \(x\) and \(y\) lie in \(L\), so does \(jx + ky\) for
all integers \(j\) and \(k\). When we form \(n\)-dimensional Euclidean
space ``modulo a lattice'', we decree two points \(x\) and \(y\) to be
the same if \(x-y\) is in \(L\). For example, all the points labelled
\(x\) in the figure above count as the same when we ``mod out by the
lattice''... so in this case, we get a \(2\)-dimensional torus.

For more on \(2\)-dimensional tori and their relation to complex
analysis, you can read \protect\hyperlink{week13}{``Week 13''}. Here we
are going to be macho and plunge right into talking about lattices and
tori in arbitrary dimensions.

To get our \(26\)-dimensional string theory to work out nicely when we
curl up \(16\)-dimensional space to a \(16\)-dimensional torus, it turns
out that we need the lattice \(L\) that we're modding out by to have
some nice properties. First of all, it needs to be an ``integral''
lattice, meaning that for any vectors \(x\) and \(y\) in \(L\) the dot
product \(x\cdot y\) must be an integer. This is no big deal --- there
are gadzillions of integral lattices. In fact, sometimes when people say
``lattice'' they really mean ``integral lattice''. What's more of a big
deal is that \(L\) must be ``even'', that is, for any \(x\) in \(L\) the
inner product \(x\cdot x\) is even. This implies that \(L\) is integral,
by the identity
\[(x + y)\cdot (x + y) = x\cdot x + 2x\cdot y + y\cdot y.\] But what's
really a big deal is that \(L\) must also be ``unimodular''. There are
different ways to define this concept. Perhaps the easiest to grok is
that the volume of each lattice cell --- e.g., each parallelogram in the
picture above --- is \(1\). Another way to say it is this. Take any
basis of \(L\), that is, a bunch of vectors in \(L\) such that any
vector in \(L\) can be uniquely expressed as an integer linear
combination of these vectors. Then make a matrix with the components of
these vectors as rows. Then take its determinant. That should equal plus
or minus \(1\). Still another way to say it is this. We can define the
``dual'' of \(L\), say \(L^*\), to be all the vectors \(x\) such that
\(x\cdot y\) is an integer for all \(y\) in \(L\). An integer lattice is
one that's contained in its dual, but \(L\) is unimodular if and only if
\(L = L^*\). So people also call unimodular lattices ``self-dual''. It's
a fun little exercise in linear algebra to show that all these
definitions are equivalent.

Why does \(L\) have to be an even unimodular lattice? Well, one can
begin to understand this a litle by thinking about what a closed string
vibrating in a torus is like. If you've ever studied the quantum
mechanics of a particle on a torus (e.g.~a circle!) you may know that
its momentum is quantized, and must be an element of \(L^*\). So the
momentum of the center of mass of the string lies in \(L^*\).

On the other hand, the string can also wrap around the torus in various
topologically different ways. Since two points in Euclidean space
correspond to the same point in the torus if they differ by a vector in
\(L\), if we imagine the string as living up in Euclidean space, and
trace our finger all around it, we don't necesarily come back to the
same point in Euclidean space: the same point \emph{plus} any vector in
\(L\) will do. So the way the string wraps around the torus is described
by a vector in \(L\). If you've heard of the ``winding number'', this is
just a generalization of that.

So both \(L\) and \(L^*\) are really important here (which has to do
with the fashionable subject of ``string duality''), and a bunch more
work shows that they \emph{both} need to be even, which implies that
\(L\) is even and unimodular.

Now something cool happens: even unimodular lattices are only possible
in certain dimensions --- namely, dimensions divisible by 8. So we luck
out, since we're in dimension 16.

In dimension 8 there is only \emph{one} even unimodular lattice (up to
isometry), namely the wonderful lattice \(\mathrm{E}_8\)! The easiest
way to think about this lattice is as follows. Say you are packing
spheres in \(n\) dimensions in a checkerboard lattice --- in other words,
you color the cubes of an \(n\)-dimensional checkerboard alternately red
and black, and you put spheres centered at the center of every red cube,
using the biggest spheres that will fit. There are some little hole left
over where you could put smaller spheres if you wanted. And as you go up
to higher dimensions, these little holes get bigger! By the time you
get up to dimension 8, there's enough room to put another sphere \emph{of the same size as the rest} in each hole! If you do that, you get the lattice
\(\mathrm{E}_8\). (I explained this and a bunch of other lattices in
\protect\hyperlink{week65}{``Week 65''}, so for more info take a look at
that.)

In dimension 16 there are only \emph{two} even unimodular lattices. One
is \(\mathrm{E}_8\oplus\mathrm{E}_8\). A vector in this is just a pair
of vectors in \(\mathrm{E}_8\). The other is called
\(\mathrm{D}_{16}^+\), which we get the same way as we got
\(\mathrm{E}_8\): we take a checkerboard lattice in 16 dimensions and
stick in extra spheres in all the holes. More mathematically, to get
\(\mathrm{E}_8\) or \(\mathrm{D}_{16}^+\), we take all vectors in
\(\mathbb{R}^8\) or \(\mathbb{R}^{16}\), respectively, whose coordinates
are either \emph{all} integers or \emph{all} half-integers, for which
the coordinates add up to an even integer. (A ``half-integer'' is an
integer plus \(1/2\).)

So \(\mathrm{E}_8\oplus\mathrm{E}_8\) and \(\mathrm{D}_{16}^+\) give us
the two kinds of heterotic string theory! They are often called the
\(\mathrm{E}_8\oplus\mathrm{E}_8\) and \(\mathrm{SO}(32)\) heterotic
theories.

In \protect\hyperlink{week63}{``Week 63''} and
\protect\hyperlink{week64}{``Week 64''} I explained a bit about lattices
and Lie groups, and if you know about that stuff, I can explain why the
second sort of string theory is called ``\(\mathrm{SO}(32)\)''. Any
compact Lie group has a maximal torus, which we can think of as some
Euclidean space modulo a lattice. There's a group called
\(\mathrm{E}_8\), described in \protect\hyperlink{week90}{``Week 90''},
which gives us the \(\mathrm{E}_8\) lattice this way, and the product of
two copies of this group gives us \(\mathrm{E}_8\oplus\mathrm{E}_8\). On
the other hand, we can also get a lattice this way from the group
\(\mathrm{SO}(32)\) of rotations in 32 dimensions, and after a little
finagling this gives us the \(\mathrm{D}_{16}^+\) lattice (technically,
we need to use the lattice generated by the weights of the adjoint
representation and one of the spinor representations, according to
Gross). In any event, it turns out that these two versions of heterotic
string theory act, at low energies, like gauge field theories with gauge
group \(\mathrm{E}_8\times\mathrm{E}_8\) and \(\mathrm{SO}(32)\),
respectively! People seem especially optimistic that the
\(\mathrm{E}_8\times\mathrm{E}_8\) theory is relevant to real-world
particle physics; see for example:

\begin{enumerate}
\def\labelenumi{\arabic{enumi})}
\setcounter{enumi}{5}
\tightlist
\item
  Edward Witten, ``Unification in ten dimensions'', in \emph{Workshop on
  Unified String Theories}, eds.~M. Green and D. Gross, World
  Scientific, Singapore, 1986, pp.~438--456.

Edward Witten, ``Topological tools in ten dimensional physics'', with an
appendix by R. E. Stong, in \emph{Workshop on Unified String Theories},
eds.~M. Green and D. Gross, World Scientific, Singapore, 1986,
pp.~400--437.
\end{enumerate}

The first paper listed here is about particle physics; I mention the
second here just because \(\mathrm{E}_8\) fans should enjoy it --- it
discusses the classification of bundles with \(\mathrm{E}_8\) as gauge
group.

Anyway, what does all this have to do with Lucas and his stack of cannon
balls?

Well, in dimension 24, there are \emph{24} even unimodular lattices,
which were classified by Niemeier. A few of these are obvious, like
\(\mathrm{E}_8\oplus\mathrm{E}_8\oplus\mathrm{E}_8\) and
\(\mathrm{E}_8\oplus\mathrm{D}_{16}^+\), but the coolest one is the
``Leech lattice'', which is the only one having no vectors of length 2.
This is related to a whole \emph{world} of bizarre and perversely fascinating
mathematics, like the ``Monster group'', the largest sporadic finite
simple group --- and also to string theory. I said a bit about this
stuff in \protect\hyperlink{week66}{``Week 66''}, and I will say more in
the future, but for now let me just describe how to get the Leech
lattice.

First of all, let's think about Lorentzian lattices, that is, lattices
in Minkowski spacetime instead of Euclidean space. The difference is
just that now the dot product is defined by
\[(x_1,\ldots,x_n)\cdot(y_1,\ldots,y_n) = -x_1y_1+x_2 y_2+\cdots+x_ny_n\]
with the first coordinate representing time. It turns out that the only
even unimodular Lorentzian lattices occur in dimensions of the form
\(8k + 2\). There is only \emph{one} in each of those dimensions, and it
is very easy to describe: it consists of all vectors whose coordinates
are either all integers or all half-integers, and whose coordinates add
up to an even number.

Note that the dimensions of this form: 2, 10, 18, 26, etc., are
precisely the dimensions I said were specially important in
\protect\hyperlink{week93}{``Week 93''} for some \emph{other}
string-theoretic reason. Is this a ``coincidence''? Well, all I can say
is that I don't understand it.

Anyway, the \(10\)-dimensional even unimodular Lorentzian lattice is
pretty neat and has attracted some attention in string theory:

\begin{enumerate}
\def\labelenumi{\arabic{enumi})}
\setcounter{enumi}{6}
\tightlist
\item
  Reinhold W. Gebert and Hermann Nicolai, ``\(\mathrm{E}_{10}\) for
  beginners'', available as
  \href{https://arxiv.org/abs/hep-th/9411188}{\texttt{hep-th/9411188}}
\end{enumerate}
\noindent
but the \(26\)-dimensional one is even more neat. In particular, thanks
to the cannonball trick of Lucas, the vector
\[v = (70,0,1,2,3,4,\ldots,24)\] is ``lightlike''. In other words,
\[v\cdot v=0.\] What this implies is that if we let \(T\) be the set of
all integer multiples of \(v\), and let \(S\) be the set of all vectors
\(x\) in our lattice with \(x\cdot v = 0\), then \(T\) is contained in
\(S\), and \(S/T\) is a \(24\)-dimensional lattice --- the Leech
lattice!

Now \emph{that} has all sorts of ramifications that I'm just barely
beginning to understand. For one, it means that if we do bosonic string
theory in 26 dimensions on \(\mathbb{R}^{26}\) modulo the
\(26\)-dimensional even unimodular lattice, we get a theory having lots
of symmetries related to those of the Leech lattice. In some sense this
is a ``maximally symmetric'' approach to \(26\)-dimensional bosonic
string theory:

\begin{enumerate}
\def\labelenumi{\arabic{enumi})}
\setcounter{enumi}{7}
\tightlist
\item
  Gregory Moore, ``Finite in all directions'', available as
  \href{https://arxiv.org/abs/hep-th/9305139}{\texttt{hep-th/9305139}}.
\end{enumerate}
\noindent
Indeed, the Monster group is lurking around as a symmetry group here!
For a physics-flavored introduction to that aspect, try:

\begin{enumerate}
\def\labelenumi{\arabic{enumi})}
\setcounter{enumi}{8}
\tightlist
\item
  Reinhold W. Gebert, ``Introduction to vertex algebras, Borcherds
  algebras, and the Monster Lie algebra'', \emph{Int. Jour. Mod. Phys. A}
  \textbf{8} (1993), 5441--5503.   Also available as
  \href{https://arxiv.org/abs/hep-th/9308151}{\texttt{hep-th/9308151}}.
\end{enumerate}
\noindent
and for a detailed mathematical tour see:

\begin{enumerate}
\def\labelenumi{\arabic{enumi})}
\setcounter{enumi}{9}
\tightlist
\item
  Igor Frenkel, James Lepowsky, and Arne Meurman, \emph{Vertex Operator
  Algebras and the Monster}, Academic Press, New York, 1988.
\end{enumerate}
\noindent
Also try the very readable review articles by Richard Borcherds, who
came up with a lot of this business:

\begin{enumerate}
\def\labelenumi{\arabic{enumi})}
\setcounter{enumi}{10}
\item
  Richard Borcherds, ``Automorphic forms and Lie algebras'', available at \hfill \break
  \href{https://math.berkeley.edu/~reb/papers/index.html#preprints}{\texttt{https://math.berkeley.edu/\(\sim\)reb/papers/index.html\#preprints}}.
  
  Richard Borcherds, ``Sporadic groups and string theory'', available at \hfill \break
  \href{https://math.berkeley.edu/~reb/papers/index.html#preprints}{\texttt{https://math.berkeley.edu/\(\sim\)reb/papers/index.html\#preprints}}.
\end{enumerate}

Well, there is a lot more to say, but I need to go home and pack for my
Thanksgiving travels. Let me conclude by answering a natural followup
question: how many even unimodular lattices are there in 32 dimensions?
Well, one can show that there are \emph{at least 80 million!}

Some of you may have wondered what's happened to the ``Tale of
\(n\)-Categories''. I haven't forgotten that! In fact, earlier this fall
I finished writing a big fat paper on 2-Hilbert spaces (which are to
Hilbert spaces as categories are to sets), and since then I have been
struggling to finish another big fat paper with James Dolan, on the
general definition of ``weak \(n\)-categories''. I want to talk about
this sort of thing, and other progress on \(n\)-categories and physics,
but I've been so busy \emph{working} on it that I haven't had time to
\emph{chat} about it on This Week's Finds. Maybe soon I'll find time.

\begin{center}\rule{0.5\linewidth}{0.5pt}\end{center}

\textbf{Addenda:} Robin Chapman pointed out that Anglin's proof also
appears in the \emph{American Mathematical Monthly}, February 1990,
pp.~120--124, and that another elementary proof has subsequently
appeared in the \emph{Journal of Number Theory}. David Morrison pointed out in
email that since the sum \[1^2 + 2^2 + \cdots + n^2 = m^2\] is
\(n(n+1)(2n+1)/6\), this problem can be solved by finding all the
rational points \((n,m)\) on the elliptic curve
\[\frac{n^3}{3} + \frac{n^2}2 + \frac{n}{6} = m^2\] which is the sort of
thing folks know how to do.

Also, here's something Michael Thayer wrote on one of the newsgroups,
and my reply:

\begin{quote}
John Baez wrote:

\begin{quote} In particular, thanks to the cannonball trick of Lucas, the vector
\[          v = (70,0,1,2,3,4,\ldots,24) \]
is ``lightlike".  In other words,  \(v \cdot v = 0\).
\end{quote}

I don't see what is so significant about the vector v.
For instance, the 10 dimensional vector
(3,1,1,1,1,1,1,1,1,1) is also light like, and you make
no big deal about that. Is there some reason why the
ascending values in \(v\) are important?
\end{quote}

\noindent
Yikes! Thanks for catching that massive hole in the exposition.

You're right that there's no shortage of lightlike vectors in the even
unimodular Lorentzian lattices of other dimensions \(8n+2\); there are
also lots of other lightlike vectors in the \(26\)-dimensional one. Any
one of these gives us a lattice in \(8n\)-dimensional Euclidean space.
In fact, we can get all 24 even unimodular lattices in
\(24\)-dimensional Euclidean space by suitable choices of lightlike
vector. The lightlike vector you wrote down happens to give us the
\(\mathrm{E}_8\) lattice in 8 dimensions.

So what's so special about I wrote, which gives the Leech lattice? Of
course the Leech lattice is itself special, but what does this have to
do with the nicely ascending values of the components of \(v\)?

Alas, I don't know the real answer. I'm not an expert on this stuff; I'm
just explaining it in order to try to learn it. Let me just say what I
know, which all comes from Chap. 27 of Conway and Sloane's book \emph{Sphere
Packings, Lattices, and Groups}.

If we have a lattice, we say a vector \(r\) in it is a ``root'' if the
reflection through \(r\) is a symmetry of the lattice. Corresponding to
each root is a hyperplane consisting of all vectors perpendicular to
that root. These chop space into a bunch of ``fundamental regions''. If
we pick a fundamental region, the roots corresponding to the hyperplanes
that form the walls of this region are called ``fundamental roots''. The
nice thing about the fundamental roots is that the reflection through
any root is a product of reflections through these fundamental roots.

(For more stuff on reflection groups and lattices see
\protect\hyperlink{week62}{``Week 62''} and the following weeks.)

In 1983 John Conway published a paper where he showed various amazing
things; this is now Chapter 27 of the above book. First, he shows that
the fundamental roots of the even unimodular Lorentzian lattices in
dimensions 10, 18, and 26 are the vectors \(r\) with \(r\cdot r = 2\)
and \(r\cdot v = -1\), where the ``Weyl vector'' \(v\) is
\[(28,0,1,2,3,4,5,6,7,8)\] \[(46,0,1,2,3,\ldots,16)\] and
\[(70,0,1,2,3,\ldots,70)\] respectively.

They all have this nice ascending form but only in 26 dimensions is the
Weyl vector lightlike!

Howerver, Conway doesn't seem to explain \emph{why} the Weyl vectors
have this ascending form. So I'm afraid I really don't understand how
all the pieces fit together. All I can say is that for some reason the
Weyl vectors have this ascending form, and the fact that the Weyl vector
is also lightlike makes a lot of magic happen in 26 dimensions. For
example, it turns out that in 26 dimensions there are \emph{infinitely
many} fundamental roots, unlike in the two lower dimensional cases.

Just to add mystery upon mystery, Conway notes that in higher dimensions
there is no vector \(v\) for which all the fundamental roots \(r\) have
\(r\cdot v\) equal to some constant. So the pattern above does not
continue.

I find this stuff fascinating, but it would drive me nuts to try to work
on it. It's as if God had a day off and was seeing how many strange
features he could build into mathematics without actually making it
inconsistent.

\begin{center}\rule{0.5\linewidth}{0.5pt}\end{center}

\textbf{Yet another addendum (August 2001):} now, with the rise of
interest in \(11\)-dimensional physics, there is even a paper on
\(\mathrm{E}_{11}\):

\begin{enumerate}
\def\labelenumi{\arabic{enumi})}
\setcounter{enumi}{11}
\tightlist
\item
  P. West, \(\mathrm{E}_{11}\) and M-theory, \emph{Class. Quant. Grav.}
  \textbf{18} (2001), 4443--4460. Also available as
  \href{https://arxiv.org/abs/hep-th/0104081}{\texttt{hep-th/0104081}}.
\end{enumerate}

\hypertarget{week96}{%
\section{December 16, 1996}\label{week96}}

Lots of cool papers have been appearing which I've been neglecting in my
attempts to write expository stuff about string theory, lattices,
category theory, and all that. It's time to start catching up!

Let me start with the following book:

\begin{enumerate}
\def\labelenumi{\arabic{enumi})}
\tightlist
\item
  J. Scott Carter, Daniel E. Flath and Masahico Saito, \emph{The
  Classical and Quantum \(6j\)-Symbols}, Princeton U.\ Press,
  Princeton, 1995.  
\end{enumerate}

Ever since Jones discovered the Jones polynomial invariant of knots, an
amazing story has been unfolding about the relation between algebra and
3-dimensional topology. Some key players in this story are the ``quantum
groups'': certain noncommutative algebras analogous to the commutative
algebras of functions on groups. In fact, not merely are they analogous,
they depend on a parameter, usually called Planck's constant or
\(\hbar\), and in the classical limit where \(\hbar\to0\) they actually
reduce to algebras of functions on familiar groups. The simplest case is
``quantum \(\mathrm{SU}(2)\)'', which reduces in the classical limit to
the group \(\mathrm{SU}(2)\) of \(2\times2\) unitary matrices with
determinant \(1\). Ironically, it's good old ``classical
\(\mathrm{SU}(2)\)'' that governs the quantum mechanical theory of
angular momentum. Quantum \(\mathrm{SU}(2)\) was first discovered by
people working on physics in 2-dimensional spacetime, where when you
quantize certain systems you also need to quantize their group of
symmetries!

Nowadays, mathematicians find it simpler to work with the closely
related ``quantum \(\mathrm{SL}(2)\)'', a quantization of the the group
\(\mathrm{SL}(2)\) of all \(2\times2\) complex matrices with determinant
\(1\). The above book is largely about quantum \(\mathrm{SL}(2)\) and
its applications to topology.

All quantum groups give rise to invariants of knots, links, and tangles.
They also give rise to \(3\)-dimensional topological quantum field
theories of ``Turaev--Viro type''. This is a kind of quantum field theory
you can define on a \(3\)-dimensional spacetime that you've
triangulated, i.e., chopped up into tetrahedra. One of the main things
you want to compute in a quantum field theory is the ``partition
function'', and we say the Turaev--Viro theories are ``topological''
because you get the same answer for the partition function no matter how
you triangulate the 3-dimensional manifold corresponding to your
spacetime: the partition function only depends on the topology of the
manifold. The \(\mathrm{SU}(2)\) Turaev--Viro theory, the first one to be
discovered, is also one of the most interesting because, modulo a few
subtle points, this theory is just quantum gravity in 3 dimensions (see
\protect\hyperlink{week16}{``Week 16''}). The basic idea, though, is
that you compute the partition function by summing over all ways of
labelling the edges of your tetrahedra by ``spins''
\(j = 0, 1/2, 1, 3/2,\ldots\). Ponzano and Regge had tried to set up
\(3\)-dimensional quantum gravity this way previously, but there were
problems getting the sum to converge. The neat thing about the quantum
group is that you only sum over spins less than some fixed spin
depending on the value of \(\hbar\). Since the sums are finite, they
automatically converge.

It turns out that in these Turaev--Viro theories you are not actually
taking advantage of all the structure of the quantum group. Using the
extra structure, you can also use quantum groups to define certain
\emph{4-dimensional} topological quantum field theories, those of
``Crane--Yetter--Broda'' type. Here you triangulate a \(4\)-dimensional
manifold and, in the \(\mathrm{SU}(2)\) case, you label both the 2d
faces the 3d tetrahedra with spins. Actually, lots of people think the
Crane--Yetter--Broda theories are boring, because they look sort of boring
if you only examine their implications for \(4\)-dimensional topology.
However, they become interesting when you realize that, like all
topological quantum field theories defined using triangulations, they
are ``extended topological quantum field theories''. Roughly speaking
this means that they have implications for all dimensions below the
dimension they live in.

In particular, the Crane--Yetter--Broda theories spawn \(3\)-dimensional
topological quantum field theories of
``Chern--Simons-Reshetikhin--Turaev'' type, and most people agree that
\emph{these} are interesting. I like to emphasize, however, that a deep
understanding of these \(3\)-dimensional progeny requires an
understanding of their seemingly innocuous \(4\)-dimensional ancestors.
Also, there are a lot of interesting relationships between the
\(\mathrm{SU}(2)\) Crane--Yetter--Broda model and quantum gravity in 4
dimensions, which we are just beginning to understand. See
\protect\hyperlink{week56}{``Week 56''} for a bit about this.

If you haven't yet joined the fun, Carter, Saito, and Flath's book is a
great place to start learning about the marvelous interplay between
algebra, topology, and physics in 3 and 4 dimensions. Needless to say,
it doesn't cover all the ground I've sketched above. Instead, it focuses
on a rather specific and concrete aspect: the \(6j\) symbols. This
should make it especially handy for beginners who aren't familiar with
category theory, path integrals, and all that jazz.

What are the \(6j\) symbols, anyway? Here I need to get a wee bit more
technical. The ``classical'' \(6j\) symbols are important in the
representation theory of plain old classical \(\mathrm{SU}(2)\), while
the ``quantum'' ones are analogous gadgets applicable to quantum
\(\mathrm{SU}(2)\). In either case the idea is the same.
\(\mathrm{SU}(2)\), classical or quantum, has different representations
corresponding to different spins \(j = 0, 1/2, 1, 3/2,\ldots\). (If you
don't know what I mean by this, try \protect\hyperlink{week5}{``Week
5''}.) If we take three representations \(j_1\), \(j_2\), and \(j_3\),
we can tensor them either like this: \[(j_1\otimes j_2)\otimes j_3\] or
like this \[j_1\otimes (j_2\otimes j_3)\] The tensor product is
associative, but that doesn't mean that the above two representations
are \emph{equal}. They are only \emph{isomorphic}. This
\emph{isomorphism} can be thought of as just a big fat matrix, and the
entries in this matrix are a bunch of numbers, the \(6j\) symbols.

Turaev and Viro used the quantum \(6j\) symbols to define the original
Turaev--Viro model. It goes like this: first you chop your
\(3\)-dimensional manifold up into tetrahedra, and then you consider all
possible ways of labelling the edges with spins. Each tetrahedron gets
labelled with 6 spins since it has 6 edges, and from these spins we can
compute a number: the \(6j\) symbol. Then we multiply all these
together, one for each tetrahedron, and finally we sum over labellings
to get the partition function. Marvelously, the identities satisfied by
the \(6j\) symbols are precisely what's needed to make the result
independent of the triangulation! See \protect\hyperlink{week38}{``Week
38''} for an explanation of this seeming miracle: it's actually no
miracle at all.

\begin{enumerate}
\def\labelenumi{\arabic{enumi})}
\setcounter{enumi}{1}
\tightlist
\item
  E. Guadagnini, L. Pilo, ``Three-manifold invariants and their relation
  with the fundamental group'', 22 pages in LaTeX available as
  \href{https://arxiv.org/ps/hep-th/9612090}{\texttt{hep-th/9612090}}.
\end{enumerate}

Fans of topological field theory may like this one, though I must admit
I haven't gotten around to doing more than reading the abstract yet. In
this paper the authors give evidence for the conjecture that among
3-manifolds \(M\) for which the Chern--Simons invariant
\(\mathrm{CS}(M)\) is nonzero, the absolute value \(|\mathrm{CS}(M)|\)
only depends on the fundamental group of \(M\). Chern--Simons theory
depends on a choice of group; they prove the conjecture for certain
manifolds (``lens spaces'') when the group is \(\mathrm{SU}(2)\), and
give numerical evidence when the gauge group is \(\mathrm{SU}(3)\).

What's interesting about this to me is that \(|\mathrm{CS}(M)|^2\) is
just the Turaev--Viro theory partition function, so this conjecture is
saying that the Turaev--Viro theories discussed above have a tendency to
notice only the fundamental group.

\begin{enumerate}
\def\labelenumi{\arabic{enumi})}
\setcounter{enumi}{2}
\tightlist
\item
  Michael Reisenberger and Carlo Rovelli, ``\,`Sum over surfaces' form
  of loop quantum gravity'', available as
  \href{https://arxiv.org/ps/gr-qc/9612035}{\texttt{gr-qc/9612035}}.
\end{enumerate}
\noindent
This wonderful paper should really push forwards our understanding of
the loop representation of quantum gravity. I talked a little bit about
the basic idea in \protect\hyperlink{week86}{``Week 86''}. In the loop
representation, a state of quantum gravity at a given moment is
represented by a bunch of knotted loops or ``spin networks'' in space.
What's the spacetime picture? Well, if you have a surface in spacetime
and look at it at one moment of time, it typically looks like a bunch of
loops... so maybe the spacetime picture of quantum gravity is that
spacetime is packed with \(2\)-dimensional surfaces, all tangled up.
Interestingly, this is also very reminiscent of the picture of quantum
gravity in string theory!

I've been working on this sort of idea ever since I wrote a paper
suggesting that the loop representation and string theory might be two
faces of the same ideas (see \protect\hyperlink{week18}{``Week 18''}).
Since then, most of the time I've been trying to understand how these
ideas relate to the Crane--Yetter--Broda theories, and trying to set up
the necessary \emph{algebra} (\(n\)-category theory) to deal nicely with
surfaces in 4-dimensional spacetime.

But there are many other angles from which one can attack this problem,
and one of the best is to start directly from Einstein's equations for
general relativity, try to quantize them using the path-integral
approach, and see how the path integral can be written as a sum over
surfaces. Reisenberger has already begun work on this in the context of
``simplicial quantum gravity'' --- where you chop spacetime up into the
4-dimensional analog of tetrahedra. But during the Vienna workshop on
canonical quantum gravity this summer, we talked about a different,
still more direct approach (see \protect\hyperlink{week89}{``Week
89''}). The idea is to copy standard quantum field theory, write the
propagator describing time evolution as a time-ordered exponential, and
interpret the terms in the resulting sum as surfaces in spacetime. It's
all very analogous to traditional Feynman diagrams, where you write the
propagator as a sum over diagrams, but now the ``Feynman diagrams'' are
\(2\)-dimensional surfaces. (Again, this is reminiscent of string theory
--- but with many important differences.)

There is much more to say, but I think I'll leave it at that\ldots. Over
in the world of \(n\)-categories there is also some very interesting stuff
happening, which I will discuss more next week. I'm almost done writing
a paper with James Dolan on the definition of \(n\)-categories, but in
the meantime some other folks have been coming up with other definitions
of \(n\)-categories, so we will soon be in the position to compare
definitions and see how similar or different they are, and start
erecting the formalism needed to deal with all these topological quantum
field theories and ``sums over surfaces'' in a really elegant way!
Everything looks like its fitting together. At least, that's my
momentary optimistic feeling. Perhaps it's just the fact that classes
are over that is making me so happy. Yes, it's probably just that.

\hypertarget{week97}{%
\section{February 8, 1997}\label{week97}}

I've taken a break from writing This Week's Finds in order to finish up
a paper with James Dolan in which we give a definition of ``weak
\(n\)-categories'' for all \(n\). This paper is now available on my
website, and I'm immodestly eager to talk about it, and I will, but a
lot of stuff has accumulated in the meantime which I want to discuss
first.

First, I'm sure you remember a while back when atoms were first coaxed
to form true Bose--Einstein condensates. The basic idea is that particles
come in two basic kinds, fermions and bosons, and while the fermions
have half-integer spin and obey the Pauli exclusion principle saying
that no two identical fermions can be in the same state at the same
time, bosons have integer spin and are gregarious: they \emph{love} to
be in the same state at the same time.

Why is spin related to what happens when you try to put a bunch of
particles in the same state? Well, it all has to do with the relation
between twisting something around: \[
  \begin{tikzpicture}
    \begin{knot}[clip width=5]
      \strand[thick] (0,0)
        to (0,-0.5)
        to [out=down,in=left] (0.3,-1)
        to [out=right,in=right,looseness=2] (0.3,-0.5);
      \strand[thick] (0.3,-0.5)
        to [out=left,in=up] (0,-1)
        to (0,-1.5);
    \flipcrossings{1}
    \end{knot}
  \end{tikzpicture}
\] and switching two things: \[
  \begin{tikzpicture}
    \begin{knot}[clip width=7]
      \strand[thick] (1,0) to (0,-2);
      \strand[thick] (0,0) to (1,-2);
    \end{knot}
  \end{tikzpicture}
\] To understand this, try

\begin{enumerate}
\def\labelenumi{\arabic{enumi})}
\tightlist
\item
  John Baez, ``Spin, statistics, CPT and all that jazz'',
  \href{http://math.ucr.edu/home/baez/spin.stat.html}{\texttt{http://math.ucr.edu/home/}}  \href{http://math.ucr.edu/home/baez/spin.stat.html}{\texttt{baez/spin.stat.html}}
\end{enumerate}

But let's consider some examples. Since photons have spin 1 they are
bosons. In laser light one has a bunch of photons all in the same state.
Thanks to the Heisenberg uncertainty principle, of course, we can't know
both their position and momentum. In a laser we don't know the position
of the photons: each photon is all over the laser beam in a spread-out
sort of way. However, we do know the momentum of the photons and they
all have the same momentum. This means that we have ``coherent light''
in which all the photons are like waves wiggling perfectly in phase. One
can demonstrate this by interfering two beams of laser light and seeing
beautifully perfect interference fringes, bright and dark stripes in
places where the two beams are either in phase with each other and
adding up, or out of phase and cancelling out.

Now, other particles are bosons as well, and they can do similar tricks.
Bose and Einstein predicted that when any gas of noninteracting bosons
gets sufficiently cold, all --- or at least a sizeable fraction --- of
them will be found in the same state: the state of least possible
energy. Unfortunately, when things get cold they are usually liquids or
solids rather than a gas, and the particles in a liquid or gas interact
a lot, so true Bose--Einstein condensation is hard to achieve.

Some related things have been studied for decades. If you get an even
number of fermions together they act approximately like a boson, at
least if the density is not too high. Helium stays liquid at
temperatures arbitrarily close to absolute zero, when the pressure is
low enough. Since helium 4 has 2 protons, 2 neutrons, and 2 electrons,
and all these particles are fermions, helium 4 acts like a boson. At
really low temperatures, helium 4 becomes ``superfluid'' --- a
substantial fraction of the atoms fall into the same state and the
liquid acquires shocking properties, like zero viscosity. Similarly, in
certain metals at low temperatures electrons will, by a subtle
mechanism, form ``Cooper pairs'', and these act like bosons. When a
bunch of these fall into the same state, you have a ``superconductor''.

But neither of these is a Bose--Einstein condensate in the technical
sense of the term, because the helium atoms interact a lot in superfluid
helium, and the Cooper pairs interact a lot in a superconductor. Only
recently have people been able to get dilute gases of bosonic atoms cold
enough to study true Bose--Einstein condensation.

The first team to do it, the ``JILA'' team in Boulder, Colorado got a
Bose--Einstein condensate of about 2000 rubidium atoms to form in a
magnetic trap at less than \(2 \times 10^{-7}\) degrees above absolute
zero. A team at Rice University did it with lithium soon after, followed
by a team at MIT, who did it with sodium.

Check out:

\begin{enumerate}
\def\labelenumi{\arabic{enumi})}
\setcounter{enumi}{1}
\item
   NIST, ``Physicists create new state of matter'', July 13, 1995.  Available at
  \href{https://www.nist.gov/news-events/news/1995/07/physicists-create-new-state-matter-record-low-temperature}{\texttt{https://www.nist.gov/news-events/news/1995/07/physicists-create-new-state-}}
 \hfill \break \href{https://www.nist.gov/news-events/news/1995/07/physicists-create-new-state-matter-record-low-temperature}{\texttt{matter-record-low-temperature}}

  Atomcool home page, \href{http://atomcool.rice.edu/}{\texttt{http://atomcool.rice.edu/}}

  ``Bose--Einstein condensation of ultracold sodium''
  \href{https://web.archive.org/web/19990203201237/http://amo.mit.edu/\%7Ebec/news/news.html}{\texttt{https://web.archive.org/web/}} \href{https://web.archive.org/web/19990203201237/http://amo.mit.edu/\%7Ebec/news/news.html}{\texttt{19990203201237/http://amo.mit.edu/\%7Ebec/news/news.html}}
\end{enumerate}

So what's the news? Well, recently the team at MIT, led by Wolfgang
Ketterle, made two blobs of Bose--Einstein condensate out of sodium
atoms. Ramming these into each other, they were able to see interference
fringes just as in a laser! In other words, they are seeing interference
of matter waves, just as quantum mechanics predicts, but involving lots
of atoms in a coherent state rather than a single electron as in the
famous double slit experiment. 

In honor of this event, I hereby present the following limerick composed
by the poet Lisa Raphals, with myself serving as science consultant. It
may aid your appreciation if I note first that ``Squantum'' is an actual
town in Massachusetts. With no further ado:

\begin{quote}
A metaphysician from Squantum \\
Was asked, what's the state of the
quantum? \\ It's all reciprocity: Position, velocity --- \\
 They're never both
there when you want 'em!
\end{quote}
\noindent
Now on to some more technical stuff\ldots.

I am now visiting the Center for Gravitational Physics and Geometry here
at Penn State, which is a delightful place for people interested in the
loop representation of quantum gravity (see
\protect\hyperlink{week77}{``Week 77''}). Right now everyone is working
on using the loop representation to derive Hawking's formula which says
that the entropy of a black hole is proportional to the surface area of
its event horizon.

When I arrived, Jorge Pullin handed me a copy of his book:

\begin{enumerate}
\def\labelenumi{\arabic{enumi})}
\setcounter{enumi}{2}
\tightlist
\item
  Rodolfo Gambini and Jorge Pullin, \emph{Loops, Knots, Gauge Theories,
  and Quantum Gravity}, Cambridge U. Press, Cambridge, 1996.
\end{enumerate}

Here is the table of contents:

\begin{enumerate}
\def\labelenumi{\arabic{enumi}.}
\tightlist
\item
  Holonomies and the group of loops
\item
  Loop coordinates and the extended group of loops
\item
  The loop representation
\item
  Maxwell theory
\item
  Yang--Mills theories
\item
  Lattice techniques
\item
  Quantum gravity
\item
  The loop representation of quantum gravity
\item
  Loop representation: further developments
\item
  Knot theory and physical states of quantum gravity
\item
  The extended loop representation of quantum gravity
\item
  Conclusions, present status and outlook
\end{enumerate}

This is presently the most complete introduction available to the ``loop
representation'' concept, as applied to electromagnetism, Yang--Mills
theory, and quantum gravity. Gambini was one of the original inventors
of this notion, and this book covers the whole sweep of its
ramifications, with a special emphasis on a particular technical form,
the ``extended loop representation'', which he has been developing with
Pullin and other collaborators.

What the heck is the loop representation, anyway? Well, all the forces
we know are described by gauge theories, and gauge theories all describe
the ``phase'', or generalization thereof, that a particle acquires when
you carry it around a loop. In the case of electromagnetism, for
example, a charged quantum particle carried around a loop in space
acquires a phase equal to \[\exp(-iqB/\hbar)\] where \(q\) is the
particle's charge, \(\hbar\) is Planck's constant, and \(B\) is the
magnetic flux through the loop: i.e., the integral of the magnetic field
over any surface spanning the loop. Knowing these phases for all loops is
the same as knowing the magnetic field. Similarly, if we knew the phase
for all loops in \emph{spacetime} instead of just space, we would know both the
electric and magnetic fields throughout spacetime.

General relativity is similar except that instead of a phase one gets a
rotation, or more generally a Lorentz transformation, when one parallel
transports a little arrow around a loop.

The theories of the electroweak and strong forces are similar but the
analog of the ``phase'' is a bit more abstract: an element of the group
\(\mathrm{SU}(2)\times\mathrm{U}(1)\) or \(\mathrm{SU}(3)\),
respectively.

The idea of the loop representation is to take these ``phases acquired
around loops'' as basic variables for describing the laws of physics.

That's the idea in a nutshell. It turns out, not surprisingly, that
there are many interesting relationships with such topics involving
loops, such as string theory and knot theory.

Gambini and Pullin's book develops this theme in many directions. Let me
say a bit about one fascinating topic that they mention, which I would
like to understand better: Gerard 't Hooft's work on confinement in
chromodynamics using his ``order and disorder operators''.

I explained some basic ideas about confinement and asymptotic freedom in
\protect\hyperlink{week84}{``Week 84''} and \protect\hyperlink{week94}{``Week 94''},
so I'll assume you've read that stuff. Remember, the basic idea of
confinement is that if you take a meson and try to pull the quark and
antiquark it contains apart, the force required does not decrease with
distance like \(1/r^2\), because the chromoelectric field --- the strong
force analog of the electric field --- does not spread out in all
directions like an ordinary electric field does. Instead, all the field
lines are confined to a ``flux tube'', so the force is roughly
independent of the distance.

This means that the energy is roughly proportional to the distance.
Since action has dimensions of energy times time, this means that if we
consider the creation and subsequent annihilation of a virtual
quark-antiquark pair: \[
  \begin{tikzpicture}
    \begin{knot}
      \strand[thick] (0,0)
        to [out=up,in=up,looseness=2] (1.5,0);
      \strand[thick] (0,0)
        to [out=down,in=down,looseness=2] (1.5,0);
    \end{knot}
    \node[fill=white] at (0.2,0.5) {q};
    \node[fill=white] at (1.3,0.5) {$\overline{\mbox{q}}$};
  \end{tikzpicture}
\] the total action is proportional to the \emph{area} of the loop
traced out in spacetime. Here I am neglecting the action due to the
kinetic energy of the quark and antiquark, and only worrying about the
potential energy due to the flux tube joining them. This amounts to
treating the quark and antiquark as ``test particles'' to study the
behavior of the strong force.

Now, when we study quantum physics using Euclidean path integrals the
basic principle is that the probability of the occurence of any process
is proportional to \[\exp(-S/\hbar)\] where \(S\) is the action of that
process and \(\hbar\) is Planck's constant again. So in this framework
the \emph{probability} of a particular virtual quark-antiquark pair
tracing out a loop like the above one is proportional to \[\exp(-cA)\]
where \(c > 0\) is some constant and \(A\) is the area of the loop. This
``area law'' was first proposed by Kenneth Wilson in his pioneering work
on confinement; he proposed it as a way to tell, mathematically, if
confinement was happening in some theory. Just compute the probability
of a virtual quark-antiquark pair tracing out a particular loop and see
if it decreases exponentially with the area!

Deriving confinement from chromodynamics is something that people have
worked on for quite a while, and it's not easy: there is still no
rigorous proof, even though there are a bunch of heuristic arguments for
it, and computer simulations seem to demonstrate that it's bound to
occur. One approach to studying the puzzle is due to 't Hooft and
involves ``order'' and ``disorder'' operators.

I'll explain what these are, and what they have to do with knot theory,
but not how 't Hooft actually uses them in his argument for confinement.
For the actual argument, try Gambini and Pullin's book, or else 't
Hooft's paper:

\begin{enumerate}
\def\labelenumi{\arabic{enumi})}
\setcounter{enumi}{4}
\tightlist
\item
  Gerard 't Hooft, On the phase transition towards permanent quark confinement, \emph{Nucl. Phys. B} \textbf{138} (1978) 1--25.
\end{enumerate}

Let us work in space at a given time, rather than in the Euclidean path
integral approach. We'll do ``canonical quantization'', meaning that now
observables will be operators on some Hilbert space.

If we have any loop \(g\) in space, there is an observable called the
``Wilson loop'' \(W(g)\), which is the trace of the holonomy of the
connection around \(g\). The precise way of stating Wilson's area law
for confinement in this context is that
\[\langle W(g) \rangle \sim \exp(-cA)\] where \(\langle W(g) \rangle\)
is the vacuum expectation value of the Wilson loop, and \(A\) is the
area spanned by the loop \(g\). The point is that
\(\langle W(g) \rangle\) is the same as what I was (a bit sloppily)
calling the probability of the quark-antiquark pair tracing out the loop
\(g\).

't Hooft calls the Wilson loops ``order operators''. We don't really
need to worry why he calls them this, but if you know how physicists
think, you may know that the Wilson loops are keeping track of a kind of
``order parameter'' of the vacuum state. Anyway, his idea was to study
the Wilson loops by introducing some other operators he called
``disorder operators''.

Chromodynamics is an \(\mathrm{SU}(3)\) gauge theory but it's a little
clearer if we work with any \(\mathrm{SU}(N)\) gauge theory. Notice that
the center of the group \(\mathrm{SU}(N)\) consists of the matrices of
the form \[\exp(2\pi in/N)\]
where \(n\) is an integer. So if we have a loop \(h\), we can imagine an
operator that does the following thing: it modifies the connection, or
vector potential, in such a way that if you do parallel transport around
a tiny loop linking \(h\) once, the holonomy changes to
\(\exp(2\pi i/N)\) times what it had been. Note: this is a
gauge-invariant thing to do, because that \(\exp(2\pi i/N)\) is in the
center of \(\mathrm{SU}(N)\)! So just as the Wilson loop observables are
gauge-invariant, we can hope for some some ``disorder operators''
\(V(h)\) that modify the connection in this way.

If you think about it, what this means is that the following commutation
relations hold: \[W(g) V(h) = V(h) W(g) \exp(2\pi i L(g,h)/N)\] where
\(L(g,h)\) is the linking number of the loops \(g\) and \(h\), which
counts how many times \(g\) wraps around \(h\).

There is an obvious symmetry or ``duality'' between the \(V\)'s and the
\(W\)'s going on here. In fact, just as \(W\)'s satisfy an area law
where there is confinement of chromoelectric field lines into flux
tubes, I believe the \(V\)'s satisfy an area law when there is
confinement of chromomagnetic field lines into flux tubes. The simplest
case of this kind of thing occurs in plain old electromagnetism, where
plain old magnetic field lines are confined into flux tubes in type II
superconductors. For this reason confinement of electric field lines is
sometimes called ``dual superconductivity''.

Perhaps the simplest way of beginning to understand this stuff more
deeply is to understand the wonderful formula proved by Ashtekar and
Corichi in the following paper:

\begin{enumerate}
\def\labelenumi{\arabic{enumi})}
\setcounter{enumi}{5}
\tightlist
\item
  Abhay Ashtekar and Alejandro Corichi, ``Gauss linking number and
  electro-magnetic uncertainty principle'', \emph{Phys. Rev. D} \textbf{56}
  (1997), 2073.  Also
  available as
  \href{https://arxiv.org/ps/hep-th/9701136}{\texttt{hep-th/9701136}}.
\end{enumerate}
\noindent
This formula applies to plain old electromagnetism, or more precisely,
quantum electrodynamics. If we work in units where \(\hbar = 1\), and
consider a particle of charge \(1\), the Wilson loop operator \(W(g)\)
in electromagnetism is just \[W(g) = \exp(-iB(g))\] where \(B\) is the
magnetic flux flowing through the loop \(g\). But instead we can just
work with \(B(g)\) directly. Similarly, instead of \(V(h)\)'s we can
work with the operator \(E(h)\) corresponding to the electric flux
through the loop \(h\). Then we have
\[B(g) E(h) - E(h) B(g) = i L(g,h).\] In other words, the electric and
magnetic fields don't commute in quantum electrodynamics, and the
Heisenberg uncertainty of the electric field flowing through a loop \(g\)
and the magnetic field flowing through a loop \(h\) is proportional to
the linking number of \(g\) and \(h\)!

Quantum mechanics, electromagnetism, and knot theory are clearly quite
tangled up here. Since the linking number was first discovered by Gauss
in his work on magnetism, it's all quite fitting.

And that leads me to the last paper I want to mention this week. It
should be of great interest to Vassiliev invariant fans; see
\protect\hyperlink{week3}{``Week 3''} if you don't know what a Vassiliev
invariant is.

\begin{enumerate}
\def\labelenumi{\arabic{enumi})}
\setcounter{enumi}{6}
\tightlist
\item
  Dror Bar-Natan and Alexander Stoimenow, ``The fundamental theorem of
  Vassiliev invariants'', in \emph{Geometry and Physics}, CRC Press,  Boca Raton,
  2021,  pp\ 101--134.  Also available as
  \href{https://arxiv.org/ps/q-alg/9702009}{\texttt{q-alg/9702009}}.
\end{enumerate}
\noindent
Let me just quote the abstract here:

\begin{quote}
The ``fundamental theorem of Vassiliev invariants'' says that every
weight system can be integrated to a knot invariant. We discuss four
different approaches to the proof of this theorem: a
topological/combinatorial approach following M. Hutchings, a geometrical
approach following Kontsevich, an algebraic approach following
Drinfel'd's theory of associators, and a physical approach coming from
the Chern--Simons quantum field theory. Each of these approaches is
unsatisfactory in one way or another, and hence we argue that we still
don't really understand the fundamental theorem of Vassiliev invariants.
\end{quote}

\hypertarget{week98}{%
\section{February 27, 1997}\label{week98}}

I feel guilty for slacking off on This Week's Finds, so I should explain
the reason. Lots of things have been building up that I'm dying to talk
about, but new ones keep coming in at such a rapid rate that I never
feel I have time!

I will have to be ruthless, and face up to the fact that a quick and
imperfect exposition is better than none.

First of all, here at the Center for Gravitational Physics and Geometry
there are a lot of interesting attempts going on to compute the entropy
of black holes from first principles. Bekenstein, Hawking and many
others have used a wide variety of semiclassical arguments to argue that
black holes satisfy \[S = A/4\] where \(S\) is the entropy and \(A\) is
the area of the event horizon, both measured in Planck's units, where
\(G = c = \hbar = 1\).

For example, using purely classical reasoning (general relativity, but
no quantum theory) one can prove the ``2nd law of black hole
thermodynamics'', which says that \(A\) always increases. As Bekenstein
noted, this suggests that the area of the event horizon is somehow
analogous to entropy. However, by itself this does not determine the
magic number \(1/4\), which can only be derived using quantum theory (as
one can see by simple dimensional analysis).

By semiclassical reasoning --- studying quantum electrodynamics in the
Schwarzschild metric used to describe black holes --- Hawking showed
that black holes should radiate as if they had a temperature inversely
proportional to their mass: \[T = \frac{1}{8\pi M}.\] This made the
analogy between entropy and event horizon area much more than an
analogy, because it meant that one could assign a temperature to black
holes and see if they satisfy the laws of thermodynamics. It turns out
that if you consider \(A/4\) to be the entropy of a black hole, you can
eliminate seeming violations of the 2nd law that otherwise arise in
thought experiments where you get rid of entropy by throwing it into a
black hole. In other words, if you throw something with entropy \(S\)
into a black hole, calculations seem to show that the area of the event
horizon always increases by at least \(4S\)!

So far nothing I've said is related to full-fledged quantum gravity,
because in the semiclassical arguments one is still working in the
approximation where the gravitational field is treated classically. An
interesting test of any theory of quantum gravity is whether can use it
to derive \(S = A/4\). In a subject with no real experimental evidence,
this is the closest we have to an ``experimental result'' that our
theory should predict.

Recently the string theorists have done some calculations claiming to
show that string theory predicts \(S = A/4\). Personally I feel that
while these calculations are interesting they are far from definitive.
For example, they all involve taking calculations done using
perturbative string theory on \emph{flat} spacetime and extrapolating
them drastically to the regime in which string theory approximates
general relativity. One typically uses ideas from supersymmetry to
justify such extrapolations; however, these ideas only seem to apply to
``extremal black holes'', having the maximum possible charge for a black
hole of a given mass and angular momentum. Realistic black holes are far
from extremal. In short, while exciting, these calculations still need
to be taken with a grain of salt.

Of course, I am biased because I am interested in another approach to
quantum gravity, the loop representation of quantum gravity, which folks
are working on here at the CGPG, among other places. This is in many
ways a more conservative approach. The idea is to simply take Einstein's
equation for general relativity and quantize it, rather than trying to
develop a theory of \emph{all} particles and forces as in string theory.
Of course, for various reasons it is not so easy to quantize Einstein's
equation. String theorists think it's \emph{impossible} without dragging
in all sorts of other forces and particles, but folks working on the
loop representation are more optimistic. This is an ongoing argument,
but certainly a good test of the loop representation would be to try to
use it to derive Hawking's formula S = A/4. If the loop representation
is really any good, this should be possible, because many different
lines of reasoning using only general relativity and quantum theory lead
to this formula.

I've already mentioned a few attempts to do this in
\protect\hyperlink{week56}{"Week 56}``, \protect\hyperlink{week57}{"Week
57}'', and \protect\hyperlink{week87}{``Week 87''}. These were
promising, but they didn't get the magical number \(1/4\). Also, they
are rather rough, in that they do computations on some region with
boundary, but don't use anything that ensures the boundary is an event
horizon.

Recently Kirill Krasnov has made some progress:

\begin{enumerate}
\def\labelenumi{\arabic{enumi})}
\tightlist
\item
  Kirill Krasnov, ``On quantum statistical mechanics of Schwarzschild black
  hole'', \emph{Gen. Rel. Grav.} \textbf{30} (1998), 53--68.  Also available as
  \href{https://arxiv.org/ps/gr-qc/9605047}{\texttt{gr-qc/9605047}}.
\end{enumerate}

This paper still doesn't get the magic number \(1/4\), and Krasnov later
realized it has a few mistakes in it, but it does something very cool.
It notes that the boundary conditions holding on the event horizon of a
Schwarzschild black hole are closely related to Chern--Simons theory. Now
is not the time for me to go into Chern--Simons theory, but basically, it
lets you apply a lot of neat mathematics to calculate everything to your
heart's content, very much as Carlip did on his work on the toy model of
a 2+1-dimensional black hole (see \protect\hyperlink{week41}{``Week
41''}). Also, it sheds new light on the relationship between topological
quantum field theory and quantum gravity, something I am always trying
to understand better.

While I'm at it, I should note the existence of a paper that reworks
Carlip's calculation from a slightly different angle:

\begin{enumerate}
\def\labelenumi{\arabic{enumi})}
\setcounter{enumi}{1}
\tightlist
\item
  Maximo Banados and Andres Gomberoff, ``Black hole entropy in the
  Chern--Simons formulation of 2+1 gravity'', \emph{Phys. Rev. D} \textbf{55} (1997)     
  6162--6167.  Also available as
  \href{https://arxiv.org/ps/gr-qc/9611044}{\texttt{gr-qc/9611044}}.
\end{enumerate}
\noindent
2+1-dimensional quantum gravity is very simple compared to the
3+1-dimensional quantum gravity we'd really like to understand: in a
sense it's ``exactly solvable''. But there are still some puzzling
things about Carlip's computation of the entropy of a black hole in 2+1
dimensions which need figuring out, so every paper on the subject is
worth looking at, if you're interested in black hole entropy.

Speaking of topological quantum field theory and quantum gravity, I just
finished a paper on these topics:

\begin{enumerate}
\def\labelenumi{\arabic{enumi})}
\setcounter{enumi}{2}
\tightlist
\item
  John Baez, ``Degenerate solutions of general relativity from
  topological field theory'', \emph{Commun. Math. Phys.}
   \textbf{193} (1998) 219--231.   Also available as 
  \href{https://arxiv.org/abs/gr-qc/9702051}{\texttt{gr-qc/9702051}}.
\end{enumerate}

Let me just summarize the basic idea, resisting the temptation to become
insanely technical.

A while ago Rovelli and Smolin introduced Penrose's notion of ``spin
network'' into the loop representation of quantum gravity. I described
spin networks pretty carefully in \protect\hyperlink{week43}{``Week
43''}, but here let me just say that they are graphs embedded in space
with edges labelled by spins \(j = 0, 1/2, 1, 3/2,\ldots\), just as in
the quantum mechanics of angular momentum, and with vertices labelled by
``intertwining operators'', which are other gadgets that come up in the
study of angular momentum. In the loop representation these spin
networks form a basis of states. Geometrical observables like the area
of surfaces and the volumes of regions have been quantized and their
matrix elements computed in the spin network basis, giving us a nice
picture of ``quantum 3-geometries'', that is, the possible geometries of
space in the context of quantum gravity. In this picture, the edges of
spin networks play the role of quantized flux tubes of area, much as the
magnetic field comes in quantized flux tubes in a type II
superconductor. To work out the area of a surface in some spin network
state, you just total up contributions from each edge of the spin
network that pokes through the surface. An edge labelled with spin \(j\)
carries an area equal to \(\sqrt{j(j+1)}\) times the Planck length
squared. What's cool is that this is not merely postulated, it's derived
from fairly standard ideas about how you turn observables into operators
in quantum mechanics.

However, the dynamics of quantum gravity is more obscure. Technical
issues aside, the main problem is that while we have a nice picture of
quantum 3-geometries, we don't have a similar picture of the
\emph{4-dimensional}, or \emph{spacetime}, aspects of the theory. To
represent a physical state of quantum gravity, a spin network state (or
linear combination thereof) has to satisfy something called the
Wheeler--DeWitt equation. This is sort of the quantum gravity analog of
the Schrodinger equation. There is a lot of controversy over the
Wheeler--DeWitt equation and what's the right way to write it down in the
loop representation. The really annoying thing, however, is that even if
you feel you know how to write it down --- for example, Thomas Thiemann
has worked out one way (see \protect\hyperlink{week85}{``Week 85''}) ---
and can find solutions, you still don't necessarily have a good
intuition as to what the solutions \emph{mean}. For example, almost
everyone seems to agree that spin networks with no vertices should
satisfy the Wheeler--DeWitt equation. These are just knots or links with
edges are labelled by spins. We know these states are supposed to
represent ``quantum 4-geometries'' satisfying the quantized Einstein
equations. But how should we visualize these states in \(4\)-dimensional
terms?

In search of some insight into the \(4\)-dimensional interpretation of
these states, I turn to classical general relativity. In my paper, I
construct solutions of the equations of general relativity which at a
typical fixed time look like ``flux tubes of area'' reminiscent of the
loop states of quantum gravity. These are ``degenerate solutions'',
meaning that the ``3-metric'', the tensor you use to measure distances
in 3-dimensional space, is zero in lots of regions of space. Here I
should warn you that ordinary general relativity doesn't allow
degenerate metrics like this. The loop representation works with an
extension of general relativity that covers the case of degenerate
metrics; for more on this, see \protect\hyperlink{week88}{``Week 88''}.

More precisely, if you look at these ``flux tube'' solutions at a
typical time, the 3-metric vanishes outside a collection of solid tori
embedded in space, while inside any of these solid tori the metric is
degenerate in the longitudinal direction, but nondegenerate in the two
transverse directions.

Now since these are classical solutions --- no quantum theory in sight!
--- there is no problem with understanding what they do as time passes.
We can solve Einstein's equation and get a 4-metric, a metric on
spacetime. The \(4\)-dimensional picture is as follows: given any
surface \(\Sigma\) embedded in spacetime, I get solutions for which the
4-metric vanishes outside a neighborhood of \(\Sigma\). Inside this
neighborhood, the 4-metric is zero in the two directions tangent to
\(\Sigma\) but nondegenerate in the two transverse directions. In the
4-geometry defined by one of these solutions, the area of a typical
surface \(\Sigma'\) intersecting \(\Sigma\) in some isolated points is a
sum of contributions from the points where \(\Sigma\) and \(\Sigma'\)
intersect.

The solutions I study are inspired by the work of Mike Reisenberger, who
studied a solution for which the metric vanishes outside a neighborhood
of a sphere embedded in \(\mathbb{R}^4\). I consider more general
surfaces embedded in more general 4-manifolds, so I need to worry a lot
more about topological issues. Also, I allow the possibility of a
nonzero cosmological constant (this being a parameter in Einstein's
equation that determines the energy density of the vacuum). A lot of the
most interesting stuff happens for nonzero cosmological constant, and
this case actually helps one understand the case of vanishing
cosmological constant as a kind of limiting case.

It turns out that the interesting degrees of freedom of the metric
living on the surface \(\Sigma\) in spacetime are described by fields
living on this surface. In fact, these fields are solutions of a
\(2\)-dimensional topological field theory called \(BF\) theory. To
prove this, I take advantage of the relation between general relativity
and \(BF\) theory in 4 dimensions, together with the fact that \(BF\)
theory behaves in a simple manner under dimensional reduction.

Another neat thing is that to get a solution of general relativity this
way, we need to pick a ``framing'' of \(\Sigma\). Roughly speaking, this
means we need to pick a way of thickening up the surface \(\Sigma\) to a
neighborhood that looks like \(\Sigma\times D^2\), where \(D^2\) is the
\(2\)-dimensional disc. This is precisely the \(4\)-dimensional analog
of a framing of a knot or link in 3-dimensions. People who know about
topological quantum field theory know that framings are very important.
In fact, I can show that my solutions of general relativity are closely
related to Chern--Simons theory, a \(3\)-dimensional topological field
theory famous for giving invariants of framed knots and links. What's
beginning to emerge is a picture that makes the \emph{spacetime} aspects
of framings easier to understand.

Now before I plunge into some even more esoteric stuff, let me briefly
return to reality and answer the question you've all been secretly dying
to ask: how does general relativity impact the world of big business?

In plain terms: is all this fancy physics just an excuse to have fun
visualizing evolving spin networks in terms of quantum field theories on
surfaces embedded in \(4\)-dimensional spacetime, etcetera
etcetera... or does it actually contribute to the well-being of the
corporations upon which we depend?

Well, you may be surprised to know that general relativity plays an
significant role in a \$200-million business. Surprised? Read on! What
follows is taken from the latest issue of ``Matters of Gravity'', the
newsletter put out by Jorge Pullin. More precisely, it's from:

\begin{enumerate}
\def\labelenumi{\arabic{enumi})}
\setcounter{enumi}{3}
\tightlist
\item
  Neil Ashby, ``General relativity in the global positioning system'',
  in \emph{Matters of Gravity}, ed.~Jorge Pullin, no. \textbf{9},
  available at \href{http://www.phys.lsu.edu//mog/mog9/node9.html}{\texttt{http://www.phys.lsu.edu//mog/mog9/}} \href{http://www.phys.lsu.edu//mog/mog9/node9.html}{\texttt{node9.html}}.
\end{enumerate}

I will simply quote some excerpts from this fascinating article:

\begin{quote}
"The Global Position System (GPS) consists of 24 earth-orbiting
satellites, each carrying accurate, stable atomic clocks. Four
satellites are in each of six different orbital planes, of inclination
55 degrees with respect to earth's equator. Orbital periods are 12 hours
(sidereal), so that the apparent position of a satellite against the
background of stars repeats in 12 hours. Clock-driven transmitters send
out synchronous time signals, tagged with the position and time of the
transmission event, so that a receiver near the earth can determine its
position and time by decoding navigation messages from four satellites
to find the transmission event coordinates, and then solving four
simultaneous one-way signal propagation equations. Conversely,
\(\gamma\)-ray detectors on the satellites could determine the
space-time coordinates of a nuclear event by measuring signal arrival
times and solving four one-way propagation delay equations.
\end{quote}

\begin{quote}
Apart possibly from high-energy accelerators, there are no other
engineering systems in existence today in which both special and general
relativity have so many applications. The system is based on the
principle of the constancy of \(c\) in a local inertial frame: the
Earth-Centered Inertial or ECI frame. Time dilation of moving clocks is
significant for clocks in the satellites as well as clocks at rest on
earth. The weak principle of equivalence finds expression in the
presence of several sources of large gravitational frequency shifts.
Also, because the earth and its satellites are in free fall,
gravitational frequency shifts arising from the tidal potentials of the
moon and sun are only a few parts in 10\^{}16 and can be neglected.
\end{quote}

\begin{quote}
{[}\ldots{]}
\end{quote}

\begin{quote}
At the time of launch of the first NTS-2 satellite (June 1977), which
contained the first Cesium clock to be placed in orbit, there were some
who doubted that relativistic effects were real. A frequency synthesizer
was built into the satellite clock system so that after launch, if in
fact the rate of the clock in its final orbit was that predicted by GR,
then the synthesizer could be turned on bringing the clock to the
coordinate rate necessary for operation. The atomic clock was first
operated for about 20 days to measure its clock rate before turning on
the synthesizer. The frequency measured during that interval was +442.5
parts in \(10^{12}\) faster than clocks on the ground; if left
uncorrected this would have resulted in timing errors of about 38,000
nanoseconds per day. The difference between predicted and measured
values of the frequency shift was only 3.97 parts in \(10^{12}\), well
within the accuracy capabilities of the orbiting clock. This then gave
about a 1\% validation of the combined motional and gravitational shifts
for a clock at 4.2 earth radii.
\end{quote}

\begin{quote}
{[}\ldots{]}
\end{quote}

\begin{quote}
This system was intended primarily for navigation by military users
having access to encrypted satellite transmissions which are not
available to civilian users. Uncertainty of position determination in
real time by using the Precise Positioning code is now about 2.4 meters.
Averaging over time and over many satellites reduces this uncertainty to
the point where some users are currently interested in modelling many
effects down to the millimeter level. Even without this impetus, the GPS
provides a rich source of examples for the applications of the concepts
of relativity.
\end{quote}

\begin{quote}
New and surprising applications of position determination and time
transfer based on GPS are continually being invented. Civilian
applications include for example, tracking elephants in Africa, studies
of crustal plate movements, surveying, mapping, exploration, salvage in
the open ocean, vehicle fleet tracking, search and rescue, power line
fault location, and synchronization of telecommunications nodes. About
60 manufacturers now produce over 350 different commercial GPS products.
Millions of receivers are being made each year; prices of receivers at
local hardware stores start in the neighborhood of \$200."
\end{quote}

Pretty cool, eh?

Okay, now for something completely different --- homotopy theory! Well,
everything I write about is actually secretly part of my grand plan to
see how everything interesting is related to everything else, but let me
not delve into how homotopy theory is related to topological quantum
field theory and thus quantum gravity. Let me simply mention the
existence of this great book:

\begin{enumerate}
\def\labelenumi{\arabic{enumi})}
\setcounter{enumi}{4}
\tightlist
\item
  \emph{Handbook of Algebraic Topology}, ed.~I. M. James, North-Holland,
  Amsterdam, 1995, 1324 pages.
\end{enumerate}
\noindent
Occasionally you come across a book that you wish you just download into
your brain; for me this is one of those books. It is probably not a good
idea to read it if you are just wanting to get started on algebraic
topology; it assumes you are pretty familiar with the basic ideas
already, and it goes into a lot of depth, mainly in hardcore homotopy
theory. A lot of it is too technical for me to appreciate, but let me
list a few chapters that I can understand and like.

\begin{itemize}
\item
  Chapter 1, ``Homotopy types'' by Hans-Joachim Baues, is a great survey
  of different models of homotopy types. Remember, we say two
  topological spaces \(X\) and \(Y\) are homotopy equivalent if there
  are continuous functions \(f\colon X\to Y\) and \(g\colon Y\to X\)
  that are inverses ``up to homotopy''. In other words, we don't require
  that \(fg\) and \(gf\) are \emph{equal} to identity functions, but
  merely that they can both be \emph{continuously deformed} to identity
  functions. So for example the circle and an annulus are homotopy
  equivalent, and we say therefore that they represent the same
  ``homotopy type''.

  The cool thing is that there turn out to be very elegant algebraic and
  combinatorial ways of describing homotopy types that don't mention
  topology at all. Perhaps the most beautiful way of all is a way that
  in a sense hasn't been fully worked out yet: namely, thinking of
  homotopy types as ``\(\omega\)-groupoids''. The idea is this. An
  ``\(\omega\)-category'' is something that has

  \begin{itemize}
  \tightlist
  \item
    objects like \(x\)
  \item
    morphisms between objects like \(f\colon x\to y\)
  \item
    \(2\)-morphisms between morphisms like \(F\colon f\to g\)
  \item
    \(3\)-morphisms between \(2\)-morphisms like \(T\colon F\to G\)
  \item
    \ldots{}
  \end{itemize}

  and so on ad infinitum. There should be some ways of composing these,
  and these should satisfy some axioms, and that of course is the tricky
  part. But the basic idea is that if you hand me a topological space
  \(X\), I can cook up an \(\omega\)-category whose

  \begin{itemize}
  \tightlist
  \item
    objects are points in \(X\)
  \item
    morphisms are paths between points in \(X\)
  \item
    \(2\)-morphisms are continuous 1-parameter families of paths in
    \(X\), i.e.~``paths of paths'' in \(X\)
  \item
    \(3\)-morphisms are ``paths of paths of paths'' in \(X\)
  \item
    \ldots{}
  \end{itemize}

  and so on. This is better than your garden-variety \(\omega\)-category
  because all the morphisms and \(2\)-morphisms and \(3\)-morphisms and
  so on have inverses, at least ``up to homotopy''. We call it an
  ``\(\omega\)-groupoid''. This \(\omega\)-groupoid keeps track of the
  homotopy type of \(X\) in a very nice way. (If this ``\(\omega\)''
  stuff is too mind-boggling, you may want to start by reading a bit
  about plain old categories and groupoids in
  \protect\hyperlink{week74}{``Week 74''}.)

  Conversely, given any \(\omega\)-groupoid there should be a nice way
  to cook up a homotopy corresponding to it. This is just the
  infinite-dimensional generalization of something I described in
  \protect\hyperlink{week75}{``Week 75''}. There, I showed how you could
  get a groupoid from a ``homotopy 1-type'' and vice versa. Here there are
  \(1\)-morphisms but no interesting \(2\)-morphisms, \(3\)-morphisms,
  and so on, because the topology of a ``homotopy 1-type'' is boring in
  dimensions greater than 1. (In case any experts are reading this, what
  I mean is that its higher homotopy groups are trivial; its higher
  homology and cohomology groups can be very interesting.)

  So we can --- and should --- think of homotopy theory as, among other
  things, the study of \(\omega\)-groupoids, and thus a very useful
  warmup to the study of \(\omega\)-categories. In my occasional series
  on This Week's Finds called ``the Tale of \(n\)-Categories'', I have
  tried to explain why \(n\)-categories, and ultimately
  \(\omega\)-categories, should serve as a powerful unifying approach to
  lots of mathematics and physics. In trying to understand this subject,
  I find time and time again that homotopy theorists are the ones to
  listen to.
\item
  Chapter 2, ``Homotopy theories and model categories'', by W. G. Dwyer
  and J. Spalinski, is a nice introduction to the formal idea of using
  different ``models'' for homotopy types. For example, above I was
  sketching how one might do homotopy theory using the ``model
  category'' of \(\omega\)-groupoids. Other model categories include
  gadgets like Kan complexes, CW complexes, simplicial complexes, and so
  on.
\item
  Chapter 6, ``Modern foundations for stable homotopy theory'', by A. D.
  Elmendorf, I. Kriz, M. Mandell and J. P. May describes a very nice
  approach to spectra. Loosely speaking, we can think of a spectrum as a
  \(\mathbb{Z}\)-groupoid, where \(\mathbb{Z}\) denotes the integers. In
  other words, in addition to \(j\)-morphisms for all natural numbers
  \(j\), we also have \(j\)-morphisms for negative \(j\)! This may seem
  bizarre, but it's a lot like how in homology theory one is interested
  in chain complexes that extend in both the positive and negative
  directions. In fact, we can think of a chain complex as a very special
  sort of \(\mathbb{Z}\)-groupoid or spectrum. The study of spectra is
  called stable homotopy theory.
\item
  Chapter 13, ``Stable homotopy and iterated loop spaces'', by G.
  Carlsson and R. J. Milgram, is packed with handy information about
  stable homotopy theory.
\item
  Chapter 21, ``Classifying spaces of compact Lie groups and finite loop
  spaces'', by D. Notbohm, is a good source of heavy-duty information on
  classifying spaces of your favorite Lie groups. To study vector
  bundles and the like one really needs to become comfortable with
  classifying spaces, and I'm finally doing this, and I hope eventually
  I'll be comfortable enough with them to really understand all these
  results.
\end{itemize}

There is a lot more, but I will stop here.

\hypertarget{week99}{%
\section{March 15, 1997}\label{week99}}

Life here at the Center for Gravitational Physics and Geometry is
tremendously exciting. In two weeks I have to return to U. C. Riverside
and my mundane life as a teacher of calculus, but right now I'm still
living it up. I'm working with Ashtekar, Corichi, and Krasnov on
computing the entropy of black holes using the loop representation of
quantum gravity, and also I'm talking to lots of people about an
interesting \(4\)-dimensional formulation of the loop representation in
terms of ``spin foams'' --- roughly speaking, soap-bubble-like
structures with faces labelled by spins.

Here are some papers I've come across while here:

\begin{enumerate}
\def\labelenumi{\arabic{enumi})}
\tightlist
\item
  Lee Smolin, ``The future of spin networks'', in \emph{The Geometric
  Universe: Science, Geometry, and the Work of Roger Penrose}, eds.~S.
  Hugget, Paul Tod, and L.J. Mason, Oxford U.\ Press,  Oxford, 1998. Also
  available as
  \href{https://arxiv.org/abs/gr-qc/9702030}{\texttt{gr-qc/9702030}}.
\end{enumerate}

I've spoken a lot about spin networks here on This Week's Finds. They
were first invented by Penrose as a radical alternative to the usual way
of thinking of space as a smooth manifold. For him, they were purely
discrete, purely combinatorial structures: graphs with edges labelled by
``spins'' \(j = 0, 1/2, 1, 3/2, \ldots\), and with three edges meeting
at each vertex. He showed how when these spin networks become
sufficiently large and complicated, they begin in certain ways to mimic
ordinary 3-dimensional Euclidean space. Interestingly, he never got
around to publishing his original paper on the subject:

\begin{enumerate}
\def\labelenumi{\arabic{enumi})}
\setcounter{enumi}{1}
\tightlist
\item
  Roger Penrose, ``Theory of quantized directions'', unpublished
  manuscript, available at \href{https://math.ucr.edu/home/baez/penrose/}{\texttt{https://math.ucr.edu/home/baez/penrose/}}.
\end{enumerate}
\noindent
Anyway, Smolin's paper is a kind of tribute to Penrose, and it traces
the curiously twisting history of spin networks from their origin up to
the present day, where they play a major role in topological quantum
field theory and the loop representation --- now more appropriately
called the spin network representation! --- of quantum gravity. (See
\protect\hyperlink{week55}{``Week 55''} for more on spin networks.)

Note however that the title of the paper refers to the \emph{future} of
spin networks. Smolin argues that spin networks are a major clue about
the future of physics, and he paints a picture of what this future might
be... which I urge you to look at.

For more on this, try:

\begin{enumerate}
\def\labelenumi{\arabic{enumi})}
\setcounter{enumi}{2}
\tightlist
\item
  Fotini Markopoulou and Lee Smolin, ``Causal evolution of spin
  networks'', available as
  \href{http://arxiv.org/abs/gr-qc/9702025}{\texttt{gr-qc/9702025}}.
\end{enumerate}
\noindent
Fotini Markopoulou is a student of Chris Isham at Imperial College, but
now she's visiting the CGPG and working with Lee Smolin on spin
networks. In this paper they describe some theories in which spin
networks evolve in time in discrete steps. The evolution is ``local'' in
the sense that in a given step, any vertex of the spin network changes
in a way that only depends on its immediate neighbors --- vertices
connected to it by an edge. It is also ``causal'' in the sense that
history of spin network evolving according to their rules gives a causal
set, i.e.~a set equipped with a partial ordering which we think of as
saying which points come ``before'' which other points. This ties their
work to the work of Rafael Sorkin on causal sets, e.g.:

\begin{enumerate}
\def\labelenumi{\arabic{enumi})}
\setcounter{enumi}{3}
\tightlist
\item
  Luca Bombelli, Joohan Lee, David Meyer and Rafael D. Sorkin,
  ``Space-time as a causal set'', \emph{Phys. Rev.~Lett.} \textbf{59}
  (1987), 521.
\end{enumerate}

Unlike the related work of Reisenberger and Rovelli (see
\protect\hyperlink{week96}{``Week 96''}), Markopolou and Smolin do not
attempt to ``derive'' their rules from general relativity by standard
quantization techniques. Instead, they hope that some theory of the sort
they consider will approximate general relativity in the large-scale
limit. To check this will require some new techniques akin to the
``renormalization group'' approach to studying the large-scale limits of
statistical mechanical systems defined on a lattice. This is a bit
daunting, but it seems likely that no matter how one proceeds to pursue
a spin-network-based theory of quantum gravity, one will need to develop
such techniques at some point.

\begin{center}\rule{0.5\linewidth}{0.5pt}\end{center}

\hypertarget{week99_tale}{
Now I'd like to switch gears and return to...}

\begin{center}
THE TALE OF \(n\)-CATEGORIES!
\end{center}

Recall that in our last episode, in \protect\hyperlink{week92_tale}{``Week
92''}, we had worked our way up to \(2\)-categories, and we were
beginning to see what they had to do with \(2\)-dimensional physics and
topology. I described how to get monads from adjunctions, and what this
has to do with matrix multiplication, Yang--Mills theory, and the 4-color
theorem.

Next week I want to get serious and start talking about \(n\)-categories
for arbitrary \(n\). One reason is that at the end of this month there's
a conference on \(n\)-categories and physics that I want to report on:

\begin{enumerate}
\def\labelenumi{\arabic{enumi})}
\setcounter{enumi}{4}
\tightlist
\item
  \emph{Workshop on Higher Category Theory and Physics}, March 28-30,
  1997, Northwestern University, Evanston, Illinois. Organized by Ezra
  Getzler and Mikhail Kapranov.
\end{enumerate}
\noindent
But before doing this, I want to say a bit about what category theory
has to do with quantum mechanics!

First remember the big picture: \(n\)-category theory is a language to
talk about processes that turn processes into other processes. Roughly
speaking, an \(n\)-category is some sort of structure with objects,
morphisms between objects, \(2\)-morphisms between morphisms, and so on
up to \(n\)-morphisms. A \(0\)-category is just a set, with its objects usually
being called ``elements''. Things get trickier as \(n\) increases. For a
precise definition of \(n\)-categories for \(n = 1\) and \(2\), see
\protect\hyperlink{week73}{``Week 73''} and
\protect\hyperlink{week80}{``Week 80''}, respectively.

Most familiar mathematical gadgets are sets equipped with extra bells
and whistles: groups, vector spaces, Hilbert spaces, and so on all have
underlying sets. This is why set theory plays an important role in
mathematics. However, we can also consider fancier gadgets that are
\emph{categories} equipped with extra bells and whistles. Some of the
most interesting examples are just ``categorifications'' of well-known
gadgets.

For example, a ``monoid'' is a simple gadget, just a set equipped with
an associative product and multiplicative identity. An example we all
know and love is the complex numbers: the product is just ordinary
multiplication, and the multiplicative identity is the number \(1\).

We may categorify the notion of ``monoid'' and define a ``monoidal
category'' to be a \emph{category} equipped with an associative product
and multiplicative identity. I gave the precise definition back in
\protect\hyperlink{week89}{``Week 89''}; the point here is that while
they may sound scary, monoidal categories are actually very familiar.
For example, the category of Hilbert spaces is a monoidal category where
the product of Hilbert spaces is the tensor product and the
multiplicative identity is \(\mathbb{C}\), the complex numbers.

If one systematically studies categorification one discovers an amazing
fact: many deep-sounding results in mathematics are just
categorifications of stuff we all learned in high school. There is a
good reason for this, I believe. All along, mathematicians have been
unwittingly ``decategorifying'' mathematics by pretending that
categories are just sets. We ``decategorify'' a category by forgetting
about the morphisms and pretending that isomorphic objects are equal. We
are left with a mere set: the set of isomorphism classes of objects.

I gave an example in \protect\hyperlink{week73}{``Week 73''}. There is a
category FinSet whose objects are finite sets and whose morphisms are
functions. If we decategorify this, we get the set of natural numbers!
Why? Well, two finite sets are isomorphic if they have the same number
of elements. ``Counting'' is thus the primordial example of
decategorification.

I like to think of it in terms of the following fairy tale. Long ago, if
you were a shepherd and wanted to see if two finite sets of sheep were
isomorphic, the most obvious way would be to look for an isomorphism. In
other words, you would try to match each sheep in herd \(A\) with a
sheep in herd \(B\). But one day, along came a shepherd who invented
decategorification. This person realized you could take each set and
``count'' it, setting up an isomorphism between it and some set of
``numbers'', which were nonsense words like ``one, two, three,
four,...'' specially designed for this purpose. By comparing the
resulting numbers, you could see if two herds were isomorphic without
explicitly establishing an isomorphism!

According to this fairy tale, decategorification started out as the
ultimate stroke of mathematical genius. Only later did it become a
matter of dumb habit, which we are now struggling to overcome through
the process of ``categorification''.

Okay, so what does this have to do with quantum mechanics?

Well, a Hilbert space is a set with extra bells and whistles, so maybe
there is some gadget called a ``2-Hilbert space'' which is a
\emph{category} with analogous extra bells and whistles. And maybe if we
figure out this notion we will learn something about quantum mechanics.

Actually the notion of 2-Hilbert space didn't arise in this
simple-minded way. It arose in some work by Daniel Freed on topological
quantum field theory:

\begin{enumerate}
\def\labelenumi{\arabic{enumi})}
\setcounter{enumi}{4}
\tightlist
\item
  Daniel Freed, ``Higher algebraic structures and quantization'',
  \emph{Commun. Math. Phys.} \textbf{159} (1994), 343--398.  Also
  available as
  \href{https://arxiv.org/ps/hep-th/9212115}{\texttt{hep-th/9212115}};
  see also \protect\hyperlink{week48}{``Week 48''}.
\end{enumerate}
\noindent
Later, Louis Crane adopted this notion as part of his program to reduce
quantum gravity to \(n\)-category theory:

\begin{enumerate}
\def\labelenumi{\arabic{enumi})}
\setcounter{enumi}{5}
\tightlist
\item
  Louis Crane: ``Clock and category: is quantum gravity algebraic?'',
  \emph{Jour. Math. Phys.} \textbf{36} (1995), 6180--6193.  Also
  available as
  \href{https://arxiv.org/ps/gr-qc/9504038}{\texttt{gr-qc/9504038}}.
\end{enumerate}
\noindent
These papers made is clear that 2-Hilbert spaces are interesting things
and that one should go further and think about ``\(n\)-Hilbert spaces''
for higher \(n\), too. However, neither of them gave a precise
definition of 2-Hilbert space, so at some point I decided to do this. It
took a while for me to learn enough category theory, but eventually I
wrote something about it:

\begin{enumerate}
\def\labelenumi{\arabic{enumi})}
\setcounter{enumi}{6}
\tightlist
\item
  John Baez, ``Higher-dimensional algebra II: 2-Hilbert spaces'', 
  \textsl{Adv.\ Math.\ }\textbf{127} (1997), 125--189. Also available as
  \href{https://arxiv.org/abs/q-alg/9609018}{\texttt{q-alg/9609018}}.
\end{enumerate}
\noindent
To understand this requires a little category theory, so let me explain
the basic ideas here.

I'll concentrate on finite-dimensional Hilbert spaces, since the
infinite-dimensional case introduces extra complications. To define
2-Hilbert spaces we need to start by categorifying the various
ingredients in the definition of Hilbert space. These are:

\begin{enumerate}
\def\labelenumi{\arabic{enumi}.}
\tightlist
\item
  the zero element,
\item
  addition,
\item
  subtraction,
\item
  scalar multiplication, and
\item
  the inner product.
\end{enumerate}
\noindent
The first four have well-known categorical analogs. The fifth one, which
is really the essence of a Hilbert space, may seem a bit more mysterious
at first, but as we shall see, it's really the key to the whole
business.

\begin{enumerate}
\def\labelenumi{\arabic{enumi})}
\item
  The analog of the zero vector is a `zero object'. A zero object in a
  category is an object that is both initial and terminal. That is,
  there is exactly one morphism from it to any object, and exactly one
  morphism to it from any object. Consider for example the category
  \(\mathsf{Hilb}\) having finite-dimensional Hilbert spaces as objects,
  and linear maps between them as morphisms. In \(\mathsf{Hilb}\), any
  zero-dimensional Hilbert space is a zero object.

  Note: there isn't really a unique zero object in the ``strict'' sense
  of the term. Instead, any two zero objects are canonically isomorphic.
  The reason is that if you have two zero objects, say \(0\) and \(0'\),
  there is a unique morphism \(f\colon 0\to 0'\) and a unique morphism
  \(g\colon 0'\to 0\). These morphisms are inverses of each other so
  they are isomorphisms. Why are they inverses? Well,
  \(fg\colon 0\to 0'\) must be the identity morphism
  \(1_0\colon 0 \to 0\), because there is only one morphism from \(0\)
  to \(0\)! Similarly, \(gf\) is the identity on \(0'\). (Note that I am
  using category theorist's notation, where the composite of
  \(f\colon x\to y\) and \(g\colon y\to z\) is denoted
  \(fg\colon x\to z\).)

  This is typical in category theory. We don't expect stuff to be
  unique; it should only be unique up to a canonical isomorphism.
\item
  The analog of adding two vectors is forming the ``coproduct'' of two
  objects. Coproducts are just a fancy way of talking about direct sums.
  Any decent quantum mechanic knows about the direct sum of Hilbert
  spaces. But in fact, we can define this notion very generally in any
  category, where it goes under the name of a ``coproduct''. (I give the
  definition below; if I gave it here it would scare people away.) As
  with zero objects, coproducts are typically not unique, but they are
  always unique up to canonical isomorphism, which is what matters. It's
  a good little exercise to show this.
\item
  The analog of subtracting vectors is forming the ``cokernel'' of a
  morphism \(f\colon x\to y\). If \(x\) and \(y\) are Hilbert spaces,
  the cokernel of \(f\) is just the orthogonal complement of the range
  of \(f\). In other words, for Hilbert spaces we have ``direct
  differences'' as well as direct sums. However, the notion of cokernel
  makes sense in any category with a zero object. I won't burden you
  with the precise definition here.
\end{enumerate}

An important difference between zero, addition and subtraction and their
categorical analogs is that these operations represent extra
\emph{structure} on a set, while having a zero object, coproducts of two
objects, or cokernels of morphisms is merely a \emph{property} of a
category. Thus these concepts are in some sense more intrinsic to
categories than to sets. On the other hand, we've seen one pays a price
for this: while the zero element, sums, and differences are unique in a
Hilbert space, the zero object, coproducts, and cokernels are typically
unique only up to canonical isomorphism.

\begin{enumerate}
\def\labelenumi{\arabic{enumi})}
\setcounter{enumi}{3}
\tightlist
\item
  The analog of multiplying a vector by a complex number is tensoring an
  object by a Hilbert space. Besides its additive properties (zero
  object, binary coproducts, and cokernels), \(\mathsf{Hilb}\) is also a
  monoidal category: we can multiply Hilbert spaces by tensoring them,
  and there is a multiplicative identity, namely the complex numbers
  \(\mathbb{C}\). In fact, \(\mathsf{Hilb}\) is a ``ring category'', as
  defined by Laplaza and Kelly.
\end{enumerate}

We expect \(\mathsf{Hilb}\) to play a role in 2-Hilbert space theory
analogous to the role played by the ring \(\mathbb{C}\) of complex
numbers in Hilbert space theory. Thus we expect 2-Hilbert spaces to be
``module categories'' over \(\mathsf{Hilb}\), as defined by Kapranov and
Voevodsky.

An important part of our philosophy here is that \(\mathbb{C}\) is the
primordial Hilbert space: the simplest one, upon which the rest are
modelled. By analogy, we expect \(\mathsf{Hilb}\) to be the primordial
2-Hilbert space. This is part of a general pattern pervading
higher-dimensional algebra; for example, there is a sense in which the
\((n+1)\)-category of all (small) \(n\)-categories, \(n\mathsf{Cat}\),
is the primordial \((n+1)\)-category. The real significance of this
pattern remains mysterious.

\begin{enumerate}
\def\labelenumi{\arabic{enumi})}
\setcounter{enumi}{4}
\tightlist
\item
  Finally, what is the categorification of the inner product in a
  Hilbert space? It's the `\(\operatorname{Hom}\) functor'! The inner
  product in a Hilbert space eats two vectors \(v\) and \(w\) and spits
  out a complex number \[\langle v,w \rangle\] Similarly, given two
  objects \(v\) and \(w\) in a category, the \(\operatorname{Hom}\)
  functor gives a \emph{set} \[\operatorname{Hom}(x,y)\] namely the set
  of morphisms from \(x\) to \(y\). Note that the inner product
  \(\langle v,w \rangle\) is linear in \(w\) and conjugate-linear in
  \(y\), and similarly, the \(\operatorname{Hom}\) functor
  \(\operatorname{Hom}(x,y)\) is covariant in \(y\) and contravariant in
  \(x\). This hints at the category theory secretly underlying quantum
  mechanics. In quantum theory the inner product \(\langle v,w \rangle\)
  represents the \emph{amplitude} to pass from \(v\) to \(w\), while in
  category theory \(\operatorname{Hom}(x,y)\) is the \emph{set} of ways
  to get from \(x\) to \(y\). In Feynman path integrals, we do an
  integral over the set of ways to get from here to there, and get a
  number, the amplitude to get from here to there. So when physicists do
  Feynman path integration --- just like a shepherd counting sheep ---
  they are engaged in a process of decategorification!
\end{enumerate}

To understand this analogy better, note that any morphism
\(f\colon x\to y\) in \(\mathsf{Hilb}\) can be turned around or
``dualized'' to obtain a morphism \(f^*\colon y\to x\). This is usually
called the ``adjoint'' of \(f\), and it satisfies
\[\langle fv,w \rangle = \langle v,f^*w \rangle\] for all \(v\) in
\(x\), and \(w\) in \(y\). This ability to dualize morphisms is crucial
to quantum theory. For example, observables are represented by
self-adjoint morphisms, while symmetries are represented by unitary
morphisms, whose adjoint equals their inverse.

However, it should now be clear --- at least to the categorically minded
--- that this sort of adjoint is just a decategorified version of the
``adjoint functors'' so important in category theory. As I explained in
\protect\hyperlink{week79}{``Week 79''}, a functor
\(F^*\colon\mathcal{D}\to\mathcal{C}\) is a ``right adjoint'' of
\(F\colon\mathcal{C}\to\mathcal{D}\) if there is, not an equation, but a
natural isomorphism
\[\operatorname{Hom}(Fc,d) \cong \operatorname{Hom}(c,F^*d)\] for all
objects \(c\) in \(\mathcal{C}\), and \(d \in \mathcal{D}\).

Anyway, in the paper I proceed to use these ideas to give a precise
definition of 2-Hilbert spaces, and then I prove all sorts of stuff
which I won't describe here.

Let me wrap up by explaining the definition of ``coproduct''. This is
one of those things they should teach all math grad students, but for
some reason they don't. It's a bit dry but it'll be good for you. A
coproduct of the objects \(x\) and \(y\) is an object \(x+y\) equipped
with morphisms \[i\colon x \to x+y\] and \[j\colon y \to x+y\] that is
universal with respect to this property. In other words, if we have any
\emph{other} object, say \(z\), and morphisms \[i'\colon x \to z\] and
\[j'\colon y \to z\] then there is a unique morphism
\(f\colon x+y \to z\) such that \[i' = if\] and \[j' = jf.\] This kind
of definition automatically implies that the coproduct is unique up to
canonical isomorphism. To understand this abstract nonsense, it's good
to check that the coproduct of sets or topological spaces is just their
disjoint union, while the coproduct of vector spaces or Hilbert spaces
is their direct sum.

To continue reading the ``Tale of \(n\)-Categories'', see
\protect\hyperlink{week100_tale}{``Week 100''}.

\hypertarget{week100}{%
\section{March 23, 1997}\label{week100}}

\hypertarget{week100_tale}{Pretty much ever since I started writing 
``This Week's Finds'' I've been
trying to get folks interested in \(n\)-categories and other aspects of
higher-dimensional algebra.} There is really an enormous world out there
that only becomes visible when you break out of ``linear thinking'' ---
the mental habits that go along with doing math by writing strings of
symbols in a line. For example, when we write things in a line, the sums
\(a+b\) and \(b+a\) look very different. Then we introduce a profound
and mysterious equation, the ``commutative law'': \[a + b = b + a\]
which says that actually they are the same. But in real life, we prove
this equation using higher-dimensional reasoning:
\[a+b = {}^{\mbox{\normalfont $a$}}+{}_{\mbox{\normalfont $b$}} = \underset{{\mbox{\normalfont $b$}}}{\overset{{\mbox{\normalfont $a$}}}{+}} = {}_{{\mbox{\normalfont $b$}}}+{}^{{\mbox{\normalfont $a$}}} = b+a\]
If this seems silly, think about explaining to a kid what \(9+17\)
means, and how you could prove that \(9+17 = 17+9\). You might take a
pile of 9 rocks and set it to the left of a pile of 17 rocks, and say
``this is 9+17 rocks''. Alternatively, you might put the pile of 9 rocks
to the right of the pile of 17 rocks, and say ``this is 17+9 rocks''.
Thus to prove that \(9+17=17+9\), you would simply need to \emph{switch}
the two piles by moving one around the other.

This is all very simple. Historically, however, it took people a long time to
really understand. It's one of those things that's too simple to take
seriously until it turns out to have complicated ramifications. Now it
goes by the name of the ``Eckmann--Hilton theorem'', which says that ``a
monoid object in the category of monoids is a commutative monoid''. You
practically need a PhD in math to understand \emph{that}! However, lest
you think that Eckmann and Hilton were merely dressing up the obvious in
fancy jargon, it's important to note that what they did was to figure
out a \emph{framework} that turns the above ``picture proof'' that
\(a+b = b+a\) into an actual rigorous proof! This is one of the goals of
higher-dimensional algebra.

The above proof that \(a+b = b+a\) uses \(2\)-dimensional space, but if
you really think about it also uses a 3rd dimension, namely time: the
time that passes as you move ``\(a\)'' around ``\(b\)''. If we draw this
3rd dimension as space rather than time we can visualize the process of
moving \(a\) around \(b\) as follows: \[
  \begin{tikzpicture}
    \begin{knot}[clip width=7]
      \strand[thick] (1,0) to (0,-2);
      \strand[thick] (0,0) to (1,-2);
    \end{knot}
    \node[label=above:{$a$}] at (0,0) {};
    \node[label=above:{$b$}] at (1,0) {};
    \node[label=below:{$a$}] at (1,-2) {};
    \node[label=below:{$b$}] at (0,-1.9) {};
  \end{tikzpicture}
\] This picture is an example of what mathematicians call a ``braid''.
This particular one is a boring little braid with only two strands and
one place where the two strands cross. It illustrates another major idea
behind higher-dimensional algebra: equations are best thought of as
summarizing ``processes'' (or technically, ``isomorphisms''). The
equation \(a+b = b+a\) is a summary of the process of switching \(a\)
and \(b\). There is more information in the process than in the mere
equation \(a+b = b+a\), because in fact there are two \emph{different}
ways to switch \(a\) and \(b\): the above way and \[
  \begin{tikzpicture}
    \begin{knot}[clip width=7]
      \strand[thick] (0,0) to (1,-2);
      \strand[thick] (1,0) to (0,-2);
    \end{knot}
    \node[label=above:{$a$}] at (0,0) {};
    \node[label=above:{$b$}] at (1,0) {};
    \node[label=below:{$a$}] at (1,-2) {};
    \node[label=below:{$b$}] at (0,-1.9) {};
  \end{tikzpicture}
\] If one has a bunch of objects one can switch them around in a lot of
ways, getting lots of different braids.

In fact, the mathematics of braids, and related things like knots, is
crucially important for understanding quantum gravity in
\(3\)-dimensional spacetime. Spacetime is really \(4\)-dimensional, of
course, but quantum gravity in \(4\)-dimensional spacetime is awfully
difficult, so in the late 1980s people got serious about studying
\(3\)-dimensional quantum gravity as a kind of warmup exercise. It
turned out that the math required was closely related to some mysterious
new mathematics related to knots and ``braidings''. At first this must
sound bizarre: a deep relationship between knots and \(3\)-dimensional
quantum gravity! However, after you fight your way through the
sophisticated mathematical physics that's involved, it becomes clear why
they are related: both rely crucially on ``3-dimensional algebra'', the
algebra describing how you can move things around in \(3\)-dimensional
spacetime.

However, there is more to the story, because knot theory also seems
deeply related to \emph{\(4\)-dimensional} quantum gravity. Here the knots
arise as ``flux tubes of area'' living in \(3\)-dimensional space at a
given time. Recent work on quantum gravity suggests that as time passes
these knots (or more generally, ``spin networks'') move around and
change topology as time passes.

To really understand this, we probably need to understand
``4-dimensional algebra''. Unfortunately, not enough is known about
4-dimensional algebra. The problem is that we don't know much about
4-categories! To do \(n\)-dimensional algebra in a really nice way, you
need to know about \(n\)-categories. Roughly speaking, an \(n\)-category
is an algebraic structure that has a bunch of things called ``objects'',
a bunch of things called ``morphisms'' that go between objects, and
similarly \(2\)-morphisms going between morphisms, \(3\)-morphisms going
between\(2\)-morphisms, and so on up to the number \(n\). You can think of the
objects as ``things'' of whatever sort you like, the morphisms as
processes going from one thing to another, the \(2\)-morphisms as
meta-processes going from one process to another, and so on. Depending
on how you play the \(n\)-category game, there are either no morphisms
after level \(n\), or only simple and bland ones playing the role of
``equations''. The idea is that in the world of \(n\)-categories, one
keeps track of things, processes, meta-processes, and so on to the
\(n\)th level, but after that one calls it quits and uses equations.

So what is the definition of \(4\)-categories? Well, Eilenberg and Mac
Lane defined \(1\)-categories, or simply ``categories'', in a paper that
was published in 1945:

\begin{enumerate}
\def\labelenumi{\arabic{enumi})}
\tightlist
\item
  Samuel Eilenberg and Saunders Mac Lane, ``General theory of natural
  equivalences'', \emph{Trans. Amer. Math. Soc.} \textbf{58} (1945),
  231--294.
\end{enumerate}
\noindent
B\'enabou defined \(2\)-categories --- though actually he called them
``bicategories'' --- in a 1967 paper:

\begin{enumerate}
\def\labelenumi{\arabic{enumi})}
\setcounter{enumi}{1}
\tightlist
\item
  Jean B\'enabou, \emph{Introduction to Bicategories}, Springer Lecture
  Notes in Mathematics \textbf{47}, New York, 1967, pp.~1--77.
\end{enumerate}
\noindent
Gordon, Power, and Street defined \(3\)-categories --- or actually
``tricategories'' --- in a paper that came out in 1995:

\begin{enumerate}
\def\labelenumi{\arabic{enumi})}
\setcounter{enumi}{2}
\tightlist
\item
  Rorbert Gordon, A. John Power, and Ross Street, ``Coherence for
  tricategories'', \emph{Memoirs Amer. Math. Soc.} \textbf{117} (1995).
\end{enumerate}
\noindent
This step took a long time in part because it took a long time for
people to understand deeply where \emph{braidings} fit into the picture.

But what about \(4\)-categories and higher \(n\)? Well, the history is
complicated and I won't get it right, but let me say a bit anyway. First
of all, there are some things called ``strict \(n\)-categories'' that
people have known how to define for arbitrarily high \(n\) for quite a
while. In fact, people know how to go up to infinity and define ``strict
\(\omega\)-categories''; see for example:

\begin{enumerate}
\def\labelenumi{\arabic{enumi})}
\setcounter{enumi}{3}
\tightlist
\item
  Sjoerd E. Crans, \emph{On combinatorial models for higher dimensional
  homotopies}, Ph.D.~thesis, University of Utrecht, Utrecht, 1991.
\end{enumerate}
\noindent
Strict \(n\)-categories are quite interesting and important, but I'm
mainly mentioning them here to emphasize that they are \emph{not} what
I'm talking about. People sometimes often call strict \(n\)-categories
simply ``\(n\)-categories'', and call the more general \(n\)-categories
I'm talking about above ``weak \(n\)-categories''. However, I think the
weak \(n\)-categories will eventually be called simply
``\(n\)-categories'', because they are far more interesting and
important than the strict ones. Anyway, that's what I'm doing here.

Secondly, when you define \(n\)-categories you have to make some choice
about the ``shapes'' of your \(j\)-morphisms. In general they should be
some \(j\)-dimensional things, but they could be simplices, or cubes, or
other shapes. In some ways the simplest shapes are ``globes'', a
\(j\)-dimensional globe being a \(j\)-dimensional ball with its boundary
divided into two hemispheres, the ``inface'' and ``outface'', which are
themselves \((j-1)\)-dimensional globes. This corresponds to a picture
where each ``process'' has one input and one output, which are
themselves processes having the same input and output. The definitions
of category, bicategory, and tricategory work this way. In fact, Ross
Street came up with a very nice definition of \(n\)-categories for all
\(n\) using simplices in 1987:

\begin{enumerate}
\def\labelenumi{\arabic{enumi})}
\setcounter{enumi}{4}
\tightlist
\item
  Ross Street, ``The algebra of oriented simplexes'', \emph{Jour. Pure
  Appl. Alg.} \textbf{49} (1987), 283--335.
\end{enumerate}
\noindent
Since then, however, he and his students and collaborators seem to have
been working to translate this definition into the ``globular''
formalism... while also making some other important adjustments too
technical to discuss here. In particular, Dominic Verity and Todd
Trimble have done a lot of work on getting the definition of
\(n\)-category worked out, and a while ago I learned that Trimble came
up with a definition of ``tetracategory'' (or what I'm calling simply
``4-category'') in August of 1995. I don't think this has been
published, however.

James Dolan came to U.\ C.\ Riverside in the fall of 1993, and ever since
then, he and I have been talking about \(n\)-categories and their role in
physics. Most of the category theory I know, I learned in this process.
It soon became clear that we needed a nice definition of \(n\)-category
for all \(n\) in order to turn our hopes and dreams into theorems. After
a while we started working pretty hard on this. His job was to come up
with all the bright ideas, and mine was to get him to explain them, to
try to poke holes in them, and to figure out rigorous proofs of all the
things that were so obvious to him that he couldn't figure out how (or
why) to prove them. We sent a summarized version of our definition to
Ross Street at the end of 1995:

\begin{enumerate}
\def\labelenumi{\arabic{enumi})}
\setcounter{enumi}{5}
\tightlist
\item
  John Baez and James Dolan, ``\(n\)-Categories --- sketch of a definition'',
  letter to Ross Street, Nov.~29, 1995, available at
  \href{http://math.ucr.edu/home/baez/ncat.def.html}{\texttt{http://math.ucr.edu/home/baez/ncat.def.html}}.
\end{enumerate}
\noindent
and then for a year I worked on trying to write up a longer, clearer
version, while all the meantime Dolan kept coming up with new ways of
looking at everything. I finished in February of this year:

\begin{enumerate}
\def\labelenumi{\arabic{enumi})}
\setcounter{enumi}{6}
\tightlist
\item
  John Baez and James Dolan, ``Higher-dimensional algebra III:
  \(n\)-Categories and the algebra of opetopes'', to appear in
  \emph{Adv. Math.}  Also available as
  \href{https://arxiv.org/ps/q-alg/9702014}{\texttt{q-alg/9702014}}.
\end{enumerate}

The key feature of this definition is that it uses ``\(j\)-dimensional
opetopes'' as the shapes for \(j\)-morphisms. These shapes are very
handy because the \((j+1)\)-dimensional opetopes describe all the legal
ways of sticking together a bunch of \(j\)-dimensional opetopes to form
another \(j\)-dimensional opetope! They are related to the theory of
``operads'', which is part of the reason for their name. (By the way,
the first two syllables are pronounced exactly as in ``operation''.)

In the meantime, Michael Makkai and John Power had begun work using our
definition. Also, other definitions of ``\(n\)-category'' have appeared
on the scene! Zouhair Tamsamani came up with one in terms of
``multi-simplicial sets'':

\begin{enumerate}
\def\labelenumi{\arabic{enumi})}
\setcounter{enumi}{7}
\tightlist
\item
  Zouhair Tamsamani, \emph{Sur des notions de \(\infty\)-categorie et
  \(\infty\)-groupoide non-strictes via des ensembles
  multi-simpliciaux}, Ph.D.~thesis, Universite Paul Sabatier, Toulouse,
  France, 1995.
\end{enumerate}

Michael Batanin also has a definition of \(\omega\)-categories, of the
``globular'' sort:

\begin{enumerate}
\def\labelenumi{\arabic{enumi})}
\setcounter{enumi}{8}
\tightlist
\item
  Michael A. Batanin, ``On the definition of weak \(\omega\)-category'',
  \emph{Macquarie Mathematics Report} number \textbf{96/207}.
\end{enumerate}
\noindent
Now the fun will begin! These different definitions of (weak)
\(n\)-category should be equivalent, albeit in a rather subtle sense, so
we should check to see if they really are. Also, we need to develop many
more tools for working with \(n\)-categories. Then we can really start
using them as a tool.

When I started writing this I thought I was going to explain the
definition that Dolan and I came up with. Now I'm too tired! It takes a
while to explain, so I think I'll stop here and save that for some other
week or weeks. Perhaps I'll mix it in with my report on the Workshop on
Higher Category Theory and Physics, which is taking place next weekend
at Northwestern University.

This is the end of the ``Tale of \(n\)-Categories''. If you want more,
try \href{https://arxiv.org/abs/q-alg/9705009}{``An Introduction to
\(n\)-Categories''}, or else read the above papers.

\end{document}